\documentclass[12pt]{amsbook}
\usepackage{amssymb}
\usepackage[all]{xy}

\usepackage[colorlinks=true,linkcolor=magenta,citecolor=magenta]{hyperref}
\makeindex

\newtheorem{theorem}{Theorem}[chapter]
\newtheorem{advertisement}[theorem]{Advertisement}
\newtheorem{advice}[theorem]{Advice}
\newtheorem{algorithm}[theorem]{Algorithm}
\newtheorem{answer}[theorem]{Answer}
\newtheorem{cat}[theorem]{Cat}
\newtheorem{claim}[theorem]{Claim}
\newtheorem{conclusion}[theorem]{Conclusion}
\newtheorem{conjecture}[theorem]{Conjecture}
\newtheorem{convention}[theorem]{Convention}
\newtheorem{definition}[theorem]{Definition}
\newtheorem{dilemma}[theorem]{Dilemma}
\newtheorem{examples}[theorem]{Examples}
\newtheorem{exercise}[theorem]{Exercise}

\newtheorem{fact}[theorem]{Fact}

\newtheorem{interpretation}[theorem]{Interpretation}
\newtheorem{method}[theorem]{Method}
\newtheorem{principle}[theorem]{Principle}
\newtheorem{problem}[theorem]{Problem}
\newtheorem{proposition}[theorem]{Proposition}
\newtheorem{puzzle}[theorem]{Puzzle}
\newtheorem{question}[theorem]{Question}
\newtheorem{remark}[theorem]{Remark}
\newtheorem{thought}[theorem]{Thought}
\newtheorem{warning}[theorem]{Warning}

\begin{document}

\title{Calculus and applications}

\author{Teo Banica}
\address{Department of Mathematics, University of Cergy-Pontoise, F-95000 Cergy-Pontoise, France. {\tt teo.banica@gmail.com}}


\begin{abstract}
This is an introduction to calculus, and its applications to basic questions from physics. We first discuss the theory of functions $f:\mathbb R\to\mathbb R$, with the notion of continuity, and the construction of the derivative $f'(x)$ and of the integral $\int_a^bf(x)dx$. Then we investigate the case of the complex functions $f:\mathbb C\to\mathbb C$, and notably the holomorphic functions, and harmonic functions. Then, we discuss the multivariable functions, $f:\mathbb R^N\to\mathbb R^M$ or $f:\mathbb R^N\to\mathbb C^M$ or $f:\mathbb C^N\to\mathbb C^M$, with general theory, integration results, maximization questions, and basic applications to physics.
\end{abstract}

\maketitle

\chapter*{Preface}

Understanding what happens in the real life surrounding us, in phenomena involving physics, chemistry, biology and so on, is not an easy task. What we can do as humans is to come up with some machinery, and perform measurements, recording quantities such as length, volume, temperature, pressure and so on, and then see how these quantities, called ``variables'', and denoted $x,y,z,\ldots$ depend on each other, and change in time. 

\bigskip

Calculus is the study of the correspondences $x\to y$ between such variables. Such correspondences are called ``functions'', and are denoted $y=f(x)$, with $f$ standing for the abstract machinery, or mathematical formula, producing $y$ out of $x$.

\bigskip

The basics of calculus were developed by Newton, Leibnitz and others, a long time ago. The idea is very simple. The simplest functions $f:\mathbb R\to\mathbb R$ are the linear ones, $f(x)=a+bx$ with $a,b\in\mathbb R$, but of course not any function is linear. Miraculously, however, most functions $f:\mathbb R\to\mathbb R$ are ``locally linear'', in the sense that around any given point $c\in\mathbb R$, we have a formula of type $f(c+x)\simeq a+bx$, for $x$ small. Why? Obviously, $a\in\mathbb R$ can only be the value of our function at that point, $a=f(c)$. As for the number $b\in\mathbb R$, this can be taken to be the rate of change of $f$ around that point, called derivative of the function at that point, and denoted $b=f'(c)$. 

\bigskip

So, this was the main idea of calculus, ``functions are locally linear''. This idea applies as well to more complicated functions, such as the ``multivariable'' ones $f:\mathbb R^N\to\mathbb R^M$, relating vector variables $x\in\mathbb R^N$ to vector variables $y\in\mathbb R^M$, with the linear approximation formula $f(c+x)\simeq a+bx$ needing this time as parameters a vector $a=f(c)\in\mathbb R^M$, and a linear map, or beast called rectangular matrix, $b=f'(c)\in M_{M\times N}(\mathbb R)$.

\bigskip

Further ideas of calculus, which are more advanced, include the facts that: (1) the remainder $\varepsilon(x)$ given by $f(c+x)=a+bx+\varepsilon(x)$ can studied by using again derivatives, (2) in several variables, the geometric understanding of the derivatives $f'(c)\in M_{M\times N}(\mathbb R)$ is best done by using complex numbers, (3) in fact, the use of complex numbers is useful even for one-variable functions $f:\mathbb R\to\mathbb R$, and (4) in one variable at least, there is a magic relation between derivatives and weighted averages, called integrals and denoted $\int_a^bf(x)dx$, the idea being that ``the derivative of the integral is the function itself''.

\bigskip

Calculus can be learned from many places, with this being mostly a matter of taste. Personally as a student I read the books of Rudin \cite{ru1}, \cite{ru2}, and this was a very good investment, never had any trouble with calculus since, be that for research, or teaching. And these are still the books that I recommend to my students, although in the present modern age there are so many alternative resources, for having the basics learned.

\bigskip

The present book is an introduction to calculus, based on lecture notes from various classes that I taught at Cergy, and previously at Toulouse. The material inside claims of course no originality, basically going back to Newton, Leibnitz and others. But in what regards the presentation, there are a few ideas behind it, none of these claiming of course originality either, but their combination being something original, I hope:

\bigskip

(1) One complex variable comes before several real variables. This is perhaps not that standard, but as a quantum physicist, I just love complex numbers.

\bigskip

(2) Applications to probability everywhere, scattered throughout the book. Again, coming from experience with mathematics, physics, and science in general.

\bigskip

(3) Combinatorics, binomials and factorials all over the place, with joy. With this being a quite popular approach, who in mathematics does not love binomials.

\bigskip

(4) Applications to physics too, including even the hydrogen atom, at the end. In short, read this book, and you'll understand how hydrogen $_1{\rm H}$ works.

\bigskip

In the hope that you will like this book. The presentation will be quite quick, and for more on 1-variable calculus, you have here my book on functions \cite{ba1}, which normally comes first. If you would like to know more about several variables, you have here my book on measure and integration \cite{ba2}. And finally, in case you find the present book too mathematical and rigorous, go with my physics book \cite{ba3}, you won't be deceived.

\bigskip

As already mentioned, the present book is based on lecture notes from classes at Toulouse and Cergy, and I would like to thank my students. Many thanks go as well to my cats, for useful pieces of advice, often complementary to the pieces of advice of my colleagues, and for some help with the underlying PDE and physics.

\bigskip

\

{\em Cergy, January 2026}

\smallskip

{\em Teo Banica}\

\baselineskip=15.95pt
\tableofcontents
\baselineskip=14pt

\part{Basic calculus}

\ \vskip50mm

\begin{center}
{\em I've got to stand and fight

In this creation

Vanity I know

Can't guide I alone}
\end{center}

\chapter{Sequences, series}

\section*{1a. Binomials, factorials}

We denote by $\mathbb N$ the set of positive integers, $\mathbb N=\{0,1,2,3,\ldots\}$, with $\mathbb N$ standing for ``natural''. Quite often in computations we will need negative numbers too, and we denote by $\mathbb Z$ the set of all integers, $\mathbb Z=\{\ldots,-2,-1,0,1,2,\ldots\}$, with $\mathbb Z$ standing from ``zahlen'', which is German for ``numbers''. Finally, there are many questions in mathematics involving fractions, or quotients, which are called rational numbers:

\index{fractions}

\begin{definition}
The rational numbers are the quotients of type
$$r=\frac{a}{b}$$
with $a,b\in\mathbb Z$, and $b\neq0$, identified according to the usual rule for quotients, namely:
$$\frac{a}{b}=\frac{c}{d}\iff ad=bc$$
We denote the set of rational numbers by $\mathbb Q$, standing for ``quotients''.
\end{definition}

Observe that we have inclusions $\mathbb N\subset\mathbb Z\subset\mathbb Q$. The integers add and multiply according to the rules that you know well. As for the rational numbers, these add according to the usual rule for quotients, which is as follows, and death penalty for forgetting it: 
$$\frac{a}{b}+\frac{c}{d}=\frac{ad+bc}{bd}$$

Also, the rational numbers multiply according to the usual rule for quotients, namely:
$$\frac{a}{b}\cdot\frac{c}{d}=\frac{ac}{bd}$$

Beyond rationals, we have the real numbers, whose set is denoted $\mathbb R$, and which include beasts such as $\sqrt{3}=1.73205\ldots$ or $\pi=3.14159\ldots$ But more on these later. For the moment, let us see what can be done with integers, and their quotients. As a first theorem, solving a problem which often appears in real life, we have:

\index{binomial coefficient}
\index{factorials}

\begin{theorem}
The number of possibilities of choosing $k$ objects among $n$ objects is
$$\binom{n}{k}=\frac{n!}{k!(n-k)!}$$
called binomial number, where $n!=1\cdot2\cdot3\ldots(n-2)(n-1)n$, called ``factorial $n$''.
\end{theorem}

\begin{proof}
Imagine a set consisting of $n$ objects. We have $n$ possibilities for choosing our 1st object, then $n-1$ possibilities for choosing our 2nd object, out of the $n-1$ objects left, and so on up to $n-k+1$ possibilities for choosing our $k$-th object, out of the $n-k+1$ objects left. Since the possibilities multiply, the total number of choices is:
\begin{eqnarray*}
N
&=&n(n-1)\ldots(n-k+1)\\
&=&n(n-1)\ldots(n-k+1)\cdot\frac{(n-k)(n-k-1)\ldots2\cdot1}{(n-k)(n-k-1)\ldots 2\cdot1}\\
&=&\frac{n(n-1)\ldots2\cdot 1}{(n-k)(n-k-1)\ldots 2\cdot1}\\
&=&\frac{n!}{(n-k)!}
\end{eqnarray*}

But is this correct. Normally a mathematical theorem coming with mathematical proof is guaranteed to be $100\%$ correct, and if in addition the proof is truly clever, like the above proof was, with that fraction trick, the confidence rate jumps up to $200\%$. 

\medskip

This being said, never knows, so let us doublecheck, by taking for instance $n=3,k=2$. Here we have to choose 2 objects among 3 objects, and this is something easily done, because what we have to do is to dismiss one of the objects, and $N=3$ choices here, and keep the 2 objects left. Thus, we have $N=3$ choices. On the other hand our genius math computation gives $N=3!/1!=6$, which is obviously the wrong answer.

\medskip

So, where is the mistake? Thinking a bit, the number $N$ that we computed is in fact the number of possibilities of choosing $k$ ordered objects among $n$ objects. Thus, we must divide everything by the number $M$ of orderings of the $k$ objects that we chose:
$$\binom{n}{k}=\frac{N}{M}$$

In order to compute now the missing number $M$, imagine a set consisting of $k$ objects. There are $k$ choices for the object to be designated $\#1$, then $k-1$ choices for the object to be designated $\#2$, and so on up to 1 choice for the object to be designated $\#k$. We conclude that we have $M=k(k-1)\ldots 2\cdot 1=k!$, and so:
$$\binom{n}{k}=\frac{n!/(n-k)!}{k!}=\frac{n!}{k!(n-k)!}$$

And this is the correct answer, because, well, that is how things are. In case you doubt, at $n=3,k=2$ for instance we obtain $3!/2!1!=3$, which is correct.
\end{proof}

All this is quite interesting, and in addition to having some exciting mathematics going on, and more on this in a moment, we have as well some philosophical conclusions. Formulae can be right or wrong, and as the above shows, good-looking, formal mathematical proofs can be right or wrong too. So, what to do? Here is my advice:

\begin{advice}
Always doublecheck what you're doing, regularly, and definitely at the end, either with an alternative proof, or with some numerics.
\end{advice}

This is something very serious. Unless you're doing something very familiar, that you're used to for at least 5-10 years or so, like doing additions and multiplications for you, or some easy calculus for me, formulae and proofs that you can come upon are by default wrong. In order to make them correct, and ready to use, you must check and doublecheck and correct them, helped by alternative methods, or numerics.

\bigskip

Which brings us into the question on whether mathematics is an exact science or not. Not clear. Chemistry for instance is an exact science, because findings of type ``a mixture of water and salt cannot explode'' look rock-solid. Same for biology, with findings of type ``crocodiles eat fish'' being rock-solid too. In what regards mathematics however, and theoretical physics too, things are always prone to human mistake.

\bigskip

And for ending this discussion, you might ask then, what about engineering? After all, this is mathematics and physics, which is usually $100\%$ correct, because most of the bridges, buildings and other things built by engineers don't collapse. Well, this is because engineers follow, and in a truly maniac way, the above Advice 1.3. You won't declare a project for a bridge, building, engine and so on final and correct, ready for production, until you checked and doublechecked it with 10 different methods or so, won't you.

\bigskip

Back to work now, as an important adding to Theorem 1.2, we have:

\begin{convention}
By definition, $0!=1$.
\end{convention}

This convention comes, and no surprise here, from Advice 1.3. Indeed, we obviously have $\binom{n}{n}=1$, but if we want to recover this formula via Theorem 1.2 we are a bit in trouble, and so we must declare that $0!=1$, as for the following computation to work:
$$\binom{n}{n}=\frac{n!}{n!0!}=\frac{n!}{n!\times1}=1$$

Going ahead now with more mathematics and less philosophy, with Theorem 1.2 complemented by Convention 1.4 being in final form (trust me), we have:

\index{binomial formula}

\begin{theorem}
We have the binomial formula
$$(a+b)^n=\sum_{k=0}^n\binom{n}{k}a^kb^{n-k}$$
valid for any two numbers $a,b\in\mathbb Q$. 
\end{theorem}

\begin{proof}
We have to compute the following quantity, with $n$ terms in the product:
$$(a+b)^n=(a+b)(a+b)\ldots(a+b)$$

When expanding, we obtain a certain sum of products of $a,b$ variables, with each such product being a quantity of type $a^kb^{n-k}$. Thus, we have a formula as follows:
$$(a+b)^n=\sum_{k=0}^nC_ka^kb^{n-k}$$

In order to finish, it remains to compute the coefficients $C_k$. But, according to our product formula, $C_k$ is the number of choices for the $k$ needed $a$ variables among the $n$ available $a$ variables. Thus, according to Theorem 1.2, we have:
$$C_k=\binom{n}{k}$$

We are therefore led to the formula in the statement.
\end{proof}

Theorem 1.5 is something quite interesting, so let us doublecheck it with some numerics. At small values of $n$ we obtain the following formulae, which are all correct:
$$(a+b)^0=1$$
$$(a+b)^1=a+b$$
$$(a+b)^2=a^2+2ab+b^2$$
$$(a+b)^3=a^3+3a^2b+3ab^2+b^3$$
$$(a+b)^4=a^4+4a^3b+6a^2b^2+4ab^3+b^4$$
$$(a+b)^5=a^5+5a^4b+10a^3b^2+10a^2b^3+5ab^4+b^5$$
$$\vdots$$

Now observe that in these formulae, say for memorization purposes, the powers of the $a,b$ variables are something very simple, that can be recovered right away. What matters are the coefficients, which are the binomial coefficients $\binom{n}{k}$, which form a triangle. So, it is enough to memorize this triangle, and this can be done by using:

\index{Pascal triangle}

\begin{theorem}
The Pascal triangle, formed by the binomial coefficients $\binom{n}{k}$,
$$1$$
$$1\ \ ,\ \ 1$$
$$1\ \ ,\ \ 2\ \ ,\ \ 1$$
$$1\ \ ,\ \ 3\ \ ,\ \ 3\ \ ,\ \ 1$$
$$1\ \ ,\ \ 4\ \ ,\ \ 6\ \ ,\ \ 4\ \ ,\ \ 1$$
$$1\ \ ,\ \ 5\ \ ,\ \ 10\ \ ,\ \ 10\ \ ,\ \ 5\ \ ,\ \ 1$$
$$\vdots$$
has the property that each entry is the sum of the two entries above it.
\end{theorem}

\begin{proof}
In practice, the theorem states that the following formula holds:
$$\binom{n}{k}=\binom{n-1}{k-1}+\binom{n-1}{k}$$

There are many ways of proving this formula, all instructive, as follows:

\medskip

(1) Brute-force computation. We have indeed, as desired:
\begin{eqnarray*}
\binom{n-1}{k-1}+\binom{n-1}{k}
&=&\frac{(n-1)!}{(k-1)!(n-k)!}+\frac{(n-1)!}{k!(n-k-1)!}\\
&=&\frac{(n-1)!}{(k-1)!(n-k-1)!}\left(\frac{1}{n-k}+\frac{1}{k}\right)\\
&=&\frac{(n-1)!}{(k-1)!(n-k-1)!}\cdot\frac{n}{k(n-k)}\\
&=&\binom{n}{k}
\end{eqnarray*}

(2) Algebraic proof. We have the following formula, to start with:
$$(a+b)^n=(a+b)^{n-1}(a+b)$$

By using the binomial formula, this formula becomes:
$$\sum_{k=0}^n\binom{n}{k}a^kb^{n-k}=\left[\sum_{r=0}^{n-1}\binom{n-1}{r}a^rb^{n-1-r}\right](a+b)$$

Now let us perform the multiplication on the right. We obtain a certain sum of terms of type $a^kb^{n-k}$, and to be more precise, each such $a^kb^{n-k}$ term can either come from the $\binom{n-1}{k-1}$ terms $a^{k-1}b^{n-k}$ multiplied by $a$, or from the $\binom{n-1}{k}$ terms $a^kb^{n-1-k}$ multiplied by $b$. Thus, the coefficient of $a^kb^{n-k}$ on the right is $\binom{n-1}{k-1}+\binom{n-1}{k}$, as desired.

\medskip

(3) Combinatorics. Let us count $k$ objects among $n$ objects, with one of the $n$ objects having a hat on top. Obviously, the hat has nothing to do with the count, and we obtain $\binom{n}{k}$. On the other hand, we can say that there are two possibilities. Either the object with hat is counted, and we have $\binom{n-1}{k-1}$ possibilities here, or the object with hat is not counted, and we have $\binom{n-1}{k}$ possibilities here. Thus $\binom{n}{k}=\binom{n-1}{k-1}+\binom{n-1}{k}$, as desired.
\end{proof}

There are many more things that can be said about binomial coefficients, with all sorts of interesting formulae, but the idea is always the same, namely that in order to find such formulae you have a choice between algebra and combinatorics, and that when it comes to proofs, the brute-force computation method is useful too. In practice, the best is to master all 3 techniques. Among others, because of Advice 1.3. You will have in this way 3 different methods, for making sure that your formulae are correct indeed.

\section*{1b. Real numbers, analysis}

All the above was very nice, but remember that we are here for doing science and physics, and more specifically for mathematically understanding the numeric variables $x,y,z,\ldots$ coming from real life. Such variables can be lengths, volumes, pressures and so on, which vary continuously with time, and common sense dictates that there is little to no chance for our variables to be rational, $x,y,z,\ldots\notin\mathbb Q$. In fact, we will even see soon a theorem, stating that the probability for such a variable to be rational is exactly 0. Or, to put it in a dramatic way, ``rational numbers don't exist in real life''.

\bigskip

You are certainly familiar with the real numbers, but let us review now their definition, which is something quite tricky. As a first goal, we would like to construct a number $x=\sqrt{2}$ having the property $x^2=2$. But how to do this? Let us start with:

\index{square root}

\begin{proposition}
There is no number $r\in\mathbb Q_+$ satisfying $r^2=2$. In fact, we have
$$\mathbb Q_+=\left\{p\in\mathbb Q_+\Big|p^2<2\right\}\bigsqcup\left\{q\in\mathbb Q_+\Big|q^2>2\right\}$$
with this being a disjoint union.
\end{proposition}

\begin{proof}
In what regards the first assertion, assuming that $r=a/b$ with $a,b\in\mathbb N$ prime to each other satisfies $r^2=2$, we have $a^2=2b^2$, so $a\in2\mathbb N$. But by using again $a^2=2b^2$ we obtain $b\in2\mathbb N$, contradiction. As for the second assertion, this is obvious.
\end{proof}

It looks like we are a bit stuck. We can't really tell who $\sqrt{2}$ is, and the only piece of information about $\sqrt{2}$ that we have comes from the knowledge of the rational numbers satisfying $p^2<2$ or $q^2>2$. To be more precise, the picture that emerges is:

\begin{conclusion}
The number $\sqrt{2}$ is the abstract beast which is bigger than all rationals satisfying $p^2<2$, and smaller than all positive rationals satisfying $q^2>2$.
\end{conclusion}

This does not look very good, but you know what, instead of looking for more clever solutions to our problem, what about relaxing, or being lazy, or coward, or you name it, and taking Conclusion 1.8 as a definition for $\sqrt{2}$. This is actually something not that bad, and leads to the following ``lazy'' definition for the real numbers:

\index{real number}
\index{Dedekind cut}

\begin{definition}
The real numbers $x\in\mathbb R$ are formal cuts in the set of rationals,
$$\mathbb Q=A_x\sqcup B_x$$
with such a cut being by definition subject to the following conditions:
$$p\in A_x\ ,\ q\in B_x\implies p<q\qquad,\qquad\inf B_x\notin B_x$$
These numbers add and multiply by adding and multiplying the corresponding cuts.
\end{definition}

As a first observation, we have an inclusion $\mathbb Q\subset\mathbb R$, obtained by identifying each rational number $r\in\mathbb Q$ with the obvious cut that it produces, namely:
$$A_r=\left\{p\in\mathbb Q\Big|p\leq r\right\}
\quad,\quad
B_r=\left\{q\in\mathbb Q\Big|q>r\right\}$$

As a second observation, the addition and multiplication of real numbers, obtained by adding and multiplying the corresponding cuts, in the obvious way, is something very simple. To be more precise, in what regards the addition, the formula is as follows:
$$A_{x+y}=A_x+A_y$$ 

As for the multiplication, the formula here is similar, namely $A_{xy}=A_xA_y$, up to some mess with positives and negatives, which is quite easy to untangle, and with this being a good exercise. We can also talk about order between real numbers, as follows:
$$x\leq y\iff A_x\subset A_y$$

But let us perhaps leave more abstractions for later, and go back to more concrete things. As a first success of our theory, we can formulate the following theorem:

\index{square root}

\begin{theorem}
The equation $x^2=2$ has two solutions over the real numbers, namely the positive solution, denoted $\sqrt{2}$, and its negative counterpart, which is $-\sqrt{2}$.
\end{theorem}

\begin{proof}
By using $x\to-x$, it is enough to prove that $x^2=2$ has exactly one positive solution $\sqrt{2}$. But this is clear, because $\sqrt{2}$ can only come from the following cut:
$$A_{\sqrt{2}}=\mathbb Q_-\bigsqcup\left\{p\in\mathbb Q_+\Big|p^2<2\right\}\quad,\quad B_{\sqrt{2}}=\left\{q\in\mathbb Q_+\Big|q^2>2\right\}$$

Thus, we are led to the conclusion in the statement.
\end{proof}

More generally, the same method works in order to extract the square root $\sqrt{r}$ of any number $r\in\mathbb Q_+$, or even of any number $r\in\mathbb R_+$, and we have the following result:

\index{degree 2 equation}

\begin{theorem}
The solutions of $ax^2+bx+c=0$ with $a,b,c\in\mathbb R$ are
$$x_{1,2}=\frac{-b\pm\sqrt{b^2-4ac}}{2a}$$
provided that $b^2-4ac\geq0$. In the case $b^2-4ac<0$, there are no solutions.
\end{theorem}

\begin{proof}
We can write our equation in the following way:
\begin{eqnarray*}
ax^2+bx+c=0
&\iff&x^2+\frac{b}{a}x+\frac{c}{a}=0\\
&\iff&\left(x+\frac{b}{2a}\right)^2-\frac{b^2}{4a^2}+\frac{c}{a}=0\\
&\iff&x+\frac{b}{2a}=\pm\frac{\sqrt{b^2-4ac}}{2a}
\end{eqnarray*}

Thus, we are led to the conclusion in the statement.
\end{proof}

Summarizing, we have a nice definition for the real numbers, that we can certainly do some math with. However, for anything more advanced we are in need of the decimal writing for the real numbers. The result here is as follows:

\index{decimal writing}

\begin{theorem}
The real numbers $x\in\mathbb R$ can be written in decimal form,
$$x=\pm a_1\ldots a_n.b_1b_2b_3\ldots\ldots$$
with $a_i,b_i\in\{0,1,\ldots,9\}$, with the convention $\ldots b999\ldots=\ldots(b+1)000\ldots$
\end{theorem}

\begin{proof}
This is something non-trivial, even for the rationals $x\in\mathbb Q$ themselves, which require some work in order to be put in decimal form, the idea being as follows:

\medskip

(1) First of all, our precise claim is that any $x\in\mathbb R$ can be written in the form in the statement, with the integer $\pm a_1\ldots a_n$ and then each of the digits $b_1,b_2,b_3,\ldots$ providing the best approximation of $x$, at that stage of the approximation. 

\medskip

(2) Moreover, we have a second claim as well, namely that any expression of type $x=\pm a_1\ldots a_n.b_1b_2b_3\ldots\ldots$ corresponds to a real number $x\in\mathbb R$, and that with the convention $\ldots b999\ldots=\ldots(b+1)000\ldots\,$, the correspondence is bijective.

\medskip

(3) In order to prove now these two assertions, our first claim is that we can restrict the attention to the case $x\in[0,1)$, and with this meaning of course $0\leq x<1$, with respect to the order relation for the reals discussed in the above.

\medskip

(4) Getting started now, let $x\in\mathbb R$, coming from a cut $\mathbb Q=A_x\sqcup B_x$. Since the set $A_x\cap\mathbb Z$ consists of integers, and is bounded from above by any element $q\in B_x$ of your choice, this set has a maximal element, that we can denote $[x]$:
$$[x]=\max\left(A_x\cap\mathbb Z\right)$$

It follows from definitions that $[x]$ has the usual properties of the integer part, namely:
$$[x]\leq x<[x]+1$$

Thus we have $x=[x]+y$ with $[x]\in\mathbb Z$ and $y\in[0,1)$, and getting back now to what we want to prove, namely (1,2) above, it is clear that it is enough to prove these assertions for the remainder $y\in[0,1)$. Thus, we have proved (3), and we can assume $x\in[0,1)$.

\medskip

(5) So, assume $x\in[0,1)$. We are first looking for a best approximation from below of type $0.b_1$, with $b_1\in\{0,\ldots,9\}$, and it is clear that such an approximation exists, simply by comparing $x$ with the numbers $0.0,0.1,\ldots,0.9$. Thus, we have our first digit $b_1$, and then we can construct the second digit $b_2$ as well, by comparing $x$ with the numbers $0.b_10,0.b_11,\ldots,0.b_19$. And so on, which finishes the proof of our claim (1).

\medskip

(6) In order to prove now the remaining claim (2), let us restrict again the attention, as explained in (4), to the  case $x\in[0,1)$. First, it is clear that any expression of type $x=0.b_1b_2b_3\ldots$ defines a real number $x\in[0,1]$, simply by declaring that the corresponding cut $\mathbb Q=A_x\sqcup B_x$ comes from the following set, and its complement:
$$A_x=\bigcup_{n\geq1}\left\{p\in\mathbb Q\Big|p\leq 0.b_1\ldots b_n\right\}$$

(7) Thus, we have our correspondence between real numbers as cuts, and real numbers as decimal expressions, and we are left with the question of investigating the bijectivity of this correspondence. But here, the only bug that happens is that numbers of type $x=\ldots b999\ldots$, which produce reals $x\in\mathbb R$ via (6), do not come from reals $x\in\mathbb R$ via (5). So, in order to finish our proof, we must investigate such numbers.

\medskip

(8) So, consider an expression of type $\ldots b999\ldots$ Going back to the construction in (6), we are led to the conclusion that we have the following equality:
$$A_{b999\ldots}=B_{(b+1)000\ldots}$$

Thus, at the level of the real numbers defined as cuts, we have:
$$\ldots b999\ldots=\ldots(b+1)000\ldots$$

But this solves our problem, because by identifying $\ldots b999\ldots=\ldots(b+1)000\ldots$ the bijectivity issue of our correspondence is fixed, and we are done.
\end{proof}

The above theorem was of course quite difficult, but this is how things are. Let us record as well the following result, coming as a useful complement to the above:

\index{probability 0}

\begin{theorem}
The following happen, in relation with $\mathbb Q\to\mathbb R$:
\begin{enumerate}
\item $\mathbb Q$ is countable, while $\mathbb R$ is not countable.

\item $r\in\mathbb R$ is rational precisely when its decimal writing is periodic.

\item The probability for a randomly picked $x\in\mathbb R$ to be rational is $0$.
\end{enumerate}
\end{theorem}

\begin{proof}
We have several things to be proved, the idea being as follows:

\medskip

(1) We can count the positive rationals, with some redundancies, as follows:
$$\xymatrix@R=10pt@C=10pt{
1/1\ar[r]&1/2\ar[dl]&1/3\ar[r]&1/4\ar[dl]&\ldots\\
2/1\ar[d]&2/2\ar[ur]&2/3\ar[dl]&2/4&\ldots\\
3/1\ar[ur]&3/2\ar[dl]&3/3&3/4&\ldots\\
4/1&4/2&4/3&4/4&\ldots\\
\vdots&\vdots&\vdots&\vdots&\ddots}$$

Now after eliminating the redundancies, and then adding the negatives, which must be countable too, say via an alternating $+/-$ scheme, countability of $\mathbb Q$ proved.
 
\medskip

(2) Regarding now the reals, assume by contradiction that $[0,1]$ is countable, listed as follows, and with the convention that the writings of type $\ldots999\ldots$ are avoided:
$$x_1=0.a_1a_2a_3\ldots$$
$$x_2=0.b_1b_2b_3\ldots$$
$$x_3=0.c_1c_2c_3\ldots$$
$$\ldots$$

Now pick digits $\sigma_1\neq a_1$, $\sigma_2\neq b_2$, $\sigma_3\neq c_3$ and so on, again with a technical convention here, that these are different from 9, and define $x\in\mathbb R$ as follows:
$$x=0.\sigma_1\sigma_2\sigma_3\ldots$$

We have then $x\in[0,1]$, and since $x$ is obviously not on the above list, this is a contradiction. Thus $[0,1]$ is not countable, and it follows that $\mathbb R$ is not countable either.

\medskip

(3) Assuming that $r\in\mathbb R$ has periodic decimal writing, the following computation, based on $(10^p-1)(\sum_k10^{-kp})=1$, obtained by multiplying, gives $r\in\mathbb Q$:
\begin{eqnarray*}
r
&=&\pm a_1\ldots a_m.b_1\ldots b_nc_1\ldots c_pc_1\ldots c_p\ldots\\
&=&\pm\frac{1}{10^n}\left(a_1\ldots a_mb_1\ldots b_n+c_1\ldots c_p\left(\frac{1}{10^p}+\frac{1}{10^{2p}}+\ldots\right)\right)\\
&=&\pm\frac{1}{10^n}\left(a_1\ldots a_mb_1\ldots b_n+\frac{c_1\ldots c_p}{10^p-1}\right)
\end{eqnarray*}

(4) Conversely, given a rational number $r=k/l$, we can find its decimal writing by performing the usual division algorithm, $k$ divided by $l$. But this algorithm will be surely periodic, after some time, so the decimal writing of $r$ is indeed periodic.

\medskip

(5) Finally, regarding the probability assertion, in order to avoid some troubles, we will prove that the probability for a real number $x\in[0,1]$ to be rational is 0. So, let us write the rational numbers $r\in[0,1]$ in the form of a sequence $r_1,r_2,r_3\ldots\,$ as follows:
$$\mathbb Q\cap[0,1]=\big\{r_1,r_2,r_3,\ldots\big\}$$

Let us also pick a number $c>0$. Since the probability of having $x=r_1$ is certainly smaller than $c/2$, then the probability of having $x=r_2$ is certainly smaller than $c/4$, then the probability of having $x=r_3$ is certainly smaller than $c/8$ and so on, the probability for $x$ to be rational satisfies the following inequality:
$$P\leq\frac{c}{2}+\frac{c}{4}+\frac{c}{8}+\ldots=c\left(\frac{1}{2}+\frac{1}{4}+\frac{1}{8}+\ldots\right)=c$$

Here we have used the well-known formula $\frac{1}{2}+\frac{1}{4}+\frac{1}{8}+\ldots=1$, which comes by dividing $[0,1]$ into half, and then one of the halves into half again, and so on, and then saying in the end that the pieces that we have must sum up to 1. Thus, we have indeed $P\leq c$, and since the number $c>0$ was arbitrary, we obtain $P=0$, as desired.
\end{proof}

Moving ahead now, let us construct now some more real numbers. We already know about $\sqrt{2}$ and other numbers of the same type, namely roots of polynomials, and our knowledge here being quite decent, no hurry with this, we will be back to it later. So, let us get now into $\pi$ and trigonometry. To start with, we have the following result:

\index{pi}

\begin{theorem}
The following two definitions of $\pi$ are equivalent:
\begin{enumerate}
\item The length of the unit circle is $L=2\pi$.

\item The area of the unit disk is $A=\pi$.
\end{enumerate}
\end{theorem}

\begin{proof}
In order to prove this theorem let us cut the unit disk as a pizza, into $N$ slices, and forgetting about gastronomy, leave aside the rounded parts:
$$
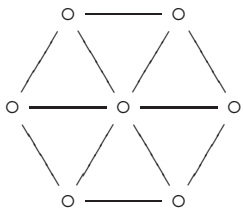
\xymatrix@R=23pt@C=8pt{
&\circ\ar@{-}[rr]\ar@{-}[dl]\ar@{-}[dr]&&\circ\ar@{-}[dl]\ar@{-}[dr]\\
\circ\ar@{-}[rr]&&\circ\ar@{-}[rr]&&\circ\\
&\circ\ar@{-}[rr]\ar@{-}[ul]\ar@{-}[ur]&&\circ\ar@{-}[ul]\ar@{-}[ur]
}$$

The area to be eaten can be then computed as follows, where $H$ is the height of the slices, $S$ is the length of their sides, and $P=NS$ is the total length of the sides:
$$A
=N\times \frac{HS}{2}
=\frac{HP}{2}
\simeq\frac{1\times L}{2}$$

Thus, with $N\to\infty$ we obtain that we have $A=L/2$, as desired.
\end{proof}

In what regards now the precise value of $\pi$, the above picture at $N=6$ shows that we have $\pi>3$, but not by much. The precise figure is $\pi=3.14159\ldots\,$, but we will come back to this later, once we will have appropriate tools for dealing with such questions. It is also possible to prove that $\pi$ is irrational, $\pi\notin\mathbb Q$, but this is not trivial either.

\bigskip

Let us end this discussion about real numbers with some trigonometry. There are many things that can be said, that you certainly know, the basics being as follows:

\index{cos}
\index{sin}
\index{Pythagoras' theorem}

\begin{theorem}
The following happen:
\begin{enumerate}
\item We can talk about angles $t\in\mathbb R$, by using the unit circle, in the usual way, and in this correspondence, the right angle has a value of $\pi/2$.

\item Associated to any $t\in\mathbb R$ are numbers $\sin t,\cos t\in\mathbb R$, constructed in the usual way, by using a triangle. These numbers satisfy $\sin^2t+\cos^2t=1$.
\end{enumerate}
\end{theorem}

\begin{proof}
There are certainly things that you know, the idea being as follows:

\medskip

(1) The formula $L=2\pi$ from Theorem 1.14 shows that the length of a quarter of the unit circle is $l=\pi/2$, and so the right angle has indeed this value, $\pi/2$.

\medskip

(2) As for $\sin^2t+\cos^2t=1$, called Pythagoras' theorem, this comes from the following picture, consisting of two squares and four identical triangles, as indicated:
$$
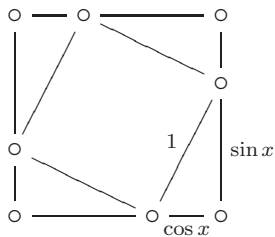
\xymatrix@R=13pt@C=13pt{
\circ\ar@{-}[r]\ar@{-}[dd]&\circ\ar@{-}[rr]\ar@{-}[drr]&&\circ\ar@{-}[d]\\
&&&\circ\ar@{-}[dd]^{\sin t}&\\
\circ\ar@{-}[d]\ar@{-}[uur]\ar@{-}[drr]&&&&\\
\circ\ar@{-}[rr]&&\circ\ar@{-}[r]_{\cos t}\ar@{-}[uur]^1&\circ
}$$

Indeed, when computing the area of the outer square, in two ways, we obtain:
$$(\sin t+\cos t)^2=1+4\times\frac{\sin t\cos t}{2}$$

Now when expanding we obtain from this $\sin^2t+\cos^2t=1$, as claimed.
\end{proof}

It is possible to say many more things about angles and $\sin t$, $\cos t$, and also talk about some supplementary quantities, such as $\tan t=\sin t/\cos t$. But more on this later, once we will have some appropriate tools, beyond basic geometry, in order to discuss this.

\section*{1c. Sequences, convergence} 

We already met, on several occasions, infinite sequences or sums, and their limits. Time now to clarify all this. Let us start with the following definition:

\index{sequence}
\index{convergent sequence}
\index{limit of sequence}

\begin{definition}
We say that a sequence $\{x_n\}_{n\in\mathbb N}\subset\mathbb R$ converges to $x\in\mathbb R$ when:
$$\forall\varepsilon>0,\exists N\in\mathbb N,\forall n\geq N,|x_n-x|<\varepsilon$$
In this case, we write $\lim_{n\to\infty}x_n=x$, or simply $x_n\to x$.
\end{definition}

This might look quite scary, at a first glance, but when thinking a bit, there is nothing scary about it. Indeed, let us try to understand, how shall we translate $x_n\to x$ into mathematical language. The condition $x_n\to x$ tells us that ``when $n$ is big, $x_n$ is close to $x$'', and to be more precise, it tells us that ``when $n$ is big enough, $x_n$ gets arbitrarily close to $x$''. But $n$ big enough means $n\geq N$, for some $N\in\mathbb N$, and $x_n$ arbitrarily close to $x$ means $|x_n-x|<\varepsilon$, for some $\varepsilon>0$. Thus, we are led to the above definition. 

\bigskip

As a basic example for all this, we have:

\begin{proposition}
We have $1/n\to0$.
\end{proposition}

\begin{proof}
This is obvious, but let us prove it by using Definition 1.16. We have:
$$\left|\frac{1}{n}-0\right|<\varepsilon
\iff\frac{1}{n}<\varepsilon
\iff\frac{1}{\varepsilon}<n$$

Thus we can take $N=[1/\varepsilon]+1$ in Definition 1.16, and we are done.
\end{proof}

There are many other examples, and more on this in a moment. Going ahead with more theory, let us complement Definition 1.16 with:

\begin{definition}
We write $x_n\to\infty$ when the following condition is satisfied:
$$\forall K>0,\exists N\in\mathbb N,\forall n\geq N,x_n>K$$
Similarly, we write $x_n\to-\infty$ when the same happens, with $x_n<-K$ at the end.
\end{definition}

Again, this is something very intuitive, coming from the fact that $x_n\to\infty$ can only mean that $x_n$ is arbitrarily big, for $n$ big enough. As a basic illustration, we have:

\begin{proposition}
We have $n^2\to\infty$.
\end{proposition}

\begin{proof}
As before, this is obvious, but let us prove it using Definition 1.18. We have:
$$n^2>K\iff n>\sqrt{K}$$

Thus we can take $N=[\sqrt{K}]+1$ in Definition 1.18, and we are done.
\end{proof}

We can unify and generalize Proposition 1.17 and Proposition 1.19, as follows:

\begin{proposition}
We have the following convergence, with $n\to\infty$:
$$n^a\to\begin{cases}
0&(a<0)\\
1&(a=0)\\
\infty&(a>0)
\end{cases}$$
\end{proposition}

\begin{proof}
This follows indeed by using the same method as in the proof of Proposition 1.17 and Proposition 1.19, first for $a$ rational, and then for $a$ real as well.
\end{proof}

We have some general results about limits, summarized as follows:

\index{increasing sequence}
\index{decreasing sequence}
\index{subsequence}
\index{bounded sequence}

\begin{theorem}
The following happen:
\begin{enumerate}
\item The limit $\lim_{n\to\infty}x_n$, if it exists, is unique.

\item If $x_n\to x$, with $x\in(-\infty,\infty)$, then $x_n$ is bounded.

\item If $x_n$ is increasing or descreasing, then it converges.

\item Assuming $x_n\to x$, any subsequence of $x_n$ converges to $x$.
\end{enumerate}
\end{theorem}

\begin{proof}
All this is elementary, coming from definitions:

\medskip

(1) Assuming $x_n\to x$, $x_n\to y$ we have indeed, for any $\varepsilon>0$, for $n$ big enough:
$$|x-y|\leq|x-x_n|+|x_n-y|<2\varepsilon$$

(2) Assuming $x_n\to x$, we have $|x_n-x|<1$ for $n\geq N$, and so, for any $k\in\mathbb N$:
$$|x_k|<1+|x|+\sup\left(|x_1|,\ldots,|x_{n-1}|\right)$$

(3) By using $x\to-x$, it is enough to prove the result for increasing sequences. But here we can construct the limit $x\in(-\infty,\infty]$ in the following way:
$$\bigcup_{n\in\mathbb N}(-\infty,x_n)=(-\infty,x)$$

(4) This is clear from definitions.
\end{proof}

Here are as well some general rules for computing limits:

\begin{theorem}
The following happen, with the conventions $\infty+\infty=\infty$, $\infty\cdot\infty=\infty$, $1/\infty=0$, and with the conventions that $\infty-\infty$ and $\infty\cdot0$ are undefined:
\begin{enumerate}
\item $x_n\to x$ implies $\lambda x_n\to\lambda x$.

\item $x_n\to x$, $y_n\to y$ implies $x_n+y_n\to x+y$.

\item $x_n\to x$, $y_n\to y$ implies $x_ny_n\to xy$.

\item $x_n\to x$ with $x\neq0$ implies $1/x_n\to 1/x$.
\end{enumerate}
\end{theorem}

\begin{proof}
All this is again elementary, coming from definitions:

\medskip

(1) This is something which is obvious from definitions.

\medskip

(2) This follows indeed from the following estimate:
$$|x_n+y_n-x-y|\leq|x_n-x|+|y_n-y|$$

(3) This follows indeed from the following estimate:
\begin{eqnarray*}
|x_ny_n-xy|
&=&|(x_n-x)y_n+x(y_n-y)|\\
&\leq&|x_n-x|\cdot|y_n|+|x|\cdot|y_n-y|
\end{eqnarray*}

(4) This is again clear, by estimating $1/x_n-1/x$, in the obvious way.
\end{proof}

As an application of the above rules, we have the following useful result:

\begin{proposition}
The $n\to\infty$ limits of quotients of polynomials are given by
$$\lim_{n\to\infty}\frac{a_pn^p+a_{p-1}n^{p-1}+\ldots+a_0}{b_qn^q+b_{q-1}n^{q-1}+\ldots+b_0}=\lim_{n\to\infty}\frac{a_pn^p}{b_qn^q}$$
with the limit on the right being $\pm\infty$, $0$, $a_p/b_q$, depending on the values of $p,q$.
\end{proposition}

\begin{proof}
The first assertion comes from the following computation:
\begin{eqnarray*}
\lim_{n\to\infty}\frac{a_pn^p+a_{p-1}n^{p-1}+\ldots+a_0}{b_qn^q+b_{q-1}n^{q-1}+\ldots+b_0}
&=&\lim_{n\to\infty}\frac{n^p}{n^q}\cdot\frac{a_p+a_{p-1}n^{-1}+\ldots+a_0n^{-p}}{b_q+b_{q-1}n^{-1}+\ldots+b_0n^{-q}}\\
&=&\lim_{n\to\infty}\frac{a_pn^p}{b_qn^q}
\end{eqnarray*}

As for the second assertion, this comes from Proposition 1.20.
\end{proof}

Getting back now to theory, some sequences which obviously do not converge, like for instance $x_n=(-1)^n$, have however ``2 limits instead of 1''. So let us formulate:

\index{lim sup}
\index{lim inf}
\index{subsequence}

\begin{definition}
Given a sequence $\{x_n\}_{n\in\mathbb N}\subset\mathbb R$, we let
$$\liminf_{n\to\infty}x_n\in[-\infty,\infty]\quad,\quad\limsup_{n\to\infty}x_n\in[-\infty,\infty]$$
to be the smallest and biggest limit of a subsequence of $(x_n)$.
\end{definition}

Observe that the above quantities are defined indeed for any sequence $x_n$. For instance, for $x_n=(-1)^n$ we obtain $-1$ and $1$. Also, for $x_n=n$ we obtain $\infty$ and $\infty$. And so on. Of course, and generalizing the $x_n=n$ example, if $x_n\to x$ we obtain $x$ and $x$.

\bigskip

Going ahead with more theory, here is a key result:

\index{Cauchy sequence}

\begin{theorem}
A sequence $x_n$ converges, with finite limit $x\in\mathbb R$, precisely when
$$\forall\varepsilon>0,\exists N\in\mathbb N,\forall m,n\geq N,|x_m-x_n|<\varepsilon$$
called Cauchy condition.
\end{theorem}

\begin{proof}
In one sense, this is clear. In the other sense, we can say for instance that the Cauchy condition forces the decimal writings of our numbers $x_n$ to coincide more and more, with $n\to\infty$, and so we can construct a limit $x=\lim_{n\to\infty}x_n$, as desired.
\end{proof}

The above result is quite interesting, and as an application, we have:

\index{complete space}

\begin{theorem}
$\mathbb R$ is the completion of $\mathbb Q$, in the sense that it is the space of Cauchy sequences over $\mathbb Q$, identified when the virtual limit is the same, in the sense that: 
$$x_n\sim y_n\iff |x_n-y_n|\to0$$
Moreover, $\mathbb R$ is complete, in the sense that it equals its own completion.
\end{theorem}

\begin{proof}
Let us denote the completion operation by $X\to\bar{X}=C_X/\sim$, where $C_X$ is the space of Cauchy sequences over $X$, and $\sim$ is the above equivalence relation. Since by Theorem 1.25 any Cauchy sequence $(x_n)\in C_\mathbb Q$ has a limit $x\in\mathbb R$, we obtain $\bar{\mathbb Q}=\mathbb R$. As for the equality $\bar{\mathbb R}=\mathbb R$, this is clear again by using Theorem 1.25.
\end{proof}

\section*{1d. Series, the number e}

With the above understood, we are now ready to get into some truly interesting mathematics. Let us start with the following definition:

\index{series}
\index{limit of series}
\index{convergent series}

\begin{definition}
Given numbers $x_0,x_1,x_2,\ldots\in\mathbb R$, we write
$$\sum_{n=0}^\infty x_n=x$$
with $x\in[-\infty,\infty]$ when $\lim_{k\to\infty}\sum_{n=0}^kx_n=x$.
\end{definition}

As before with the sequences, there is some general theory that can be developed for the series, and more on this in a moment. As a first, basic example, we have:

\index{geometric series}

\begin{theorem}
We have the ``geometric series'' formula
$$\sum_{n=0}^\infty x^n=\frac{1}{1-x}$$
valid for any $|x|<1$. For $|x|\geq1$, the series diverges.
\end{theorem}

\begin{proof}
Our first claim, which comes by multiplying and simplifying, is that:
$$\sum_{n=0}^kx^n=\frac{1-x^{k+1}}{1-x}$$

But this proves the first assertion, because with $k\to\infty$ we get:
$$\sum_{n=0}^kx^n\to\frac{1}{1-x}$$

As for the second assertion, this is clear as well from our formula above.
\end{proof}

Less trivial now is the following result, due to Riemann:

\index{Riemann series}

\begin{theorem}
We have the following formula:
$$1+\frac{1}{2}+\frac{1}{3}+\frac{1}{4}+\ldots=\infty$$
In fact, $\sum_n1/n^a$ converges for $a>1$, and diverges for $a\leq1$.
\end{theorem}

\begin{proof}
We have to prove several things, the idea being as follows:

\medskip

(1) The first assertion comes from the following computation:
\begin{eqnarray*}
1+\frac{1}{2}+\frac{1}{3}+\frac{1}{4}+\ldots
&=&1+\frac{1}{2}+\left(\frac{1}{3}+\frac{1}{4}\right)
+\left(\frac{1}{5}+\frac{1}{6}+\frac{1}{7}+\frac{1}{8}\right)+\ldots\\
&\geq&1+\frac{1}{2}+\left(\frac{1}{4}+\frac{1}{4}\right)
+\left(\frac{1}{8}+\frac{1}{8}+\frac{1}{8}+\frac{1}{8}\right)+\ldots\\
&=&1+\frac{1}{2}+\frac{1}{2}+\frac{1}{2}+\ldots\\
&=&\infty
\end{eqnarray*}

(2) Regarding now the second assertion, we have that at $a=1$, and so at any $a\leq1$. Thus, it remains to prove that at $a>1$ the series converges. Let us first discuss the case $a=2$, which will prove the convergence at any $a\geq2$. The trick here is as follows:
\begin{eqnarray*}
1+\frac{1}{4}+\frac{1}{9}+\frac{1}{16}+\ldots
&\leq&1+\frac{1}{3}+\frac{1}{6}+\frac{1}{10}+\ldots\\
&=&2\left(\frac{1}{2}+\frac{1}{6}+\frac{1}{12}+\frac{1}{20}+\ldots\right)\\
&=&2\left[\left(1-\frac{1}{2}\right)+\left(\frac{1}{2}-\frac{1}{3}\right)+\left(\frac{1}{3}-\frac{1}{4}\right)+\left(\frac{1}{4}-\frac{1}{5}\right)\ldots\right]\\
&=&2
\end{eqnarray*}

(3) It remains to prove that the series converges at $a\in(1,2)$, and here it is enough to deal with the case of the exponents $a=1+1/p$ with $p\in\mathbb N$. We already know how to do this at $p=1$, and the proof at $p\in\mathbb N$ will be based on a similar trick. We have:
$$\sum_{n=0}^\infty\frac{1}{n^{1/p}}-\frac{1}{(n+1)^{1/p}}=1$$

Let us compute, or rather estimate, the generic term of this series. By using the formula $a^p-b^p=(a-b)(a^{p-1}+a^{p-2}b+\ldots+ab^{p-2}+b^{p-1})$, we have:
\begin{eqnarray*}
\frac{1}{n^{1/p}}-\frac{1}{(n+1)^{1/p}}
&=&\frac{(n+1)^{1/p}-n^{1/p}}{n^{1/p}(n+1)^{1/p}}\\
&=&\frac{1}{n^{1/p}(n+1)^{1/p}[(n+1)^{1-1/p}+\ldots+n^{1-1/p}]}\\
&\geq&\frac{1}{n^{1/p}(n+1)^{1/p}\cdot p(n+1)^{1-1/p}}\\
&=&\frac{1}{pn^{1/p}(n+1)}\\
&\geq&\frac{1}{p(n+1)^{1+1/p}}
\end{eqnarray*}

We therefore obtain the following estimate for the Riemann sum:
\begin{eqnarray*}
\sum_{n=0}^\infty\frac{1}{n^{1+1/p}}
&=&1+\sum_{n=0}^\infty\frac{1}{(n+1)^{1+1/p}}\\
&\leq&1+p\sum_{n=0}^\infty\left(\frac{1}{n^{1/p}}-\frac{1}{(n+1)^{1/p}}\right)\\
&=&1+p
\end{eqnarray*}

Thus, we are done with the case $a=1+1/p$, which finishes the proof.
\end{proof}

Here is another tricky result, this time about alternating sums:

\index{alternating series}

\begin{theorem}
We have the following convergence result:
$$1-\frac{1}{2}+\frac{1}{3}-\frac{1}{4}+\ldots<\infty$$
However, when rearranging terms, we can obtain any $x\in[-\infty,\infty]$ as limit.
\end{theorem}

\begin{proof}
Both the assertions follow from Theorem 1.29, as follows:

\medskip

(1) We have the following computation, using the Riemann criterion at $a=2$:
\begin{eqnarray*}
1-\frac{1}{2}+\frac{1}{3}-\frac{1}{4}+\ldots
&=&\left(1-\frac{1}{2}\right)+\left(\frac{1}{3}-\frac{1}{4}\right)+\ldots\\
&=&\frac{1}{2}+\frac{1}{12}+\frac{1}{30}+\ldots\\
&<&\frac{1}{1^2}+\frac{1}{2^2}+\frac{1}{3^2}+\ldots\\
&<&\infty
\end{eqnarray*}

(2) We have the following formulae, coming from the Riemann criterion at $a=1$:
$$\frac{1}{2}+\frac{1}{4}+\frac{1}{6}+\frac{1}{8}+\ldots=\frac{1}{2}\left(1+\frac{1}{2}+\frac{1}{3}+\frac{1}{4}+\ldots\right)=\infty$$
$$1+\frac{1}{3}+\frac{1}{5}+\frac{1}{7}+\ldots\geq\frac{1}{2}+\frac{1}{4}+\frac{1}{6}+\frac{1}{8}+\ldots=\infty$$

Thus, both these series diverge. The point now is that, by using this, when rearranging terms in the alternating series in the statement, we can arrange for the partial sums to go arbitrarily high, or arbitrarily low, and we can obtain any $x\in[-\infty,\infty]$ as limit.
\end{proof}

Back now to the general case, we first have the following statement:

\begin{theorem}
The following hold, with the converses of $(1)$ and $(2)$ being wrong, and with $(3)$ not holding when the assumption $x_n\geq0$ is removed:
\begin{enumerate}
\item If $\sum_nx_n$ converges then $x_n\to0$.

\item If $\sum_n|x_n|$ converges then $\sum_nx_n$ converges.

\item If $\sum_nx_n$ converges, $x_n\geq0$ and $x_n/y_n\to1$ then $\sum_ny_n$ converges.
\end{enumerate}
\end{theorem}

\begin{proof}
This is a mixture of trivial and non-trivial results, as follows:

\medskip

(1) We know that $\sum_nx_n$ converges when $S_k=\sum_{n=0}^kx_n$ converges. Thus by Cauchy we have $x_k=S_k-S_{k-1}\to0$, and this gives the result. As for the simplest counterexample for the converse, this is $1+\frac{1}{2}+\frac{1}{3}+\frac{1}{4}+\ldots=\infty$, coming from Theorem 1.29.

\medskip

(2) This follows again from the Cauchy criterion, by using:
$$|x_n+x_{n+1}+\ldots+x_{n+k}|\leq|x_n|+|x_{n+1}|+\ldots+|x_{n+k}|$$

As for the simplest counterexample for the converse, this is $1-\frac{1}{2}+\frac{1}{3}-\frac{1}{4}+\ldots<\infty$, coming from Theorem 1.30, coupled with $1+\frac{1}{2}+\frac{1}{3}+\frac{1}{4}+\ldots=\infty$ from (1).

\medskip

(3) Again, the main assertion here is clear, coming from, for $n$ big:
$$(1-\varepsilon)x_n\leq y_n\leq(1+\varepsilon)x_n$$

In what regards now the failure of the result, when the assumption $x_n\geq0$ is removed, this is something quite tricky, the simplest counterexample being as follows:
$$x_n=\frac{(-1)^n}{\sqrt{n}}
\quad,\quad 
y_n=\frac{1}{n}+\frac{(-1)^n}{\sqrt{n}}$$

To be more precise, we have $y_n/x_n\to1$, so $x_n/y_n\to1$ too, but according to the above-mentioned results from (1,2), modified a bit, $\sum_nx_n$ converges, while $\sum_ny_n$ diverges.  
\end{proof} 

Summarizing, we have some useful positive results about series, which are however quite trivial, along with various counterexamples to their possible modifications, which are non-trivial. Staying positive, here are some more positive results:

\begin{theorem}
The following happen, and in all cases, the situtation where $c=1$ is indeterminate, in the sense that the series can converge or diverge:
\begin{enumerate}
\item If $|x_{n+1}/x_n|\to c$, the series $\sum_nx_n$ converges if $c<1$, and diverges if $c>1$.

\item If $\sqrt[n]{|x_n|}\to c$, the series $\sum_nx_n$ converges if $c<1$, and diverges if $c>1$.

\item With $c=\limsup_{n\to\infty}\sqrt[n]{|x_n|}$, $\sum_nx_n$ converges if $c<1$, and diverges if $c>1$.
\end{enumerate}
\end{theorem}

\begin{proof}
Again, this is a mixture of trivial and non-trivial results, as follows:

\medskip

(1) Here the main assertions, regarding the cases $c<1$ and $c>1$, are both clear by comparing with the geometric series $\sum_nc^n$. As for the case $c=1$, this is what happens for the Riemann series $\sum_n1/n^a$, so we can have both convergent and divergent series.

\medskip

(2) Again, the main assertions, where $c<1$ or $c>1$, are clear by comparing with the geometric series $\sum_nc^n$, and the $c=1$ examples come from the Riemann series.

\medskip

(3) Here the case $c<1$ is dealt with as in (2), and the same goes for the examples at $c=1$. As for the case $c>1$, this is clear too, because here $x_n\to0$ fails.
\end{proof} 

Finally, generalizing the first assertion in Theorem 1.30, we have:

\index{alternating series}

\begin{theorem}
If $x_n\searrow0$ then $\sum_n(-1)^nx_n$ converges.
\end{theorem}

\begin{proof}
We have the $\sum_n(-1)^nx_n=\sum_ky_k$, where:
$$y_k=x_{2k}-x_{2k+1}$$

But, by drawing for instance the numbers $x_i$ on the real line, we see that $y_k$ are positive numbers, and that $\sum_ky_k$ is the sum of lengths of certain disjoint intervals, included in the interval $[0,x_0]$. Thus we have $\sum_ky_k\leq x_0$, and this gives the result.
\end{proof}

In order to formulate our next result, we will need the following key fact:

\begin{theorem}
We have the following inequality, for any $a_1,\ldots,a_n\geq0$,
$$\frac{a_1+\ldots+a_n}{n}\geq\sqrt[n]{a_1\ldots a_n}$$
telling us that the arithmetic mean is bigger than the geometric mean.
\end{theorem}

\begin{proof}
This can be done in several steps, as follows:

\medskip

(1) To start with, the result holds indeed at $n=2$, with this coming from:
$$\frac{a+b}{2}\geq\sqrt{ab}\iff(\sqrt{a}-\sqrt{b})^2\geq0$$

(2) But with this, we can prove our inequality at $n=4$ too, as follows:
$$\frac{a+b+c+d}{4}
\geq\frac{\sqrt{ab}+\sqrt{cd}}{2}
\geq\sqrt[4]{abcd}$$

(3) Next, we can prove our inequality, in the same way, at $n=8$, then at $n=16$, and so on. Thus, as a conclusion, we know how to prove the result at any $n=2^s$.

\medskip

(4) In general now, given numbers $a_1,\ldots,a_n\geq0$, consider their arithmetic mean:
$$m=\frac{a_1+\ldots+a_n}{n}$$

Now pick $s\in\mathbb N$ such that $n\leq 2^s$, and let us complete our series $a_1,\ldots,a_n$ with $2^s-n$ copies of $m$. The arithmetic mean stays the same, and by using (3) we obtain:
\begin{eqnarray*}
m\geq\sqrt[2^s]{a_1\ldots a_nm^{2^s-n}}
&\implies&m^{2^s}\geq a_1\ldots a_nm^{2^s-n}\\
&\implies&m^n\geq a_1\ldots a_n\\
&\implies&m\geq\sqrt[n]{a_1\ldots a_n}
\end{eqnarray*}

Thus, we are led to the conclusion in the statement. 
\end{proof}

Good news, we can talk now about a very interesting convergence, as follows:

\index{e}

\begin{theorem}
We have the following convergence
$$\left(1+\frac{1}{n}\right)^n\to e$$
where $e=2.71828\ldots$ is a certain number.
\end{theorem}

\begin{proof}
This is something quite tricky, as follows:

\medskip

(1) Our first claim is that the following sequence is increasing:
$$x_n=\left(1+\frac{1}{n}\right)^n$$ 

In order to prove this, we use the following arithmetic-geometric inequality:
$$\frac{1+\sum_{i=1}^n\left(1+\frac{1}{n}\right)}{n+1}\geq\sqrt[n+1]{1\cdot\prod_{i=1}^n\left(1+\frac{1}{n}\right)}$$

In practice, this gives the following inequality:
$$1+\frac{1}{n+1}\geq\left(1+\frac{1}{n}\right)^{n/(n+1)}$$

Now by raising to the power $n+1$ we obtain, as desired:
$$\left(1+\frac{1}{n+1}\right)^{n+1}\geq\left(1+\frac{1}{n}\right)^n$$

(2) Normally we are left with proving that $x_n$ is bounded from above, but this is non-trivial, and we have to use a trick. Consider the following sequence:
$$y_n=\left(1+\frac{1}{n}\right)^{n+1}$$ 

We will prove that this sequence $y_n$ is decreasing, and together with the fact that we have $x_n/y_n\to1$, this will give the result. So, this will be our plan. 

\medskip

(3) In order to prove now that $y_n$ is decreasing, we use, a bit as before:
$$\frac{1+\sum_{i=1}^n\left(1-\frac{1}{n}\right)}{n+1}\geq\sqrt[n+1]{1\cdot\prod_{i=1}^n\left(1-\frac{1}{n}\right)}$$

In practice, this gives the following inequality:
$$1-\frac{1}{n+1}\geq\left(1-\frac{1}{n}\right)^{n/(n+1)}$$

Now by raising to the power $n+1$ we obtain from this:
$$\left(1-\frac{1}{n+1}\right)^{n+1}\geq\left(1-\frac{1}{n}\right)^n$$

And, by inverting this inequality that we found, we obtain, as desired:
$$\left(1+\frac{1}{n}\right)^{n+1}\leq\left(1+\frac{1}{n-1}\right)^n$$

(4) But with this, we can now finish. Indeed, the sequence $x_n$ is increasing, the sequence $y_n$ is decreasing, and we have $x_n<y_n$, as well as:
$$\frac{y_n}{x_n}=1+\frac{1}{n}\to1$$

Thus, both sequences $x_n,y_n$ converge to a certain number $e$, as desired. 

\medskip

(5) Finally, regarding the numerics for our limiting number $e$, we know from the above that we have $x_n<e<y_n$ for any $n\in\mathbb N$, which reads:
$$\left(1+\frac{1}{n}\right)^n<e<\left(1+\frac{1}{n}\right)^{n+1}$$

Thus $e\in[2,3]$, and with a bit of patience, or a computer, we obtain $e=2.71828\ldots$ We will actually come back to this question later, with better methods.
\end{proof}

We should mention that there are many other ways of getting into $e$. For instance it is possible to prove that we have the following interesting formula:
$$\sum_{n=0}^\infty\frac{1}{n!}=e$$

Importantly, this is not the end of the story with $e$, because we have as well:
$$\lim_{n\to\infty}\left(1+\frac{x}{n}\right)^n=\sum_{n=0}^\infty\frac{x^n}{n!}=e^x$$

However, all this is quite tricky, requiring a good knowledge of real functions. And, good news, real functions will be what we will be doing in chapters 2-4 below.

\section*{1e. Exercises}

This opening chapter was a bit special, containing a lot of material in need to be known, and compacted to the maximum. As exercises, again compacted, we have:

\begin{exercise}
Prove by recurrence on $a$, using $(a+1)^p=\sum_{k=0}^p\binom{p}{k}a^k$, the little Fermat theorem, stating that $a^p=a(p)$, for $p$ prime.
\end{exercise}

\begin{exercise}
Find geometric proofs, using triangles in the plane, for the well-known formulae for $\sin(x+y)$ and $\cos(x+y)$.
\end{exercise}

\begin{exercise}
Develop some convergence theory for $x_n=a^n$ with $a>0$, notably by proving that $a^n/n^k\to\infty$ for any $a>1$, and any $k\in\mathbb N$.
\end{exercise}

\begin{exercise}
Prove that $\sum_{n=0}^\infty\frac{1}{n!}=e$. Also, prove that $\left(1+\frac{x}{n}\right)^n\to e^x$, and that $\sum_{n=0}^\infty\frac{x^n}{n!}=e^x$, for $x=-1$, then for $x\in\mathbb Z$, then for $x\in\mathbb R$.
\end{exercise}

These exercises are probably quite difficult, unless you are already a bit familiar with all this. If this is not the case, a good idea at this point is to pick a random entry-level calculus book, and work out a few dozen exercises from there, as a warm-up.

\chapter{Functions, continuity}

\section*{2a. Continuous functions}

We are now ready to talk about functions, which are the main topic of this book. A function $f:\mathbb R\to\mathbb R$ is a correspondence $x\to f(x)$, which to each real number $x\in\mathbb R$ associates a real number $f(x)\in\mathbb R$. As examples, we have $f(x)=x^2$, $f(x)=2^x$ and so on. This suggests that any function $f:\mathbb R\to\mathbb R$ should be given by some kind of ``mathematical formula'', but unfortunately this is not correct, because, with suitable definitions of course, there are more functions than mathematical formulae.

\bigskip

However, we will see that under suitable regularity assumptions on $f:\mathbb R\to\mathbb R$, we have indeed a mathematical formula for $f(x)$ in terms of $x$, at least approximately, and locally. And with this being actually the main idea of calculus, that will take some time to be developed. But more on this later, once we will know more about functions.

\bigskip

Getting started now, let us keep from the above discussion the idea that we should focus our study on the functions $f:\mathbb R\to\mathbb R$ having suitable regularity properties. In what regards these regularity properties, the most basic of them is continuity:

\index{function}
\index{continuous function}

\begin{definition}
A function $f:\mathbb R\to\mathbb R$, or more generally $f:X\to\mathbb R$, with $X\subset\mathbb R$ being a subset, is called continuous when, for any $x_n,x\in X$:
$$x_n\to x\implies f(x_n)\to f(x)$$
Also, we say that $f:X\to\mathbb R$ is continuous at a given point $x\in X$ when the above condition is satisfied, for that point $x$.
\end{definition}

Observe that a function $f:X\to\mathbb R$ is continuous precisely when it is continuous at any point $x\in X$. We will see examples in a moment. Still speaking theory, there are many equivalent formulations of the notion of continuity, with a well-known one, coming by reminding in the above definition what convergence of a sequence means, twice, for both the convergences $x_n\to x$ and $f(x_n)\to f(x)$, being as follows: 
$$\forall x\in X,\forall\varepsilon>0,\exists\delta>0,|x-y|<\delta\implies|f(x)-f(y)|<\varepsilon$$

At the level of examples, basically all the functions that you know, including powers $x^a$, exponentials $a^x$, and more advanced functions like $\sin,\cos,\exp,\log$, are continuous. However, proving this will take some time. Let us start with:

\begin{theorem}
If $f,g$ are continuous, then so are:
\begin{enumerate}
\item $f+g$.

\item $fg$.

\item $f/g$.

\item $f\circ g$.
\end{enumerate}
\end{theorem}

\begin{proof}
Before anything, we should mention that the claim is that (1-4) hold indeed, provided that at the level of domains and ranges, the statement makes sense. For instance in (1,2,3) we are talking about functions having the same domain, and with $g(x)\neq0$ for the needs of (3), and there is a similar discussion regarding (4).

\medskip

(1) The claim here is that if both $f,g$ are continuous at a point $x$, then so is the sum $f+g$. But this is clear from the similar result for sequences, namely:
$$\lim_{n\to\infty}(x_n+y_n)=\lim_{n\to\infty}x_n+\lim_{n\to\infty}y_n$$

(2) Again, the statement here is similar, and the result follows from:
$$\lim_{n\to\infty}x_ny_n=\lim_{n\to\infty}x_n\lim_{n\to\infty}y_n$$

(3) Here the claim is that if both $f,g$ are continuous at $x$, with $g(x)\neq0$, then $f/g$ is continuous at $x$. In order to prove this, observe that by continuity, $g(x)\neq0$ shows that $g(y)\neq0$ for $|x-y|$ small enough. Thus we can assume $g\neq0$, and with this assumption made, the result follows from the similar result for sequences, namely:
$$\lim_{n\to\infty}x_n/y_n=\lim_{n\to\infty}x_n/\lim_{n\to\infty}y_n$$

(4) Here the claim is that if $g$ is continuous at $x$, and $f$ is continuous at $g(x)$, then $f\circ g$ is continuous at $x$. But this is clear, coming from:
\begin{eqnarray*}
x_n\to x
&\implies&g(x_n)\to g(x)\\
&\implies&f(g(x_n))\to f(g(x))
\end{eqnarray*}

Alternatively, let us prove this as well by using that scary $\varepsilon,\delta$ condition given after Definition 2.1. So, let us pick $\varepsilon>0$. We want in the end to have something of type $|f(g(x))-f(g(y))|<\varepsilon$, so we must first use that $\varepsilon,\delta$ condition for the function $f$. So, let us start in this way. Since $f$ is continuous at $g(x)$, we can find $\delta>0$ such that:
$$|g(x)-z|<\delta\implies|f(g(x))-f(z)|<\varepsilon$$

On the other hand, since $g$ is continuous at $x$, we can find $\gamma>0$ such that:
$$|x-y|<\gamma\implies|g(x)-g(y)|<\delta$$

Now by combining the above two inequalities, with $z=g(y)$, we obtain:
$$|x-y|<\gamma\implies|f(g(x))-f(g(y))|<\varepsilon$$

Thus, the composition $f\circ g$ is continuous at $x$, as desired.
\end{proof}

At the level of concrete examples of continuous functions, we have:

\index{sin}
\index{cos}

\begin{theorem}
The following functions are continuous:
\begin{enumerate}
\item $x^n$, with $n\in\mathbb Z$.

\item $P(x)$, with $P\in\mathbb R[X]$.

\item $P(x)/Q(x)$, with $P,Q\in\mathbb R[X]$.

\item $\sin x$, $\cos x$, $\tan x$.

\item $\sec x$, $\csc x$, $\cot x$.
\end{enumerate}
\end{theorem}

\begin{proof}
This is a mixture of trivial and non-trivial results, as follows:

\medskip

(1) Since $f(x)=x$ is continuous, by using Theorem 2.2 we obtain the result for exponents $n\in\mathbb N$, and then for general exponents $n\in\mathbb Z$ too.

\medskip

(2) This follows from (1) with $n\in\mathbb N$, and from the trivial fact, that I forgot to mention in Theorem 2.2, that if $f$ is continuous, then so is $\lambda f$, for any $\lambda\in\mathbb R$.

\medskip

(3) The statement here, which generalizes (1,2), follows exactly as (1,2), by using the various findings from Theorem 2.2, plus the extra rule regarding $f\to\lambda f$.

\medskip

(4) We must first prove here that $x_n\to x$ implies $\sin x_n\to\sin x$, which in practice amounts in proving that $\sin(x+y)\simeq\sin x$ for $y$ small. But this follows from:
$$\sin(x+y)=\sin x\cos y+\cos x\sin y$$ 

To be more precise, let us first establish this formula. In order to do so, consider the following picture, consisting of a length 1 line segment, with angles $x,y$ drawn on each side, and with everything being completed, and lengths computed, as indicated:
$$
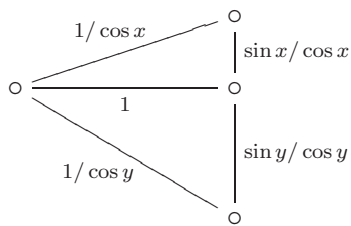
\xymatrix@R=15pt@C=70pt{
&\circ\ar@{-}[d]^{\sin x/\cos x}\\
\circ\ar@{-}[ur]^{1/\cos x}\ar@{-}[r]_1\ar@{-}[ddr]_{1/\cos y}&\circ\ar@{-}[dd]^{\sin y/\cos y}\\
\\
&\circ
}$$

Now let us compute the area of the big triangle, or rather the double of that area. We can do this in two ways, either directly, with a formula involving $\sin(x+y)$, or by using the two small triangles, involving functions of $x,y$. We obtain in this way:
$$\frac{1}{\cos x}\cdot\frac{1}{\cos y}\cdot\sin(x+y)=\frac{\sin x}{\cos x}\cdot 1+\frac{\sin y}{\cos y}\cdot 1$$

But this gives the formula for $\sin(x+y)$ claimed above.

\medskip

(5) Now with this formula in hand, we can establish the continuity of $\sin x$, as follows, with the limits at 0 which are used being both clear on pictures:
\begin{eqnarray*}
\lim_{y\to0}\sin(x+y)
&=&\lim_{y\to0}\left(\sin x\cos y+\cos x\sin y\right)\\
&=&\sin x\lim_{y\to0}\cos y+\cos x\lim_{y\to0}\sin y\\
&=&\sin x\cdot1+\cos x\cdot0\\
&=&\sin x
\end{eqnarray*} 

(6) Moving ahead now with $\cos x$, here the continuity follows from the continuity of $\sin x$, by using the following formula, which is obvious from definitions:
$$\cos x=\sin\left(\frac{\pi}{2}-x\right)$$

(7) Alternatively, and let us do this because we will need later the formula, by using the formula for $\sin(x+y)$ we can deduce a formula for $\cos(x+y)$, as follows:
\begin{eqnarray*}
\cos(x+y)
&=&\sin\left(\frac{\pi}{2}-x-y\right)\\
&=&\sin\left[\left(\frac{\pi}{2}-x\right)+(-y)\right]\\
&=&\sin\left(\frac{\pi}{2}-x\right)\cos(-y)+\cos\left(\frac{\pi}{2}-x\right)\sin(-y)\\
&=&\cos x\cos y-\sin x\sin y
\end{eqnarray*}

(8) But with this, we can use the same method as in (5), and we get, as desired:
\begin{eqnarray*}
\lim_{y\to0}\cos(x+y)
&=&\lim_{y\to0}\left(\cos x\cos y-\sin x\sin y\right)\\
&=&\cos x\lim_{y\to0}\cos y-\sin x\lim_{y\to0}\sin y\\
&=&\cos x\cdot1-\sin x\cdot0\\
&=&\cos x
\end{eqnarray*} 

(9) Regarding $\tan x=\sin x/\cos x$, its continuity follows from that of $\sin x,\cos x$, or via a computation as in (5,8), based on the following formula, coming from (4,7):
$$\tan(x+y)=\frac{\tan x+\tan y}{1-\tan x\tan y}$$

Finally, we can talk about the secondary trigonometric functions, as follows:
$$\sec x=\frac{1}{\cos x}\quad,\quad \csc x=\frac{1}{\sin x}\quad,\quad \cot x=\frac{1}{\tan x}$$

And these are continuous on their domain, again by using Theorem 2.2.
\end{proof}

We will be back to more examples later, and in particular to functions of type $x^a$ and $a^x$ with $a\in\mathbb R$, which are more tricky to define. Also, we will talk as well about inverse functions $f^{-1}$, with as particular cases the basic inverse trigonometric functions.

\bigskip

Going ahead with more theory, when looking at continuous functions with a good miscroscope, some appear to be ``more continuous'' than other, as shown by:

\begin{theorem}
Consider the following properties, regarding $f:X\to\mathbb R$ with $X\subset\mathbb R$:
\begin{enumerate}
\item $f$ has the following property, with $K>0$, called Lipschitz property:
$$|f(x)-f(y)|\leq K|x-y|$$

\item $f$ is uniformly continuous, in the sense that the following happens:
$$\forall\varepsilon>0,\exists\delta>0,|x-y|<\delta\implies|f(x)-f(y)|<\varepsilon$$

\item $f$ is continuous in the usual sense, namely:
$$\forall x\in X,\forall\varepsilon>0,\exists\delta>0,|x-y|<\delta\implies|f(x)-f(y)|<\varepsilon$$
\end{enumerate}
We have then $(1)\implies(2)\implies(3)$. Also, the converse implications do not hold.
\end{theorem}

\begin{proof}
This is something quite self-explanatory, the idea being as follows:

\medskip

$(1)\implies(2)$ This is clear, coming by talking $\delta=\varepsilon/K$.

\medskip

$(2)\implies(3)$ This is something which is plainly trivial. 

\medskip

$(3)\ \,\not\!\!\!\!\implies\!(2)$ The simplest counterexample here is $f(x)=x^2$. Indeed, this function is continuous, and its uniform continuity property, applied with $\varepsilon=1$, would lead to the existence of $\delta>0$ such that the following happens, which is wrong:
$$|x-y|<\delta\implies|x^2-y^2|<1$$

$(2)\ \,\not\!\!\!\!\implies\!(1)$ The simplest counterexample here is $f(x)=\sqrt{x}$. Indeed, observe first that we have the following estimate, valid for any $x>y>0$:
$$(\sqrt{x}-\sqrt{y})^2\leq(\sqrt{x}-\sqrt{y})(\sqrt{x}+\sqrt{y})=x-y$$

Thus our function is indeed uniformly continuous, and this because we have:
$$|x-y|<\varepsilon^2\implies|\sqrt{x}-\sqrt{y}|<\varepsilon$$

In what regards now the Lipschitz property, observe that we have:
$$\frac{\sqrt{x}-\sqrt{y}}{x-y}=\frac{1}{\sqrt{x}+\sqrt{y}}$$

Now since this can be arbitrarily big, when $x,y$ are small, Lipschitz fails indeed.
\end{proof}

Generally speaking, the usual functions are quite often Lipschitz, on suitable domains, and up to you to do some computations here. As a main application of the Lipschitz property, we have the following result, mixing old and modern mathematics:

\begin{theorem}
The following happen, in relation with iterations:
\begin{enumerate}
\item A function $f:X\to X$ with $X\subset\mathbb R$ closed, which is a contraction, meaning Lipschitz with $K<1$, has a unique fixed point, obtained by iterating $f$.

\item In particular, we can extract the square root of $a>0$ by iterating the function $f(x)=(x+a/x)/2$, with this being known as the Babylonian method.
\end{enumerate}
\end{theorem}

\begin{proof}
As mentioned, this is mixture of old and modern mathematics, as follows:

\medskip

(1) Consider indeed a function $f:X\to X$, then pick $x\in X$, and start iterating $f$:
$$x_0=x\quad,\quad x_{n+1}=f(x_n)$$

Now when assuming that $f$ is a contraction, $\{x_n\}_{n\in\mathbb N}$ will be Cauchy, due to:
$$|x_{n+s}-x_n|\leq\sum_{t=n}^{n+s-1}|x_{t+1}-x_t|\leq\sum_{t=n}^{n+s-1}K^t|x_1-x_0|\leq\frac{K^n}{1-K}\,|x_1-x_0|$$

Thus we have a limit $x_n\to z\in\mathbb R$, and with $X$ assumed to be closed, $z\in X$. But with this we are done, because $x_{n+1}=f(x_n)$ gives in the limit $z=f(z)$. As for the uniqueness of the fixed point $z$, this is plainly clear, from the contraction property of $f$. 

\medskip

(2) This is something quite self-explanatory, based on the following computation:
$$\frac{x+a/x}{2}=x\iff a/x=x\iff x=\sqrt{a}$$

So, let us see if (1) applies. Assuming $a>1$, the graph of our function is as follows:
$$\xymatrix@R=4pt@C=4pt{
&\\
\\
&a\ar[uu]\ar@{.}[rrrrr]&&&&&&\\
&\\
&\frac{a+1}{2}\ar@{-}[uu]\ar@{.}[rrrrr]&\bullet&&&&\bullet\\
&\sqrt{a}\ar@{-}[u]\ar@{.}[rrrrr]&&\bullet\ar@{-}@/^/[ul]\ar@{-}@/_/[urrr]&&&\\
&1\ar@{-}[u]\\
\ar@{-}[rr]&\ar@{.}[uuuuurrrrr]&1\ar@{.}[uuu]\ar@{-}[r]&\sqrt{a}\ar@{.}[uu]\ar@{-}[r]&\frac{a+1}{2}\ar@{-}[rr]\ar@{.}[uuu]&&a\ar[rr]\ar@{.}[uuuuu]&&\\
&\ar@{-}[uu]}$$

Which makes it quite clear that we will have to work a bit, in order to see if (1) applies indeed. But here, I did all the computations for us, with the following conclusions:

\smallskip

-- For $a<3$ we have a contraction on $[1,a]$, with constant $K=(a-1)/2$.

\smallskip

-- For $a>3$ we have a contraction on $[\sqrt{a},a/3]$, with constant $K=(9/a-1)/2$.

\smallskip

And up to you to fill in the details, in relation with all this, and don't forget to work out some numerics too, say for $a=2$, I mean if the ancient Babylonians, or an old man like me, were able to fully do the math here, so should you, young reader. 
\end{proof}

In relation now with the notion of uniform continuity, which is something more abstract, we have the following remarkable result, due to Heine and Cantor:

\begin{theorem}
Any continuous function defined on a closed, bounded interval
$$f:[a,b]\to\mathbb R$$
is automatically uniformly continuous.
\end{theorem}

\begin{proof}
This is something quite subtle, and we are punching here a bit above our weight, but here is the proof, with everything or almost included:

\medskip

(1) Given $\varepsilon>0$, for any $x\in[a,b]$ we know that we have a $\delta_x>0$ such that:
$$|x-y|<\delta_x\implies|f(x)-f(y)|<\frac{\varepsilon}{2}$$

So, consider the following open intervals, centered at the various points $x\in[a,b]$:
$$I_x=\left(x-\frac{\delta_x}{2}\,,\,x+\frac{\delta_x}{2}\right)$$

These intervals then obviously cover $[a,b]$, in the sense that we have:
$$[a,b]\subset\bigcup_{x\in[a,b]}I_x$$

Now assume that we managed to prove that this cover has a finite subcover. Then we can most likely choose our $\delta>0$ to be the smallest of the $\delta_x>0$ involved, or perhaps half of that, and then get our uniform continuity condition, via the triangle inequality.

\medskip

(2) So, let us prove that the cover in (1) has a finite subcover. For this purpose, we proceed by contradiction. So, assume that $[a,b]$ has no finite subcover, and let us cut this interval in half. Then one of the halves must have no finite subcover either, and we can repeat the procedure, by cutting this smaller interval in half. And so on. But this leads to a contradiction, because the limiting point $x\in[a,b]$ that we obtain in this way, as the intersection of these smaller and smaller intervals, must be covered by something, and so one of these small intervals leading to it must be covered too, contradiction. 

\medskip

(3) With this done, we are ready to finish, as announced in (1). Indeed, let us denote by $[a,b]\subset\bigcup_iI_{x_i}$ the finite subcover found in (2), and let us set:
$$\delta=\min_i\frac{\delta_{x_i}}{2}$$

Now assume $|x-y|<\delta$, and pick $i$ such that $x\in I_{x_i}$. By the triangle inequality we have then $|x_i-y|<\delta_{x_i}$, which shows that we have $y\in I_{x_i}$ as well. But by applying now $f$, this gives as desired $|f(x)-f(y)|<\varepsilon$, again via the triangle inequality.
\end{proof}

\section*{2b. Intermediate values}

Moving ahead with more theory, we would like to explain now an alternative formulation of the notion of continuity, which is quite abstract, and a bit difficult to understand and master when you are a beginner, but which is definitely worth learning, because it is quite powerful, solving some of the questions that we have left. Let us start with:

\index{open set}
\index{closed set}

\begin{definition}
The open and closed sets are defined as follows:
\begin{enumerate}
\item Open means that there is a small interval around each point.

\item Closed means that our set is closed under taking limits.
\end{enumerate}
\end{definition} 

As basic examples, the open intervals $(a,b)$ are open, and the closed intervals $[a,b]$ are closed. Observe also that $\mathbb R$ itself is open and closed at the same time. Further examples, or rather results which are easy to establish, include the fact that the finite unions or intersections of open or closed sets are open or closed. We will be back to all this later, with some precise results in this sense. For the moment, we will only need:

\begin{proposition}
A set $O\subset\mathbb R$ is open precisely when its complement $C\subset\mathbb R$ is closed, and vice versa.
\end{proposition}

\begin{proof}
It is enough to prove the first assertion, since the ``vice versa'' part will follow from it, by taking complements. But this can be done as follows:

\medskip

``$\implies$'' Assume that $O\subset\mathbb R$ is open, and let $C=\mathbb R-O$. In order to prove that $C$ is closed, assume that $\{x_n\}_{n\in\mathbb N}\subset C$ converges to $x\in\mathbb R$. We must prove that $x\in C$, and we will do this by contradiction. So, assume  $x\notin C$. Thus $x\in O$, and since $O$ is open we can find a small interval $(x-\varepsilon,x+\varepsilon)\subset O$. But since $x_n\to x$ this shows that $x_n\in O$ for $n$ big enough, which contradicts $x_n\in C$ for all $n$, and we are done.

\medskip

``$\Longleftarrow$'' Assume that $C\subset\mathbb R$ is open, and let $O=\mathbb R-C$. In order to prove that $O$ is open, let $x\in O$, and consider the intervals $(x-1/n,x+1/n)$, with $n\in\mathbb N$. If one of these intervals lies in $O$, we are done. Otherwise, this would mean that for any $n\in\mathbb N$ we have at least one point $x_n\in(x-1/n,x+1/n)$ satisfying $x_n\notin O$, and so $x_n\in C$. But since $C$ is closed and $x_n\to x$, we get $x\in C$, and so $x\notin O$, contradiction, and we are done. 
\end{proof}

As basic illustrations for the above result, $\mathbb R-(a,b)=(-\infty,a]\cup[b,\infty)$ is closed, and $\mathbb R-[a,b]=(-\infty,a)\cup(b,\infty)$ is open. Getting now back to functions, we have:

\index{continuous function}

\begin{theorem}
A function is continuous precisely when $f^{-1}(O)$ is open, for any $O$ open. Equivalently, $f^{-1}(C)$ must be closed, for any $C$ closed.
\end{theorem}

\begin{proof}
Here the first assertion follows from definitions, and more specifically from the $\varepsilon,\delta$ definition of continuity, which was as follows:
$$\forall x\in X,\forall\varepsilon>0,\exists\delta>0,|x-y|<\delta\implies|f(x)-f(y)|<\varepsilon$$

Indeed, if $f$ satisfies this condition, it is clear that if $O$ is open, then $f^{-1}(O)$ is open, and the converse holds too. As for the second assertion, this can be proved either directly, by using the $f(x_n)\to f(x)$ definition of continuity, or by taking complements.
\end{proof}

As a test for the above criterion, let us reprove the fact, that we know from Theorem 2.2, that if $f,g$ are continuous, so is $f\circ g$. But this is clear, coming from:
$$(f\circ g)^{-1}(O)=g^{-1}(f^{-1}(O))$$

Not bad. In order to reach now to true applications of Theorem 2.9, we need to know more about the open and closed sets. Let us begin with a useful result, as follows:

\index{open set}
\index{closed set}

\begin{proposition}
The following happen:
\begin{enumerate}
\item Union of open sets is open.

\item Intersection of closed sets is closed.

\item Finite intersection of open sets is open.

\item Finite union of closed sets is closed.
\end{enumerate}
\end{proposition}

\begin{proof}
Here (1) is clear from definitions, (3) is clear from definitions too, and (2,4) follow from (1,3) by taking complements $E\to E^c$, using the following formulae:
$$\left(\bigcup_iE_i\right)^c=\bigcap_iE_i^c
\quad,\quad
\left(\bigcap_iE_i\right)^c=\bigcup_iE_i^c$$

Thus, we are led to the conclusions in the statement.  
\end{proof}

As an important comment, (3,4) above do not hold when removing the finiteness assumption, with the simplest counterexamples being as follows:
$$\bigcap_{n\in\mathbb N}\left(-\frac{1}{n}\,,\,\frac{1}{n}\right)=\{0\}
\quad,\quad 
\bigcup_{n\in\mathbb N}\left[0\,,\,1-\frac{1}{n}\right]=[0,1)$$

All this is quite interesting, and leads us to the question about what the open and closed sets really are. And fortunately, this question can be answered, as follows:

\index{union of intervals}

\begin{theorem}
The open and closed sets are as follows:
\begin{enumerate}
\item The open sets are the disjoint unions of open intervals.

\item The closed sets are the complements of these unions.
\end{enumerate}
\end{theorem}

\begin{proof}
We have two assertions to be proved, the idea being as follows:

\medskip

(1) We know that the open intervals are those of type $(a,b)$ with $a<b$, with the values $a,b=\pm\infty$ allowed, and by Proposition 2.10 a union of such intervals is open. 

\medskip

(2) Conversely, given $O\subset\mathbb R$ open, we can cover each point $x\in O$ with an open interval $I_x\subset O$, and we have $O=\cup_xI_x$, so $O$ is a union of open intervals. 

\medskip

(3) In order to finish the proof of the first assertion, it remains to prove that the union $O=\cup_xI_x$ in (2) can be taken to be disjoint. For this purpose, our first observation is that, by approximating points $x\in O$ by rationals $y\in\mathbb Q\cap O$, we can make our union to be countable. But once our union is countable, we can start merging intervals, whenever they meet, and we are left in the end with a countable, disjoint union, as desired.

\medskip

(4) Finally, the second assertion comes from Proposition 2.8.
\end{proof}

Moving towards more concrete things, and applications, let us formulate:

\index{compact set}
\index{connected set}
\index{cover}
\index{subcover}

\begin{definition}
The compact and connected sets are defined as follows:
\begin{enumerate}
\item A subset $K\subset\mathbb R$ is called compact when any open cover of it, $K\subset\cup_xO_x$, with all $O_x\subset\mathbb R$ being open sets, has a finite subcover.

\item A subset $E\subset\mathbb R$ is called connected when it cannot be broken into two parts, in the sense that we cannot have $E\subset A\sqcup B$, with $A,B\neq\emptyset$ open.
\end{enumerate}
\end{definition}

As basic examples, the closed bounded intervals $[a,b]$ are compact, as we know from the proof of Theorem 2.6, and so are the finite unions of such intervals. As for connected sets, the basic examples here are the various types of intervals, namely $(a,b)$, $(a,b]$, $[a,b)$, $[a,b]$, and it looks impossible to come up with more examples. In fact, we have:

\index{closed and bounded}

\begin{theorem}
The compact and connected sets are as follows:
\begin{enumerate}
\item The compact sets are those which are closed and bounded.

\item The connected sets are the various types of intervals.
\end{enumerate}
\end{theorem}

\begin{proof}
This is something quite intuitive, the idea being as follows:

\medskip

(1) Compact implies closed, because assuming the contrary, given $\{x_n\}\subset K$ converging to $x\notin K$, we have the following open cover, having no finite subcover:
$$K\subset\bigcup_{n\in\mathbb N}\left[x-\frac{1}{n},x+\frac{1}{n}\right]^c$$

(2) Similarly, compact implies bounded, and this due to the following open cover:
$$K\subset\bigcup_{n\in\mathbb N}(-n,n)$$

(3) As for the converse, stating that closed and bounded implies compact, assume that we have a bounded set, $K\subset[a,b]$, which is closed, and consider an open cover of it:
$$K\subset\bigcup_xO_x$$

By adding on the right $K^c$, we obtain an open cover of $[a,b]$. But, we know from the proof of Theorem 2.6 that $[a,b]$ is compact, so this cover must have a finite subcover. And by removing $K^c$ from this latter cover, we obtain a finite cover of $K$, as desired.

\medskip

(4) Finally, regarding the second assertion, this is something quite obvious, because $E\subset\mathbb R$ being connected means $a,b\in E\implies[a,b]\subset E$, and this gives the result.
\end{proof}

With this discussed, let us go back to the continuous functions. We have:

\begin{theorem}
Assuming that $f$ is continuous:
\begin{enumerate}
\item If $K$ is compact, then $f(K)$ is compact.

\item If $E$ is connected, then $f(E)$ is connected.
\end{enumerate}
\end{theorem}

\begin{proof}
This is something very standard, the idea being as follows:

\medskip

(1) Given an open cover $f(K)\subset\cup_xO_x$ we have $K\subset\cup_xf^{-1}(O_x)$, open cover too, and if $K\subset\cup_yf^{-1}(O_y)$ is a finite subcover of this, then $f(K)\subset\cup_yO_y$, as desired.

\medskip

(2) This is clear too, because assuming $f(E)\subset A\sqcup B$ with $A,B\neq\emptyset$ open, we have $E\subset f^{-1}(A)\sqcup f^{-1}(B)$ with $f^{-1}(A),f^{-1}(B)\neq\emptyset$ open, contradiction.
\end{proof}

Let us record as well the following useful generalization of Theorem 2.6:

\begin{theorem}
Any continuous function defined on a compact set $f:X\to\mathbb R$ is automatically uniformly continuous.
\end{theorem}

\begin{proof}
We can prove this exactly as Theorem 2.6, by using the compactness of $X$.
\end{proof}

Time for some applications? Here is a beautiful theorem, coming from this:

\index{intermediate value}
\index{maximum}
\index{minimum}

\begin{theorem}
The following happen for a continuous function $f:[a,b]\to\mathbb R$:
\begin{enumerate}
\item Its image is a compact interval, $Im(f)=[c,d]$.

\item $f$ takes all intermediate values between $f(a),f(b)$.

\item $f$ has a minimum and maximum on $[a,b]$.

\item If $f(a),f(b)$ have different signs, $f(x)=0$ has a solution.
\end{enumerate}
\end{theorem}

\begin{proof}
All these statements are related, and are called altogether ``intermediate value theorem''. Regarding now the proof, (1), which gives everything, says that:
$$X={\rm compact,\,connected}\ \implies\ f(X)={\rm compact,\,connected}$$

But this follows from Theorem 2.14, making the whole thing trivial. Good.
\end{proof}

Along the same lines, we have as well the following result:

\index{monotone function}

\begin{theorem}
A continuous surjective function $f$ is injective, and so invertible, precisely when it is monotone, and in this case, the inverse function $f^{-1}$ must be monotone and continuous too. Moreover, this statement holds both locally, and globally.
\end{theorem}

\begin{proof}
The first assertion follows from Theorem 2.16, and the fact that $f^{-1}$ is monotone is clear. Regarding now the continuity of $f^{-1}$, we want to prove that we have:
$$x_n\to x\implies f^{-1}(x_n)\to f^{-1}(x)$$

But with $x_n=f(y_n)$ and $x=f(y)$, this condition becomes:
$$f(y_n)\to f(y)\implies y_n\to y$$

And this latter condition being true since $f$ is monotone, we are done.
\end{proof}

As a basic application now of Theorem 2.17, we have:

\index{inverse function}

\begin{proposition}
The various usual inverse functions, such as the inverse trigonometric functions $\arcsin$, $\arccos$, $\arctan$, ${\rm arcsec}$, ${\rm arccsc}$, ${\rm arccot}$, are all continuous.
\end{proposition}

\begin{proof}
This is indeed a very standard application of Theorem 2.17.
\end{proof}

As another basic application of our technology, we have:

\index{roots of polynomial}
\index{square root}

\begin{proposition}
The following happen:
\begin{enumerate}
\item Any polynomial $P\in\mathbb R[X]$ of odd degree has a root.

\item Given $n\in2\mathbb N+1$, we can extract $\sqrt[n]{x}$, for any $x\in\mathbb R$.

\item Given $n\in\mathbb N$, we can extract $\sqrt[n]{x}$, for any $x\in[0,\infty)$.
\end{enumerate}
\end{proposition}

\begin{proof}
All these results come as applications of Theorem 2.16, as follows:

\medskip

(1) This is clear from Theorem 2.16 (3), applied on $[-\infty,\infty]$.

\medskip

(2) This follows from (1), by using the polynomial $P(z)=z^n-x$.

\medskip 

(3) This follows as well by applying Theorem 2.16 (3) to the polynomial $P(z)=z^n-x$, but this time on $[0,\infty)$.
\end{proof}

As a concrete application, in relation with powers, we have the following result, completing our series of results regarding the basic mathematical functions:

\index{power function}

\begin{theorem}
The function $x^a$ is defined and continuous on $(0,\infty)$, for any $a\in\mathbb R$. Moreover, when trying to extend it to $\mathbb R$, we have $4$ cases, as follows,
\begin{enumerate}
\item For $a\in\mathbb Q_{odd}$, $a>0$, the maximal domain is $\mathbb R$.

\item For $a\in\mathbb Q_{odd}$, $a\leq0$, the maximal domain is $\mathbb R-\{0\}$.

\item For $a\in\mathbb R-\mathbb Q$ or $a\in\mathbb Q_{even}$, $a>0$, the maximal domain is $[0,\infty)$.

\item For $a\in\mathbb R-\mathbb Q$ or $a\in\mathbb Q_{even}$, $a\leq0$, the maximal domain is $(0,\infty)$.
\end{enumerate}
where $\mathbb Q_{odd}$ is the set of rationals $r=p/q$ with $q$ odd, and $\mathbb Q_{even}=\mathbb Q-\mathbb Q_{odd}$.
\end{theorem}

\begin{proof}
The idea is that we know how to extract roots by using Proposition 2.19, and all the rest follows by continuity. To be more precise:

\medskip

(1) Assume $a=p/q$, with $p,q\in\mathbb N$, $p\neq0$ and $q$ odd. Given a number $x\in\mathbb R$, we can construct the power $x^a$ in the following way, by using Proposition 2.19:
$$x^a=\sqrt[q]{x^p}$$

Then, it is straightforward to prove that $x^a$ is indeed continuous on $\mathbb R$.

\medskip

(2) In the case $a=-p/q$, with $p,q\in\mathbb N$ and $q$ odd, the same discussion applies, with the only change coming from the fact that $x^a$ cannot be applied to $x=0$. 

\medskip

(3) Assume first $a\in\mathbb Q_{even}$, $a>0$. This means $a=p/q$ with $p,q\in\mathbb N$, $p\neq0$ and $q$ even, and as before in (1), we can set $x^a=\sqrt[q]{x^p}$ for $x\geq0$, by using Proposition 2.19. It is then straightforward to prove that $x^a$ is indeed continuous on $[0,\infty)$, and not extendable either to the negatives. Thus, we are done with the case $a\in\mathbb Q_{even}$, $a>0$, and the case left, namely $a\in\mathbb R-\mathbb Q$, $a>0$, follows as well by continuity.

\medskip

(4) In the cases $a\in\mathbb Q_{even}$, $a\leq0$ and $a\in\mathbb R-\mathbb Q$, $a\leq0$, the same discussion applies, with the only change coming from the fact that $x^a$ cannot be applied to $x=0$. 
\end{proof}

Let us record as well a result about the function $a^x$, as follows:

\begin{theorem}
The function $a^x$ is as follows:
\begin{enumerate}
\item For $a>0$, this function is defined and continuous on $\mathbb R$.

\item For $a=0$, this function is defined and continuous on $(0,\infty)$.

\item For $a<0$, the domain of this function contains no interval.
\end{enumerate}
\end{theorem}

\begin{proof}
This is a sort of reformulation of Theorem 2.20, by exchanging the variables, $x\leftrightarrow a$. To be more precise, the situation is as follows:

\medskip

(1) We know from Theorem 2.20 that things fine with $x^a$ for $x>0$, no matter what $a\in\mathbb R$ is. But this means that things fine with $a^x$ for $a>0$, no matter what $x\in\mathbb R$ is.

\medskip

(2) This is something trivial, and we have of course $0^x=0$, for any $x>0$. As for the powers $0^x$ with $x\leq0$, these are impossible to define, for obvious reasons.

\medskip

(3) Given $a<0$, we know from Theorem 2.20 that we cannot define $a^x$ for $x\in\mathbb Q_{even}$. But since $\mathbb Q_{even}$ is dense in $\mathbb R$, this gives the result.
\end{proof}

\section*{2c. Sequences and series} 

Our goal now is to extend the material from chapter 1 regarding the numeric sequences and series, to the case of the sequences and series of functions. To start with, we can talk about the convergence of sequences of functions, $f_n\to f$, as follows:

\index{sequence of functions}
\index{pointwise convergence}

\begin{definition}
We say that $f_n$ converges pointwise to $f$, and write $f_n\to f$, if
$$f_n(x)\to f(x)$$
for any $x$. Equivalently, $\forall x,\forall\varepsilon>0,\exists N\in\mathbb N,\forall n\geq N,|f_n(x)-f(x)|<\varepsilon$.
\end{definition}

The question is now, assuming that $f_n$ are continuous, does it follow that $f$ is continuous? I am pretty much sure that you think that the answer is yes, based on:
\begin{eqnarray*}
\lim_{y\to x}f(y)
&=&\lim_{y\to x}\lim_{n\to\infty}f_n(y)\\
&=&\lim_{n\to\infty}\lim_{y\to x}f_n(y)\\
&=&\lim_{n\to\infty}f_n(x)\\
&=&f(x)
\end{eqnarray*}

However, this proof is wrong, because we know well from chapter 1 that we cannot intervert limits, with this being a common beginner mistake. In fact, the result itself is wrong in general, because if we consider the functions $f_n:[0,1]\to\mathbb R$ given by $f_n(x)=x^n$, which are obviously continuous, their limit is discontinuous, given by:
$$\lim_{n\to\infty}x^n
=\begin{cases}
0&,\quad x\in[0,1)\\
1&,\quad x=1
\end{cases}$$

Of course, you might say here that allowing $x=1$ in all this might be a bit unnatural, for whatever reasons, but there is an answer to this too. We can do worse, as follows:

\index{step function}
\index{totally discontinuous}

\begin{proposition}
The basic step function, namely the sign function
$$sgn(x)=\begin{cases}
-1&,\quad x<0\\
0&,\quad x=0\\
1&,\quad x>0
\end{cases}$$
can be approximated by suitable modifications of $\arctan(x)$. Even worse, there are examples of $f_n\to f$ with each $f_n$ continuous, and with $f$ totally discontinuous. 
\end{proposition}

\begin{proof}
To start with, $\arctan(x)$ looks a bit like $sgn(x)$, so to say, but one problem comes from the fact that its image is $[-\pi/2,\pi/2]$, instead of the desired $[-1,1]$. Thus, we must first rescale $\arctan(x)$ by $\pi/2$. Now with this done, we can further stretch the variable $x$, as to get our function closer and closer to $sgn(x)$, as desired. This proves the first assertion, and regarding the second assertion, which is a bit more technical, and that we will not really need in what follows, we will leave this as an exercise for you.
\end{proof}

Sumarizing, we are a bit in trouble, because we would like to have in our bag of theorems something saying that $f_n\to f$ with $f_n$ continuous implies $f$ continuous. Fortunately, this can be done, with a suitable refinement of the notion of convergence, as follows:

\index{uniform convergence}

\begin{definition}
We say that $f_n$ converges uniformly to $f$, and write $f_n\to_uf$, if:
$$\forall\varepsilon>0,\exists N\in\mathbb N,\forall n\geq N,|f_n(x)-f(x)|<\varepsilon,\forall x$$
That is, the same condition as for $f_n\to f$ must be satisfied, but with the $\forall x$ at the end. 
\end{definition}

And it is this ``$\forall x$ at the end'' which makes the difference, and will make our theory work. In order to understand this, which is something quite subtle, let us compare Definition 2.22 and Definition 2.24. As a first observation, we have:

\index{pointwise convergence}

\begin{proposition}
Uniform convergence implies pointwise convergence,
$$f_n\to_uf\implies f_n\to f$$
but the converse is not true, in general.
\end{proposition}

\begin{proof}
Here the first assertion is clear from definitions, just by thinking at what is going on, with no computations needed. As for the second assertion, the simplest counterexamples here are the functions $f_n:[0,1]\to\mathbb R$ given by $f_n(x)=x^n$, that we met before in Proposition 2.23. Indeed, uniform convergence on $[0,1)$ would mean:
$$\forall\varepsilon>0,\exists N\in\mathbb N,\forall n\geq N,x^n<\varepsilon,\forall x\in[0,1)$$

But this is wrong, because no matter how big $N$ is, we have $\lim_{x\to1}x^N=1$, and so we can find $ x\in[0,1)$ such that $x^N>\varepsilon$. Thus, we have our counterexample.
\end{proof}

Moving ahead now, let us state our main theorem on uniform convergence, as follows:

\index{limit of continuous functions}

\begin{theorem}
Assuming that $f_n$ are continuous, and that
$$f_n\to_uf$$
then $f$ is continuous. That is, uniform limit of continuous functions is continuous.
\end{theorem}

\begin{proof}
As previously advertised, it is the ``$\forall x$ at the end'' in Definition 2.24 that will make this work. Indeed, let us try to prove that the limit $f$ is continuous at some point $x$. For this, we pick a number $\varepsilon>0$. Since $f_n\to_uf$, we can find $N\in\mathbb N$ such that:
$$|f_N(z)-f(z)|<\frac{\varepsilon}{3}\quad,\quad\forall z$$

On the other hand, since $f_N$ is continuous at $x$, we can find $\delta>0$ such that:
$$|x-y|<\delta\implies|f_N(x)-f_N(y)|<\frac{\varepsilon}{3}$$

But with this, we are done. Indeed, for $|x-y|<\delta$ we have:
\begin{eqnarray*}
|f(x)-f(y)|
&\leq&|f(x)-f_N(x)|+|f_N(x)-f_N(y)|+|f_N(y)-f(y)|\\
&\leq&\frac{\varepsilon}{3}+\frac{\varepsilon}{3}+\frac{\varepsilon}{3}\\
&=&\varepsilon
\end{eqnarray*}

Thus, the limit function $f$ is continuous at $x$, and we are done.
\end{proof}

Obviously, the notion of uniform convergence in Definition 2.24 is something quite interesting, worth some more study. As a first result, we have:

\begin{proposition}
The following happen, regarding uniform limits:
\begin{enumerate}
\item $f_n\to_uf$, $g_n\to_ug$ imply $f_n+g_n\to_uf+g$.

\item $f_n\to_uf$, $g_n\to_ug$ imply $f_ng_n\to_ufg$.

\item $f_n\to_uf$, $f\neq0$ imply $1/f_n\to_u1/f$.

\item $f_n\to_uf$, $g$ continuous imply $f_n\circ g\to_uf\circ g$.

\item $f_n\to_uf$, $g$ continuous imply $g\circ f_n\to_ug\circ f$.
\end{enumerate}
\end{proposition}

\begin{proof}
All this is routine, exactly as for the results for numeric sequences from chapter 1, that we know well, with no difficulties or tricks involved, exercise for you.
\end{proof}

Finally, there is some abstract mathematics to be done as well. Indeed, observe that the notion of uniform convergence, as formulated in Definition 2.24, means that:
$$\sup_x\big|f_n(x)-f(x)\big|\ \longrightarrow_{n\to\infty}\ 0$$

This suggests measuring the distance between functions via a supremum as above, and in relation with this, we have the following result:

\index{distance}

\begin{theorem}
The uniform convergence, $f_n\to_uf$, means that we have $f_n\to f$ with respect to the following distance,
$$d(f,g)=\sup_x\big|f(x)-g(x)\big|$$
which is indeed a distance function.
\end{theorem}

\begin{proof}
This is indeed something quite self-explanatory, and we will leave some thinking at all this, distance axioms and their verification, as an instructive exercise.
\end{proof}

\section*{2d. Basic functions}

With the above theory in hand, let us get now to interesting things, namely computations. Regarding the power functions $x^a$, we first have the following result:

\index{binomial formula}
\index{generalized binomial formula}
\index{generalized binomial numbers}

\begin{theorem}
We have the generalized binomial formula
$$(1+x)^a=\sum_{k=0}^\infty\binom{a}{k}x^k$$
with the generalized binomial coefficients being given by
$$\binom{a}{k}=\frac{a(a-1)\ldots(a-k+1)}{k!}$$
valid for any exponent $a\in\mathbb Z$, and any $|x|<1$.
\end{theorem}

\begin{proof}
This is something quite tricky, the idea being as follows:

\medskip

(1) For exponents $a\in\mathbb N$, this is something that we know well from chapter 1, and which is valid for any $x\in\mathbb R$, coming from the usual binomial formula, namely:
$$(1+x)^n=\sum_{k=0}^n\binom{n}{k}x^k$$

(2) For the exponent $a=-1$ this is something that we know from chapter 1 too, coming from the following formula, valid for any $|x|<1$:
$$\frac{1}{1+x}=1-x+x^2-x^3+\ldots$$

Indeed, this is exactly our generalized binomial formula at $a=-1$, because:
$$\binom{-1}{k}
=\frac{(-1)(-2)\ldots(-k)}{k!}
=(-1)^k$$

(3) Let us discuss now the general case $a\in-\mathbb N$. With $a=-n$, and $n\in\mathbb N$, the generalized binomial coefficients are given by the following formula:
\begin{eqnarray*}
\binom{-n}{k}
&=&\frac{(-n)(-n-1)\ldots(-n-k+1)}{k!}\\
&=&(-1)^k\frac{n(n+1)\ldots(n+k-1)}{k!}\\
&=&(-1)^k\frac{(n+k-1)!}{(n-1)!k!}\\
&=&(-1)^k\binom{n+k-1}{n-1}
\end{eqnarray*}

Thus, our generalized binomial formula at $a=-n$, and $n\in\mathbb N$, reads:
$$\frac{1}{(1+x)^n}=\sum_{k=0}^\infty(-1)^k\binom{n+k-1}{n-1}x^k$$

(4) In order to prove this formula, it is convenient to write it with $-x$ instead of $x$, in order to get rid of signs. The formula to be proved becomes:
$$\frac{1}{(1-x)^n}=\sum_{k=0}^\infty\binom{n+k-1}{n-1}x^k$$

We prove this by recurrence on $n$. At $n=1$ this formula definitely holds, as explained in (2) above. So, assume that the formula holds at $n\in\mathbb N$. We have then:
\begin{eqnarray*}
\frac{1}{(1-x)^{n+1}}
&=&\frac{1}{1-x}\cdot\frac{1}{(1-x)^n}\\
&=&\sum_{k=0}^\infty x^k\sum_{l=0}^\infty\binom{n+l-1}{n-1}x^l\\
&=&\sum_{s=0}^\infty x^s\sum_{l=0}^s\binom{n+l-1}{n-1}
\end{eqnarray*}

Thus, in order to finish, we must prove the following formula:
$$\sum_{l=0}^s\binom{n+l-1}{n-1}=\binom{n+s}{n}$$

(5) In order to prove this latter formula, we proceed by recurrence on $s\in\mathbb N$. At $s=0$ the formula is trivial, $1=1$. So, assume that the formula holds at $s\in\mathbb N$. In order to prove the formula at $s+1$, we are in need of the following formula:
$$\binom{n+s}{n}+\binom{n+s}{n-1}=\binom{n+s+1}{n}$$

But this is the Pascal formula, that we know from chapter 1, and we are done.
\end{proof}

Quite interestingly, the formula in Theorem 2.29 holds in fact at any $a\in\mathbb R$, but this is something non-trivial, whose proof will have to wait until chapter 3 below. However, in the meantime, let us investigate the case $a=\pm 1/2$. We first have here:

\begin{proposition}
The generalized binomial coefficients at $a=\pm1/2$ are
$$\binom{1/2}{k}=-2\left(\frac{-1}{4}\right)^kC_{k-1}\quad,\quad
\binom{-1/2}{k}=\left(\frac{-1}{4}\right)^kD_k$$
with $D_k=\binom{2k}{k}$ being the central binomial coefficients, and $C_k=\frac{1}{k+1}\binom{2k}{k}$.
\end{proposition}

\begin{proof}
At $a=1/2$, the generalized binomial coefficients are as follows:
\begin{eqnarray*}
\binom{1/2}{k}
&=&\frac{1/2(-1/2)\ldots(3/2-k)}{k!}\\
&=&(-1)^{k-1}\frac{1\cdot 3\cdot 5\ldots(2k-3)}{2^kk!}\\
&=&(-1)^{k-1}\frac{(2k-2)!}{2^{k-1}(k-1)!2^kk!}\\
&=&-2\left(\frac{-1}{4}\right)^kC_{k-1}
\end{eqnarray*}

As for the case $a=-1/2$, here the binomial coefficients are as follows:
\begin{eqnarray*}
\binom{-1/2}{k}
&=&\frac{-1/2(-3/2)\ldots(1/2-k)}{k!}\\
&=&(-1)^k\frac{1\cdot 3\cdot 5\ldots(2k-1)}{2^kk!}\\
&=&(-1)^k\frac{(2k)!}{2^kk!2^kk!}\\
&=&\left(\frac{-1}{4}\right)^kD_k
\end{eqnarray*}

Thus, we are led to the conclusions in the statement.
\end{proof}

Which brings us into the numbers $C_k=\frac{1}{k+1}\binom{2k}{k}$, whose theory is as follows:

\begin{theorem}
The Catalan numbers $C_k=\frac{1}{k+1}\binom{2k}{k}$ are integers, and count:
\begin{enumerate}
\item The length $2k$ loops on $\mathbb N$, based at $0$.

\item The noncrossing pairings of $1,\ldots,2k$.

\item The noncrossing partitions of $1,\ldots,k$.

\item The length $2k$ Dyck paths in the plane.
\end{enumerate}
\end{theorem}

\begin{proof}
All this is standard combinatorics, the idea being as follows:

\medskip

(1) To start with, in what regards the various objects involved, the length $2k$ loops on $\mathbb N$ are the length $2k$ loops on $\mathbb N$ that we know, and the same goes for the noncrossing pairings of $1,\ldots,2k$, and for the noncrossing partitions of $1,\ldots,k$, the idea here being that you must be able to draw the pairing or partition in a noncrossing way. 

\medskip

(2) As for the length $2k$ Dyck paths in the plane, these are by definition the paths from $(0,0)$ to $(k,k)$, marching North-East over the integer lattice $\mathbb Z^2\subset\mathbb R^2$, by staying inside the square $[0,k]\times[0,k]$, and staying as well under the diagonal of this square. 

\medskip

(3) As an illustration for this, here are the 5 possible Dyck paths at $n=3$:
$$\xymatrix@R=4pt@C=4pt
{\circ&\circ&\circ&\circ\\
\circ&\circ&\circ&\circ\ar@{-}[u]\\
\circ&\circ&\circ&\circ\ar@{-}[u]\\
\circ\ar@{-}[r]&\circ\ar@{-}[r]&\circ\ar@{-}[r]&\circ\ar@{-}[u]}
\qquad
\xymatrix@R=4pt@C=4pt
{\circ&\circ&\circ&\circ\\
\circ&\circ&\circ&\circ\ar@{-}[u]\\
\circ&\circ&\circ\ar@{-}[r]&\circ\ar@{-}[u]\\
\circ\ar@{-}[r]&\circ\ar@{-}[r]&\circ\ar@{-}[u]&\circ}
\qquad
\xymatrix@R=4pt@C=4pt
{\circ&\circ&\circ&\circ\\
\circ&\circ&\circ\ar@{-}[r]&\circ\ar@{-}[u]\\
\circ&\circ&\circ\ar@{-}[u]&\circ\\
\circ\ar@{-}[r]&\circ\ar@{-}[r]&\circ\ar@{-}[u]&\circ}
\qquad
\xymatrix@R=4pt@C=4pt
{\circ&\circ&\circ&\circ\\
\circ&\circ&\circ&\circ\ar@{-}[u]\\
\circ&\circ\ar@{-}[r]&\circ\ar@{-}[r]&\circ\ar@{-}[u]\\
\circ\ar@{-}[r]&\circ\ar@{-}[u]&\circ&\circ}
\qquad
\xymatrix@R=4pt@C=4pt
{\circ&\circ&\circ&\circ\\
\circ&\circ&\circ\ar@{-}[r]&\circ\ar@{-}[u]\\
\circ&\circ\ar@{-}[r]&\circ\ar@{-}[u]&\circ\\
\circ\ar@{-}[r]&\circ\ar@{-}[u]&\circ&\circ}$$

(4) Thus, we have definitions for all objects involved, and in each case, if you start counting them, you always end up with the same sequence of numbers, namely:
$$1,2,5,14,42,132,429,1430,4862,16796,58786,\ldots$$

(5) In order to prove now that (1-4) produce indeed the same numbers, many things can be said, with a total of $\binom{4}{2}=6$ bijective proofs being possible. However, as a matter of having our claim proved, here is a quick argument. The point is that, in each of the cases (1-4) under consideration, the numbers $C_k$ that we get are easily seen to satisfy:
$$C_{k+1}=\sum_{a+b=k}C_aC_b$$ 

Now the initial data being the same, namely $C_1=1$ and $C_2=2$, in each of the cases (1-4) under consideration, we get indeed the same numbers, as desired. 

\medskip

(6) What is next? In view of what we already have, it remains to pick one of the objects (1-4), skilfully do the count, and conclude that we have indeed:
$$C_k=\frac{1}{k+1}\binom{2k}{k}$$

(7) The most convenient is to count the Dyck paths. For this purpose, we can use a trick. Indeed, if we ignore the assumption that our path must stay under the diagonal of the square, we have $\binom{2k}{k}$ such paths. And among these, we have the ``good'' ones, those that we want to count, and then the ``bad'' ones, those that we want to ignore.

\medskip

(8) So, let us count the bad paths, those crossing the diagonal of the square, and reaching the higher diagonal next to it, the one joining $(0,1)$ and $(k,k+1)$. In order to count these, the trick is to ``flip'' their bad part over that higher diagonal, as follows:
$$
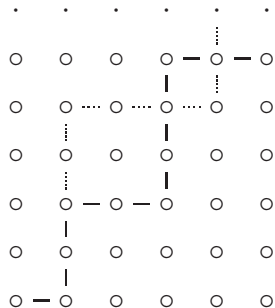
\xymatrix@R=4pt@C=4pt
{\cdot&\cdot&\cdot&\cdot&\cdot&\cdot\\
\circ&\circ&\circ&\circ\ar@{-}[r]&\circ\ar@{-}[r]\ar@{.}[u]&\circ\\
\circ&\circ\ar@{.}[r]&\circ\ar@{.}[r]&\circ\ar@{.}[r]\ar@{-}[u]&\circ\ar@{.}[u]&\circ\\
\circ&\circ\ar@{.}[u]&\circ&\circ\ar@{-}[u]&\circ&\circ\\
\circ&\circ\ar@{-}[r]\ar@{.}[u]&\circ\ar@{-}[r]&\circ\ar@{-}[u]&\circ&\circ\\
\circ&\circ\ar@{-}[u]&\circ&\circ&\circ&\circ\\
\circ\ar@{-}[r]&\circ\ar@{-}[u]&\circ&\circ&\circ&\circ}$$

(9) Now observe that, as it is obvious on the above picture, due to the flipping, the flipped bad path will no longer end in $(k,k)$, but rather in $(k-1,k+1)$. Moreover, more is true, in the sense that, by thinking a bit, we see that the flipped bad paths are precisely those ending in $(k-1,k+1)$. Thus, we can count these flipped bad paths, and so the bad paths, and so the good paths too, and so good news, we are done.

\medskip

(10) To finish now, by putting everything together, we have:
$$C_k
=\binom{2k}{k}-\binom{2k}{k-1}
=\binom{2k}{k}-\frac{k}{k+1}\binom{2k}{k}
=\frac{1}{k+1}\binom{2k}{k}$$

Thus, we are led to the various conclusions in the statement.
\end{proof}

We can go back now to the generalized binomial formula, and we have:

\index{square root}
\index{Catalan numbers}
\index{central binomial coefficients}

\begin{theorem}
The generalized binomial formula at $a=1/2,-1/2$ reads
$$\sqrt{1-4t}=1-2\sum_{k=1}^\infty C_{k-1}t^k\quad,\quad 
\frac{1}{\sqrt{1-4t}}=\sum_{k=0}^\infty D_kt^k$$
with $C_k=\frac{1}{k+1}\binom{2k}{k}$ and $D_k=\binom{2k}{k}$, and these formulae hold indeed, for $|t|<1/4$.
\end{theorem}

\begin{proof}
This is quite standard, based on what we have, as follows:

\medskip

(1) To start with, the formulae in Proposition 2.30 suggest to make the change of variables $x=-4t$ in the generalized binomial formula, and with this change made, that binomial formula at $a=1/2,-1/2$ corresponds precisely to the formulae above.

\medskip

(2) In order to prove our two formulae, we must establish the following identities:
$$\left(1-2\sum_{k=1}^\infty C_{k-1}t^k\right)^2=1-4t\quad,\quad 
\left(\sum_{k=0}^\infty D_kt^k\right)^2=\frac{1}{1-4t}$$

(3) But the first formula is equivalent to the following identity for the Catalan numbers, that we know well to hold, as explained in the proof of Theorem 2.31:
$$\sum_{k+l=n}C_kC_l=C_{n+1}$$

(4) As for the second formula, by using the standard series for $1/(1-4t)$, this is equivalent to the following formula, involving the central binomial coefficients:
$$\sum_{k+l=n}D_kD_l=4^n$$

(5) Now instead of doing again some combinatorics, this time for the numbers $D_k$, let us pull an analysis trick. With $t\to t-\varepsilon$ and $\varepsilon\simeq0$, our first formula becomes:
\begin{eqnarray*}
\sqrt{1-4t+4\varepsilon}
&=&1-2\sum_{k=1}^\infty C_{k-1}(t-\varepsilon)^k\\
&\simeq&1-2\sum_{k=1}^\infty C_{k-1}(t^k-kt^{k-1}\varepsilon)\\
&=&\sqrt{1-4t}+2\varepsilon\sum_{k=1}^\infty D_{k-1}t^{k-1}
\end{eqnarray*}

(6) On the other hand, again with $\varepsilon\simeq0$, we have the following estimate:
$$\sqrt{1-4t+4\varepsilon}-\sqrt{1-4t}=\frac{4\varepsilon}{\sqrt{1-4t+4\varepsilon}+\sqrt{1-4t}}\simeq\frac{2\varepsilon}{\sqrt{1-4t}}$$

We conclude from this that we have the following formula, as desired:
$$\frac{1}{\sqrt{1-4t}}=\sum_{k=1}^\infty D_{k-1}t^{k-1}$$

Summarizing, both formulae in the statement proved, one way or another.
\end{proof}

With this discussed, let us get now into exp and log. We first have:

\index{e}
\index{exponential}

\begin{proposition}
We have the following formula,
$$\left(1+\frac{x}{n}\right)^n\to e^x$$
valid for any $x\in\mathbb R$.
\end{proposition}

\begin{proof}
We know that this holds at $x=1$, by definition of $e$, and by inverting, we have it at $x=-1$ too. But then, when $x\in\mathbb R$ is arbitrary, we can proceed as follows:
$$\left(1+\frac{x}{n}\right)^n=\left[\left(1+\frac{x}{n}\right)^{n/x}\right]^x
\to e^x$$

Thus, we are led to the conclusion in the statement.
\end{proof}

Next, we have the following result, which is something quite far-reaching:

\begin{theorem}
We have the following formula,
$$e^x=\sum_{k=0}^\infty\frac{x^k}{k!}$$
valid for any $x\in\mathbb R$.
\end{theorem}

\begin{proof}
This can be done in several steps, as follows:

\medskip

(1) At $x=1$, we want to prove that we have the following equality:
$$\sum_{k=0}^\infty\frac{1}{k!}=\lim_{n\to\infty}\left(1+\frac{1}{n}\right)^n$$

For this purpose, the first observation is that we have the following estimate:
$$2<\sum_{k=0}^\infty\frac{1}{k!}<\sum_{k=0}^\infty\frac{1}{2^{k-1}}=3$$

In order to prove now that this limit is indeed $e$, observe that we have:
$$\left(1+\frac{1}{n}\right)^n
=\sum_{k=0}^n\binom{n}{k}\cdot\frac{1}{n^k}
\leq\sum_{k=0}^n\frac{1}{k!}$$

Thus, with $n\to\infty$, we get that the limit of the series $\sum_{k=0}^\infty\frac{1}{k!}$ belongs to $[e,3)$.

\medskip

(2) For the reverse inequality, we have the following computation:
\begin{eqnarray*}
\sum_{k=0}^n\frac{1}{k!}-\left(1+\frac{1}{n}\right)^n
&=&\sum_{k=2}^n\frac{n^k-n(n-1)\ldots(n-k+1)}{n^kk!}\\
&\leq&\sum_{k=2}^n\frac{n^k-(n-k)^k}{n^kk!}\\
&=&\sum_{k=2}^n\frac{1-\left(1-\frac{k}{n}\right)^k}{k!}
\end{eqnarray*}

Next, we can use the following trivial inequality, valid for any number $x\in(0,1)$:
$$1-x^k=(1-x)(1+x+x^2+\ldots+x^{k-1})\leq(1-x)k$$

Indeed, we can use this with $x=1-k/n$, and we obtain in this way:
$$\sum_{k=0}^n\frac{1}{k!}-\left(1+\frac{1}{n}\right)^n
\leq\sum_{k=2}^n\frac{\frac{k}{n}\cdot k}{k!}
\leq\frac{1}{n}\sum_{k=2}^n\frac{2}{2^{k-2}}
<\frac{4}{n}$$

Thus, we have our needed estimate, and so done with the case $x=1$.

\medskip

(3) In order to deal now with the general case, consider the following function:
$$f(x)=\sum_{k=0}^\infty\frac{x^k}{k!}$$

Our claim, which is the key one, is that we have the following formula:
$$f(x+y)=f(x)f(y)$$

Indeed, by using the binomial formula, we have the following computation:
\begin{eqnarray*}
f(x+y)
&=&\sum_{k=0}^\infty\frac{(x+y)^k}{k!}\\
&=&\sum_{k=0}^\infty\sum_{s=0}^k\binom{k}{s}\cdot\frac{x^sy^{k-s}}{k!}\\
&=&\sum_{k=0}^\infty\sum_{s=0}^k\frac{x^sy^{k-s}}{s!(k-s)!}\\
&=&f(x)f(y)
\end{eqnarray*}

(4) As a first observation, this shows that $f$ is continuous. Indeed, at $x=0$ we have:
$$\lim_{t\to0}f(t)=\lim_{t\to0}\left(1+t\sum_{k=1}^\infty\frac{t^{k-1}}{k!}\right)=1$$

But from this, we get $f(x+t)=f(x)f(t)\to f(x)$ with $t\to0$, at any $x$. Thus, as a conclusion, our function $f$ is continuous, and satisfies the following conditions:
$$f(x+y)=f(x)f(y)\quad,\quad f(1)=e$$

But with this, we can finish. Indeed, by iterating, we have $f(nx)=f(x)^n$ for any $n\in\mathbb N$. Then, by extracting roots, we have $f(rx)=f(x)^r$ for any $r\in\mathbb Q$. Thus $f(r)=e^r$ for any $r\in\mathbb Q$, and by continuity we obtain $f(x)=e^x$ for any $x\in\mathbb R$, as desired.
\end{proof}

Next, we can talk about the logarithm function, which appears as follows:

\begin{theorem}
The exponential function, as constructed before,
$$\exp:\mathbb R\to(0,\infty)\quad,\quad x\to e^x$$
is invertible, with its inverse being the logarithm function, denoted as follows:
$$\log:(0,\infty)\to\mathbb R\quad,\quad \log=\exp^{-1}$$
This logarithm function is continuous, increasing, and bijective.
\end{theorem}

\begin{proof}
This is indeed something self-explanatory, based on Theorem 2.17.
\end{proof}

In practice, the logarithm is a very useful function, and more on this later. In the meantime, how to compute it? And here, we have the following collection of results:

\begin{theorem}
The logarithm can be computed via the equivalent formulae
$$e^{\log x}=x\quad,\quad \log(e^x)=x$$
and in practice, we have the following useful rules:
\begin{enumerate}
\item $\log(xy)=\log x+\log y$.

\item $\log(1/y)=-\log y$.

\item $\log(x/y)=\log x-\log y$.

\item $\log(x^p)=p\log x$.
\end{enumerate}
\end{theorem}

\begin{proof}
The first two formulae come from the fact that the inverse of a function $f$ can be defined either via $f(f^{-1}(x))=x$, or via $f^{-1}(f(x))=x$. As for the rest:

\medskip

(1) This comes indeed from the following computation:
$$e^{\log(xy)}=xy=e^{\log x}e^{\log y}=e^{\log x+\log y}$$

(2) This comes from (1), by setting $x=1/y$, and using $\log 1=0$.

\medskip

(3) This comes also from (1), by replacing $y\to1/y$, and using (2).

\medskip

(4) This formula, generalizing (1) with $x=y$, and also (2), comes as follows:
$$e^{\log(x^p)}=x^p=(e^{\log x})^p=e^{p\log x}$$

Thus, we are led to the conclusions in the statement.
\end{proof}

We have as well the following formula, for $|x|<1$, which can be proved by showing, using Theorem 2.29, that the series on the right satisfies $f(xy)=f(x)+f(y)$:
$$\log(1+x)=\sum_{k=1}^\infty(-1)^{k+1}\frac{x^k}{k}$$

And more on this later. We will be actually back to powers, exponentials, logarithms and trigonometric functions on numerous occasions, in what follows.

\section*{2e. Exercises}

There are many possible exercises on the above, and here are a few of them:

\begin{exercise}
Compute Lipschitz constants for all functions that you know.
\end{exercise}

\begin{exercise}
Learn some alternative proofs of the intermediate value theorem.
\end{exercise}

\begin{exercise}
Prove the generalized binomial formula, for exponents $a\in\mathbb Z/2$.
\end{exercise}

\begin{exercise}
Rewrite the theory of $e$, with $e=\sum_k1/k!$ as definition.
\end{exercise}

For the rest, business as usual, more exercises are easy to find. Find and solve them.

\chapter{Derivatives}

\section*{3a. Derivatives, rules}

Welcome to calculus. In this chapter we go for the real thing, namely development of modern calculus, following some amazing ideas of Newton, Leibnitz and others. The material will be quite difficult, mixing geometry and intuition with formal mathematics and computations, and needing some time to be understood. But we will survive.

\bigskip

The basic idea of calculus is very simple. We are interested in functions $f:\mathbb R\to\mathbb R$, and we already know that when $f$ is continuous at a point $x$, we can write an approximation formula as follows, for the values of our function $f$ around that point $x$:
$$f(x+t)\simeq f(x)$$

The problem is now, how to improve this? And a bit of thinking at all this suggests to look at the slope of $f$ at the point $x$, according to the following picture:
$$\xymatrix@R=4pt@C=15pt{
&&&&&&\\
&&&&&&\\
&&&&&&\\
&&&&&&\\
&&&&&&\\
&f(x+t)\ar[uuuuu]&&&&\circ\ar@{--}[dddd]\ar@{--}[llll]\ar@{.}[uur]\ar@{-}@/_/[uuuur]^f&\\
&&&&&&\\
&f(x)\ar@{-}[uu]&&&\circ\ar@{--}[dd]\ar@{--}[lll]\ar@{--}[r]\ar@{.}[uur]\ar@{-}@/_/[uur]&&\\
&&\ar@{-}@/_/[urr]&&&\\
\ar@{-}[rrrr]&&&\ar@{.}[uur]\ar@{-}[r]^\alpha&x\ar@{-}[r]&x+t\ar[rr]&&\\
&\\
&\ar@{-}[uuuu]&&&&\\}$$

Indeed, the slope of $f$ at the point $x$ appears as the $t\to0$ limit of the quantities $\tan\alpha=(f(x+t)-f(x))/t$, which leads us into the following notion:

\index{differentiable function}
\index{derivative}

\begin{definition}
A function $f:\mathbb R\to\mathbb R$ is called differentiable at $x$ when
$$f'(x)=\lim_{t\to0}\frac{f(x+t)-f(x)}{t}$$
called derivative, or slope of $f$ at that point $x$, exists.
\end{definition}

As a first remark, in order for $f$ to be differentiable at $x$, that is to say, in order for the above limit to converge, the numerator must go to $0$, as the denominator $t$ does:
$$\lim_{t\to0}\left[f(x+t)-f(x)\right]=0$$

Thus, $f$ must be continuous at $x$. However, the converse is not true, a basic counterexample being $f(x)=|x|$ at $x=0$. Let us summarize these findings as follows:

\index{modulus}

\begin{proposition}
If $f$ is differentiable at $x$, then $f$ must be continuous at $x$. However, the converse is not true, a basic counterexample being $f(x)=|x|$, at $x=0$.
\end{proposition}

\begin{proof}
The first assertion is something that we already know, from the above. As for the second assertion, regarding $f(x)=|x|$, this is something quite clear on the picture of $f$, but let us prove this mathematically, based on Definition 3.1. We have:
$$\lim_{t\searrow 0}\frac{|0+t|-|0|}{t}=\lim_{t\searrow 0}\frac{t-0}{t}=1$$

On the other hand, we have as well the following computation:
$$\lim_{t\nearrow 0}\frac{|0+t|-|0|}{t}=\lim_{t\nearrow 0}\frac{-t-0}{t}=-1$$

Thus, the limit in Definition 3.1 does not converge, so we have our counterexample. 
\end{proof}

Generally speaking, the last assertion in Proposition 3.2 should not bother us much, because most of the basic continuous functions are differentiable, and we will see examples in a moment. Before that, however, let us recall why we are here, namely improving the basic estimate $f(x+t)\simeq f(x)$. We can now do this, using the derivative, as follows:

\index{derivative}
\index{locally affine}

\begin{theorem}
Assuming that $f$ is differentiable at $x$, we have:
$$f(x+t)\simeq f(x)+f'(x)t$$
In other words, $f$ is, approximately, locally affine at $x$.
\end{theorem}

\begin{proof}
Assume indeed that $f$ is differentiable at $x$, and let us set, as before:
$$f'(x)=\lim_{t\to0}\frac{f(x+t)-f(x)}{t}$$

By multiplying by $t$, we obtain that we have, once again in the $t\to0$ limit:
$$f(x+t)-f(x)\simeq f'(x)t$$

Thus, we are led to the conclusion in the statement.
\end{proof}

Still talking theory, let us make as well the following remark:

\begin{remark}
The differentiability is, exactly as continuity, a local notion. That is, we can talk about the differentiability of $f:X\to\mathbb R$, with $X\subset\mathbb R$, at any point $x\in X$, provided that $X$ is open, or at least contains a small interval, around $x$.
\end{remark}

Now with this book being an introduction to calculus, as opposed to a treatise on calculus, or to a life and death matter, in general, we will agree to be a bit sloppy on this, by keeping talking about functions $f:\mathbb R\to\mathbb R$, which is quite practical, and makes the basic concepts more visible. These functions will be understood to be as in Remark 3.4, or even more general, because we can talk for instance about the differentiability $f:[a,b]\to\mathbb R$ at the endpoints $a,b$, using right and left limits, in the obvious way.

\bigskip

Time for some examples? Here is a first key computation of derivatives:

\index{power function}

\begin{proposition}
We have the differentiation formula
$$(x^p)'=px^{p-1}$$
valid for any exponent $p\in\mathbb R$.
\end{proposition}

\begin{proof}
We can do this in three steps, as follows:

\medskip

(1) In the case $p\in\mathbb N$ we can use the binomial formula, which gives, as desired:
\begin{eqnarray*}
(x+t)^p
&=&\sum_{k=0}^n\binom{p}{k}x^{p-k}t^k\\
&=&x^p+px^{p-1}t+\ldots+t^p\\
&\simeq&x^p+px^{p-1}t
\end{eqnarray*}

(2) Let us discuss now the general case $p\in\mathbb Q$. We write $p=m/n$, with $m\in\mathbb Z$ and $n\in\mathbb N$. In order to do the computation, we use the following formula:
$$a^n-b^n=(a-b)(a^{n-1}+a^{n-2}b+\ldots+b^{n-1})$$

We set in this formula $a=(x+t)^{m/n}$ and $b=x^{m/n}$. We obtain, as desired: 
\begin{eqnarray*}
(x+t)^{m/n}-x^{m/n}
&=&\frac{(x+t)^m-x^m}{(x+t)^{m(n-1)/n}+\ldots+x^{m(n-1)/n}}\\
&\simeq&\frac{(x+t)^m-x^m}{nx^{m(n-1)/n}}\\
&\simeq&\frac{mx^{m-1}t}{nx^{m(n-1)/n}}\\
&=&\frac{m}{n}\cdot x^{m-1-m+m/n}\cdot t\\
&=&\frac{m}{n}\cdot x^{m/n-1}\cdot t
\end{eqnarray*}

(3) In the general case now, where $p\in\mathbb R$ is real, we can use a similar argument. Indeed, given any integer $n\in\mathbb N$, we have the following computation:
$$(x+t)^p-x^p
=\frac{(x+t)^{pn}-x^{pn}}{(x+t)^{p(n-1)}+\ldots+x^{p(n-1)}}
\simeq\frac{(x+t)^{pn}-x^{pn}}{nx^{p(n-1)}}$$

Now observe that we have the following estimate, with $[.]$ being the integer part:
$$(x+t)^{[pn]}\leq (x+t)^{pn}\leq (x+t)^{[pn]+1}$$

By using the binomial formula on both sides, for the integer exponents $[pn]$ and $[pn]+1$ there, we deduce that with $n>>0$ we have the following estimate:
$$(x+t)^{pn}\simeq x^{pn}+pnx^{pn-1}t$$

Thus, we can finish our computation started above as follows:
$$(x+t)^p-x^p
\simeq\frac{pnx^{pn-1}t}{nx^{pn-p}}
=px^{p-1}t$$

But this gives $(x^p)'=px^{p-1}$, which finishes the proof.
\end{proof}

Here are some further computations, for other basic functions that we know:

\index{sin}
\index{cos}
\index{exp}
\index{log}

\begin{proposition}
We have the following results:
\begin{enumerate}
\item $(\sin x)'=\cos x$.

\item $(\cos x)'=-\sin x$.

\item $(e^x)'=e^x$.

\item $(\log x)'=x^{-1}$.
\end{enumerate}
\end{proposition}

\begin{proof}
This is quite tricky, as always when computing derivatives, as follows:

\medskip

(1) Regarding $\sin$, the computation here goes as follows:
\begin{eqnarray*}
(\sin x)'
&=&\lim_{t\to0}\frac{\sin(x+t)-\sin x}{t}\\
&=&\lim_{t\to0}\frac{\sin x\cos t+\cos x\sin t-\sin x}{t}\\
&=&\lim_{t\to0}\sin x\cdot\frac{\cos t-1}{t}+\cos x\cdot\frac{\sin t}{t}\\
&=&\cos x
\end{eqnarray*}

Here we have used the fact, which is clear on pictures, by drawing the trigonometric circle, that we have $\sin t\simeq t$ for $t\simeq 0$, plus the fact, which follows from this and from Pythagoras, $\sin^2+\cos^2=1$, that we have as well $\cos t\simeq 1-t^2/2$, for $t\simeq 0$.

\medskip

(2) The computation for $\cos$ is similar, as follows:
\begin{eqnarray*}
(\cos x)'
&=&\lim_{t\to0}\frac{\cos(x+t)-\cos x}{t}\\
&=&\lim_{t\to0}\frac{\cos x\cos t-\sin x\sin t-\cos x}{t}\\
&=&\lim_{t\to0}\cos x\cdot\frac{\cos t-1}{t}-\sin x\cdot\frac{\sin t}{t}\\
&=&-\sin x
\end{eqnarray*}

(3) For the exponential, the derivative can be computed as follows:
$$(e^x)'
=\left(\sum_{k=0}^\infty\frac{x^k}{k!}\right)'
=\sum_{k=0}^\infty\frac{kx^{k-1}}{k!}
=e^x$$

(4) As for the logarithm, the computation here is as follows, using $\log(1+t)\simeq t$ for $t\simeq 0$, which follows from the basic estimate $e^t\simeq 1+t$, by taking the logarithm:
\begin{eqnarray*}
(\log x)'
&=&\lim_{t\to0}\frac{\log(x+t)-\log x}{t}\\
&=&\lim_{t\to0}\frac{\log(1+t/x)}{t}\\
&=&\frac{1}{x}
\end{eqnarray*}

Thus, we are led to the formulae in the statement.
\end{proof}

Speaking exponentials, we can now formulate a nice result about them:

\index{exponential}

\begin{theorem}
The following happen, regarding $e$ and $\exp$:
\begin{enumerate}
\item $e$ is the unique number satisfying $e^t\simeq 1+t$.

\item $\exp$ is the unique power series satisfying $f'=f$ and $f(0)=1$.
\end{enumerate}
\end{theorem}

\begin{proof}
The first assertion follows from the following computation, with $a>0$:
$$a^t=e^{t\log a}\simeq 1+t\log a$$

Regarding the second assertion, consider a power series satisfying $f'=f$ and $f(0)=1$. Due to $f(0)=1$, the first term must be 1, so our series must look as follows:
$$f(x)=1+\sum_{k=1}^\infty c_kx^k$$

According to our differentiation rules, the derivative of this series is given by:
$$f(x)=\sum_{k=1}^\infty kc_kx^{k-1}$$

Thus, the equation $f'=f$ is equivalent to the following equalities:
$$c_1=1\quad,\quad 2c_2=c_1\quad,\quad 3c_3=c_2\quad,\quad 4c_4=c_3\quad,\quad\ldots$$

But this system of equations can be solved by recurrence, as follows:
$$c_1=1\quad,\quad c_2=\frac{1}{2}\quad,\quad c_3=\frac{1}{2\times 3}\quad,\quad c_4=\frac{1}{2\times 3\times 4}\quad,\quad\ldots$$

We therefore have $c_k=1/k!$, leading to $f(x)=e^x$, as claimed.
\end{proof}

Back now to computations and rules, we have here the following statement:

\index{chain rule}

\begin{theorem}
We have the following formulae:
\begin{enumerate}
\item $(f+g)'=f'+g'$.

\item $(fg)'=f'g+fg'$.

\item $(f\circ g)'=(f'\circ g)\cdot g'$.
\end{enumerate}
\end{theorem}

\begin{proof}
All these formulae are elementary, the idea being as follows:

\medskip

(1) This follows indeed from definitions, the computation being as follows:
\begin{eqnarray*}
(f+g)'(x)
&=&\lim_{t\to0}\frac{(f+g)(x+t)-(f+g)(x)}{t}\\
&=&\lim_{t\to0}\left(\frac{f(x+t)-f(x)}{t}+\frac{g(x+t)-g(x)}{t}\right)\\
&=&\lim_{t\to0}\frac{f(x+t)-f(x)}{t}+\lim_{t\to0}\frac{g(x+t)-g(x)}{t}\\
&=&f'(x)+g'(x)
\end{eqnarray*}

(2) This follows from definitions too, the computation, by using the more convenient formula $f(x+t)\simeq f(x)+f'(x)t$ as a definition for the derivative, being as follows:
\begin{eqnarray*}
(fg)(x+t)
&=&f(x+t)g(x+t)\\
&\simeq&(f(x)+f'(x)t)(g(x)+g'(x)t)\\
&\simeq&f(x)g(x)+(f'(x)g(x)+f(x)g'(x))t
\end{eqnarray*}

Indeed, we obtain from this that the derivative is the coefficient of $t$, namely:
$$(fg)'(x)=f'(x)g(x)+f(x)g'(x)$$

(3) Regarding compositions, the computation here is as follows, again by using the more convenient formula $f(x+t)\simeq f(x)+f'(x)t$ as a definition for the derivative:
\begin{eqnarray*}
(f\circ g)(x+t)
&=&f(g(x+t))\\
&\simeq&f(g(x)+g'(x)t)\\
&\simeq&f(g(x))+f'(g(x))g'(x)t
\end{eqnarray*}

Indeed, we obtain from this that the derivative is the coefficient of $t$, namely:
$$(f\circ g)'(x)=f'(g(x))g'(x)$$

Thus, we are led to the conclusions in the statement.
\end{proof}

As a basic application of the formula (3) above, called chain rule, we have:
$$(a^x)'=(e^{x\log a})'=e^{x\log a}\cdot\log a=a^x\log a$$

We can of course combine the above formulae, and we obtain for instance:

\index{fraction}

\begin{proposition}
The derivatives of fractions are given by:
$$\left(\frac{f}{g}\right)'=\frac{f'g-fg'}{g^2}$$
In particular, we have the following formula, for the derivative of inverses:
$$\left(\frac{1}{f}\right)'=-\frac{f'}{f^2}$$
In fact, we have $(f^p)'=pf^{p-1}$, for any exponent $p\in\mathbb R$.
\end{proposition}

\begin{proof}
This statement is written a bit upside down, and for the proof it is better to proceed backwards. To be more precise, by using $(x^p)'=px^{p-1}$ and Theorem 3.8 (3), we obtain the third formula. Then, with $p=-1$, we obtain from this the second formula. And finally, by using this second formula and Theorem 3.8 (2), we obtain:
$$\left(\frac{f}{g}\right)'
=f'\cdot\frac{1}{g}+f\left(\frac{1}{g}\right)'
=\frac{f'}{g}-\frac{fg'}{g^2}
=\frac{f'g-fg'}{g^2}$$

Thus, we are led to the formulae in the statement.
\end{proof}

And with this, good news, we have all needed formulae in our bag. Sky is the limit, and in what regards for instance the various trigonometric functions, we have:
\index{tan}
\index{arctan}

\begin{proposition}
We have the following formulae,
$$(\tan x)'=\frac{1}{\cos^2x}\quad,\quad (\arctan x)'=\frac{1}{1+x^2}$$
and the derivatives of the remaining trigonometric functions can be computed as well.
\end{proposition}

\begin{proof}
For $\tan=\sin/\cos$, we have indeed the following computation:
\begin{eqnarray*}
(\tan x)'
&=&\frac{\sin'x\cos x-\sin x\cos'x}{\cos^2x}\\
&=&\frac{\cos^2x+\sin^2x}{\cos^2x}\\
&=&\frac{1}{\cos^2x}
\end{eqnarray*}

As for $\arctan=\tan^{-1}$, we can use here the following computation:
\begin{eqnarray*}
(\tan\circ\arctan)'(x)
&=&\tan'(\arctan x)\arctan'(x)\\
&=&\frac{1}{\cos^2(\arctan x)}\arctan'(x)
\end{eqnarray*}

Indeed, since the term on the left is simply $x'=1$, we obtain from this:
$$\arctan'(x)=\cos^2(\arctan x)$$

On the other hand, with $t=\arctan x$ we know that we have $\tan t=x$, and so:
$$\cos^2(\arctan x)=\cos^2t=\frac{1}{1+\tan^2t}=\frac{1}{1+x^2}$$

Thus, we are led to the formula in the statement, namely:
$$(\arctan x)'=\frac{1}{1+x^2}$$

As for the last assertion, we will leave this as an exercise. Of course, in case of an emergency, you can find the formulae in any 1-variable book, including mine \cite{ba1}.
\end{proof}

At the theoretical level now, further building on Theorem 3.3, we have:

\index{local minimum}
\index{local maximum}
\index{minimum}
\index{maximum}

\begin{theorem}
The local minima and maxima of a differentiable function $f:\mathbb R\to\mathbb R$ appear at the points $x\in\mathbb R$ where:
$$f'(x)=0$$
However, the converse of this fact is not true in general.
\end{theorem}

\begin{proof}
The first assertion follows from the formula in Theorem 3.3, namely:
$$f(x+t)\simeq f(x)+f'(x)t$$

Indeed, saying that our function $f$ has a local maximum at $x\in\mathbb R$ means that there exists a number $\varepsilon>0$ such that the following happens:
$$f(x+t)\geq f(x)\quad,\quad\forall t\in[-\varepsilon,\varepsilon]$$

Thus we must have $f'(x)t\geq0$ for sufficiently small $t$, and since this small $t$ can be both positive or negative, this gives the following condition, as desired:
$$f'(x)=0$$

The discussion for the local minima is similar. Finally, in what regards the converse, the simplest counterexample here is the function $f(x)=x^3$, at $x=0$.
\end{proof}

In practice, Theorem 3.11 can be used in order to find the maximum and minimum of any differentiable function, and this method is best recalled as follows:

\begin{algorithm}
In order to find the minimum and maximum of $f:[a,b]\to\mathbb R$:
\begin{enumerate}
\item Compute the derivative $f'$.

\item Solve the equation $f'(x)=0$.

\item Add $a,b$ to your set of solutions.

\item Compute $f(x)$, for all your solutions.

\item Compute the min/max of all these $f(x)$ values.

\item Then this is the min/max of your function.
\end{enumerate}
\end{algorithm}

As another important consequence of Theorem 3.11, we have:

\begin{theorem}
Assuming that $f:[a,b]\to\mathbb R$ is differentiable, we have
$$\frac{f(b)-f(a)}{b-a}=f'(c)$$
for some $c\in(a,b)$, called mean value property of $f$.
\end{theorem}

\begin{proof}
In the case $f(a)=f(b)$, the result, called Rolle theorem, states that we have $f'(c)=0$ for some $c\in(a,b)$, and follows from Theorem 3.11. Now in what regards our statement, due to Lagrange, this follows from Rolle, applied to the following function:
$$g(x)=f(x)-\frac{f(b)-f(a)}{b-a}\cdot x$$

Indeed, we have $g(a)=g(b)$, due to our choice of the constant on the right, so we get $g'(c)=0$ for some $c\in(a,b)$, which translates into the formula in the statement.
\end{proof}

As a key consequence of Theorem 3.13, of great practical interest, we have:

\begin{theorem}
For a differentiable function we have
$$f'=0\quad\implies\quad f={\rm constant}$$
and with the converse of this being of course true too.
\end{theorem}

\begin{proof}
This is indeed something self-explanatory, coming from Theorem 3.13.
\end{proof}

Time for some applications? Remember from chapter 2 the difficulties with the generalized binomial formula. We can now nuke all that problematics, as follows:

\begin{theorem}
We have the generalized binomial formula
$$(1+x)^p=\sum_{k=0}^\infty\binom{p}{k}x^k$$
with the generalized binomial coefficients being given by
$$\binom{p}{k}=\frac{p(p-1)\ldots(p-k+1)}{k!}$$
valid for any exponent $p\in\mathbb R$, and any $|x|<1$.
\end{theorem}

\begin{proof}
The series in the statement $f$ converges indeed at $|x|<1$, thanks to our convergence results from chapter 1, which apply. Also, we have the following formula:
$$(1+x)f'(x)=pf(x)$$

Now by using this formula, we have the following computation:
$$\left((1+x)^{-p}f(x)\right)'=-p(1+x)^{-p-1}f(x)+(1+x)^{-p}f'(x)=0$$

Thus we have $f(x)=c(1+x)^p$, with $c=f(0)=1$, as desired.
\end{proof}

As another application of Theorem 3.14, we can improve Theorem 3.7, as follows:

\begin{theorem}
The exponential function is the unique solution of
$$f'=f\quad,\quad f(0)=1$$
and as a consequence, $e=f(1)$, with $f$ being this unique solution.
\end{theorem}

\begin{proof}
Since we have $f(0)=1$ and $f'=f$ we conclude that we have $f\geq1$ for $x\geq0$, and a similar backwards argument shows that we have as well $f>0$, for $x<0$. In short, we have $f>0$ over the whole $\mathbb R$, and in particular $f\neq 0$. But with this, we have:
\begin{eqnarray*}
f'=f
&\implies&\frac{f'}{f}=1\\
&\implies&(\log f)'=1\\
&\implies&\log f=x+c\\
&\implies&f=\lambda e^x
\end{eqnarray*}

Now by using $f(0)=1$ we conclude that we have $f(x)=e^x$, as desired.
\end{proof}

Quite interesting, this latter result. Is analysis part of algebra? To be seen.

\section*{3b. Second derivatives}

The derivative theory that we have is already quite powerful, and can be used in order to solve all sorts of interesting questions, but with a bit more effort, we can do better. Indeed, at a more advanced level, we can come up with the following notion:

\begin{definition}
We say that $f:\mathbb R\to\mathbb R$ is twice differentiable if it is differentiable, and its derivative $f':\mathbb R\to\mathbb R$ is differentiable too. The derivative of $f'$ is denoted
$$f'':\mathbb R\to\mathbb R$$ 
and is called second derivative of $f$.
\end{definition}

You might probably wonder why coming with this definition, which looks a bit abstract and complicated, instead of further developing the theory of the first derivative, which looks like something very reasonable and useful. Good point, and answer to this coming in a moment. But before that, let us get a bit familiar with $f''$. We first have:

\begin{interpretation}
The second derivative $f''(x)\in\mathbb R$ is the number which:
\begin{enumerate}
\item Expresses the growth rate of the slope $f'(z)$ at the point $x$.

\item Gives us the acceleration of the function $f$ at the point $x$.

\item Computes how much different is $f(x)$, compared to $f(z)$ with $z\simeq x$.

\item Tells us how much convex or concave is $f$, around the point $x$.
\end{enumerate}
\end{interpretation}

So, this is the truth about the second derivative, making it clear that what we have here is a very interesting notion. In practice now, (1) follows from the usual interpretation of the derivative, as both a growth rate, and a slope. Regarding (2), this is some sort of reformulation of (1), using the intuitive meaning of the word ``acceleration'', with the relevant physics equations, due to Newton, being as follows: 
$$v=\dot{x}\quad,\quad a=\dot{v}$$

To be more precise, here $x,v,a$ are the position, speed and acceleration, and the dot denotes the time derivative, and according to these equations, we have $a=\ddot{x}$, second derivative. We will be back to these equations at the end of the present chapter.

\bigskip

Regarding now (3) in the above, this is something more subtle, of statistical nature, that we will clarify with some mathematics, in a moment. As for (4), this is something quite subtle too, that we will again clarify with some mathematics, in a moment. 

\bigskip

In practice now, let us first compute the second derivatives of the functions that we are familiar with, see what we get. The result here, which is perhaps not very enlightening at this stage of things, but which certainly looks technically useful, is as follows:

\begin{proposition}
The second derivatives of the basic functions are as follows:
\begin{enumerate}
\item $(x^p)''=p(p-1)x^{p-2}$.

\item $\sin''=-\sin$.

\item $\cos''=-\cos$.

\item $\exp'=\exp$.

\item $\log'(x)=-1/x^2$.
\end{enumerate}
Also, there are functions which are differentiable, but not twice differentiable.
\end{proposition}

\begin{proof}
We have several assertions here, the idea being as follows:

\medskip

(1) Regarding the various formulae in the statement, these all follow from the various formulae for the derivatives established before, as follows:
$$(x^p)''=(px^{p-1})'=p(p-1)x^{p-2}$$
$$(\sin x)''=(\cos x)'=-\sin x$$
$$(\cos x)''=(-\sin x)'=-\cos x$$
$$(e^x)''=(e^x)'=e^x$$
$$(\log x)''=(-1/x)'=-1/x^2$$

Of course, this is not the end of the story, because these formulae remain quite opaque, and must be examined in view of Interpretation 3.18, in order to see what exactly is going on. Also, we have $\tan$ and the other trigonometric functions too. In short, plenty of good exercises here, for you, and the more you solve, the better your calculus will be.

\medskip

(2) Regarding now the counterexample, recall first that the simplest example of a function which is continuous, but not differentiable, was $f(x)=|x|$, the idea behind this being to use a ``piecewise linear function whose branches do not fit well''. In connection now with our question, piecewise linear will not do, but we can use a similar idea, namely ``piecewise quadratic function whose branches do not fit well''. So, let us set:
$$f(x)=\begin{cases}
-x^2& (x\leq0)\\
x^2& (x\geq 0)
\end{cases}$$

The derivative is then $f'(x)=2|x|$, which is not differentiable, as desired.
\end{proof}

Getting back now to theory, we will need the following standard result:

\index{L'H\^opital's rule}

\begin{proposition}
The $0/0$ type limits can be computed according to the formula
$$\frac{f(x)}{g(x)}\simeq\frac{f'(x)}{g'(x)}$$
called L'H\^opital's rule.
\end{proposition}

\begin{proof}
The above formula holds indeed, as an application of the general first derivative theory from before, which gives, in the situation from the statement:
$$\frac{f(x+t)}{g(x+t)}
\simeq\frac{f(x)+f'(x)t}{g(x)+g'(x)t}
=\frac{f'(x)t}{g'(x)t}
=\frac{f'(x)}{g'(x)}$$

Thus, we are led to the conclusion in the statement.
\end{proof}

We can now formulate the following key result, improving Theorem 3.3:

\index{Taylor formula}
\index{second derivative}

\begin{theorem}
Given a twice differentiable function $f:\mathbb R\to\mathbb R$, we have:
$$f(x+t)\simeq f(x)+f'(x)t+\frac{f''(x)}{2}\,t^2$$
That is, $f$ is approximately locally quadratic.
\end{theorem}

\begin{proof}
Assume indeed that $f$ is twice differentiable at $x$, and let us try to construct an approximation of $f$ around $x$ by a quadratic function, as follows:
$$f(x+t)\simeq a+bt+ct^2$$

We must have $a=f(x)$, and we also know from Theorem 3.3 that $b=f'(x)$ is the correct choice for the coefficient of $t$. Thus, our approximation must be as follows:
$$f(x+t)\simeq f(x)+f'(x)t+ct^2$$

In order to find the correct choice for $c\in\mathbb R$, observe that the function $t\to f(x+t)$ matches with $t\to f(x)+f'(x)t+ct^2$ in what regards the value at $t=0$, and also in what regards the value of the derivative at $t=0$. Thus, the correct choice of $c\in\mathbb R$ should be the one making match the second derivatives at $t=0$, and this gives:
$$f''(x)=2c$$

We are therefore led to the formula in the statement, namely:
$$f(x+t)\simeq f(x)+f'(x)t+\frac{f''(x)}{2}\,t^2$$

In order to prove now that this formula holds indeed, we can use L'H\^opital's rule. Indeed, by using it, if we denote by $\varphi(t)\simeq P(t)$ the formula to be proved, we have:
\begin{eqnarray*}
\frac{\varphi(t)-P(t)}{t^2}
&\simeq&\frac{\varphi'(t)-P'(t)}{2t}\\
&\simeq&\frac{\varphi''(t)-P''(t)}{2}\\
&=&\frac{f''(x)-f''(x)}{2}\\
&=&0
\end{eqnarray*}

Thus, we are led to the conclusion in the statement.
\end{proof}

The above result substantially improves Theorem 3.3, and there are many applications of it. As a first such application, justifying Interpretation 3.18 (3), we have the following statement, which is a bit heuristic, but we will call it however Proposition:

\begin{proposition}
Intuitively speaking, the second derivative $f''(x)\in\mathbb R$ computes how much different is $f(x)$, compared to the average of $f(z)$, with $z\simeq x$.
\end{proposition}

\begin{proof}
As already mentioned, this is something a bit heuristic, but which is good to know. Let us write the formula in Theorem 3.21 as such, and with $t\to-t$ too:
$$f(x+t)\simeq f(x)+f'(x)t+\frac{f''(x)}{2}\,t^2$$
$$f(x-t)\simeq f(x)-f'(x)t+\frac{f''(x)}{2}\,t^2$$

By making the average, we obtain from this the following formula:
$$\frac{f(x+t)+f(x-t)}{2}\simeq f(x)+\frac{f''(x)}{2}\,t^2$$

Now assume that we have found a way of averaging things over $t\in[-\varepsilon,\varepsilon]$, with the corresponding averages being denoted $I$. We obtain from the above:
$$I(f)\simeq f(x)+f''(x)I\left(\frac{t^2}{2}\right)$$

But this is what our statement says, save for some uncertainties regarding the averaging method, and the precise value of $I(t^2/2)$. We will leave this for later.
\end{proof}

Back to rigorous mathematics now, after this physics intermezzo, as a second application of Theorem 3.21, we can improve as well Theorem 3.11, as follows:

\index{local minimum}
\index{local maximum}

\begin{theorem}
The local minima and local maxima of a twice differentiable function $f:\mathbb R\to\mathbb R$ appear at the points $x\in\mathbb R$ where
$$f'(x)=0$$
with the local minima corresponding to the case $f'(x)\geq0$, and with the local maxima corresponding to the case $f''(x)\leq0$.
\end{theorem}

\begin{proof}
The first assertion is something that we already know. As for the second assertion, we can use the formula in Theorem 3.21, which in the case $f'(x)=0$ reads:
$$f(x+t)\simeq f(x)+\frac{f''(x)}{2}\,t^2$$

Indeed, assuming $f''(x)\neq 0$, it is clear that the condition $f''(x)>0$ will produce a local minimum, and that the condition $f''(x)<0$ will produce a local maximum.
\end{proof}

As before with Theorem 3.11, the above result is not the end of the story with the local minima and maxima, because things are undetermined when:
$$f'(x)=f''(x)=0$$

For instance the functions $\pm x^n$ with $n\in\mathbb N$ all satisfy this condition at $x=0$, which is a minimum for the functions of type $x^{2m}$, a maximum for the functions of type $-x^{2m}$, and not a local minimum or local maximum for the functions of type $\pm x^{2m+1}$.

\bigskip

There are some comments to be made in relation with Algorithm 3.12 as well. Normally that algorithm stays strong, because Theorem 3.23 can only help in relation with the final steps, and is it worth it to compute the second derivative $f''$, just for getting rid of roughly $1/2$ of the $f(x)$ values to be compared. However, in certain cases, this method proves to be useful, so Theorem 3.23 is good to know, when applying that algorithm.

\bigskip

As a main concrete application now of the second derivative, which is something very useful in practice, and related to Interpretation 3.18 (4), we have the following result:

\index{convex function}
\index{concave function}
\index{Jensen inequality}

\begin{theorem}
Given a convex function $f:\mathbb R\to\mathbb R$, we have the following Jensen inequality, for any $x_1,\ldots,x_N\in\mathbb R$, and any $\lambda_1,\ldots,\lambda_N>0$ summing up to $1$,
$$f(\lambda_1x_1+\ldots+\lambda_Nx_N)\leq\lambda_1f(x_1)+\ldots+\lambda_Nx_N$$
with equality when $x_1=\ldots=x_N$. In particular, by taking the weights $\lambda_i$ to be all equal, we obtain the following Jensen inequality, valid for any $x_1,\ldots,x_N\in\mathbb R$,
$$f\left(\frac{x_1+\ldots+x_N}{N}\right)\leq\frac{f(x_1)+\ldots+f(x_N)}{N}$$
and once again with equality when $x_1=\ldots=x_N$. A similar statement holds for the concave functions, with all the inequalities being reversed.
\end{theorem}

\begin{proof}
This is indeed something quite routine, the idea being as follows:

\medskip

(1) First, we can talk about convex functions in a usual, intuitive way, with this meaning by definition that the following inequality must be satisfied:
$$f\left(\frac{x+y}{2}\right)\leq\frac{f(x)+f(y)}{2}$$

(2) But this means, via a simple argument, by approximating numbers $t\in[0,1]$ by sums of powers $2^{-k}$, that for any $t\in[0,1]$ we must have:
$$f(tx+(1-t)y)\leq tf(x)+(1-t)f(y)$$

Alternatively, via yet another simple argument, this time by doing some geometry with triangles, this means that we must have:
$$f\left(\frac{x_1+\ldots+x_N}{N}\right)\leq\frac{f(x_1)+\ldots+f(x_N)}{N}$$

But then, again alternatively, by combining the above two simple arguments, the following must happen, for any $\lambda_1,\ldots,\lambda_N>0$ summing up to $1$:
$$f(\lambda_1x_1+\ldots+\lambda_Nx_N)\leq\lambda_1f(x_1)+\ldots+\lambda_Nx_N$$

(3) Summarizing, all our Jensen inequalities, at $N=2$ and at $N\in\mathbb N$ arbitrary, are equivalent. The point now is that, if we look at what the first Jensen inequality, that we took as definition for the convexity, exactly means, this is simply equivalent to:
$$f''(x)\geq0$$

(4) Thus, we are led to the conclusions in the statement, regarding the convex functions. As for the concave functions, the proof here is similar. Alternatively, we can say that $f$ is concave precisely when $-f$ is convex, and get the results from what we have.
\end{proof}

As a basic application of the Jensen inequality, which is very classical, we have:

\index{Cauchy-Schwarz}

\begin{theorem}
For any $p\in(1,\infty)$ we have the following inequality,
$$\left|\frac{x_1+\ldots+x_N}{N}\right|^p\leq\frac{|x_1|^p+\ldots+|x_N|^p}{N}$$
and for any $p\in(0,1)$ we have the following inequality,
$$\left|\frac{x_1+\ldots+x_N}{N}\right|^p\geq\frac{|x_1|^p+\ldots+|x_N|^p}{N}$$
with in both cases equality precisely when $|x_1|=\ldots=|x_N|$.
\end{theorem}

\begin{proof}
This follows indeed from Theorem 3.24, because we have:
$$(x^p)''=p(p-1)x^{p-2}$$

Thus $x^p$ is convex for $p>1$ and concave for $p<1$, which gives the results.
\end{proof}

Observe that at $p=2$ we obtain as particular case of the above inequality the Cauchy-Schwarz inequality, or rather something equivalent to it, namely:
$$\left(\frac{x_1+\ldots+x_N}{N}\right)^2\leq\frac{x_1^2+\ldots+x_N^2}{N}$$

We will be back to this, Cauchy-Schwarz and versions, later on in this book. Finally, as yet another important application of the Jensen inequality, we have:

\begin{theorem}
We have the following Young inequality,
$$ab\leq \frac{a^p}{p}+\frac{b^q}{q}$$
valid for any $a,b\geq0$, and any exponents $p,q>1$ satisfying $\frac{1}{p}+\frac{1}{q}=1$. 
\end{theorem}

\begin{proof}
We use the logarithm function, which is concave on $(0,\infty)$, due to:
$$(\log x)''=\left(-\frac{1}{x}\right)'=-\frac{1}{x^2}$$

Thus we can apply the Jensen inequality, and we obtain in this way:
\begin{eqnarray*}
\log\left(\frac{a^p}{p}+\frac{b^q}{q}\right)
&\geq&\frac{\log(a^p)}{p}+\frac{\log(b^q)}{q}\\
&=&\log(a)+\log(b)\\
&=&\log(ab)
\end{eqnarray*}

Now by exponentiating, we obtain the Young inequality.
\end{proof}

Observe that for the simplest exponents, namely $p=q=2$, the Young inequality reads $ab\leq(a^2+b^2)/2$. In general, what we have is a useful refinement of this, that can be used for many purposes. We will be back to applications later, in chapter 7.

\section*{3c. The Taylor formula}

Back now to the general theory of the derivatives, and their theoretical applications, we can further develop our basic approximation method, at order 3, at order 4, and so on, the ultimate result on the subject, called Taylor formula, being as follows:

\index{Taylor formula}
\index{L'H\^opital's rule}
\index{higher derivatives}

\begin{theorem}
Assuming that $f:\mathbb R\to\mathbb R$ is differentiable $n$ times, we have
$$f(x+t)\simeq\sum_{k=0}^n\frac{f^{(k)}(x)}{k!}\,t^k$$
where $f^{(k)}(x)$ are the higher derivatives of $f$ at the point $x$.
\end{theorem}

\begin{proof}
Consider indeed the function to be approximated, namely:
$$\varphi(t)=f(x+t)$$

Let us try to best approximate this function at a given order $n\in\mathbb N$. We are therefore looking for a certain polynomial in $t$, of the following type:
$$P(t)=a_0+a_1t+\ldots+a_nt^n$$

The natural conditions to be imposed are those stating that $P$ and $\varphi$ should match at $t=0$, at the level of the actual value, of the derivative, second derivative, and so on up the $n$-th derivative. Thus, we are led to the approximation in the statement:
$$f(x+t)\simeq\sum_{k=0}^n\frac{f^{(k)}(x)}{k!}\,t^k$$

Next, we can prove this approximation $\varphi(t)\simeq P(t)$ by using L'H\^opital's rule:
\begin{eqnarray*}
\frac{\varphi(t)-P(t)}{t^n}
&\simeq&\frac{\varphi'(t)-P'(t)}{nt^{n-1}}\\
&\simeq&\frac{\varphi''(t)-P''(t)}{n(n-1)t^{n-2}}\\
&\vdots&\\
&\simeq&\frac{\varphi^{(n)}(t)-P^{(n)}(t)}{n!}\\
&=&\frac{f^{(n)}(x)-f^{(n)}(x)}{n!}\\
&=&0
\end{eqnarray*}

Thus, we are led to the conclusion in the statement.
\end{proof}

Here is a related interesting statement, inspired from the above proof:

\index{Taylor formula}
\index{polynomial}

\begin{proposition}
For a polynomial of degree $n$, the Taylor approximation
$$f(x+t)\simeq\sum_{k=0}^n\frac{f^{(k)}(x)}{k!}\,t^k$$
is an equality. The converse of this statement holds too.
\end{proposition}

\begin{proof}
By linearity, it is enough to check the equality in question for the monomials $f(x)=x^p$, with $p\leq n$. But here, the formula to be proved is as follows:
$$(x+t)^p\simeq\sum_{k=0}^p\frac{p(p-1)\ldots(p-k+1)}{k!}\,x^{p-k}t^k$$

We recognize the binomial formula, so our result holds indeed. As for the converse, this is clear, because the Taylor approximation is a polynomial of degree $n$.
\end{proof}

In order to further comment now on Theorem 3.27, which remains something quite subtle, it is perhaps time to clarify our meaning of $\simeq$. We have been using this sign, since the beginning of this book, for approximation in a general, intuitive sense. However, since in Theorem 3.27 we have all kinds of infinitesimals appearing, namely $t,t^2,\ldots,t^n$, we must invent something better. And here, the best is to state things as follows:

\begin{theorem}
Assuming that $f:\mathbb R\to\mathbb R$ is differentiable $n$ times we have
$$f(x+t)=\sum_{k=0}^n\frac{f^{(k)}(x)}{k!}\,t^k+o(t^n)$$
with $o$ being the Landau symbol, meaning $o(s)/s\to0$, with $s\to0$.
\end{theorem}

\begin{proof}
This is indeed something self-explanatory, based on Theorem 3.27.
\end{proof}

As a last comment about Theorem 3.27, an interesting situation, which appears quite often, is that when $f$ is infinitely differentiable. Here the result is as follows:

\begin{theorem}
Assuming that $f:\mathbb R\to\mathbb R$ is infinitely differentiable, we have
$$f(x+t)=\sum_{k=0}^n\frac{f^{(k)}(x)}{k!}\,t^k+o(t^n)$$
for any $n\in\mathbb N$, according to the Taylor theorem. However, the asymptotic formula
$$f(x+t)=\sum_{k=0}^\infty\frac{f^{(k)}(x)}{k!}\,t^k$$
might hold or not, depending on $f$, and generically, does not hold.
\end{theorem}

\begin{proof}
This is something quite tricky, the idea being as follows:

\medskip

(1) To start with, the first assertion is something that we know well.

\medskip

(2) The second assertion is something more subtle, with the examples there abounding. However, we have as well counterexamples, with a standard counterexample being:
$$f(x)=\begin{cases}
e^{-1/x^2}&(x\neq0)\\
0&(x=0)
\end{cases}$$

Indeed, for this function we have the following estimate, valid for any $n\in\mathbb N$:
$$f(t)=(e^{1/t^2})^{-1}\leq\left(\frac{1/t^{2n}}{n!}\right)^{-1}=
 n!t^{2n}=o(t^n)$$

Thus $f$ is infinitely differentiable at $0$, with all its derivatives vanishing there, and so its Taylor series at $0$ is the null series, which cannot be equal to $f$ itself. That is, $f$ is designed not to take off from $0$, but it manages however to take off, very slowly.

\medskip

(3) In what regards now the very last claim, this is something more technical, an intuitive explanation here being that there should be more functions $f:\mathbb R\to\mathbb R$, even taken infinitely differentiable, than series $\psi=\sum c_kt^k$. And in practice, up to you to learn here how to count such beasts, as an exercise, and reach to the above conclusion.
\end{proof}

In relation now with the local extrema, and getting back to our usual, informal $\simeq$ convention, in order to quickly explain what happens, we have the following result:

\begin{theorem}
Assuming that $f:\mathbb R\to\mathbb R$ is $n$ times differentiable, and
$$f(x+t)\simeq f(x)+\frac{f^{(n)}(x)}{n!}\,t^n$$
with $f^{(n)}(x)\neq0$, this tells us if $x$ is a local minimum or maximum of $f$.
\end{theorem} 

\begin{proof}
This is a quite compact statement, coming from the Taylor formula, the idea in practice being that we have an algorithm here, as follows:

\medskip

(1) We can start with $n=1$, and with the following formula, that we know well:
$$f(x+t)\simeq f(x)+f'(x)t$$

Indeed, this formula tells us that when $f'(x)\neq0$, the point $x$ cannot be a local minimum or maximum, due to the fact that $t\to-t$ will invert the growth. 

\medskip

(2) In the case left, $f'(x)=0$, we switch to $n=2$, where the Taylor formula is:
$$f(x+t)\simeq f(x)+\frac{f''(x)}{2}\,t^2$$

And here, when $f''(x)<0$ we have a local maximum, and when $f''(x)>0$ we have a local minimum. As for the remaining case, $f''(x)=0$, things here remain open.

\medskip

(3) In the case left, $f''(x)=0$, we switch to $n=3$, where the Taylor formula is:
$$f(x+t)\simeq f(x)+\frac{f'''(x)}{6}\,t^3$$

But this solves the problem in the case $f'''(x)\neq0$, because here we cannot have a local minimum or maximum, due to $t\to-t$, which switches growth. As for the remaining case, $f'''(x)=0$, things here remain open, and we have to go at higher order.

\medskip

(4) Summarizing, we have a recurrence method for solving our problem. In order to comment now on what happens at the $n$-th step, let us write, as in the statement:
$$f(x+t)\simeq f(x)+\frac{f^{(n)}(x)}{n!}\,t^n$$

Then, when $n$ is even, if $f^{(n)}(x)<0$ we have a local maximum, and if $f^{(n)}(x)>0$ we have a local minimum. As for the case where $n$ is odd, here with $f^{(n)}(x)\neq0$ we cannot have a local minimum or maximum, due to $t\to-t$ which switches growth.

\medskip

(5) And so on, until the algorithm stops, either due to $f^{(n)}(x)\neq0$, as it would be ideal, solving our problem, or due to the fact that $f^{(n)}$ is no longer differentiable at $x$, telling us that we have to use some alternative methods, or due the fact that we are facing a tricky function, resisting our algorithm until the very end, after $n=\infty$ steps, such as the function $f(x)=e^{-1/x^2}$, that we met in the proof of Theorem 3.30.
\end{proof}

Getting now to more concrete things, let us compute the Taylor series of the basic functions that we know. We first have here the following result:

\index{Taylor formula}
\index{exp}
\index{log}
\index{sin}
\index{cos}

\begin{theorem}
We have the following formulae,
$$e^x=\sum_{k=0}^\infty\frac{x^k}{k!}\quad,\quad 
\log(1+x)=\sum_{k=0}^\infty(-1)^{k+1}\frac{x^k}{k}$$
as well as the following formulae,
$$\sin x=\sum_{l=0}^\infty(-1)^l\frac{x^{2l+1}}{(2l+1)!}\quad,\quad 
\cos x=\sum_{l=0}^\infty(-1)^l\frac{x^{2l}}{(2l)!}$$
as Taylor series, and in general as well, with $|x|<1$ needed for $\log$.
\end{theorem}

\begin{proof}
There are several statements here, the proofs being as follows:

\medskip

(1) Regarding $\exp$ and $\log$, the needed differentiation formulae are as follows:
$$(e^x)'=e^x\quad,\quad 
(\log x)'=x^{-1}\quad,\quad
(x^p)'=px^{p-1}$$

As for $\sin$ and $\cos$, here the needed differentiation formulae are as follows:
$$(\sin x)'=\cos x\quad,\quad 
(\cos x)'=-\sin x$$

Thus, we can compute, and we obtain the Taylor series in the statement.

\medskip

(2) Regarding now the fact that the formulae in the statement are in fact exact, this is something more subtle. For $\exp$, this is something that we know. In order to deal with $\sin$ and $\cos$, we have the following trick, using a formal number satisfying $i^2=-1$:
$$(\cos x+i\sin x)'=-\sin x+i\cos x=i(\cos x+i\sin x)$$

Thus $f(x)=\cos i+i\sin x$ satisfies $f'=if$, and by arguing like in the proof of Theorem 3.16, we conclude that we have $f(x)=e^{ix}$. But this does the job, because:
\begin{eqnarray*}
e^{ix}
&=&\sum_{k=0}^\infty\frac{(ix)^k}{k!}\\
&=&\sum_{l=0}^\infty\frac{(ix)^{2l}}{(2l)!}+\sum_{l=0}^\infty\frac{(ix)^{2l+1}}{(2l+1)!}\\
&=&\sum_{l=0}^\infty(-1)^l\frac{x^{2l}}{(2l)!}+i\sum_{l=0}^\infty(-1)^l\frac{x^{2l+1}}{(2l+1)!}
\end{eqnarray*}

Of course, I can hear you screaming that this is not rigorous. In answer, this is in fact rigorous, when assuming a better knowledge of $i^2=-1$. More on this in chapter 5.

\medskip

(3) Finally, regarding $\log$, if we set $f(1+x)=\sum_k(-1)^{k+1}x^k/k$, we have the following computation, using the generalized binomial formula, with exponent $-r$:
\begin{eqnarray*}
f((1+x)(1+y))
&=&\sum_{k=1}^\infty(-1)^{k+1}\frac{(x+y+xy)^k}{k}\\
&=&\sum_{k=1}^\infty\frac{(-1)^{k+1}}{k}\sum_{r=0}^k\binom{k}{r}(x+xy)^ry^{k-r}\\
&=&\sum_{r+s\geq1}\frac{(-1)^{r+s+1}}{r+s}\binom{r+s}{r}x^r(1+y)^ry^s\\
&=&\sum_{r>1}(-1)^{r+1}\frac{x^r(1+y)^r}{r}\cdot\frac{1}{(1+y)^r}+\sum_{s\geq1}(-1)^{s+1}\frac{y^s}{s}\\
&=&f(1+x)+f(1+y)
\end{eqnarray*}

Thus $f$ satisfies the functional equation of $\log$, and since the slope at 1 is the correct one, we can proceed as in chapter 2 for the exponential, and we get $f=\log$, as stated.
\end{proof}

Next, let us record the result for the arctangent, which is as follows:

\begin{theorem}
The Taylor series of $\arctan$, which is exact too, is
$$\arctan x=\sum_{k=0}^\infty\frac{(-1)^k}{2k+1}\,x^{2k+1}$$
and in particular, with $x=1$ we obtain the Leibnitz formula
$$\frac{\pi}{4}=1-\frac{1}{3}+\frac{1}{5}-\frac{1}{7}+\frac{1}{9}-\frac{1}{11}+\ldots$$
which can be used in order to find the decimals of $\pi$.
\end{theorem}

\begin{proof}
We have indeed the following computation, based on Proposition 3.10:
$$(\arctan x)'=\frac{1}{1+x^2}=\sum_{k=0}^\infty(-1)^kx^{2k}$$

But this gives the formula in the statement, and the Leibnitz formula follows too.
\end{proof}

As for the remaining 9 trigonometric functions, the Taylor series can be computed too, and are exact too, but with the difficulty varying, depending on the function. You can find formulae and details in any dedicated 1-variable book, including mine \cite{ba1}.

\section*{3d. Differential equations}

Good news, with the calculus that we know we can do some physics, in 1 dimension. Let us start with something immensely important, in the history of science:

\begin{fact}
Newton invented calculus for formulating the laws of motion as
$$v=\dot{x}\quad,\quad a=\dot{v}$$
where $x,v,a$ are the position, speed and acceleration, and the dots are time derivatives.
\end{fact}

To be more precise, the variable in Newton's physics is time $t\in\mathbb R$, playing the role of the variable $x\in\mathbb R$ that we have used in the above. And we are looking at a particle whose position is described by a function $x=x(t)$. Then, it is quite clear that the speed of this particle should be described by the first derivative $v=x'(t)$, and that the acceleration of the particle should be described by the second derivative $a=v'(t)=x''(t)$.

\bigskip

Summarizing, with Newton's theory of derivatives, as we learned it in this chapter, we can certainly do some mathematics for the motion of bodies. But, for these bodies to move, we need them to be acted upon by some forces, right? The simplest such force is gravity, and in our present, modest 1 dimensional setting, we have:

\index{differential equation}
\index{gravity}

\begin{theorem}
The equation of a gravitational free fall, in $1$ dimension, is
$$\ddot{x}=-\frac{GM}{x^2}$$
with $M$ being the attracting mass, and $G\simeq 6.674\times 10^{-11}$ being a constant.
\end{theorem}

\begin{proof}
Assume indeed that we have a free falling object, in $1$ dimension:
$$\xymatrix@R=20pt@C=10pt{
\circ_m\ar[d]\\
\bullet_M}$$

In order to reach to calculus as we know it, we must peform a rotation, as to have all this happening on the $Ox$ axis. By doing this, and assuming that $M$ is fixed at $0$, our picture becomes as follows, with the attached numbers being now the coordinates:
$$\xymatrix@R=20pt@C=20pt{
\bullet_0&\circ_x\ar[l]}$$

Now comes the physics. The gravitational force exterted by $M$, which is fixed in our formalism, on the object $m$ which moves, is subject to the following equations:
$$F=-G\cdot\frac{Mm}{x^2}\quad,\quad F=ma\quad,\quad a=\dot{v}\quad,\quad v=\dot{x}$$

We are therefore led to the equation of motion in the statement.
\end{proof}

As more physics, we can talk as well about waves in 1 dimension, as follows:

\index{wave equation}
\index{lattice model}
\index{Hooke law}
\index{Newton law}

\begin{theorem}
The wave equation in $1$ dimension is
$$\ddot{f}=v^2f''$$
with the dot denoting time derivatives, and $v>0$ being the propagation speed.
\end{theorem}

\begin{proof}
In order to understand the propagation of the waves, let us model the space, which is $\mathbb R$ for us, as a network of balls, with springs between them, as follows:
$$\cdots\times\!\!\!\times\!\!\!\times\bullet\times\!\!\!\times\!\!\!\times\bullet\times\!\!\!\times\!\!\!\times\bullet\times\!\!\!\times\!\!\!\times\bullet\times\!\!\!\times\!\!\!\times\bullet\times\!\!\!\times\!\!\!\times\cdots$$

Now let us send an impulse, and see how balls will be moving. For this purpose, we zoom on one ball. The situation here is as follows, $l$ being the spring length:
$$\cdots\cdots\bullet_{f(x-l)}\times\!\!\!\times\!\!\!\times\bullet_{f(x)}\times\!\!\!\times\!\!\!\times\bullet_{f(x+l)}\cdots\cdots$$

We have two forces acting at $x$. First is the Newton motion force, mass times acceleration, which is as follows, with $m$ being the mass of each ball:
$$F_n=m\cdot\ddot{f}(x)$$

And second is the Hooke force, displacement of the spring, times spring constant. Since we have two springs at $x$, this is as follows, $k$ being the spring constant:
\begin{eqnarray*}
F_h
&=&F_h^r-F_h^l\\
&=&k(f(x+l)-f(x))-k(f(x)-f(x-l))\\
&=&k(f(x+l)-2f(x)+f(x-l))
\end{eqnarray*}

We conclude that the equation of motion, in our model, is as follows:
$$m\cdot\ddot{f}(x)=k(f(x+l)-2f(x)+f(x-l))$$

Now let us take the limit of our model, as to reach to continuum. For this purpose we will assume that our system consists of $N>>0$ balls, having a total mass $M$, and spanning a total distance $L$. Thus, our previous infinitesimal parameters are as follows, with $K$ being the spring constant of the total system, which is of course lower than $k$:
$$m=\frac{M}{N}\quad,\quad k=KN\quad,\quad l=\frac{L}{N}$$

With these changes, our equation of motion found in (1) reads:
$$\ddot{f}(x)=\frac{KN^2}{M}(f(x+l)-2f(x)+f(x-l))$$

Now with $N\to\infty$, and therefore with $l\to0$, we obtain in this way:
$$\ddot{f}(x)=\frac{KL^2}{M}\cdot f''(x)$$

Thus, we are led to the conclusion in the statement.
\end{proof}

Along the same lines, we can talk as well about heat in 1D, as follows:

\index{heat equation}
\index{lattice model}

\begin{theorem}
The heat diffusion equation in $1$ dimension is
$$\dot{f}=\alpha f''$$
where $\alpha>0$ is the thermal diffusivity of the medium.
\end{theorem}

\begin{proof}
As before with the wave equation, this is not exactly a theorem, but rather what comes out of experiments, but we can justify this mathematically, as follows:

\medskip

(1) As an intuitive explanation for this equation, since the second derivative $f''$ computes the average value of a function $f$ around a point, minus the value of $f$ at that point, as we know from Proposition 3.22, the heat equation as formulated above tells us that the rate of change $\dot{f}$ of the temperature of the material at any given point must be proportional, with proportionality factor $\alpha>0$, to the average difference of temperature between that given point and the surrounding material. Which sounds reasonable.

\medskip

(2) In practice now, we can use, a bit like before for the wave equation, a lattice model as follows, with distance $l>0$ between the neighbors:
$$\xymatrix@R=10pt@C=20pt{
\ar@{-}[r]&\circ_{x-l}\ar@{-}[r]^l&\circ_x\ar@{-}[r]^l&\circ_{x+l}\ar@{-}[r]&
}$$

In order to model now heat diffusion, we have to implement the intuitive mechanism explained above, and in practice, this leads to a condition as follows, expressing the change of the temperature $f$, over a small period of time $\delta>0$:
$$f(x,t+\delta)=f(x,t)+\frac{\alpha\delta}{l^2}\sum_{x\sim y}\left[f(y,t)-f(x,t)\right]$$

But this leads, via manipulations as before, to $\dot{f}(x,t)=\alpha f''(x,t)$, as claimed.
\end{proof}

All this is very nice, so with the calculus that we know, we can certainly talk about physics. We will see later in this book how to deal with the above equations.

\section*{3e. Exercises}

Here are some basic exercises, in relation with the above:

\begin{exercise}
Work out the basic estimates for $\sin,\cos,\tan$ at $t=0$.
\end{exercise}

\begin{exercise}
Work out all details in the proof of the Jensen inequality.
\end{exercise}

\begin{exercise}
Learn a bit about analytic and non-analytic functions.
\end{exercise}

\begin{exercise}
$\tan$, $\sec$, $\csc$, $\cot$, $\arcsin$, $\arccos$, ${\rm arcsec}$, ${\rm arccsc}$, ${\rm arccot}$.
\end{exercise}

As a bonus exercise, try solving the gravity, wave and heat equations in 1D.

\chapter{Integration}

\section*{4a. Integration theory}

We have seen so far the foundations of calculus, with lots of interesting results regarding the functions $f:\mathbb R\to\mathbb R$, and their derivatives $f':\mathbb R\to\mathbb R$. The general idea was that in order to understand $f$, we first need to compute its derivative $f'$. The overall conclusion, coming from the Taylor formula, was that if we are able to compute $f'$, but then also $f''$, and $f'''$ and so on, we will have a good understanding of $f$ itself.

\bigskip 

However, the story is not over here, and there is one more twist to the plot. Which will be a major twist, of similar magnitude to that of the Taylor formula. For reasons which are quite tricky, that will become clear later on, we will be interested in the integration of the functions $f:\mathbb R\to\mathbb R$. With the claim that this is related to calculus.

\bigskip

There are several possible viewpoints on the integral, which are all useful, and good to know. To start with, we have something very simple, as follows:

\index{integral}
\index{area}

\begin{definition}
The integral of a continuous function $f:[a,b]\to\mathbb R$, denoted
$$\int_a^bf(x)dx$$
is the area below the graph of $f$, signed $+$ where $f\geq0$, and signed $-$ where $f\leq0$.
\end{definition}

Here it is of course understood that the area in question can be computed, and with this being something quite subtle, that we will get into later. For the moment, let us just trust our intuition, our function $f$ being continuous, the area in question can ``obviously'' be computed. More on this later, but for being rigorous, however, let us formulate:

\index{volume}

\begin{method}
In practice, the integral of $f\geq 0$ can be computed as follows,
\begin{enumerate}
\item Cut the graph of $f$ from 3mm plywood,

\item Plunge that graph into a square container of water,

\item Measure the water displacement, as to have the volume of the graph,

\item Divide by $3\times 10^{-3}$ that volume, as to have the area,
\end{enumerate} 
and for general $f$, we can use this plus $f=f_+-f_-$, with $f_+,f_-\geq0$.
\end{method}

So far, so good, we have a rigorous definition, so let us do now some computations. In order to compute areas, and so integrals of functions, without wasting precious water, we can use our geometric knowledge. Here are some basic results of this type:

\index{piecewise linear}

\begin{proposition}
We have the following results:
\begin{enumerate}
\item When $f$ is linear, we have the following formula:
$$\int_a^bf(x)dx=(b-a)\cdot\frac{f(a)+f(b)}{2}$$

\item In fact, when $f$ is piecewise linear on $[a=a_1,a_2,\ldots,a_n=b]$, we have:
$$\int_a^bf(x)dx=\sum_{i=1}^{n-1}(a_{i+1}-a_i)\cdot\frac{f(a_i)+f(a_{i+1})}{2}$$

\item We have as well the formula $\int_{-1}^1\sqrt{1-x^2}\,dx=\pi/2$.
\end{enumerate}
\end{proposition}

\begin{proof}
These results all follow from basic geometry, as follows:

\medskip

(1) Assuming $f\geq0$, we must compute the area of a trapezoid having sides $f(a),f(b)$, and height $b-a$. But this is the same as the area of a rectangle having side $(f(a)+f(b))/2$ and height $b-a$, and we obtain $(b-a)(f(a)+f(b))/2$, as claimed.

\medskip

(2) This is clear indeed from the formula found in (1), by additivity.

\medskip

(3) The integral in the statement is by definition the area of the upper unit half-disc. But since the area of the whole unit disc is $\pi$, this half-disc area is $\pi/2$.
\end{proof}

As an interesting observation, (2) in the above result makes it quite clear that $f$ does not necessarily need to be continuous, in order to talk about its integral. Indeed, assuming that $f$ is piecewise linear on $[a=a_1,a_2,\ldots,a_n=b]$, but not necessarily continuous, we can still talk about its integral, in the obvious way, exactly as in Definition 4.1, and we have an explicit formula for this integral, generalizing the one found in (2), namely:
$$\int_a^bf(x)dx=\sum_{i=1}^{n-1}(a_{i+1}-a_i)\cdot\frac{f(a_i^+)+f(a_{i+1}^-)}{2}$$

Based on this observation, let us upgrade our formalism, as follows:

\index{integrable function}

\begin{definition}
We say that a function  $f:[a,b]\to\mathbb R$ is integrable when the area below its graph is computable. In this case we denote by
$$\int_a^bf(x)dx$$
this area, signed $+$ where $f\geq0$, and signed $-$ where $f\leq0$.
\end{definition}

As basic examples of integrable functions, we have the continuous ones, provided indeed that our intuition, or that Method 4.2, works indeed for any such function. We will soon see that this is indeed true, coming with mathematical proof. As further examples, we have the functions which are piecewise linear, or more generally piecewise continuous. We will also see, later, as another class of examples, that the piecewise monotone functions are integrable. But more on this later, let us not bother for the moment with all this.

\bigskip

This being said, one more thing regarding theory, that you surely have in mind: is any function integrable? Not clear. I would say that if the Devil comes with some sort of nasty, totally discontinuous function $f:\mathbb R\to\mathbb R$, then you will have big troubles in cutting its graph from 3mm plywood, as required by Method 4.2. More on this later.

\bigskip

Back to work now, here are some general results regarding the integrals:

\begin{proposition}
We have the following formulae,
$$\int_a^bf(x)+g(x)dx=\int_a^bf(x)dx+\int_a^bg(x)dx$$
$$\int_a^b\lambda f(x)=\lambda\int_a^bf(x)$$
valid for any functions $f,g$ and any scalar $\lambda\in\mathbb R$.
\end{proposition}

\begin{proof}
Both these formulae are indeed clear from definitions.
\end{proof}

Moving ahead now, passed the above results, which are of purely algebraic and geometric nature, and perhaps a few more of the same type, which are all quite trivial and that we we will not get into here, we must do some analysis, in order to compute integrals. This is something quite tricky, and we have here the following result:

\index{Riemann integration}

\begin{theorem}
We have the Riemann integration formula,
$$\int_a^bf(x)dx=(b-a)\times\lim_{N\to\infty}\frac{1}{N}\sum_{k=1}^Nf\left(a+\frac{b-a}{N}\cdot k\right)$$
which can serve as a definition for the integral.
\end{theorem}

\begin{proof}
This is standard, by drawing rectangles. We have indeed the following formula, which can stand as a definition for the signed area below the graph of $f$:
$$\int_a^bf(x)dx=\lim_{N\to\infty}\sum_{k=1}^N\frac{b-a}{N}\cdot f\left(a+\frac{b-a}{N}\cdot k\right)$$

Thus, we are led to the formula in the statement.
\end{proof}

Observe that the above formula suggests that $\int_a^bf(x)dx$ is the length of the interval $[a,b]$, namely $b-a$, times the average of $f$ on the interval $[a,b]$. Thinking a bit, this is indeed something true, with no need for Riemann sums, coming directly from Definition 4.1, because area means side times average height. Thus, we can formulate:

\index{average of function}

\begin{theorem}
The integral of a function $f:[a,b]\to\mathbb R$ is given by
$$\int_a^bf(x)dx=(b-a)\times A(f)$$
where $A(f)$ is the average of $f$ over the interval $[a,b]$.
\end{theorem}

\begin{proof}
As explained above, this is clear from Definition 4.1, via some geometric thinking. Alternatively, this is something which certainly comes from Theorem 4.6.
\end{proof}

The point of view in Theorem 4.7 is something quite useful, and as an illustration for this, let us review the results that we already have, by using this interpretation. First, we have the formula for linear functions from Proposition 4.3, namely:
$$\int_a^bf(x)dx=(b-a)\cdot\frac{f(a)+f(b)}{2}$$

But this formula is totally obvious with our new viewpoint, from Theorem 4.7. The same goes for the results in Proposition 4.5, which become even more obvious with the viewpoint from Theorem 4.7. However, not everything trivializes in this way, and the result which is left, namely the formula $\int_{-1}^1\sqrt{1-x^2}\,dx=\pi/2$ from Proposition 4.3 (3), not only does not trivialize, but becomes quite opaque with our new philosophy.

\bigskip

In short, modesty. Integration is a quite delicate business, and we have several equivalent points of view on what an integral means, and all these points of view are useful, and must be learned, with none of them being clearly better than the others.

\bigskip

Going ahead with more interpretations of the integral, we have:

\index{Monte Carlo integration}
\index{random number}

\begin{theorem}
We have the Monte Carlo integration formula,
$$\int_a^bf(x)dx=(b-a)\times\lim_{N\to\infty}\frac{1}{N}\sum_{k=1}^Nf(x_i)$$
with $x_1,\ldots,x_N\in[a,b]$ being random.
\end{theorem}

\begin{proof}
We recall from Theorem 4.6 that the idea is that we have a formula as follows, with the points $x_1,\ldots,x_N\in[a,b]$ being uniformly distributed: 
$$\int_a^bf(x)dx=(b-a)\times\lim_{N\to\infty}\frac{1}{N}\sum_{k=1}^Nf(x_i)$$

But this works as well when the points $x_1,\ldots,x_N\in[a,b]$ are randomly distributed, for somewhat obvious reasons, and this gives the result.
\end{proof}

Observe that Monte Carlo integration works better than Riemann integration, for instance when trying to improve the estimate, via $N\to N+1$. Indeed, in the context of Riemann integration, assume that we managed to find an estimate as follows, which in practice requires computing $N$ values of our function $f$, and making their average:
$$\int_a^bf(x)dx\simeq\frac{b-a}{N}\sum_{k=1}^Nf\left(a+\frac{b-a}{N}\cdot k\right)$$

In order to improve this estimate, any extra computed value of our function $f(y)$ will be unuseful. For improving our formula, what we need are $N$ extra values of our function, $f(y_1),\ldots,f(y_N)$, with the points $y_1,\ldots,y_N$ being the midpoints of the previous division of $[a,b]$, so that we can write an improvement of our formula, as follows:
$$\int_a^bf(x)dx\simeq\frac{b-a}{2N}\sum_{k=1}^{2N}f\left(a+\frac{b-a}{2N}\cdot k\right)$$

With Monte Carlo, things are far more flexible. Assume indeed that we managed to find an estimate as follows, which again requires computing $N$ values of our function:
$$\int_a^bf(x)dx\simeq\frac{b-a}{N}\sum_{k=1}^Nf(x_i)$$

Now if we want to improve this, any extra computed value of our function $f(y)$ will be helpful, because we can set $x_{n+1}=y$, and improve our estimate as follows:
$$\int_a^bf(x)dx\simeq\frac{b-a}{N+1}\sum_{k=1}^{N+1}f(x_i)$$

And isn't this potentially useful, and powerful, when thinking at practically computing integrals, either by hand, or by using a computer. Let us record this finding as follows:

\begin{conclusion}
Monte Carlo integration works better than Riemann integration, when it comes to computing as usual, by estimating, and refining the estimate.
\end{conclusion}

As another interesting feature of Monte Carlo integration, this works better than Riemann integration, for functions having various symmetries, because Riemann integration can get ``fooled'' by these symmetries, while Monte Carlo remains strong. 

\bigskip

As an example for this phenomeon, chosen to be quite drastic, let us attempt to integrate, via both Riemann and Monte Carlo, the following function $f:[0,\pi]\to\mathbb R$:
$$f(x)=\Big|\sin(120x)\Big|$$

The first few Riemann sums for this function are then as follows:
$$I_2(f)=\frac{\pi}{2}(|\sin 0|+|\sin 60\pi|)=0$$
$$I_3(f)=\frac{\pi}{3}(|\sin 0|+|\sin 40\pi|+|\sin 80\pi|)=0$$
$$I_4(f)=\frac{\pi}{4}(|\sin 0|+|\sin 30\pi|+|\sin 60\pi|+|\sin 90\pi|)=0$$
$$I_5(f)=\frac{\pi}{5}(|\sin 0|+|\sin 24\pi|+|\sin 48\pi|+|\sin 72\pi|+|\sin 96\pi|)=0$$
$$I_6(f)=\frac{\pi}{6}(|\sin 0|+|\sin 20\pi|+|\sin 40\pi|+|\sin 60\pi|+|\sin 80\pi|+|\sin 100\pi|)=0$$
$$\vdots$$

Based on this evidence, we will conclude, obviously, that we have:
$$\int_0^\pi f(x)dx=0$$

With Monte Carlo, however, such things cannot happen. Indeed, since there are finitely many points $x\in[0,\pi]$ having the property $\sin(120 x)=0$, a random point $x\in[0,\pi]$ will have the property $|\sin(120 x)|>0$, so Monte Carlo will give, at any $N\in\mathbb N$:
$$\int_0^\pi f(x)dx\simeq\frac{b-a}{N}\sum_{k=1}^Nf(x_i)>0$$

Again, this is something interesting, when practically computing integrals, either by hand, or by using a computer. So, let us record, as a complement to Conclusion 4.9:

\begin{conclusion}
Monte Carlo integration is smarter than Riemann integration, because the symmetries of the function can fool Riemann, but not Monte Carlo.
\end{conclusion}

All this is good to know, when computing integrals in practice, especially with a computer. Finally, here is one more useful interpretation of the integral:

\index{random variable}
\index{expectation}

\begin{theorem}
The integral of a function $f:[a,b]\to\mathbb R$ is given by
$$\int_a^bf(x)dx=(b-a)\times E(f)$$
where $E(f)$ is the expectation of $f$, regarded as random variable.
\end{theorem}

\begin{proof}
This is just some sort of fancy reformulation of Theorem 4.7, the idea being that what we can ``expect'' from a random variable is of course its average. We will be back to this later in this book, when systematically discussing probability theory.
\end{proof}

\section*{4b. Riemann sums}

Our purpose now will be to understand which functions $f:\mathbb R\to\mathbb R$ are integrable, and how to compute their integrals. For this purpose, the Riemann formula in Theorem 4.6 will be our favorite tool. Let us begin with some theory. We first have:

\index{piecewise continuous}
\index{piecewise monotone}

\begin{theorem}
The following functions are integrable:
\begin{enumerate}
\item The piecewise continuous functions.

\item The piecewise monotone functions.
\end{enumerate}
\end{theorem}

\begin{proof}
This is indeed something quite standard, as follows:

\medskip

(1) It is enough to prove the first assertion for a function $f:[a,b]\to\mathbb R$ which is continuous, and our claim here is that this follows from the uniform continuity of $f$. To be more precise, given $\varepsilon>0$, let us choose $\delta>0$ such that the following happens:
$$|x-y|<\delta\implies|f(x)-f(y)|<\varepsilon$$

In order to prove the result, let us pick two divisions of $[a,b]$, as follows:
$$I=[a=a_1<a_2<\ldots<a_n=b]$$
$$I'=[a=a_1'<a_2'<\ldots<a_m'=b]$$

Our claim, which will prove the result, is that if these divisions are sharp enough, of resolution $<\delta/2$, then the associated Riemann sums $\Sigma_I(f),\Sigma_{I'}(f)$ are close within $\varepsilon$:
$$a_{i+1}-a_i<\frac{\delta}{2}\ ,\ a_{i+1}'-a_i'<\delta_2\ \implies\ \big|\Sigma_I(f)-\Sigma_{I'}(f)\big|<\varepsilon$$

(2) In order to prove this claim, let us denote by $l$ the length of the intervals on the real line. Our assumption is that the lengths of the divisions $I,I'$ satisfy:
$$l\big([a_i,a_{i+1}]\big)<\frac{\delta}{2}\quad,\quad l\big([a_i',a_{i+1}']\big)<\frac{\delta}{2}$$

Now let us intersect the intervals of our divisions $I,I'$, and set:
$$l_{ij}=l\big([a_i,a_{i+1}]\cap[a_j',a_{j+1}']\big)$$

The difference of Riemann sums that we are interested in is then given by:
\begin{eqnarray*}
\Big|\Sigma_I(f)-\Sigma_{I'}(f)\Big|
&=&\left|\sum_{ij}l_{ij}f(a_i)-\sum_{ij}l_{ij}f(a_j')\right|\\
&=&\left|\sum_{ij}l_{ij}(f(a_i)-f(a_j'))\right|
\end{eqnarray*}

(3) Now let us estimate $f(a_i)-f(a_j')$. Since in the case $l_{ij}=0$ we do not need this estimate, we can assume $l_{ij}>0$. Now by remembering what the definition of the numbers $l_{ij}$ was, we conclude that we have at least one point $x\in\mathbb R$ satisfying:
$$x\in [a_i,a_{i+1}]\cap[a_j',a_{j+1}']$$

But then, by using this point $x$ and our assumption on $I,I'$ involving $\delta$, we get:
\begin{eqnarray*}
|a_i-a_j'|
&\leq&|a_i-x|+|x-a_j'|\\
&\leq&\frac{\delta}{2}+\frac{\delta}{2}\\
&=&\delta
\end{eqnarray*}

Thus, according to our definition of $\delta$ from (1), in relation to $\varepsilon$, we get:
$$|f(a_i)-f(a_j')|<\varepsilon$$

(4) But this is what we need, in order to finish. Indeed, with the estimate that we found, we can finish the computation started in (2), as follows:
\begin{eqnarray*}
\Big|\Sigma_I(f)-\Sigma_{I'}(f)\Big|
&=&\left|\sum_{ij}l_{ij}(f(a_i)-f(a_j'))\right|\\
&\leq&\varepsilon\sum_{ij}l_{ij}\\
&=&\varepsilon(b-a)
\end{eqnarray*}

Thus our two Riemann sums are close enough, provided that they are both chosen to be fine enough, and this finishes the proof of the first assertion.

\medskip

(5) Regarding now the second assertion, this is something more technical, that we will not really need in what follows. We will leave the proof here, which uses similar ideas to those in the proof of (1) above, namely subdivisions and estimates, as an exercise.
\end{proof}

Going ahead with more theory, let us establish some abstract properties of the integration operation. We already know from Proposition 4.5 that the integrals behave well with respect to sums and multiplication by scalars. Along the same lines, we have:

\begin{theorem}
The integrals behave well with respect to taking limits,
$$\int_a^b\left(\lim_{n\to\infty}f_n(x)\right)dx=\lim_{n\to\infty}\int_a^bf_n(x)dx$$
and with respect to taking infinite sums as well,
$$\int_a^b\left(\sum_{n=0}^\infty f_n(x)\right)dx=\sum_{n=0}^\infty\int_a^bf_n(x)dx$$
with both these formulae being valid, under mild assumptions.
\end{theorem}

\begin{proof}
This is something quite standard, by using the general theory developed in chapter 3 for the sequences and series of functions. To be more precise, (1) follows by using the material there, via Riemann sums, and then (2) follows as a particular case of (1). We will leave the clarification of all this as an instructive exercise.
\end{proof}

Finally, still at the general level, let us record as well the following result:

\index{mean value property}

\begin{theorem}
Given a continuous function $f:[a,b]\to\mathbb R$, we have
$$\exists c\in[a,b]\quad,\quad \int_a^bf(x)dx=(b-a)f(c)$$
with this being called mean value property.
\end{theorem}

\begin{proof}
Our claim is that this follows from the following trivial estimate:
$$\min(f)\leq f\leq\max(f)$$

Indeed, by integrating this over $[a,b]$, we obtain the following estimate:
$$(b-a)\min(f)\leq\int_a^bf(x)dx\leq(b-a)\max(f)$$

Now observe that this latter estimate can be written as follows:
$$\min(f)\leq\frac{\int_a^bf(x)dx}{b-a}\leq\max(f)$$

Since $f$ must takes all values on $[\min(f),\max(f)]$, we get a $c\in[a,b]$ such that:
$$ \frac{\int_a^bf(x)dx}{b-a}=f(c)$$

Thus, we are led to the conclusion in the statement.
\end{proof}

At the level of examples now, let us first look at the simplest functions that we know, namely the power functions $f(x)=x^p$. However, things here are tricky, as follows:

\index{Riemann sum}

\begin{theorem}
We have the integration formula
$$\int_a^bx^pdx=\frac{b^{p+1}-a^{p+1}}{p+1}$$
valid for any $p\in\mathbb N$.
\end{theorem}

\begin{proof}
This is something quite tricky, the idea being as follows:

\medskip

(1) By linearity we can assume that our interval $[a,b]$ is of the form $[0,c]$, and the formula that we want to establish is as follows:
$$\int_0^cx^pdx=\frac{c^{p+1}}{p+1}$$

(2) We can further assume $c=1$, and by expressing the left term as a Riemann sum, we are in need of the following estimate, in the $N\to\infty$ limit:
$$1^p+2^p+\ldots+N^p\simeq \frac{N^{p+1}}{p+1}$$

(3) So, let us try to prove this. At $p=0$, obviously nothing to do, because we have the following formula, which is exact, and which proves our estimate:
$$1^0+2^0+\ldots+N^0=N$$

(4) At $p=1$ now, we are confronted with a well-known question, namely the computation of $1+2+\ldots+N$. But this is simplest done by arguing that the average of the numbers $1,2,\ldots,N$ being the number in the middle, we have:
$$\frac{1+2+\ldots+N}{N}=\frac{N+1}{2}$$

Thus, we obtain the following formula, which again solves our question:
$$1+2+\ldots+N=\frac{N(N+1)}{2}\simeq\frac{N^2}{2}$$

(5) At $p=2$ now, go compute $1^2+2^2+\ldots+N^2$. This is not obvious, so as a preliminary, let us go back to the case $p=1$, and try to find a new proof there, which might have some chances to extend at $p=2$. And here, we have the following trick:
$$\xymatrix@R=10pt@C=10pt{
\ar@{-}[ddddd]_N\ar@{-}[rrrrrr]&&&&&\ar@{-}[d]&\ar@{-}[ddddd]\\
&&&&\ar@{-}[d]&\ar@{-}[l]\\
&&&\ar@{-}[d]&\ar@{-}[l]&\\
&&\ar@{-}[d]&\ar@{-}[l]&\\
&\ar@{-}[d]&\ar@{-}[l]&\\
\ar@{-}[rrrrrr]_{N+1}&&&&&&}$$

Now the point is that this trick works at $p=2$ too. Indeed, if we consider the 3D shape $P$ formed by a succession of solid squares, having sizes $1\times 1$, $2\times2$, $3\times3$, and so on up to $N\times N$, if we stack 6 copies of $P$ we get a parallelepiped, which gives:
$$1^2+2^2+\ldots+N^2=\frac{N(N+1)(2N+1)}{6}\simeq\frac{N^3}{3}$$

Alternatively, you can get this by computing $1^2+2^2+\ldots+N^2$ for small values of $N$, then conjecturing the above formula, and proving your conjecture by recurrence.

\medskip

(6) At $p=3$ now, the legend has it that by deeply thinking in 4D we are led to the following formula, a bit as in the cases $p=1,2$, explained above:
$$1^3+2^3+\ldots+N^3=\frac{N^2(N+1)^2}{4}\simeq\frac{N^4}{4}$$

Alternatively, assuming that the gods of combinatorics are with us, we can see right away the following formula, which coupled with (2) gives the result:
$$1^3+2^3+\ldots+N^3=(1+2+\ldots+N)^2$$

In any case, in practice, the above formula holds indeed, and you can of course come upon it via some numerics too, and then check it by recurrence.

\medskip

(7) All this is very nice, but coming now as bad news, at $p=4,5,6,\ldots$ the situation is more complicated, with the formulae, provable by recurrence, being as follows:
$$1^4+2^4+\ldots+N^4=\frac{N(N+1)(2N+1)(3N^2+3N-1)}{30}\simeq\frac{N^5}{5}$$
$$1^5+2^5+\ldots+N^5=\frac{N^2(N+1)^2(2N^2+2N+1)}{12}\simeq\frac{N^6}{6}$$
$$1^6+2^6+\ldots+N^6=\frac{N(N+1)(2N+1)(3N^4+6N^3-3N+1)}{42}\simeq\frac{N^7}{7}$$
$$\vdots$$

(8) In fact, we have the following formula, valid at any $p\in\mathbb N$, provable by recurrence, making appear certain mysterious numbers $B_k\in\mathbb Q$, called Bernoulli numbers:
$$1^p+2^p+\ldots+N^p=\frac{1}{p+1}\sum_{k=0}^p(-1)^k\binom{p+1}{k}B_kN^{p+1-k}$$

To be more precise, this formula is indeed compatible with those above at $p=0,\ldots,6$, with the numeric data for the first Bernoulli numbers being as follows:
$$B_0=1\ ,\quad B_1=-\frac{1}{2}\ ,\quad B_2=\frac{1}{6}\ ,\quad B_3=0\ ,\quad 
B_4=-\frac{1}{30}\ ,\quad B_5=0\ ,\quad B_6=\frac{1}{42}$$

As for the proof, in general, when trying to establish this by recurrence, we are led into the recurrence formula for the Bernoulli numbers, which is as follows:
$$\sum_{k=0}^m\binom{m+1}{k}B_k=\delta_{m0}$$

And we will leave clarifying this, and some further learning here, as an exercise. 

\medskip

(9) Now the point is that the above formula does the job for our integration purposes, because we obtain right away the following estimate, exactly as needed:
$$1^p+2^p+\ldots+N^p\simeq\frac{N^{p+1}}{p+1}$$

Summarizing, theorem proved, with the cases $p=0,1,2,3$ being elementary, then $p=4,5,6,\ldots$ being doable too, and $p\in\mathbb N$ general requiring Bernoulli numbers. 
\end{proof}

Quite interesting all this, and based on what we have, let us formulate:

\begin{conjecture}
We have the following integration formula,
$$\int_a^bx^pdx=\frac{b^{p+1}-a^{p+1}}{p+1}$$
valid for any exponent $p\in\mathbb R$.
\end{conjecture}

Now, instead of struggling with this conjecture, let us look at some other functions, which are not polynomial. And here, as good news, we have: 

\begin{theorem}
We have the following integration formula,
$$\int_a^be^xdx=e^b-e^a$$
valid for any two real numbers $a<b$.
\end{theorem}

\begin{proof}
This follows indeed from the Riemann integration formula, because:
\begin{eqnarray*}
\int_a^be^xdx
&=&\lim_{N\to\infty}\frac{e^a+e^{a+(b-a)/N}+e^{a+2(b-a)/N}+\ldots+e^{a+(N-1)(b-a)/N}}{N}\\
&=&\lim_{N\to\infty}\frac{e^a}{N}\cdot\left(1+e^{(b-a)/N}+e^{2(b-a)/N}+\ldots+e^{(N-1)(b-a)/N}\right)\\
&=&\lim_{N\to\infty}\frac{e^a}{N}\cdot\frac{e^{b-a}-1}{e^{(b-a)/N}-1}\\
&=&(e^b-e^a)\lim_{N\to\infty}\frac{1}{N(e^{(b-a)/N}-1)}\\
&=&e^b-e^a
\end{eqnarray*}

Thus, we are led to the conclusion in the statement.
\end{proof}

\section*{4c. Advanced results}

The problem is now, what to do with what we have, namely Conjecture 4.16 and Theorem 4.17. Not obvious, so stuck, and time to ask the cat. And cat says:

\index{infinitesimal}

\begin{cat}
Summing the infinitesimals of the rate of change of the function should give you the global change of the function. Obvious.
\end{cat}

Which is quite puzzling, usually my cat is quite helpful. Guess he must be either a reincarnation of Newton or Leibnitz, these gentlemen used to talk like that, or that I should take care at some point of my garden, remove catnip and other weeds.

\bigskip

This being said, wait. There is suggestion to connect integrals and derivatives, and this is in fact what we have, coming from Conjecture 4.16 and Theorem 4.17, due to:
$$\left(\frac{x^{p+1}}{p+1}\right)'=x^p\quad,\quad (e^x)'=e^x$$

So, eureka, we have our idea, thanks cat. Moving ahead now, following this idea, we first have the following result, called fundamental theorem of calculus:

\index{fundamental theorem of calculus}

\begin{theorem}
Given a continuous function $f:[a,b]\to\mathbb R$, if we set
$$F(x)=\int_a^xf(s)ds$$
then $F'=f$. That is, the derivative of the integral is the function itself. 
\end{theorem}

\begin{proof}
This follows from the Riemann integration picture, and more specifically, from the mean value property from Theorem 4.14. Indeed, we have:
$$\frac{F(x+t)-F(x)}{t}=\frac{1}{t}\int_x^{x+t}f(x)dx$$

On the other hand, our function $f$ being continuous, by using the mean value property from Theorem 4.14, we can find a number $c\in[x,x+t]$ such that:
$$\frac{1}{t}\int_x^{x+t}f(x)dx=f(c)$$

Thus, putting our formulae together, we conclude that we have:
$$\frac{F(x+t)-F(x)}{t}=f(c)$$

Now with $t\to0$, no matter how the number $c\in[x,x+t]$ varies, one thing that we can be sure about is that we have $c\to x$. Thus, by continuity of $f$, we obtain:
$$\lim_{t\to0}\frac{F(x+t)-F(x)}{t}=f(x)$$

But this means exactly that we have $F'=f$, and we are done.
\end{proof}

We have as well the following result, also called fundamental theorem of calculus:

\index{fundamental theorem of calculus}

\begin{theorem}
Given a differentiable function $F:\mathbb R\to\mathbb R$, we have
$$\int_a^bF'(x)dx=F(b)-F(a)$$
for any interval $[a,b]$.
\end{theorem}

\begin{proof}
This is something which follows from Theorem 4.19, and is in fact equivalent to it. Indeed, given $F:\mathbb R\to\mathbb R$ as above, consider the following function:
$$G(s)=\int_a^sF'(x)dx$$

By using Theorem 4.19 we have $G'=F'$, and so our functions $F,G$ differ by a constant. But with $s=a$ we have $G(a)=0$, and so the constant is $F(a)$, and we get:
$$F(s)=G(s)+F(a)$$

Now with $s=b$ this gives $F(b)=G(b)+F(a)$, which reads:
$$F(b)=\int_a^bF'(x)dx+F(a)$$

Thus, we are led to the conclusion in the statement.
\end{proof}

As a first illustration for all this, solving our previous problems, we have:

\begin{theorem}
We have the following integration formulae,
$$\int_a^bx^pdx=\frac{b^{p+1}-a^{p+1}}{p+1}\quad,\quad 
\int_a^b\frac{1}{x}\,dx=\log\left(\frac{b}{a}\right)$$
$$\int_a^b\sin x\,dx=\cos a-\cos b\quad,\quad
\int_a^b\cos x\,dx=\sin b-\sin a$$
$$\int_a^be^xdx=e^b-e^a\quad,\quad 
\int_a^b\log x\,dx=b\log b-a\log a-b+a$$
all obtained, in case you ever forget them, via the fundamental theorem of calculus.
\end{theorem}

\begin{proof}
We already know some of these formulae, but the best is to do everything, using the fundamental theorem of calculus. The computations go as follows:

\medskip

(1) With $F(x)=x^{p+1}$ we have $F'(x)=px^p$, and we get, as desired:
$$\int_a^bpx^p\,dx=b^{p+1}-a^{p+1}$$

(2) Observe first that the formula (1) does not work at $p=-1$. However, here we can use $F(x)=\log x$, having as derivative $F'(x)=1/x$, which gives, as desired:
$$\int_a^b\frac{1}{x}\,dx=\log b-\log a=\log\left(\frac{b}{a}\right)$$

(3) With $F(x)=\cos x$ we have $F'(x)=-\sin x$, and we get, as desired:
$$\int_a^b-\sin x\,dx=\cos b-\cos a$$

(4) With $F(x)=\sin x$ we have $F'(x)=\cos x$, and we get, as desired:
$$\int_a^b\cos x\,dx=\sin b-\sin a$$

(5) With $F(x)=e^x$ we have $F'(x)=e^x$, and we get, as desired:
$$\int_a^be^x\,dx=e^b-e^a$$

(6) This is something more tricky. We are looking for a function satisfying:
$$F'(x)=\log x$$

This does not look doable, but fortunately the answer to such things can be found on the internet. But, what if the internet connection is down? So, let us think a bit, and try to solve our problem. Speaking logarithm and derivatives, what we know is:
$$(\log x)'=\frac{1}{x}$$

But then, in order to make appear $\log$ on the right, the idea is quite clear, namely multiplying on the left by $x$. We obtain in this way the following formula:
$$(x\log x)'=1\cdot\log x+x\cdot\frac{1}{x}=\log x+1$$

We are almost there, all we have to do now is to substract $x$ from the left, as to get:
$$(x\log x-x)'=\log x$$

But this this formula in hand, we can go back to our problem, and we get the result.
\end{proof}

Getting back now to theory, inspired by the above, let us formulate:

\begin{definition}
Given a function $f$, we call primitive of $f$ any function $F$ satisfying:
$$F'=f$$
We denote such primitives by $\int f$, and also call them indefinite integrals.
\end{definition}

Observe that the primitives are unique up to an additive constant, in the sense that if $F$ is a primitive, then  so is $F+c$, for any $c\in\mathbb R$, and conversely, if $F,G$ are two primitives, then we must have $G=F+c$, for some $c\in\mathbb R$, with this latter fact coming from a result from chapter 3, saying that the derivative vanishes when the function is constant.

\bigskip

As for the convention at the end, $F=\int f$, this comes from the fundamental theorem of calculus, which can be written as follows, by using this convention:
$$\int_a^bf(x)dx=\left(\int f\right)(b)-\left(\int f\right)(a)$$

By the way, observe that there is no contradiction here, coming from the indeterminacy of $\int f$. Indeed, when adding a constant $c\in\mathbb R$ to the chosen primitive $\int f$, when conputing the above difference the $c$ quantities will cancel, and we will obtain the same result.

\bigskip

As an application, we can reformulate Theorem 4.21 in a more digest form, as follows:

\begin{theorem}
We have the following formulae for primitives,
$$\int x^p=\frac{x^{p+1}}{p+1}\quad,\quad\int\frac{1}{x}=\log x$$
$$\int\sin x=-\cos x\quad,\quad \int\cos x=\sin x$$
$$\int e^x=e^x\quad,\quad\int\log x=x\log x-x$$
allowing us to compute the corresponding definite integrals too.
\end{theorem}

\begin{proof}
Here the various formulae in the statement follow from Theorem 4.21, or rather from the proof of Theorem 4.21, or even from chapter 3, for most of them, and the last assertion comes from the integration formula given after Definition 4.22.
\end{proof}

Getting back now to theory, we have the following key result:

\index{integration by parts}

\begin{theorem}
We have the formula
$$\int f'g+\int fg'=fg$$
called integration by parts.
\end{theorem}

\begin{proof}
This follows by integrating the Leibnitz differentiation formula, namely:
$$(fg)'=f'g+fg'$$

Indeed, with our convention for primitives, this gives the formula in the statement.
\end{proof}

It is then possible to pass to usual integrals, and we obtain a formula here as well, as follows, also called integration by parts, with the convention $[\varphi]_a^b=\varphi(b)-\varphi(a)$:
$$\int_a^bf'g+\int_a^bfg'=\Big[fg\Big]_a^b$$

In practice, the most interesting case is that when $fg$ vanishes on the boundary $\{a,b\}$ of our interval, leading to the following formula:
$$\int_a^bf'g=-\int_a^bfg'$$

Examples of this usually come with $[a,b]=[-\infty,\infty]$, and more on this later. Now still at the theoretical level, we have as well the following result:

\index{change of variable}
\index{chain rule}

\begin{theorem}
We have the change of variable formula
$$\int_a^bf(x)dx=\int_c^df(g(t))g'(t)dt$$
where $c=g^{-1}(a)$ and $d=g^{-1}(b)$.
\end{theorem}

\begin{proof}
This follows with $f=F'$, from the following differentiation rule:
$$(Fg)'(t)=F'(g(t))g'(t)$$

Indeed, by integrating between $c$ and $d$, we obtain the result.
\end{proof}

As a theoretical application now of our integration theory, we have:

\index{Taylor formula}
\index{remainder}
\index{higher derivative}

\begin{theorem}
Given a $n$ times differentiable function $f:\mathbb R\to\mathbb R$, we have
$$f(x+t)=\sum_{k=0}^{n-1}\frac{f^{(k)}(x)}{k!}\,t^k+\int_x^{x+t}\frac{f^{(n)}(s)}{n!}(x+t-s)^{n-1}\,ds$$
called Taylor formula with integral formula for the remainder.
\end{theorem}

\begin{proof}
This is something quite standard, the idea being as follows:

\medskip

(1) At $n=1$ the formula in the statement is as follows, and certainly holds, due to the fundamental theorem of calculus, which gives $\int_x^{x+t}f'(s)ds=f(x+t)-f(x)$:
$$f(x+t)=f(x)+\int_x^{x+t}f'(s)ds$$

(2) At $n=2$, the formula in the statement becomes more complicated, as follows:
$$f(x+t)=f(x)+f'(x)t+\int_x^{x+t}f''(s)(x+t-s)ds$$

As a first observation, this formula holds indeed for the linear functions, where we have $f(x+t)=f(x)+f'(x)t$, and $f''=0$. So, let us try $f(x)=x^2$. Here we have:
$$f(x+t)-f(x)-f'(x)t=(x+t)^2-x^2-2xt=t^2$$

On the other hand, the integral remainder is given by the same formula, namely:
\begin{eqnarray*}
\int_x^{x+t}f''(s)(x+t-s)ds
&=&2\int_x^{x+t}(x+t-s)ds\\
&=&2t(x+t)-2\int_x^{x+t}sds\\
&=&2t(x+t)-((x+t)^2-x^2)\\
&=&2tx+2t^2-2tx-t^2\\
&=&t^2
\end{eqnarray*}

(3) Still at $n=2$, let us try now to prove the formula in the statement, in general. Since what we have to prove is an equality, this cannot be that hard, and the first thought goes towards differentiating. But this method works indeed, and we obtain the result.

\medskip

(4) In general, the proof is similar, by differentiating, the computations being similar to those at $n=2$, and we will leave this as an instructive exercise.
\end{proof}

So long for basic integration theory. As a first concrete application now, we can compute all sorts of areas and volumes. Normally such computations are the business of multivariable calculus, and we will be back to this later, but with the technology that we have so far, we can do a number of things. As a first such computation, we have:

\begin{proposition}
The area of an ellipse, given by the equation
$$\left(\frac{x}{a}\right)^2+\left(\frac{y}{b}\right)^2=1$$
with $a,b>0$ being half the size of a box containing the ellipse, is $A=\pi ab$.
\end{proposition}

\begin{proof}
The idea is that of cutting the ellipse into vertical slices. First observe that, according to our equation $(x/a)^2+(y/b)^2=1$, the $x$ coordinate can range as follows:
$$x\in[-a,a]$$

For any such $x$, the other coordinate $y$, satisfying $(x/a)^2+(y/b)^2=1$, is given by:
$$y=\pm b\sqrt{1-\frac{x^2}{a^2}}$$

Thus the length of the vertical ellipse slice at $x$ is given by the following formula:
$$l(x)=2b\sqrt{1-\frac{x^2}{a^2}}$$

We conclude from this discussion that the area of the ellipse is given by:
\begin{eqnarray*}
A
&=&2b\int_{-a}^a\sqrt{1-\frac{x^2}{a^2}}\,dx\\
&=&\frac{4b}{a}\int_0^a\sqrt{a^2-x^2}\,dx\\
&=&4ab\int_0^1\sqrt{1-y^2}\,dy\\
&=&4ab\cdot\frac{\pi}{4}\\
&=&\pi ab
\end{eqnarray*}

Finally, as a verification, for $a=b=1$ we get $A=\pi$, as we should.
\end{proof}

Moving now to 3D, as an obvious challenge here, we can try to compute the volume of the sphere. This can be done a bit as for the ellipse, the answer being as follows:

\begin{theorem}
The volume of the unit sphere is given by:
$$V=\frac{4\pi}{3}$$
More generally, the volume of the sphere of radius $R$ is $V=4\pi R^3/3$.
\end{theorem}

\begin{proof}
We proceed a bit as for the ellipse. The equation of the sphere is:
$$x^2+y^2+z^2=1$$

Now when the first coordinate $x\in[-1,1]$ is fixed, the other coordinates $y,z$ vary on a circle, given by the equation $y^2+z^2=1-x^2$, and so having radius as follows:
$$r(x)=\sqrt{1-x^2}$$

Thus, the vertical slice of our sphere at $x$ has area as follows:
$$a(x)=\pi r(x)^2=\pi(1-x^2)$$

We conclude from this discussion that the volume of the sphere is given by:
\begin{eqnarray*}
V
&=&\pi\int_{-1}^11-x^2\,dx\\
&=&\pi\int_{-1}^1\left(x-\frac{x^3}{3}\right)'dx\\
&=&\pi\left[\left(1-\frac{1}{3}\right)-\left(-1+\frac{1}{3}\right)\right]\\
&=&\pi\left(\frac{2}{3}+\frac{2}{3}\right)\\
&=&\frac{4\pi}{3}
\end{eqnarray*}

Finally, the last assertion is clear too, by multiplying everything by $R$, which amounts in multiplying the final result of our volume computation by $R^3$.
\end{proof}

As another application of our integration methods, we can now study the 1D wave equation. In order to explain this, we will need a standard calculus result, as follows:

\begin{proposition}
The derivative of a function of type
$$f(x)=\int_{a(x)}^{b(x)}g(s)ds$$
is given by the formula $f'(x)=g(b(x))b'(x)-g(a(x))a'(x)$.
\end{proposition}

\begin{proof}
Consider a primitive of the function that we integrate, $G'=g$. We have:
\begin{eqnarray*}
f(x)
&=&\int_{a(x)}^{b(x)}g(s)ds\\
&=&\int_{a(x)}^{b(x)}G'(s)ds\\
&=&G(b(x))-G(a(x))
\end{eqnarray*}

By using now the chain rule for derivatives, we obtain from this:
\begin{eqnarray*}
f'(x)
&=&G'(b(x))b'(x)-G'(a(x))a'(x)\\
&=&g(b(x))b'(x)-g(a(x))a'(x)
\end{eqnarray*}

Thus, we are led to the formula in the statement.
\end{proof}

Now back to the 1D waves, the result here, due to d'Alembert, is as follows:

\begin{theorem}
The solution of the 1D wave equation $\ddot{f}=v^2f''$ with initial value conditions $f(x,0)=g(x)$ and $\dot{f}(x,0)=h(x)$ is given by the d'Alembert formula:
$$f(x,t)=\frac{g(x-vt)+g(x+vt)}{2}+\frac{1}{2v}\int_{x-vt}^{x+vt}h(s)ds$$
Moreover, in the context of our previous lattice model discretizations, what happens is more or less that the above d'Alembert integral gets computed via Riemann sums.
\end{theorem}

\begin{proof}
There are several things going on here, the idea being as follows:

\medskip

(1) Let us first check that the d'Alembert solution is indeed a solution of the wave equation $\ddot{f}=v^2f''$. The first time derivative is computed as follows:
$$\dot{f}(x,t)=\frac{-vg'(x-vt)+vg'(x+vt)}{2}+\frac{1}{2v}(vh(x+vt)+vh(x-vt))$$

The second time derivative is computed as follows:
$$\ddot{f}(x,t)=\frac{v^2g''(x-vt)+v^2g''(x+vt)}{2}+\frac{vh'(x+vt)-vh'(x-vt)}{2}$$

Regarding now space derivatives, the first one is computed as follows:
$$f'(x,t)=\frac{g'(x-vt)+g'(x+vt)}{2}+\frac{1}{2v}(h(x+vt)-h(x-vt))$$

As for the second space derivative, this is computed as follows:
$$f''(x,t)=\frac{g''(x-vt)+g''(x+vt)}{2}+\frac{h'(x+vt)-h'(x-vt)}{2v}$$

Thus we have indeed $\ddot{f}=v^2f''$. As for the initial conditions, $f(x,0)=g(x)$ is clear from our definition of $f$, and $\dot{f}(x,0)=h(x)$ is clear from our above formula of $\dot{f}$.

\medskip

(2) Conversely, it is possible to prove that the solution is unique, but this is something quite tricky, requiring a better calculus knowledge, and we will be back to this later. As for the last assertion, we will leave the study here as an instructive exercise.
\end{proof}

Regarding the other equations from chapter 3, namely gravity and heat, things there are quite tricky too. We will be back to them later, when we will know more things.

\section*{4d. Some probability}

As another application of the integration theory developed above, let us develop now some theoretical probability theory. You probably know, from real life, what probability is. But in practice, when trying to axiomatize this, in mathematical terms, things can be quite tricky. So, here comes our point, the definition saving us is as follows:

\index{probability density}

\begin{definition}
A probability density is a function $\varphi:\mathbb R\to\mathbb R$ satisfying
$$\varphi\geq0\qquad,\qquad \int_\mathbb R\varphi(x)dx=1$$
with the convention that we allow Dirac masses, $\delta_x$ with $x\in\mathbb R$, as components of $\varphi$.
\end{definition}

To be more precise, in what regards the convention at the end, which is something of physics flavor, this states that our density function $\varphi:\mathbb R\to\mathbb R$ must be a combination as follows, with $\psi:\mathbb R\to\mathbb R$ being a usual function, and with $\alpha_i,x_i\in\mathbb R$:
$$\varphi=\psi+\sum_i\alpha_i\delta_{x_i}$$

Assuming that $x_i$ are distinct, and with the usual convention that the Dirac masses integrate up to 1, the conditions on our density function $\varphi:\mathbb R\to\mathbb R$ are as follows:
$$\psi\geq 0\quad,\quad \alpha_i\geq0\quad,\quad \int_\mathbb R\psi(x)dx+\sum_i\alpha_i=1$$

Observe the obvious relation with intuitive probability theory, where the probability for something to happen is always positive, $P\geq0$, and where the overall probability for something to happen, with this meaning for one of the possible events to happen, is of course $\Sigma P=1$, and this because life goes on, and something must happen, right.

\bigskip

In short, what we are proposing with Definition 4.31 is some sort of continuous generalization of basic probability theory, coming from coins, dice and cards, that you surely know. Moving now ahead, let us formulate, as a continuation of Definition 4.31:

\begin{definition}
We say that a random variable $f$ follows the density $\varphi$ if
$$P(f\in[a,b])=\int_a^b\varphi(x)dx$$
holds, for any interval $[a,b]\subset\mathbb R$.
\end{definition}

With this, we are now one step closer to what we know from coins, dice, cards and so on. For instance when rolling a die, the corresponding density is as follows:
$$\varphi=\frac{1}{6}\left(\delta_1+\delta_2+\delta_3+\delta_4+\delta_5+\delta_6\right)$$

\bigskip

In what regards now the random variables $f$, described as above by densities $\varphi$, the first questions regard their mean and variance, constructed as follows:

\index{mean}
\index{moment}
\index{variance}

\begin{definition}
Given a random variable $f$, with probability density $\varphi$:
\begin{enumerate}
\item Its mean is the quantity $M=\int_\mathbb Rx\varphi(x)\,dx$.

\item More generally, its $k$-th moment is $M_k=\int_\mathbb Rx^k\varphi(x)\,dx$.

\item Its variance is the quantity $V=M_2-M_1^2$.
\end{enumerate}
\end{definition}

Before going further, with more theory and examples, let us observe that, in both Definition 4.32 and Definition 4.33, what really matters is not the density $\varphi$ itself, but rather the related quantity $\mu=\varphi(x)dx$. So, let us upgrade our formalism, as follows:

\begin{definition}[upgrade]
A basic real probability measure is a quantity of the following type, with $\psi\geq 0$, $\alpha_i\geq0$ and $x_i\in\mathbb R$, satisfying $\int_\mathbb R\psi(x)dx+\sum_i\alpha_i=1$:
$$\mu=\psi(x)dx+\sum_i\alpha_i\delta_{x_i}$$
We say that a random variable $f$ follows $\mu$ when $P(f\in[a,b])=\int_a^bd\mu(x)$. In this case
$$M_k=\int_\mathbb Rx^k d\mu(x)$$
are called moments of $f$, and $M=M_1$ and $V=M_2-M_1^2$ are called mean, and variance.
\end{definition}

In practice now, let us look for some illustrations for this. The simplest random variables are those following discrete laws, $\psi=0$, and as a basic example here, when flipping a coin and being rewarded $\$1$ for heads, and $\$0$ for tails, the corresponding law is $\mu=\frac{1}{2}(\delta_0+\delta_1)$. More generally, playing the same game with a biased coin, which lands on heads with probability $p\in(0,1)$, leads to the following law, called Bernoulli law:
$$\mu_p=(1-p)\delta_0+p\delta_1$$

Many more things can be said here, notably with a study of what happens when you play the game $n$ times in a row, leading to some sort of powers of the Bernoulli laws, called binomial laws. Skipping some discussion here, and getting straight to the point, the most important laws in discrete probability are the Poisson laws, constructed as follows:

\begin{definition}
The Poisson law of parameter $1$ is the following measure,
$$p_1=\frac{1}{e}\sum_{n\geq0}\frac{\delta_n}{n!}$$
and more generally, the Poisson law of parameter $t>0$ is the following measure,
$$p_t=e^{-t}\sum_{n\geq0}\frac{t^n}{n!}\,\delta_n$$
with the letter ``p'' standing for Poisson.
\end{definition}

In general, the idea with the Poisson laws is that these appear a bit everywhere, in the real life, the reasons for this coming from the Poisson Limit Theorem (PLT). However, this theorem uses more advanced calculus, and we will leave it for later. In the meantime, however, we can have some fun with moments, the result here being as follows:

\begin{theorem}
The moments of $p_1$ are the Bell numbers,
$$M_k(p_1)=|P(k)|$$
where $P(k)$ is the set of partitions of $\{1,\ldots,k\}$. More generally, we have
$$M_k(p_t)=\sum_{\pi\in P(k)}t^{|\pi|}$$
for any $t>0$, where $|.|$ is the number of blocks.
\end{theorem}

\begin{proof}
The moments of $p_1$ satisfy the following recurrence formula:
\begin{eqnarray*}
M_{k+1}
&=&\frac{1}{e}\sum_{n\geq1}\frac{n^{k+1}}{n!}\\
&=&\frac{1}{e}\sum_{m\geq0}\frac{(m+1)^k}{m!}\\
&=&\frac{1}{e}\sum_{m\geq0}\frac{m^k}{m!}\left(1+\frac{1}{m}\right)^k\\
&=&\frac{1}{e}\sum_{m\geq0}\frac{m^k}{m!}\sum_{s=0}^k\binom{k}{s}m^{-s}\\
&=&\sum_{s=0}^k\binom{k}{s}\cdot\frac{1}{e}\sum_{m\geq0}\frac{m^{k-s}}{m!}\\
&=&\sum_{s=0}^k\binom{k}{s}M_{k-s}
\end{eqnarray*}

With this done, let us try now to find a recurrence for the Bell numbers, $B_k=|P(k)|$. Since a partition of $\{1,\ldots,k+1\}$ appears by choosing $s$ buddies for $1$, among the $k$ numbers available, and then partitioning the $k-s$ elements left, we have:
$$B_{k+1}=\sum_{s=0}^k\binom{k}{s}B_{k-s}$$

Since the initial values coincide, $M_1=B_1=1$ and $M_2=B_2=2$, we obtain by recurrence $M_k=B_k$, as claimed. Regarding now the law $p_t$ with $t>0$, we have here a similar recurrence formula for the moments, as follows:
\begin{eqnarray*}
M_{k+1}
&=&e^{-t}\sum_{m\geq0}\frac{t^{m+1}(m+1)^k}{m!}\\
&=&e^{-t}\sum_{m\geq0}\frac{t^{m+1}m^k}{m!}\sum_{s=0}^k\binom{k}{s}m^{-s}\\
&=&\sum_{s=0}^k\binom{k}{s}\cdot e^{-t}\sum_{m\geq0}\frac{t^{m+1}m^{k-s}}{m!}\\
&=&t\sum_{s=0}^k\binom{k}{s}M_{k-s}
\end{eqnarray*}

Regarding the initial values, the first moment of $p_t$ is given by:
$$M_1
=e^{-t}\sum_{n\geq0}\frac{t^nn}{n!}
=e^{-t}\sum_{m\geq0}\frac{t^{m+1}}{m!}
=t$$

Now by using the above recurrence we obtain from this:
$$M_2
=t\sum_{s=0}^1\binom{1}{s}M_{k-s}
=t(1+t)
=t+t^2$$

On the other hand, some standard combinatorics, a bit as before at $t=1$, shows that the numbers in the statement $S_k=\sum_{\pi\in P(k)}t^{|\pi|}$ satisfy the same recurrence relation, and with the same initial values. Thus we have $M_k=S_k$, as claimed.
\end{proof}

The above result looks quite exciting, making a link between subtle probability, and subtle combinatorics. We will be back to it, after learning more calculus.

\section*{4e. Exercises}

We had a lot of interesting theory in this chapter, and as exercises, we have:

\begin{exercise}
Clarify the integrability property of piecewise monotone functions.
\end{exercise}

\begin{exercise}
Find a geometric proof for $1^3+\ldots+N^3=(1+\ldots+N)^2$.
\end{exercise}

\begin{exercise}
Find estimates for the remainder, in the Taylor formula.
\end{exercise}

\begin{exercise}
Work out with full details the formula of $M_k(p_t)$, at $t>0$. 
\end{exercise}

As a bonus exercise, compute areas and volumes. As many as you can.

\part{Complex functions}

\ \vskip50mm

\begin{center}
{\em How can I change the world

If I can't even change myself

How can I change the way I am

I don't know, I don't know}
\end{center}

\chapter{Complex numbers}

\section*{5a. Complex numbers}

In this second part of the present book we discuss the functions of one complex variable $f:\mathbb C\to\mathbb C$. This is certainly something to be done before going to several variables, be them real or complex, because we have some unfinished business with $\sin,\cos,\exp,\log$, which remain quite mysterious objects. We will see here that these functions are much better understood as complex functions $f:\mathbb C\to\mathbb C$. In fact, even a dumb polynomial $P\in\mathbb R[X]$ is better understood as complex polynomial, $P\in\mathbb C[X]$, because its roots, that might not always exist as real numbers, always exist as complex numbers.

\bigskip

In short, many interesting things to be discussed, and this not necessarily for the sake of doing complex functions, but also with the goal of better understanding the real functions themselves. Let us begin with the complex numbers. There is a lot of magic here, and we will carefully explain this material. Their definition is as follows:

\index{complex number}
\index{i}

\begin{definition}
The complex numbers are variables of the form
$$x=a+ib$$
with $a,b\in\mathbb R$, which add in the obvious way, and multiply according to the following rule:
$$i^2=-1$$
Each real number can be regarded as a complex number, $a=a+i\cdot 0$.
\end{definition}

In other words, we consider variables as above, without bothering for the moment with their precise meaning. Now consider two such complex numbers:
$$x=a+ib\quad,\quad 
y=c+id$$

The formula for the sum is then the obvious one, as follows:
$$x+y=(a+c)+i(b+d)$$

As for the formula of the product, by using the rule $i^2=-1$, we obtain:
\begin{eqnarray*}
xy
&=&(a+ib)(c+id)\\
&=&ac+iad+ibc+i^2bd\\
&=&ac+iad+ibc-bd\\
&=&(ac-bd)+i(ad+bc)
\end{eqnarray*}

The advantage of using the complex numbers comes from the fact that the equation $x^2=1$ has now a solution, $x=i$. In fact, this equation has two solutions, namely:
$$x=\pm i$$

This is of course very good news. More generally, we have the following result, regarding the arbitrary degree 2 equations, with real coefficients:

\index{degree 2 equation}
\index{square root}

\begin{theorem}
The complex solutions of $ax^2+bx+c=0$ with $a,b,c\in\mathbb R$ are
$$x_{1,2}=\frac{-b\pm\sqrt{b^2-4ac}}{2a}$$
with the square root of negative real numbers being defined as
$$\sqrt{-m}=\pm i\sqrt{m}$$
and with the square root of positive real numbers being the usual one.
\end{theorem}

\begin{proof}
We can write indeed our equation in the following way:
\begin{eqnarray*}
ax^2+bx+c=0
&\iff&\left(x+\frac{b}{2a}\right)^2=\frac{b^2-4ac}{4a^2}\\
&\iff&x+\frac{b}{2a}=\pm\frac{\sqrt{b^2-4ac}}{2a}
\end{eqnarray*}

Thus, we are led to the conclusion in the statement.
\end{proof}

We will see later that any degree 2 complex equation has solutions as well, and that more generally, any polynomial equation, real or complex, has solutions. Moving ahead now, we can represent the complex numbers in the plane, in the following way:

\index{complex number}
\index{vector}
\index{sum of vectors}
\index{parallelogram rule}

\begin{proposition}
The complex numbers, written as usual
$$x=a+ib$$
can be represented in the plane, according to the following identification:
$$x=\binom{a}{b}$$
With this convention, the sum of complex numbers is the usual sum of vectors. 
\end{proposition}

\begin{proof}
Consider indeed two arbitrary complex numbers:
$$x=a+ib\quad,\quad 
y=c+id$$

Their sum is then by definition the following complex number:
$$x+y=(a+c)+i(b+d)$$

Now let us represent $x,y$ in the plane, as in the statement:
$$x=\binom{a}{b}\quad,\quad y=\binom{c}{d}$$

In this picture, their sum is given by the following formula:
$$x+y=\binom{a+c}{b+d}$$

But this is indeed the vector corresponding to $x+y$, so we are done.
\end{proof}

Here we have assumed that you are a bit familiar with vector calculus. If not, no problem, the idea is simply that vectors add by forming a parallelogram, as follows:
$$
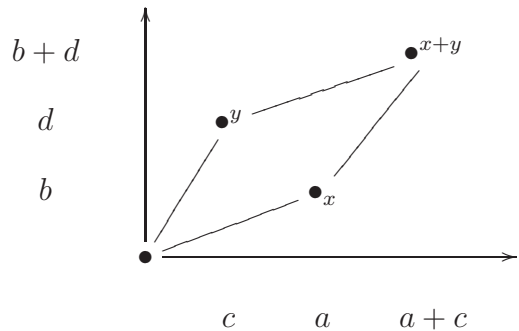
\xymatrix@R=8pt@C=15pt{
&&&\\
b+d&&&&\bullet^{x+y}\\
d&&\bullet^y\ar@{-}[urr]&&\\
b&&&\bullet_x\ar@{-}[uur]&\\
&\bullet\ar@{-}[urr]\ar[rrrr]\ar[uuuu]\ar@{-}[uur]&&&&\\
&&\ c\ &\ a\ &a+c}$$

Observe that in our geometric picture from Proposition 5.3, the real numbers correspond to the numbers on the $Ox$ axis. As for the purely imaginary numbers, these lie on the $Oy$ axis. As an illustration for this, we have the following basic picture:
$$\xymatrix@R=15pt@C=13pt{
&&\\
&&\bullet^i\ar@{-}[d]\ar@{.}@/^/[dr]\ar[u]\\
\ar@{-}[r]&\bullet^{-1}\ar@{-}[r]\ar@{.}@/^/[ur]&\ar@{-}[r]\ar@{-}[d]&\bullet^1\ar@{.}@/^/[dl]\ar[r]&\\
&&\bullet_{-i}\ar@{.}@/^/[ul]\ar@{-}[d]\\
&&&}$$

You might perhaps wonder why I chose to draw that circle, connecting the numbers $1,i,-1,-i$, which does not look very useful. More on this in a moment, the idea being that that circle can be indeed very useful, and coming in advance, some advice:

\begin{advice}
When drawing complex numbers, always begin with the coordinate axes $Ox,Oy$, and with a copy of the unit circle.
\end{advice}

We have so far a quite good understanding of their complex numbers, and their addition. In order to understand now the multiplication operation, we must do something more complicated, namely using polar coordinates. Let us start with:

\index{polar coordinates}
\index{modulus of complex number}
\index{argument of complex number}

\begin{definition}
The complex numbers $x=a+ib$ can be written in polar coordinates,
$$x=r(\cos t+i\sin t)$$
with the connecting formulae being as follows,
$$a=r\cos t\quad,\quad 
b=r\sin t$$
and in the other sense being as follows,
$$r=\sqrt{a^2+b^2}\quad,\quad 
\tan t=\frac{b}{a}$$
and with $r,t$ being called modulus, and argument.
\end{definition}

There is a clear relation here with the vector notation from Proposition 5.3, because $r$ is the length of the vector, and $t$ is the angle made by the vector with the $Ox$ axis. To be more precise, the picture for what is going on in Definition 5.5 is as follows:
$$
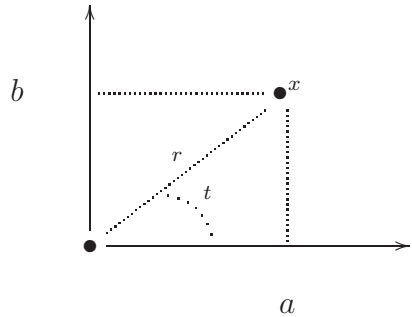
\xymatrix@R=10pt@C=15pt{
&&&\\
&&&\\
b&\ar@{.}[rrr]&&&\bullet^x\ar@{.}[ddd]\\
&&&&\\
&&&&\\
&\bullet\ar[rrrrr]\ar[uuuuu]\ar@{.}[uuurrr]^r&&\ar@{.}@/_/[ul]_t&&&\\
&&&&a}$$

As a basic example here, the number $i$ takes the following form:
$$i=\cos\left(\frac{\pi}{2}\right)+i\sin\left(\frac{\pi}{2}\right)$$

The point now is that in polar coordinates, the multiplication formula for the complex numbers, which was so far something quite opaque, takes a very simple form:

\begin{theorem}
Two complex numbers written in polar coordinates,
$$x=r(\cos s+i\sin s)\quad,\quad 
y=p(\cos t+i\sin t)$$
multiply according to the following formula:
$$xy=rp(\cos(s+t)+i\sin(s+t))$$
In other words, the moduli multiply, and the arguments sum up.
\end{theorem}

\begin{proof}
This follows from the following formulae, that we know well:
$$\cos(s+t)=\cos s\cos t-\sin s\sin t$$
$$\sin(s+t)=\cos s\sin t+\sin s\cos t$$

Indeed, we can assume that we have $r=p=1$, by dividing everything by these numbers. Now with this assumption made, we have the following computation:
\begin{eqnarray*}
xy
&=&(\cos s+i\sin s)(\cos t+i\sin t)\\
&=&(\cos s\cos t-\sin s\sin t)+i(\cos s\sin t+\sin s\cos t)\\
&=&\cos(s+t)+i\sin(s+t)
\end{eqnarray*}

Thus, we are led to the conclusion in the statement.
\end{proof}

The above result, which is based on some non-trivial trigonometry, is quite powerful. As a basic application of it, we can now compute powers, as follows:

\index{powers of complex number}

\begin{theorem}
The powers of a complex number, written in polar form,
$$x=r(\cos t+i\sin t)$$
are given by the following formula, valid for any exponent $k\in\mathbb N$:
$$x^k=r^k(\cos kt+i\sin kt)$$
Moreover, this formula holds in fact for any $k\in\mathbb Z$, and even for any $k\in\mathbb Q$.
\end{theorem}

\begin{proof}
Given a complex number $x$, written in polar form as above, and an exponent $k\in\mathbb N$, we have indeed the following computation, with $k$ terms everywhere:
\begin{eqnarray*}
x^k
&=&x\ldots x\\
&=&r(\cos t+i\sin t)\ldots r(\cos t+i\sin t)\\
&=&r^k([\cos(t+\ldots+t)+i\sin(t+\ldots+t))\\
&=&r^k(\cos kt+i\sin kt)
\end{eqnarray*}

Thus, we are done with the case $k\in\mathbb N$. Regarding now the generalization to the case $k\in\mathbb Z$, it is enough here to do the verification for $k=-1$, where the formula is:
$$x^{-1}=r^{-1}(\cos(-t)+i\sin(-t))$$

But this number $x^{-1}$ is indeed the inverse of $x$, as shown by:
\begin{eqnarray*}
xx^{-1}
&=&r(\cos t+i\sin t)\cdot r^{-1}(\cos(-t)+i\sin(-t))\\
&=&\cos(t-t)+i\sin(t-t)\\
&=&\cos 0+i\sin 0\\
&=&1
\end{eqnarray*}

Finally, regarding the generalization to the case $k\in\mathbb Q$, it is enough to do the verification for exponents of type $k=1/n$, with $n\in\mathbb N$. The claim here is that:
$$x^{1/n}=r^{1/n}\left[\cos\left(\frac{t}{n}\right)+i\sin\left(\frac{t}{n}\right)\right]$$

In order to prove this, let us compute the $n$-th power of this number. We can use the power formula for the exponent $n\in\mathbb N$, that we already established, and we obtain:
\begin{eqnarray*}
(x^{1/n})^n
&=&(r^{1/n})^n\left[\cos\left(n\cdot\frac{t}{n}\right)+i\sin\left(n\cdot\frac{t}{n}\right)\right]\\
&=&r(\cos t+i\sin t)\\
&=&x
\end{eqnarray*}

Thus, we have indeed a $n$-th root of $x$, and our proof is now complete.
\end{proof}

We should mention that there is a bit of ambiguity in the above, in the case of the exponents $k\in\mathbb Q$, due to the fact that the square roots, and the higher roots as well, can take multiple values, in the complex number setting. We will be back to this.

\bigskip

As a basic application of Theorem 5.7, we have the following result:

\index{square root}

\begin{proposition}
Each complex number, written in polar form,
$$x=r(\cos t+i\sin t)$$
has two square roots, given by the following formula:
$$\sqrt{x}=\pm\sqrt{r}\left[\cos\left(\frac{t}{2}\right)+i\sin\left(\frac{t}{2}\right)\right]$$
When $x>0$, these roots are $\pm\sqrt{x}$. When $x<0$, these roots are $\pm i\sqrt{-x}$.  
\end{proposition}

\begin{proof}
The first assertion is clear indeed from the general formula in Theorem 5.7, at $k=1/2$. As for its particular cases with $x\in\mathbb R$, these are clear from it.
\end{proof}

As a comment here, for $x>0$ we are very used to call the usual $\sqrt{x}$ square root of $x$. However, for $x<0$, or more generally for $x\in\mathbb C-\mathbb R_+$, there is less interest in choosing one of the possible $\sqrt{x}$ and calling it ``the'' square root of $x$, because all this is based on our convention that $i$ comes up, instead of down, which is something rather arbitrary. Actually, clocks turning clockwise, $i$ should be rather coming down. All this is a matter of taste, but in any case, for our math, the best is to keep some ambiguity, as above. 

\bigskip

With the above results in hand, and notably with the square root formula from Proposition 5.8, we can now go back to the degree 2 equations, and we have:

\index{degree 2 equation}
\index{square root}

\begin{theorem}
The complex solutions of $ax^2+bx+c=0$ with $a,b,c\in\mathbb C$ are
$$x_{1,2}=\frac{-b\pm\sqrt{b^2-4ac}}{2a}$$
with the square root of complex numbers being defined as above.
\end{theorem}

\begin{proof}
This is clear, the computations being the same as in the real case. To be more precise, our degree 2 equation can be written as follows:
$$\left(x+\frac{b}{2a}\right)^2=\frac{b^2-4ac}{4a^2}$$

Now since we know from Proposition 5.8 that any complex number has a square root, we are led to the conclusion in the statement.
\end{proof}

As a last general topic regarding the complex numbers, let us discuss conjugation. This is something quite tricky, complex number specific, as follows:

\index{complex conjugate}
\index{reflection}

\begin{definition}
The complex conjugate of $x=a+ib$ is the following number,
$$\bar{x}=a-ib$$
obtained by making a reflection with respect to the $Ox$ axis.
\end{definition}

As before with other such operations on complex numbers, a quick picture says it all. Here is the picture, with the numbers $x,\bar{x},-x,-\bar{x}$ being all represented:
$$
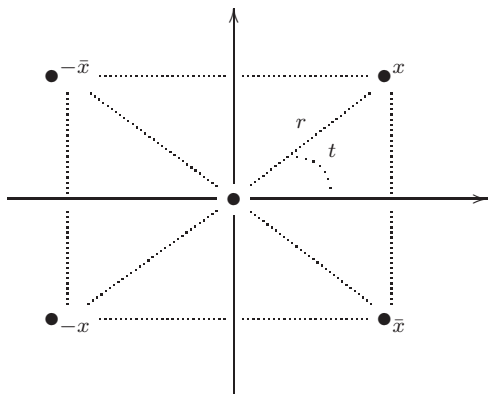
\xymatrix@R=6pt@C=10pt{
&&&&&&\\
&&&&&&\\
&\bullet^{-\bar{x}}\ar@{.}[ddd]\ar@{.}[rrr]&&&\ar@{.}[rrr]&&&\bullet^x\ar@{.}[ddd]\\
&&&&&&&\\
&&&&&&&\\
\ar@{-}[rrrr]&&&&\bullet\ar[rrrrr]\ar[uuuuu]\ar@{.}[dddlll]\ar@{.}[dddrrr]\ar@{.}[uuulll]\ar@{-}[ddddd]\ar@{.}[uuurrr]^r&&\ar@{.}@/_/[ul]_t&&&\\
&&&&&&&\\
&&&&&&&\\
&\bullet_{-x}\ar@{.}[rrr]\ar@{.}[uuu]&&&\ar@{.}[rrr]&&&\bullet_{\bar{x}}\ar@{.}[uuu]\\
&&&&&&\\
&&&&&&
}$$

Observe that the conjugate of a real number $x\in\mathbb R$ is the number itself, $x=\bar{x}$. In fact, the equation $x=\bar{x}$ characterizes the real numbers, among the complex numbers. At the level of non-trivial examples now, we have the following formula:
$$\bar{i}=-i$$

There are many things that can be said about the conjugation of the complex numbers, and here is a summary of basic such things that can be said:

\begin{theorem}
The conjugation operation $x\to\bar{x}$ has the following properties:
\begin{enumerate}
\item $x=\bar{x}$ precisely when $x$ is real.

\item $x=-\bar{x}$ precisely when $x$ is purely imaginary.

\item $x\bar{x}=|x|^2$, with $|x|=r$ being as usual the modulus.

\item With $x=r(\cos t+i\sin t)$, we have $\bar{x}=r(\cos t-i\sin t)$.

\item We have the formula $\overline{xy}=\bar{x}\bar{y}$, for any $x,y\in\mathbb C$.

\item The solutions of $ax^2+bx+c=0$ with $a,b,c\in\mathbb R$ are conjugate.
\end{enumerate}
\end{theorem}

\begin{proof}
These results are all elementary, the idea being as follows:

\medskip

(1) This is something that we already know, coming from definitions.

\medskip

(2) This is something clear too, because with $x=a+ib$ our equation $x=-\bar{x}$ reads $a+ib=-a+ib$, and so $a=0$, which amounts in saying that $x$ is purely imaginary.

\medskip

(3) This is a key formula, which can be proved as follows, with $x=a+ib$:
$$x\bar{x}
=(a+ib)(a-ib)
=a^2+b^2
=|x|^2$$

(4) This is clear indeed from the picture following Definition 5.10.

\medskip

(5) This is something quite magic, which can be proved as follows:
\begin{eqnarray*}
\overline{(a+ib)(c+id)}
&=&\overline{(ac-bd)+i(ad+bc)}\\
&=&(ac-bd)-i(ad+bc)\\
&=&(a-ib)(c-id)
\end{eqnarray*}

(6) This comes from the formula of the solutions, that we know from Theorem 5.2, but we can deduce this as well directly, without computations. Indeed, by using our assumption that the coefficients are real, $a,b,c\in\mathbb R$, we have:
\begin{eqnarray*}
ax^2+bx+c=0
&\implies&\overline{ax^2+bx+c}=0\\
&\implies&\bar{a}\bar{x}^2+\bar{b}\bar{x}+\bar{c}=0\\
&\implies&a\bar{x}^2+b\bar{x}+c=0
\end{eqnarray*}

Thus, we are led to the conclusion in the statement.
\end{proof}

\section*{5b. Exponential writing}

Let us discuss now the theory of complex functions $f:\mathbb C\to\mathbb C$, in analogy with the theory of the real functions $f:\mathbb R\to\mathbb R$. We will see that many results that we know from the real setting extend to the complex setting. Before starting, two remarks on this:

\bigskip

(1) Most of the real functions $f:\mathbb R\to\mathbb R$ that we know, such as $\sin,\cos,\exp,\log$, extend into complex functions $f:\mathbb C\to\mathbb C$, and the study of these latter extensions will bring some new light on the original real functions. Thus, what we will be doing here will be, in a certain sense, a refinement of the theory developed in Part I.

\bigskip

(2) On the other hand, since we have $\mathbb C\simeq\mathbb R^2$, the complex functions $f:\mathbb C\to\mathbb C$ that we will study here can be regarded as functions $f:\mathbb R^2\to\mathbb R^2$. This is something quite subtle, but in any case, what we will be doing here will stand as well as an introduction to the functions of type $f:\mathbb R^N\to\mathbb R^M$, that we will study in Parts III-IV.

\bigskip

In short, one complex variable is something in between one real variable, and two or more real variables, and we can only expect to end up with a mysterious mixture of surprising and unsurprising results. Welcome to complex analysis. Let us start with:

\index{complex function}
\index{continuous function}

\begin{definition}
A complex function $f:\mathbb C\to\mathbb C$, or more generally $f:X\to\mathbb C$, with $X\subset\mathbb C$ being a subset, is called continuous when, for any $x_n,x\in X$
$$x_n\to x\implies f(x_n)\to f(x)$$
where the convergence of the sequences of complex numbers, $x_n\to x$, means by definition that for $n$ big enough, the quantity $|x_n-x|$ becomes arbitrarily small.
\end{definition}

Observe that in real coordinates, $x=(a,b)$, the distances appearing in the definition of the convergence $x_n\to x$ are given by the following formula:
$$|x_n-x|=\sqrt{(a_n-a)^2+(b_n-b)^2}$$

Thus $x_n\to x$ in the complex sense means that $(a_n,b_n)\to(a,b)$ in the usual, intuitive sense, with respect to the usual distance in the plane $\mathbb R^2$, and as a consequence, a function $f:\mathbb C\to\mathbb C$ is continuous precisely when it is continuous, in an intuitive sense, when regarded as function $f:\mathbb R^2\to\mathbb R^2$. But more on this later, in Part III.

\bigskip

At the level of basic examples now, we first have the following result:

\index{exponential}

\begin{theorem}
We can exponentiate the complex numbers, according to the formula
$$e^x=\sum_{k=0}^\infty\frac{x^k}{k!}$$
and the function $x\to e^x$ is continuous, and satisfies $e^{x+y}=e^xe^y$.
\end{theorem}

\begin{proof}
We must first prove that the series converges. But this follows from:
$$|e^x|
=\left|\sum_{k=0}^\infty\frac{x^k}{k!}\right|
\leq\sum_{k=0}^\infty\left|\frac{x^k}{k!}\right|
=\sum_{k=0}^\infty\frac{|x|^k}{k!}
=e^{|x|}<\infty$$

Regarding the formula $e^{x+y}=e^xe^y$, this follows too as in the real case, as follows:
$$e^{x+y}
=\sum_{k=0}^\infty\frac{(x+y)^k}{k!}
=\sum_{k=0}^\infty\sum_{s=0}^k\frac{x^sy^{k-s}}{s!(k-s)!}
=e^xe^y$$

Finally, the continuity of $x\to e^x$ comes at $x=0$ from the following computation:
$$|e^t-1|
=\left|\sum_{k=1}^\infty\frac{t^k}{k!}\right|
\leq\sum_{k=1}^\infty\left|\frac{t^k}{k!}\right|
=\sum_{k=1}^\infty\frac{|t|^k}{k!}
=e^{|t|}-1$$

As for the continuity of $x\to e^x$ in general, this can be deduced now as follows:
$$\lim_{t\to0}e^{x+t}=\lim_{t\to0}e^xe^t=e^x\lim_{t\to0}e^t=e^x\cdot 1=e^x$$

Thus, we are led to the conclusions in the statement.
\end{proof}

As an application of this, let us discuss now the final and most convenient writing of the complex numbers, $x=re^{it}$. The point with this formula comes from:

\index{polar writing}

\begin{theorem}
We have the following formula, 
$$e^{it}=\cos t+i\sin t$$
due to Euler, valid for any $t\in\mathbb R$.
\end{theorem}

\begin{proof}
Many things can be said here, the idea being as follows:

\medskip

(1) To start with, the Euler formula shows that $t\to e^{it}$ maps $\mathbb R\to\mathbb T$, which might seem quite surprising, so let us try to understand this. We have, for any $x\in\mathbb C$:
$$e^{\bar{x}}=\sum_{k=0}^\infty\frac{\bar{x}^k}{k!}=\overline{\sum_{k=0}^\infty\frac{x^k}{k!}}=\overline{e^x}$$

Also, we have as well the following computation, again valid for any $x\in\mathbb C$:
$$e^xe^{-x}=e^{x-x}=e^0=1\implies (e^x)^{-1}=e^{-x}$$

(2) Now these two formulae, applied with $x=it$, with $t\in\mathbb R$, give:
$$e^{-it}=\overline{e^{it}}\quad,\quad (e^{it})^{-1}=e^{-it}$$

We conclude that the complex number $z=e^{it}$ has the following property:
$$z^{-1}=\bar{z}$$

But this is exactly the equation of the unit circle $\mathbb T\subset\mathbb C$, as desired.

\medskip

(3) Getting now to the proof of the Euler formula, we know that $t\to e^{it}$ maps $\mathbb R\to\mathbb T$, and also that this is a group morphism, meaning mapping sums in $\mathbb R$ into products in $\mathbb T$. But in view of this, barring any pathologies, this operation can only appear by ``wrapping". That is, we must have a formula as follows, for a certain $\alpha\in\mathbb R$:
$$e^{it}=\cos(\alpha t)+i\sin(\alpha t)$$

(4) In order now to find the parameter $\alpha\in\mathbb R$, let us look at what happens around $t=0$. And here, we have the following basic estimate, obtained by truncating $\exp$:
$$e^{it}\simeq 1+it$$

On the other hand, according to the basic trigonometry estimates for $\sin$ and $\cos$, from plane geometry, we have as well the following estimate, again around $t=0$:
$$\cos(\alpha t)+i\sin(\alpha t)\simeq 1+i\alpha t$$

We conclude that we must have $\alpha=1$, which gives the Euler formula. 

\medskip

(5) In practice now, the above proof is not exactly complete, because we still have to show that the pathologies evoked in (3) cannot appear, that is, that our continuous group morphism $\mathbb R\to\mathbb T$ must appear indeed via wrapping. And we will leave this, turning what we have into a full, rigorous proof of the Euler formula, as an instructive exercise.

\medskip

(6) As an alternative to all this, short-circuiting all the above discussion, we can kill the problem with calculus. Consider indeed the following function $f:\mathbb R\to\mathbb C$:
$$f(t)=\frac{\cos t+i\sin t}{e^{it}}$$

The point now is that we can compute the derivative of this function $f$ by using our first derivative formulae for $\exp$, $\sin$, $\cos$, and we obtain in this way:
\begin{eqnarray*}
f'(t)
&=&(e^{-it}(\cos t+i\sin t))'\\
&=&-ie^{-it}(\cos t+i\sin t)+e^{-it}(-\sin t+i\cos t)\\
&=&e^{-it}(-i\cos t+\sin t)+e^{-it}(-\sin t+i\cos t)\\
&=&0
\end{eqnarray*}

We conclude that our function $f:\mathbb R\to\mathbb C$ is constant, and the constant in question can be found by setting $t=0$, where we obtain, as desired:
$$f(0)=\frac{\cos 0+i\sin 0}{e^{i0}}=\frac{1}{1}=1$$

(7) Summarizing, Euler formula proved, one way or another. There are of course many other proofs, as for instance by using $(\cos t+i\sin t)'=i(\cos t+i\sin t)$, and then solving $g'=ig$, via standard calculus methods, as suggested in chapter 3.

\medskip

(8) Talking calculus, and what we saw in chapter 3, now that we proved the Euler formula $e^{it}=\cos t+i\sin t$, we can safely use the following consequences of it:
$$\cos t=\sum_{l=0}^\infty(-1)^l\frac{t^{2l}}{(2l)!}\quad,\quad 
\sin t=\sum_{l=0}^\infty(-1)^l\frac{t^{2l+1}}{(2l+1)!}$$

(9) And with the comment here that, conversely, these formulae can potentially lead to yet one more proof of the Euler formula, provided that we can find a good argument proving that $\sin$, $\cos$ are analytic, that is, equal to their own Taylor series. And exercise of course for you to explore a bit more all this, and the Euler formula in general.
\end{proof}

As an interesting application of the Euler formula, we have:

\begin{theorem}
We have the following formula,
$$e^{\pi i}=-1$$
and we have $E=mc^2$ as well.
\end{theorem}

\begin{proof}
We have two assertions here, the idea being as follows:

\medskip

(1) The first formula, $e^{\pi i}=-1$, which is actually the main formula in mathematics, comes from Theorem 5.14, by setting $t=\pi$. Indeed, we obtain:
\begin{eqnarray*}
e^{\pi i}
&=&\cos\pi+i\sin\pi\\
&=&-1+i\cdot 0\\
&=&-1
\end{eqnarray*}

(2) As for $E=mc^2$, which is the main formula in physics, this is something deep too. Although we will not really need it here, we recommend learning it as well, for symmetry reasons between math and physics, say from Feynman \cite{fe1}, \cite{fe2}, \cite{fe3}.
\end{proof}

Now back to our $x=re^{it}$ objectives, with the above theory in hand we can indeed use from now on this notation, the complete statement being as follows:

\index{polar coordinates}

\begin{theorem}
The complex numbers $x=a+ib$ can be written in polar coordinates,
$$x=re^{it}$$
with the connecting formulae being
$$a=r\cos t\quad,\quad 
b=r\sin t$$
and in the other sense being
$$r=\sqrt{a^2+b^2}\quad,\quad 
\tan t=\frac{b}{a}$$
and with $r,t$ being called modulus, and argument.
\end{theorem}

\begin{proof}
This is a reformulation of our previous Definition 5.5, by using the formula $e^{it}=\cos t+i\sin t$ from Theorem 5.14, and multiplying everything by $r$.
\end{proof}

With this in hand, we can now go back to the basics, namely the addition and multiplication of the complex numbers. We have the following result:

\index{multiplication of complex numbers}

\begin{theorem}
In polar coordinates, the complex numbers multiply as
$$re^{is}\cdot pe^{it}=rp\,e^{i(s+t)}$$
with the arguments $s,t$ being taken modulo $2\pi$.
\end{theorem}

\begin{proof}
This is something that we already know, from Theorem 5.6, reformulated by using the notations from Theorem 5.16. Observe that this follows as well directly, from the fact that we have $e^{x+y}=e^xe^y$, that we know from Theorem 5.13.
\end{proof}

The above formula is obviously very powerful. However, in polar coordinates we do not have a simple formula for the sum. Thus, this formalism has its limitations.

\bigskip

We can investigate as well more complicated operations, as follows:

\begin{theorem}
We have the following operations on the complex numbers, written in polar form, as above:
\begin{enumerate}
\item Inversion: $(re^{it})^{-1}=r^{-1}e^{-it}$.

\item Square roots: $\sqrt{re^{it}}=\pm\sqrt{r}e^{it/2}$.

\item Powers: $(re^{it})^a=r^ae^{ita}$.

\item Conjugation: $\overline{re^{it}}=re^{-it}$.
\end{enumerate}
\end{theorem}

\begin{proof}
This is something that we already know, from Theorem 5.7, but we can now discuss all this, from a more conceptual viewpoint, the idea being as follows:

\medskip

(1) We have indeed the following computation, using Theorem 5.17:
$$(re^{it})(r^{-1}e^{-it})
=rr^{-1}\cdot e^{i(t-t)}
=1$$

(2) Once again by using Theorem 5.17, we have:
$$(\pm\sqrt{r}e^{it/2})^2
=(\sqrt{r})^2e^{i(t/2+t/2)}
=re^{it}$$

(3) Given an arbitrary number $a\in\mathbb R$, we can define, as stated:
$$(re^{it})^a=r^ae^{ita}$$

Due to Theorem 5.17, this operation $x\to x^a$ is indeed the correct one.

\medskip

(4) This comes from the fact, that we know from Theorem 5.11, that the conjugation operation $x\to\bar{x}$ keeps the modulus, and switches the sign of the argument.
\end{proof}

\section*{5c. Equations, roots}

Getting back to algebra, recall from Theorem 5.9 that any degree 2 equation has 2 complex roots. We can in fact prove that any polynomial equation, of arbitrary degree $N\in\mathbb N$, has exactly $N$ complex solutions, counted with multiplicities:

\index{roots of polynomial}
\index{complex roots}

\begin{theorem}
Any polynomial $P\in\mathbb C[X]$ decomposes as
$$P=c(X-a_1)\ldots (X-a_N)$$
with $c\in\mathbb C$ and with $a_1,\ldots,a_N\in\mathbb C$.
\end{theorem}

\begin{proof}
The problem is that of proving that our polynomial has at least one root, because afterwards we can proceed by recurrence. We prove this by contradiction. So, assume that $P$ has no roots, and pick a number $z\in\mathbb C$ where $|P|$ attains its minimum:
$$|P(z)|=\min_{x\in\mathbb C}|P(x)|>0$$ 

Since $Q(t)=P(z+t)-P(z)$ is a polynomial which vanishes at $t=0$, this polynomial must be of the form $ct^k$ + higher terms, with $c\neq0$, and with $k\geq1$ being an integer. We obtain from this that, with $t\in\mathbb C$ small, we have the following estimate:
$$P(z+t)\simeq P(z)+ct^k$$

Now let us write $t=rw$, with $r>0$ small, and with $|w|=1$. Our estimate becomes:
$$P(z+rw)\simeq P(z)+cr^kw^k$$

Now recall that we assumed $P(z)\neq0$. We can therefore choose $w\in\mathbb T$ such that $cw^k$ points in the opposite direction to that of $P(z)$, and we obtain in this way:
\begin{eqnarray*}
|P(z+rw)|
&\simeq&|P(z)+cr^kw^k|\\
&=&|P(z)|(1-|c|r^k)
\end{eqnarray*}

Now by choosing $r>0$ small enough, as for the error in the first estimate to be small, and overcame by the negative quantity $-|c|r^k$, we obtain from this:
$$|P(z+rw)|<|P(z)|$$

But this contradicts our definition of $z\in\mathbb C$, as a point where $|P|$ attains its minimum. Thus $P$ has a root, and by recurrence it has $N$ roots, as stated.
\end{proof}

All this is very nice, and we will see applications in a moment. As a word of warning, however, we should mention that the above result remains something quite theoretical. Indeed, the proof is by contradiction, and there is no way of recycling the material there into something explicit, that can be used for effectively computing the roots.

\bigskip

Still talking polynomials and their roots, let us try however to understand what the analogue of $\Delta=b^2-4ac$ is, for an arbitrary polynomial $P\in\mathbb C[X]$. We will need:

\index{resultant}
\index{common roots}

\begin{theorem}
Given two polynomials $P,Q\in\mathbb C[X]$, written as follows,
$$P=c(X-a_1)\ldots(X-a_k)\quad,\quad 
Q=d(X-b_1)\ldots(X-b_l)$$
the following quantity, which is called resultant of $P,Q$,
$$R(P,Q)=c^ld^k\prod_{ij}(a_i-b_j)$$
is a polynomial in the coefficients of $P,Q$, with integer coefficients, and we have
$$R(P,Q)=0$$
precisely when $P,Q$ have a common root.
\end{theorem}

\begin{proof}
Given $P,Q\in\mathbb C[X]$, we can certainly construct the quantity $R(P,Q)$ in the statement, and we have then $R(P,Q)=0$ precisely when $P,Q$ have a common root. The whole point is that of proving that $R(P,Q)$ is a polynomial in the coefficients of $P,Q$, with integer coefficients. But this can be checked as follows:

\medskip

(1) We can expand the formula of $R(P,Q)$, and in what regards $a_1,\ldots,a_k$, which are the roots of $P$, we obtain in this way certain symmetric functions in these variables, which will be therefore polynomials in the coefficients of $P$, with integer coefficients.

\medskip

(2) We can then look what happens with respect to the remaining variables $b_1,\ldots,b_l$, which are the roots of $Q$. Once again what we have here are certain symmetric functions, and so polynomials in the coefficients of $Q$, with integer coefficients.

\medskip

(3) Thus, we are led to the conclusion in the statement, that $R(P,Q)$ is a polynomial in the coefficients of $P,Q$, with integer coefficients, and with the remark that the $c^ld^k$ factor is there for these latter coefficients to be indeed integers, instead of rationals.
\end{proof}

All this might seem a bit complicated, and as an illustration, let us work out an example. Consider the case of a polynomial of degree 2, and a polynomial of degree 1:
$$P=ax^2+bx+c\quad,\quad 
Q=dx+e$$

In order to compute the resultant, let us factorize our polynomials:
$$P=a(x-p)(x-q)\quad,\quad 
Q=d(x-r)$$

The resultant can be then computed as follows, by using the method above:
\begin{eqnarray*}
R(P,Q)
&=&ad^2(p-r)(q-r)\\
&=&ad^2(pq-(p+q)r+r^2)\\
&=&cd^2+bd^2r+ad^2r^2\\
&=&cd^2-bde+ae^2
\end{eqnarray*}

Finally, observe that $R(P,Q)=0$ corresponds indeed to the fact that $P,Q$ have a common root. Indeed, the root of $Q$ is $r=-e/d$, and we have:
$$P(r)
=\frac{ae^2}{d^2}-\frac{be}{d}+c
=\frac{R(P,Q)}{d^2}$$

Thus $P(r)=0$ precisely when $R(P,Q)=0$, as predicted by Theorem 5.20.

\bigskip

With this, we can now talk about the discriminant of any polynomial, as follows:

\index{discriminant}
\index{double root}
\index{single roots}

\begin{theorem}
Given a polynomial $P\in\mathbb C[X]$, written as
$$P(X)=cX^N+dX^{N-1}+\ldots$$
its discriminant, defined as being the following quantity,
$$\Delta(P)=\frac{(-1)^{\binom{N}{2}}}{c}R(P,P')$$
is a polynomial in the coefficients of $P$, with integer coefficients, and
$$\Delta(P)=0$$
happens precisely when $P$ has a double root.
\end{theorem}

\begin{proof}
This follows from Theorem 5.20, applied with $P=Q$, with the division by $c$ being indeed possible, under $\mathbb Z$, and with the sign being there for various reasons, including the compatibility with some well-known formulae, at small values of $N\in\mathbb N$.
\end{proof}

As an illustration, let us see what happens in degree 2. Here we have:
$$P=aX^2+bX+c\quad,\quad 
P'=2aX+b$$

Thus, the resultant is given by the following formula:
\begin{eqnarray*}
R(P,P')
&=&ab^2-b(2a)b+c(2a)^2\\
&=&4a^2c-ab^2\\
&=&-a(b^2-4ac)
\end{eqnarray*}

With the normalizations in Theorem 5.21 made, we obtain, as we should:
$$\Delta(P)=b^2-4ac$$

As another illustration, let us work out what happens in degree 3. Here the result, which is useful and interesting, and is probably new to you, is as follows:

\begin{theorem}
The discriminant of a degree $3$ polynomial,
$$P=aX^3+bX^2+cX+d$$
is given by $\Delta(P)=b^2c^2-4ac^3-4b^3d-27a^2d^2+18abcd$.
\end{theorem}

\begin{proof}
We need to do some tough computations here. Let us first compute resultants. Consider two polynomials, of degree 3 and degree 2, written as follows:
$$P=aX^3+bX^2+cX+d=a(X-p)(X-q)(X-r)$$
$$Q=eX^2+fX+g=e(X-s)(X-t)$$

The resultant of these two polynomials is then given by:
\begin{eqnarray*}
R(P,Q)
&=&a^2e^3(p-s)(p-t)(q-s)(q-t)(r-s)(r-t)\\
&=&a^2\cdot e(p-s)(p-t)\cdot e(q-s)(q-t)\cdot e(r-s)(r-t)\\
&=&a^2Q(p)Q(q)Q(r)\\
&=&a^2(ep^2+fp+g)(eq^2+fq+g)(er^2+fr+g)
\end{eqnarray*}

By expanding, we obtain the following formula for this resultant:
\begin{eqnarray*}
\frac{R(P,Q)}{a^2}
&=&e^3p^2q^2r^2+e^2f(p^2q^2r+p^2qr^2+pq^2r^2)\\
&+&e^2g(p^2q^2+p^2r^2+q^2r^2)+ef^2(p^2qr+pq^2r+pqr^2)\\
&+&efg(p^2q+pq^2+p^2r+pr^2+q^2r+qr^2)+f^3pqr\\
&+&eg^2(p^2+q^2+r^2)+f^2g(pq+pr+qr)\\
&+&fg^2(p+q+r)+g^3
\end{eqnarray*}

Note in passing that we have 27 terms on the right, as we should, and with this kind of check being mandatory, when doing such computations. Next, we have:
$$p+q+r=-\frac{b}{a}\quad,\quad
pq+pr+qr=\frac{c}{a}\quad,\quad 
pqr=-\frac{d}{a}$$

By using these formulae, we can produce some more, as follows:
$$p^2+q^2+r^2=(p+q+r)^2-2(pq+pr+qr)=\frac{b^2}{a^2}-\frac{2c}{a}$$
$$p^2q+pq^2+p^2r+pr^2+q^2r+qr^2=(p+q+r)(pq+pr+qr)-3pqr=-\frac{bc}{a^2}+\frac{3d}{a}$$
$$p^2q^2+p^2r^2+q^2r^2=(pq+pr+qr)^2-2pqr(p+q+r)=\frac{c^2}{a^2}-\frac{2bd}{a^2}$$

By plugging now this data into the formula of $R(P,Q)$, we obtain:
\begin{eqnarray*}
R(P,Q)
&=&a^2e^3\cdot\frac{d^2}{a^2}-a^2e^2f\cdot\frac{cd}{a^2}+a^2e^2g\left(\frac{c^2}{a^2}-\frac{2bd}{a^2}\right)+a^2ef^2\cdot\frac{bd}{a^2}\\
&+&a^2efg\left(-\frac{bc}{a^2}+\frac{3d}{a}\right)-a^2f^3\cdot\frac{d}{a}\\
&+&a^2eg^2\left(\frac{b^2}{a^2}-\frac{2c}{a}\right)+a^2f^2g\cdot\frac{c}{a}-a^2fg^2\cdot\frac{b}{a}+a^2g^3
\end{eqnarray*}

Thus, we have the following formula for the resultant:
\begin{eqnarray*}
R(P,Q)
&=&d^2e^3-cde^2f+c^2e^2g-2bde^2g+bdef^2-bcefg+3adefg\\
&-&adf^3+b^2eg^2-2aceg^2+acf^2g-abfg^2+a^2g^3
\end{eqnarray*}

Getting back now to our discriminant problem, with $Q=P'$, which corresponds to $e=3a$, $f=2b$, $g=c$, we obtain the following formula:
\begin{eqnarray*}
R(P,P')
&=&27a^3d^2-18a^2bcd+9a^2c^3-18a^2bcd+12ab^3d-6ab^2c^2+18a^2bcd\\
&-&8ab^3d+3ab^2c^2-6a^2c^3+4ab^2c^2-2ab^2c^2+a^2c^3
\end{eqnarray*}

By simplifying terms, and dividing by $a$, we obtain the following formula:
$$-\Delta(P)=27a^2d^2-18abcd+4ac^3+4b^3d-b^2c^2$$

But this gives the formula in the statement, and we are done.
\end{proof}

Still talking degree 3 equations, let us try to solve $P=0$, with $P=aX^3+bX^2+cX+d$ as above. By linear transformations we can assume $a=1,b=0$, and then it is convenient to write $c=3p,d=2q$. Thus, our equation becomes $x^3+3px+2q=0$, and regarding such equations, we have the following famous result, due to Cardano:

\begin{theorem}
For a normalized degree $3$ equation, namely 
$$x^3+3px+2q=0$$
the discriminant is $\Delta=-108(p^3+q^2)$. Assuming $p,q\in\mathbb R$ and $\Delta<0$, the number
$$x=\sqrt[3]{-q+\sqrt{p^3+q^2}}+\sqrt[3]{-q-\sqrt{p^3+q^2}}$$
is a real solution of our equation.
\end{theorem}

\begin{proof}
The formula of $\Delta$ is clear from definitions, and with $108=4\times 27$. Now with $x$ as in the statement, by using $(a+b)^3=a^3+b^3+3ab(a+b)$, we have:
\begin{eqnarray*}
x^3
&=&\left(\sqrt[3]{-q+\sqrt{p^3+q^2}}+\sqrt[3]{-q-\sqrt{p^3+q^2}}\right)^3\\
&=&-2q+3\sqrt[3]{-q+\sqrt{p^3+q^2}}\cdot\sqrt[3]{-q-\sqrt{p^3+q^2}}\cdot x\\
&=&-2q+3\sqrt[3]{q^2-p^3-q^2}\cdot x\\
&=&-2q-3px
\end{eqnarray*}

Thus, we are led to the conclusion in the statement.
\end{proof}

There are many more things that can be said about degree 3 equations, along these lines, and we will certainly have an exercise about this, at the end of this chapter.

\section*{5d. Roots of unity}

We kept the best for the end. As a last topic regarding the complex numbers, which is something really beautiful, we have the roots of unity. Let us start with:

\index{roots of unity}

\begin{theorem}
The equation $x^N=1$ has $N$ complex solutions, namely
$$\left\{w^k\Big|k=0,1,\ldots,N-1\right\}\quad,\quad w=e^{2\pi i/N}$$
which are called roots of unity of order $N$.
\end{theorem}

\begin{proof}
This follows from the general multiplication formula for complex numbers from Theorem 5.17. Indeed, with $x=re^{it}$ our equation reads:
$$r^Ne^{itN}=1$$

Thus $r=1$, and $t\in[0,2\pi)$ must be a multiple of $2\pi/N$, as stated.
\end{proof}

As an illustration here, the roots of unity of small order, along with some of their basic properties, which are very useful for computations, are as follows:

\medskip

\underline{$N=1$}. Here the unique root of unity is 1.

\medskip

\underline{$N=2$}. Here we have two roots of unity, namely 1 and $-1$.

\medskip

\underline{$N=3$}. Here we have 1, then $w=e^{2\pi i/3}$, and then $w^2=\bar{w}=e^{4\pi i/3}$.

\medskip

\underline{$N=4$}. Here the roots of unity, read as usual counterclockwise, are $1,i,-1,-i$. 

\medskip

\underline{$N=5$}. Here, with $w=e^{2\pi i/5}$, the roots of unity are $1,w,w^2,w^3,w^4$.

\medskip

\underline{$N=6$}. Here a useful alternative writing is $\{\pm1,\pm w,\pm w^2\}$, with $w=e^{2\pi i/3}$.

\medskip

\underline{$N=7$}. Here, with $w=e^{2\pi i/7}$, the roots of unity are $1,w,w^2,w^3,w^4,w^5,w^6$.

\medskip

\underline{$N=8$}. Here the roots of unity, read as usual counterclockwise, are the numbers $1,w,i,iw,-1,-w,-i,-iw$, with $w=e^{\pi i/4}$, which is also given by $w=(1+i)/\sqrt{2}$. 

\medskip

The roots of unity are very useful variables, and have many interesting properties. As a first application, we can now solve the ambiguity questions related to the extraction of $N$-th roots, from Theorem 5.7 and Theorem 5.18, the statement being as follows:

\index{roots of unity}

\begin{theorem}
Any nonzero $x=re^{it}$ has $N$ roots of order $N$, which appear as
$$y=r^{1/N}e^{it/N}$$
multiplied by the $N$ roots of unity of order $N$.
\end{theorem}

\begin{proof}
We must solve the equation $z^N=x$, over the complex numbers. Since the number $y$ in the statement clearly satisfies $y^N=x$, our equation reformulates as: 
$$z^N=x\iff z^N=y^N\iff\left(\frac{z}{y}\right)^N=1$$

Thus, we are led to the conclusion in the statement.
\end{proof}

Here is now a very useful formula, regarding the roots of unity:

\index{roots of unity}
\index{barycenter}

\begin{theorem}
The roots of unity, $\{w^k\}$ with $w=e^{2\pi i/N}$, have the property
$$\sum_{k=0}^{N-1}(w^k)^s=N\delta_{N|s}$$
for any exponent $s\in\mathbb N$, where on the right we have a Kronecker symbol.
\end{theorem}

\begin{proof}
The numbers in the statement, when written more conveniently as $(w^s)^k$ with $k=0,\ldots,N-1$, form a certain regular polygon in the plane $P_s$. Thus, if we denote by $C_s$ the barycenter of this polygon, we have the following formula:
$$\frac{1}{N}\sum_{k=0}^{N-1}w^{ks}=C_s$$

Now observe that in the case $N\slash\hskip-1.6mm|\,s$ our polygon $P_s$ is non-degenerate, circling around the unit circle, and having center $C_s=0$. As for the case $N|s$, here the polygon is degenerate, lying at 1, and having center $C_s=1$. We conclude that we have:
$$C_s=\delta_{N|s}$$

Thus, we have indeed the formula in the statement.
\end{proof}

Time for some applications? Following Napoleon, and no joke here, we first have:

\begin{theorem}
In the context of the Napoleon configuration, namely
$$\xymatrix@R=9pt@C=0pt{
&&&&&&&&&&&&&&&&&&E\\
F&&&&&&&&A\ar@{-}[dddllll]\ar@{-}[dddrrrrrr]\ar@{--}[llllllll]\ar@{--}[urrrrrrrrrr]\\
&&&&R&&&&&&&&&Q\\
\\
&&&&B\ar@{-}[rrrrrrrrrr]\ar@{--}[uuullll]\ar@{--}[dddrrrrr]&&&&&&&&&&C\ar@{--}[dddlllll]\ar@{--}[uuuurrrr]\\
&&&&&&&&&P\\
\\
&&&&&&&&&D}$$
with equilateral triangles, and their barycenters drawn, the triangle $PQR$ is equilateral.
\end{theorem}

\begin{proof}
This can be proved by using the root of unity $w=e^{2\pi i/3}$, as follows:

\medskip

(1) Our first claim is that a triangle $IJK$, with $I,J,K$ appearing counterclockwise, is equilateral precisely when its vertices, regarded as complex numbers, satisfy:
$$I+wJ+w^2K=0$$

Indeed, the clockwise rotation by $60^\circ$ is $P\to -wP$, and by using this, along with $1+w+w^2$, coming from $w^3=1$, the condition for $IJK$ to be equilateral reads:
\begin{eqnarray*}
I-K=-w(J-K)
&\iff&I+wJ-(1+w)K=0\\
&\iff&I+wJ+w^2K=0
\end{eqnarray*}

(2) But with this, the Napoleon theorem can be proved, majestically, as follows:
\begin{eqnarray*}
P+wQ+w^2R
&=&\frac{B+C+D}{3}+w\cdot\frac{A+C+E}{3}+w^2\cdot\frac{A+B+F}{3}\\
&=&\frac{B+wA+w^2F}{3}+\frac{C+wE+w^2A}{3}+\frac{D+wC+w^2B}{3}\\
&=&0+0+0\\
&=&0
\end{eqnarray*}

(3) And exercise for you to learn more about the above configuration, and its relevance to the original triangle $ABC$, following Napoleon, Torricelli, Fermat and others.
\end{proof}

As a second application, this time in relation with degree 3 equations, we have:

\begin{theorem}
For a normalized degree $3$ equation, namely 
$$x^3+3px+2q=0$$
the discriminant is $\Delta=-108(p^3+q^2)$. Assuming $p,q\in\mathbb R$ and $\Delta<0$, the numbers
$$z=w\sqrt[3]{-q+\sqrt{p^3+q^2}}+w^2\sqrt[3]{-q-\sqrt{p^3+q^2}}$$
with $w=1,e^{2\pi i/3},e^{4\pi i/3}$ are the solutions of our equation.
\end{theorem}

\begin{proof}
As before in Theorem 5.23, using $(a+b)^3=a^3+b^3+3ab(a+b)$, we get:
\begin{eqnarray*}
z^3
&=&\left(w\sqrt[3]{-q+\sqrt{p^3+q^2}}+w^2\sqrt[3]{-q-\sqrt{p^3+q^2}}\right)^3\\
&=&-2q+3\sqrt[3]{-q+\sqrt{p^3+q^2}}\cdot\sqrt[3]{-q-\sqrt{p^3+q^2}}\cdot z\\
&=&-2q+3\sqrt[3]{q^2-p^3-q^2}\cdot z\\
&=&-2q-3pz
\end{eqnarray*}

Thus, we are led to the conclusion in the statement.
\end{proof}

As a comment here, the above formula holds in the case $\Delta>0$ too, and also when the coefficients are complex numbers, $p,q\in\mathbb C$. However, these extensions are quite often not very useful, because when it comes to extract the above square and cubic roots, for complex numbers, you can end up with the initial question, the one you started with.

\bigskip

Next, as a wrong application of the roots of unity, which is however sweet, we have:

\begin{theorem}
For a three times differentiable function we have, with $w=e^{2\pi i/3}$, 
$$\frac{f(x+t)+f(x+wt)+f(x+w^2t)}{3}\simeq f(x)+\frac{f'''(x)}{6}\,t^3$$
suggesting that $f'''(x)$ measures how far is $f(z)$ with $z\simeq x$ from $f(x)$, over $\mathbb C$.
\end{theorem}

\begin{proof}
We have indeed, with $a=f'(x)$, $b=f''(x)/2$ and $c=f'''(x)/6$:
\begin{eqnarray*}
f(x+t)+f(x+wt)+f(x+w^2t)
&\simeq&f(x)+at+bt^2+ct^3\\
&+&f(x)+wat+w^2bt^2+ct^3\\
&+&f(x)+w^2at+wbt^2+ct^3\\
&=&3f(x)+3ct^3
\end{eqnarray*}

Thus, formula proved, but the last conclusion is wrong, because the same trick will work with roots of unity of arbitrary order $N$, ruining our conclusion. Indeed, in case you don't believe me, here is what we get at $N=4$, with $d=f''''(x)/24$:
\begin{eqnarray*}
f(x+t)+f(x+it)+f(x-t)+f(x-it)
&\simeq&f(x)+at+bt^2+ct^3+dt^4\\
&+&f(x)+iat-bt^2-ict^3+dt^4\\
&+&f(x)-at+bt^2-ct^3+dt^4\\
&+&f(x)-iat-bt^2+ict^3+dt^4\\
&=&4f(x)+4dt^4
\end{eqnarray*}

In other words, assuming that $f$ is differentiable four times, we have:
$$\frac{f(x+t)+f(x+it)+f(x-t)+f(x-it)}{4}\simeq f(x)+\frac{f''''(x)}{24}\,t^4$$

Thus $f''''(x)$ measures too how far is $f(z)$ with $z\simeq x$ from $f(x)$, over $\mathbb C$. Damn.
\end{proof}

Finally, in case you are a bit familiar with linear algebra, you will certainly appreciate:

\begin{theorem}
The all-one matrix diagonalizes as follows,
$$\begin{pmatrix}
1&\ldots&1\\
\vdots&&\vdots\\
1&\ldots&1\end{pmatrix}
=\frac{1}{N}\,F_N\begin{pmatrix}
N\\
&0\\
&&\ddots\\
&&&0
\end{pmatrix}F_N^*$$
with $F_N=(w^{ij})_{ij}$ with $w=e^{2\pi i/N}$, with $i,j=0,1,\ldots,N-1$, being the Fourier matrix.
\end{theorem}

\begin{proof}
The all-one matrix being $N$ times the projection on the all-one vector, the diagonal form is the one in the statement. In order to find now the explicit diagonalization formula, with passage matrix and its inverse, we must solve the following equation:
$$x_0+\ldots+x_{N-1}=0$$

And this is not easy, if we want a nice basis for the space of solutions. Fortunately, the complex numbers come to the rescue, via Theorem 5.26, and solve the problem.
\end{proof}

\section*{5e. Exercises}

This was all beautiful material, and our exercises will be beautiful as well:

\begin{exercise}
Learn more about triangles, from Napoleon, Torricelli, Fermat.
\end{exercise}

\begin{exercise}
Learn more about degree $3$ equations, and about degree $4$ too.
\end{exercise}

\begin{exercise}
Further meditate on the possible interpretations of $f'''$ and $f''''$.
\end{exercise}

\begin{exercise}
Prove that the Fourier matrix $F_N$ diagonalizes any circulant matrix.
\end{exercise}

As bonus exercise, learn some physics, of all types, involving waves, particles and more, as much as you can. This is where the complex numbers really shine.

\chapter{Complex functions}

\section*{6a. Complex functions}

With the complex numbers understood, we would like to discuss now the theory of complex functions $f:\mathbb C\to\mathbb C$, in analogy with the theory of the real functions $f:\mathbb R\to\mathbb R$. We will see that many results that we know from the real setting extend to the complex setting, but there will be quite a number of new phenomena too. Let us start with:

\index{distance}
\index{complex plane}
\index{convergent sequence}

\begin{definition}
The distance between two complex numbers is the usual distance in the plane between them, namely: 
$$d(x,y)=|x-y|$$
With this, we can talk about convergence, by saying that $x_n\to x$ when $d(x_n,x)\to 0$.
\end{definition}

Observe that in real coordinates, the distance formula is quite complicated:
\begin{eqnarray*}
d(a+ib,c+id)
&=&|(a+ib)-(c+id)|\\
&=&|(a-c)+i(b-d)|\\
&=&\sqrt{(a-c)^2+(b-d)^2}
\end{eqnarray*}

However, for most computations, we will not need this formula, and we can get away with the various tricks regarding complex numbers that we know. As a first result now, regarding $\mathbb C$ and its distance, that we will need in what follows, we have:

\index{complete space}
\index{Cauchy sequence}

\begin{proposition}
The complex plane $\mathbb C$ is complete, in the sense that any Cauchy sequence converges.
\end{proposition}

\begin{proof}
Consider indeed a Cauchy sequence $\{x_n\}_{n\in\mathbb N}\subset\mathbb C$. If we write $x_n=a_n+ib_n$ for any $n\in\mathbb N$, then we have the following estimates:
$$|a_n-a_m|\leq\sqrt{(a_n-a_m)^2+(b_n-b_m)^2}=|x_n-x_m|$$
$$|b_n-b_m|\leq\sqrt{(a_n-a_m)^2+(b_n-b_m)^2}=|x_n-x_m|$$

Thus both the sequences $\{a_n\}_{n\in\mathbb N}\subset\mathbb R$ and $\{b_n\}_{n\in\mathbb N}\subset\mathbb R$ are Cauchy, and with $a_n\to a$ and $b_n\to b$ we can set $x=a+ib$, and we have our limit $x_n\to x$, as desired.
\end{proof}

Talking complex functions now, we have here the following definition:

\index{pointwise convergence}
\index{uniform convergence}

\begin{definition}
A complex function $f:\mathbb C\to\mathbb C$, or more generally $f:X\to\mathbb C$, with $X\subset\mathbb C$ being a subset, is called continuous when, for any $x_n,x\in X$:
$$x_n\to x\implies f(x_n)\to f(x)$$
Also, we can talk about pointwise convergence of functions, $f_n\to f$, and about uniform convergence too, $f_n\to_uf$, exactly as for the real functions.
\end{definition}

Observe that, since $x_n\to x$ in the complex sense means that $(a_n,b_n)\to(a,b)$ in the usual, real plane sense, a function $f:\mathbb C\to\mathbb C$ is continuous precisely when it is continuous when regarded as real function, $f:\mathbb R^2\to\mathbb R^2$. But more on this later in this book. At the level of examples now, we first have the polynomials, $P\in\mathbb C[X]$. We already met such polynomials in chapter 5, so let us recall from there that we have:

\index{polynomial}
\index{roots of polynomial}

\begin{theorem}
Each polynomial $P\in\mathbb C[X]$ can be regarded as a continuous function $P:\mathbb C\to\mathbb C$. Moreover, we have the formula
$$P(x)=a(x-r_1)\ldots(x-r_n)$$
with $a\in\mathbb C$, and with the numbers $r_1,\ldots,r_n\in\mathbb C$ being the roots of $P$.
\end{theorem}

\begin{proof}
This is something that we know from chapter 5, the idea being that one root can be always constructed, by reasoning by contradiction, and doing some analysis around the minimum of $|P|$, and then a recurrence on the degree $n\in\mathbb N$ does the rest.
\end{proof}

Next in line, we have the rational functions, which are defined as follows:

\index{rational function}
\index{pole}
\index{field}

\begin{theorem}
The quotients of complex polynomials $f=P/Q$ are called rational functions. When written in reduced form, with $P,Q$ prime to each other,
$$f=\frac{P}{Q}$$
is well-defined and continuous outside the zeroes $P_f\subset\mathbb C$ of $Q$, called poles of $f$:
$$f:\mathbb C-P_f\to\mathbb C$$
In addition, the rational functions, regarded as algebraic expressions, are stable under summing, making products and taking inverses.
\end{theorem}

\begin{proof}
There are several things going on here, the idea being as follows:

\medskip

(1) First of all, we can surely talk about quotients of polynomials, $f=P/Q$, regarded as abstract algebraic expressions. Also, the last assertion is clear, because we can indeed perform sums, products, and take inverses, by using the following formulae:
$$\frac{P}{Q}+\frac{R}{S}=\frac{PS+QR}{QS}\quad,\quad 
\frac{P}{Q}\cdot\frac{R}{S}=\frac{PR}{QS}\quad,\quad
\left(\frac{P}{Q}\right)^{-1}=\frac{Q}{P}$$

(2) The question is now, given a rational function $f$, can we regard it as a complex function? In general, we cannot say that we have $f:\mathbb C\to\mathbb C$, for instance because $f(x)=x^{-1}$ is not defined at $x=0$. More generally, assuming $f=P/Q$ with $P,Q\in\mathbb C[X]$, we cannot talk about $f(x)$ when $x$ is a root of $Q$, unless of course we are in the special situation where $x$ is a root of $P$ too, and we can simplify the fraction.

\medskip

(3) In view of this discussion, in order to solve our question, we must avoid the situation where the polynomials $P,Q$ have common roots. But this can be done by writing our rational function $f$ in reduced form, as follows, with $P,Q\in\mathbb C[X]$ prime to each other:
$$f=\frac{P}{Q}$$

(4) Now with this convention made, it is clear that $f$ is well-defined, and continuous too, outside of the zeroes of $f$. Now since these zeroes can be obviously recovered from the knowledge of $f$ itself, as being the points where ``$f$ explodes'', we can call them poles of $f$, and so we have a function $f:\mathbb C-P_f\to\mathbb C$, as in the statement.
\end{proof}

Obviously, what we have in Theorem 6.5 is a quite subtle mixture of algebra and analysis. So, let us first clarify the algebra part. In relation with this, we have:

\begin{definition}
A field is a set $F$ with a sum operation $+$ and a product operation $\times$, subject to the following conditions:
\begin{enumerate}
\item $a+b=b+a$, $a+(b+c)=(a+b)+c$, there exists $0\in F$ such that $a+0=0$, and any $a\in F$ has an inverse $-a\in F$, satisfying $a+(-a)=0$.

\item $ab=ba$, $a(bc)=(ab)c$, there exists $1\in F$ such that $a1=a$, and any $a\neq0$ has a multiplicative inverse $a^{-1}\in F$, satisfying $aa^{-1}=1$.

\item The sum and product are compatible via $a(b+c)=ab+ac$.
\end{enumerate}
\end{definition}

As basic examples of fields, we have the rational numbers $\mathbb Q$, the real numbers $\mathbb R$, and the complex numbers $\mathbb C$. Further examples come from arithmetic. In view of this, it is useful to think of any field $F$ as being a ``field of numbers'', and this because the elements $a,b,c,\ldots\in F$ behave under the operations $+$ and $\times$ exactly as the usual numbers do.

\bigskip

In what regards the various spaces of functions, such as the polynomials $\mathbb C[X]$, or the continuous functions $C(\mathbb R)$, these certainly have sum and product operations $+$ and $\times$, but are in general not fields, because they do not satisfy the following field axiom:
$$f\neq 0\implies\exists f^{-1}$$

However, and here comes our point, Theorem 6.5 tells us that the rational functions form a field. Which is quite interesting, so let us record this finding, as follows:

\begin{theorem}
The rational functions, $f=P/Q$ with $P,Q\in\mathbb C[X]$, are as follows:
\begin{enumerate}
\item They form a field, denoted $\mathbb C(X)$.

\item In addition, they are stable by composition.
\end{enumerate}
\end{theorem}

\begin{proof}
This is something self-explanatory, with (1) coming from Theorem 6.5 and its proof, and with (2) being something which is clear too.
\end{proof}

Moving on, with some analysis this time, we have the following result, which is something very useful, when dealing with computations with rational functions:

\begin{theorem}
The complex rational functions can be written as follows,
$$f(x)=\sum_i\frac{A_i(x)}{(r_i-x)^{n_i}}$$
with $A_i\in\mathbb C[X]$, and $r_i\in\mathbb C$ being the poles. Also, we have the following formula, 
$$\frac{1}{(r-x)^n}=\frac{1}{r^n}\sum_{k=0}^\infty\binom{n+k-1}{n-1}\left(\frac{x}{r}\right)^k$$
valid for $|x|<r$, which computes these rational functions, in practice. 
\end{theorem}

\begin{proof}
Consider indeed a rational function $f=P/Q$, with $P,Q\in\mathbb C[X]$ chosen prime to each other. By factorizing $Q$ we have a formula as follows, with $P(r_i)\neq0$:
$$f(x)=\frac{P(x)}{(r_1-x)^{n_1}\ldots(r_k-x)^{n_k}}$$

Now comes the trick. Assuming that $S,T\in\mathbb C[X]$ are prime to each other, we can find, a bit as for the usual numbers, by performing successive divisions, polynomials $A,B\in\mathbb C[X]$ such that $AT+BS=1$, and so that the following happens:
$$\frac{1}{ST}=\frac{A}{S}+\frac{B}{T}$$

Thus, we are led to the first formula in the statement. As for the second formula, this is something that we know since chapter 2, binomial formula with negative exponent.
\end{proof}
 
Next, let us point out that, contrary to what the above might suggest, everything does not always extend easily from the real to the complex case. For instance, we have:

\index{geometric series}
\index{spiral}

\begin{proposition}
We have the following formula, valid for any $|x|<1$,
$$\frac{1}{1-x}=1+x+x^2+\ldots$$
but, for $x\in\mathbb C-\mathbb R$, the geometric meaning of this formula is quite unclear.
\end{proposition}

\begin{proof}
Here the formula in the statement holds indeed, by multiplying and cancelling terms, and with the convergence being justified by the following estimate:
$$\left|\sum_{n=0}^\infty x^n\right|\leq\sum_{n=0}^\infty|x|^n=\frac{1}{1-|x|}$$

As for the last assertion, this is something quite informal. To be more precise, for $x=1/2$ our formula is clear, by cutting the interval $[0,2]$ into half, and so on:
$$1+\frac{1}{2}+\frac{1}{4}+\frac{1}{8}+\ldots=2$$

More generally, for $x\in(-1,1)$ the meaning of the formula in the statement is something quite clear and intuitive, geometrically speaking, by using a similar argument. Next, let us consider the case $x=rw$, with $r\in[0,1)$ and $w^N=1$. We have:
\begin{eqnarray*}
1+rw+r^2w^2+\ldots
&=&(1+rw+\ldots+r^{N-1}w^{N-1})\\
&+&(r^N+r^{N+1}w\ldots+r^{2N-1}w^{N-1})\\
&+&(r^{2N}+r^{2N+1}w\ldots+r^{3N-1}w^{N-1})\\
&+&\ldots
\end{eqnarray*}

Thus, by grouping the terms with the same argument, our infinite sum is:
\begin{eqnarray*}
1+rw+r^2w^2+\ldots
&=&(1+r^N+r^{2N}+\ldots)\\
&+&(r+r^{N+1}+r^{2N+1}+\ldots)w\\
&+&\ldots\\
&+&(r^{N-1}+r^{2N-1}+r^{3N-1}+\ldots)w^{N-1}
\end{eqnarray*}

But the sums of each ray can be computed with the real formula for geometric series, that we know and understand well, and with an extra bit of algebra, we get:
\begin{eqnarray*}
1+rw+r^2w^2+\ldots
&=&\frac{1}{1-r^N}+\frac{rw}{1-r^N}+\ldots+\frac{r^{N-1}w^{N-1}}{1-r^N}\\
&=&\frac{1}{1-r^N}\cdot\frac{1-(rw)^N}{1-rw}\\
&=&\frac{1}{1-rw}
\end{eqnarray*}

Summarizing, the geometric series formula can still be understood, in a purely geometric way, for variables of type $x=rw$, with $r\in[0,1)$, and with $w$ being a root of unity. In general, however, this formula tells us that the numbers on a certain infinite spiral sum up to a certain number, which remains something quite mysterious.
\end{proof}

Getting now to more complicated functions, such as $\sin$, $\cos$, $\exp$, $\log$, here things extend well from real to complex, with the basic theory being as follows:

\index{sin}
\index{cos}
\index{exp}
\index{log}

\begin{theorem}
The functions $\sin,\cos,\exp,\log$ have complex extensions, given by
$$\sin x=\sum_{l=0}^\infty(-1)^l\frac{x^{2l+1}}{(2l+1)!}\quad,\quad  
\cos x=\sum_{l=0}^\infty(-1)^l\frac{x^{2l}}{(2l)!}$$
$$e^x=\sum_{k=0}^\infty\frac{x^k}{k!}\quad,\quad 
\log(1+x)=\sum_{k=1}^\infty(-1)^{k+1}\frac{x^k}{k}$$
with $|x|<1$ needed for $\log$, which are continuous over their domain, and satisfy the formulae $e^{x+y}=e^xe^y$ and $e^{ix}=\cos x+i\sin x$.
\end{theorem}

\begin{proof}
This is a mixture of trivial and non-trivial results, as follows:

\medskip

(1) We already know about $e^x$ from chapter 5, the idea being that the convergence of the series, and then the continuity of $e^x$, come from the following estimate:
$$|e^x|\leq\sum_{k=0}^\infty\frac{|x|^k}{k!}=e^{|x|}<\infty$$

But the same can be said about $\sin$ and $\cos$, with the estimates being as follows:
$$|\sin x|\leq e^{|x|}\quad,\quad |\cos x|\leq e^{|x|}$$

(2) Regarding now the formulae satisfied by $\sin,\cos,\exp$, we already know from chapter 5 that the exponential has the following property, exactly as in the real case:
$$e^{x+y}=e^xe^y$$

We also have the following formula, connecting $\sin,\cos,\exp$, again as before:
$$e^{ix}
=\sum_{l=0}^\infty\frac{(ix)^{2l}}{(2l)!}+\sum_{l=0}^\infty\frac{(ix)^{2l+1}}{(2l+1)!}
=\cos x+i\sin x$$

(3) In order to discuss now the complex logarithm function $\log$, let us first study some more the complex exponential function $\exp$. By using $e^{x+y}=e^xe^y$ we obtain $e^x\neq0$ for any $x\in\mathbb C$, so the complex exponential function is as follows:
$$\exp:\mathbb C\to\mathbb C-\{0\}$$

Now since we have $e^{x+iy}=e^xe^{iy}$ for $x,y\in\mathbb R$, with $e^x$ being surjective onto $(0,\infty)$, and with $e^{iy}$ being surjective onto the unit circle $\mathbb T$, we deduce that $\exp:\mathbb C\to\mathbb C-\{0\}$ is surjective. Also, again by using $e^{x+iy}=e^xe^{iy}$, we deduce that we have:
$$e^x=e^y\iff x-y\in 2\pi i\mathbb Z$$

(4) With these ingredients in hand, we can now talk about $\log$. Indeed, let us fix a horizontal strip in the complex plane, having width $2\pi$: 
$$S=\left\{x+iy\Big|x\in\mathbb R,y\in[a,a+2\pi)\right\}$$

We know from the above that the restriction map $\exp:S\to\mathbb C-\{0\}$ is bijective, so we can define $\log$ as to be the inverse of this map:
$$\log=\exp^{-1}:\mathbb C-\{0\}\to S$$

(5) In practice now, the best is to choose for instance $a=0$, or $a=-\pi$, as to have the whole real line included in our strip, $\mathbb R\subset S$. In this case on $\mathbb R_+$ we recover the usual logarithm, while on $\mathbb R_-$ we obtain complex values, as for instance $\log(-1)=\pi i$ in the case $a=0$, or $\log(-1)=-\pi i$ in the case $a=-\pi$, coming from $e^{\pi i}=-1$. 

\medskip

(6) Finally, the formula for $\log$ follows as in the real case, via the computation in chapter 3, using the generalized binomial formula, with negative integer exponent.
\end{proof}

As an interesting consequence of the above result, which is of great practical interest, we have the following useful method, for remembering the basic math formulae:

\begin{method}\
Knowing $e^x=\sum_kx^k/k!$ and $e^{ix}=\cos x+i\sin x$ gives you
$$\sin(x+y)=\sin x\cos y+\cos x\sin y$$
$$\cos(x+y)=\cos x\cos y-\sin x\sin y$$
right away, in case you forgot these formulae, as well as
$$\sin x=\sum_{l=0}^\infty(-1)^l\frac{x^{2l+1}}{(2l+1)!}\quad,\quad 
\cos x=\sum_{l=0}^\infty(-1)^l\frac{x^{2l}}{(2l)!}$$
again, right away, in case you forgot these formulae.
\end{method}

To be more precise, knowing $e^{ix}=\cos x+i\sin x$, we can get right away:
\begin{eqnarray*}
e^{i(x+y)}=e^{ix}e^{iy}
&\implies&\cos(x+y)+i\sin(x+y)=(\cos x+i\sin x)(\cos y+i\sin y)\\
&\implies&\begin{cases}
\cos(x+y)=\cos x\cos y-\sin x\sin y\\
\sin(x+y)=\sin x\cos y+\cos x\sin y
\end{cases}
\end{eqnarray*}

Moreover, by further remembering $e^x=\sum_kx^k/k!$, we get right away, as well:
\begin{eqnarray*}
e^{ix}=\sum_k\frac{(ix)^k}{k!}
&\implies&\cos x+i\sin x=\sum_k\frac{(ix)^k}{k!}\\
&\implies&\begin{cases}
\cos x=\sum_{l=0}^\infty(-1)^l\frac{x^{2l}}{(2l)!}\\
\sin x=\sum_{l=0}^\infty(-1)^l\frac{x^{2l+1}}{(2l+1)!}
\end{cases}
\end{eqnarray*}

Moving ahead, Theorem 6.10 leads us into the question on whether the other formulae that we know about $\sin,\cos$, such as the values of these functions on sums $x+y$, or on doubles $2x$, extend to the complex setting. Things are quite tricky here, and in relation with this, we have the following result, which is something of general interest: 

\index{sinh}
\index{cosh}

\begin{theorem}
The following functions, called hyperbolic sine and cosine,
$$\sinh x=\frac{e^x-e^{-x}}{2}
\quad,\quad
\cosh x=\frac{e^x+e^{-x}}{2}$$
are subject to the following formulae:
\begin{enumerate}
\item $e^x=\cosh x+\sinh x$.

\item $\sinh(ix)=i\sin x$, $\cosh(ix)=\cos x$, for $x\in\mathbb R$.

\item $\sinh(x+y)=\sinh x\cosh y+\cosh x\sinh y$.

\item $\cosh(x+y)=\cosh x\cosh y+\sinh x\sinh y$.

\item $\sinh x=\sum_l\frac{x^{2l+1}}{(2l+1)!}$, $\cosh x=\sum_l\frac{x^{2l}}{(2l)!}$.
\end{enumerate}
\end{theorem}

\begin{proof}
The formula (1) follows from definitions. As for (2), this follows from:
$$\sinh(ix)=\frac{e^{ix}-e^{-ix}}{2}=\frac{\cos x+i\sin x}{2}-\frac{\cos x-i\sin x}{2}=i\sin x$$
$$\cosh(ix)=\frac{e^{ix}+e^{-ix}}{2}=\frac{\cos x+i\sin x}{2}+\frac{\cos x-i\sin x}{2}=\cos x$$

Regarding now (3,4), observe first that the formula $e^{x+y}=e^xe^y$ reads:
$$\cosh(x+y)+\sinh(x+y)=(\cosh x+\sinh x)(\cosh y+\sinh y)$$

Thus, we have some good explanation for (3,4), and in practice, these formulae can be checked by direct computation, as follows:
$$\frac{e^{x+y}-e^{-x-y}}{2}=
\frac{e^x-e^{-x}}{2}\cdot\frac{e^y+e^{-y}}{2}+
\frac{e^x+e^{-x}}{2}\cdot\frac{e^y-e^{-y}}{2}$$
$$\frac{e^{x+y}+e^{-x-y}}{2}=
\frac{e^x+e^{-x}}{2}\cdot\frac{e^y+e^{-y}}{2}+
\frac{e^x-e^{-x}}{2}\cdot\frac{e^y-e^{-y}}{2}$$

Finally, (5) is clear from the definition of $\sinh$, $\cosh$, and from $e^x=\sum_k\frac{x^k}{k!}$.
\end{proof}

Getting now to the Taylor series for our new beasts, many things can be said here. To start with, with obvious definitions for the remaining hyperbolic functions, namely $\tanh=\sinh/\cosh$ and so on, we have a total of 24 trigonometric functions, as follows:
$$\begin{matrix}
\sin&\cos&\tan&\sec&\csc&\cot\\
\arcsin&\arccos&\arctan&{\rm arcsec}&{\rm arccsc}&{\rm arccot}\\
\sinh&\cosh&\tanh&{\rm sech}&{\rm csch}&{\rm coth}\\
{\rm arcsinh}&{\rm arccosh}&{\rm arctanh}&{\rm arcsech}&{\rm arccsch}&{\rm arccoth}
\end{matrix}$$

Thus, a lot of work to be done. In analogy with what we did in chapter 3, where we only talked about $\arctan$, passed $\sin$ and $\cos$, we will do something similar here. Passed $\sinh$ and $\cosh$, an interesting hyperbolic function is $\coth$, and we have, about it:

\begin{theorem}
With the Bernoulli numbers $B_n\in\mathbb Q$ defined via
$$\frac{x}{e^x-1}=\sum_{n=0}^\infty\frac{B_n}{n!}\,x^n$$
the Taylor series of ${\rm coth}$ is given by the following formula:
$${\rm coth}\,x=\sum_{k=0}^\infty\frac{4^kB_{2k}}{(2k)!}\,x^{2k-1}$$
Moreover, we have similar formulae for $\tan$, $\csc$, $\cot$ and for $\tanh$, ${\rm csch}$. 
\end{theorem}

\begin{proof}
Many things going on here, the idea being as follows:

\medskip

(1) To start with, with some patience we can compute $x/(e^x-1)$:
$$\frac{x}{e^x-1}=\frac{1}{1+\frac{x}{2}+\frac{x^2}{6}+\frac{x^3}{24}+\ldots}=1-\frac{x}{2}+\frac{x^2}{12}-\frac{x^4}{720}+\frac{x^6}{30240}-\ldots$$

Which seems to lead nowhere, but by some magic, the coefficients are $c_n=B_n/n!$, with $B_n\in\mathbb Q$ being the Bernoulli numbers, that we met in chapter 4.

\medskip

(2) But with this, we can finish. Indeed, $c_n=B_n/n!$ can be shown to hold indeed, by recurrence, and then the passage to ${\rm coth}$ is straightforward, coming from:
$$x({\rm coth}\,x-1)
=x\left(\frac{e^x+e^{-x}}{e^x-e^{-x}}-1\right)
=x\left(\frac{e^{2x}+1}{e^{2x}-1}-1\right)
=\frac{2x}{e^{2x}-1}$$

(3) As in what regards $\tan$, $\csc$, $\cot$ and $\tanh$, ${\rm csch}$, good exercise for you. Finally, for our story to be complete, the 24 trigonometric functions fall in 5 classes:

\smallskip

-- $\sin$, $\cos$, $\sinh$, $\cosh$, whose Taylor series involve factorials $k!$.

\smallskip

-- $\arctan$, ${\rm arccot}$, ${\rm arctanh}$, ${\rm arccoth}$, which only need inverses $1/k$.

\smallskip

-- $\arcsin$, $\arccos$, ${\rm arcsinh}$, ${\rm arccosh}$, ${\rm arcsec}$, ${\rm arccsc}$, ${\rm arcsech}$, ${\rm arccsch}$, needing $D_k=\binom{2k}{k}$.

\smallskip

-- $\tan$, $\csc$, $\cot$, $\tanh$, ${\rm csch}$, ${\rm coth}$, needing the Bernoulli numbers $B_k$.

\smallskip

-- $\sec$, ${\rm sech}$, needing yet another type of numbers, the Euler ones $E_k$.

\smallskip

And exercise of course for you, to have the second and third classes fully worked out, and then to learn more about the fourth and fifth classes, say from my book \cite{ba1}.
\end{proof}

Moving on, we can talk as well about complex powers, in the following way:

\index{power function}

\begin{fact}
Under suitable assumptions, we can talk about 
$$x^y=e^{y\log x}$$
with $x,y\in\mathbb C$, and in particular about the functions $a^x$ and $x^a$, with $a\in\mathbb C$.
\end{fact}

However, all this is quite complicated, and we will leave some study or learning here as an exercise. Finally, as something easier and refreshing, let us record the following general result, standing as an intermediate value theorem, for the complex functions:

\index{connected set}
\index{intermediate value}

\begin{theorem}
Assuming that $f:X\to\mathbb C$ with $X\subset\mathbb C$ is continuous, if the domain $X$ is compact and connected, then so is its image $f(X)$.
\end{theorem}

\begin{proof}
This follows exactly as in the real case, with just a bit of discussion being needed, in relation with the open, closed, compact and connected sets, inside $\mathbb C$. We will leave this as an exercise, and come back to it later, directly in the general $\mathbb R^N$ setting.
\end{proof}
 
\section*{6b. Holomorphic functions}

Let us study now the differentiability of the complex functions $f:\mathbb C\to\mathbb C$. Things here are quite tricky, but let us start with a straightforward definition, as follows:

\index{differentiable function}
\index{complex function}
\index{holomorphic function}

\begin{definition}
We say that a function $f:X\to\mathbb C$ is differentiable in the complex sense when the following limit is defined for any $x\in X$:
$$f'(x)=\lim_{t\to0}\frac{f(x+t)-f(x)}{t}$$
In this case, we also say that $f$ is holomorphic, and we write $f\in H(X)$.
\end{definition}

As basic examples, we have the power functions $f(x)=x^n$. Indeed, the derivative of such a power function can be computed exactly as in the real case, and we get:
\begin{eqnarray*}
(x^n)'
&=&\lim_{t\to0}\frac{(x+t)^n-x^n}{t}\\
&=&\lim_{t\to0}\frac{nx^{n-1}t+\binom{n}{2}x^{n-2}t^2+\ldots+t^n}{t}\\
&=&\lim_{t\to0}\frac{nx^{n-1}t}{t}\\
&=&nx^{n-1}
\end{eqnarray*}

We will see later more computations of this type, similar to those from the real case. To summarize, our definition of differentiability seems to work nicely, so let us start developing some theory. The general results from the real case extend well, as follows:

\begin{proposition}
We have the following results:
\begin{enumerate}
\item $(f+g)'=f'+g'$.

\item $(\lambda f)'=\lambda f'$.

\item $(fg)'=f'g+fg'$.

\item $(f\circ g)'=f'(g)g'$.
\end{enumerate}
\end{proposition}

\begin{proof}
These formulae are all clear from definitions, following exactly as in the real case. Thus, we are led to the conclusions in the statement.
\end{proof}

As an obvious consequence of (1,2) above, any poynomial $P\in\mathbb C[X]$ is differentiable, with its derivative being given by the same formula as in the real case, namely:
$$P(x)=\sum_{k=0}^nc_kx^k\implies P'(x)=\sum_{k=1}^nkc_kx^{k-1}$$

More generally, any rational function $f\in\mathbb C(X)$ is differentiable on its domain, that is, outside its poles, because if we write $f=P/Q$ with $P,Q\in\mathbb C[X]$, we have:
$$f'=\left(\frac{P}{Q}\right)'=\frac{P'Q-PQ'}{Q^2}$$

Let us record these conclusions in a statement, as follows:

\index{rational function}
\index{infinitely differentiable}

\begin{proposition}
The following happen:
\begin{enumerate}
\item Any polynomial $P\in\mathbb C[X]$ is holomorphic, and in fact infinitely differentiable in the complex sense, with all its derivatives being polynomials.

\item Any rational function $f\in\mathbb C(X)$ is holomorphic, and in fact infinitely differentiable, with all its derivatives being rational functions.
\end{enumerate}
\end{proposition}

\begin{proof}
This follows indeed from the above discussion.
\end{proof}

Regarding now more complicated complex functions that we know, we have here:

\index{complex conjugate}
\index{modulus}

\begin{theorem}
The following happen:
\begin{enumerate}
\item $\sin,\cos,\exp,\log$ are holomorphic, and in fact are infinitely differentiable, with their derivatives being given by the same formulae as in the real case. 

\item However, functions like $\bar{x}$ or $|x|$ are not holomorphic, and this because the limit defining $f'(x)$ depends on the way we choose $t\to0$.
\end{enumerate}
\end{theorem}

\begin{proof}
There are several things going on here, the idea being as follows:

\medskip

(1) Here the first assertion is standard, because our functions $\sin,\cos,\exp,\log$ have Taylor series that we know, and the derivative can be therefore computed by using the same rule as in the real case, similar to the one for polynomials, namely:
$$f(x)=\sum_{k=0}^\infty c_kx^k\implies f'(x)=\sum_{k=1}^\infty kc_kx^{k-1}$$

(2) Regarding now the function $f(x)=\bar{x}$, the point here is that we have:
$$\frac{f(x+t)-f(x)}{t}=\frac{\bar{x}+\bar{t}-\bar{x}}{t}=\frac{\bar{t}}{t}$$

But this limit does not converge with $t\to0$, for instance because with $t\in\mathbb R$ we obtain 1 as limit, while with $t\in i\mathbb R$ we obtain $-1$ as limit. In fact, with $t=rw$ with $|w|=1$ fixed and $r\in\mathbb R$, $r\to0$, we can obtain as limit any number on the unit circle:
$$\lim_{r\to0}\frac{f(x+rw)-f(x)}{rw}=\lim_{r\to0}\frac{r\bar{w}}{rw}=\bar{w}^2$$

(3) The situation for the function $f(x)=|x|$ is similar. To be more precise, we have:
$$\frac{f(x+rw)-f(x)}{rw}=\frac{|x+rw|-|x|}{r}\cdot\bar{w}$$

Thus with $|w|=1$ fixed and $r\to0$ we obtain a certain multiple of $\bar{w}$, with the multiplication factor being computed as follows:
$$\frac{|x+rw|-|x|}{r}
=\frac{|x+rw|^2-|x|^2}{(|x+rw|+|x|)r}
\simeq\frac{xr\bar{w}+\bar{x}rw}{2|x|r}
=Re\left(\frac{x\bar{w}}{|x|}\right)$$

Thus, as before in (1), we are led to the conclusion in the statement.
\end{proof}

The above result is quite surprising, because we are so used, from the real case, to the notion of differentiability to correspond to some form of ``smoothness'' of the function, and to be more precise, to ``smoothness at first order''. Or, if you prefer, to correspond to the ``non-bumpiness'' of the function. So, we are led to the following dilemma:

\begin{dilemma}
It's either that $\bar{x}$ and $|x|$ are smooth, as the intuition suggests, and we are wrong with our definition of differentiability. Or that $\bar{x}$ and $|x|$ are bumpy, while this being not very intuitive, and we are right with our definition of differentiability.
\end{dilemma}

And we won't get discouraged by this. After all, this is just some empty talking, and if there is something to rely upon, mathematics and computations, these are the computations from the proof of Theorem 6.19. So, moving ahead now, based on that computations, let us formulate the following definition, coming as a complement to Definition 6.16:

\index{radial limit}

\begin{definition}
A function $f:X\to\mathbb C$ is called differentiable:
\begin{enumerate}
\item In the real sense, if the following two limits converge, for any $x\in X$:
$$f_1'(x)=\lim_{t\in\mathbb R\to0}\frac{f(x+t)-f(x)}{t}\quad,\quad 
f_i'(x)=\lim_{t\in i\mathbb R\to0}\frac{f(x+t)-f(x)}{t}$$

\item In a radial sense, if the following limit converges, for any $x\in X$, and $w\in\mathbb T$:
$$f_w'(x)=\lim_{t\in w\mathbb R\to0}\frac{f(x+t)-f(x)}{t}$$

\item In the complex sense, if the following limit converges, for any $x\in X$:
$$f'(x)=\lim_{t\to0}\frac{f(x+t)-f(x)}{t}$$
\end{enumerate} 
If $f$ is differentiable in the complex sense, we also say that $f$ is holomorphic.
\end{definition}

With this, we have $(3)\implies(2)\implies(1)$, and most of the functions that we know, namely the polynomials, the rational functions, and $\sin,\cos,\exp,\log$, satisfy (3). However, the functions $\bar{x},|x|$ satisfy (2), which is not bad, but do not satisfy (3).

\bigskip

Moving on, all the examples of holomorphic functions that we have are infinitely differentiable, and this raises the question of finding a function such that $f'$ exists, while $f''$ does not exist. Quite surprisingly, we will see that such functions do not exist. In order to get into this latter phenomenon, let us start with:

\begin{theorem}
Each power series $f(x)=\sum_nc_nx^n$ has a radius of convergence 
$$R\in[0,\infty]$$
which is such that $f$ converges for $|x|<R$, and diverges for $|x|>R$. We have:
$$R=\frac{1}{C}\quad,\quad C=\limsup_{n\to\infty}\sqrt[n]{|c_n|}$$
Also, in the case $|x|=R$ the function $f$ can either converge, or diverge.
\end{theorem}

\begin{proof}
This follows from the Cauchy criterion for series, from chapter 1, which says that a series $\sum_nx_n$ converges if $c<1$, and diverges if $c>1$, where:
$$c=\limsup_{n\to\infty}\sqrt[n]{|x_n|}$$

Indeed, with $x_n=|c_nx^n|$ we obtain that the convergence radius $R\in[0,\infty]$ exists, and is given by the formula in the statement. Finally, for the examples and counterexamples at the end, when $|x|=R$, the simplest here is to use $f(x)=\sum_nx^n/n$, where $R=1$.
\end{proof}

Back now to our questions regarding higher derivatives, we have:

\index{analytic function}

\begin{theorem}
Assuming that a function $f:X\to\mathbb C$ is analytic, in the sense that it is a series, around each point $x\in X$, 
$$f(x+t)=\sum_{n=0}^\infty c_nt^n$$
it follows that $f$ is infinitely differentiable, in the complex sense. In particular, $f'$ exists, and so $f$ is holomorphic in our sense.
\end{theorem}

\begin{proof}
Assuming that $f$ is analytic, as in the statement, we have:
$$f'(x+t)=\sum_{n=1}^\infty nc_nt^{n-1}$$

Moreover, the radius of convergence is the same, as shown by the following computation, using the Cauchy formula for the convergence radius, and $\sqrt[n]{n}\to1$:
$$\frac{1}{R'}
=\limsup_{n\to\infty}\sqrt[n]{|nc_n|}
=\limsup_{n\to\infty}\sqrt[n]{|c_n|}
=\frac{1}{R}$$

Thus $f'$ exists and is analytic, on the same domain, and this gives the result.
\end{proof}

\section*{6c. Cauchy formula} 

Quite interesting all the above, and our goal in what follows will be that of proving that any holomorphic function is analytic. However, this is something quite subtle, which cannot be proved with bare hands, and requires substantial preliminaries. 

\bigskip

Getting to these preliminaries, let us start with some heuristics. Our claim is that a lot of useful knowledge, in order to deal with the holomorphic functions, can be gained by further studying the analytic functions, and even the usual polynomials $P\in\mathbb C[X]$.

\bigskip

So, let us further study the polynomials $P\in\mathbb C[X]$, and other analytic functions. We already know from chapter 5 that in the polynomial case, $P\in\mathbb C[X]$, some interesting things happen, because any such polynomial has a root, and even $\deg(P)$ roots, after a recurrence. Keeping looking at polynomials, with the same methods, we are led to:

\index{boundary of domain}
\index{maximum principle}

\begin{theorem}
Any polynomial $P\in\mathbb C[X]$ satisfies the maximum principle, in the sense that given a disk $D$, with boundary $\gamma$, we have:
$$\exists x\in\gamma\quad,\quad |P(x)|=\max_{y\in D}|P(y)|$$
In fact, the maximum of $|P|$ over any domain is attained on the boundary.
\end{theorem}

\begin{proof}
This can be proved by contradiction, by using the same arguments as in the proof of the existence of a root, from chapter 5. To be more precise, assume $\deg P\geq 1$, and that the maximum of $|P|$ is attained at the center of a disk $D=D(z,r)$:
$$|P(z)|=\max_{x\in D}|P(x)|$$

We can write then $P(z+t)\simeq P(z)+ct^k$ with $c\neq0$, for $t$ small, and by suitably choosing the argument of $t$ on the unit circle we conclude, exactly as in chapter 5, that the function $|P|$ cannot have a local maximum at $z$, as stated.
\end{proof}

A good explanation for the fact that the maximum principle holds for polynomials $P\in\mathbb C[X]$ could be that the values of such a polynomial inside a disk can be recovered from its values on the boundary. And fortunately, this is indeed the case:

\index{main value formula}

\begin{theorem}
Given a polynomial $P\in\mathbb C[X]$, and a disk $D$, with boundary $\gamma$, we have the following formulae, with the integrations being the normalized, mass $1$ ones:
\begin{enumerate}
\item $P$ satisfies the plain mean value formula $P(x)=\int_DP(y)dy$.

\item $P$ satisfies the boundary mean value formula $P(x)=\int_\gamma P(y)dy$.
\end{enumerate}
\end{theorem}

\begin{proof}
As a first observation, the two mean value formulae in the statement are equivalent, by restricting the attention to disks $D$, having as boundaries circles $\gamma$, and using annuli and polar coordinates for the proof of the equivalence. As for the formulae themselves, these can be checked by direct computation for a disk $D$, with the formulation in (2) being the most convenient. Indeed, for a monomial $P(x)=x^n$ we have:
\begin{eqnarray*}
\int_\gamma y^ndy
&=&\frac{1}{2\pi}\int_0^{2\pi}(x+re^{it})^ndt\\
&=&\frac{1}{2\pi}\int_0^{2\pi}\sum_{k=0}^n\binom{n}{k}x^k(re^{it})^{n-k}dt\\
&=&\sum_{k=0}^n\binom{n}{k}x^kr^{n-k}\frac{1}{2\pi}\int_0^{2\pi}e^{i(n-k)t}dt\\
&=&\sum_{k=0}^n\binom{n}{k}x^kr^{n-k}\delta_{kn}\\
&=&x^n
\end{eqnarray*}

Thus, we have the result for monomials, and the general case follows by linearity.
\end{proof}

Finally, in relation with some previous speculations from chapter 5, we have:

\begin{theorem}
Assuming that $f:X\to\mathbb C$ is $n$ times differentiable, we have
$$\frac{\sum_{s=1}^nf(x+tw^s)}{n}\simeq f(x)+\frac{f^{(n)}(x)}{n!}\,t^n$$
with $w=e^{2\pi i/n}$, suggesting that $f$ should satisfy the mean value formula.
\end{theorem}

\begin{proof}
We have indeed the following computation, using the Taylor formula:
\begin{eqnarray*}
\frac{\sum_{s=1}^nf(x+tw^s)}{n}
&\simeq&\frac{1}{n}\sum_{s=1}^n\sum_{k=0}^n\frac{f^{(k)}(x)}{k!}\,(tw^s)^k\\
&=&\sum_{k=0}^n\frac{f^{(k)}(x)}{k!}\,t^k\cdot\frac{1}{n}\sum_{s=1}^n(w^k)^s\\
&=&f(x)+\frac{f^{(n)}(x)}{n!}\,t^n
\end{eqnarray*}

As for the last assertion, this is obviously just a speculation, based on this.
\end{proof}

All the above is very nice, but we can in fact do even better. We first have:

\index{integral over curve}

\begin{proposition}
We can integrate functions $f$ over curves $\gamma$ by setting
$$\int_\gamma f(x)dx=\int_a^bf(\gamma(t))\gamma'(t)dt$$
with this quantity being independent on the parametrization $\gamma:[a,b]\to\mathbb C$.
\end{proposition}

\begin{proof}
We must prove that the quantity in the statement is independent on the parametrization. In other words, we must prove that if we pick two different parametrizations $\gamma,\eta:[a,b]\to\mathbb C$ of our curve, then we have the following formula:
$$\int_a^bf(\gamma(t))\gamma'(t)dt=\int_a^bf(\eta(t))\eta'(t)dt$$

But for this purpose, let us write $\gamma=\eta\phi$, with $\phi:[a,b]\to[a,b]$ being a certain function, that we can assume to be bijective, via an elementary cut-and-paste argument. By using the chain rule for derivatives, and the change of variable formula, we have:
\begin{eqnarray*}
\int_a^bf(\gamma(t))\gamma'(t)dt
&=&\int_a^bf(\eta\phi(t))(\eta\phi)'(t)dt\\
&=&\int_a^bf(\eta\phi(t))\eta'(\phi(t))\phi'(t)dt\\
&=&\int_a^bf(\eta(t))\eta'(t)dt
\end{eqnarray*}

Thus, we are led to the conclusions in the statement.
\end{proof}

The main properties of the above integration method are as follows:

\begin{proposition}
We have the following formula, for a union of paths:
$$\int_{\gamma\cup\eta}f(x)dx=\int_\gamma f(x)dx+\int_\eta f(x)dx$$
Also, when reversing the path, the integral changes its sign.
\end{proposition}

\begin{proof}
Here the first assertion is clear from definitions, and the second assertion comes from the change of variable formula, by using Proposition 6.27.
\end{proof}

Now by getting back to polynomials, we have the following result:

\index{Cauchy formula}

\begin{theorem}
Any polynomial $P\in\mathbb C[X]$ satisfies the Cauchy formula
$$P(x)=\frac{1}{2\pi i}\int_\gamma\frac{P(y)}{y-x}\,dy$$
with the integration over $\gamma$ being constructed as above.
\end{theorem}

\begin{proof}
This follows by using abstract arguments and computations similar to those in the proof of Theorem 6.24. Indeed, by linearity we can assume $P(x)=x^n$. Also, by using a cut-and-paste argument, we can assume that we are on a circle:
$$\gamma:[0,2\pi]\to\mathbb C\quad,\quad \gamma(t)=x+re^{it}$$

By using now the computation from the proof of Theorem 6.24, we obtain:
\begin{eqnarray*}
\int_\gamma\frac{y^n}{y-x}\,dy
&=&\int_0^{2\pi}\frac{(x+re^{it})^n}{re^{it}}\,rie^{it}dt\\
&=&i\int_0^{2\pi}(x+re^{it})^ndt\\
&=&i\cdot 2\pi x^n
\end{eqnarray*}

Thus, we are led to the formula in the statement.
\end{proof}

All this is quite interesting, and obviously, we are now into some serious mathematics. Importantly, Theorem 6.24, Theorem 6.25 and Theorem 6.29 provide us with a path for proving the converse of Theorem 6.23. Indeed, if we manage to prove the Cauchy formula for any holomorphic function $f:X\to\mathbb C$, then it will follow that our function is in fact analytic, and so infinitely differentiable. So, let us start with the following result:

\index{Cauchy formula}

\begin{theorem}
The Cauchy formula, namely
$$f(x)=\frac{1}{2\pi i}\int_\gamma\frac{f(y)}{y-x}\,dy$$
holds for any holomorphic function $f:X\to\mathbb C$.
\end{theorem}

\begin{proof}
This is something standard, which can be proved as follows:

\medskip

(1) Our first claim is that given $f\in H(X)$, with $f'\in C(X)$, the integral of $f'$ vanishes on any path. Indeed, by using the change of variable formula, we have:
\begin{eqnarray*}
\int_\gamma f'(x)dx
&=&\int_a^bf'(\gamma(t))\gamma'(t)dt\\
&=&f(\gamma(b))-f(\gamma(a))\\
&=&0
\end{eqnarray*}

(2) Our second claim is that given $f\in H(X)$ and a triangle $\Delta\subset X$, we have:
$$\int_\Delta f(x)dx=0$$

Indeed, let us call $\Delta=ABC$ our triangle. Now consider the midpoints $A',B',C'$ of the edges $BC,CA,AB$, and then consider the following smaller triangles:
$$\Delta_1=AC'B'\quad,\quad\Delta_2=BA'C'\quad,\quad\Delta_3=CB'A'\quad,\quad\Delta_4=A'B'C'$$

These smaller triangles partition then $\Delta$, and due to our above conventions for the vertex ordering, which produce cancellations when integrating over them, we have:
$$\int_\Delta f(x)dx=\sum_{i=1}^4\int_{\Delta_i}f(x)dx$$

Thus we can pick, among the triangles $\Delta_i$, a triangle $\Delta^{(1)}$ such that:
$$\left|\int_\Delta f(x)dx\right|\leq 4\left|\int_{\Delta^{(1)}}f(x)dx\right|$$

In fact, by repeating the procedure, we obtain triangles $\Delta^{(n)}$ such that:
$$\left|\int_\Delta f(x)dx\right|\leq 4^n\left|\int_{\Delta^{(n)}}f(x)dx\right|$$

(3) Now let $z$ be the limiting point of these triangles $\Delta^{(n)}$, and fix $\varepsilon>0$. By using the fact that the functions $1,x$ integrate over paths up to 0, coming from (1), we obtain the following estimate, with $n\in\mathbb N$ being big enough, and $L$ being the perimeter of $\Delta$:
\begin{eqnarray*}
\left|\int_{\Delta^{(n)}}f(x)dx\right|
&=&\left|\int_{\Delta^{(n)}}f(x)-f(z)-f'(z)(x-z)dx\right|\\
&\leq&\int_{\Delta^{(n)}}\left|f(x)-f(z)-f'(z)(x-z)\right|dx\\
&\leq&\int_{\Delta^{(n)}}\varepsilon|x-z|dx\\
&\leq&\varepsilon(2^{-n}L)^2
\end{eqnarray*}

Now by combining this with the estimate in (2), this proves our claim.

\medskip

(4) The rest is quite routine. First, we can pass from triangles to boundaries of convex sets, in a straightforward way, with the same conclusion as in (2), namely:
$$\int_\gamma f(x)dx=0$$

Getting back to what we want to prove, namely the Cauchy formula for an arbitrary holomorphic function $f\in H(X)$, let $x\in X$, and consider the following function:
$$g(y)=\begin{cases}
\frac{f(y)-f(x)}{y-x}&(y\neq x)\\
f'(x)&(y=x)
\end{cases}$$

Now assuming that $\gamma$ encloses a convex set, we can apply what we found, namely vanishing of the integral, to this function $g$, and we obtain the Cauchy formula for $f$.

\medskip

(5) Finally, the extension to general curves is standard, and standard as well is the discussion of what exactly happens at $x$, in the above proof. See Rudin \cite{ru2}.
\end{proof}

As a main application of the Cauchy formula, we have:

\index{holomorphic function}
\index{infinitely differentiable}
\index{analytic function}
\index{Cauchy formula}

\begin{theorem}
The following conditions are equivalent, for a function $f:X\to\mathbb C$:
\begin{enumerate}
\item $f$ is holomorphic.

\item $f$ is infinitely differentiable.

\item $f$ is analytic.

\item The Cauchy formula holds for $f$.
\end{enumerate}
\end{theorem}

\begin{proof}
This is routine from what we have, the idea being as follows:

\medskip

$(1)\implies(4)$ is non-trivial, but we know this from Theorem 6.30.

\medskip

$(4)\implies(3)$ is something trivial, because we can expand the series in the Cauchy formula, and we conclude that our function is indeed analytic.

\medskip

$(3)\implies(2)\implies(1)$ are both elementary, known from Theorem 6.23.
\end{proof}

As another application of the Cauchy formula, we have:

\index{maximum principle}

\begin{theorem}
Any holomorphic function $f:X\to\mathbb C$ satisfies the maximum principle, in the sense that given a domain $D$, with boundary $\gamma$, we have:
$$\exists x\in\gamma\quad,\quad |f(x)|=\max_{y\in D}|f(y)|$$
That is, the maximum of $|f|$ over a domain is attained on the boundary.
\end{theorem}

\begin{proof}
This follows indeed from the Cauchy formula. Observe that the converse is not true, for instance because functions like $\bar{x}$ satisfy too the maximum principle. We will be back to this later, in chapter 8, when talking about harmonic functions.
\end{proof}

As before with polynomials, a good explanation for the fact that the maximum principle holds could be that the values of our function inside a disk can be recovered from its values on the boundary. And fortunately, this is indeed the case, and we have:

\index{main value formula}

\begin{theorem}
Given an holomorphic function $f:X\to\mathbb C$, and a disk $D$, with boundary $\gamma$, the following happen:
\begin{enumerate}
\item $f$ satisfies the plain mean value formula $f(x)=\int_Df(y)dy$.

\item $f$ satisfies the boundary mean value formula $f(x)=\int_\gamma f(y)dy$.
\end{enumerate}
\end{theorem}

\begin{proof}
As usual, this follows from the Cauchy formula, with of course some care in passing from integrals constructed as in Proposition 6.27 to integrals viewed as averages, which are those that we refer to, in the present statement.
\end{proof}

Finally, as yet another application of the Cauchy formula, due to Liouville, we have:

\index{Liouville theorem}

\begin{theorem}
An entire, bounded holomorphic function
$$f:\mathbb C\to\mathbb C\quad,\quad |f|\leq M$$
must be constant. In particular, knowing $f\to0$ with $z\to\infty$ gives $f=0$.
\end{theorem}

\begin{proof}
This follows as usual from the Cauchy formula. Alternatively, we can view this as a consequence of Theorem 6.33, because given two points $x\neq y$, we can view the values of $f$ at these points as averages over big disks centered at these points, say $D=D_x(R)$ and $E=D_y(R)$, with $R>>0$, with the formulae being as follows:
$$f(x)=\int_Df(z)dz\quad,\quad f(y)=\int_Ef(z)dz$$

Indeed, the point is that when the radius goes to $\infty$, these averages tend to be equal, and so we have $f(x)\simeq f(y)$, which gives $f(x)=f(y)$ in the limit.
\end{proof}

\section*{6d. Stieltjes inversion}

We would like to end this chapter with an interesting application of the complex functions to probability theory. We learned some basic probability in chapter 4, and in view of the material there, an interesting question is how to recover a probability measure out of its moments. And the answer here, which is non-trivial, is as follows:

\begin{theorem}
The density of a real probability measure $\mu$ can be recaptured from the sequence of moments $\{M_k\}_{k\geq0}$ via the Stieltjes inversion formula
$$d\mu (x)=\lim_{t\searrow 0}-\frac{1}{\pi}\,Im\left(G(x+it)\right)\cdot dx$$
where the function on the right, given in terms of moments by
$$G(\xi)=\xi^{-1}+M_1\xi^{-2}+M_2\xi^{-3}+\ldots$$
is the Cauchy transform of the measure $\mu$.
\end{theorem}

\begin{proof}
The Cauchy transform of our measure $\mu$ is given by:
\begin{eqnarray*}
G(\xi)
&=&\xi^{-1}\sum_{k=0}^\infty M_k\xi^{-k}\\\
&=&\int_\mathbb R\frac{\xi^{-1}}{1-\xi^{-1}y}\,d\mu(y)\\
&=&\int_\mathbb R\frac{1}{\xi-y}\,d\mu(y)
\end{eqnarray*}

Now with $\xi=x+it$, we obtain the following formula:
\begin{eqnarray*}
Im(G(x+it))
&=&\int_\mathbb RIm\left(\frac{1}{x-y+it}\right)d\mu(y)\\
&=&\int_\mathbb R\frac{1}{2i}\left(\frac{1}{x-y+it}-\frac{1}{x-y-it}\right)d\mu(y)\\
&=&-\int_\mathbb R\frac{t}{(x-y)^2+t^2}\,d\mu(y)
\end{eqnarray*}

By integrating over $[a,b]$ we obtain, with the change of variables $x=y+tz$:
\begin{eqnarray*}
\int_a^bIm(G(x+it))dx
&=&-\int_\mathbb R\int_a^b\frac{t}{(x-y)^2+t^2}\,dx\,d\mu(y)\\
&=&-\int_\mathbb R\int_{(a-y)/t}^{(b-y)/t}\frac{t}{(tz)^2+t^2}\,t\,dz\,d\mu(y)\\
&=&-\int_\mathbb R\int_{(a-y)/t}^{(b-y)/t}\frac{1}{1+z^2}\,dz\,d\mu(y)\\
&=&-\int_\mathbb R\left(\arctan\frac{b-y}{t}-\arctan\frac{a-y}{t}\right)d\mu(y)
\end{eqnarray*}

Now observe that with $t\searrow0$ we have:
$$\lim_{t\searrow0}\left(\arctan\frac{b-y}{t}-\arctan\frac{a-y}{t}\right)
=\begin{cases}
\frac{\pi}{2}-\frac{\pi}{2}=0& (y<a)\\
\frac{\pi}{2}-0=\frac{\pi}{2}& (y=a)\\
\frac{\pi}{2}-(-\frac{\pi}{2})=\pi& (a<y<b)\\
0-(-\frac{\pi}{2})=\frac{\pi}{2}& (y=b)\\
-\frac{\pi}{2}-(-\frac{\pi}{2})=0& (y>b)
\end{cases}$$

We therefore obtain the following formula:
$$\lim_{t\searrow0}\int_a^bIm(G(x+it))dx=-\pi\left(\mu(a,b)+\frac{\mu(a)+\mu(b)}{2}\right)$$

Thus, we are led to the conclusion in the statement.
\end{proof}

Before getting further, let us mention that the above result does not fully solve the moment problem, because we still have the question of understanding when a sequence of numbers $M_1,M_2,M_3,\ldots$ can be the moments of a measure $\mu$.  We have here:

\begin{theorem}
A sequence of numbers $M_0,M_1,M_2,M_3,\ldots\in\mathbb R$, with $M_0=1$, is the series of moments of a real probability measure $\mu$ precisely when:
$$\begin{vmatrix}M_0\end{vmatrix}\geq0\quad,\quad 
\begin{vmatrix}
M_0&M_1\\
M_1&M_2
\end{vmatrix}\geq0\quad,\quad 
\begin{vmatrix}
M_0&M_1&M_2\\
M_1&M_2&M_3\\
M_2&M_3&M_4\\
\end{vmatrix}\geq0\quad,\quad 
\ldots$$
That is, the associated Hankel determinants must be all positive.
\end{theorem}

\begin{proof}
This is something a bit more advanced, the idea being as follows:

\medskip

(1) As a first observation, the positivity conditions in the statement tell us that the following associated linear forms must be positive:
$$\sum_{i,j=1}^nc_i\bar{c}_jM_{i+j}\geq0$$

(2) But this is something very classical, in one sense the result being elementary, coming from the following computation, which shows that we have positivity indeed:
$$\int_\mathbb R\left|\sum_{i=1}^nc_ix^i\right|^2d\mu(x)
=\int_\mathbb R\sum_{i,j=1}^nc_i\bar{c}_jx^{i+j}d\mu(x)
=\sum_{i,j=1}^nc_i\bar{c}_jM_{i+j}$$

(3) As for the other sense, here the result comes once again from the above formula, this time via some standard functional analysis.
\end{proof}

Getting back now to more concrete things, the point is that we have:

\begin{fact}
Given a graph $X$, with distinguished vertex $*$, we can talk about the probability measure $\mu$ having as $k$-th moment the number of length $k$ loops based at $*$:
$$M_k=\Big\{\!*=i_0-i_1-i_2-\ldots-i_{k-1}-i_k=*\Big\}$$
As basic examples, for the graph $\mathbb N$ the even moments are the Catalan numbers $C_k$, and for the graph $\mathbb Z$, the even moments are the central binomial coefficients $D_k$.
\end{fact}

To be more precise, the first assertion, regarding the existence and uniqueness of $\mu$, follows from a basic linear algebra computation, by diagonalizing the adjacency matrix of $X$. As for the examples, for the graph $\mathbb N$ we more or less already know this, from our various Catalan number considerations from chapter 2, and for the graph $\mathbb Z$ this is something elementary, that we will leave here as an instructive exercise.

\bigskip

Needless to say, counting loops on graphs, as in Fact 6.37, is something important in applied mathematics, and physics. So, back to our business now, motivated by all this, as a basic application of the Stieltjes formula, let us solve the moment problem for the Catalan numbers $C_k$, and for the central binomial coefficients $D_k$. We first have:

\begin{theorem}
The real measure having as even moments the Catalan numbers, $C_k=\frac{1}{k+1}\binom{2k}{k}$, and having all odd moments $0$ is the measure
$$\gamma_1=\frac{1}{2\pi}\sqrt{4-x^2}dx$$
called Wigner semicircle law on $[-2,2]$.
\end{theorem}

\begin{proof}
In order to apply the inversion formula, our starting point will be the formula from chapter 2 for the generating series of the Catalan numbers, namely:
$$\sum_{k=0}^\infty C_kz^k=\frac{1-\sqrt{1-4z}}{2z}$$

By using this formula with $z=\xi^{-2}$, we obtain the following formula:
\begin{eqnarray*}
G(\xi)
&=&\xi^{-1}\sum_{k=0}^\infty C_k\xi^{-2k}\\
&=&\xi^{-1}\cdot\frac{1-\sqrt{1-4\xi^{-2}}}{2\xi^{-2}}\\
&=&\frac{\xi}{2}\left(1-\sqrt{1-4\xi^{-2}}\right)\\
&=&\frac{\xi}{2}-\frac{1}{2}\sqrt{\xi^2-4}
\end{eqnarray*}

Now let us apply Theorem 6.35. The study here goes as follows:

\medskip

(1) According to the general philosophy of the Stieltjes formula, the first term, namely $\xi/2$, which is ``trivial'', will not contribute to the density. 

\medskip

(2) As for the second term, which is something non-trivial, this will contribute to the density, the rule here being that the square root $\sqrt{\xi^2-4}$ will be replaced by the ``dual'' square root $\sqrt{4-x^2}\,dx$, and that we have to multiply everything by $-1/\pi$. 

\medskip

(3) As a conclusion, by Stieltjes inversion we obtain the following density:
$$d\mu(x)
=-\frac{1}{\pi}\cdot-\frac{1}{2}\sqrt{4-x^2}\,dx
=\frac{1}{2\pi}\sqrt{4-x^2}dx$$

Thus, we have obtained the mesure in the statement, and we are done.
\end{proof}

We have as well the following interesting version of the above result:

\begin{theorem}
The real measure having as sequence of moments the Catalan numbers, $C_k=\frac{1}{k+1}\binom{2k}{k}$, is the measure
$$\pi_1=\frac{1}{2\pi}\sqrt{4x^{-1}-1}\,dx$$
called Marchenko-Pastur law on $[0,4]$.
\end{theorem}

\begin{proof}
As before, we use the standard formula for the generating series of the Catalan numbers. With $z=\xi^{-1}$ in that formula, we obtain the following formula:
\begin{eqnarray*}
G(\xi)
&=&\xi^{-1}\sum_{k=0}^\infty C_k\xi^{-k}\\
&=&\xi^{-1}\cdot\frac{1-\sqrt{1-4\xi^{-1}}}{2\xi^{-1}}\\
&=&\frac{1}{2}-\frac{1}{2}\sqrt{1-4\xi^{-1}}
\end{eqnarray*}

With this in hand, let us apply now the Stieltjes inversion formula, from Theorem 6.35. We obtain, a bit as before in Theorem 6.38, the following density:
$$d\mu(x)
=-\frac{1}{\pi}\cdot-\frac{1}{2}\sqrt{4x^{-1}-1}\,dx
=\frac{1}{2\pi}\sqrt{4x^{-1}-1}\,dx$$

Thus, we are led to the conclusion in the statement.
\end{proof}

Regarding now the central binomial coefficients, we have here:

\begin{theorem}
The real probability measure having as moments the central binomial coefficients, $D_k=\binom{2k}{k}$, is the measure
$$\alpha_1=\frac{1}{\pi\sqrt{x(4-x)}}\,dx$$
called arcsine law on $[0,4]$.
\end{theorem}

\begin{proof}
We have the following computation, again using formulae from chapter 2:
\begin{eqnarray*}
G(\xi)
&=&\xi^{-1}\sum_{k=0}^\infty D_k\xi^{-k}\\
&=&\frac{1}{\xi}\cdot\frac{1}{\sqrt{1-4/\xi}}\\
&=&\frac{1}{\sqrt{\xi(\xi-4)}} 
\end{eqnarray*}

But this gives the density in the statement, via Theorem 6.35. 
\end{proof}

Finally, we have the following version of the above result:

\begin{theorem}
The real probability measure having as moments the middle binomial coefficients, $E_k=\binom{k}{[k/2]}$, is the following law on $[-2,2]$,
$$\sigma_1=\frac{1}{2\pi}\sqrt{\frac{2+x}{2-x}}\,dx$$
called modified the arcsine law on $[-2,2]$.
\end{theorem}

\begin{proof}
In terms of the central binomial coefficients $D_k$, we have:
$$E_{2k}=\binom{2k}{k}=\frac{(2k)!}{k!k!}=D_k$$
$$E_{2k-1}=\binom{2k-1}{k}=\frac{(2k-1)!}{k!(k-1)!}=\frac{D_k}{2}$$

Standard calculus based on the Taylor formula for $(1+t)^{-1/2}$ gives:
$$\frac{1}{2x}\left(\sqrt{\frac{1+2x}{1-2x}}-1\right)=\sum_{k=0}^\infty E_kx^k$$

With $x=\xi^{-1}$ we obtain the following formula for the Cauchy transform:
\begin{eqnarray*}
G(\xi)
&=&\xi^{-1}\sum_{k=0}^\infty E_k\xi^{-k}\\
&=&\frac{1}{\xi}\left(\sqrt{\frac{1+2/\xi}{1-2/\xi}}-1\right)\\
&=&\frac{1}{\xi}\left(\sqrt{\frac{\xi+2}{\xi-2}}-1\right)
\end{eqnarray*}

By Stieltjes inversion we obtain the density in the statement.
\end{proof}

All this is very nice, and we are obviously building here, as this book goes by, some good knowledge in probability theory. We will be back to all this later.

\section*{6e. Exercises}

Here are a few exercises, in relation with the material from this chapter:

\begin{exercise}
Find a geometric proof for $1+x+x^2+\ldots=1/(1-x)$.
\end{exercise}

\begin{exercise}
Clarify the definition and properties of $x^y$, with $x,y\in\mathbb C$.
\end{exercise}

\begin{exercise}
Complete our proof of the Cauchy formula, in the general case.
\end{exercise}

\begin{exercise}
Learn more complex analysis, including the residue formula.
\end{exercise}

As a bonus exercise, learn some more probability, and random matrices.

\chapter{Fourier analysis}

\section*{7a. Function spaces}

In this chapter we go back to the functions of one real variable, $f:\mathbb R\to\mathbb R$ or $f:\mathbb R\to\mathbb C$, with some applications, by using our complex function technology. We will be mainly interested in constructing the Fourier transform, which is an operation $f\to\widehat{f}$ on such functions, which can solve various questions. Some advertisement first:

\begin{advertisement}
There is only one serious tool in mathematics, and this is the Fourier transform.
\end{advertisement}

Or at least, that's what some people say. And as fans of calculus, we should of course side with them. We will see for instance how, via Fourier transform, light can be decomposed into color components, distant galaxies can be observed, atoms can be distinguished from each other, the future of the universe can be predicted, and so on.

\bigskip

Before doing that, however, let us study the spaces that the functions $f:\mathbb R\to\mathbb C$ can form. These functions can be continuous, differentiable, infinitely differentiable, and so on, but there are many more properties that these functions can have, that we will investigate now. This will lead to various spaces of functions $f:\mathbb R\to\mathbb C$, that can be used, among others, in order to well-define the Fourier transform operation $f\to\widehat{f}$.

\bigskip

Let us start with some well-known and useful inequalities, as follows:

\index{Minkowski inequality}
\index{H\"older inequality}

\begin{theorem}
Given two functions $f,g:\mathbb R\to\mathbb C$, and $p\in[1,\infty]$, we have
$$\left(\int_\mathbb R|f+g|^p\right)^{1/p}\leq\left(\int_\mathbb R|f|^p\right)^{1/p}+\left(\int_\mathbb R|g|^p\right)^{1/p}$$
called Minkowski inequality. Also, if $p,q\in[1,\infty]$ satisfy $1/p+1/q=1$, we have
$$\int_\mathbb R|fg|\leq \left(\int_\mathbb R|f|^p\right)^{1/p}\left(\int_\mathbb R|g|^q\right)^{1/q}$$
called H\"older inequality. The same can be said about  $f,g:X\to\mathbb C$, with $X\subset\mathbb R$.
\end{theorem}

\begin{proof}
All this is very standard, the idea being as follows:

\medskip

(1) As a first observation, at $p=2$, which is a special exponent, we have $q=2$ as well, and the Minkowski and H\"older inequalities are as follows:
$$\left(\int_\mathbb R|f+g|^2\right)^{1/2}\leq\left(\int_\mathbb R|f|^2\right)^{1/2}+\left(\int_\mathbb R|g|^2\right)^{1/2}$$
$$\int_\mathbb R|fg|\leq \left(\int_\mathbb R|f|^2\right)^{1/2}\left(\int_\mathbb R|g|^2\right)^{1/2}$$

But the proof of the H\"older inequality, called Cauchy-Schwarz inequality in this case, is something elementary, coming from the fact that $I(t)=\int_\mathbb R|f+twg|^2$ with $|w|=1$ is a positive degree 2 polynomial in $t\in\mathbb R$, and so its discriminant must be negative. As for the Minkowski inequality, this follows from this, by taking squares and simplifying.

\medskip

(2) In general, let us first prove H\"older, in the case of finite exponents, $p,q\in(1,\infty)$. By linearity we can assume that $f,g$ are normalized, in the following way:
$$\int_\mathbb R|f|^p=\int_\mathbb R|g|^q=1$$

In this case, we want to prove that we have $\int_\mathbb R|fg|\leq1$. And for this purpose, we can use the Young inequality from chapter 3, which gives, for any $x\in\mathbb R$:
$$|f(x)g(x)|\leq\frac{|f(x)|^p}{p}+\frac{|g(x)|^q}{q}$$

By integrating now over $x\in\mathbb R$, we obtain from this, as desired:
$$\int_\mathbb R|fg|
\leq\int_\mathbb R\frac{|f|^p}{p}+\frac{|g|^q}{q}
=\frac{1}{p}+\frac{1}{q}
=1$$

(3) Let us prove now Minkowski, again in the finite exponent case, $p\in(1,\infty)$. We have the following estimate, using the H\"older inequality, and the conjugate exponent:
\begin{eqnarray*}
\int_\mathbb R|f+g|^p
&=&\int_\mathbb R|f+g|\cdot|f+g|^{p-1}\\
&\leq&\int_\mathbb R|f|\cdot|f+g|^{p-1}+\int_\mathbb R|g|\cdot|f+g|^{p-1}\\
&\leq&\left(\int_\mathbb R|f|^p\right)^{1/p}\left(\int_\mathbb R|f+g|^{(p-1)q}\right)^{1/q}\\
&+&\left(\int_\mathbb R|g|^p\right)^{1/p}\left(\int_\mathbb R|f+g|^{(p-1)q}\right)^{1/q}\\
&=&\left[\left(\int_\mathbb R|f|^p\right)^{1/p}+\left(\int_\mathbb R|g|^p\right)^{1/p}\right]\left(\int_\mathbb R|f+g|^p\right)^{1-1/p}
\end{eqnarray*}

Thus, we are led to the Minkowski inequality in the statement.

\medskip

(4) In order to deal now with infinite exponents, observe that we have the following formula, with the essential supremum being computed as the usual supremum, but by neglecting sets of measure 0, which will not contribute to the integrals on the left:
$$\lim_{p\to\infty}\left(\int_\mathbb R|f|^p\right)^{1/p}
=\ {\rm ess\;sup}\,|f|$$

As an example here, a function satisfying $f=0$ at any $x\neq0$ has essential supremum 0, and this no matter what the value $f(0)$ is. Now with the above formula in hand, it makes sense to say that $(\int_\mathbb R|f|^p)^{1/p}$ takes as value at $p=\infty$ the essential supremum of $f$, and with this convention, the Minkowski inequality holds as well at $p=\infty$, trivially:
$${\rm ess\;sup}\,\big|f+g\big|\leq{\rm ess\;sup}\,\big|f\big|+{\rm ess\;sup}\,\big|g\big|$$

The same goes for the H\"older inequality at $p=\infty,q=1$, which is simply:
$$\int_\mathbb R|fg|\leq{\rm ess\;sup}\,\big|f\big|\times\int_\mathbb R|g|$$

And the same goes for the H\"older inequality at $p=1,q=\infty$, which is trivial too.

\medskip

(5) Finally, as mentioned in the statement, all the above extends, in the obvious way, to the case of the functions $f,g:X\to\mathbb C$ with $X\subset\mathbb R$, and with no new proofs needed, because we can simply declare $f,g$ to be zero on $\mathbb R-X$, and we get the results.
\end{proof}

As a consequence of the above result, we can now formulate:

\index{equal almost everywhere}
\index{function space}
\index{p-norm}
\index{normed space}

\begin{theorem}
Given a subset $X\subset\mathbb R$ and an exponent $p\in[1,\infty]$, the following space, with the convention that functions are identified up to equality almost everywhere,
$$L^p(X)=\left\{f:X\to\mathbb C\,\Big|\int_X|f(x)|^pdx<\infty\right\}$$
is a vector space, and the following quantity,
$$||f||_p=\left(\int_X|f(x)|^p\right)^{1/p}$$
is a norm on it, in the sense that it satisfies the usual conditions for a vector space norm. Moreover, $L^p(X)$ is complete with respect to the distance $d(f,g)=||f-g||_p$.
\end{theorem}

\begin{proof}
This basically follows from Theorem 7.2, the idea being as follows:

\medskip

(1) To start with, the vector space norm axioms mentioned above, inspired by the properties satisfied by the length of vectors $||x||=\sqrt{\sum_i|x_i|^2}$, are as follows:

\smallskip

-- $||x||\geq0$, and $||x||=0\iff x=0$.

\smallskip

-- $||\lambda x||=|\lambda|\cdot||x||$, for any $x$, and $\lambda\in\mathbb C$.

\smallskip

-- $||x+y||\leq||x||+||y||$, for any $x,y$.

\medskip

(2) But in our case, the first two norm axioms are clear, once we make the convention for $f=g$ in the statement, and the third axiom comes from Minkowski. 

\medskip

(3) Next, given a normed space as in (1), the function $d(x,y)=||x-y||$ is a distance, in the sense that it satisfies the following properties, a bit as the usual vectors:

\smallskip

-- $d(x,y)\geq0$, and $d(x,y)=0\iff x-y$.

\smallskip

-- $d(x,y)=d(y,x)$, for any $x,y$.

\smallskip

-- $d(x,y)\leq d(x,z)+d(y,z)$, for any $x,y,z$.

\medskip

(4) Thus, it makes sense to talk about complete normed vector spaces, also called Banach spaces, and in our case, $L^p(X)$ is indeed complete, with this being something quite elementary, and exercise here for you, to learn more about all this.
\end{proof}

The above is quite interesting, and as a conclusion to this, let us formulate:

\begin{conclusion}
We can have some geometric understanding of the functions 
$$f:X\to\mathbb C$$
by thinking of them as vectors, having as length $||f||_p$, with $p\in[1,\infty]$ suitably chosen.
\end{conclusion}

Now before exploting this, let us have some further look at the case of the exponent $p=2$, which quite often is the most useful one. We have here, refining Theorem 7.3:

\begin{theorem}
Given $X\subset\mathbb R$, the space of square-summable functions on it,
$$L^2(X)=\left\{f:X\to\mathbb C\,\Big|\int_X|f(x)|^2dx<\infty\right\}$$
with the functions being identified up to equality almost everywhere, has
$$<f,g>=\int_Xf(x)\overline{g(x)}\,dx$$
as scalar product on it, and $||f||=\sqrt{<f,f>}$. Thus, $L^2(X)$ is a Hilbert space.
\end{theorem}

\begin{proof}
As before with Theorem 7.3, this is something quite compact:

\medskip

(1) To start with, we can talk about scalar products on complex vector spaces, with the axioms, inspired from the properties of $<x,y>=\sum_ix_i\bar{y}_i$, being as follows:

\smallskip

-- $<x,y>$ is linear in $x$, and antilinear in $y$.

\smallskip

-- $\overline{<x,y>}=<y,x>$, for any $x,y$.

\smallskip

-- $<x,x>>0$, for any $x\neq0$.

\medskip

(2) But in our case, all these axioms are clear. Also, $||f||=\sqrt{<f,f>}$ is the previous 2-norm, from Theorem 7.3, and finally Hilbert space means complete with respect to the norm, which is something that we already know for $L^2(X)$, from Theorem 7.3.
\end{proof}

Quite nice all this, so let us update our previous conclusion as follows:

\begin{conclusion}[update]
The best functions are the square-summable ones
$$f\in L^2(X)$$
because we can think of them as vectors, subject to scalar products $<f,g>$.
\end{conclusion}

Time now for some exciting geometry, for the functions, based on this? I would say so, but for this purpose, we first need to have a look at what can be done with the Hilbert spaces, in general, and select from there results that we would like to apply to our $L^2$ spaces. And here, among the countless possible choices, all good, we have:

\begin{theorem}
Any Hilbert space $H$ has an orthonormal basis $\{e_i\}_{i\in I}$, which is by definition a set of vectors whose span is dense in $H$, and which satisfy
$$<e_i,e_j>=\delta_{ij}$$
with $\delta$ being a Kronecker symbol. The cardinality $|I|$ of the index set, which can be finite, countable, or uncountable, depends only on $H$, and is called dimension of $H$. We have
$$H\simeq l^2(I)\quad,\quad \sum_i\lambda_ie_i\to(\lambda_i)_{i\in I}$$
with $l^2(I)$ being the space of square-summable sequences $(\lambda_i)_{i\in I}$. The Hilbert spaces with $\dim H=|I|$ countable, such as $l^2(\mathbb N)$, are all isomorphic, and are called separable.
\end{theorem}

\begin{proof}
Many things going on here, the idea being as follows:

\medskip

(1) In finite dimensions, we can turn any vector space basis $\{f_i\}_{i\in I}$ into an orthogonal basis $\{e_i\}_{i\in I}$, by using the Gram-Schmidt procedure, as follows, with $\alpha_i,\beta_i,\gamma_i,\ldots$ being uniquely determined by the fact at each step, $e_k$ must be orthogonal to $f_1,\ldots,f_{k-1}$:
$$e_1=f_1$$
$$e_2=f_2+\alpha_1f_1$$
$$e_3=f_3+\beta_1f_1+\beta_2f_2$$
$$e_4=f_4+\gamma_1f_1+\gamma_2f_2+\gamma_3f_3$$
$$\vdots$$

And then, by replacing $e_i\to e_i/||e_i||$, we have our orthonormal basis, as desired.

\medskip

(2) In general, the same method works, namely Gram-Schmidt, with a subtlety coming from the fact that the basis $\{e_i\}_{i\in I}$ will not span in general the whole $H$, but just a dense subspace of it, as it is in fact obvious by looking at the standard basis of $l^2(\mathbb N)$. 

\medskip

(3) And there is a second subtlety as well, coming from the fact that the recurrence procedure needed for Gram-Schmidt must be replaced by some sort of ``transfinite recurrence'', using standard tools from logic, and more specifically the Zorn lemma.

\medskip

(4) Finally, everything at the end, regarding our notion of separability for the Hilbert spaces,  is clear from definitions, and from our various results above.
\end{proof}

So, can we apply Theorem 7.7 to the spaces of functions $L^2(X)$, and what do we get in this way? As a first famous result here, due to Weierstrass, we have:

\begin{theorem}
The following happen, regarding the functions $f:[a,b]\to\mathbb C$:
\begin{enumerate}
\item Any continuous function $f:[a,b]\to\mathbb C$ can be uniformly approximated by polynomials. Thus, $\{x^n\}_{n\in\mathbb N}$ is an algebraic basis of the space $L^2[a,b]$.

\item By applying Gram-Schmidt we obtain certain polynomials $\{L_n\}_{n\in\mathbb N}$, called Legendre polynomials, which give an explicit isomorphism $L^2[a,b]\simeq l^2(\mathbb N)$.
\end{enumerate}
\end{theorem}

\begin{proof}
We can assume by linearity that $[a,b]=[0,1]$. Let us set:
$$b_{kn}(x)=\binom{n}{k}x^k(1-x)^{n-k}$$

The claim is that $f:[0,1]\to\mathbb C$ is approximated by the following polynomials:
$$f_n(x)=\sum_{k=0}^nf\left(\frac{k}{n}\right)b_{kn}(x)$$

As for the proof, this is something a bit technical, using some probability know-how, based on the fact that $b_{kn}$ are the densities of the binomial laws. And exercise here for you, to learn more about this, and about the Legendre polynomials too.
\end{proof}

As a second application now of Theorem 7.7, which is something famous too, and that we will discuss this time in detail, following Fourier, we have:

\begin{theorem}
The space of $2\pi$-periodic square-summable functions on $\mathbb R$,
$$L^2(\mathbb R)_{per}=\left\{f:\mathbb R\to\mathbb C\,\Big|\,f(t)=f(t+2\pi),\,\int_{-\pi}^\pi|f(t)|^2dt<\infty\right\}$$
has $\{e^{int}\}_{n\in\mathbb Z}$ as orthonormal basis, with respect to the normalized mass $1$ measure.
\end{theorem}

\begin{proof}
This is something quite tricky, which came as a big surprise at the time of its discovery by Fourier, the idea with all this being as follows:

\medskip

(1) To start with, according to what the statement says at the end, the scalar product on $L^2(\mathbb R)_{per}$ is by definition given by the following formula:
$$<f,g>=\frac{1}{2\pi}\int_{-\pi}^\pi f(t)\overline{g(t)}\,dt$$

As for the corresponding norm on $L^2(\mathbb R)_{per}$, this is given by the following formula:
$$||f||=\frac{1}{2\pi}\int_{-\pi}^\pi|f(t)|^2dt$$

Observe that we could have used $[0,2\pi]$ for integrating, or more generally, any interval $I\subset\mathbb R$ having length $2\pi$. For certain technical reasons, we prefer to use $I=[-\pi,\pi]$.

\medskip

(2) Getting now to the proof, as a first basic computation that we can do, coming from the $2\pi$-periodicity of $e^{it}=\cos t+i\sin t$, we have the following formula:
$$<e^{int},e^{imt}>
=\frac{1}{2\pi}\int_{-\pi}^\pi e^{i(n-m)t}dt
=\delta_{nm}$$

Thus $\{e^{int}\}_{n\in\mathbb Z}$ is as orthonormal system, and our theorem says that this is a basis.

\medskip

(3) But, how to prove this? For reasons that will become clear in a moment, consider the following trigonometric polynomials, with $c_k\in\mathbb R$ being chosen for having mass $1$:
$$E_k(t)=c_k\left(\frac{1+\cos t}{2}\right)^k$$

Our claim is that $E_k\to 0$ uniformly on any $[-\pi,-\delta]\cup[\delta,\pi]$, meaning that:
$$\lim_{k\to\infty}\ \sup_{\delta<|t|<\pi}E_k(t)=0$$

(4) So, let us prove this claim. We have the following elementary estimate:
\begin{eqnarray*}
1
&=&\frac{c_k}{\pi}\int_0^\pi\left(\frac{1+\cos t}{2}\right)^kdt\\
&>&\frac{c_k}{\pi}\int_0^\pi\left(\frac{1+\cos t}{2}\right)^k\sin tdt\\
&=&\frac{c_k}{\pi}\left[-\frac{2}{k+1}\left(\frac{1+\cos t}{2}\right)^{k+1}\right]_0^\pi\\
&=&\frac{2c_k}{\pi(k+1)}
\end{eqnarray*}

Now since $E_k$ is decreasing on $[0,\pi]$, we obtain from this, for $\delta<|t|<\pi$:
$$E_k(t)<E_k(\delta)<\frac{\pi(k+1)}{2}\left(\frac{1+\cos\delta}{2}\right)^k$$

But this proves our claim, because for $\delta>0$ we have $(1+\cos\delta)/2<1$, as needed.

\medskip

(5) Getting now to what we wanted to do, we must prove that $\{e^{int}\}_{n\in\mathbb Z}$ spans a dense subset of $L^2(\mathbb R)_{per}$. Since $C(\mathbb R)_{per}\subset L^2(\mathbb R)_{per}$ is dense, it is enough to prove that any $f\in C(\mathbb R)_{per}$ can be approximated by trigonometric polynomials $\sum_nc_ne^{int}$. Moreover, since $||.||_2\leq||.||_\infty$, it is enough to prove our approximation with respect to $||.||_\infty$.

\medskip

(6) All in all, it remains to prove that given a function $f\in C(\mathbb R)_{per}$ and a number $\varepsilon>0$, we can always come with a trigonometric polynomial $\sum_nc_ne^{int}$, such that:
$$\left|f(t)-\sum_nc_ne^{int}\right|<\varepsilon\quad,\quad\forall t\in[-\pi,\pi]$$

But for this, we can use the polynomials $E_k$ from (3). Let us set indeed:
$$Q_k(t)=\frac{1}{2\pi}\int_{-\pi}^\pi f(t-s)E_k(s)\,ds$$

As a first observation, with the change of variables $s\to t-s$, we have the following alternative formula, which shows that $Q_k(t)$ are indeed trigonometric polynomials:
$$Q_k(t)=\frac{1}{2\pi}\int_{-\pi}^\pi f(s)E_k(t-s)\,ds$$

(7) Now given $\varepsilon>0$, let us prove that the estimate in (6) holds indeed, with the trigonometric polynomial there being $Q_k(t)$, for $k>>0$ large enough. For this purpose, we use the uniform continuity of $f$, which tells us that we can find $\delta>0$ such that:
$$|s-t|<\delta\implies|f(s)-f(t)|<\varepsilon$$

Indeed, by using this, we have the following estimate, for the error in (6):
\begin{eqnarray*}
|Q_k(t)-f(t)|
&=&\frac{1}{2\pi}\left|\int_{-\pi}^\pi(f(t-s)-f(t))E_k(s)\,ds\right|\\
&\leq&\frac{1}{2\pi}\int_{-\pi}^\pi|f(t-s)-f(t)|E_k(s)\,ds
\end{eqnarray*}

(8) Now let us split the last integral into three parts, according to:
$$[-\pi,\pi]=[-\pi,-\delta]\cup[-\delta,\delta]\cup[\delta,\pi]$$

On the middle part the integrand is $<\varepsilon$, so the middle integral is $<\varepsilon$. As for the other two integrals, on $[-\pi,-\delta]\cup[\delta,\pi]$, we can use here (4), telling us that $E_k(t)\to0$ uniformly, on that domain. Indeed, with $k>>0$ big enough the other two integrals are $<\varepsilon$ too, so we have obtained (6) as desired, with $\varepsilon\to3\varepsilon$, which finishes the proof. 
\end{proof}

In practice now, as a useful reformulation of Theorem 7.9, we have:

\begin{theorem}
We have an isomorphism $L^2(\mathbb R)_{per}\simeq l^2(\mathbb Z)$, as follows:
\begin{enumerate}
\item Associated to $f\in L^2(\mathbb R)_{per}$ are its Fourier coefficients, given by:
$$\widehat{f}(n)=\frac{1}{2\pi}\int_{-\pi}^\pi f(t)e^{int}dt$$

\item Associated to $g\in l^2(\mathbb Z)$ is the series $S_g(t)=\sum_{n\in\mathbb Z}g(n)e^{-int}$.
\end{enumerate}
\end{theorem}

\begin{proof}
This is something self-explanatory, based on the general orthonormal basis theory from Theorem 7.7, with the Fourier coefficients of $f\in L^2(\mathbb R)_{per}$ being its coefficients $\widehat{f}(n)=<f,e^{-int}>$ with respect to the basis $\{e^{int}\}$ from Theorem 7.9.
\end{proof}

As yet another reformulation of what we have, we have the Parseval formula:

\begin{theorem}
The Fourier coefficients $\widehat{f}(n)=\frac{1}{2\pi}\int_{-\pi}^\pi f(t)e^{int}dt$ satisfy
$$\sum_{n\in\mathbb Z}\widehat{f}(n)\overline{\widehat{g}(n)}=\frac{1}{2\pi}\int_{-\pi}^\pi f(t)\overline{g(t)}\,dt$$
for any $f,g\in L^2(\mathbb R)_{per}$, and in particular satisfy the formula
$$\sum_{n\in\mathbb Z}|\widehat{f}(n)|^2=\frac{1}{2\pi}\int_{-\pi}^\pi|f(t)|^2dt$$
for any $f\in L^2(\mathbb R)_{per}$, with this being called Parseval formula.
\end{theorem}

\begin{proof}
This is indeed yet another reformulation of what we have, coming from the fact that the scalar products and norms are invariant under $L^2(\mathbb R)_{per}\simeq l^2(\mathbb Z)$.
\end{proof}

With this discussed, time perhaps for some applications? And here, there are countless of them, because the above technology can be used in order to decompose various signals, such as mechanical, electromagnetic, seismic or acoustic waves, or even solutions of more complicated differential equations, somewhat of wave type, into sinusoids. 

\bigskip

Thus, many things to be learned here, and for more, have a look at any advanced physics book. In what concerns us, let us just present an application to arithmetic:

\begin{theorem}
We have the following formula of Euler,
$$1+\frac{1}{4}+\frac{1}{9}+\frac{1}{16}+\frac{1}{25}+\ldots=\frac{\pi^2}{6}$$
solving the Basel problem, asking for the computation of this sum.
\end{theorem}

\begin{proof}
This is something quite tricky, the idea being as follows:

\medskip

(1) To start with, we have the following computation of Euler, obtained by factorizing the function $\sin x$, having zeroes at $\mathbb Z/2$, a bit like a polynomial, which gives the result:
\begin{eqnarray*}
\frac{\sin x}{x}
&=&1-\frac{x^2}{3!}+\frac{x^4}{5!}-\frac{x^6}{7!}+\ldots\\
&=&\left(1-\frac{x}{\pi}\right)\left(1+\frac{x}{\pi}\right)\left(1-\frac{x}{2\pi}\right)\left(1+\frac{x}{2\pi}\right)\ldots\\
&=&\left(1-\frac{x^2}{\pi^2}\right)\left(1-\frac{x^2}{4\pi^2}\right)\left(1-\frac{x^2}{9\pi^2}\right)\ldots\\
&=&1-\frac{1}{\pi^2}\sum_{n=1}^\infty\frac{1}{n^2}\,x^2+\ldots
\end{eqnarray*}

Of course this is far from being rigorous, but following Weierstrass, it is possible to fix this, with the above formula for $\sin x/x$ being indeed true. And, exercise of course for you, to learn more about all this, with this being first-class mathematics, for sure.

\medskip

(2) As an alternative approach, much quicker, we can use our Fourier series knowledge. Indeed, the nonzero Fourier coefficients of the function $f(t)=t$ on $[-\pi,\pi]$ are:
$$\widehat{f}(n)
=\frac{1}{2\pi}\int_{-\pi}^\pi te^{int}\,dt
=\frac{1}{2\pi}\left[\frac{1-int}{n^2}\,e^{int}\right]_{-\pi}^\pi
=\frac{(-1)^{n+1}}{n}\,i$$

Thus, the Parseval formula for the function $f(t)=t$ gives:
$$\sum_{n\in\mathbb Z^*}\frac{1}{n^2}=\frac{1}{2\pi}\int_{-\pi}^\pi t^2dt=\frac{1}{2\pi}\cdot\frac{2\pi^3}{3}=\frac{\pi^2}{3}$$

And we therefore solved the Basel problem, just like that. Amazing.
\end{proof}

\section*{7b. Fourier transform}

Going ahead with our study of functions $f:\mathbb R\to\mathbb C$, let us define a key operation on such functions, called convolution, which can be useful for many purposes:

\begin{definition}
The convolution of two functions $f,g:\mathbb R\to\mathbb C$ is the function
$$(f*g)(x)=\int_\mathbb R f(x-y)g(y)dy$$
provided that the function $y\to f(x-y)g(y)$ is indeed integrable, for any $x$.
\end{definition}

There are many reasons for introducing this operation, that we will gradually discover, in what follows. As a basic example, let us take $g=\chi_{[0,1]}$. We have then:
$$(f*g)(x)=\int_0^1f(x-y)dy$$

Thus, with this choice of $g$, the operation $f\to f*g$ has some sort of ``regularizing effect'', that can be useful for many purposes. We will be back to this, later.

\bigskip

Getting now to some theory, let us first try to understand when the convolution operation is well-defined. We have here the following basic result:

\begin{proposition}
The convolution operation is well-defined on the space
$$C_c(\mathbb R)=\left\{f\in C(\mathbb R)\Big|supp(f)=\ {\rm compact}\right\}$$
of continuous functions $f:\mathbb R\to\mathbb C$ having compact support.
\end{proposition}

\begin{proof}
We have several things to be proved, the idea being as follows:

\medskip

(1) First we must show that given two functions $f,g\in C_c(\mathbb R)$, their convolution $f*g$ is well-defined, as a function $f*g:\mathbb R\to\mathbb C$. But this follows from the following estimate, where $l$ denotes the length of the compact subsets of $\mathbb R$:
\begin{eqnarray*}
\int_\mathbb R|f(x-y)g(y)|dy
&=&\int_{supp(g)}|f(x-y)g(y)|dy\\
&\leq&\max(g)\int_{supp(g)}|f(x-y)|dy\\
&\leq&\max(g)\cdot l(supp(g))\cdot\max(f)\\
&<&\infty
\end{eqnarray*}

(2) Next, we must show that the function $f*g:\mathbb R\to\mathbb C$ that we constructed is indeed continuous. But this follows from the following estimate, where $K_f$ is the constant of uniform continuity for the function $f\in C_c(\mathbb R)$:
\begin{eqnarray*}
|(f*g)(x+\varepsilon)-(f*g)(x)|
&=&\left|\int_\mathbb R f(x+\varepsilon-y)g(y)dy-\int_\mathbb R f(x-y)g(y)dy\right|\\
&=&\left|\int_\mathbb R\left(f(x+\varepsilon-y)-f(x-y)\right)g(y)dy\right|\\
&\leq&\int_\mathbb R|f(x+\varepsilon-y)-f(x-y)|\cdot|g(y)|dy\\
&\leq&K_f\cdot\varepsilon\cdot\int_\mathbb R|g|
\end{eqnarray*}

(3) Finally, we must show that the function $f*g\in C(\mathbb R)$ that we constructed has indeed compact support. For this purpose, our claim is that we have:
$$supp(f+g)\subset supp(f)+supp(g)$$

In order to prove this claim, observe that we have, by definition of $f*g$:
$$(f*g)(x)
=\int_\mathbb R f(x-y)g(y)dy
=\int_{supp(g)}f(x-y)g(y)dy$$

But this latter quantity being 0 for $x\notin supp(f)+supp(g)$, this gives the result.
\end{proof}

In relation with derivatives, and with the ``regularizing effect'' of the convolution operation mentioned after Definition 7.13, we have the following result:

\begin{theorem}
Given two functions $f,g\in C_c(\mathbb R)$, assuming that $g$ is differentiable, then so is $f*g$, with derivative given by the following formula:
$$(f*g)'=f*g'$$
More generally, given $f,g\in C_c(\mathbb R)$, and assuming that $g$ is $k$ times differentiable, then so is $f*g$, with $k$-th derivative given by $(f*g)^{(k)}=f*g^{(k)}$.
\end{theorem}

\begin{proof}
In what regards the first assertion, with $y=x-t$, then $t=x-y$, we get:
\begin{eqnarray*}
(f*g)'(x)
&=&\frac{d}{dx}\int_\mathbb R f(x-y)g(y)dy\\
&=&\frac{d}{dx}\int_\mathbb R f(t)g(x-t)dt\\
&=&\int_\mathbb R f(t)g'(x-t)dt\\
&=&\int_\mathbb R f(x-y)g'(y)dy\\
&=&(f*g')(x)
\end{eqnarray*}

As for the second assertion, this follows form the first one, by recurrence.
\end{proof}

Finally, getting beyond the compactly supported continuous functions, we have the following result, which is of particular theoretical importance:

\begin{theorem}
The convolution operation is well-defined on $L^1(\mathbb R)$, and we have:
$$||f*g||_1\leq||f||_1||g||_1$$
Thus, if $f\in L^1(\mathbb R)$ and $g\in C_c^k(\mathbb R)$, then $f*g$ is well-defined, and $f*g\in C_c^k(\mathbb R)$.
\end{theorem}

\begin{proof}
In what regards the first assertion, this follows from:
\begin{eqnarray*}
\int_\mathbb R|(f*g)(x)|dx
&\leq&\int_\mathbb R\int_\mathbb R|f(x-y)g(y)|dydx\\
&=&\int_\mathbb R\int_\mathbb R|f(x-y)g(y)|dxdy\\
&=&\int_\mathbb R|f(x)|dx\int_\mathbb R|g(y)|dy
\end{eqnarray*}

As for the second assertion, this follows from this, and from Theorem 7.15.
\end{proof}

Let us discuss now the construction and main properties of the Fourier transform, which is the main tool in analysis, and mathematics in general. We first have:

\index{Fourier transform}

\begin{definition}
Given $f\in L^1(\mathbb R)$, we define a function $\widehat{f}:\mathbb R\to\mathbb C$ by
$$\widehat{f}(x)=\int_\mathbb R f(t)e^{ixt}dt$$
and call it Fourier transform of $f$.
\end{definition}

As a first observation, even if $f$ is a real function, $\widehat{f}$ is a complex function, which is not necessarily real. Also, $\widehat{f}$ is obviously well-defined, because $f\in L^1(\mathbb R)$ and $|e^{ixt}|=1$. Also, the condition $f\in L^1(\mathbb R)$ is basically needed for construcing $\widehat{f}$, because:
$$\widehat{f}(0)=\int_\mathbb Rf(t)dt$$

Generally speaking, the Fourier transform is there for helping with various computations, with the above formula $\widehat{f}(0)=\int f$ being something quite illustrating. Here are some basic properties of the Fourier transform, all providing some good motivations:

\begin{proposition}
The Fourier transform has the following properties:
\begin{enumerate}
\item Linearity: $\widehat{f+g}=\widehat{f}+\widehat{g}$, $\widehat{\lambda f}=\lambda\widehat{f}$.

\item Regularity: $\widehat{f}$ is continuous and bounded.

\item If $f$ is even then $\widehat{f}$ is even. 

\item If $f$ is odd then $\widehat{f}$ is odd. 
\end{enumerate}
\end{proposition}

\begin{proof}
All this is very standard, and we will leave the proof as an exercise.
\end{proof}

Here are as well some basic computations of Fourier transforms:

\begin{proposition}
We have the following Fourier transform formulae,
$$f=\chi_{[-a,a]}\implies\widehat{f}(x)=\frac{2\sin(ax)}{x}$$
$$f=e^{-at}\chi_{[0,\infty]}(t)\implies\widehat{f}(x)=\frac{1}{a-ix}$$
$$f=e^{at}\chi_{[-\infty,0]}(t)\implies\widehat{f}(x)=\frac{1}{a+ix}$$
$$f=e^{-a|t|}\implies\widehat{f}(x)=\frac{2a}{a^2+x^2}$$
valid for any number $a>0$.
\end{proposition}

\begin{proof}
All this is again standard, and we will leave the proof as an exercise.
\end{proof}

Back now to theory, we have the following result, adding to the various properties in Proposition 7.18, and providing more motivations for the Fourier transform:

\index{exchange of hat}

\begin{proposition}
Given $f,g\in L^1(\mathbb R)$ we have $\widehat{f}g,f\widehat{g}\in L^1(\mathbb R)$ and
$$\int_\mathbb R f(x)\widehat{g}(x)dx=\int_\mathbb R\widehat{f}(x)g(x)dx$$
called ``exchange of hat'' formula.
\end{proposition}

\begin{proof}
Regarding the fact that we have indeed $\widehat{f}g,f\widehat{g}\in L^1(\mathbb R)$, this is actually a bit non-trivial, but we will be back to this later. Assuming this, we have:
$$\int_\mathbb R f(x)\widehat{g}(x)dx=\int_\mathbb R\int_\mathbb R f(x)g(y)e^{ixy}dxdy$$

On the other hand, we have as well the following formula:
$$\int_\mathbb R\widehat{f}(x)g(x)dx=\int_\mathbb R\int_\mathbb R f(y)e^{iyx}g(x)dydx$$

Thus, with $x\leftrightarrow y$, we are led to the formula in the statement.
\end{proof}

As a key result now, showing the power of the Fourier transform, we have:

\begin{theorem}
Given $f:\mathbb R\to\mathbb C$ such that $f,f'\in L^1(\mathbb R)$, we have:
$$\widehat{f'}(x)=-ix\widehat{f}(x)$$
More generally, assuming $f,f',f'',\ldots,f^{(n)}\in L^1(\mathbb R)$, we have
$$\widehat{f^{(k)}}(x)=(-ix)^k\widehat{f}(x)$$
for any $k=1,2,\ldots,n$.
\end{theorem}

\begin{proof}
Assuming that $f:\mathbb R\to\mathbb C$ has compact support, we have indeed:
\begin{eqnarray*}
\widehat{f'}(x)
&=&\int_\mathbb R f'(t)e^{ixt}dt\\
&=&-\int_\mathbb R f(t)\cdot ixe^{ixt}dt\\
&=&-ix\int_\mathbb R f(t)e^{ixt}dt\\
&=&-ix\widehat{f}(x)
\end{eqnarray*}

As for the higher derivatives, the formula here follows by recurrence.
\end{proof}

Importantly, we have a converse statement as well, as follows:

\begin{theorem}
Assuming that $f\in L^1(\mathbb R)$ is such that $F(t)=tf(t)$ belongs to $L^1(\mathbb R)$ too, the function $\widehat{f}$ is differentiable, with derivative given by:
$$(\widehat{f})'(x)=i\widehat{F}(x)$$
More generally, if $F_k(t)=t^kf(t)$ belongs to $L^1(\mathbb R)$, for $k=0,1,\ldots,n$, we have
$$(\widehat{f})^{(k)}(x)=i^k\widehat{F_k}(x)$$
for any $k=1,2,\ldots,n$.
\end{theorem}

\begin{proof}
Regarding the first assertion, the computation here is as follows:
\begin{eqnarray*}
(\widehat{f})'(x)
&=&\frac{d}{dx}\int_\mathbb R f(t)e^{ixt}dt\\
&=&\int_\mathbb R f(t)\cdot ite^{ixt}dt\\
&=&i\int_\mathbb R tf(t)e^{ixt}dt\\
&=&i\widehat{F}(x)
\end{eqnarray*}

As for the second assertion, this follows from the first one, by recurrence.
\end{proof}

Here is another useful result, of the same type, this time regarding convolutions:

\begin{theorem}
Assuming $f,g\in L^1(\mathbb R)$, the following happens:
$$\widehat{f*g}=\widehat{f}\cdot\widehat{g}$$
Conversely, we have $\widehat{fg}=\frac{1}{2\pi}\,\widehat{f}*\widehat{g}$, which holds almost everywhere.
\end{theorem}

\begin{proof}
The first assertion is something elementary, coming as follows:
\begin{eqnarray*}
\widehat{f*g}(x)
&=&\int_\mathbb R (f*g)(t)e^{ixt}dt\\
&=&\int_\mathbb R\int_\mathbb R f(t-s)g(s)e^{ixt}dsdt\\
&=&\int_\mathbb R \left(\int_\mathbb R f(t-s)e^{ix(t-s)}dt\right)g(s)e^{ixs}ds\\
&=&\int_\mathbb R \left(\int_\mathbb R f(r)e^{ixr}dr\right)g(s)e^{ixs}ds\\
&=&\int_\mathbb R \widehat{f}(x)g(s)e^{ixs}ds\\
&=&\widehat{f}(x)\widehat{g}(x)
\end{eqnarray*}

As for the second assertion, this is something more tricky, which follows from the first one by using the Fourier inversion formula, that we will soon learn.
\end{proof}

Let us develop now more theory for the Fourier transform. We first have:

\index{Riemann-Lebesgue property}

\begin{theorem}
Given $f\in L^1(\mathbb R)$, its Fourier transform satisfies
$$\lim_{x\to\pm\infty}\widehat{f}(x)=0$$
called Riemann-Lebesgue property of $\widehat{f}$.
\end{theorem}

\begin{proof}
This is something quite technical, as follows:

\medskip

(1) Given a function $f:\mathbb R\to\mathbb C$ and a number $r\in\mathbb R$, let us set:
$$f_r(t)=f(t-r)$$

Our claim is then is that if $f\in L^p(\mathbb R)$, then the following function is uniformly continuous, with respect to the usual $p$-norm on the right:
$$\mathbb R\to L^p(\mathbb R)\quad,\quad r\to f_r$$

(2) In order to prove this, fix $\varepsilon>0$. Since $f\in L^p(\mathbb R)$, we can find a function of type $g:[-K,K]\to\mathbb C$ which is continuous, such that:
$$||f-g||_p<\varepsilon$$

Now since $g$ is uniformly continuous, we can find $\delta\in(0,K)$ such that:
$$|u-v|<\delta\implies|g(u)-g(v)|<(3K)^{-1/p}\varepsilon$$

But this shows that we have the following estimate:
\begin{eqnarray*}
||g_r-g_s||_p
&=&\left(\int_\mathbb R\big|g(t-r)-g(t-s)\big|^pdt\right)^{1/p}\\
&<&\left[(3K)^{-1}\varepsilon^p(2k+\delta)\right]^{1/p}\\
&<&\varepsilon
\end{eqnarray*}

By using now the formula $||f||_p=||f_r||_p$, which is clear, we obtain:
\begin{eqnarray*}
||f_r-f_s||_p
&\leq&||f_r-g_r||_p+||g_r-g_s||_p+||g_s-f_s||_p\\
&<&\varepsilon+\varepsilon+\varepsilon\\
&=&3\varepsilon
\end{eqnarray*}

But this being true for any $|r-s|<\delta$, we have proved our claim.

\medskip

(3) Let us prove now the Riemann-Lebesgue property of $\widehat{f}$, as formulated in the statement. By using $e^{\pi i}=-1$, and the change of variables $t\to t-\pi/x$, we have:
\begin{eqnarray*}
\widehat{f}(x)
&=&\int_\mathbb Rf(t)e^{ixt}dt\\
&=&-\int_\mathbb Re^{\pi i}f(t)e^{ixt}dt\\
&=&-\int_\mathbb Rf(t)e^{ix(t+\pi/x)}dt\\
&=&-\int_\mathbb Rf\left(t-\frac{\pi}{x}\right)e^{ixt}dt
\end{eqnarray*}

On the other hand, we have as well the following formula:
$$\widehat{f}(x)=\int_\mathbb Rf(t)e^{ixt}dt$$

Thus by summing, we obtain the following formula:
$$2\widehat{f}(x)=\int_\mathbb R\left(f(t)-f\left(t-\frac{\pi}{x}\right)\right)e^{ixt}dt$$

But this gives the following estimate, with notations from (1):
$$2|\widehat{f}(x)|\leq||f-f_{\pi/x}||_1$$

Since by (1) this goes to $0$ with $x\to\pm\infty$, this gives the result.
\end{proof}

As a main result now, quite remarkably, $f$ can be recaptured from $\widehat{f}$, as follows:

\index{Fourier inversion}

\begin{theorem}
Assuming $f,\widehat{f}\in L^1(\mathbb R)$, we have
$$f(t)=\frac{1}{2\pi}\int_\mathbb R \widehat{f}(x)e^{-itx}dx$$
almost everywhere, called Fourier inversion formula.
\end{theorem}

\begin{proof}
Consider the following function, depending on a parameter $\lambda>0$:
$$\varphi_\lambda(s)=\frac{1}{2\pi}\int_\mathbb Re^{-isx-\lambda|x|}dx$$

We have then the following standard convolution computation:
\begin{eqnarray*}
(f*\varphi_\lambda)(t)
&=&\int_\mathbb Rf(t-s)\varphi_\lambda(s)ds\\
&=&\frac{1}{2\pi}\int_\mathbb R\int_\mathbb Rf(t-s)e^{-isx-\lambda|x|}dxds\\
&=&\frac{1}{2\pi}\int_\mathbb R\left(\int_\mathbb Rf(t-s)e^{-isx}ds\right)e^{-\lambda|x|}dx\\
&=&\frac{1}{2\pi}\int_\mathbb R\widehat{f}(x)e^{-itx}e^{-\lambda|x|}dx
\end{eqnarray*}

By letting now $\lambda\to0$, we obtain from this the following formula:
$$\lim_{\lambda\to0}(f*\varphi_\lambda)(t)=\frac{1}{2\pi}\int_\mathbb R \widehat{f}(x)e^{-itx}dx$$

On the other hand, a direct computation based on Theorem 7.24 gives:
$$\lim_{\lambda\to0}(f*\varphi_\lambda)(t)=f(t)$$

We are therefore led to the conclusion in the statement.
\end{proof}

As a first application, we have now the fix for Theorem 7.23. In general, the Fourier transform can be used a bit like the Fourier series, for dealing with all sorts of differential equations, and exercise of course for you, to learn more about all this.

\bigskip

Also, there are many more things that can be said about Fourier transform, a key result being the Plancherel formula, allowing us to talk about the Fourier transform over the space $L^2(\mathbb R)$. Also, we can talk about the Fourier transform over the space $\mathcal S$ of functions all whose derivatives are rapidly decreasing, called Schwartz space. As before, exercise for you to learn more about all this, from any functional analysis book.

\section*{7c. Groups, extensions} 

We have met Joseph Fourier and his mathematics at least 3 times, recently, namely at the end of chapter 5, when talking about the Fourier matrix $F_N\in M_N(\mathbb T)$, and then twice in this chapter, first in relation with the Fourier series isomorphism, which can be written in a compact form as $L^2(\mathbb T)\simeq l^2(\mathbb Z)$, and then in relation with the Fourier transform isomorphism, which can be basically written as $L^2(\mathbb R)\simeq L^2(\mathbb R)$. Which leads us to:

\begin{question}
What is the relation between the various Fourier transform theories that we have? And, are there any more?
\end{question}

In answer, the idea is that we can have a Fourier transform over any locally compact abelian group $G$. In order to discuss this, let us start with some algebra:

\index{abelian group}

\begin{definition}
An abelian group is a set $G$ with a multiplication operation 
$$(g,h)\to gh$$
which must satisfy the following conditions:
\begin{enumerate}
\item Commutativity: we have $gh=hg$, for any $g,h\in G$.

\item Associativity: we have $(gh)k=g(hk)$, for any $g,h,k\in G$.

\item Unit: there is an element $1\in G$ such that $g1=g$, for any $g\in G$.

\item Inverses: for any $g\in G$ there is $g^{-1}\in G$ such that $gg^{-1}=1$.
\end{enumerate}
\end{definition}

There are many examples of abelian groups, and in order to understand what is going on, and to develop our Fourier transform theory, let us first restrict the attention to the finite group case. Here we have as basic examples the cyclic groups $\mathbb Z_N$, formed by the $N$-th roots of unity, and more generally, the products of such cyclic groups:
$$G=\mathbb Z_{N_1}\times\ldots\times\mathbb Z_{N_k}$$

Our goal in what follows will be that of proving that these are in fact all the finite abelian groups, and also that a Fourier transform theory, called ``discrete Fourier transform'', can be developed over such groups. Let us start our study with:

\begin{proposition}
Given a finite abelian group $G$, the group morphisms
$$\chi:G\to\mathbb T$$
with $\mathbb T$ being the unit circle, called characters of $G$, form a finite abelian group $\widehat{G}$.
\end{proposition}

\begin{proof}
Our first claim is that $\widehat{G}$ is a group, with multiplication as follows:
$$(\chi\rho)(g)=\chi(g)\rho(g)$$

Indeed, if $\chi,\rho$ are characters, so is $\chi\rho$, and so the multiplication is well-defined on $\widehat{G}$. Regarding the unit, this is the trivial character $1:G\to\mathbb T$, mapping $g\to1$, for any $g\in G$. Finally, we have inverses, with the inverse of $\chi:G\to\mathbb T$ being its conjugate:
$$\bar{\chi}:G\to\mathbb T\quad,\quad 
g\to\overline{\chi(g)}$$

Thus, we have indeed a group $\widehat{G}$, which is by definition abelian. Our next claim is that $\widehat{G}$ is finite. Indeed, given a group element $g\in G$, we can talk about its order, which is smallest integer $k\in\mathbb N$ such that $g^k=1$. Now given a character $\chi\in\widehat{G}$, we have:
$$\chi(g)^k=1$$

Thus $\chi(g)$ must be one of the $k$-th roots of unity, and in particular there are finitely many choices for $\chi(g)$. Thus, there are finitely many choices for $\chi$, as desired.
\end{proof}

The above construction is quite interesting, and coming next, we have:

\index{cyclic group}
\index{product of cyclic groups}
\index{self-dual group}

\begin{theorem}
The character group operation $G\to\widehat{G}$ for the finite abelian groups, called Pontrjagin duality, has the following properties:
\begin{enumerate}
\item The dual of a cyclic group is the group itself, $\widehat{\mathbb Z}_N=\mathbb Z_N$.

\item The dual of a product is the product of duals, $\widehat{G\times H}=\widehat{G}\times\widehat{H}$.

\item Any product of cyclic groups $G=\mathbb Z_{N_1}\times\ldots\times\mathbb Z_{N_k}$ is self-dual, $G=\widehat{G}$.
\end{enumerate}
\end{theorem}

\begin{proof}
A cyclic group character $\chi:\mathbb Z_N\to\mathbb T$ is uniquely determined by its value $z=\chi(g)$ on the standard generator $g\in\mathbb Z_N$. But this value must satisfy:
$$z^N=1$$

Thus we must have a $N$-th root of unity, $z\in\mathbb Z_N$. Now conversely, any $N$-th root of unity $z\in\mathbb Z_N$ defines a character $\chi:\mathbb Z_N\to\mathbb T$, by setting, for any $r\in\mathbb N$:
$$\chi(g^r)=z^r$$

Thus we have $\widehat{\mathbb Z}_N=\mathbb Z_N$, as claimed in (1). As for (2,3), these are both clear.
\end{proof}

At a more advanced level now, we have the following key result:

\begin{theorem}
The finite abelian groups are the following groups,
$$G=\mathbb Z_{N_1}\times\ldots\times\mathbb Z_{N_k}$$
and these groups are all self-dual, $G=\widehat{G}$.
\end{theorem}

\begin{proof}
This is something quite tricky. Given a finite abelian group $G$, and a prime number $p$, let $G_p\subset G$ be the subset of elements having as order a power of $p$. It is then routine to check that $G_p$ is a subgroup, and that we have an isomorphism as follows:
$$\prod_pG_p\simeq G\quad,\quad (g_p)\to\prod_pg_p$$

Thus, we are left with proving that each component $G_p$ decomposes as a product of cyclic groups, having as orders powers of $p$, as follows:
$$G_p=\mathbb Z_{p^{r_1}}\times\ldots\times\mathbb Z_{p^{r_s}}$$

But this is something that can be proved by recurrence on $|G_p|$, via some routine quotient and lifting arguments, and we are led to the conclusions in the statement.
\end{proof}

In relation now with Fourier analysis, the result is as follows:

\begin{theorem}
Given a finite abelian group $G$, we have an isomorphism as follows, obtained by linearizing/delinearizing the characters,
$$C^*(G)\simeq C(\widehat{G})$$
where $C^*(G)$ is the algebra of functions $\varphi:G\to\mathbb C$, with convolution product, and $C(\widehat{G})$ is the algebra of functions $\varphi:\widehat{G}\to\mathbb C$, with usual product.
\end{theorem}

\begin{proof}
Again, this is something a bit heavy. To start with, we can talk about the algebra $C^*(G)$ formed by the functions $\varphi:G\to\mathbb C$, with convolution product, namely:
$$(\varphi*\psi)(g)=\sum_{h\in G}\varphi(gh^{-1})\psi(h)$$

Now that we know what the statement is about, let us go for the proof. The first observation is that we have a morphism of algebras as follows:
$$C^*(G)\to C(\widehat{G})\quad,\quad g\to\left[\chi\to\chi(g)\right]$$

Now since on both sides we have vector spaces of dimension $N=|G|$, it is enough to check that this morphism is injective. But this follows from Theorem 7.30.
\end{proof}

We can feel that Theorem 7.31 is related to Fourier analysis, and we have:

\begin{fact}
The following happen, regarding the locally compact abelian groups:
\begin{enumerate}
\item What we did in the finite case, namely group characters, and construction and basic properties of the dual, can be extended to them.

\item As basic examples of this, besides what we have in the finite case, and notably $\widehat{\mathbb Z}_N=\mathbb Z_N$, we have $\widehat{\mathbb Z}=\mathbb T$, $\widehat{\mathbb T}=\mathbb Z$, and also $\widehat{\mathbb R}=\mathbb R$.

\item With some care for analytic aspects, $C^*(G)\simeq C(\widehat{G})$ remains true in this setting, and with $G=\mathbb Z,\mathbb T$ and $G=\mathbb R$ we recover the previous Fourier transforms.
\end{enumerate}
\end{fact}

In practice, all this is a bit complicated, and back now to the finite groups, let us work out a softer version of all the above. The result that we need is as follows:

\begin{theorem}
Given a finite abelian group $G$, with dual group $\widehat{G}=\{\chi:G\to\mathbb T\}$, consider the corresponding Fourier coupling, namely:
$$\mathcal F_G:G\times\widehat{G}\to\mathbb T\quad,\quad 
(i,\chi)\to\chi(i)$$
\begin{enumerate}
\item Via the standard isomorphism $G\simeq\widehat{G}$, this Fourier coupling can be regarded as a square matrix, $F_G\in M_G(\mathbb T)$, which is a complex Hadamard matrix.

\item For the cyclic group $G=\mathbb Z_N$ we obtain in this way, via the standard identification $\mathbb Z_N=\{1,\ldots,N\}$, the usual Fourier matrix, $F_N=(w^{ij})$ with $w=e^{2\pi i/N}$.

\item In general, when using a decomposition $G=\mathbb Z_{N_1}\times\ldots\times\mathbb Z_{N_k}$, the corresponding Fourier matrix is given by $F_G=F_{N_1}\otimes\ldots\otimes F_{N_k}$.
\end{enumerate}
\end{theorem}

\begin{proof}
Regarding (1), with the identification $G\simeq\widehat{G}$ made our matrix is given by $(F_G)_{i\chi}=\chi(i)$, and the scalar products between the rows are given by:
$$<R_i,R_j>
=\sum_\chi\chi(i)\overline{\chi(j)}
=\sum_\chi\chi(i-j)
=|G|\cdot\delta_{ij}$$

Thus the rows are pairwise orthogonal, telling us that $F_G$ is Hadamard. Regarding (2), with $\mathbb Z_N=\widehat{\mathbb Z}_N=\{1,\ldots,N\}$, the Fourier coupling is as follows, with $w=e^{2\pi i/N}$:
$$(i,j)\to w^{ij}$$

Finally, (3) follows from $F_{H\times K}=F_H\otimes F_K$, which is something elementary.
\end{proof}

As a nice application of discrete Fourier analysis, generalizing some things that we met before in chapter 5, when first talking about the Fourier matrix $F_N$, we have:

\begin{theorem}
For a matrix $A\in M_N(\mathbb C)$, the following are equivalent:
\begin{enumerate}
\item $A$ is circulant, $A_{ij}=\xi_{j-i}$, for a certain vector $\xi\in\mathbb C^N$.

\item $A$ is Fourier-diagonal, $A=F_NQF_N^*$, for a certain diagonal matrix $Q$.
\end{enumerate}
Moreover, if these conditions hold, then $\xi=F_N^*q$, where $q=(Q_{11},\ldots,Q_{NN})$.
\end{theorem}

\begin{proof}
In one sense, assuming $A_{ij}=\xi_{j-i}$, the matrix $Q=F_N^*AF_N$ is indeed diagonal, as shown by the following computation:
$$Q_{ij}
=\sum_{kl}w^{jl-ik}\xi_{l-k}
=\sum_{kr}w^{j(k+r)-ik}\xi_r
=N\delta_{ij}\sum_rw^{jr}\xi_r$$

Conversely, assuming $Q=diag(q_1,\ldots,q_N)$, the matrix $A=F_NQF_N^*$ is indeed circulant, as shown by the following computation:
$$A_{ij}
=\sum_kw^{ik}Q_{kk}w^{-jk}
=\sum_kw^{(i-j)k}q_k$$

Indeed, since the last term depends only on $j-i$, we have $A_{ij}=\xi_{j-i}$, with $\xi_i=\sum_kw^{-ik}q_k=(F_N^*q)_i$. Thus, we are led to the conclusions in the statement.
\end{proof}

\section*{7d. Independence, limits}

With the above discussed, let us go back now to probability. Our claim is that things are very interesting here, with the real-life notion of independence corresponding to our mathematical notion of convolution, and with the Fourier transform being something very useful in probability, in order to understand the independence. Let us start with:

\index{independence}

\begin{definition}
Two variables $f,g\in L^\infty(X)$ are called independent when
$$E(f^kg^l)=E(f^k)\, E(g^l)$$
happens, for any $k,l\in\mathbb N$.
\end{definition}

This definition hides some non-trivial things. Indeed, by linearity, we would like to have a formula as follows, valid for any polynomials $P,Q\in\mathbb R[X]$:
$$E[P(f)Q(g)]=E[P(f)]\,E[Q(g)]$$

By a continuity argument, it is enough to have this formula for the characteristic functions $\chi_I,\chi_J$ of the arbitrary measurable sets of real numbers $I,J\subset\mathbb R$:
$$E[\chi_I(f)\chi_J(g)]=E[\chi_I(f)]\,E[\chi_J(g)]$$

Thus, we are led to the usual definition of independence, namely:
$$P(f\in I,g\in J)=P(f\in I)\,P(g\in J)$$

All this might seem a bit abstract, but in practice, the idea is of course that $f,g$ must be independent, in an intuitive, real-life sense. As a first result now, we have:

\index{convolution}

\begin{theorem}
Assuming that $f,g\in L^\infty(X)$ are independent, we have
$$\mu_{f+g}=\mu_f*\mu_g$$
where $*$ is the convolution of real probability measures.
\end{theorem}

\begin{proof}
We have the following computation, using the independence of $f,g$:
\begin{eqnarray*}
M_k(f+g)
&=&\sum_r\binom{k}{r}E(f^rg^{k-r})\\
&=&\sum_r\binom{k}{r}M_r(f)M_{k-r}(g)
\end{eqnarray*}

On the other hand, we have as well the following computation:
\begin{eqnarray*}
\int_\mathbb Rx^kd(\mu_f*\mu_g)(x)
&=&\sum_r\binom{k}{r}\int_\mathbb Rx^rd\mu_f(x)\int_\mathbb Ry^{k-r}d\mu_g(y)\\
&=&\sum_r\binom{k}{r}M_r(f)M_{k-r}(g)
\end{eqnarray*}

Thus $\mu_{f+g}$ and $\mu_f*\mu_g$ have the same moments, so they coincide, as claimed.
\end{proof}

Here is now a second result on independence, which is something more advanced:

\index{independence}
\index{Fourier transform}

\begin{theorem}
Assuming that $f,g\in L^\infty(X)$ are independent, we have
$$F_{f+g}=F_fF_g$$
where $F_f(x)=E(e^{ixf})$ is the Fourier transform.
\end{theorem}

\begin{proof}
We have the following computation, using Theorem 7.36:
\begin{eqnarray*}
F_{f+g}(x)
&=&\int_\mathbb Re^{ixz}d\mu_{f+g}(z)\\
&=&\int_\mathbb Re^{ixz}d(\mu_f*\mu_g)(z)\\
&=&\int_{\mathbb R\times\mathbb R}e^{ix(z+t)}d\mu_f(z)d\mu_g(t)\\
&=&\int_\mathbb Re^{ixz}d\mu_f(z)\int_\mathbb Re^{ixt}d\mu_g(t)\\
&=&F_f(x)F_g(x)
\end{eqnarray*}

Thus, we are led to the conclusion in the statement.
\end{proof}

In order to work out some basic illustrations, let us go back to the Poisson laws $p_t$, from chapter 4. We can now say more about them, first with the following result:

\begin{theorem}
The Fourier transform of $p_t$ with $t>0$ is given by:
$$F_{p_t}(y)=\exp\left((e^{iy}-1)t\right)$$
As a consequence, we have convolution semigroup property $p_s*p_t=p_{s+t}$.
\end{theorem}

\begin{proof}
We have indeed the following computation, for the Fourier transform:
\begin{eqnarray*}
F_{p_t}(y)
&=&e^{-t}\sum_k\frac{t^k}{k!}F_{\delta_k}(y)\\
&=&e^{-t}\sum_k\frac{t^k}{k!}\,e^{iky}\\
&=&e^{-t}\sum_k\frac{(e^{iy}t)^k}{k!}\\
&=&\exp(-t)\exp(e^{iy}t)\\
&=&\exp\left((e^{iy}-1)t\right)
\end{eqnarray*}

As for the second assertion, which can be proved of course by a direct computation too, say exercise for you, this follows from this, $\log F_{p_t}$ being linear in $t$.
\end{proof}

Good news, we can now establish the Poisson Limit Theorem, as follows:

\index{PLT}
\index{Poisson Limit Theorem}
\index{Bernoulli laws}

\begin{theorem}[PLT]
We have the following convergence, in moments,
$$\left(\left(1-\frac{t}{n}\right)\delta_0+\frac{t}{n}\,\delta_1\right)^{*n}\to p_t$$
for any $t>0$.
\end{theorem}

\begin{proof}
Let us denote by $\nu_n$ the measure under the convolution sign, namely:
$$\nu_n=\left(1-\frac{t}{n}\right)\delta_0+\frac{t}{n}\,\delta_1$$

We have the following computation, for the Fourier transform of the limit: 
\begin{eqnarray*}
F_{\delta_r}(y)=e^{iry}
&\implies&F_{\nu_n}(y)=\left(1-\frac{t}{n}\right)+\frac{t}{n}\,e^{iy}\\
&\implies&F_{\nu_n^{*n}}(y)=\left(\left(1-\frac{t}{n}\right)+\frac{t}{n}\,e^{iy}\right)^n\\
&\implies&F_{\nu_n^{*n}}(y)=\left(1+\frac{(e^{iy}-1)t}{n}\right)^n\\
&\implies&F(y)=\exp\left((e^{iy}-1)t\right)
\end{eqnarray*}

Thus, we obtain indeed the Fourier transform of $p_t$, as desired. Of course, there are far more things that can be said here, mixing Bernoulli, binomial and Poisson laws, and we will leave some further learning about all this as an interesting exercise.
\end{proof}

Quite nice all this, we have accumulated so far a lot of probability knowledge, going along with our learning of calculus, in chapter 4, chapter 6, and then here. All this certainly needs more discussion, and we will be back to it once we will have the full calculus tools that are needed, in chapter 14, which will be dedicated to probability.

\section*{7e. Exercises}

This was a dense, introductory chapter to Fourier analysis, with several technical things missing, due to some lack of space, and here are some exercises on all this:

\begin{exercise}
Learn more about Gram-Schmidt, Weierstrass and Legendre.
\end{exercise}

\begin{exercise}
Compute the values of $\zeta(s)=\sum_{n\geq1}\frac{1}{n^s}$, at $s\in 2\mathbb N$.
\end{exercise}

\begin{exercise}
Learn about the Plancherel formula, and Fourier over $L^2(\mathbb R)$.
\end{exercise}

\begin{exercise}
Learn also about the Schwartz space $\mathcal S$, and Fourier over it.
\end{exercise}

As bonus exercise, learn some probability. Indeed, as mentioned above, we already know some, but we will take a long break from this, until chapter 14 below.

\chapter{Harmonic functions}

\section*{8a. Laplace operator}

Welcome to theoretical physics. We discuss in this chapter some applications of the theory of real and complex functions developed so far, to some questions from physics, along with some more theory, needed in order to reach to these applications. We will be mainly interested in the two fundamental equations of physics, namely the wave equation, and the heat equation. We will have a look as well into basic quantum mechanics.

\bigskip

All this sounds exciting, doesn't it. However, in practice, technically speaking, we are not exactly ready yet for all these things. But the temptation remains high to talk about them, a bit like physicists do. So, what to do. Ask the cat of course, who says:

\begin{cat}
The more you know, the better that is.
\end{cat}

Thanks cat. So, we will talk indeed about all this, a bit like physicists do, by taking some freedom with surfing on difficult mathematics. In short, this will be a physics chapter, and for a more rigorous treatment of all this, which will be more general too, do not worry, that will come, later in this book, once we will know more calculus.

\bigskip

Getting started, we talked about physics and equations in chapter 3, and time to review that material. In relation with our first equation there, 1D gravity, we have:

\begin{theorem}
The equation of a gravitational free fall in $1$ dimension,
$$\ddot{x}=-\frac{k}{x^2}$$
can be successfully studied, by computing $t=t(x)$ instead of $x=x(t)$, and we have
$$t=\sqrt{\frac{x_0^3}{2k}}\left(\sqrt{\frac{x}{x_0}\left(1-\frac{x}{x_0}\right)}+\arccos\sqrt{\frac{x}{x_0}}\right)$$
with $x_0$ being the initial position, at launching.
\end{theorem}

\begin{proof}
Many things can be said here, the idea being as follows:

\medskip

(1) To start with, the equation in the statement, $\ddot{x}=-k/x^2$, is not really solvable. In order to say however something on the subject, we can trick as in the statement, by computing $t=t(x)$ instead of $x=x(t)$. Now in order to do the inversion, we will need the following standard computation, coming from the chain rule, applied twice:
\begin{eqnarray*}
f(g(x))=x
&\implies&f'(g(x))g'(x)=1\\
&\implies&f'(g(x))=\frac{1}{g'(x)}\\
&\implies&f''(g(x))g'(x)=-\frac{g''(x)}{g'(x)^2}\\
&\implies&f''(g(x))=-\frac{g''(x)}{g'(x)^3}
\end{eqnarray*}

(2) So, consider our equation, written as follows, with $f(t)$ being the position:
$$f''(t)=-\frac{k}{f(t)^2}$$

When setting $t=g(x)$, with $f(g(x))=x$ as above, our equation becomes:
\begin{eqnarray*}
-\frac{g''(x)}{g'(x)^3}=-\frac{k}{x^2}
&\implies&\left(\frac{1}{g'(x)^2}\right)'=\left(\frac{2k}{x}\right)'\\
&\implies&\frac{1}{g'(x)^2}=\frac{2k}{x}+c\\
&\implies&g'(x)=-\frac{1}{\sqrt{2k/x+c}}
\end{eqnarray*}

(3) Next, at the initial position $x_0$ the time must vanish,  $g(x_0)=0$, and we can expect its derivative to explode there, $g'(x_0)=-\infty$. Thus the constant $c$ appearing in the above must be $c=-2k/x_0$, and our equation takes the following form:
$$g'(x)=-\frac{1}{\sqrt{2k}}\cdot\frac{1}{\sqrt{1/x-1/x_0}}$$

Even better, with the change of variables $x=x_0y$, this equation becomes:
$$g'(x_0y)=-\frac{x_0}{\sqrt{2k}}\cdot\frac{1}{\sqrt{1/(x_0y)-1/x_0}}=-\sqrt{\frac{x_0^3}{2k}}\cdot\frac{1}{\sqrt{1/y-1}}$$

(4) Now in order to solve this latter equation, observe that we have:
\begin{eqnarray*}
\left(\sqrt{y-y^2}+\arccos\sqrt{y}\right)'
&=&\frac{1-2y}{2\sqrt{y-y^2}}-\frac{1}{2\sqrt{y}}\cdot\frac{1}{\sqrt{1-y}}\\
&=&-\frac{y}{\sqrt{y-y^2}}\\
&=&-\frac{1}{\sqrt{1/y-1}}
\end{eqnarray*}

Thus, the solution of our equation in (2), with initial data $g(x_0)=0$, is as follows:
$$g(x_0y)=\sqrt{\frac{x_0^3}{2k}}\left(\sqrt{y-y^2}+\arccos\sqrt{y}\right)$$

Now with $x=x_0y$ we are led to the formula in the statement, namely:
$$g(x)=\sqrt{\frac{x_0^3}{2k}}\left(\sqrt{\frac{x}{x_0}\left(1-\frac{x}{x_0}\right)}+\arccos\sqrt{\frac{x}{x_0}}\right)$$

(5) What is next? Many possible things, based on this, and as a first application, we can compute the stopping time as function of the initial position $x_0$, as follows:
$$t_{final}=g(0)=\sqrt{\frac{x_0^3}{2k}}\cdot\arccos(0)=\pi\sqrt{\frac{x_0^3}{8k}}$$

And we will leave some further exploration here as an exercise. Among others, with the above formula of $g=f^{-1}$ in hand, you can do some numerics for the power series expansion of $f$, with the conclusion that this leads nowhere, as stated in (1).
\end{proof}

Good to have this basic gravity question solved. Coming next, we have 1D waves:

\index{wave equation}

\begin{theorem}
The solution of the 1D wave equation $\ddot{f}=v^2f''$ with initial value conditions $f(x,0)=g(x)$ and $\dot{f}(x,0)=h(x)$ is given by the d'Alembert formula:
$$f(x,t)=\frac{g(x-vt)+g(x+vt)}{2}+\frac{1}{2v}\int_{x-vt}^{x+vt}h(s)ds$$
Moreover, in the context of our previous lattice model discretizations, what happens is more or less that the above d'Alembert integral gets computed via Riemann sums.
\end{theorem}

\begin{proof}
This is something that we know from chapter 4, reproduced here for convenience. And with the comment that, as explained in chapter 4, the uniqueness is still to be proved. However, although we have now some good tools for studying such questions, such as the Fourier transform, as apprentice physicists, we will not bother with this. We will come back to this issue later, in chapter 10, when doing mathematics.
\end{proof}

And finally, still in 1D, we can talk as well about heat diffusion, as follows:

\begin{theorem}
The heat diffusion equation, $\dot{f}=\alpha f''$ with $\alpha>0$, with initial condition $f(x,0)=g(x)$, has as solution the function
$$f(x,t)=\int_\mathbb RK_t(x-y)g(y)dy$$
where the function $K_t:\mathbb R\to\mathbb R$, called heat kernel, given by
$$K_t(x)=\frac{1}{\sqrt{4\pi\alpha t}}\,e^{-x^2/4\alpha t}$$
is the standard solution, coming from the initial data $g=\delta_0$, Dirac mass at $0$.
\end{theorem}

\begin{proof}
There are several things going on here, the idea being as follows:

\medskip

(1) Let us first check that $K_t$ is indeed a solution. The time derivative is:
$$\dot{K_t}=-\frac{1}{2t\sqrt{4\pi\alpha t}}\,e^{-x^2/4\alpha t}+\frac{x^2}{4\alpha t^2\sqrt{4\pi\alpha t}}\,e^{-x^2/4\alpha t}$$

Regarding the first space derivative, this is given by the following formula:
$$K_t'=-\frac{x}{2\alpha t\sqrt{4\pi\alpha t}}\,e^{-x^2/4\alpha t}$$

As for the second space derivative, this is given by the following formula:
$$K_t''=-\frac{1}{2\alpha t\sqrt{4\pi\alpha t}}\,e^{-x^2/4\alpha t}+\frac{x^2}{4\alpha^2t^2\sqrt{4\pi\alpha t}}\,e^{-x^2/4\alpha t}$$

Thus, we can see that the heat equation $\dot{f}=\alpha f''$ is indeed satisfied.

\medskip

(2) Now the point is that when convolving the heat kernel $K_t$ with an arbitrary function $g$, all these computations will perturb well, according to the following formulae:
$$\dot{f}(x,t)=\int_\mathbb R\dot{K}_t(x-y)g(y)dy$$
$$f''(x,t)=\int_\mathbb RK''_t(x-y)g(y)dy$$

Thus, we can see that the heat equation $\dot{f}=\alpha f''$ is again satisfied. 

\medskip

(3) Which looks suspiciously neat and mathematical, and you might wonder here, are we doing physics as advertised, or what. In answer, no worries, because there is a bug indeed with this, coming from the correspondence $\delta_0\to K_t$ mentioned at the end, saying that $K_t$ comes from the simplest situation, that of a radiating point body placed at $0$.

\medskip

(4) To be more precise, when getting into this, and we will leave the details here for later in this book, we reach to the conclusion that the mysterious $\sqrt{4\pi\alpha}$ factor in the definition of $K_t$ is indeed needed. By the way, all this is related as well to the mysterious $2\pi$ factor appearing in chapter 7, at Fourier inversion. More on this later.

\medskip

(5) Finally, in what regards uniqueness issues, we will follow here our usual physics policy, as for the waves. That is, our solutions look good, and unique, and job for mathematicians to formally prove this. Back to this later, when doing math.
\end{proof}

With this discussed, what is next? Obviously, 2D. And here, leaving aside gravity, which truly requires multivariable calculus, and that we will discuss in chapter 11, we are left with waves and heat. So, here is the question that we would like to solve:

\begin{question}
What are the $2D$ analogues of the wave and heat equations
$$\ddot{f}=v^2f''\quad,\quad 
\dot{f}=\alpha f''$$ 
plus of course, can we solve these $2D$ wave and heat equations, once found?
\end{question}

In answer now, we certainly have techniques, for dealing with this. The above 1D equations were indeed established in chapter 3 by using 1D lattice models, and we can certainly consider similar 2D models, do the computations, and see what we get. 

\bigskip

Before doing this, however, let us have some mathematical thinking at our problem. It looks reasonable to conjecture that the 2D equations are similar to the 1D ones, provided that we have a reasonable definition for $f''$, as a function $\mathbb R^2\to\mathbb R$. Moreover, with the heat discussion from chapter 3 in mind, we also know that this mysterious $f'':\mathbb R^2\to\mathbb R$ function should measure how far is the average of $f(z)$ with $z\simeq x$ from $f(x)$.

\bigskip

So, let us try to differentiate the functions $f:\mathbb R^2\to\mathbb R$. We first have:

\begin{definition}
We can talk about the partial derivatives of $f:\mathbb R^2\to\mathbb R$,
$$\frac{df}{dx}(x,y)=\lim_{t\to0}\frac{f(x+t,y)-f(x,y)}{t}$$
$$\frac{df}{dy}(x,y)=\lim_{t\to0}\frac{f(x,y+t)-f(x,y)}{t}$$
provided that these derivatives exist indeed.
\end{definition}

Observe the relation with the material from chapter 6, where we already met such things, radial derivatives, in the context of the complex functions $f:\mathbb C\to\mathbb C$. However, the theory developed there is not of much use, right now, in the present setting.

\bigskip

Moving now to higher order, as a natural continuation of Definition 8.6, we have:

\index{double derivatives}

\begin{definition}
We can talk about the second partial derivatives of $f:\mathbb R^2\to\mathbb R$,
$$\frac{d^2f}{dx^2}=\frac{d}{dx}\left(\frac{df}{dx}\right)
\quad,\quad \frac{d^2f}{dxdy}=\frac{d}{dx}\left(\frac{df}{dy}\right)$$
$$\frac{d^2f}{dydx}=\frac{d}{dy}\left(\frac{df}{dx}\right)
\quad,\quad \frac{d^2f}{dy^2}=\frac{d}{dy}\left(\frac{df}{dy}\right)$$
provided that these exist indeed. We can talk about higher order derivatives, as well.
\end{definition}

As before with Definition 8.6, many things that can be said here, mathematically speaking, but we will leave this discussion for later, starting with chapter 9 below, when systematically investigating functions of several variables. In fact, the material in this chapter will be, among others, an introduction to that. However, we will need:

\begin{theorem}
The double derivatives satisfy the formula
$$\frac{d^2f}{dxdy}=\frac{d^2f}{dydx}$$
called Clairaut formula, provided that they are continuous.
\end{theorem}

\begin{proof}
This is something very standard, the idea being as follows:

\medskip

(1) As an intuitive proof for this, we know from chapter 7 that the 1D continuous functions can be approximated by polynomials. Now by extrapolating a bit, we can say that the same should happen for 2D functions, and with the approximation taking care of second derivatives too. But with this, done, because for $f(x,y)=x^py^q$ we have:
$$\frac{d^2f}{dxdy}=\frac{d^2f}{dydx}=pqx^{p-1}y^{q-1}$$

(2) Getting now to more standard techniques, in view of a rigorous proof, given a point $(x,y)$, consider the following functions, depending on $h,k\in\mathbb R$ small:
$$u(h,k)=f(x+h,y+k)-f(x+h,y)$$
$$v(h,k)=f(x+h,y+k)-f(x,y+k)$$
$$w(h,k)=f(x+h,y+k)-f(x+h,y)-f(x,y+k)+f(x,y)$$

(3) By the mean value theorem, for $h,k\neq0$ we can find $\alpha,\beta\in\mathbb R$ such that:
\begin{eqnarray*}
w(h,k)
&=&u(h,k)-u(0,k)\\
&=&h\cdot\frac{d}{dx}\,u(\alpha h,k)\\
&=&h\left(\frac{d}{dx}\,f(x+\alpha h,y+k)-\frac{d}{dx}\,f(x+\alpha h,y)\right)\\
&=&hk\cdot\frac{d}{dy}\cdot\frac{d}{dx}\,f(x+\alpha h,y+\beta k)
\end{eqnarray*}

(4) Similarly, again for $h,k\neq0$, we can find $\gamma,\delta\in\mathbb R$ such that:
\begin{eqnarray*}
w(h,k)
&=&v(h,k)-v(h,0)\\
&=&k\cdot\frac{d}{dy}\,v(h,\delta k)\\
&=&k\left(\frac{d}{dy}\,f(x+h,y+\delta k)-\frac{d}{dy}\,f(x,y+\delta k)\right)\\
&=&hk\cdot\frac{d}{dx}\cdot\frac{d}{dy}\,f(x+\gamma h,y+\delta k)
\end{eqnarray*}

(5) Now by dividing everything by $hk\neq0$, we conclude from this that the following equality holds, with the numbers $\alpha,\beta,\gamma,\delta\in\mathbb R$ being found as above:
$$\frac{d}{dy}\cdot\frac{d}{dx}\,f(x+\alpha h,y+\beta k)
=\frac{d}{dx}\cdot\frac{d}{dy}\,f(x+\gamma h,y+\delta k)$$

But with $h,k\to0$ we get from this the Clairaut formula at $(x,y)$, as desired.
\end{proof}

Back now to physics, in analogy with what we saw in chapter 3, in 1D, we have: 

\index{wave equation}
\index{Laplace operator}
\index{lattice model}
\index{Hooke law}
\index{Newton law}

\begin{theorem}
The wave equation in the plane $\mathbb R^2$ is
$$\ddot{f}=v^2\Delta f$$
where dots denote time derivatives, $\Delta$ is the Laplace operator, given by
$$\Delta f=\frac{d^2f}{dx^2}+\frac{d^2f}{dy^2}$$
and $v>0$ is the propagation speed.
\end{theorem}

\begin{proof}
We have already met this equation in chapter 3, in one dimension, and in 2 dimensions the study is similar, by using a lattice model, as follows:

\medskip

(1) In order to understand the propagation of waves in 2 dimensions, we can model the whole space $\mathbb R^2$ as a network of balls, with springs between them, as follows:
$$
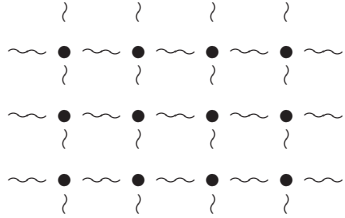
\xymatrix@R=12pt@C=15pt{
&\ar@{~}[d]&\ar@{~}[d]&\ar@{~}[d]&\ar@{~}[d]\\
\ar@{~}[r]&\bullet\ar@{~}[r]\ar@{~}[d]&\bullet\ar@{~}[r]\ar@{~}[d]&\bullet\ar@{~}[r]\ar@{~}[d]&\bullet\ar@{~}[r]\ar@{~}[d]&\\
\ar@{~}[r]&\bullet\ar@{~}[r]\ar@{~}[d]&\bullet\ar@{~}[r]\ar@{~}[d]&\bullet\ar@{~}[r]\ar@{~}[d]&\bullet\ar@{~}[r]\ar@{~}[d]&\\
\ar@{~}[r]&\bullet\ar@{~}[r]\ar@{~}[d]&\bullet\ar@{~}[r]\ar@{~}[d]&\bullet\ar@{~}[r]\ar@{~}[d]&\bullet\ar@{~}[r]\ar@{~}[d]&\\
&&&&&}$$

As before in chapter 3 one dimension, let us send now an impulse, and zoom on one ball. The situation here is as follows, with $l$ being the spring length:
$$
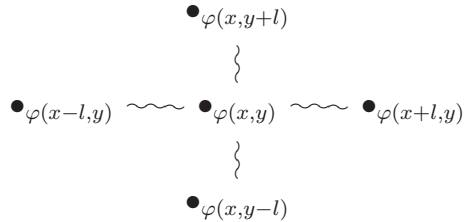
\xymatrix@R=20pt@C=20pt{
&\bullet_{f(x,y+l)}\ar@{~}[d]&\\
\bullet_{f(x-l,y)}\ar@{~}[r]&\bullet_{f(x,y)}\ar@{~}[r]\ar@{~}[d]&\bullet_{f(x+l,y)}\\
&\bullet_{f(x,y-l)}}$$

We have two forces acting at $(x,y)$. First is the Newton motion force, mass times acceleration, which is as follows, with $m$ being the mass of each ball:
$$F_n=m\cdot\ddot{f}(x,y)$$

And second is the Hooke force, displacement of the spring, times spring constant. Since we have four springs at $(x,y)$, this is as follows, $k$ being the spring constant:
\begin{eqnarray*}
F_h
&=&F_h^r-F_h^l+F_h^u-F_h^d\\
&=&k(f(x+l,y)-f(x,y))-k(f(x,y)-f(x-l,y))\\
&+&k(f(x,y+l)-f(x,y))-k(f(x,y)-f(x,y-l))\\
&=&k(f(x+l,y)-2f(x,y)+f(x-l,y))\\
&+&k(f(x,y+l)-2f(x,y)+f(x,y-l))
\end{eqnarray*}

We conclude that the equation of motion, in our model, is as follows:
\begin{eqnarray*}
m\cdot\ddot{f}(x,y)
&=&k(f(x+l,y)-2f(x,y)+f(x-l,y))\\
&+&k(f(x,y+l)-2f(x,y)+f(x,y-l))
\end{eqnarray*}

(2) Now let us take the limit of our model, as to reach to continuum. For this purpose we will assume that our system consists of $B^2>>0$ balls, having a total mass $M$, and spanning a total area $L^2$. Thus, our previous infinitesimal parameters are as follows, with $K$ being the spring constant of the total system, taken to be equal to $k$:
$$m=\frac{M}{B^2}\quad,\quad k=K\quad,\quad l=\frac{L}{B}$$

With these changes, our equation of motion found in (1) reads:
\begin{eqnarray*}
\ddot{f}(x,y)
&=&\frac{KB^2}{M}(f(x+l,y)-2f(x,y)+f(x-l,y))\\
&+&\frac{KB^2}{M}(f(x,y+l)-2f(x,y)+f(x,y-l))
\end{eqnarray*}

Now observe that this equation can be written, more conveniently, as follows:
\begin{eqnarray*}
\ddot{f}(x,y)
&=&\frac{KL^2}{M}\times\frac{f(x+l,y)-2f(x,y)+f(x-l,y)}{l^2}\\
&+&\frac{KL^2}{M}\times\frac{f(x,y+l)-2f(x,y)+f(x,y-l)}{l^2}
\end{eqnarray*}

With $N\to\infty$, and therefore $l\to0$, we obtain in this way:
$$\ddot{f}(x,y)=\frac{KL^2}{M}\left(\frac{d^2f}{dx^2}+\frac{d^2f}{dy^2}\right)(x,y)$$

Thus, we are led to the conclusion in the statement.
\end{proof}

Regarding now heat diffusion, we have here a similar equation, as follows:

\index{heat equation}
\index{Laplace operator}

\begin{theorem}
Heat diffusion in $\mathbb R^2$ is described by the heat equation
$$\dot{f}=\alpha\Delta f$$
where $\alpha>0$ is the thermal diffusivity of the medium, and $\Delta$ is the Laplace operator.
\end{theorem}

\begin{proof}
We have already met this equation in chapter 3, in one dimension, and in 2 dimensions the study is similar, by using a lattice model, as follows:

\medskip

(1) In order to understand the propagation of heat in 2 dimensions, we can model the whole space $\mathbb R^2$ as a network as follows, with all lengths being $l>0$:
$$
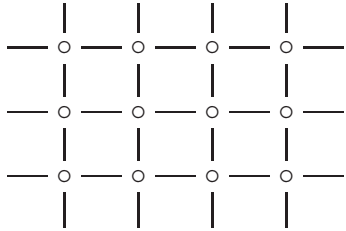
\xymatrix@R=12pt@C=15pt{
&\ar@{-}[d]&\ar@{-}[d]&\ar@{-}[d]&\ar@{-}[d]\\
\ar@{-}[r]&\circ\ar@{-}[r]\ar@{-}[d]&\circ\ar@{-}[r]\ar@{-}[d]&\circ\ar@{-}[r]\ar@{-}[d]&\circ\ar@{-}[r]\ar@{-}[d]&\\
\ar@{-}[r]&\circ\ar@{-}[r]\ar@{-}[d]&\circ\ar@{-}[r]\ar@{-}[d]&\circ\ar@{-}[r]\ar@{-}[d]&\circ\ar@{-}[r]\ar@{-}[d]&\\
\ar@{-}[r]&\circ\ar@{-}[r]\ar@{-}[d]&\circ\ar@{-}[r]\ar@{-}[d]&\circ\ar@{-}[r]\ar@{-}[d]&\circ\ar@{-}[r]\ar@{-}[d]&\\
&&&&&
}$$

We have to implement now the physical heat diffusion mechanism, namely ``the rate of change of the temperature of the material at any given point must be proportional, with proportionality factor $\alpha>0$, to the average difference of temperature between that given point and the surrounding material''. In practice, this leads to a condition as follows, expressing the change of the temperature $f$, over a small period of time $\delta>0$:
$$f(x,y,t+\delta)=f(x,y,t)+\frac{\alpha\delta}{l^2}\sum_{(x,y)\sim(u,v)}\left[f(u,v,t)-f(x,y,t)\right]$$

In fact, we can rewrite our equation as follows, making it clear that we have here an equation regarding the rate of change of temperature at $x$:
$$\frac{f(x,y,t+\delta)-f(x,y,t)}{\delta}=\frac{\alpha}{l^2}\sum_{(x,y)\sim(u,v)}\left[f(u,v,t)-f(x,y,t)\right]$$

(2) Now, let us do the math. In the context of our 2D model the neighbors of $x$ are the points $(x\pm l,y\pm l)$, so the equation above takes the following form:
\begin{eqnarray*}
&&\frac{f(x,y,t+\delta)-f(x,y,t)}{\delta}\\
&=&\frac{\alpha}{l^2}\Big[(f(x+l,y,t)-f(x,y,t))+(f(x-l,y,t)-f(x,y,t))\Big]\\
&+&\frac{\alpha}{l^2}\Big[(f(x,y+l,t)-f(x,y,t))+(f(x,y-l,t)-f(x,y,t))\Big]
\end{eqnarray*}

Now observe that we can write this equation as follows:
\begin{eqnarray*}
\frac{f(x,y,t+\delta)-f(x,y,t)}{\delta}
&=&\alpha\cdot\frac{f(x+l,y,t)-2f(x,y,t)+f(x-l,y,t)}{l^2}\\
&+&\alpha\cdot\frac{f(x,y+l,t)-2f(x,y,t)+f(x,y-l,t)}{l^2}
\end{eqnarray*}

As it was the case before when modeling the wave equation, we recognize on the right the usual approximation of the second derivative, coming from calculus. Thus, when taking the continuous limit of our model, $l\to 0$, we obtain the following equation:
$$\frac{f(x,y,t+\delta)-f(x,y,t)}{\delta}
=\alpha\left(\frac{d^2f}{dx^2}+\frac{d^2f}{dy^2}\right)(x,y,t)$$

Now with $t\to0$, we are led in this way to the heat equation in the statement.
\end{proof}

\section*{8b. Harmonic functions} 

To summarize, we have been doing some physics, with the conclusion that both the wave and heat equations involve the Laplace operator. So, let us formulate:

\index{Laplace operator}
\index{harmonic function}

\begin{definition}
The Laplace operator in $2$ dimensions is:
$$\Delta f=\frac{d^2f}{dx^2}+\frac{d^2f}{dy^2}$$
A function $f:\mathbb R^2\to\mathbb C$ satisfying $\Delta f=0$ will be called harmonic.
\end{definition}

Here the formula of $\Delta$ is the one coming from Theorems 8.9 and 8.10. As for the notion of harmonic function, this is something quite natural. Indeed, we can think of $\Delta$ as being a linear operator on the space of functions $f:\mathbb R^2\to\mathbb C$, and this suggests looking at the eigenvectors of $\Delta$. But the simplest such eigenvectors are those corresponding to the eigenvalue $\lambda=0$, and these are exactly our harmonic functions, satisfying:
$$\Delta f=0$$

Getting now to more concrete things, let us try to find the functions $f:\mathbb R^2\to\mathbb C$ which are harmonic. And here, coming as a good surprise, we have:

\index{complex conjugate}

\begin{theorem}
Any holomorphic function $f:\mathbb C\to\mathbb C$, when regarded as function
$$f:\mathbb R^2\to\mathbb C$$
is harmonic. Moreover, the conjugates $\bar{f}$ of holomorphic functions are harmonic too.
\end{theorem}

\begin{proof}
The first assertion follows from the following computation, for the power functions $f(z)=z^n$, with the usual notation $z=x+iy$:
\begin{eqnarray*}
\Delta z^n
&=&\frac{d^2z^n}{dx^2}+\frac{d^2z^n}{dy^2}\\
&=&\frac{d(nz^{n-1})}{dx}+\frac{d(inz^{n-1})}{dy}\\
&=&n(n-1)z^{n-2}-n(n-1)z^{n-2}\\
&=&0
\end{eqnarray*}

As for the second assertion, this follows from the equality $\Delta\bar{f}=\overline{\Delta f}$.
\end{proof}

As a next goal, in order to understand the harmonic functions, we can try to find the homogeneous polynomials $P\in\mathbb R[x,y]$ which are harmonic. In order to do so, the most convenient is to use the variable $z=x+iy$, and think of these polynomials as being homogeneous polynomials $P\in\mathbb R[z,\bar{z}]$. With this convention, the result is as follows:

\index{homogeneous polynomial}

\begin{theorem}
The degree $n$ homogeneous polynomials $P\in\mathbb R[x,y]$ which are harmonic are precisely the linear combinations of
$$P=z^n\quad,\quad P=\bar{z}^n$$
with the usual identification $z=x+iy$.
\end{theorem}

\begin{proof}
As explained above, any homogeneous polynomial $P\in\mathbb R[x,y]$ can be regarded as an homogeneous polynomial $P\in\mathbb R[z,\bar{z}]$, with the change of variables $z=x+iy$, and in this picture, the degree $n$ homogeneous polynomials are as follows:
$$P(z)=\sum_{k+l=n}c_{kl}z^k\bar{z}^l$$

In oder to solve now the Laplace equation $\Delta P=0$, we must compute the quantities $\Delta(z^k\bar{z}^l)$, for any $k,l$. But the computation here is routine. We first have the following formula, with the derivatives being computed with respect to the variable $x$:
\begin{eqnarray*}
\frac{d(z^k\bar{z}^l)}{dx}
&=&(z^k)'\bar{z}^l+z^k(\bar{z}^l)'\\
&=&kz^{k-1}\bar{z}^l+lz^k\bar{z}^{l-1}
\end{eqnarray*}

By taking one more time the derivative with respect to $x$, we obtain:
\begin{eqnarray*}
\frac{d^2(z^k\bar{z}^l)}{dx^2}
&=&k(z^{k-1}\bar{z}^l)'+l(z^k\bar{z}^{l-1})'\\
&=&k\left[(z^{k-1})'\bar{z}^l+z^{k-1}(\bar{z}^l)'\right]+l\left[(z^k)'\bar{z}^{l-1}+z^k(\bar{z}^{l-1})'\right]\\
&=&k\left[(k-1)z^{k-2}\bar{z}^l+lz^{k-1}\bar{z}^{l-1}\right]+l\left[kz^{k-1}\bar{z}^{l-1}+(l-1)z^k\bar{z}^{l-2}\right]\\
&=&k(k-1)z^{k-2}\bar{z}^l+2klz^{k-1}\bar{z}^{l-1}+l(l-1)z^k\bar{z}^{l-2}
\end{eqnarray*}

With respect to the variable $y$, the computations are similar, but some $\pm i$ factors appear, due to $z'=i$ and $\bar{z}'=-i$, coming from $z=x+iy$. We first have:
\begin{eqnarray*}
\frac{d(z^k\bar{z}^l)}{dy}
&=&(z^k)'\bar{z}^l+z^k(\bar{z}^l)'\\
&=&ikz^{k-1}\bar{z}^l-ilz^k\bar{z}^{l-1}
\end{eqnarray*}

By taking one more time the derivative with respect to $y$, we obtain:
\begin{eqnarray*}
\frac{d^2(z^k\bar{z}^l)}{dy^2}
&=&ik(z^{k-1}\bar{z}^l)'-il(z^k\bar{z}^{l-1})'\\
&=&ik\left[(z^{k-1})'\bar{z}^l+z^{k-1}(\bar{z}^l)'\right]-il\left[(z^k)'\bar{z}^{l-1}+z^k(\bar{z}^{l-1})'\right]\\
&=&ik\left[i(k-1)z^{k-2}\bar{z}^l-ilz^{k-1}\bar{z}^{l-1}\right]-il\left[ikz^{k-1}\bar{z}^{l-1}-i(l-1)z^k\bar{z}^{l-2}\right]\\
&=&-k(k-1)z^{k-2}\bar{z}^l+2klz^{k-1}\bar{z}^{l-1}-l(l-1)z^k\bar{z}^{l-2}
\end{eqnarray*}

We can now sum the formulae that we found, and we obtain:
\begin{eqnarray*}
\Delta(z^k\bar{z}^l)
&=&\frac{d^2(z^k\bar{z}^l)}{dx^2}+\frac{d^2(z^k\bar{z}^l)}{dy^2}\\
&=&k(k-1)z^{k-2}\bar{z}^l+2klz^{k-1}\bar{z}^{l-1}+l(l-1)z^k\bar{z}^{l-2}\\
&&-k(k-1)z^{k-2}\bar{z}^l+2klz^{k-1}\bar{z}^{l-1}-l(l-1)z^k\bar{z}^{l-2}\\
&=&4klz^{k-1}\bar{z}^{l-1}
\end{eqnarray*}

In other words, we have reached to the following interesting formula:
$$f=z^k\bar{z}^l\implies\Delta f=\frac{4klf}{|z|^2}$$

Now let us get back to our homogeneous polynomial $P$, written as follows:
$$P(z)=\sum_{k+l=n}c_{kl}z^k\bar{z}^l$$

By using the above formula, the Laplacian of $P$ is given by:
$$\Delta P(z)=\frac{4}{|z|^2}\sum_{k+l=n}klc_{kl}z^k\bar{z}^l$$

We conclude that the Laplace equation for $P$ takes the following form:
\begin{eqnarray*}
\Delta P=0
&\iff&klc_{kl}=0,\forall k,l\\
&\iff&[k,l\neq0\implies c_{kl}=0]\\
&\iff&P=c_{n0}z^n+c_{0n}\bar{z}^n
\end{eqnarray*}

Thus, we are led to the conclusion in the statement.
\end{proof}

At the general level now, let us formulate the following definition:

\begin{definition}
The Cauchy-Riemann operators are
$$\partial=\frac{1}{2}\left(\frac{d}{dx}-i\frac{d}{dy}\right)
\quad,\quad
\bar{\partial}=\frac{1}{2}\left(\frac{d}{dx}+i\frac{d}{dy}\right)$$
where $\frac{d}{dx}$ and $\frac{d}{dy}$ are the usual partial derivatives for complex functions. 
\end{definition}

There are many things that can be said about the Cauchy-Riemann operators $\partial,\bar{\partial}$, the idea being that in many contexts, these are better to use than the usual partial derivatives $\frac{d}{dx},\frac{d}{dy}$, and with this being a bit like the usage of the variables $z,\bar{z}$, instead of the decomposition $z=x+iy$, for many questions regarding the complex numbers. 

\bigskip

We have already seen in fact some instances of this, in our computations above. At the general level, the main properties of $\partial,\bar{\partial}$ can be summarized as follows:

\begin{theorem}
Assume that $f:X\to\mathbb C$ is differentiable in the real sense.
\begin{enumerate}
\item $f$ is holomorphic precisely when $\bar{\partial}f=0$.

\item In this case, its derivative is $f'=\partial f$.

\item The Laplace operator is given by $\Delta=4\partial\bar{\partial}$.

\item $f$ is harmonic precisely when $\partial\bar{\partial}f=0$.
\end{enumerate}
\end{theorem}

\begin{proof}
We can assume by linearity that we are dealing with differentiability questions at $0$. Since our function $f:X\to\mathbb C$ is differentiable in the real sense, we have a formula as follows, with $z=x+iy$, and with $a,b\in\mathbb C$ being the partial derivatives at $0$:
$$f(z)=ax+by+o(z)$$

Now observe that we can write this formula in the following way:
\begin{eqnarray*}
f(z)
&=&a\cdot\frac{z+\bar{z}}{2}+b\cdot\frac{z-\bar{z}}{2i}+o(z)\\
&=&a\cdot\frac{z+\bar{z}}{2}+b\cdot\frac{i\bar{z}-iz}{2}+o(z)\\
&=&\frac{a-ib}{2}\cdot z+\frac{a+ib}{2}\cdot\bar{z}+o(z)
\end{eqnarray*}

Now by dividing by $z$, we obtain from this the following formula:
\begin{eqnarray*}
\frac{f(z)}{z}
&=&\frac{a-ib}{2}+\frac{a+ib}{2}\cdot\frac{\bar{z}}{z}+o(1)\\
&=&\partial f(0)+\bar{\partial}f(0)\cdot\frac{\bar{z}}{z}+o(1)
\end{eqnarray*}

But this gives the first two assertions, because in order for the derivative $f'(0)$ to exist, appearing as the $z\to0$ limit of the above quantity, the coefficient of $\bar{z}/z$, which does not converge, must vanish. Regarding now the third assertion, this follows from:
\begin{eqnarray*}
\Delta
&=&\frac{d^2}{dx^2}+\frac{d^2}{dy^2}\\
&=&\left(\frac{d}{dx}-i\frac{d}{dy}\right)
\left(\frac{d}{dx}+i\frac{d}{dy}\right)\\
&=&4\partial\bar{\partial}
\end{eqnarray*}

As for the last assertion, this is clear from this latter formula of $\Delta$.
\end{proof}

In analogy now with the theory of the holomorphic functions, we have:

\index{maximum modulus}
\index{mean value formula}
\index{Liouville theorem}

\begin{theorem}
The harmonic functions obey to the same general principles as the holomorphic functions, namely:
\begin{enumerate}
\item The maximum modulus principle.

\item The plain mean value formula. 

\item The boundary mean value formula.

\item The Liouville theorem.
\end{enumerate}
Also, locally, the real harmonic functions are the real parts of holomorphic functions. 
\end{theorem}

\begin{proof}
This is something quite tricky, the idea being as follows:

\medskip

(1) Regarding the plain mean value formula, here the statement is that given an harmonic function $f:X\to\mathbb C$, and a disk $D$, the following happens:
$$f(x)=\int_Df(y)dy$$

In order to prove this, we can assume that $D$ is centered at $0$, of radius $r>0$. If we denote by $\chi_r$ the characteristic function of this disk, normalized as to integrate up to 1, in terms of the convolution operation from chapter 7, we want to prove that we have:
$$f=f*\chi_r$$

For doing so, let us pick a number $0<s<r$, and then a solution $w$ of the following equation on $D$, which can be constructed explicitly, say as a radial function:
$$\Delta w=\chi_r-\chi_s$$

By using the properties of the convolution operation $*$ from chapter 7, we have:
\begin{eqnarray*}
f*\chi_r-f*\chi_s
&=&f*(\chi_r-\chi_s)\\
&=&f*\Delta w\\
&=&\Delta f*w\\
&=&0
\end{eqnarray*}

Thus $f*\chi_r=f*\chi_s$, and by letting now $s\to0$, we get $f*\chi_r=f$, as desired.

\medskip

(2) Regarding the boundary mean value formula, here the statement is that given an harmonic function $f:X\to\mathbb C$, and a disk $D$, with boundary $\gamma$, the following happens:
$$f(x)=\int_\gamma f(y)dy$$

But this follows as a consequence of the plain mean value formula in (1), with our two mean value formulae, the one there and the one here, being in fact equivalent, by using annuli and radial integration for the proof of the equivalence, in the obvious way. 

\medskip

(3) Regarding the maximum modulus principle, the statement here is that any holomorphic function $f:X\to\mathbb C$ has the property that the maximum of $|f|$ over a domain is attained on its boundary. That is, given a domain $D$, with boundary $\gamma$, we have:
$$\exists x\in\gamma\quad,\quad |f(x)|=\max_{y\in D}|f(y)|$$

But this is something which follows again from the mean value formula in (1), first for the disks, and then in general, by using a standard division argument.

\medskip

(4) Regarding the Liouville theorem, the statement here is that an entire, bounded harmonic function must be constant: 
$$f:\mathbb R^2\to\mathbb C\quad,\quad\Delta f=0\quad,\quad |f|\leq M\quad\implies\quad f={\rm constant}$$

As a slightly weaker statement, again called Liouville theorem, we have the fact that an entire harmonic function which vanishes at $\infty$ must vanish globally: 
$$f:\mathbb R^2\to\mathbb C\quad,\quad\Delta f=0\quad,\quad\lim_{x\to\infty}f(x)=0\quad\implies\quad f=0$$

But we can view these as a consequence of the mean value formula in (1), because given two points $x\neq y$, we can view the values of $f$ at these points as averages over big disks centered at these points, say $B=D_x(R)$ and $C=D_y(R)$, with $R>>0$:
$$f(x)=\int_Bf(z)dz\quad,\quad f(y)=\int_Cf(z)dz$$

Indeed, the point is that when the radius goes to $\infty$, these averages tend to be equal, and so we have $f(x)\simeq f(y)$, which gives $f(x)=f(y)$ in the limit, as desired.

\medskip

(5) Finally, the last assertion, which certainly agrees with our findings so far, is something more subtle, the claim here being that a real harmonic function on the unit disk $f:D\to\mathbb R$ must equal the real part of the following function, which is holomorphic:
$$g(z)=\frac{1}{2\pi}\int_{-\pi}^\pi\frac{e^{it}+z}{e^{it}-z}\,f(e^{it})dt$$

Indeed, observe that this holds for $f(z)=Re(z^n)$, and by taking linear combinations, it holds for the real parts of holomorphic functions. However, the proof in general is non-trivial, using techniques and estimates similar to those from chapter 7, for the Fourier series, and we refer to Rudin \cite{ru2} for this. And with the promise that we will be back to this, harmonic functions, later in this book, once we will know more calculus.
\end{proof}

Finally, no discussion about harmonic functions would be complete without a word on radial harmonics, meaning harmonic functions of the following special type: 
$$f(z)=\varphi(|z|)$$

With some physics in mind, namely waves or heat radiating from the origin, it is reasonable to call such functions fundamental solutions of the Laplace equation $\Delta f=0$, and also to allow for singularities at 0, that is, to only assume $\varphi:(0,\infty)\to\mathbb C$. And here, surprise, here is what we get, messing up everything that we know:

\index{Laplace equation}
\index{radial function}

\begin{theorem}
We can talk about harmonic functions in $N$ dimensions, in the obvious way, and the fundamental radial solutions of the Laplace equation $\Delta f=0$ are
$$f(x)=\begin{cases}
||x||^{2-N}&(N\neq 2)\\
\log||x||&(N=2)
\end{cases}$$
with the $\log$ at $N=2$ basically coming from $\log'=1/x$. Thus, we have a blowup phenomenon at the dimension value $N=2$, which is therefore the bad dimension.
\end{theorem}

\begin{proof}
This is something quite tricky, the idea being as follows:

\medskip

(1) We can certainly talk about the Laplace operator and harmonic functions in $N$ dimensions, in the obvious way, and pretty much everything that we know at $N=2$ extends well, with this including the wave and heat equations, and the bulk of Theorem 8.16, meaning everything, except for the last assertion. Which is certainly very nice, and we will come back to this, with details, in chapter 12, when talking $N$ dimensions.

\medskip

(2) In relation with radial harmonics, again this is something that we will discuss in chapter 12, but the computations are fairly straightforward. To be more precise, with $f(x)=\varphi(r)$, and $r=||x||$, the Laplace equation $\Delta f=0$ reformulates as follows:
$$r\varphi''+(N-1)\varphi'=0$$

But the fundamental solutions of this latter equation are easy to find, as follows:
$$\varphi(r)=\begin{cases}
r^{2-N}&(N\neq 2)\\
\log r&(N=2)
\end{cases}$$

Thus, we are led to the above conclusions. And more on this later, in chapter 12.
\end{proof}

\section*{8c. Light, spectroscopy} 

Getting back now to physics, a natural question would be that of going back to the equations that we started with, namely the wave and heat equations in 2D, and see if our accumulated knowledge about harmonic functions, which is after all not that bad, can help there. Unfortunately, this is not exactly the case, and we are still a long way to go, from solving that equations. So, obviously, time to ask the cat. And cat says:

\begin{cat}
If you don't understand physics, do more physics.
\end{cat}

Thanks cat, and this seems wise indeed, let's have some fun with more physics, still related of course to the Laplace operator $\Delta$, and for difficult mathematics, which is still to be solved, we will see later, towards the end of the present book. Getting started now, we will need for what we want to talk about something quite heavy, namely:

\index{Maxwell equations}

\begin{theorem}
Light is an electromagnetic wave, subject to the wave equation
$$\ddot{f}=v^2\Delta f$$
traveling in vacuum at speed $v=c$, numerically given by
$$c=299,792,458$$
with this figure being exact, and in non-vacuum at a lower speed $v<c$.
\end{theorem}

\begin{proof}
This is something quite tricky, involving a lot of physics, as follows:

\medskip

(1) To start with, moving charges produce magnetism. Thus, in contrast with classical mechanics, where static and dynamic problems are described by a unique field, the gravitational one, in electrodynamics we have two fields, namely the electric field $E$, and the magnetic field $B$. And these are subject to the Maxwell equations, namely:
$$<\nabla,E>=\frac{\rho}{\varepsilon_0}\quad,\quad<\nabla,B>=0$$
$$\nabla\times E=-\dot{B}\quad,\quad\nabla\times B=\mu_0J+\mu_0\varepsilon_0\dot{E}$$

(2) To be more precise, regarding the math, the dots denote derivatives with respect to time, $<\,,>$ are usual scalar products, and $\nabla$ and $\times$ are given by:
$$\nabla=\begin{pmatrix}
\frac{d}{dx}\\
\frac{d}{dy}\\
\frac{d}{dz}
\end{pmatrix}\qquad,\qquad 
\begin{pmatrix}x_1\\ x_2\\x_3\end{pmatrix}
\times\begin{pmatrix}y_1\\ y_2\\y_3\end{pmatrix}
=\begin{pmatrix}x_2y_3-x_3y_2\\ x_3y_1-x_1y_3\\x_1y_2-x_2y_1\end{pmatrix}$$

As the the physics, $\rho$ is the charge, $J$ is the volume current density, and the constants $\mu_0,\varepsilon_0$, are known to be related to the speed of light by the following formula:
$$\mu_0\varepsilon_0=\frac{1}{c^2}$$ 

\medskip

(3) The point now is that, in regions of space where there is no charge or current present, the Maxwell equations take a quite simple form, as follows: 
$$<\nabla,E>=<\nabla,B>=0$$
$$\nabla\times E=-\dot{B}\quad,\quad 
\nabla\times B=\mu_0\varepsilon_0\dot{E}$$

But this shows, via some math, that both $E$ and $B$ satisfy the wave equation:
$$\ddot{f}=c^2\Delta f$$

(4) Next, and again coming from Maxwell, the point is that accelerating or decelerating charges produce electromagnetic waves. These are the waves predicted by (3), traveling at speed $c$, and with the important extra property that they depend on a real positive parameter, that can be called, upon taste, frequency, wavelength, or color.

\medskip

(5) In practice, this phenomenon can be observed is a variety of situations, such as the usual light bulbs, where electrons get decelerated by the filament, acting as a resistor, or in usual fire, which is a chemical reaction, with the electrons moving around, as they do in any chemical reaction, or in more complicated machinery like nuclear plants, particle accelerators, and so on, leading there to all sorts of eerie glows, of various colors. 

\medskip

(6) So, this was for the story of light, and more on this later, in chapter 15. In the meantime, for simplifying, we can say that light being certainly not a particle, and this because it does not hurt your eyes, when reaching them, it must be a wave, as stated.
\end{proof}

Next, in order to understand light, we will need a piece of math, as follows:

\begin{theorem}
The 1D wave equation with speed $v$, namely
$$\ddot{f}=v^2f''$$
has as basic solutions the following functions,
$$f(x,t)=A\cos(kx-wt+\delta)$$
with $A$ being called amplitude, $kx-wt+\delta$ being called the phase, $k$ being the wave number, $w$ being the angular frequency, and $\delta$ being the phase constant. We have
$$\lambda=\frac{2\pi}{k}\quad,\quad T=\frac{2\pi}{kv}\quad,\quad \nu=\frac{1}{T}\quad,\quad w=2\pi\nu$$
relating the wavelength $\lambda$, period $T$, frequency $\nu$, and angular frequency $w$. Moreover, any solution of the wave equation appears as a linear combination of such basic solutions.
\end{theorem}

\begin{proof}
There are several things going on here, the idea being as follows:

\medskip

(1) Our first claim is that the function $f$ in the statement satisfies indeed the wave equation, with speed $v=w/k$. But this is clear, coming from:
$$\ddot{f}=-w^2f\quad,\quad f''=-k^2f$$

(2) Regarding now the other things in the statement, all this is basically terminology, which is very natural, when thinking how $f(x,t)=A\cos(kx-wt+\delta)$ propagates.

\medskip

(3) Finally, the last assertion is standard, coming either from d'Alembert, or from Fourier analysis, and we will leave clarifying this as an instructive exercise.
\end{proof}

Moving ahead now towards electromagnetism and 3D, let us formulate:

\begin{definition}
A monochromatic plane wave is a solution of the 3D wave equation which moves in only 1 direction, making it in practice a solution of the 1D wave equation, and which is of the special from found in Theorem 8.20, with no frequencies mixed.
\end{definition}

In other words, we are making here two assumptions on our wave. First is the 1-dimensionality assumption, which gets us into the framework of Theorem 8.20. And second is the assumption, in connection with the Fourier decomposition result from the end of Theorem 8.20, that our solution is of pure type, meaning a wave having a well-defined wavelenght and frequency, instead of being a packet of such pure waves.

\bigskip

Summarizing, we have now a decent intuition about what light is, and more on this later. Let us discuss now the examples. The idea is that we have various types of light, depending on frequency and wavelength. These are normally referred to as electromagnetic waves, but for keeping things simple and luminous, we will keep using the familiar term ``light''. The classification, in a rough form, is as follows:
$$\begin{matrix}
{\rm Frequency}&\ \ \ \ \ \ \ \ \ \ {\rm Type}\ \ \ \ \ \ \ \ \ \ &{\rm Wavelength}\\
&-\\
10^{18}-10^{20}&\gamma\ {\rm rays}&10^{-12}-10^{-10}\\
10^{16}-10^{18}&{\rm X-rays}&10^{-10}-10^{-8}\\
10^{15}-10^{16}&{\rm UV}&10^{-8}-10^{-7}\\
&-\\
10^{14}-10^{15}&{\rm blue}&10^{-7}-10^{-6}\\
10^{14}-10^{15}&{\rm yellow}&10^{-7}-10^{-6}\\
10^{14}-10^{15}&{\rm red}&10^{-7}-10^{-6}\\
&-\\
10^{11}-10^{14}&{\rm IR}&10^{-6}-10^{-3}\\
10^9-10^{11}&{\rm microwave}&10^{-3}-10^{-1}\\
1-10^9&{\rm radio}&10^{-1}-10^8\\
\end{matrix}$$

Observe the tiny space occupied by the visible light, all colors there, and the many more missing, being squeezed under the $10^{14}-10^{15}$ frequency banner. Here is a zoom on that part, with of course the comment that all this, colors, is something subjective:
$$\begin{matrix}
{\rm Frequency\ THz}=10^{12}\ {\rm Hz}&\ \ \ {\rm Color}\ \ \ &{\rm Wavelength\ nm}=10^{-9}\ {\rm m}\\
&-\\
670-790&{\rm violet}&380-450\\
620-670&{\rm blue}&450-485\\
600-620&{\rm cyan}&485-500\\
530-600&{\rm green}&500-565\\
510-530&{\rm yellow}&565-590\\
480-510&{\rm orange}&590-625\\
400-480&{\rm red}&625-750
\end{matrix}$$

With this in hand, we can now do some basic optics. Light usually comes in bundles, with waves of several wavelenghts coming at the same time, from the same source, and the first challenge is that of separating these wavelenghts. In order to discuss this, from a practical perspective, let us start with the following fact:

\begin{fact}
When traveling through a material, and hitting a new material, some of the light gets reflected, at the same angle, and some of it gets refracted, at a different angle, depending both on the old and the new material, and on the wavelength.
\end{fact}

Again, this is something deep, and there are many things that can be said here, ranging from elementary to fairly advanced. As a basic formula, we have the famous Snell law, which relates the incidence angle $\theta_1$ to the refraction angle $\theta_2$, as follows:
$$\frac{\sin\theta_2}{\sin\theta_1}=\frac{n_1(\lambda)}{n_2(\lambda)}$$

Here $n_i(\lambda)$ are the refraction indices of the two materials, adjusted for the wavelength, and with this adjustment for wavelength being the whole point, which is something quite complicated. Now as a simple consequence of the above, we have:

\begin{theorem}
Light can be decomposed, by using a prism.
\end{theorem}

\begin{proof}
This follows from Fact 8.22. Indeed, when hitting a piece of glass, provided that the hitting angle is not $90^\circ$, the light will decompose over the wavelenghts present, with the corresponding refraction angles depending on these wavelengths. And we can capture these split components at the exit from the piece of glass, again deviated a bit, provided that the exit surface is not parallel to the entry surface. And the simplest device doing the job, that is, having two non-parallel faces, is a prism.
\end{proof}

As an application of this, we can study events via spectroscopy, by capturing the light the event has produced, decomposing it with a prism, carefully recording its spectral signature, consisting of the wavelenghts present, and their density, and then doing some reverse engineering, consisting in reconstructing the event out of its spectral signature. 

\section*{8d. Atomic spectrum}

Getting now to some truly exciting applications of light and spectroscopy, let us discuss the beginnings of the atomic theory. There is a long story here, involving many discoveries, around 1890-1900, focusing on hydrogen ${\rm H}$. We will present here things a bit retrospectively. First on our list is the following discovery, by Lyman in 1906:

\begin{fact}[Lyman]
The hydrogen atom has spectral lines given by the formula
$$\frac{1}{\lambda}=R\left(1-\frac{1}{n^2}\right)$$
where $R\simeq 1.097\times 10^7$ and $n\geq2$, which are as follows,
$$\begin{matrix}
n&{\rm Name}&{\rm Wavelength}&{\rm Color}\\
&-&-\\
2&\alpha&121.567&{\rm UV}\\
3&\beta&102.572&{\rm UV}\\
4&\gamma&97.254&{\rm UV}\\
\vdots&\vdots&\vdots&\vdots\\
\infty&{\rm limit}&91.175&{\rm UV}
\end{matrix}$$
called Lyman series of the hydrogen atom.
\end{fact}

Observe that all the Lyman series lies in UV, which is invisible to the naked eye. Due to this fact, this series, while theoretically being the most important, was discovered only second. The first discovery, which was the big one, and the breakthrough, was by Balmer, the founding father of all this, back in 1885, in the visible range, as follows:

\begin{fact}[Balmer]
The hydrogen atom has spectral lines given by the formula
$$\frac{1}{\lambda}=R\left(\frac{1}{4}-\frac{1}{n^2}\right)$$
where $R\simeq 1.097\times 10^7$ and $n\geq3$, which are as follows,
$$\begin{matrix}
n&{\rm Name}&{\rm Wavelength}&{\rm Color}\\
&-&-\\
3&\alpha&656.279&{\rm red}\\
4&\beta&486.135&{\rm aqua}\\
5&\gamma&434.047&{\rm blue}\\
6&\delta&410.173&{\rm violet}\\
7&\varepsilon&397.007&{\rm UV}\\
\vdots&\vdots&\vdots&\vdots\\
\infty&{\rm limit}&346.600&{\rm UV}
\end{matrix}$$
called Balmer series of the hydrogen atom.
\end{fact}

So, this was Balmer's original result, which started everything. As a third main result now, this time in IR, due to Paschen in 1908, we have:

\begin{fact}[Paschen]
The hydrogen atom has spectral lines given by the formula
$$\frac{1}{\lambda}=R\left(\frac{1}{9}-\frac{1}{n^2}\right)$$
where $R\simeq 1.097\times 10^7$ and $n\geq4$, which are as follows,
$$\begin{matrix}
n&{\rm Name}&{\rm Wavelength}&{\rm Color}\\
&-&-\\
4&\alpha&1875&{\rm IR}\\
5&\beta&1282&{\rm IR}\\
6&\gamma&1094&{\rm IR}\\
\vdots&\vdots&\vdots&\vdots\\
\infty&{\rm limit}&820.4&{\rm IR}
\end{matrix}$$
called Paschen series of the hydrogen atom.
\end{fact}

Observe the striking similarity between the above three results. In fact, we have here the following fundamental, grand result, due to Rydberg in 1888, based on the Balmer series, and with later contributions by Ritz in 1908, using the Lyman series as well:

\begin{conclusion}[Rydberg, Ritz]
The spectral lines of the hydrogen atom are given by the Rydberg formula, depending on integer parameters $n_1<n_2$,
$$\frac{1}{\lambda_{n_1n_2}}=R\left(\frac{1}{n_1^2}-\frac{1}{n_2^2}\right)$$
with $R$ being the Rydberg constant for hydrogen, which is as follows:
$$R\simeq1.096\ 775\ 83\times 10^7$$
These spectral lines combine according to the Ritz-Rydberg principle, as follows:
$$\frac{1}{\lambda_{n_1n_2}}+\frac{1}{\lambda_{n_2n_3}}=\frac{1}{\lambda_{n_1n_3}}$$
Similar formulae hold for other atoms, with suitable fine-tunings of $R$.
\end{conclusion}

Here the first part, the Rydberg formula, generalizes the results of Lyman, Balmer, Paschen, which appear at $n_1=1,2,3$, at least retrospectively. The Rydberg formula predicts further spectral lines, appearing at $n_1=4,5,6,\ldots\,$, and these were discovered later, by Brackett in 1922, Pfund in 1924, Humphreys in 1953, and others afterwards, with all these extra lines being in far IR. The simplified complete table is as follows:
$$\begin{matrix}
n_1&n_2&{\rm Series\ name}&{\rm Wavelength}\ n_2=\infty&{\rm Color}\ n_2=\infty\\
&&-&-\\
1&2-\infty&{\rm Lyman}&91.13\ {\rm nm}&{\rm UV}\\
2&3-\infty&{\rm Balmer}&364.51\ {\rm nm}&{\rm UV}\\
3&4-\infty&{\rm Paschen}&820.14\ {\rm nm}&{\rm IR}\\
&&-&-\\
4&5-\infty&{\rm Brackett}&1458.03\ {\rm nm}&{\rm far\ IR}\\
5&6-\infty&{\rm Pfund}&2278.17\ {\rm nm}&{\rm far\ IR}\\
6&7-\infty&{\rm Humphreys}&3280.56\ {\rm nm}&{\rm far\ IR}\\
\vdots&\vdots&\vdots&\vdots&\vdots\\
\end{matrix}$$

Regarding the last assertion, concerning other elements, this was something conjectured and partly verified by Ritz, and fully verified and clarified later, via many experiments, the fine-tuning of $R$ being basically $R\to RZ^2$, where $Z$ is the atomic number. 

\bigskip

From a theoretical physics viewpoint, the main result remains the middle assertion, called Ritz-Rydberg combination principle. This is something at the same time extremely simple, and completely puzzling, the informal conclusion being as follows:

\begin{thought}
The simplest observables of the hydrogen atom, combining via
$$\frac{1}{\lambda_{n_1n_2}}+\frac{1}{\lambda_{n_2n_3}}=\frac{1}{\lambda_{n_1n_3}}$$
look like quite weird quantities. Why wouldn't they just sum normally.
\end{thought}

Fortunately, mathematics comes to the rescue. Indeed, the Ritz-Rydberg combination principle reminds the formula $e_{n_1n_2}e_{n_2n_3}=e_{n_1n_3}$ for the usual matrix units $e_{ij}:e_j\to e_i$. In short, we are in familiar territory here, and we can start dreaming of:

\begin{thought}
Observables in quantum mechanics should be some sort of infinite matrices, generalizing  the Lyman, Balmer, Paschen lines of the hydrogen atom, and multiplying between them as the matrices do, as to produce further observables.
\end{thought}

Good news, time to put everything together. As a main problem that we would like to solve, we have the understanding the intimate structure of matter, at the atomic level. There is of course a long story here, regarding the intimate structure of matter, going back centuries and even millennia ago, and our presentation here will be quite simplified. As a starting point, since we need a starting point, let us agree on:

\begin{claim}
Ordinary matter is made of small particles called atoms, with each atom appearing as a mix of even smaller particles, namely protons $+$, neutrons $0$ and electrons $-$, with the same number of protons $+$ and electrons $-$.
\end{claim}

As a first observation, this is something which does not look obvious at all, with probably lots of work, by many people, being involved, as to lead to this claim. And so it is. The story goes back to the discovery of charges and electricity, which were attributed to a small particle, the electron $-$. Now since matter is by default neutral, this naturally leads to the consideration to the proton $+$, having the same charge as the electron.

\bigskip

But, as a natural question, why should be these electrons $-$ and protons $+$ that small? And also, what about the neutron 0? These are not easy questions, and the fact that it is so came from several clever experiments. Let us first recall that careful experiments with tiny particles are practically impossible. However, all sorts of brutal experiments, such as bombarding matter with other pieces of matter, accelerated to the extremes, or submitting it to huge electric and magnetic fields, do work. And it is such kind of experiments, due to Thomson, Rutherford and others, ``peeling off'' protons $+$, neutrons $0$ and electrons $-$ from matter, and observing them, that led to the conclusion that these small beasts $+,0,-$ exist indeed, in agreement with Claim 8.30.

\bigskip

Of particular importance here was as well the radioactivity theory of Becquerel and Pierre and Marie Curie, involving this time such small beasts, or perhaps some related radiation, peeling off by themselves, in heavy elements such as uranium $_{92}{\rm U}$, polonium $_{84}{\rm Po}$ and radium $_{88}{\rm Ra}$. And there was also Einstein's work on the photoelectric effect, light interacting with matter, suggesting that even light itself might have associated to it some kind of particle, called photon. All this goes of course beyond Claim 8.30, with further particles involved, and more on this later, but as a general idea, all this deluge of small particle findings, all coming around 1900-1910, further solidified Claim 8.30.

\bigskip

So, taking now Claim 8.30 for granted, how are then the atoms organized, as mixtures of protons $+$, neutrons $0$ and electrons $-$? The answer here lies again in the above-mentioned ``brutal'' experiments of Thomson, Rutherford and others, which not only proved Claim 8.30, but led to an improved version of it, as follows:

\index{atom}

\begin{claim}
The atoms are formed by a core of protons $+$ and neutrons $0$, surrounded by a cloud of electrons $-$, gravitating around the core.
\end{claim}

This is a considerable advance, because we are now into familiar territory, namely some kind of mechanics. And with this in mind, all the pieces of our puzzle start fitting together, and we are led to the following grand conclusion:

\begin{claim}[Bohr and others]
The atoms are formed by a core of protons and neutrons, surrounded by a cloud of electrons, basically obeying to a modified version of electromagnetism. And with a fine mechanism involved, as follows:
\begin{enumerate}
\item The electrons are free to move only on certain specified elliptic orbits, labeled $1,2,3,\ldots\,$, situated at certain specific heights.

\item The electrons can jump or fall between orbits $n_1<n_2$, absorbing or emitting light and heat, that is, electromagnetic waves, as accelerating charges.

\item The energy of such a wave, coming from $n_1\to n_2$ or $n_2\to n_1$, is given, via the Planck viewpoint, by the Rydberg formula, applied with $n_1<n_2$. 

\item The simplest such jumps are those observed by Lyman, Balmer, Paschen. And multiple jumps explain the Ritz-Rydberg formula.
\end{enumerate}
\end{claim}

And isn't this beautiful. Moreover, some further claims, also by Bohr and others, are that the theory can be further extended and fine-tuned as to explain many other phenomena, such as the above-mentioned findings of Einstein, and of Becquerel and Pierre and Marie Curie, and generally speaking, all the physics and chemistry known. 

\smallskip

And the story is not over here. Following now Heisenberg, the next claim is that the underlying mathematics in all the above can lead to a beautiful axiomatization of quantum mechanics, as a ``matrix mechanics'', along the lines of Thought 8.29. We will be back to all this in chapter 16 below, at the end of the present book.

\section*{8e. Exercises}

We had a tough physics chapter here, and here are some physics exercises:

\begin{exercise}
Learn more about the Laplace equation, and about Maxwell too.
\end{exercise}

\begin{exercise}
Learn about soap films, and their relation with harmonic functions.
\end{exercise}

\begin{exercise}
Explore a bit the wave and heat equations, in $2$ dimensions.
\end{exercise}

\begin{exercise}
Work out the details for the decomposition of wave packets.
\end{exercise}

As bonus exercise, try developing quantum mechanics, based on the above. With this being not a joke, Heisenberg did not know much more than that, when he did it.

\part{Several variables}

\ \vskip50mm

\begin{center}
{\em This is my church

This is where I heal my hurts

For tonight

God is a DJ}
\end{center}

\chapter{Linear maps}

\section*{9a. Linear maps}

Welcome to several variables. We have already met them, in chapter 8, and obviously this is a quite tricky business. But we have a whole 200 pages for clarifying this. As a main question, that we would like to solve, over all these pages to follow, we have:

\begin{question}
The main idea of calculus was that the functions $f:\mathbb R\to\mathbb R$ are locally approximately linear. In view of this, when looking for generalizations:
\begin{enumerate}
\item What can we say about the linear maps $f:\mathbb R^N\to\mathbb R^M$?

\item Then, what can we say about the arbitrary functions $f:\mathbb R^N\to\mathbb R^M$?
\end{enumerate}
\end{question}

Which sounds good, now we have a serious plan, and time to develop it. Regarding the first question, about the linear maps, we first have the following result:

\index{linear map}
\index{linear function}
\index{matrix}
\index{rectangular matrix}

\begin{theorem}
The linear maps $f:\mathbb R^N\to\mathbb R^M$ are in correspondence with the matrices $A\in M_{M\times N}(\mathbb R)$, with the linear map associated to such a matrix being
$$f(x)=Ax$$
and with the matrix associated to a linear map being $A_{ij}=<f(e_j),e_i>$. Similarly, the linear maps $f:\mathbb C^N\to\mathbb C^N$ are in correspondence with the matrices $A\in M_{M\times N}(\mathbb C)$.
\end{theorem}

\begin{proof}
The first assertion is clear, because a linear map $f:\mathbb R^N\to\mathbb R^M$ must send a vector $x\in\mathbb R^N$ to a certain vector $f(x)\in\mathbb R^M$, all whose components are linear combinations of the components of $x$. Thus, we can write, for certain numbers $A_{ij}\in\mathbb R$:
$$f\begin{pmatrix}
x_1\\
\vdots\\
x_N
\end{pmatrix}
=\begin{pmatrix}
A_{11}x_1+\ldots+A_{1N}x_N\\
\vdots\\
A_{M1}x_1+\ldots+A_{MN}x_N
\end{pmatrix}$$

Now the parameters $A_{ij}\in\mathbb R$ can be regarded as being the entries of a rectangular matrix $A\in M_{M\times N}(\mathbb R)$, and with the usual convention for the rectangular matrix multiplication, the above formula is precisely the one in the statement, namely:
$$f(x)=Ax$$

Regarding the second assertion, with $f(x)=Ax$ as above, if we denote by $e_1,\ldots,e_N$ the standard basis of $\mathbb R^N$, then we have the following formula:
$$f(e_j)
=\begin{pmatrix}
A_{1j}\\
\vdots\\
A_{Mj}
\end{pmatrix}$$

But this gives $<f(e_j),e_i>=A_{ij}$, as desired. As for the last assertion, regarding the complex maps and matrices, the proof here is similar, and with the complex scalar products being by definition given, here and in what follows, by $<x,y>=\sum_ix_i\bar{y}_i$.
\end{proof}

At the level of examples, let us first discuss the linear maps $f:\mathbb R^2\to\mathbb R^2$. We have:

\index{rotation}
\index{symmetry}

\begin{proposition}
The rotation of angle $t\in\mathbb R$, and the symmetry with respect to the $Ox$ axis rotated by an angle $t/2\in\mathbb R$, are given by the matrices
$$R_t=\begin{pmatrix}\cos t&-\sin t\\ \sin t&\cos t\end{pmatrix}\quad,\quad 
S_t=\begin{pmatrix}\cos t&\sin t\\ \sin t&-\cos t\end{pmatrix}$$
both depending on $t\in\mathbb R$ taken modulo $2\pi$.
\end{proposition}

\begin{proof}
The rotation being linear, it must correspond to a certain matrix:
$$R_t=\begin{pmatrix}a&b\\ c&d\end{pmatrix}$$

We can guess this matrix, via its action on the basic coordinate vectors $\binom{1}{0}$ and $\binom{0}{1}$. Indeed, a quick picture in the plane shows that we must have:
$$\begin{pmatrix}a&b\\ c&d\end{pmatrix}\begin{pmatrix}1\\ 0\end{pmatrix}=
\begin{pmatrix}\cos t\\ \sin t\end{pmatrix}\quad,\quad 
\begin{pmatrix}a&b\\ c&d\end{pmatrix}\begin{pmatrix}0\\ 1\end{pmatrix}=
\begin{pmatrix}-\sin t\\ \cos t\end{pmatrix}$$

Guessing now the matrix is not complicated, because the first equality gives us the first column, and the second equality gives us the second column:
$$\binom{a}{c}=\begin{pmatrix}\cos t\\ \sin t\end{pmatrix}\quad,\quad 
\binom{b}{d}=\begin{pmatrix}-\sin t\\ \cos t\end{pmatrix}$$

Thus, we can just put together these two vectors, and we obtain our matrix $R_t$. As for the symmetry, the proof here is similar, again by computing $S_t\binom{1}{0}$ and $S_t\binom{0}{1}$.
\end{proof}

Let us record as well a result regarding the projections, as follows:

\index{projection}

\begin{proposition}
The projection on the $Ox$ axis rotated by an angle $t/2\in\mathbb R$ is
$$P_t=\frac{1}{2}\begin{pmatrix}1+\cos t&\sin t\\ \sin t&1-\cos t\end{pmatrix}$$
depending on $t\in\mathbb R$ taken modulo $2\pi$.
\end{proposition}

\begin{proof}
A quick picture in the plane, using similarity of triangles, and the basic trigonometry formulae for the duplication of angles, show that we must have:
$$P_t\begin{pmatrix}1\\ 0\end{pmatrix}
=\cos\frac{t}{2}\binom{\cos\frac{t}{2}}{\sin\frac{t}{2}}
=\frac{1}{2}\begin{pmatrix}1+\cos t\\ \sin t\end{pmatrix}$$

Similarly, another quick picture plus trigonometry show that we must have:
$$P_t\begin{pmatrix}0\\ 1\end{pmatrix}
=\sin\frac{t}{2}\binom{\cos\frac{t}{2}}{\sin\frac{t}{2}}
=\frac{1}{2}\begin{pmatrix}\sin t\\1-\cos t\end{pmatrix}$$

Now by putting together these two vectors, and we obtain our matrix.
\end{proof}

Back to Theorem 9.2, our claim is that, no matter what we want to do with $f$ or $A$, we will run at some point into their adjoints $f^*$ and $A^*$, constructed as follows:

\index{adjoint matrix}
\index{adjoint map}

\begin{theorem}
The adjoint linear map $f^*:\mathbb C^N\to\mathbb C^N$, which is given by
$$<f(x),y>=<x,f^*(y)>$$
corresponds to the adjoint matrix $A^*\in M_N(\mathbb C)$, given by
$$(A^*)_{ij}=\bar{A}_{ji}$$
via the correspondence between linear maps and matrices constructed above.
\end{theorem}

\begin{proof}
Given a linear map $f:\mathbb C^N\to\mathbb C^N$, fix $y\in\mathbb C^N$, and consider the linear form $\varphi(x)=<f(x),y>$. This form must be as follows, for a certain vector $f^*(y)\in\mathbb C^N$:
$$\varphi(x)=<x,f^*(y)>$$

Thus, we have constructed a map $y\to f^*(y)$ as in the statement, which is obviously linear, and that we can call $f^*$. Now by taking the vectors $x,y\in\mathbb C^N$ to be elements of the standard basis of $\mathbb C^N$, our defining formula for $f^*$ reads:
$$<f(e_i),e_j>=<e_i,f^*(e_j)>$$

By reversing the scalar product on the right, this formula can be written as:
$$<f^*(e_j),e_i>=\overline{<f(e_i),e_j>}$$

But this means that the matrix of $f^*$ is given by $(A^*)_{ij}=\bar{A}_{ji}$, as desired.
\end{proof}

Getting back to our claim, the adjoints $*$ are indeed ubiquitous, as shown by:

\index{scalar product}
\index{unitary}
\index{projection}
\index{polarization formula}

\begin{theorem}
The following happen:
\begin{enumerate}
\item $f(x)=Ux$ with $U\in M_N(\mathbb C)$ is an isometry precisely when $U^*=U^{-1}$.

\item $f(x)=Px$ with $P\in M_N(\mathbb C)$ is a projection precisely when $P^2=P^*=P$.
\end{enumerate}
\end{theorem}

\begin{proof}
Let us first recall that the lengths, or norms, of the vectors $x\in\mathbb C^N$ can be recovered from the knowledge of the scalar products, as follows:
$$||x||=\sqrt{<x,x>}$$

Conversely, we can recover the scalar products out of norms, by using the following difficult to remember formula, called complex polarization identity:
\begin{eqnarray*}
&&||x+y||^2-||x-y||^2+i||x+iy||^2-i||x-iy||^2\\
&=&||x||^2+||y||^2-||x||^2-||y||^2+i||x||^2+i||y||^2-i||x||^2-i||y||^2\\
&&+2Re(<x,y>)+2Re(<x,y>)+2iIm(<x,y>)+2iIm(<x,y>)\\
&=&4<x,y>
\end{eqnarray*}

Finally, we will use Theorem 9.5, and more specifically the following formula coming from there, valid for any matrix $A\in M_N(\mathbb C)$ and any two vectors $x,y\in\mathbb C^N$:
$$<Ax,y>=<x,A^*y>$$

(1) Given a matrix $U\in M_N(\mathbb C)$, we have indeed the following equivalences, with the first one coming from the polarization identity, and with the other ones being clear:
\begin{eqnarray*}
||Ux||=||x||
&\iff&<Ux,Uy>=<x,y>\\
&\iff&<x,U^*Uy>=<x,y>\\
&\iff&U^*Uy=y\\
&\iff&U^*U=1\\
&\iff&U^*=U^{-1}
\end{eqnarray*}

(2) Given a matrix $P\in M_N(\mathbb C)$, in order for $x\to Px$ to be an oblique projection, we must have $P^2=P$. Now observe that this projection is orthogonal when:
\begin{eqnarray*}
<Px-x,Py>=0
&\iff&<P^*Px-P^*x,y>=0\\
&\iff&P^*Px-P^*x=0\\
&\iff&P^*P-P^*=0\\
&\iff&P^*P=P^*
\end{eqnarray*}

The point now is that by conjugating the last formula, we obtain $P^*P=P$. Thus we must have $P=P^*$, and this gives the result. 
\end{proof}

Summarizing, the linear operators come in pairs $T,T^*$, and the associated matrices come as well in pairs $A,A^*$. We will keep this in mind, and come back to it later.

\section*{9b. Matrix inversion} 

We have seen so far that most of the interesting maps $f:\mathbb R^N\to\mathbb R^N$ that we know, such as the rotations, symmetries and projections, are linear, and can be written in the following form, with $A\in M_N(\mathbb R)$ being a square matrix:
$$f(v)=Av$$

We develop now more general theory for such linear maps. We will be interested in the question of inverting the linear maps $f:\mathbb R^N\to\mathbb R^N$. And the point is that this is the same question as inverting the corresponding matrices $A\in M_N(\mathbb R)$, due to:

\index{invertible matrix}

\begin{theorem}
A linear map $f:\mathbb R^N\to\mathbb R^N$, written as
$$f(v)=Av$$
is invertible precisely when $A$ is invertible, and in this case we have $f^{-1}(v)=A^{-1}v$.
\end{theorem}

\begin{proof}
This comes indeed from the fact that, with the notation $f_A(v)=Av$, we have the formula $f_Af_B=f_{AB}$. Thus, we are led to the conclusion in the statement.
\end{proof}

In order to study invertibility questions, for matrices or linear maps, let us begin with some examples. In the simplest case, in 2 dimensions, the result is as follows:

\index{inversion formula}

\begin{theorem}
We have the following inversion formula, for the $2\times2$ matrices:
$$\begin{pmatrix}a&b\\ c&d\end{pmatrix}^{-1}
=\frac{1}{ad-bc}\begin{pmatrix}d&-b\\ -c&a\end{pmatrix}$$
When $ad-bc=0$, the matrix is not invertible.
\end{theorem}

\begin{proof}
We have two assertions to be proved, the idea being as follows:

\medskip

(1) As a first observation, when $ad-bc=0$ we must have, for some $\lambda\in\mathbb R$:
$$b=\lambda a\quad,\quad 
d=\lambda c$$

Thus our matrix must be of the following special type:
$$\begin{pmatrix}a&b\\ c&d\end{pmatrix}=\begin{pmatrix}a&\lambda a\\ a&\lambda c\end{pmatrix}$$

But in this case the columns are proportional, so the linear map associated to the matrix is not invertible, and so the matrix itself is not invertible either.

\medskip

(2) When $ad-bc\neq 0$, let us look for an inversion formula of the following type:
$$\begin{pmatrix}a&b\\ c&d\end{pmatrix}^{-1}
=\frac{1}{ad-bc}\begin{pmatrix}*&*\\ *&*\end{pmatrix}$$

We must therefore solve the following system of equations:
$$\begin{pmatrix}a&b\\ c&d\end{pmatrix}
\begin{pmatrix}*&*\\ *&*\end{pmatrix}=
\begin{pmatrix}ad-bc&0\\ 0&ad-bc\end{pmatrix}$$

But the solution to these equations is obvious, is as follows:
$$\begin{pmatrix}a&b\\ c&d\end{pmatrix}
\begin{pmatrix}d&-b\\ -c&a\end{pmatrix}=
\begin{pmatrix}ad-bc&0\\ 0&ad-bc\end{pmatrix}$$

Thus, we are led to the formula in the statement.
\end{proof}

In order to deal now with the inversion problem in general, for the arbitrary matrices $A\in M_N(\mathbb R)$, we will use the same method as the one above, at $N=2$. Let us write indeed our matrix as follows, with $v_1,\ldots,v_N\in\mathbb R^N$ being its column vectors:
$$A=[v_1,\ldots,v_N]$$

We know from the above that, in order for the matrix $A$ to be invertible, the vectors $v_1,\ldots,v_N$ must be linearly independent. Thus, we are led into the question of understanding when a family of vectors $v_1,\ldots,v_N\in\mathbb R^N$ are linearly independent. In order to deal with this latter question, let us introduce the following notion:

\index{volume of parallelepiped}

\begin{definition}
Associated to any vectors $v_1,\ldots,v_N\in\mathbb R^N$ is the volume
$${\rm det}^+(v_1\ldots v_N)=vol<v_1,\ldots,v_N>$$
of the parallelepiped made by these vectors.
\end{definition}

Here the volume is taken in the standard $N$-dimensional sense. At $N=1$ this volume is a length, at $N=2$ this volume is an area, at $N=3$ this is the usual 3D volume, and so on. In general, the volume of a body $X\subset\mathbb R^N$ is by definition the number $vol(X)\in[0,\infty]$ of copies of the unit cube $C\subset\mathbb R^N$ which are needed for filling $X$. Now with this notion in hand, in relation with our inversion problem, we have the following statement:

\index{invertible matrix}

\begin{proposition}
The quantity ${\rm det}^+$ that we constructed, regarded as a function of the corresponding square matrices, formed by column vectors,
$${\rm det}^+:M_N(\mathbb R)\to\mathbb R_+$$
has the property that a matrix $A\in M_N(\mathbb R)$ is invertible precisely when ${\rm det}^+(A)>0$.
\end{proposition}

\begin{proof}
This follows from the fact that a matrix $A\in M_N(\mathbb R)$ is invertible precisely when its column vectors $v_1,\ldots,v_N\in\mathbb R^N$ are linearly independent. But this latter condition is equivalent to the fact that we must have the following strict inequality: 
$$vol<v_1,\ldots,v_N>>0$$

Thus, we are led to the conclusion in the statement.
\end{proof}

Summarizing, all this leads us into the explicit computation of ${\rm det}^+$. As a first observation, in 1 dimension we obtain the absolute value of the real numbers:
$${\rm det}^+(a)=|a|$$

In 2 dimensions now, the computation is non-trivial, and we have the following result, making the link with our main inversion result so far, namely Theorem 9.8:

\begin{theorem}
In $2$ dimensions we have the following formula,
$${\rm det}^+\begin{pmatrix}a&b\\ c&d\end{pmatrix}=|ad-bc|$$
with ${\rm det}^+:M_2(\mathbb R)\to\mathbb R_+$ being the function constructed above.
\end{theorem}

\begin{proof}
We must show that the area of the parallelogram formed by $\binom{a}{c},\binom{b}{d}$ equals $|ad-bc|$. We can assume $a,b,c,d>0$ for simplifying, the proof in general being similar. Moreover, by switching if needed the vectors $\binom{a}{c},\binom{b}{d}$, we can assume that we have:
$$\frac{a}{c}>\frac{b}{d}$$

Now let us slide the upper side of the parallelogram downwards left, until we reach $Oy$. Our parallelogram, which has not changed its area in this process, becomes:
$$
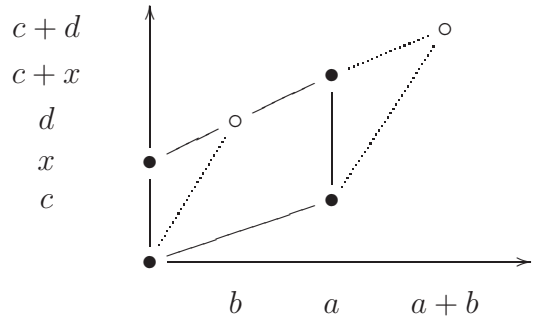
\xymatrix@R=1pt@C=16pt{
&&&\\
c+d&&&&\circ\\
c+x&&&\bullet\ar@{.}[ur]&\\
d&&\circ\ar@{-}[ur]&&&&&\\
x&\bullet\ar@{-}[ur]\ar[uuuu]&&&\\
c&&&\bullet\ar@{-}[uuu]\ar@{.}[uuuur]&\\
&&&&&&\\
&\bullet\ar@{-}[uurr]\ar@{-}[uuu]\ar[rrrr]\ar@{.}[uuuur]&&&&\\
&&\ b\ &\ a\ &a+b}$$

Moreover, we can further modify this parallelogram, once again by not altering its area, by sliding the right side downwards, until we reach the $Ox$ axis:
$$
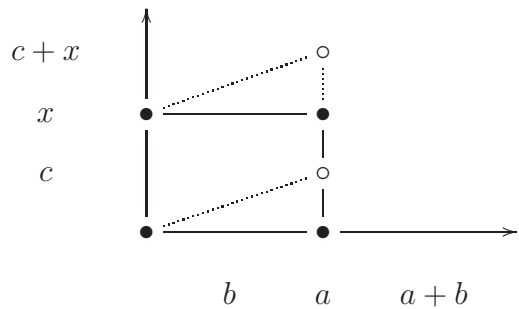
\xymatrix@R=10pt@C=15pt{
&&&\\
c+x&&&\circ&\\
x&\bullet\ar@{.}[urr]\ar[uu]\ar@{-}[rr]&&\bullet\ar@{.}[u]&\\
c&&&\circ\ar@{-}[u]&\\
&\bullet\ar@{.}[urr]\ar@{-}[uu]\ar@{-}[rr]&&\bullet\ar@{-}[u]\ar[rr]&&\\
&&\ b\ &\ a\ &a+b}$$

Let us compute now the area. Since our two sliding operations have not changed the area of the original parallelogram, this area is given by the following formula:
$$A=ax$$

In order to compute the quantity $x$, observe that in the context of the first move, we have two similar triangles, according to the following picture:
$$
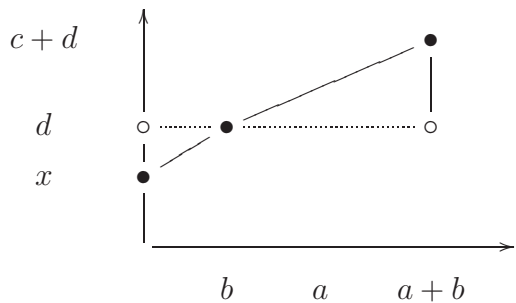
\xymatrix@R=4pt@C=15pt{
&&&&\\
c+d&&&&\bullet\\
&&&&&&\\
d&\circ\ar@{.}[r]\ar[uuu]&\bullet\ar@{.}[rr]\ar@{-}[uurr]&&\circ\ar@{-}[uu]\\
x&\bullet\ar@{-}[u]\ar@{-}[ur]&&&\\
&&&&\\
&\ar@{-}[uu]\ar[rrrr]&&&&\\
&&\ b\ &\ a\ &a+b}$$

Thus, we are led to the following equation for the number $x$:
$$\frac{d-x}{b}=\frac{c}{a}$$

But this gives $x=d-bc/a$, so the area of our parallelogram, or rather of the final rectangle obtained from it, which has the same area as the original parallelogram, is:
$$A=ax=ad-bc$$

We are therefore led to the conclusion in the statement.
\end{proof}

All this is very nice, and obviously we have a beginning of theory here. However, when looking carefully, we can see that our theory has a weakness, because:

\medskip

\begin{enumerate}
\item In 1 dimension the number $a$, which is the simplest function of $a$ itself, is certainly a better quantity than the number $|a|$.

\medskip

\item In 2 dimensions the number $ad-bc$, which is linear in $a,b,c,d$, is certainly a better quantity than the number $|ad-bc|$. 
\end{enumerate}

\medskip

So, let us upgrade now our theory, by constructing a better function, which takes signed values. In order to do this, we must come up with a way of splitting the systems of vectors $v_1,\ldots,v_N\in\mathbb R^N$ into two classes, call them positive and negative. And here, the answer is quite clear, because a bit of thinking leads to the following definition:

\index{oriented system of vectors}
\index{unoriented system of vectors}
\index{sign of system of vectors}

\begin{definition}
A system of vectors $v_1,\ldots,v_N\in\mathbb R^N$ is called:
\begin{enumerate}
\item Oriented, if one can continuously pass from the standard basis to it.

\item Unoriented, otherwise.
\end{enumerate}
The associated sign is $+$ in the oriented case, and $-$ in the unoriented case. 
\end{definition}

As a first example, in 1 dimension the basis consists of the single vector $e=1$, which can be continuously deformed into any vector $a>0$. Thus, the sign is the usual one:
$$sgn(a)=
\begin{cases}
+&{\rm if}\ a>0\\
-&{\rm if}\ a<0
\end{cases}$$

Thus, in connection with our original question, we are definitely on the good track, because when multiplying $|a|$ by this sign we obtain $a$ itself, as desired:
$$a=sgn(a)|a|$$

In 2 dimensions now, the explicit formula of the sign is as follows:

\begin{proposition}
We have the following formula, valid for any $2$ vectors in $\mathbb R^2$,
$$sgn\left[\binom{a}{c},\binom{b}{d}\right]=sgn(ad-bc)$$
with the sign function on the right being the usual one, in $1$ dimension.
\end{proposition}

\begin{proof}
According to our conventions, the sign of $\binom{a}{c},\binom{b}{d}$ is as follows:

\medskip

(1) The sign is $+$ when these vectors come in this precise order with respect to the counterclockwise rotation in the plane, around 0.

\medskip

(2) The sign is $-$ otherwise, meaning when these vectors come in this order with respect to the clockwise rotation in the plane, around 0. 

\medskip

If we assume now $a,b,c,d>0$ for simplifying, we are left with comparing the angles having the numbers $c/a$ and $d/b$ as tangents, and we obtain in this way:
$$sgn\left[\binom{a}{c},\binom{b}{d}\right]=
\begin{cases}
+&{\rm if}\ \frac{c}{a}<\frac{d}{b}\\
-&{\rm if}\ \frac{c}{a}>\frac{d}{b}
\end{cases}$$

But this gives the formula in the statement. The proof in general is similar.
\end{proof}

Once again, in connection with our original question, we are on the good track, because when multiplying $|ad-bc|$ by this sign we obtain $ad-bc$ itself, as desired:
$$ad-bc=sgn(ad-bc)|ad-bc|$$

At the level of the general results now, we have:

\begin{proposition}
The orientation of a system of vectors changes as follows:
\begin{enumerate}
\item If we switch the sign of a vector, the associated sign switches.

\item If we permute two vectors, the associated sign switches as well.
\end{enumerate}
\end{proposition}

\begin{proof}
Both these assertions are clear from the definition of the sign, because the two operations in question change the orientation of the system of vectors.
\end{proof}

With the above notion in hand, we can now formulate:

\index{determinant}
\index{signed volume}

\begin{definition}
The determinant of $v_1,\ldots,v_N\in\mathbb R^N$ is the signed volume
$$\det(v_1\ldots v_N)=\pm vol<v_1,\ldots,v_N>$$
of the parallelepiped made by these vectors.
\end{definition}

In other words, we are upgrading here Definition 9.9, by adding a sign to the quantity ${\rm det}^+$ constructed there, as to potentially reach to good additivity properties:
$$\det(v_1\ldots v_N)=\pm {\rm det}^+(v_1\ldots v_N)$$

In relation with our original inversion problem for the square matrices, this upgrade does not change what we have so far, and we have the following statement:

\index{invertible matrix}

\begin{theorem}
The quantity $\det$ that we constructed, regarded as a function of the corresponding square matrices, formed by column vectors,
$$\det:M_N(\mathbb R)\to\mathbb R$$
has the property that a matrix $A\in M_N(\mathbb R)$ is invertible precisely when $\det(A)\neq 0$.
\end{theorem}

\begin{proof}
We know from the above that a matrix $A\in M_N(\mathbb R)$ is invertible precisely when ${\rm det}^+(A)=|\det A|$ is strictly positive, and this gives the result.
\end{proof}

Let us try now to compute the determinant. In 1 dimension we have of course the formula $\det(a)=a$, because the absolute value fits, and so does the sign:
$$\det(a)
=sgn(a)\times|a|
=a$$

In 2 dimensions now, we have the following result:

\begin{theorem}
In $2$ dimensions we have the following formula,
$$\begin{vmatrix}a&b\\ c&d\end{vmatrix}=ad-bc$$
with $|\,.\,|=\det$ being the determinant function constructed above.
\end{theorem}

\begin{proof}
According to our definition, to the computation in Theorem 9.11, and to the sign formula from Proposition 9.13, the determinant of a $2\times2$ matrix is given by:
\begin{eqnarray*}
\det\begin{pmatrix}a&b\\ c&d\end{pmatrix}
&=&sgn\left[\binom{a}{c},\binom{b}{d}\right]\times {\rm det}^+\begin{pmatrix}a&b\\ c&d\end{pmatrix}\\
&=&sgn\left[\binom{a}{c},\binom{b}{d}\right]\times|ad-bc|\\
&=&sgn(ad-bc)\times|ad-bc|\\
&=&ad-bc
\end{eqnarray*}

Thus, we have obtained the formula in the statement.
\end{proof}

\section*{9c. The determinant} 

In order to discuss now arbitrary dimensions, we will need a number of theoretical results. Here is a first series of formulae, coming straight from definitions:

\begin{theorem}
The determinant has the following properties:
\begin{enumerate}
\item When multiplying by scalars, the determinant gets multiplied as well:
$$\det(\lambda_1v_1,\ldots,\lambda_Nv_N)=\lambda_1\ldots\lambda_N\det(v_1,\ldots,v_N)$$

\item When permuting two columns, the determinant changes the sign:
$$\det(\ldots,u,\ldots,v,\ldots)=-\det(\ldots,v,\ldots,u,\ldots)$$

\item The determinant $\det(e_1,\ldots,e_N)$ of the standard basis of $\mathbb R^N$ is $1$.
\end{enumerate}
\end{theorem}

\begin{proof}
All this is clear from definitions, as follows:

\medskip

(1) This follows from definitions, and from Proposition 9.14 (1).

\medskip

(2) This follows as well from definitions, and from Proposition 9.14 (2).

\medskip

(3) This is clear from our definition of the determinant.
\end{proof}

As an application of the above result, we have:

\begin{theorem}
The determinant of a diagonal matrix is given by:
$$\begin{vmatrix}
\lambda_1\\ 
&\ddots\\
&&\lambda_N\end{vmatrix}=\lambda_1\ldots\lambda_N$$
That is, we obtain the product of diagonal entries.
\end{theorem}

\begin{proof}
The formula in the statement is clear by using Theorem 9.18, which gives:
$$\begin{vmatrix}
\lambda_1\\ 
&\ddots\\
&&\lambda_N\end{vmatrix}
=\lambda_1\ldots\lambda_N
\begin{vmatrix}
1\\ 
&\ddots\\
&&1\end{vmatrix}
=\lambda_1\ldots\lambda_N$$

As for the last assertion, this is rather a remark.
\end{proof}

In order to reach to a more advanced theory, let us adopt now the linear map point of view. In this setting, the definition of the determinant reformulates as follows: 

\index{inflation coefficient}

\begin{theorem}
Given a linear map, written as $f(v)=Av$, its ``inflation coefficient'', obtained as the signed volume of the image of the unit cube, is given by:
$$I_f=\det A$$
More generally, $I_f$ is the inflation ratio of any parallelepiped in $\mathbb R^N$, via the transformation $f$. In particular $f$ is invertible precisely when $\det A\neq0$.
\end{theorem}

\begin{proof}
The only non-trivial thing in all this is the fact that the inflation coefficient $I_f$, as defined above, is independent of the choice of the parallelepiped. But this is a generalization of the Thales theorem, which follows from the Thales theorem itself.
\end{proof}

As a first application of the above linear map viewpoint, we have:

\index{determinant of products}

\begin{theorem}
We have the following formula, valid for any matrices $A,B$:
$$\det(AB)=\det A\cdot\det B$$
In particular, we have $\det(AB)=\det(BA)$.
\end{theorem}

\begin{proof}
The first formula follows from the formula $f_{AB}=f_Af_B$ for the associated linear maps. As for $\det(AB)=\det(BA)$, this is clear from the first formula.
\end{proof}

In general now, still at the theoretical level, we have the following key result:

\begin{theorem}
The determinant has the additivity property
$$\det(\ldots,u+v,\ldots)
=\det(\ldots,u,\ldots)
+\det(\ldots,v,\ldots)$$
valid for any choice of the vectors involved.
\end{theorem}

\begin{proof}
This follows indeed by doing some elementary geometry, in the spirit of the computations in the proof of Theorem 9.11, by using the Thales theorem. 
\end{proof}

As a basic application of the above result, we have:

\index{upper triangular matrix}
\index{lower triangular matrix}

\begin{theorem}
We have the following results:
\begin{enumerate}
\item The determinant of a diagonal matrix is the product of diagonal entries.

\item The same is true for the upper triangular matrices.

\item The same is true for the lower triangular matrices.
\end{enumerate}
\end{theorem}

\begin{proof}
Here (1) is something that we already know, and (2) comes as follows:
\begin{eqnarray*}
\begin{vmatrix}
\lambda_1&&&*\\ 
&\lambda_2\\
&&\ddots\\
0&&&\lambda_N\end{vmatrix}
&=&\begin{vmatrix}
\lambda_1&0&&*\\ 
&\lambda_2\\
&&\ddots\\
0&&&\lambda_N\end{vmatrix}\\
&&\vdots\\
&&\vdots\\
&=&\begin{vmatrix}
\lambda_1&&&0\\ 
&\lambda_2\\
&&\ddots\\
0&&&\lambda_N\end{vmatrix}\\
&=&\lambda_1\ldots\lambda_N
\end{eqnarray*}

As for (3), this comes in a similar way, by proceeding this time from right to left.
\end{proof}

As an important theoretical result now, we have:

\begin{theorem}
The determinant $\det:M_N(\mathbb R)\to\mathbb R$ is the unique map satisfying:
\begin{enumerate}
\item $\det(\ldots,u+v,\ldots)
=\det(\ldots,u,\ldots)
+\det(\ldots,v,\ldots)$.

\item $\det(\lambda_1v_1,\ldots,\lambda_Nv_N)=\lambda_1\ldots\lambda_N\det(v_1,\ldots,v_N)$.

\item $\det(\ldots,u,\ldots,v,\ldots)=-\det(\ldots,v,\ldots,u,\ldots)$.

\item The determinant $\det(e_1,\ldots,e_N)$ of the standard basis of $\mathbb R^N$ is $1$.
\end{enumerate}
\end{theorem}

\begin{proof}
The conditions in the statement are those from Theorems 9.18 and 9.22. As for the converse, it is routine to check that any map $\det':M_N(\mathbb R)\to\mathbb R$ satisfying (1-4) must coincide with $\det$ on the upper triangular matrices, and then on all matrices.
\end{proof}

Here is now another important theoretical result:

\index{row expansion}

\begin{theorem}
The determinant is subject to the row expansion formula
\begin{eqnarray*}
\begin{vmatrix}a_{11}&\ldots&a_{1N}\\
\vdots&&\vdots\\
a_{N1}&\ldots&a_{NN}\end{vmatrix}
&=&a_{11}\begin{vmatrix}a_{22}&\ldots&a_{2N}\\
\vdots&&\vdots\\
a_{N2}&\ldots&a_{NN}\end{vmatrix}
-a_{12}\begin{vmatrix}a_{21}&a_{23}&\ldots&a_{2N}\\
\vdots&\vdots&&\vdots\\
a_{N1}&a_{N3}&\ldots&a_{NN}\end{vmatrix}\\
&&+\ldots\ldots
+(-1)^{N+1}a_{1N}\begin{vmatrix}a_{21}&\ldots&a_{2,N-1}\\
\vdots&&\vdots\\
a_{N1}&\ldots&a_{N,N-1}\end{vmatrix}
\end{eqnarray*}
and this method fully computes it, by recurrence.
\end{theorem}

\begin{proof}
This follows from the fact that the formula in the statement produces a certain function $\det:M_N(\mathbb R)\to\mathbb R$, which has the 4 properties in Theorem 9.24.
\end{proof}

We can expand as well over the columns, as follows:

\index{column expansion}

\begin{theorem}
The determinant is subject to the column expansion formula
\begin{eqnarray*}
\begin{vmatrix}a_{11}&\ldots&a_{1N}\\
\vdots&&\vdots\\
a_{N1}&\ldots&a_{NN}\end{vmatrix}
&=&a_{11}\begin{vmatrix}a_{22}&\ldots&a_{2N}\\
\vdots&&\vdots\\
a_{N2}&\ldots&a_{NN}\end{vmatrix}
-a_{21}\begin{vmatrix}a_{12}&\ldots&a_{1N}\\
a_{32}&\ldots&a_{3N}\\
\vdots&&\vdots\\
a_{N2}&\ldots&a_{NN}\end{vmatrix}\\
&&+\ldots\ldots
+(-1)^{N+1}a_{N1}\begin{vmatrix}a_{12}&\ldots&a_{1N}\\
\vdots&&\vdots\\
a_{N-1,2}&\ldots&a_{N-1,N}\end{vmatrix}
\end{eqnarray*}
and this method fully computes it, by recurrence.
\end{theorem}

\begin{proof}
This follows by using the same argument as for the rows.
\end{proof}

As a first application of the above methods, we can now prove:

\index{Sarrus formula}

\begin{theorem}
The determinant of the $3\times3$ matrices is given by
$$\begin{vmatrix}a&b&c\\ d&e&f\\ g&h&i\end{vmatrix}=aei+bfg+cdh-ceg-bdi-afh$$
which can be memorized by using Sarrus' triangle method, ``triangles parallel to the diagonal, minus triangles parallel to the antidiagonal".
\end{theorem}

\begin{proof}
Here is the computation, using the above results:
\begin{eqnarray*}
\begin{vmatrix}a&b&c\\ d&e&f\\ g&h&i\end{vmatrix}
&=&a\begin{vmatrix}e&f\\h&i\end{vmatrix}
-b\begin{vmatrix}d&f\\g&i\end{vmatrix}
+c\begin{vmatrix}d&e\\g&h\end{vmatrix}\\
&=&a(ei-fh)-b(di-fg)+c(dh-eg)\\
&=&aei-afh-bdi+bfg+cdh-ceg\\
&=&aei+bfg+cdh-ceg-bdi-afh
\end{eqnarray*}

Thus, we obtain the formula in the statement.
\end{proof}

As a first application, we can now invert the $3\times3$ matrices, as follows:

\index{matrix inversion}

\begin{theorem}
The inverses of the $3\times3$ matrices are given by
$$\begin{pmatrix}a&b&c\\ d&e&f\\ g&h&i\end{pmatrix}^{-1}
=\frac{1}{D}\begin{pmatrix}ei-fh&ch-bi&bf-ce\\ fg-di&ai-cg&cd-af\\ dh-eg&bg-ah&ae-bd\end{pmatrix}$$
with $D$ being the determimant. When $D=0$, the matrix is not invertible.
\end{theorem}

\begin{proof}
As before for the $2\times 2$ matrices, in order for our matrix to be invertible, we must have $D\neq0$. The trick now is to look for solutions of the following problem:
$$\begin{pmatrix}a&b&c\\ d&e&f\\ g&h&i\end{pmatrix}
\begin{pmatrix}*&*&*\\ *&*&*\\ *&*&*\end{pmatrix}
=\begin{pmatrix}D&0&0\\ 0&D&0\\ 0&0&D\end{pmatrix}$$

We know from Theorem 9.27 that the determinant is given by:
$$D=aei+bfg+cdh-ceg-bdi-afh$$

But this leads, via some obvious choices, to the following solution:
$$\begin{pmatrix}*&*&*\\ *&*&*\\ *&*&*\end{pmatrix}
=\begin{pmatrix}ei-fh&ch-bi&bf-ce\\ fg-di&ai-cg&cd-af\\ dh-eg&bg-ah&ae-bd\end{pmatrix}$$

Thus, by rescaling, we obtain the formula in the statement.
\end{proof}

In fact, we can now fully solve the inversion problem, as follows:

\index{matrix inversion}

\begin{theorem}
The inverse of a square matrix having nonzero determinant is
$$A^{-1}=\frac{1}{\det A}
\begin{pmatrix}
\det A^{(11)}&-\det A^{(21)}&\det A^{(31)}&\ldots\\
-\det A^{(12)}&\det A^{(22)}&-\det A^{(32)}&\ldots\\
\det A^{(13)}&-\det A^{(23)}&\det A^{(33)}&\ldots\\
\vdots&\vdots&\vdots&
\end{pmatrix}$$
where $A^{(ij)}$ is the matrix $A$, with the $i$-th row and $j$-th column removed.
\end{theorem}

\begin{proof}
This follows indeed by using the row expansion formula from Theorem 9.25, which in terms of the matrix $A^{-1}$ in the statement reads $AA^{-1}=1$.
\end{proof}

Let us discuss now the general formula of the determinant, at arbitrary values $N\in\mathbb N$ of the matrix size, generalizing the formulae that we have at $N=2,3$. We will need:

\index{permutation}
\index{symmetric group}

\begin{definition}
A permutation of $\{1,\ldots,N\}$ is a bijection, as follows:
$$\sigma:\{1,\ldots,N\}\to\{1,\ldots,N\}$$
The set of such permutations is denoted $S_N$.
\end{definition}

There are many possible notations for the permutations, the simplest one consisting in writing the numbers $1,\ldots,N$, and below them, their permuted versions:
$$\sigma=\begin{pmatrix}
1&2&3&4&5\\
2&1&4&5&3
\end{pmatrix}$$

Another method, which is better for most purposes, and faster too, remember that time is money, is by denoting permutations as diagrams, going from top to bottom:
$$\xymatrix@R=3mm@C=3.5mm{
&\ar@{-}[ddr]&\ar@{-}[ddl]&\ar@{-}[ddrr]&\ar@{-}[ddl]&\ar@{-}[ddl]\\
\sigma=\\
&&&&&}$$

There are many interesting things that can be said about permutations, and check here any algebra book. In what concerns us, we will need the following key result:

\index{permutation}
\index{signature}
\index{number of inversions}
\index{transposition}

\begin{theorem}
The permutations have a signature function
$$\varepsilon:S_N\to\{\pm1\}$$
which can be defined in the following equivalent ways:
\begin{enumerate}
\item As $(-1)^c$, where $c$ is the number of inversions.

\item As $(-1)^t$, where $t$ is the number of transpositions.

\item As $(-1)^o$, where $o$ is the number of odd cycles.

\item As $(-1)^x$, where $x$ is the number of crossings.

\item As the sign of the corresponding permuted basis of $\mathbb R^N$.
\end{enumerate}
\end{theorem}

\begin{proof}
Let us begin with the precise definition of $c,t,o,x$, as numbers modulo 2:

\medskip

(1) The idea here is that given any two numbers $i<j$ among $1,\ldots,N$, the permutation  can either keep them in the same order, $\sigma(i)<\sigma(j)$, or invert them:
$$\sigma(j)>\sigma(i)$$

Now by making $i<j$ vary over all pairs of numbers in $1,\ldots,N$, we can count the number of inversions, and call it $c$. This is an integer, $c\in\mathbb N$, which is well-defined.

\medskip

(2) Here the idea, which is something quite intuitive, is that any permutation appears as a product of switches, also called transpositions: 
$$i\leftrightarrow j$$

The decomposition as a product of transpositions is not unique, but the number $t$ of the needed transpositions is unique, when considered modulo 2. This follows for instance from the equivalence of (2) with (1,3,4,5), explained below.

\medskip

(3) Here the point is that any permutation decomposes, in a unique way, as a product of cycles, which are by definition permutations of the following type:
$$i_1\to i_2\to i_3\to\ldots\ldots\to i_k\to i_1$$

Some of these cycles have even length, and some others have odd length. By counting those having odd length, we obtain a well-defined number $o\in\mathbb N$.

\medskip

(4) Here the method is that of drawing the permutation, as we usually do, and by avoiding triple crossings, and then counting the number of crossings. This number $x$ depends on the way we draw the permutations, but modulo 2, we always get the same number. Indeed, this follows from the fact that we can continuously pass from a drawing to each other, and that when doing so, the number of crossings can only jump by $\pm2$.

\medskip

Summarizing, we have 4 different definitions for the signature of the permutations, which all make sense, constructed according to (1-4) above. Regarding now the fact that we always obtain the same number, this can be established as follows:

\medskip

(1)=(2) This is clear, because any transposition inverts once, modulo 2.

\medskip

(1)=(3) This is clear as well, because the odd cycles invert once, modulo 2.

\medskip

(1)=(4) This comes from the fact that the crossings correspond to inversions.

\medskip

(2)=(3) This follows by decomposing the cycles into transpositions.

\medskip

(2)=(4) This comes from the fact that the crossings correspond to transpositions.

\medskip

(3)=(4) This follows by drawing a product of cycles, and counting the crossings.

\medskip

Finally, in what regards the equivalence of all these constructions with (5), here simplest is to use (2). Indeed, we already know that the sign of a system of vectors switches when interchanging two vectors, and so the equivalence between (2,5) is clear. 
\end{proof}

Now back to linear algebra, we can formulate a key result, as follows:

\index{determinant formula}
\index{Sarrus formula}

\begin{theorem}
We have the following formula for the determinant,
$$\det A=\sum_{\sigma\in S_N}\varepsilon(\sigma)A_{1\sigma(1)}\ldots A_{N\sigma(N)}$$
with the signature function being the one introduced above.
\end{theorem}

\begin{proof}
This follows by recurrence over $N\in\mathbb N$, as follows:

\medskip

(1) When developing the determinant over the first column, we obtain a signed sum of $N$ determinants of size $(N-1)\times(N-1)$. But each of these determinants can be computed by developing over the first column too, and so on, and we are led to the conclusion that we have a formula as in the statement, with $\varepsilon(\sigma)\in\{-1,1\}$ being certain coefficients.

\medskip

(2) But these latter coefficients $\varepsilon(\sigma)\in\{-1,1\}$ can only be the signatures of the corresponding permutations $\sigma\in S_N$, with this being something that can be viewed again by recurrence, with either of the definitions (1-5) in Theorem 9.31 for the signature.
\end{proof}

The above result is something quite tricky, and in order to get familiar with it, there is nothing better than doing some computations. As a first, basic example, in 2 dimensions we recover the usual formula of the determinant, the details being as follows:
\begin{eqnarray*}
\begin{vmatrix}a&b\\ c&d\end{vmatrix}
&=&\varepsilon(|\,|)\cdot ad+\varepsilon(\slash\hskip-2mm\backslash)\cdot cb\\
&=&1\cdot ad+(-1)\cdot cb\\
&=&ad-bc
\end{eqnarray*}

In 3 dimensions now, we recover the Sarrus formula:
$$\begin{vmatrix}a&b&c\\ d&e&f\\ g&h&i\end{vmatrix}=aei+bfg+cdh-ceg-bdi-afh$$

Let us record as well the following general formula, which was not obvious with our previous theory of the determinant, but which follows from Theorem 9.32:
$$\det(A)=\det(A^t)$$

There are countless other applications of the formula in Theorem 9.32, and we will be back to this, on several occassions. But, most importantly, that formula allows us to deal now with the complex matrices too, by formulating the following statement:

\index{complex matrix}
\index{deterrminant}

\begin{theorem}
If we define the determinant of a complex matrix $A\in M_N(\mathbb C)$ to be
$$\det A=\sum_{\sigma\in S_N}\varepsilon(\sigma)A_{1\sigma(1)}\ldots A_{N\sigma(N)}$$
then this determinant has the same properties as the determinant of the real matrices.
\end{theorem}

\begin{proof}
This follows by doing some sort of reverse engineering, with respect to what has been done in this section, and we reach to the conclusion that $\det$ has indeed all the good properties that we are familiar with. Except of course for the properties at the very beginning of this section, in relation with volumes, which don't extend well to $\mathbb C^N$.
\end{proof}

Good news, this is the end of the general theory that we wanted to develop. We have now in our bag all the needed techniques for computing the determinant.

\section*{9d. Diagonalization}

Let us discuss now the diagonalization question for linear maps and matrices. The basic diagonalization theory, formulated in terms of matrices, is as follows:

\index{eigenvalue}
\index{eigenvector}
\index{diagonalization}
\index{passage matrix}

\begin{proposition}
A vector $v\in\mathbb C^N$ is called eigenvector of $A\in M_N(\mathbb C)$, with corresponding eigenvalue $\lambda$, when $A$ multiplies by $\lambda$ in the direction of $v$:
$$Av=\lambda v$$
In the case where $\mathbb C^N$ has a basis $v_1,\ldots,v_N$ formed by eigenvectors of $A$, with corresponding eigenvalues $\lambda_1,\ldots,\lambda_N$, in this new basis $A$ becomes diagonal, as follows:
$$A\sim\begin{pmatrix}\lambda_1\\&\ddots\\&&\lambda_N\end{pmatrix}$$
Equivalently, if we denote by $D=diag(\lambda_1,\ldots,\lambda_N)$ the above diagonal matrix, and by $P=[v_1\ldots v_N]$ the square matrix formed by the eigenvectors of $A$, we have:
$$A=PDP^{-1}$$
In this case we say that the matrix $A$ is diagonalizable.
\end{proposition}

\begin{proof}
This is something which is clear, the idea being as follows:

\medskip

(1) The first assertion is clear, because the matrix which multiplies each basis element $v_i$ by a number $\lambda_i$ is precisely the diagonal matrix $D=diag(\lambda_1,\ldots,\lambda_N)$.

\medskip

(2) The second assertion follows from the first one, by changing the basis. We can prove this by a direct computation as well, because we have $Pe_i=v_i$, and so:
$$PDP^{-1}v_i
=PDe_i
=P\lambda_ie_i
=\lambda_iPe_i
=\lambda_iv_i$$

Thus, the matrices $A$ and $PDP^{-1}$ coincide, as stated.
\end{proof}

In order to study the diagonalization problem, the idea is that the eigenvectors can be grouped into linear spaces, called eigenspaces, as follows:

\begin{theorem}
Let $A\in M_N(\mathbb C)$, and for any eigenvalue $\lambda\in\mathbb C$ define the corresponding eigenspace as being the vector space formed by the corresponding eigenvectors:
$$E_\lambda=\left\{v\in\mathbb C^N\Big|Av=\lambda v\right\}$$
These eigenspaces $E_\lambda$ are then in a direct sum position, in the sense that given vectors $v_1\in E_{\lambda_1},\ldots,v_k\in E_{\lambda_k}$ corresponding to different eigenvalues $\lambda_1,\ldots,\lambda_k$, we have:
$$\sum_ic_iv_i=0\implies c_i=0$$
In particular, we have $\sum_\lambda\dim(E_\lambda)\leq N$, with the sum being over all the eigenvalues, and our matrix is diagonalizable precisely when we have equality.
\end{theorem}

\begin{proof}
We prove the first assertion by recurrence on $k\in\mathbb N$. Assume by contradiction that we have a formula as follows, with the scalars $c_1,\ldots,c_k$ being not all zero:
$$c_1v_1+\ldots+c_kv_k=0$$

By dividing by one of these scalars, we can assume that our formula is:
$$v_k=c_1v_1+\ldots+c_{k-1}v_{k-1}$$

Now let us apply $A$ to this vector. On the left we obtain:
$$Av_k
=\lambda_kv_k
=\lambda_kc_1v_1+\ldots+\lambda_kc_{k-1}v_{k-1}$$

On the right we obtain something different, as follows:
\begin{eqnarray*}
A(c_1v_1+\ldots+c_{k-1}v_{k-1})
&=&c_1Av_1+\ldots+c_{k-1}Av_{k-1}\\
&=&c_1\lambda_1v_1+\ldots+c_{k-1}\lambda_{k-1}v_{k-1}
\end{eqnarray*}

We conclude from this that the following equality must hold:
$$\lambda_kc_1v_1+\ldots+\lambda_kc_{k-1}v_{k-1}=c_1\lambda_1v_1+\ldots+c_{k-1}\lambda_{k-1}v_{k-1}$$

On the other hand, we know by recurrence that the vectors $v_1,\ldots,v_{k-1}$ must be linearly independent. Thus, the coefficients must be equal, at right and at left:
$$\lambda_kc_1=c_1\lambda_1$$
$$\vdots$$
$$\lambda_kc_{k-1}=c_{k-1}\lambda_{k-1}$$

Now since at least one of the numbers $c_i$ must be nonzero, from $\lambda_kc_i=c_i\lambda_i$ we obtain $\lambda_k=\lambda_i$, which is a contradiction. Thus our proof by recurrence of the first assertion is complete. As for the second assertion, this follows from the first one.
\end{proof}

In order to reach now to more advanced results, we can use the characteristic polynomial, which appears via the following fundamental result:

\index{characteristic polynomial}

\begin{theorem}
Given a matrix $A\in M_N(\mathbb C)$, consider its characteristic polynomial:
$$P(x)=\det(A-x1_N)$$
The eigenvalues of $A$ are then the roots of $P$. Also, we have the inequality
$$\dim(E_\lambda)\leq m_\lambda$$
where $m_\lambda$ is the multiplicity of $\lambda$, as root of $P$.
\end{theorem}

\begin{proof}
The first assertion follows from the following computation, using the fact that a linear map is bijective when the determinant of the associated matrix is nonzero:
\begin{eqnarray*}
\exists v,Av=\lambda v
&\iff&\exists v,(A-\lambda 1_N)v=0\\
&\iff&\det(A-\lambda 1_N)=0
\end{eqnarray*}

Regarding now the second assertion, given an eigenvalue $\lambda$ of our matrix $A$, consider the dimension $d_\lambda=\dim(E_\lambda)$ of the corresponding eigenspace. By changing the basis of $\mathbb C^N$, as for the eigenspace $E_\lambda$ to be spanned by the first $d_\lambda$ basis elements, our matrix becomes as follows, with $B$ being a certain smaller matrix:
$$A\sim\begin{pmatrix}\lambda 1_{d_\lambda}&0\\0&B\end{pmatrix}$$

We conclude that the characteristic polynomial of $A$ is of the following form:
$$P_A
=P_{\lambda 1_{d_\lambda}}P_B
=(\lambda-x)^{d_\lambda}P_B$$

Thus the multiplicity $m_\lambda$ of our eigenvalue $\lambda$, as a root of $P$, satisfies $m_\lambda\geq d_\lambda$, and this leads to the conclusion in the statement.
\end{proof}

Now recall that we are over $\mathbb C$. We obtain the following result:

\index{eigenvalue}
\index{eigenvector}
\index{characteristic polynomial}
\index{diagonalization}

\begin{theorem}
Given a matrix $A\in M_N(\mathbb C)$, consider its characteristic polynomial
$$P(X)=\det(A-X1_N)$$ 
then factorize this polynomial, by computing the complex roots, with multiplicities,
$$P(X)=(-1)^N(X-\lambda_1)^{n_1}\ldots(X-\lambda_k)^{n_k}$$
and finally compute the corresponding eigenspaces, for each eigenvalue found:
$$E_i=\left\{v\in\mathbb C^N\Big|Av=\lambda_iv\right\}$$
The dimensions of these eigenspaces satisfy then the following inequalities,
$$\dim(E_i)\leq n_i$$
and $A$ is diagonalizable precisely when we have equality for any $i$.
\end{theorem}

\begin{proof}
This follows by combining the above results. Indeed, by summing the inequalities $\dim(E_\lambda)\leq m_\lambda$ from Theorem 9.36, we obtain an inequality as follows:
$$\sum_\lambda\dim(E_\lambda)\leq\sum_\lambda m_\lambda\leq N$$

On the other hand, we know from Theorem 9.35 that our matrix is diagonalizable when we have global equality. Thus, we are led to the conclusion in the statement.
\end{proof}

As an illustration for all this, which is a must-know computation, we have:

\index{rotation}
\index{complex eigenvalues}

\begin{proposition}
The rotation of angle $t\in\mathbb R$ in the plane diagonalizes as:
$$\begin{pmatrix}\cos t&-\sin t\\ \sin t&\cos t\end{pmatrix}
=\frac{1}{2}\begin{pmatrix}1&1\\i&-i\end{pmatrix}
\begin{pmatrix}e^{-it}&0\\0&e^{it}\end{pmatrix}
\begin{pmatrix}1&-i\\1&i\end{pmatrix}$$
Over the reals this is impossible, unless $t=0,\pi$, where the rotation is diagonal.
\end{proposition}

\begin{proof}
Observe first that, as indicated, unlike we are in the case $t=0,\pi$, where our rotation is $\pm1_2$, our rotation is a ``true'' rotation, having no eigenvectors in the plane. Fortunately the complex numbers come to the rescue, via the following computation:
$$\begin{pmatrix}\cos t&-\sin t\\ \sin t&\cos t\end{pmatrix}\binom{1}{i}
=\binom{\cos t-i\sin t}{i\cos t+\sin t}
=e^{-it}\binom{1}{i}$$

We have as well a second complex eigenvector, coming from:
$$\begin{pmatrix}\cos t&-\sin t\\ \sin t&\cos t\end{pmatrix}\binom{1}{-i}
=\binom{\cos t+i\sin t}{-i\cos t+\sin t}
=e^{it}\binom{1}{-i}$$

Thus, we are led to the conclusion in the statement.
\end{proof}

At the level of basic examples of diagonalizable matrices, we first have the following result, which provides us with the ``generic'' examples:

\index{resultant}
\index{discriminant}
\index{diagonalizable matrix}

\begin{theorem}
For a matrix $A\in M_N(\mathbb C)$ the following conditions are equivalent,
\begin{enumerate}
\item The eigenvalues are different, $\lambda_i\neq\lambda_j$,

\item The characteristic polynomial $P$ has simple roots,

\item The characteristic polynomial satisfies $(P,P')=1$,

\item The resultant of $P,P'$ is nonzero, $R(P,P')\neq0$,

\item The discriminant of $P$ is nonzero, $\Delta(P)\neq0$,
\end{enumerate}
and in this case, the matrix is diagonalizable.
\end{theorem}

\begin{proof}
The last assertion holds indeed, due to Theorem 9.37. As for the equivalences in the statement, these are all standard, by using the theory of $R,\Delta$ from chapter 5.
\end{proof}

As already mentioned, one can prove that the matrices having distinct eigenvalues are ``generic'', and so the above result basically captures the whole situation. We have in fact the following collection of density results, which are quite advanced:

\index{diagonalizable matrix}
\index{density}

\begin{theorem}
The following happen, inside $M_N(\mathbb C)$:
\begin{enumerate}
\item The invertible matrices are dense.

\item The matrices having distinct eigenvalues are dense.

\item The diagonalizable matrices are dense.
\end{enumerate}
\end{theorem}

\begin{proof}
These are quite advanced results, which can be proved as follows:

\medskip

(1) This is clear, intuitively speaking, because the invertible matrices are given by the condition $\det A\neq 0$. Thus, the set formed by these matrices appears as the complement of the hypersurface $\det A=0$, and so must be dense inside $M_N(\mathbb C)$, as claimed. 

\medskip

(2) Here we can use a similar argument, this time by saying that the set formed by the matrices having distinct eigenvalues appears as the complement of the hypersurface given by $\Delta(P_A)=0$, and so must be dense inside $M_N(\mathbb C)$, as claimed. 

\medskip

(3) This follows from (2), via the fact that the matrices having distinct eigenvalues are diagonalizable, that we know from Theorem 9.39. There are of course some other proofs as well, for instance by putting the matrix in Jordan form.
\end{proof}

This was for the basic theory of the linear maps, that we will need in what follows. There are of course far more things that can be said, and we will be back to this gradually, when needed, and notably in chapter 12 below, when doing advanced calculus.

\section*{9e. Exercises}

This was an easy chapter, and the start of something new, at least when compared with the horrors from chapters 7-8, and here are some easy exercises for you:

\begin{exercise}
Work out all the details for the main properties of the determinant function $\det:M_N(\mathbb R)\to\mathbb R$, by using our geometric approach, and Thales.
\end{exercise}

\begin{exercise}
Prove that for $H\in M_N(\pm1)$ we have $|\det H|\leq N^{N/2}$, with equality precisely when $H$ is Hadamard, in the sense that its rows are pairwise orthogonal.
\end{exercise}

\begin{exercise}
Work out all the details for the properties of the determinant function $\det:M_N(\mathbb C)\to\mathbb C$, by using our algebraic approach, and permutations.
\end{exercise}

\begin{exercise}
Learn the Jordan form, and come up with an alternative proof for our result stating that the diagonalizable matrices are dense.
\end{exercise}

As an important bonus exercise, diagonalize $3\times3$ matrices, that you can easily find on the internet, as many as needed, as for the computation to take 15 minutes.

\chapter{Partial derivatives}

\section*{10a. Functions, continuity}

With the linear maps $f:\mathbb R^N\to\mathbb R^M$ and $f:\mathbb C^N\to\mathbb C^M$ learned and digested, it remains now to develop the theory of arbitrary maps $f:\mathbb R^N\to\mathbb R^M$ and $f:\mathbb C^N\to\mathbb C^M$, in analogy with what we know at $N=M=1$. Let us start with continuity:

\begin{problem}
What can we say about the continuity of functions 
$$f:\mathbb R^N\to\mathbb R^M\quad,\quad f:\mathbb C^N\to\mathbb C^M$$
in analogy with what we know from before, in $1$ variable?
\end{problem}

In answer to this, it is pretty much clear that things will be straightforward, by using for our analysis purposes the formula of the distance in $\mathbb R^N$ or $\mathbb C^N$, namely:
$$d(x,y)=\sqrt{\sum_{i=1}^N|x_i-y_i|^2}$$

However, in order to avoid using all the time this formula, it is convenient to relax a bit, and take an abstract point of view on all this. So, let us start our study with:

\index{distance}
\index{metric space}

\begin{definition}
A metric space is a set $X$ with a distance function $d:X\times X\to\mathbb R_+$, having the following properties:
\begin{enumerate}
\item $d(x,y)>0$ if $x\neq y$, and $d(x,x)=0$.

\item $d(x,y)=d(y,x)$.

\item $d(x,y)\leq d(x,z)+d(y,z)$.
\end{enumerate}
\end{definition}

And good definition this is, going in fact much further than needed, but with this being a good thing, with the interesting examples abounding, as follows:

\begin{examples}
The following are metric spaces:
\begin{enumerate}
\item $\mathbb R^N,\mathbb C^N$ with the usual distance, and their subsets $X\subset\mathbb R^N,\mathbb C^N$. 

\item $\mathbb R^N,\mathbb C^N$ with the $p$-distance, $p\in[1,\infty]$, and their subsets $X\subset\mathbb R^N,\mathbb C^N$.

\item $l^p(I)$, taken real or complex, with the $p$-distance, and its subsets $X\subset l^p(I)$.

\item $L^p(Y)$ with $Y\subset\mathbb R,\mathbb C$, again taken real or complex, and its subsets $X\subset L^p(Y)$.
\end{enumerate}
\end{examples}

In what follows we will be mainly interested in (1), in order to solve Problem 10.1, and then go ahead with derivatives. However, having (2,3,4) in our theory is certainly a good thing, and for more on these latter spaces, we refer to chapter 7.

\bigskip

As a further example, which is in fact not new, covered by the above, we have:

\begin{theorem}
Given a set $X$, which can be finite or not, the function
$$d(x,y)=\begin{cases}
1&{\rm if}\ x\neq y\\
0&{\rm if}\ x=y
\end{cases}$$
is a metric on it, called discrete metric, and the following happen:
\begin{enumerate}
\item When $|X|<\infty$ this space is the $N$-simplex in $\mathbb R^{N-1}$, with $N=|X|$.

\item More conveniently, we can say that we have $X=\{e_1,\ldots,e_N\}\subset \mathbb R^N$. 

\item In general, we can realize $X$ as being the subset $X=\{e_x\}_{x\in X}\subset l^2(X)$.
\end{enumerate}
\end{theorem}

\begin{proof}
There are several things going on here, the idea being as follows:

\medskip

(1) First of all, the axioms from Definition 10.2 are trivially satisfied, and with the main axiom there, namely the triangle inequality, coming from:
$$1\leq1+1$$

(2) At the level of examples, at $|X|=1$ we obtain a point, at $|X|=2$ we obtain a segment, at $|X|=3$ we obtain an equilateral triangle, at $|X|=4$ we obtain a regular tetrahedron, and so on. Thus, what we have in general, at $|X|=N$, is the arbitrary dimensional generalization of this series of geometric objects, called $N$-simplex.

\medskip

(3) In what regards now geometric realizations, we certainly have $X\subset\mathbb R^{N-1}$, but the computations here are quite complicated, as you can check yourself by studying the problem at $N=4$, that is, by parametrizing the regular tetrahedron in $\mathbb R^3$.

\medskip

(4) However, mathematics, or perhaps physics, come to the rescue, via the idea ``add a dimension, for getting smarter''. Indeed, when looking for an embedding $X\subset\mathbb R^N$ things drastically simplify, because we can simply take $X$ to be the standard basis of $\mathbb R^N$.

\medskip

(5) And finally, this latter interpretation works at $|X|=\infty$ too, as indicated, by using the standard orthonormal basis of the Hilbert space $l^2(X)$, taken real or complex.
\end{proof}

Moving ahead now with some theory, and allowing us a bit of slopiness, we have:

\begin{proposition}
We can talk about limits inside metric spaces $X$, by saying that
$$x_n\to x\iff d(x_n,x)\to0$$
and we can talk as well about continuous functions $f:X\to Y$, by requiring that
$$x_n\to x\implies f(x_n)\to f(x)$$
and with these notions in hand, all the basic results from the cases $X=\mathbb R,\mathbb C$ extend.
\end{proposition}

\begin{proof}
All this is something very standard, based on what we know for $X=\mathbb R,\mathbb C$, and we will leave the various verifications here as an instructive exercise.
\end{proof}

Next, we can talk about open and closed sets inside metric spaces $X$, as follows:

\begin{definition}
Let $X$ be a metric space.
\begin{enumerate}
\item The open balls are the sets $B_x(r)=\{y\in X|d(x,y)<r\}$.

\item The closed balls are the sets $\bar{B}_x(r)=\{y\in X|d(x,y)\leq r\}$.

\item $O\subset X$ is called open if for any $x\in O$ we have a ball $B_x(r)\subset O$.

\item $C\subset X$ is called closed if its complement $C^c\subset X$ is open.
\end{enumerate}
\end{definition}

At the level of basic examples, our notions above coincide with the usual ones, that we know well, in the cases $X=\mathbb R,\mathbb C$. In general now, we first have:

\begin{proposition}
The open balls are open, and the closed balls are closed.
\end{proposition}

\begin{proof}
This might sound like a joke, but it is not one. As for the proof:

\medskip

(1) Given an open ball $B_x(r)$ and a point $y\in B_x(r)$, by using the triangle inequality we have $B_y(r')\subset B_x(r)$, with $r'=r-d(x,y)$. Thus, $B_x(r)$ is indeed open.

\medskip

(2) Given a closed ball $\bar{B}_x(r)$ and a point $y\in B_x(r)^c$, by using the triangle inequality we have $B_y(r')\subset B_x(r)^c$, with $r'=d(x,y)-r$. Thus, $\bar{B}_x(r)$ is indeed closed.
\end{proof}

Here is now something more interesting, making the link with analysis:

\begin{theorem}
For a set $C\subset X$, the following are equivalent:
\begin{enumerate}
\item $C$ is closed in our sense, meaning that $C^c$ is open.

\item We have $x_n\to x,x_n\in C\implies x\in C$.
\end{enumerate}
\end{theorem}

\begin{proof}
We can prove this by double implication, as follows:

\medskip

$(1)\implies(2)$ Assume by contradiction $x_n\to x,x_n\in C$ with $x\notin C$. Since we have $x\in C^c$, which is open, we can pick a ball $B_x(r)\subset C^c$. But this contradicts our convergence assumption $x_n\to x$, so we are done with this implication. 

\medskip

$(2)\implies(1)$ Assume by contradiction that $C$ is not closed in our sense, meaning that $C^c$ is not open. Thus, we can find $x\in C^c$ such that there is no ball $B_x(r)\subset C^c$. But with $r=1/n$ this provides us with a point $x_n\in B_x(1/n)\cap C$, and since we have $x_n\to x$, this contradicts our assumption (2). Thus, we are done with this implication too.
\end{proof}

Here is another basic theorem about open and closed sets:

\begin{theorem}
Let $X$ be a metric space.
\begin{enumerate}
\item If $O_i$ are open, then $\cup_iO_i$ is open.

\item If $C_i$ are closed, then $\cap_iC_i$ is closed.

\item If $O_1,\ldots,O_n$ are open, then $\cap_iO_i$ is open.

\item If $C_1,\ldots,C_n$ are closed, then $\cup_iC_i$ is closed.
\end{enumerate}
\end{theorem}

\begin{proof}
The proof here is identical to the one from chapter 2, for the subsets of $\mathbb R$, and we will leave checking this, that the proof is indeed identical, as an exercise.
\end{proof}

Finally, still in relation with open and closed sets, we have as well:

\begin{definition}
Let $X$ be a metric space, and $E\subset X$ be a subset.
\begin{enumerate}
\item The interior $E^\circ\subset E$ is the set of points $x\in E$ which admit around them open balls $B_x(r)\subset E$.

\item The closure $E\subset\bar{E}$ is the set of points $x\in X$ which appear as limits of sequences $x_n\to x$, with $x\in E$.
\end{enumerate}
\end{definition}

These notions are quite interesting, because they make sense for any set $E$. That is, when $E$ is open, it is open and end of the story, and when $E$ is closed, it is closed and end of the story too. In general, however, a set $E\subset X$ is not open or closed, and what we can best do to it, in order to study it with our tools, is to ``squeeze'' it, as follows:
$$E^\circ\subset E\subset\bar{E}$$

In practice now, in order to use the above notions, we need to know a number of things, including that fact that $E$ open implies $E^\circ=E$, the fact that $E$ closed implies $\bar{E}=E$, and other such results. But all this can be done, and the useful statement here, summarizing all we need to know about interiors and closures, is as follows:

\begin{theorem}
Let $X$ be a metric space, and $E\subset X$ be a subset.
\begin{enumerate}
\item The interior $E^\circ\subset E$ is the biggest open set contained in $E$.

\item The closure $E\subset\bar{E}$ is the smallest closed set containing $E$.
\end{enumerate}
\end{theorem}

\begin{proof}
We have several things to be proved, the idea being as follows:

\medskip

(1) Let us first prove that the interior $E^\circ$ is open. For this purpose, pick $x\in E^\circ$. We know that we have a ball $B_x(r)\subset E$, and since this ball is open, it follows that we have  $B_x(r)\subset E^\circ$. Thus, the interior $E^\circ$ is open, as claimed.

\medskip

(2) Let us prove now that the closure $\bar{E}$ is closed. For this purpose, we will prove that the complement $\bar{E}^c$ is open. So, pick $x\in\bar{E}^c$. Then $x$ cannot appear as a limit of a sequence $x_n\to x$ with $x_n\in E$, so we have a ball $B_x(r)\subset \bar{E}^c$, as desired.

\medskip

(3) Finally, the maximality and minimality assertions regarding $E^\circ$ and $\bar{E}$ are both routine, coming from definitions, and we will leave them as exercises.
\end{proof}

As an application of the theory developed above, and more specifically of the notion of closure from Definition 10.10, we can talk as well about density, as follows:

\begin{definition}
We say that a subset $E\subset X$ is dense when:
$$\bar{E}=X$$
That is, any point of $X$ must appear as a limit of points of $E$.
\end{definition}

Moving on, we can talk as well about compact sets, and connected sets. But here, surprise,  things are more tricky, in the general metric space framework, deviating a bit from what we know. So, we will explain this in detail. Let us start with:

\begin{definition}
A set $K\subset X$ is called compact if any cover with open sets
$$K\subset\bigcup_iO_i$$
has a finite subcover, $K\subset(O_{i_1}\cup\ldots\cup O_{i_n})$.
\end{definition}

This definition might seem overly abstract, but our claim is that this is the correct definition, and that there is no way of doing otherwise. Indeed, we have:

\begin{proposition}
Given an infinite set $X$ with the discrete distance on it, namely $d(p,q)=1-\delta_{pq}$, this can be modeled as the basis of a suitable Hilbert space,
$$X=\{e_x\}_{x\in X}\subset l^2(X)$$
and this set is closed and bounded, but not compact.
\end{proposition}

\begin{proof}
Here the first part, regarding the explicit modeling of $X$, is something that we already know, from Theorem 10.4. Regarding now the second part:

\medskip

(1) $X$ being the total space, it is by definition closed.

\medskip

(2) $X$ is also bounded, because all distances are smaller than $1$.

\medskip

(3) However, $X$ is not compact, because its points being open, containing balls of radius $1/2$ around them, $X=\cup_{x\in X}\{x\}$ is an open cover, having no finite subcover.
\end{proof}

Quite interesting all this, obviously, infinite dimensions are no joke. Next, we have the following result, confirming that we are on the good way, with Definition 10.13:

\begin{proposition}
The following hold:
\begin{enumerate}
\item Compact implies closed.

\item Closed inside compact is compact.

\item Compact intersected with closed is compact.
\end{enumerate}
\end{proposition}

\begin{proof}
These assertions are all clear from definitions, as follows:

\medskip

(1) Assume that $K\subset X$ is compact, and let us prove that $K$ is closed. For this purpose, we will prove that $K^c$ is open. So, pick $p\in K^c$. For any $q\in K$ we set $r=d(p,q)/3$, and we consider the following balls, separating $p$ and $q$:
$$U_q=B_p(r)\quad,\quad V_q=B_q(r)$$

We have then $K\subset\cup_{q\in K}V_q$, so we can pick a finite subcover, as follows:
$$K\subset\left(V_{q_1}\cup\ldots\cup V_{q_n}\right)$$

With this done, consider now the following intersection:
$$U=U_{q_1}\cap\ldots\cap U_{q_n}$$

This intersection is then a ball around $p$, and since this ball avoids $V_{q_1},\ldots,V_{q_n}$, it avoids the whole $K$. Thus, we have proved that $K^c$ is open at $p$, as desired.

\medskip

(2) Assume indeed that $C\subset K$ is closed, with $K\subset X$ being compact. For proving our result, we can assume, by replacing $X$ with $K$, that we have $X=K$. In order to prove now that $C$ is compact, consider an open cover of it, as follows:
$$C\subset\bigcup_iO_i$$

By adding the set $C^c$, which is open, to this cover, we obtain a cover of $K$. Now since $K$ is compact, we can extract from this a finite subcover $\Omega$, and there are two cases:

\medskip

-- If $C^c\in\Omega$, by removing $C^c$ from $\Omega$ we obtain a finite cover of $C$, as desired.

\medskip

-- If $C^c\notin\Omega$, we are done too, because in this case $\Omega$ is a finite cover of $C$.

\medskip

(3) This follows from (1) and (2), because if $K\subset X$ is compact, and $C\subset X$ is closed, then $K\cap C\subset K$ is closed inside a compact, so it is compact.
\end{proof}

As a second batch of results, which are useful as well, we have:

\begin{proposition}
The following hold:
\begin{enumerate}
\item If $K_i\subset X$ are compact, satisfying $K_{i_1}\cap\ldots\cap K_{i_n}\neq\emptyset$, then $\cap_iK_i\neq\emptyset$.

\item If $K_1\supset K_2\supset K_3\supset\ldots$ are non-empty compacts, then $\cap_iK_i\neq\emptyset$.

\item If $K$ is compact, and $E\subset K$ is infinite, then $E$ has a limit point in $K$.

\item If $K$ is compact, any sequence $\{x_n\}\subset K$ has a limit point in $K$.

\item If $K$ is compact, any $\{x_n\}\subset K$ has a subsequence which converges in $K$.
\end{enumerate}
\end{proposition}

\begin{proof}
Again, these are elementary results, which can be proved as follows:

\medskip

(1) Assume by contradiction $\cap_iK_i=\emptyset$, and let us pick $K_1\in\{K_i\}$. Since any $x\in K_1$ is not in $\cap_iK_i$, there is an index $i$ such that $x\in K_i^c$, and we conclude that we have:
$$K_1\subset\bigcup_{i\neq 1}K_i^c$$

But this can be regarded as being an open cover of $K_1$, that we know to be compact, so we can extract from it a finite subcover, as follows:
$$K_1\subset\left(K_{i_1}^c\cup\ldots\cup K_{i_n}^c\right)$$

Now observe that this latter subcover tells us that we have:
$$K_1\cap K_{i_1}\cap\ldots\cap K_{i_n}=\emptyset$$

But this contradicts our intersection assumption in the statement, and we are done.

\medskip

(2) This is a particular case of (1), proved above.

\medskip

(3) We prove this by contradiction. So, assume that $E$ has no limit point in $K$. This means that any $p\in K$ can be isolated from the rest of $E$ by a certain open ball $V_p=B_p(r)$, and in both the cases that can appear, $p\in E$ or $p\notin E$, we have:
$$|V_p\cap E|=0,1$$

Now observe that these sets $V_p$ form an open cover of $K$, and so of $E$. But due to $|V_p\cap E|=0,1$ and to $|E|=\infty$, this open cover of $E$ has no finite subcover. Thus the same cover, regarded now as cover of $K$, has no finite subcover either, contradiction.

\medskip

(4) This follows from (3) that we just proved, with $E=\{x_n\}$.

\medskip

(5) This is a reformulation of (4), that we just proved.
\end{proof}

Finally, in the case that we are mostly interested in, namely $X=\mathbb R^N$, we have:

\index{compact set}
\index{closed and bounded}

\begin{theorem}
For a subset $K\subset\mathbb R^N$, the following are equivalent:
\begin{enumerate}
\item $K$ is closed and bounded.

\item $K$ is compact.

\item Any infinite subset $E\subset K$ has a limiting point in $K$.
\end{enumerate}
\end{theorem}

\begin{proof}
This is something quite tricky, the idea being as follows:

\medskip

$(1)\implies(2)$ This is something that we know from chapter 2, in the case $N=1$, and the proof there extends in a straightforward way, by replacing intervals with cubes.

\medskip

$(2)\implies(3)$ This is something that we already know, not needing $K\subset\mathbb R^N$.

\medskip

$(3)\implies(1)$ We have to prove that $K$ as in the statement is both closed and bounded, and we can do both these things by contradiction, as follows:

\medskip

-- Assume first that $K$ is not closed. But this means that we can find a point $x\notin K$ which is a limiting point of $K$. Now let us pick $x_n\in K$, with $x_n\to x$, and consider the set $E=\{x_n\}$. According to our assumption, $E$ must have a limiting point in $K$. But this limiting point can only be $x$, which is not in $K$, contradiction.

\medskip

-- Assume now that $K$ is not bounded. But this means that we can find points $x_n\in K$ satisfying $||x_n||\to\infty$, and if we consider the set $E=\{x_n\}$, then again this set must have a limiting point in $K$, which is impossible, so we have our contradiction, as desired.
\end{proof}

So long for compactness. As a last piece of general topology, we have:

\index{connected set}

\begin{definition}
We can talk about connected sets $E\subset X$, as follows:
\begin{enumerate}
\item We say that $E$ is connected if it cannot be separated as $E=E_1\cup E_2$, with the components $E_1,E_2$ satisfying $E_1\cap\bar{E}_2=\bar{E}_1\cap E_2=\emptyset$.

\item We say that $E$ is path connected if any two points $p,q\in E$ can be joined by a path, meaning a continuous $f:[0,1]\to X$, with $f(0)=p$, $f(1)=q$.
\end{enumerate}
\end{definition}

All this looks a bit technical, and indeed it is. To start with, (1) is something quite natural, but the separation condition there $E_1\cap\bar{E}_2=\bar{E}_1\cap E_2=\emptyset$ can be weakened into $E_1\cap E_2=\emptyset$, or strengthened into $\bar{E}_1\cap\bar{E}_2=\emptyset$, depending on purposes, and with our (1) as formulated above being the good compromise, for most purposes. 

\bigskip

As for (2), this condition is obviously something stronger, and we have in fact the following implications, which are both clear from definitions:
$${\rm convex}\implies{\rm path\ connected}\implies {\rm connected}$$

Anyway, leaving aside the precise axiomatic discussion here, which can be something quite technical, once all these questions clarified, the idea is that any set $E$ can be written as a disjoint union of connected components, as follows:
$$E=\bigsqcup_i E_i$$

However, the story is not over here, because when looking at the connected components $E_i$, or simply at the connected sets $E$, if you prefer, there are many things that can happen, in relation with the ``holes'' that $E$ can have or not. Thus, the classification of connected sets runs into the question of deciding how the holes of such a set $E$ can look like, and this is something quite subtle, that we will discuss in chapter 11 below.

\bigskip

Getting back now to more concrete things, functions and analysis, we have:

\index{continuous function}

\begin{theorem}
Assuming that $f:X\to Y$ is continuous, the following happen,
\begin{enumerate}
\item If $O$ is open, then $f^{-1}(O)$ is open. 

\item If $C$ is closed, then $f^{-1}(C)$ is closed. 

\item If $K$ is compact, then $f(K)$ is compact. 

\item If $E$ is connected, then $f(E)$ is connected. 
\end{enumerate}
and with $(3,4)$ standing as an abstract intermediate value theorem.
\end{theorem}

\begin{proof}
This is something fundamental, which can be proved as follows:

\medskip

(1) This is clear from the definition of continuity, written with $\varepsilon,\delta$. In fact, the converse holds too, in the sense that if $f^{-1}({\rm open})={\rm open}$, then $f$ must be continuous.

\medskip

(2) This follows from (1), by taking complements. And again, the converse holds too, in the sense that if $f^{-1}({\rm closed})={\rm closed}$, then $f$ must be continuous.

\medskip

(3) Given an open cover $f(K)\subset\cup_iO_i$, we get by (1) an open cover $K\subset\cup_if^{-1}(O_i)$, and so by compactness of $K$, a finite subcover $K\subset f^{-1}(O_{i_1})\cup\ldots\cup f^{-1}(O_{i_n})$. But this gives in turn a finite subcover $f(K)\subset O_{i_1}\cup\ldots\cup O_{i_n}$, as desired. 

\medskip

(4) This comes via the same trick as for (3). Indeed, any separation of $f(E)$ into two parts can be returned via $f^{-1}$ into a separation of $E$ into two parts, contradiction.
\end{proof}

\section*{10b. Partial derivatives}

Let us discuss now differentiation in several variables. At order 1, the situation is quite similar to the one in 1 variable, but this time involving matrices. In order to explain this material, let us start with a straightforward definition, as follows:

\index{differentiable function}
\index{derivative}

\begin{definition}
We say that a map $f:\mathbb R^N\to\mathbb R^M$ is differentiable at $x\in\mathbb R^N$ if
$$f(x+t)\simeq f(x)+f'(x)t$$
for some linear map $f'(x):\mathbb R^N\to\mathbb R^M$, called derivative of $f$ at the point $x\in\mathbb R^N$.
\end{definition}

But is this the correct definition. I can hear you thinking that we are probably going the wrong way, because for functions $f:\mathbb R\to\mathbb R$ the derivative is something much simpler, as follows, and that we should try to imitate, in our higher dimensional setting:
$$f'(x)=\lim_{t\to0}\frac{f(x+t)-f(x)}{t}$$

However, this is not possible, for a number of reasons, that are worth discussing in detail. So, here is the discussion, answering all kinds of questions that you might have:

\bigskip

(1) First of all, the above formula does not make any sense for a function $f:\mathbb R^N\to\mathbb R^M$ with $N\neq M$, because we cannot divide oranges by apples. And it doesn't make sense either at $N=M\in\mathbb N$, because here we have $\mathbb R^N$ oranges, I agree with you, but there is no way of dividing such oranges by oranges, unless we are in the special case $N=1$. 

\bigskip

(2) Next, you might say that at $N=2$ we have $\mathbb R^2\simeq\mathbb C$ oranges, which can be divided by themselves, as we saw in chapter 6, when talking about holomorphic functions. Good point, but as shown in chapter 8, there are many interesting functions $f:\mathbb R^2\to\mathbb R$ which are not holomorphic, when regarded as functions $f:\mathbb C\to\mathbb R$. So, $N=2$ remains open.

\bigskip

(3) More philosophically, now that we know that having $f'(x)$ defined as a number is difficult, the question is, do we really want to have $f'(x)$ defined as a number? And my claim here is that, this would be a pity. Think at the case where $f:\mathbb R^N\to\mathbb R^M$ is linear. Such a map is just ``perfect'', and so should equal its own derivative, $f=f'$.

\bigskip

(4) Summarizing, our Definition 10.20 is just perfection, and is waiting for some further study, and this is what we will do. And in case you are still secretly dreaming about having $f'(x)$ defined as some sort of number, wait for it. When $N=M$ at least, there is indeed a lucky number, namely $\det(f'(x))$, called Jacobian, but more on this later.

\bigskip

Getting back now to Definition 10.20 as formulated, and agreed upon, we have there a linear map $f'(x):\mathbb R^N\to\mathbb R^M$, waiting to be further understood. So, time now to use our linear algebra knowledge from chapter 9. We know from there that such linear maps correspond to rectangular matrices $A\in M_{M\times N}(\mathbb R)$, and we are led in this way to:

\begin{question}
Given a differentiable map $f:\mathbb R^N\to\mathbb R^M$, in the abstract sense of Definition 10.20, what exactly is its derivative 
$$f'(x):\mathbb R^N\to\mathbb R^M$$
regarded as a rectangular matrix, $f'(x)\in M_{M\times N}(\mathbb R)$?
\end{question}

Which must sound a bit puzzling, because we just agreed, following a long discussion, on the fact that the derivative is a linear map, and not a number, and now what, we are trying to replace this linear map by a matrix, and so by a bunch of numbers. 

\bigskip

Good point, and in answer, these are quite deep things, that took mankind several centuries to develop, and that we are now presenting in a compressed form. So yes, all this is difficult mathematics, when you first see it, I perfectly agree with you.

\bigskip

Moving on, in order to answer Question 10.21, we will need the following key notion, that we already met, earlier in this book, in chapter 8, when doing physics:

\index{partial derivatives}

\begin{definition}
We can talk about the partial derivatives of $f:\mathbb R^N\to\mathbb R$,
$$\frac{df}{dx_1}(x)=\lim_{t\to0}\frac{f(x_1+t,x_2,\ldots,x_N)-f(x_1,x_2,\ldots,x_N)}{t}$$
$$\frac{df}{dx_2}(x)=\lim_{t\to0}\frac{f(x_1,x_2+t,\ldots,x_N)-f(x_1,x_2,\ldots,x_N)}{t}$$
$$\vdots$$
$$\frac{df}{dx_N}(x)=\lim_{t\to0}\frac{f(x_1,x_2,\ldots,x_N+t)-f(x_1,x_2,\ldots,x_N)}{t}$$
provided that these derivatives exist indeed.
\end{definition}

In other words, the partial derivatives are something very simple, belonging to 1-variable calculus, computed by focusing on one variable, and keeping the other variables fixed. Here is an illustration for this, making it clear that all this is elementary:
$$\frac{d}{dx}(x^5+4e^y+2\cos(x^2y))=5x^4-4xy\sin(x^2y)$$
$$\frac{d}{dy}(x^5+4e^y+2\cos(x^2y))=4e^y-2x^2\sin(x^2y)$$

As an interesting theoretical remark, however, while the partial derivatives only need 1-variable calculus, in order to be defined and computed, their meaning in several variables can be quite different from their usual meaning in 1 variable, as shown by:

\begin{warning}
The following two-variable function $f:\mathbb R^2\to\mathbb R$,
$$f(x,y)=\begin{cases}
0&{\rm if}\ xy=0\\
1&{\rm if}\ xy\neq0
\end{cases}$$
has partial derivatives at zero $\frac{df}{dx}(0,0)=\frac{df}{dy}(0,0)=0$, but is not continuous at $(0,0)$.
\end{warning}

With this discussed, time now to answer Question 10.21? I would say that, barring pathologies like those above, by putting together all partial derivatives, of all components $f_i:\mathbb R^N\to\mathbb R$ of our function $f:\mathbb R^N\to\mathbb R^M$, we will have a win. Indeed, we first have:

\begin{theorem}
The derivative of a differentiable function $f:\mathbb R^N\to\mathbb R^M$, making
$$f(x+t)\simeq f(x)+f'(x)t$$
work, must be the matrix of partial derivatives at $x$, namely
$$f'(x)=\left(\frac{df_i}{dx_j}(x)\right)_{ij}\in M_{M\times N}(\mathbb R)$$ 
acting on the vectors $t\in\mathbb R^N$, written as usual vertically, by usual multiplication.
\end{theorem}

\begin{proof}
Assume indeed that $f:\mathbb R^N\to\mathbb R^M$ is differentiable in the sense of Definition 10.20, with this meaning that we have the following approximation formula:
$$f(x+t)\simeq f(x)+f'(x)t$$

With $f=(f_i)$, with $f_i:\mathbb R^N\to\mathbb R$ being the components, this formula reads: 
$$f_i(x+t)\simeq f_i(x)+(f'(x)t)_i$$

Now let us see what this latter formula teaches us. With the variable $t\simeq 0$ chosen in the $Ox_j$ direction, that is, with $t_k=0$ for $k\neq j$, this formula reads:
$$f_i(x+t)\simeq f_i(x)+f'(x)_{ij}t_j$$

In other words, the number $f'(x)_{ij}$ must appear as the following limit:
\begin{eqnarray*}
f'(x)_{ij}
&=&\lim_{t_j\to0}\frac{f_i(x+t)-f_i(x)}{t_j}\\
&=&\lim_{t_j\to0}\frac{f_i(x_1,\ldots,x_j+t_j,\ldots,x_N)-f_i(x_1,\ldots,x_j,\ldots,x_N)}{t}\\
&=&\frac{df_i}{dx_j}(x)
\end{eqnarray*}

We are therefore led to the conclusion in the statement.
\end{proof}

Before going further, observe that there are in fact some bugs with what we have so far, coming from a confusion between vertical and horizontal vectors. As the saying goes, errare humanum est, perseverare diabolicum, so let us make the following convention:

\begin{convention}
When doing multivariable calculus it is better not to mess up vertical and horizontal vectors, as we did above, and to use vertical vectors only,
$$\begin{pmatrix}x_1\\ \vdots\\ x_N\end{pmatrix}\to f\begin{pmatrix}x_1\\ \vdots\\ x_N\end{pmatrix}=\begin{pmatrix}y_1\\ \vdots\\ y_M\end{pmatrix}$$
as for the linear algebra to work fine. In short, horizontal vectors now banned, and everyone not complying will be sent to a Siberian labor camp, for cutting wood there.
\end{convention}

As a second comment, we have our usual remark, as in the 1-variable case, that the notion of differentiability, be that abstract, as in Definition 10.20, or in relation with the partial derivatives appearing in Theorem 10.24, is local, only requiring $f:\mathbb R^N\to\mathbb R^M$ to be defined on a small ball around the point $x\in\mathbb R^N$ in question. Thus, for being more reliable with our math, we only need to assume $f:X\to\mathbb R^M$, with $X\subset\mathbb R^N$ open.

\bigskip

So, here is a remake of Theorem 10.24, taking into account all this, and also, importantly, talking about a converse too, under suitable continuity assumptions:

\index{continuously differentiable}

\begin{theorem}
A function $f:X\to\mathbb R^M$, with $X\subset\mathbb R^N$ open, is continuously differentiable, in the sense that the following approximation formula holds,
$$f(x+t)\simeq f(x)+f'(x)t$$
and with the derivatives $x\to f'(x)$ being continuous, precisely when it has partial derivatives, which are continuous. In this case the derivative is
$$f'(x)=\left(\frac{df_i}{dx_j}(x)\right)_{ij}\in M_{M\times N}(\mathbb R)$$ 
acting on the vectors $t\in\mathbb R^N$, written as usual vertically, by usual multiplication.
\end{theorem}

\begin{proof}
This is something very standard, the idea being as follows:

\medskip

(1) As a first observation, the approximation in the statement makes sense indeed, as an equality, or rather approximation, of vectors in $\mathbb R^M$, as follows:
$$f\begin{pmatrix}x_1+t_1\\ \vdots\\ x_N+t_N\end{pmatrix}
\simeq f\begin{pmatrix}x_1\\ \vdots\\ x_N\end{pmatrix}
+\begin{pmatrix}
\frac{df_1}{dx_1}(x)&\ldots&\frac{df_1}{dx_N}(x)\\
\vdots&&\vdots\\
\frac{df_M}{dx_1}(x)&\ldots&\frac{df_M}{dx_N}(x)
\end{pmatrix}\begin{pmatrix}t_1\\ \vdots\\ t_N\end{pmatrix}$$

(2) Observe also that at $N=M=1$ what we have is a usual 1-variable function $f:\mathbb R\to\mathbb R$, and the above formula is something that we know well, namely:
$$f(x+t)\simeq f(x)+f'(x)t$$

Thus, our statement is a straightforward generalization of this latter 1-variable formula, and with a continuity assumption added, in order to avoid certain pathologies which are multivariable-specific, coming from functions as those from Warning 10.23.

\medskip

(3) Getting now to the proof, in general, assume that $f:X\to\mathbb R^M$ is continuously differentiable, with this meaning by definition that $f$ is differentiable, with the derivatives $x\to f'(x)$ being continuous with respect to the usual norm of matrices, namely:
$$||A||=\sup_{||x||=1}||Ax||$$

We know from Theorem 10.24 that $f'(x)$ must be the matrix of the partial derivatives at $x$, so these partial derivatives exist indeed. Also, for any $x,y\in X$ we have:
$$f'(x)_{ij}-f'(y)_{ij}=\frac{df_i}{dx_j}(x)-\frac{df_i}{dx_j}(y)$$

By applying now the absolute value, we obtain the following estimate, with the norm at the end being the usual norm of the square matrices, given above:
\begin{eqnarray*}
\left|\frac{df_i}{dx_j}(x)-\frac{df_i}{dx_j}(y)\right|
&=&\left|f'(x)_{ij}-f'(y)_{ij}\right|\\
&=&\left|(f'(x)-f'(y))_{ij}\right|\\
&\leq&||f'(x)-f'(y)||
\end{eqnarray*}

Here we have used at the end the following estimate, coming from Cauchy-Schwarz:
$$|A_{ij}|=|<Ae_j,e_i>|\leq||Ae_j||\cdot||e_i||\leq||A||$$

Thus, the partial derivatives must indeed exist and be continuous, as claimed.

\medskip

(4) Summarizing, one implication proved, and for the other implication, in view of Theorem 10.24, it remains to show that if $f:X\to\mathbb R^M$, with $X\subset\mathbb R^N$ open, has partial derivatives, which are continuous, then the approximation formula in (1) holds indeed.

\medskip

(5) Let us first discuss the case $N=2,M=1$. Here what we have is a function $f:\mathbb R^2\to\mathbb R$, and by using twice the basic approximation result from (2), we obtain:
\begin{eqnarray*}
f\binom{x_1+t_1}{x_2+t_2}
&\simeq&f\binom{x_1+t_1}{x_2}+\frac{df}{dx_2}(x)t_2\\
&\simeq&f\binom{x_1}{x_2}+\frac{df}{dx_1}(x)t_1+\frac{df}{dx_2}(x)t_2\\
&=&f\binom{x_1}{x_2}+\begin{pmatrix}\frac{df}{dx_1}(x)&\frac{df}{dx_2}(x)\end{pmatrix}\binom{t_1}{t_2}
\end{eqnarray*}

That is, for being fully formal, we can use the mean value theorem applied twice, in the obvious way, for each of the above estimates, and we obtain the result.

\medskip

(6) More generally now, we can deal in the same way with the general case $M=1$, with the formula here, obtained via a straightforward recurrence, being as follows:
\begin{eqnarray*}
f\begin{pmatrix}x_1+t_1\\ \vdots\\ x_N+t_N\end{pmatrix}
&\simeq&f\begin{pmatrix}x_1+t_1\\ \vdots\\ x_{N-1}+t_{N-1}\\ x_N\end{pmatrix}+\frac{df}{dx_N}(x)t_N\\
&\simeq&f\begin{pmatrix}x_1+t_1\\ \vdots\\ x_{N-2}+t_{N-2}\\ x_{N-1}\\ x_N\end{pmatrix}+\frac{df}{dx_{N-1}}(x)t_{N-1}+\frac{df}{dx_N}(x)t_N\\
&\vdots&\\
&\simeq&f\begin{pmatrix}x_1\\ \vdots\\ x_N\end{pmatrix}+\frac{df}{dx_1}(x)t_1+\ldots+\frac{df}{dx_N}(x)t_N\\
&=&f\begin{pmatrix}x_1\\ \vdots\\ x_N\end{pmatrix}+
\begin{pmatrix}\frac{df}{dx_1}(x)&\ldots&\frac{df}{dx_N}(x)\end{pmatrix}
\begin{pmatrix}t_1\\ \vdots\\ t_N\end{pmatrix}
\end{eqnarray*}

(7) But this gives the result in the case where both $N,M\in\mathbb N$ are arbitrary too. Indeed, consider a function $f:\mathbb R^N\to\mathbb R^M$, and let us write it as follows:
$$f=\begin{pmatrix}f_1\\ \vdots\\ f_M\end{pmatrix}$$

We can then apply (6) to each of the components $f_i:\mathbb R^N\to\mathbb R$, and we get:
$$f_i\begin{pmatrix}x_1+t_1\\ \vdots\\ x_N+t_N\end{pmatrix}
\simeq f_i\begin{pmatrix}x_1\\ \vdots\\ x_N\end{pmatrix}+
\begin{pmatrix}\frac{df_i}{dx_1}(x)&\ldots&\frac{df_i}{dx_N}(x)\end{pmatrix}
\begin{pmatrix}t_1\\ \vdots\\ t_N\end{pmatrix}$$

(8) But this collection of $M$ formulae tells us precisely that the following happens, as an equality, or rather approximation, of vectors in $\mathbb R^M$:
$$f\begin{pmatrix}x_1+t_1\\ \vdots\\ x_N+t_N\end{pmatrix}
\simeq f\begin{pmatrix}x_1\\ \vdots\\ x_N\end{pmatrix}
+\begin{pmatrix}
\frac{df_1}{dx_1}(x)&\ldots&\frac{df_1}{dx_N}(x)\\
\vdots&&\vdots\\
\frac{df_M}{dx_1}(x)&\ldots&\frac{df_M}{dx_N}(x)
\end{pmatrix}\begin{pmatrix}t_1\\ \vdots\\ t_N\end{pmatrix}$$

Thus, we are led to the conclusion in the statement.
\end{proof}

\section*{10c. Gradient, chain rule} 

We have now a theory of first order derivatives up and working, in the multivariable function setting, $f:\mathbb R^N\to\mathbb R^M$. However, before further developing this theory, it is worth studying some more the basics. We first have the following observation:

\begin{proposition}
Our main result, stating that for $f:\mathbb R^N\to\mathbb R^M$ we have
$$f(x+t)\simeq f(x)+f'(x)t$$
under the assumption that $f$ is continuously differentiable, with
$$f'(x)=\left(\frac{df_i}{dx_j}(x)\right)_{ij}\in M_{M\times N}(\mathbb R)$$
is more or less equivalent to its $M=1$ particular case, that when $f:\mathbb R^N\to\mathbb R$.
\end{proposition}

\begin{proof}
This is something that we used many times, in the proofs of Theorems 10.24 and 10.26, the idea being that we can always write our function as follows, with each of the components $f_1,\ldots,f_M$ being a scalar-valued function, $f_i:\mathbb R^N\to\mathbb R$:
$$f=\begin{pmatrix}f_1\\ \vdots\\ f_M\end{pmatrix}$$

With this convention, the derivative, as a rectangular matrix, is given by:
$$f'(x)=\begin{pmatrix}f_1'(x)\\ \vdots\\ f_M'(x)\end{pmatrix}$$

But now, it is clear that the formula $f(x+t)\simeq f(x)+f'(x)t$ works precisely when $f_i(x+t)\simeq f_i(x)+f'_i(x)t$ works for any $i$, so the case $M=1$ is indeed what matters. 
\end{proof}

Next, regarding the case $M=1$, we can reformulate here Theorem 10.26 in the following alternative way, making the link with the physics from chapter 8:

\begin{theorem}
The derivative of a differentiable function $f:\mathbb R^N\to\mathbb R$, making the approximation formula $f(x+t)\simeq f(x)+f'(x)t$ work, is given by 
$$f'(x)=(\nabla f(x))^t$$
where $\nabla f$, called gradient of $f$, is the column vector formed by the partial derivatives:
$$\nabla f(x)=\left(\frac{df}{dx_j}(x)\right)_j\in\mathbb R^N$$ 
Equivalently, without reference to $f'$, we have the approximation formula
$$f(x+t)\simeq f(x)+<\nabla f(x),t>$$
where $<\,,>$ is the usual scalar product on $\mathbb R^N$.
\end{theorem}

\begin{proof}
This is a reformulation of Theorem 10.26, in the particular case $M=1$. Indeed, we know from there that the derivative of our function $f:\mathbb R^N\to\mathbb R$, making the approximation formula $f(x+t)\simeq f(x)+f'(x)t$ work, is given by:
$$f'(x)=\left(\frac{df}{dx_j}(x)\right)_{1j}\in M_{1\times N}(\mathbb R)$$ 

But, what we have here is the horizontal vector formed by the partial derivatives at $x$, which is transpose to the vertical vector formed by the same partial derivatives, which is by definition $\nabla f(x)$. Thus, we are led to the formula in the statement, namely:
$$f'(x)=(\nabla f(x))^t$$

Finally, regarding the last assertion, this follows from this, as follows:
\begin{eqnarray*}
f(x+t)
&\simeq&f(x)+f'(x)t\\
&=&f(x)+(\nabla f(x))^tt\\
&=&f(x)+\sum_j(\nabla f(x))_jt_j\\
&=&f(x)+<\nabla f(x),t>
\end{eqnarray*}

Thus, we are led to the conclusions in the statement.
\end{proof}

Back now to the general case, that of the functions $f:\mathbb R^N\to\mathbb R^M$, besides $M=1$, another interesting situation appears when $M=N$. Indeed, here the derivative is a square matrix, which makes it tempting to look at the determinant of this matrix:

\index{Jacobian}
\index{volume inflation}

\begin{definition}
Given a differentiable function $f:\mathbb R^N\to\mathbb R^N$, its Jacobian is
$$J_f(x)=\det(f'(x))\in\mathbb R$$
measuring the infinitesimal rate of the volume inflation by $f$, at the given point $x$.
\end{definition}

To be more precise, the interpretation of the Jacobian as an infinitesimal rate of the volume inflation comes from our fine knowledge of the determinant, from chapter 9. Which is quite nice, and as a first remark, at $N=1$ we obtain of course the usual derivative. In general, according to our matrix formula for $f'(x)$, the Jacobian is given by:
$$J_f(x)=\det\begin{pmatrix}
\frac{df_1}{dx_1}(x)&\ldots&\frac{df_1}{dx_N}(x)\\
\vdots&&\vdots\\
\frac{df_N}{dx_1}(x)&\ldots&\frac{df_N}{dx_N}(x)
\end{pmatrix}$$

Thus, the Jacobian can be explicitly computed. However, in what regards the practical uses of the Jacobian, these are quite complicated, and this will have to wait a bit, until chapter 13 below. So, sorry for this, not yet time to enjoy Definition 10.29, and stay with me, plenty of further linear algebra, and matrices instead of numbers, to follow. 

\bigskip

Generally speaking, Theorem 10.26 is what you need to know for upgrading from calculus to multivariable calculus, of course with some help from linear algebra as well. As a first standard result here, which is something very useful, in practice, we have:

\index{chain rule}

\begin{theorem}
We have the chain derivative formula
$$(f\circ g)'(x)=f'(g(x))\cdot g'(x)$$
as an equality of matrices.
\end{theorem}

\begin{proof}
This is something that we know well in one variable, and in several variables the proof is similar, by using the notion of derivative coming from Theorem 10.26. To be more precise, consider a composition of functions, as follows:
$$f:\mathbb R^N\to\mathbb R^M\quad,\quad 
g:\mathbb R^K\to\mathbb R^N\quad,\quad 
f\circ g:\mathbb R^K\to\mathbb R^M$$

According to Theorem 10.26, the derivatives of these functions are certain linear maps, corresponding to certain rectangular matrices, as follows:
$$f'(g(x))\in M_{M\times N}(\mathbb R)\quad,\quad 
g'(x)\in M_{N\times K}(\mathbb R)\quad\quad
(f\circ g)'(x)\in M_{M\times K}(\mathbb R)$$

Thus, our formula makes sense indeed. As for proof, this comes from:
\begin{eqnarray*}
(f\circ g)(x+t)
&=&f(g(x+t))\\
&\simeq&f(g(x)+g'(x)t)\\
&\simeq&f(g(x))+f'(g(x))g'(x)t
\end{eqnarray*}

Thus, we are led to the conclusion in the statement.
\end{proof}

As an alternative formulation of the chain rule, which is useful in practice, we have:

\begin{theorem}
We have the chain derivative formula
$$\frac{d(f\circ g)_i}{dx_j}(x)=\sum_k\frac{df_i}{dx_k}(g(x))\cdot\frac{dg_k}{dx_j}(x)$$
as an equality of numbers.
\end{theorem}

\begin{proof}
This follows from the formula in Theorem 10.30, namely:
$$(f\circ g)'(x)=f'(g(x))\cdot g'(x)$$

Indeed, what we have here is an equality of matrices, and at the level of the individual entries of these matrices, by performing the multiplication on the right, we obtain:
$$(f\circ g)'(x)_{ij}=\sum_kf'(g(x))_{ik}\cdot g'(x)_{kj}$$

Thus, we are led to the formula in the statement.
\end{proof}

As a standard application now, generalizing one variable results, we have:

\begin{theorem}
Assuming that $f:X\to\mathbb R^M$ is differentiable, with $X\subset\mathbb R^N$ being convex, we have the estimate
$$||f(x)-f(y)||\leq M||x-y||$$
for any $x,y\in X$, where the quantity on the right is given by:
$$M=\sup_{x\in X}||f'(x)||$$
Moreover, this estimate can be sharp, for instance for the linear functions.
\end{theorem}

\begin{proof}
This is something quite tricky, which in several variables cannot be proved with bare hands. However, we can get it by using our chain derivative formula. Consider indeed the path $\gamma:[0,1]\to\mathbb R^M$ given by the following formula:
$$\gamma(t)=tx+(1-t)y$$

Next, let us set $g(t)=f(\gamma(t))$. We have then, according to the chain rule:
$$g'(t)
=f'(\gamma(t))\gamma'(t)
=f'(\gamma(t))(x-y)$$

But this gives the following estimate, with $M>0$ being as in the statement:
$$|g'(t)|
\leq||f'(\gamma(t))||\cdot||x-y||
\leq M||x-y||$$

Now by using one-variable results from chapter 3, we obtain from this:
$$||g(1)-g(0)||\leq||M||\cdot||x-y||$$

On the other hand, by definition of $g$ we have $g(1)=f(x),g(0)=f(y)$. Thus, we obtain the formula in the statement. Finally, the last assertion is clear.
\end{proof}

As a conclusion to all this, we have extended to the case of multivariable functions most of what we know about the 1-variable functions, at the level of continuity, and of first derivatives. And with everything being, after all, quite straightforward.

\section*{10d. Differential equations}

As an application of the methods that we recently learned, linear algebra and partial derivatives, we can have now a more systematic look at the differential equations. Generally speaking, the interesting equations are of order 2, for later. In the meantime, however, we can be kings in 1D. At the beginning of everything here, we have:

\begin{theorem}
The exponential function is the unique solution of:
$$f'=f\quad,\quad f(0)=1$$
More generally, the solutions of $f'=af$ are $f(x)=\lambda e^{ax}$.
\end{theorem}

\begin{proof}
The first assertion is something that we know since chapter 3, and the second assertion is a straightforward generalization of this, coming as follows:
\begin{eqnarray*}
f'=af
&\implies&\frac{f'}{f}=a\\
&\implies&(\log f)'=(ax)'\\
&\implies&\log f=ax+c\\
&\implies&f=\lambda e^{ax}
\end{eqnarray*}

Alternatively, $f'=af$ shows that we have $(e^{-ax}f)'=0$, which gives the result.
\end{proof}

With a bit more work, we can solve similar equations in degree 2, as follows:

\begin{theorem}
Given an equation $f''=af+bf'$, let $r,s$ be the roots of:
$$x^2=a+bx$$
\begin{enumerate}
\item In the case $r\neq s$, the solutions are $f(x)=\gamma e^{rx}+\delta e^{sx}$.

\item In the case $r= s$, the solutions are $f(x)=(\lambda x+\mu)e^{rx}$.
\end{enumerate}
\end{theorem}

\begin{proof}
This is something straightforward, the idea being as follows:

\medskip

(1) Our first goal is to put our equation in a simpler form. We have:
$$(f'-rf)'=f''-rf'=af+(b-r)f'$$ 

Now let us look for numbers $r,s$ such that this equals $s(f'-rf)$. We have:
$$af+(b-r)f'=s(f'-rf)\iff rs=-a,r+s=b$$

Thus, we certainly have numbers $r,s\in\mathbb C$ as desired, appearing as solutions of:
$$x^2-bx-a=0$$

(2) As a conclusion to this, with $r,s\in\mathbb C$ being as above, our equation reads:
$$(f'-rf)'=s(f'-rf)$$

Now with $g=f'-rf$ this equation reads $g'=sg$, which by Theorem 10.33 has as solutions the functions $g(x)=\lambda e^{sx}$. Thus, we are left with solving:
$$f'=rf+\lambda e^{sx}$$

But for $r\neq s$, by arguing like in the proof of Theorem 10.33, or modifying a bit our equation, and applying Theorem 10.33, the solutions are the following functions:
$$f(x)=\gamma e^{rx}+\frac{\lambda}{s-r}\,e^{sx}$$

As for the case $r=s$, where $f'=rf+\lambda e^{rx}$, this is similar, the solutions being:
$$f(x)=(\lambda x+\mu)e^{rx}$$

Thus, we are led to the conclusions in the statement.
\end{proof}

Along the same lines, at a more advanced level, we can use linear algebra:

\begin{theorem}
A differential equation of type $f''=af+bf'$ can be written as
$$\binom{f'}{f''}=\begin{pmatrix}0&1\\ a&b\end{pmatrix}\binom{f}{f'}$$
which in terms of $g=\binom{f}{f'}$ and $A=\binom{0\ 1}{a\ b}$ takes the following compact form,
$$g'=Ag$$
and whose solutions are as follows, with $v$ being the initial data vector,
$$g=e^{Ax}v$$
and with $e^{Ax}$ being the exponential of the $2\times2$ matrix $Ax$.
\end{theorem}

\begin{proof}
This is something more advanced, the idea being as follows:

\medskip

(1) To start with, by using the trick in the statement, and with the matrix multiplication being the usual one, ``multiply rows by columns'', our equation reads indeed:
$$g'=Ag$$ 

Now the point is that this latter equation reminds the one-variable equation $f'=af$ from Theorem 10.33, having as solutions the functions $f(x)=\lambda e^{ax}$. 

\medskip

(2) Equivalently, with some changes in the notations, we can say that this reminds the one-variable equation $g'=Ag$ with $A\in\mathbb R$, that we know how to solve, having as solutions the functions $g(x)=e^{Ax}v$, with $v\in\mathbb R$. But, based on this, we can conjecture that the solutions of our original equation are as follows, as in the statement:
$$g=e^{Ax}v$$

(3) So, can we make some sense of this? To start with, we need to know how to exponentiate the matrices of type $B=Ax$, and in answer, we can declare that:
$$e^B=\sum_{k=0}^\infty\frac{B^k}{k!}$$

(4) But, will this work. In order to talk about convergence of the above series, let us endow the space of $2\times2$ matrices with its usual norm, namely:
$$||B||=\sup_{||x||=1}||Bx||$$

To be more precise, it is easy to see that this is indeed a norm, which in addition satisfies $||AB||\leq||A||\cdot||B||$. But with this, we have convergence indeed, coming from:
$$||e^B||
\leq\sum_{k=0}^\infty\frac{||B^k||}{k!}
\leq\sum_{k=0}^\infty\frac{||B||^k}{k!}
=e^{||B||}<\infty$$

(5) Summarizing, our $g=e^{Ax}v$ conjecture above makes sense. Now regarding the proof of this conjecture, in one sense this is clear, coming from the following computation, and I will leave it to you, to check that all the algebra here works just fine:
$$(e^{Ax}v)'=(e^{Ax})'v=(Ae^{ax})v=A(e^{Ax}v)$$

(6) As for the uniqueness of our solutions, this is something a bit more complicated, but we can argue here that since by Theorem 10.34 the solutions of $g'=Ag$ depend on two parameters, these solutions can only be the functions $g=e^{Ax}v$ that we found here, depending on two parameters too, namely the entries of the vector $v\in\mathbb R^2$.

\medskip

(7) Next, for the discussion to be complete, it still makes sense to explicitly compute our solutions $g=e^{Ax}v$, see if we get indeed what we previously found in Theorem 10.34. But this can be done, in two steps. First, we must compute the powers of $A$:
$$\begin{pmatrix}0&1\\ a&b\end{pmatrix}^k=?$$

But this can be done, with some combinatorial pain, and we will leave this as an instructive exercise, and then we can compute $e^{Ax}$, and the solutions $g=e^{Ax}v$, and these solutions agree of course with what we previously found, in Theorem 10.34.

\medskip

(8) At a more advanced level now, the main problem that we have, exponentiating matrices, can be quickly solved by using diagonalization, as follows:
$$B=P\begin{pmatrix}\lambda_1&0\\0&\lambda_2\end{pmatrix}P^{-1}
\implies e^B=P\begin{pmatrix}e^{\lambda_1}&0\\0&e^{\lambda_2}\end{pmatrix}P^{-1}$$

(9) So, let us see what these advanced linear algebra methods teach us, in relation with our original problem. The eigenvalues $r,s$ of a $2\times2$ matrix $A$ are computable by using the well-known formulae $r+s=Tr(A)$ and $rs=\det A$, and in our case, we get:
$$A=\begin{pmatrix}0&1\\ a&b\end{pmatrix}\implies r+s=b,rs=-a$$

(10) Which is great, because what we have here are the roots of $x^2=a+bx$, so at least we have now a conceptual explanation for the occurrence of that roots. Good.

\medskip

(11) As for the continuation of the story, this is a bit more complicated. To start with, in the case $r\neq s$ the matrix is diagonalizable, and we know from the above that the solutions appear as linear combinations of $e^{rx}$, $e^{sx}$, as in Theorem 10.34.

\medskip

(12) As for the remaining case, $r=s$, here is where things get more complicated, because the matrix is no longer diagonalizable. However, by assuming some further linear algebra know-how, namely the Jordan form, we can again do the computations, and we reach to the linear combinations of $e^{rx}$, $xe^{rx}$, as in Theorem 10.34.
\end{proof}

More generally now, in higher degree, we have the following result:

\begin{theorem}
The equation $f^{(n)}=a_0f+a_1f'+\ldots+a_{n-1}f^{(n-1)}$ reads
$$\begin{pmatrix}f'\\ f''\\\vdots\\ f^{(n)}\end{pmatrix}
=\begin{pmatrix}
&1&&\\
&&\ddots&\\ 
&&&1\\
a_0&a_1&\ldots&a_{n-1}\end{pmatrix}
\begin{pmatrix}f\\ f'\\\vdots\\ f^{(n-1)}\end{pmatrix}$$
which in terms of the matrix $A$ and vector $g$ on the right takes the compact form
$$g'=Ag$$
and whose solutions are as follows, with $v$ being the initial data vector,
$$g=e^{Ax}v$$
and with $e^{Ax}$ being the exponential of the $n\times n$ matrix $Ax$.
\end{theorem}

\begin{proof}
This is again something quite self-explanatory, generalizing Theorem 10.35, and we will leave some further learning here as an exercise. Needless to say, there are countless things that can be said here, and the more you know, the better.
\end{proof}

As yet another application of the techniques that we recently learned, namely linear algebra and partial derivatives, this time involving partial derivatives, we can now fully solve the 1D wave equation, by using a 2D trick, as follows:

\index{d'Alembert formula}
\index{Riemann sum}
\index{discretization}

\begin{theorem}
The solution of the 1D wave equation $\ddot{f}=v^2f''$ with initial value conditions $f(x,0)=g(x)$ and $\dot{f}(x,0)=h(x)$ is given by the d'Alembert formula:
$$f(x,t)=\frac{g(x-vt)+g(x+vt)}{2}+\frac{1}{2v}\int_{x-vt}^{x+vt}h(s)ds$$
Also, in the context of our previous lattice model discretizations, what happens is more or less that the above d'Alembert integral gets computed via Riemann sums.
\end{theorem}

\begin{proof}
We already talked about waves in this book, on many occasions, but the above formula is still to be proved. Let us make the following change of variables:
$$y=x-vt\quad,\quad z=x+vt$$

(1) As a first observation, this is something quite tricky, mixing space and time variables. You would probably  even say that this is crazy, but after thinking a bit, the discussion of the basic solutions from chapter 8, involving sinusoids, makes it quite clear that the waves tend to mix space and time, so this is, after all, not that crazy. 

\medskip

(2) In any case, with this change of variables done, we have, using the chain rule:
$$\frac{df}{dt}
=\frac{df}{dy}\cdot\frac{dy}{dt}+\frac{df}{dz}\cdot\frac{dz}{dt}
=-v\,\frac{df}{dy}+v\,\frac{df}{dz}$$

By using the chain rule again, the second time derivative is given by:
\begin{eqnarray*}
\frac{d^2f}{dt^2}
&=&-v\left(\frac{d^2f}{dy^2}\cdot\frac{dy}{dt}+\frac{d^2f}{dydz}\cdot\frac{dz}{dt}\right)
+v\left(\frac{d^2f}{dzdy}\cdot\frac{dy}{dt}+\frac{d^2f}{dz^2}\cdot\frac{dz}{dt}\right)\\
&=&-v\left(-v\,\frac{d^2f}{dy^2}+v\,\frac{d^2f}{dydz}\right)
+v\left(-v\,\frac{d^2f}{dzdy}+v\,\frac{d^2f}{dz^2}\right)\\
&=&v^2\left(\frac{d^2f}{dy^2}+\frac{d^2f}{dz^2}-2\frac{d^2f}{dydz}\right)
\end{eqnarray*}

(3) Regarding now the first space derivative, this can be computed as follows:
$$\frac{df}{dx}
=\frac{df}{dy}\cdot\frac{dy}{dx}+\frac{df}{dz}\cdot\frac{dz}{dx}
=\frac{df}{dy}+\frac{df}{dz}$$

By using the chain rule again, the second space derivative is given by:
\begin{eqnarray*}
\frac{d^2f}{dt^2}
&=&\left(\frac{d^2f}{dy^2}\cdot\frac{dy}{dx}+\frac{d^2f}{dydz}\cdot\frac{dz}{dx}\right)
+\left(\frac{d^2f}{dzd\xi}\cdot\frac{dy}{dx}+\frac{d^2f}{dz^2}\cdot\frac{dz}{dx}\right)\\
&=&\left(\frac{d^2f}{dy^2}+\frac{d^2f}{dydz}\right)
+\left(\frac{d^2f}{dzdy}+\frac{d^2f}{dz^2}\right)\\
&=&\frac{d^2f}{dy^2}+\frac{d^2f}{dz^2}+2\frac{d^2f}{dydz}
\end{eqnarray*}

(4) Thus, our wave equation $\ddot{f}=v^2f''$ reformulates in a very simple way, as follows:
$$\frac{d^2f}{dydz}=0$$

But this latter equation tells us that our new $y,z$ variables get separated, and we conclude from this that the solution must be of the following special form:
$$f(x,t)=F(y)+G(z)=F(x-vt)+G(x+vt)$$

(5) In order now to finish, we must take into account the intial conditions for our problem, namely $f(x,0)=g(x)$ and $\dot{f}(x,0)=h(x)$. The first condition reads:
$$f(x,0)=g(x)\iff F(x)+G(x)=g(x)$$

As for the second condition, this can be processed as follows, with $H'=h$:
\begin{eqnarray*}
\dot{f}(x,0)=h(x)
&\iff&-vF'(x)+vG'(x)=h(x)\\
&\iff&G'(x)-F'(x)=\frac{h(x)}{v}\\
&\iff&G'(x)-F'(x)=\frac{H'(x)}{v}\\
&\iff&G(x)-F(x)=\frac{H(x)}{v}+c
\end{eqnarray*}

(6) Now by putting our two equations for $F,G$ together, we obtain:
$$F(x)=\frac{g(x)}{2}-\frac{H(x)}{2v}-\frac{c}{2}\quad,\quad 
G(x)=\frac{g(x)}{2}+\frac{H(x)}{2v}+\frac{c}{2}$$

We conclude that the solution of our equation is given by:
\begin{eqnarray*}
f(x,t)
&=&F(x-vt)+G(x+vt)\\
&=&\frac{g(x-vt)}{2}-\frac{H(x-vt)}{2v}+\frac{g(x+vt)}{2}+\frac{H(x+vt)}{2v}\\
&=&\frac{g(x-vt)+g(x+vt)}{2}+\frac{H(x+vt)-H(x-vt)}{2v}\\
&=&\frac{g(x-vt)+g(x+vt)}{2}+\frac{1}{2v}\int_{x-vt}^{x+vt}h(s)ds
\end{eqnarray*}

Thus, we are led to the d'Alembert formula in the statement.
\end{proof}

As a conclusion, we have extended to the multivariable functions most of what we know about 1-variable functions. Still waiting to be discussed are the higher derivatives, and their applications, and we will come back to this in chapter 12, after a short geometric break, and then of course integration, that we will discuss later, in chapters 13-16.

\section*{10e. Exercises}

This was a quite straightforward chapter, assuming that you have understood well Part I and chapter 9, and as exercises here, mostly straightforward, we have:

\begin{exercise}
Write a short essay about the basic theory of sequences and continuous functions over metric spaces, following the material from Part I.
\end{exercise}

\begin{exercise}
Clarify the relation between the various notions of connectedness, with the goal of decomposing each metric space into connected components.
\end{exercise}

\begin{exercise}
Enjoy defining the derivative of a function $f:\mathbb R^N\to\mathbb R^N$ as being a number, namely its Jacobian. What can you do, and what not, with this approach?
\end{exercise}

\begin{exercise}
Look up the internet for various other formulations of the chain rule for derivatives, and prove them all, quickly, by using what we know.
\end{exercise}

As bonus exercise, think a bit about higher derivatives, and what can you do with them, following the material from Part I. We will be back to this, after a break.

\chapter{Some geometry}

\section*{11a. Equations, conics}

As an application of the theory of partial derivatives that we developed, and of calculus in general, as we know it at this point, we can now talk about some beautiful things, in relation with celestial mechanics. As a starting point, you have surely noticed that the Sun moves around the Earth on a circle. However, when carefully measuring it, this circle is not exactly a circle, but rather an ellipse. Also, some further possible trajectories, of one object with respect to another, due to gravity, include parabolas and hyperbolas, that you can observe with a telescope, if you are lucky, by looking at certain asteroids.

\bigskip

So, before even starting to look at the equations of gravity, and having some fun in solving them, we need a mathematical theory of curves like ellipses, parabolas and hyperbolas, which are what we can expect to find, as trajectories, from gravity computations. And, good news, this theory exists, since the ancient Greeks. Let us start with:

\index{conic}
\index{algebraic curve}
\index{degree 2 equation}

\begin{definition}
A conic is a plane algebraic curve of the form
$$C=\left\{(x,y)\in\mathbb R^2\Big|P(x,y)=0\right\}$$
with $P\in\mathbb R[x,y]$ being of degree $\leq2$.
\end{definition}

As basic examples of conics, we have the ellipses, parabolas and hyperbolas. The simplest examples of these are as follows, with the ellipse actually being a circle:
$$x^2+y^2=1\quad,\quad x^2=y\quad,\quad xy=1$$

Observe that, due to our assumption $\deg P\leq 2$, we have as conics some degenerate curves as well, such as lines, $\emptyset$, and $\mathbb R^2$ itself, coming from $\deg P\leq1$, as follows:
$$x=0\quad,\quad 1=0\quad,\quad 0=0$$

This might suggest to replace our assumption $\deg P\leq 2$ by $\deg P=2$, but we will not do so, because $\deg P=2$ does not rule out degenerate situations, such as:
$$x^2+y^2=-1\quad,\quad x^2+y^2=0\quad,\quad x^2=0\quad,\quad xy=0$$

In fact, what we get here are $\emptyset$, points, lines, and pairs of lines, so in the end, assuming $\deg P=2$ instead of $\deg P\leq 2$ would only rule out $\mathbb R^2$ itself, which is not worth it.

\bigskip

Summarizing, our notion of conic from Definition 11.1 looks quite reasonable, so let us agree on this notion. Getting now to classification matters, we first have:

\index{ellipse}
\index{parabola}
\index{hyperbola}

\begin{proposition}
Up to non-degenerate linear transformations of the plane, 
$$\binom{x}{y}\to A\binom{x}{y}$$
with $\det A\neq0$, the conics fall into two classes, as follows:
\begin{enumerate}
\item Non-degenerate: circles, parabolas, hyperbolas.

\item Degenerate: $\emptyset$, points, lines, pairs of lines, $\mathbb R^2$.
\end{enumerate}
\end{proposition}

\begin{proof}
As a first observation, looks like we forgot the ellipses, but via linear transformations these become circles, so things fine. As for the proof, this goes as follows:

\medskip

(1) Consider an arbitrary conic, written as follows, with $a,b,c,d,e,f\in\mathbb R$:
$$ax^2+by^2+cxy+dx+ey+f=0$$

(2) Assume first $a\neq0$. By making a square out of $ax^2$, up to a linear transformation in $(x,y)$, we can get rid of the term $cxy$, and we are left with:
$$ax^2+by^2+dx+ey+f=0$$

In the case $b\neq0$ we can make two obvious squares, and again up to a linear transformation in $(x,y)$, we are left with an equation as follows:
$$x^2\pm y^2=k$$

In the case of positive sign, $x^2+y^2=k$, the solutions are the circle, when $k\geq0$, the point, when $k=0$, and $\emptyset$, when $k<0$. As for the case of negative sign, $x^2-y^2=k$, which reads $(x-y)(x+y)=k$, here once again by linearity our equation becomes $xy=l$, which is a hyperbola when $l\neq0$, and two lines when $l=0$.

\medskip

(3) In the case $b\neq0$ the study is similar, with the same solutions, so we are left with the case $a=b=0$. Here our conic is as follows, with $c,d,e,f\in\mathbb R$:
$$cxy+dx+ey+f=0$$

If $c\neq 0$, by linearity our equation becomes $xy=l$, which produces a hyperbola or two lines, as explained before. As for the remaining case, $c=0$, here our equation is:
$$dx+ey+f=0$$

But this is generically the equation of a line, unless we are in the case $d=e=0$, where our equation is $f=0$, having as solutions $\emptyset$ when $f\neq0$, and $\mathbb R^2$ when $f=0$.

\medskip

(4) So, this was the study of an arbitrary conic, and by putting now everything together, we are led to the conclusions in the statement.
\end{proof}

In order now to plainly classify the conics, without reference to a linear transformation of the plane, we just need to apply linear transformations to the curves that we found in Proposition 11.2. This leads to the following classification result:

\begin{proposition}
The conics fall into two classes, as follows:
\begin{enumerate}
\item Non-degenerate: ellipses, parabolas, hyperbolas.

\item Degenerate: $\emptyset$, points, lines, pairs of lines, $\mathbb R^2$.
\end{enumerate}
Also, the compact conics are $\emptyset$, the points, and the ellipses.
\end{proposition}

\begin{proof}
We have several assertions here, the idea being as follows:

\medskip

(1) As said above, in order to get to such a classification result, we just need to apply linear transformations to the curves that we found in Proposition 11.2. But this leaves the list there unchanged, up to the circles becoming ellipses, as stated above.

\medskip

(2) In what regards the last assertion, this is clear from the first one, but since this assertion is quite interesting, let us give it a quick, independent proof as well. Consider an arbitary conic, written as follows, with $a,b,c,d,e,f\in\mathbb R$:
$$ax^2+by^2+cxy+dx+ey+f=0$$

Compacity rules then out the case $c\neq0$, and our conic must be in fact:
$$ax^2+by^2+dx+ey+f=0$$

But then with $a,b\neq0$ we must have by compacity $a,b>0$ or $a,b<0$, and we get an ellipsis, then with $a=0,b\neq0$ or $a\neq0,b=0$ we get by compacity either $\emptyset$ or a point, and finally with $a=b=0$ the compacity rules out again everything, except for $\emptyset$.
\end{proof}

As a third main result now regarding the conics, also known since the ancient Greeks, and which justifies the name ``conics'', coming from ``cone'', we have:

\begin{proposition}
Up to some degenerate cases, the conics are exactly the curves which appear by cutting a $2$-sided cone with a plane.
\end{proposition}

\begin{proof}
This is something quite tricky, the idea being as follows:

\medskip

(1) By suitably choosing our coordinate axes $(x,y,z)$, we can assume that our 2-sided cone is given by an equation as follows, with $k>0$:
$$x^2+y^2=kz^2$$

In order to prove the result, we must intersect this cone with an arbitrary plane, which has an equation as follows, with $(a,b,c)\neq(0,0,0)$:
$$ax+by+cz=d$$

(2) However, before getting into computations, observe that what we want to find is a certain  degree 2 equation in the above plane, for the intersection. Thus, it is convenient to change the coordinates, as for our plane to be given by the following equation:
$$z=0$$

(3) But with this done, what we have to do is to see how the cone equation $x^2+y^2=kz^2$ changes, under this change of coordinates, and then set $z=0$, as to get the $(x,y)$ equation of the intersection. But this leads, via some thinking or computations, to the conclusion that the cone equation $x^2+y^2=kz^2$ becomes in this way a degree 2 equation in $(x,y)$, which can be arbitrary, and so to the final conclusion in the statement.

\medskip

(4) Alternatively, and perhaps more concretely, we can use the original coordinates, with the cone being $x^2+y^2=kz^2$, and compute the intersection, with the conclusion that what we get, depending on the slope of the cone, and modulo degenerate cases, is an ellipsis, hyperbola or parabola. So, by invoking Proposition 11.3, we obtain the result.

\medskip

(5) Summarizing, we have proved the result, modulo some details and interesting computations which are left to you, reader. Left to you as well is the full discussion concerning degree 2 curve degeneracy vs cone cutting degeneracy, with the remark that in what regards the cone cuts, the degenerate cases are very easy to identity and list, with the list consisting of $\emptyset$, the points, the lines, the pairs of lines, and $\mathbb R^2$ itself. 
\end{proof}

All this is very nice, and as a conclusion to what we have so far about conics, we have the following statement, containing all the needed essentials:

\begin{theorem}
The conics, which are the algebraic curves of degree $2$ in the plane,
$$C=\left\{(x,y)\in\mathbb R^2\Big|P(x,y)=0\right\}$$
with $\deg P\leq 2$, appear modulo degeneration by cutting a $2$-sided cone with a plane, and can be classified into ellipses, parabolas and hyperbolas.
\end{theorem}

\begin{proof}
This follows indeed by putting together the above results, and with the discussion concerning degeneration being left, as usual, as an instructive exercise.
\end{proof}

Moving ahead now, the most interesting conics, which are both compact and non-degenerate, are the ellipses. So, let us study them more in detail. As a starting point, we have the following statement, summarizing our knowledge about ellipses:

\begin{proposition}
The compact non-degenerate conics are the ellipses, which can be written, modulo rotations and translations in the plane, as
$$\left(\frac{x}{a}\right)^2+\left(\frac{y}{b}\right)^2=1$$
with the parameters $a,b>0$ measuring half the size of a box containing the ellipse. These ellipses also appear by compactly cutting a cone with a plane.
\end{proposition}

\begin{proof}
In this statement most of the mathematics is from above, and with our explanations regarding the parameters $a,b>0$ being something obvious.
\end{proof}

The above result is not the end of the story with ellipses, because we have as well, as a complement to it, or even as a rival result, which is just fine on its own:

\begin{proposition}
The ellipses appear via equations of the following type, with $p,q$ being two points in the plane, and with $l\geq d(p,q)$ being a certain length:
$$d(z,p)+d(z,q)=l$$
For an ellipse parametrized as before, $(x/a)^2+(y/b)^2=1$ with $a\geq b\geq 0$, the focal points are $p=(0,-r)$ and $q=(0,r)$, with $r=\sqrt{a^2-b^2}$, and the length is $l=2a$.
\end{proposition}

\begin{proof}
As already mentioned, it is possible to take $d(z,p)+d(z,q)=l$ as a definition for the ellipses, which is nice because all you need for drawing such an ellipse is a string and a pencil, and then work out all the theory starting from this. In what concerns us, we will rather further build on what we know from Proposition 11.6, as follows:

\medskip

(1) After some routine thinking, in order to fully prove the result, what we have to do is to take an ellipse as parametrized in Proposition 11.6, and look for the focal points:
$$
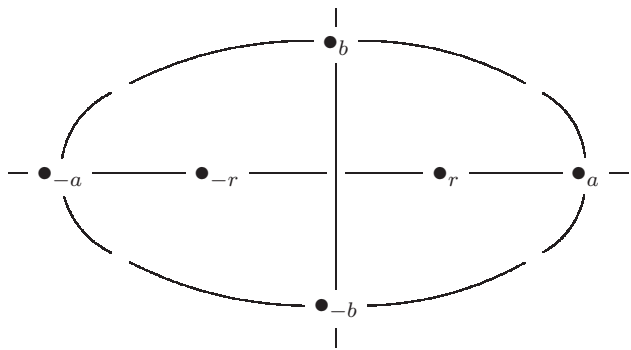
\xymatrix@R=7pt@C=7pt{
&&&&&&\ar@{-}[d]\\
&&&&&&\bullet_b\ar@{-}[dddddd]\ar@{-}@/_/[dllll]&\\
&&\ar@{-}@/_/[ddl]&&&&&&&&\ar@{-}@/_/[ullll]\\
&&&&&&&&&&\\
\ar@{-}[r]&\bullet_{-a}\ar@{-}@/_/[ddr]\ar@{-}[rrr]&&&\bullet_{-r}\ar@{-}[rr]&&\ar@{-}[rr]&&\bullet_r\ar@{-}[rrr]&&&\bullet_a\ar@{-}@/_/[uul]\ar@{-}[r]&\\
&&&&&&&&&&\\
&&\ar@{-}@/_/[drrrr]&&&&&&&&\ar@{-}@/_/[uur]\\
&&&&&&\bullet_{-b}\ar@{-}@/_/[urrrr]\ar@{-}[d]\\
&&&&&&}$$

To be more precise, we are looking for a number $r>0$, and a number $l>0$, such that our ellipse appears as $d(z,p)+d(z,q)=l$, with $p=(0,-r)$ and $q=(0,r)$.

\medskip

(2) Let us first compute these numbers $r,l>0$. Assuming that our result holds indeed as stated, by taking $z=(0,a)$, we see that the length $l$ is:
$$l=(a-r)+(a+r)=2a$$

As for the parameter $r$, by taking $z=(b,0)$, we conclude that we must have:
$$2\sqrt{b^2+r^2}=2a\implies r=\sqrt{a^2-b^2}$$

(3) With these observations made, let us prove the result. Given $l,r>0$, and setting $p=(0,-r)$ and $q=(0,r)$, we have the following computation, with $z=(x,y)$:
\begin{eqnarray*}
&&d(z,p)+d(z,q)=l\\
&\iff&\sqrt{(x+r)^2+y^2}+\sqrt{(x-r)^2+y^2}=l\\
&\iff&\sqrt{(x+r)^2+y^2}=l-\sqrt{(x-r)^2+y^2}\\
&\iff&(x+r)^2+y^2=(x-r)^2+y^2+l^2-2l\sqrt{(x-r)^2+y^2}\\
&\iff&2l\sqrt{(x-r)^2+y^2}=l^2-4xr\\
&\iff&4l^2(x^2+r^2-2xr+y^2)=l^4+16x^2r^2-8l^2xr\\
&\iff&4l^2x^2+4l^2r^2+4l^2y^2=l^4+16x^2r^2\\
&\iff&(4x^2-l^2)(4r^2-l^2)=4l^2y^2
\end{eqnarray*}

(4) Now observe that we can further process the equation that we found as follows:
\begin{eqnarray*}
(4x^2-l^2)(4r^2-l^2)=4l^2y^2
&\iff&\frac{4x^2-l^2}{l^2}=\frac{4y^2}{4r^2-l^2}\\
&\iff&\frac{4x^2-l^2}{l^2}=\frac{y^2}{r^2-l^2/4}\\
&\iff&\left(\frac{x}{2l}\right)^2-1=\left(\frac{y}{\sqrt{r^2-l^2/4}}\right)^2\\
&\iff&\left(\frac{x}{2l}\right)^2+\left(\frac{y}{\sqrt{r^2-l^2/4}}\right)^2=1
\end{eqnarray*}

(5) Thus, our result holds indeed, and with the numbers $l,r>0$ appearing, and no surprise here, via the formulae $l=2a$ and $r=\sqrt{a^2-b^2}$, found in (2) above.
\end{proof}

Finally, as an interesting analytic result regarding the ellipses, we have:

\begin{theorem}
The area of an ellipse, given as usual by
$$\left(\frac{x}{a}\right)^2+\left(\frac{y}{b}\right)^2=1$$
is $A=\pi ab$. As for the length, this is given by the formula
$$L=4\int_0^{\pi/2}\sqrt{a^2\sin^2t+b^2\cos^2t}\,dt$$
and with this integral being generically not computable.
\end{theorem}

\begin{proof}
The area formula is something that we know well, from chapter 4. As for the situation with the length, which is quite surprising, the idea here is as follows:

\medskip

(1) To start with, what is the length of a curve $\gamma:[a,b]\to\mathbb R^2$? Good question, and in answer, a physicist would say that this is the quantity obtained by integrating the magnitude of the velocity vector over the curve, with respect to time:
$$L(\gamma)=\int_a^b||\gamma'(t)||dt$$

(2) Regarding now mathematicians, these would say that the length of a curve is the following quantity, with $(t_0=a,t_1,\ldots,t_{n-1},t_n=b)$ being a uniform division of $(a,b)$:
$$L(\gamma)=\lim_{n\to\infty}\sum_{i=1}^n||\gamma(t_i)-\gamma(t_{i-1})||$$

But, by using the fundamental theorem of calculus, we can write this as follows:
$$L(\gamma)=\lim_{n\to\infty}\sum_{i=1}^n\left|\left|\int_{t_{i-1}}^{t_i}\gamma'(t)dt\right|\right|$$

And then, by doing some standard analysis, we are led to the formula in (1).

\medskip

(3) Getting back now to the ellipses, we can compute their length, as follows:
\begin{eqnarray*}
L
&=&4\int_0^{\pi/2}\sqrt{\left(\frac{dx}{dt}\right)^2+\left(\frac{dy}{dt}\right)^2}\,dt\\
&=&4\int_0^{\pi/2}\sqrt{\left(\frac{da\cos t}{dt}\right)^2+\left(\frac{db\sin t}{dt}\right)^2}\,dt\\
&=&4\int_0^{\pi/2}\sqrt{a^2\sin^2t+b^2\cos^2t}\,dt
\end{eqnarray*}

(4) And the point is that this latter integral is generically not computable, unless we are in some special cases, such as $a=b$, where we obtain of course $2\pi a$.
\end{proof}

\section*{11b. Kepler and Newton}

As a continuation of the above, or rather as a complement, let us do some physics. Theorem 11.5 was the foundational result of modern mathematics, and by some magic, the foundational result of physics is something closely related to it, as follows:

\index{conic}
\index{Kepler laws}
\index{Netwon law}
\index{gravity}

\begin{theorem}
Planets and other celestial bodies move around the Sun on conics,
$$C=\left\{(x,y)\in\mathbb R^2\Big|P(x,y)=0\right\}$$
with $P\in\mathbb R[x,y]$ being of degree $2$, which can be ellipses, parabolas or hyperbolas.
\end{theorem}

\begin{proof}
This is something quite long, due to Kepler and Newton, but no fear, we know calculus, and therefore what can resist us. The proof goes as follows:

\medskip

(1) According to observations and calculations performed over the centuries, since the ancient times, and first formalized by Newton, following some groundbreaking work of Kepler, the force of attraction between two bodies of masses $M,m$ is given by:
$$||F||=G\cdot\frac{Mm}{d^2}$$

Here $d$ is the distance between the two bodies, and $G\simeq 6.674\times 10^{-11}$ is a constant. Now assuming that $M$ is fixed at $0\in\mathbb R^3$, the force exterted on $m$ positioned at $x\in\mathbb R^3$, regarded as a vector $F\in\mathbb R^3$, is given by the following formula:
$$F
=-||F||\cdot\frac{x}{||x||}
=-\frac{GMm}{||x||^2}\cdot\frac{x}{||x||}
=-\frac{GMmx}{||x||^3}$$

But $F=ma=m\ddot{x}$, with $a=\ddot{x}$ being the acceleration, second derivative of the position, so the equation of motion of $m$, assuming that $M$ is fixed at $0$, is:
$$\ddot{x}=-\frac{GMx}{||x||^3}$$

Obviously, the problem happens in 2 dimensions, and here the most convenient is to use standard $x,y$ coordinates, and denote our point as $z=(x,y)$. With this change made, and by setting $K=GM$, the equation of motion becomes:
$$\ddot{z}=-\frac{Kz}{||z||^3}$$

(2) The idea now is that the problem can be solved via some calculus. Let us write indeed our vector $z=(x,y)$ in polar coordinates, as follows:
$$x=r\cos\theta\quad,\quad 
y=r\sin\theta$$

We have then $||z||=r$, and our equation of motion becomes:
$$\ddot{z}=-\frac{Kz}{r^3}$$

Let us differentiate now $x,y$. By using the standard calculus rules, we have:
$$\dot{x}=\dot{r}\cos\theta-r\sin\theta\cdot\dot{\theta}$$
$$\dot{y}=\dot{r}\sin\theta+r\cos\theta\cdot\dot{\theta}$$

Differentiating one more time gives the following formulae:
$$\ddot{x}=\ddot{r}\cos\theta-2\dot{r}\sin\theta\cdot\dot{\theta}-r\cos\theta\cdot\dot{\theta}^2-r\sin\theta\cdot\ddot\theta$$
$$\ddot{y}=\ddot{r}\sin\theta+2\dot{r}\cos\theta\cdot\dot{\theta}-r\sin\theta\cdot\dot{\theta}^2+r\cos\theta\cdot\ddot\theta$$

Consider now the following two quantities, appearing as coefficients in the above:
$$a=\ddot{r}-r\dot{\theta}^2\quad,\quad b=2\dot{r}\dot{\theta}+r\ddot{\theta}$$

In terms of these quantities, our second derivative formulae read:
$$\ddot{x}=a\cos\theta-b\sin\theta\quad,\quad 
\ddot{y}=a\sin\theta+b\cos\theta$$

(3) We can now solve the equation of motion from (1). Indeed, with the formulae that we found for $\ddot{x},\ddot{y}$, our equation of motion takes the following form:
$$a\cos\theta-b\sin\theta=-\frac{K}{r^2}\cos\theta\quad,\quad 
a\sin\theta+b\cos\theta=-\frac{K}{r^2}\sin\theta$$

But these two formulae can be written in the following way:
$$\left(a+\frac{K}{r^2}\right)\cos\theta=b\sin\theta\quad,\quad 
\left(a+\frac{K}{r^2}\right)\sin\theta=-b\cos\theta$$

By making now the product, and assuming that we are in a non-degenerate case, where the angle $\theta$ varies indeed, we obtain by positivity that we must have:
$$a+\frac{K}{r^2}=b=0$$

(4) Let us first examine the second equation, $b=0$. This can be solved as follows:
\begin{eqnarray*}
b=0
&\iff&2\dot{r}\dot{\theta}+r\ddot{\theta}=0\\
&\iff&\frac{\ddot{\theta}}{\dot{\theta}}=-2\frac{\dot{r}}{r}\\
&\iff&(\log\dot{\theta})'=(-2\log r)'\\
&\iff&\log\dot{\theta}=-2\log r+c\\
&\iff&\dot{\theta}=\frac{\lambda}{r^2}
\end{eqnarray*}

As for the first equation the we found, namely $a+K/r^2=0$, this becomes:
$$\ddot{r}-\frac{\lambda^2}{r^3}+\frac{K}{r^2}=0$$

As a conclusion to all this, in polar coordinates, $x=r\cos\theta$, $y=r\sin\theta$, our equations of motion are as follows, with $\lambda$ being a constant, not depending on $t$:
$$\ddot{r}=\frac{\lambda^2}{r^3}-\frac{K}{r^2}\quad,\quad 
\dot{\theta}=\frac{\lambda}{r^2}$$

Even better now, by writing $K=\lambda^2/c$, these equations read:
$$\ddot{r}=\frac{\lambda^2}{r^2}\left(\frac{1}{r}-\frac{1}{c}\right)\quad,\quad 
\dot{\theta}=\frac{\lambda}{r^2}$$

(5) In order to study the first equation, we use a trick. Let us write:
$$r(t)=\frac{1}{f(\theta(t))}$$

Abbreviated, and by reminding that $f$ takes $\theta=\theta(t)$ as variable, this reads:
$$r=\frac{1}{f}$$

With the convention that dots mean as usual derivatives with respect to $t$, and that the primes will denote derivatives with respect to $\theta=\theta(t)$, we have:
$$\dot{r}=-\frac{f'\dot{\theta}}{f^2}=-\frac{f'}{f^2}\cdot\frac{\lambda}{r^2}=-\lambda f'$$

By differentiating one more time with respect to $t$, we obtain:
$$\ddot{r}=-\lambda f''\dot{\theta}=-\lambda f''\cdot\frac{\lambda}{r^2}=-\frac{\lambda^2}{r^2}f''$$

On the other hand, our equation for $\ddot{r}$ found in (4) above reads:
$$\ddot{r}=\frac{\lambda^2}{r^2}\left(\frac{1}{r}-\frac{1}{c}\right)=\frac{\lambda^2}{r^2}\left(f-\frac{1}{c}\right)$$

Thus, in terms of $f=1/r$ as above, our equation for $\ddot{r}$ simply reads:
$$f''+f=\frac{1}{c}$$

But this latter equation is elementary to solve. Indeed, both functions $\cos t,\sin t$ satisfy $g"+g=0$, so any linear combination of them satisfies as well this equation. But the solutions of $f''+f=1/c$ being those of $g''+g=0$ shifted by $1/c$, we obtain:
$$f=\frac{1+\varepsilon\cos\theta+\delta\sin\theta}{c}$$

Now by inverting, we obtain the following formula:
$$r=\frac{c}{1+\varepsilon\cos\theta+\delta\sin\theta}$$

(6) But this leads to the conclusion that the trajectory is a conic. Indeed, in terms of the parameter $\theta$, the formulae of the coordinates are:
$$x=\frac{c\cos\theta}{1+\varepsilon\cos\theta+\delta\sin\theta}$$
$$y=\frac{c\sin\theta}{1+\varepsilon\cos\theta+\delta\sin\theta}$$

Now observe that these two functions $x,y$ satisfy the following formula:
$$x^2+y^2
=\frac{c^2(\cos^2\theta+\sin^2\theta)}{(1+\varepsilon\cos\theta+\delta\sin\theta)^2}
=\frac{c^2}{(1+\varepsilon\cos\theta+\delta\sin\theta)^2}$$

On the other hand, these two functions satisfy as well the following formula:
\begin{eqnarray*}
(\varepsilon x+\delta y-c)^2
&=&\frac{c^2\big(\varepsilon\cos\theta+\delta\sin\theta-(1+\varepsilon\cos\theta+\delta\sin\theta)\big)^2}{(1+\varepsilon\cos\theta+\delta\sin\theta)^2}\\
&=&\frac{c^2}{(1+\varepsilon\cos\theta+\delta\sin\theta)^2}
\end{eqnarray*}

We conclude that our coordinates $x,y$ satisfy the following equation:
$$x^2+y^2=(\varepsilon x+\delta y-c)^2$$

But what we have here is an equation of a conic, and we are done.
\end{proof}

From a physical perspective, that of concretely solving the gravity equation, there is a long extra discussion, and lots of additional formulae, regarding the trajectory and its parameters, as functions of the initial data. Without getting into full details here, let us record however the following result, coming as a useful version of Theorem 11.9:

\begin{theorem}
In the context of a $2$-body problem, with $M$ fixed at $0$, and $m$ starting its movement from $Ox$, the equation of motion of $m$, namely
$$\ddot{z}=-\frac{Kz}{||z||^3}$$
with $K=GM$, and $z=(x,y)$, becomes in polar coordinates, $x=r\cos\theta$,  $y=r\sin\theta$,
$$\ddot{r}=\frac{\lambda^2}{r^2}\left(\frac{1}{r}-\frac{1}{c}\right)\quad,\quad 
\dot{\theta}=\frac{\lambda}{r^2}$$
for some $\lambda,c\in\mathbb R$, related by $\lambda^2=Kc$. The value of $r$ in terms of $\theta$ is given by
$$r=\frac{c}{1+\varepsilon\cos\theta+\delta\sin\theta}$$
for some $\varepsilon,\delta\in\mathbb R$. At the level of the affine coordinates $x,y$, this means
$$x=\frac{c\cos\theta}{1+\varepsilon\cos\theta+\delta\sin\theta}\quad,\quad 
y=\frac{c\sin\theta}{1+\varepsilon\cos\theta+\delta\sin\theta}$$
with $\theta=\theta(t)$ being subject to $\dot{\theta}=\lambda^2/r$, as above. Finally, we have
$$x^2+y^2=(\varepsilon x+\delta y-c)^2$$
which is a degree $2$ equation, and so the resulting trajectory is a conic.
\end{theorem}

\begin{proof}
This is a sort of ``best of'' the formulae found in the proof of Theorem 11.9. And in the hope of course that we have not forgotten anything. Finally, let us mention that the simplest illustration for this is the circular motion, and for details on this, not included in the above, we refer to the proof of Theorem 11.9.
\end{proof}

As a concrete question now, we would like to understand how the various parameters appearing above, namely $\lambda,c,\varepsilon,\delta$, which via some basic math can only tell us more about the shape of the orbit, appear from the initial data. The formulae here are as follows:

\begin{theorem}
In the context of Theorem 11.10, and in standard polar coordinates, $x=r\cos\theta$,  $y=r\sin\theta$, the initial data is as follows, with $R=r_0$:
$$r_0=\frac{c}{1+\varepsilon}\quad,\quad\theta_0=0$$
$$\dot{r}_0=-\frac{\delta\sqrt{K}}{\sqrt{c}}\quad,\quad\dot{\theta}_0=\frac{\sqrt{Kc}}{R^2}$$
$$\ddot{r}_0=\frac{\varepsilon K}{R^2}\quad,\quad\ddot{\theta}_0=\frac{4\delta K}{R^2}$$
The corresponding formulae for the affine coordinates $x,y$ can be deduced from this. Also, the various motion parameters $c,\varepsilon,\delta$ and $\lambda=\sqrt{Kc}$ can be recovered from this data.
\end{theorem}

\begin{proof}
We have several assertions here, the idea being as follows:

\medskip

(1) As mentioned in Theorem 11.10, the object $m$ begins its movement on $Ox$. Thus we have $\theta_0=0$, and from this we get the formula of $r_0$ in the statement.

\medskip

(2) Regarding the initial speed now, the formula of $\dot{\theta}_0$ follows from:
$$\dot{\theta}=\frac{\lambda}{r^2}=\frac{\sqrt{Kc}}{r^2}$$

Also, in what concerns the radial speed, the formula of $\dot{r}_0$ follows from:
\begin{eqnarray*}
\dot{r}
&=&\frac{c(\varepsilon\sin\theta-\delta\cos\theta)\dot{\theta}}{(1+\varepsilon\cos\theta+\delta\sin\theta)^2}\\
&=&\frac{c(\varepsilon\sin\theta-\delta\cos\theta)}{c^2/r^2}\cdot\frac{\sqrt{Kc}}{r^2}\\
&=&\frac{\sqrt{K}(\varepsilon\sin\theta-\delta\cos\theta)}{\sqrt{c}}
\end{eqnarray*}

(3) Regarding now the initial acceleration, by using $\dot{\theta}=\sqrt{Kc}/r^2$ we find:
$$\ddot{\theta}=-2\sqrt{Kc}\cdot\frac{2r\dot{r}}{r^3}=-\frac{4\sqrt{Kc}\cdot\dot{r}}{r^2}$$

In particular at $t=0$ we obtain the formula in the statement, namely:
$$\ddot{\theta}_0=-\frac{4\sqrt{Kc}\cdot\dot{r}_0}{R^2}
=\frac{4\sqrt{Kc}}{R^2}\cdot\frac{\delta\sqrt{K}}{\sqrt{c}}
=\frac{4\delta K}{R^2}$$

(4) Also regarding acceleration, with $\lambda=\sqrt{Kc}$ our main motion formula reads:
$$\ddot{r}=\frac{Kc}{r^2}\left(\frac{1}{r}-\frac{1}{c}\right)$$

In particular at $t=0$ we obtain the formula in the statement, namely:
$$\ddot{r}_0=\frac{Kc}{R^2}\left(\frac{1}{R}-\frac{1}{c}\right)=\frac{Kc}{R^2}\cdot\frac{\varepsilon}{c}=\frac{\varepsilon K}{R^2}$$

(5) Finally, the last assertion is clear, and since the formulae look better anyway in polar coordinates than in affine coordinates, we will not get into details here.
\end{proof}

And we will leave some further computations here, up to the very concrete level, as exercises. Also, as further interesting exercises regarding gravity, you can try to deduce from the above the results for 1D free falls from chapter 8, and also have some work done on the uniform gravity case, where the trajectories are parabolas, quite easy to compute. Also, you can have a look at the pendulum, followed by harmonic oscillators, and damping, with the conclusion that, quite remarkably, the mathematics there coincides with the mathematics of the basic differential equations, as developed in chapter 10.

\bigskip

In a word, and never doubt about this, math is the same thing as physics.

\section*{11c. Algebraic manifolds}

We have seen so far the mathematical and physical meaning of the conics, which are the simplest curves around, with in both cases, very simple answers. The continuation of the story, however, is more complicated, because beyond conics, things ramify. A first idea, in order to generalize the conics, is to look at zeroes of arbitrary polynomials:
$$P(x,y)=0\quad,\quad P\in\mathbb R[x,y]$$

More generally, we can look at zeroes of polynomials in arbitrary $N$ dimensions:
$$P(x_1,\ldots,x_N)=0\quad,\quad P\in\mathbb R[x_1,\ldots,x_N]$$

Observe that, at $N\geq3$, what we have is not exactly a curve, but rather some sort of $(N-1)$-dimensional surface, called algebraic hypersurface. Due to this, in order to have a full collection of beasts, of all possible dimensions, we must intersect such algebraic hypersurfaces. We are led in this way to zeroes of families of polynomials, as follows:

\index{algebraic manifold}

\begin{definition}
An algebraic manifold is a space of the form
$$X=\left\{(x_1,\ldots,x_N)\in\mathbb R^N\Big|P_i(x_1,\ldots,x_N)=0,\forall i\right\}$$
with $P_i\in\mathbb R[x_1,\ldots,x_N]$ being a family of polynomials.
\end{definition}

The question is now, what can we say about such manifolds? In answer, many things, and as a natural continuation of our previous work on the conics, you can:

\bigskip

(1) Look at plane curves of higher degree, cubics, quartics and so on. And here, all sorts of interesting curves, such as sinusoidal spirals, polynomial lemniscates or stelloids, quite often appearing as field lines in physics, are waiting for your attention.

\bigskip

(2) Look at degree 2 hypersurfaces in $\mathbb R^3$ or higher, which are called quadrics. And here, again, many interesting beasts, quite often related to physics, with some appearing for instance in Einstein's relativity theory, are waiting for your attention.

\bigskip

In short, many things to be learned here. At the general level now, that of Definition 11.12 as stated, a fruitful idea is that of trying to understand the correspondence between algebraic manifolds and polynomials. And here, we first have the following result:

\begin{proposition}
The algebraic manifolds are the sets of the form
$$X=\left\{x\in\mathbb R^N\Big|P(x)=0,\forall P\in I\right\}$$
with $I\subset\mathbb R[x_1,\ldots,x_N]$ being a radical ideal, meaning a set satisfying:
\begin{enumerate}
\item Ideal property: $P_i\in I$, $Q_i\in\mathbb R[x_1,\ldots,x_N]\implies\sum_iP_iQ_i\in I$.

\item Radicality property: $P^k\in I\implies P\in I$.
\end{enumerate}
\end{proposition}

\begin{proof}
This is indeed something self-explanatory, by taking $I$ to be the set of all polynomials $P\in\mathbb R[x_1,\ldots,x_N]$ vanishing on $X$, which certainly satisfies (1,2).
\end{proof}

Summarizing, we have an interesting algebraic geometry correspondence $X\leftrightarrow I$, given in one sense by Definition 11.12, and in the other sense, by Proposition 11.13. The problem however is that this correspondence is not bijective, for instance because at $N=1$, the non-trivial ideal $I=(x^2+1)\subset\mathbb R[x]$ produces the trivial manifold $X=\emptyset$. 

\bigskip

In view of this, we must trade $\mathbb R$ for $\mathbb C$, where arbitrary polynomials have roots. And with this done, coming as good news, we have the following remarkable result:

\begin{theorem}[Nullstellensatz]
We have a bijective correspondence
$$\Big(X\subset\mathbb C^N\Big)\quad\longleftrightarrow\quad\Big(I\subset\mathbb C[x_1,\ldots,x_N]\Big)$$
between algebraic manifolds in $\mathbb C^N$, and radical ideals of $\mathbb C[x_1,\ldots,x_N]$.
\end{theorem}

\begin{proof}
This is something quite tricky, due to Hilbert, the idea being as follows:

\medskip

(1) At $N=1$ polynomials have roots, $I=(P)\implies X_I\neq\emptyset$. The point now is that, by doing some algebra, something similar happens in arbitrary $N$ dimensions, in the sense that any proper ideal $I\subset \mathbb C[x_1,\ldots,x_N]$ produces a non-empty manifold, $X_I\neq\emptyset$.

\medskip

(2) Next, what we want to prove is that given an ideal $I\subset\mathbb C[x_1,\ldots,x_N]$, any polynomial $P\in\mathbb C[x_1,\ldots,x_N]$ vanishing on $X_I$ has the property $P^k\in I$, for some $k\in\mathbb N$. For this purpose, we can add 1 dimension, and consider the following ideal:
$$J=<I,x_{N+1}P(x_1,\ldots,x_N)-1>$$

(3) Now since we have $X_J=\emptyset$, by (1) we conclude that $J$ is trivial. In order now to best interpret this finding, consider the following algebra:
$$\mathbb C[x_1,\ldots,x_N][P^{-1}]=\mathbb C[x_1,\ldots,x_{N+1}]/(x_{N+1}P-1)$$

(4) The triviality of $J$ gives then a formula of the following type, with $f_i\in I$:
$$1=f_0+f_1x_{N+1}+\ldots+f_kx_{N+1}^k$$

Now by multiplying by $P^k$, we obtain from this $P^k\in I$, as desired.
\end{proof}

Quite interesting all this, and needless to say, it is possible to build a whole skyscraper, on the Nullstellensatz. For summarizing this discussion, let us formulate:

\begin{conclusion}
For beautiful mathematics and geometry, several complex variables are what you need.
\end{conclusion}

Moving on, let us have a look at projective geometry too, which is something fun, interesting, and quite often more fun and interesting than affine geometry itself. 

\bigskip

You might have heard or not of projective geometry. In case you didn't yet, the general principle is that ``this is the wonderland where parallel lines cross''. Which might sound a bit crazy, and not very realistic, but take a picture of some railroad tracks, and look at that picture. Do these parallel railroad tracks cross, on the picture? Sure they do. So, we are certainly not into abstractions here, but rather into serious science. QED.

\bigskip

Mathematically now, here are some axioms, to start with:

\index{projective space}

\begin{definition}
A projective space is a space consisting of points and lines, subject to the following conditions:
\begin{enumerate}
\item Each $2$ points determine a line.

\item Each $2$ lines cross, on a point.
\end{enumerate}
\end{definition}

As a basic example we have the usual projective plane $P^2_\mathbb R$, which is best seen as being the space of lines in $\mathbb R^3$ passing through the origin. To be more precise, let us call each of these lines in $\mathbb R^3$ passing through the origin a ``point'' of $P^2_\mathbb R$, and let us also call each plane in $\mathbb R^3$ passing through the origin a ``line'' of $P^2_\mathbb R$. Now observe the following:

\bigskip

(1) Each $2$ points determine a line. Indeed, 2 points in our sense means 2 lines in $\mathbb R^3$ passing through the origin, and these 2 lines obviously determine a plane in $\mathbb R^3$ passing through the origin, namely the plane they belong to, which is a line in our sense.

\bigskip

(2) Each $2$ lines cross, on a point. Indeed, 2 lines in our sense means 2 planes in $\mathbb R^3$ passing through the origin, and these 2 planes obviously determine a line in $\mathbb R^3$ passing through the origin, namely their intersection, which is a point in our sense.

\bigskip

Thus, what we have is a projective space in the sense of Definition 11.16. More generally now, we have the following construction, in arbitrary dimensions:

\index{M\"obius strip}
\index{Klein bottle}

\begin{theorem}
We can define the projective space $P^{N-1}_\mathbb R$ as being the space of lines in $\mathbb R^N$ passing through the origin, and in small dimensions:
\begin{enumerate}
\item $P^1_\mathbb R$ is the usual circle.

\item $P^2_\mathbb R$ is some sort of twisted sphere.
\end{enumerate}
\end{theorem}

\begin{proof}
We have several assertions here, with all this being of course a bit informal, and self-explanatory, the idea and some further details being as follows:

\medskip

(1) To start with, the fact that the space $P^{N-1}_\mathbb R$ constructed in the statement is indeed a projective space in the sense of Definition 11.16 follows from definitions, exactly as in the discussion preceding the statement, regarding the case $N=3$.

\medskip

(2) At $N=2$ now, a line in $\mathbb R^2$ passing through the origin corresponds to 2 opposite points on the unit circle $\mathbb T\subset\mathbb R^2$, according to the following scheme:
$$
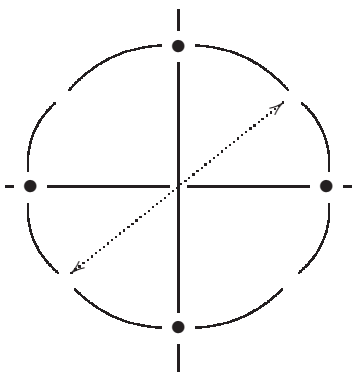
\xymatrix@R=6.7pt@C=2pt{
&&&&&&\ar@{-}[d]\\
&&&&&&\bullet\ar@{-}[dddddd]\ar@{-}@/_/[dllll]&\\
&&\ar@{-}@/_/[ddl]&&&&&&&&\ar@.[ddddllllllll]\ar@{-}@/_/[ullll]\\
&&&&&&&&&&\\
\ar@{-}[r]&\bullet\ar@{-}@/_/[ddr]\ar@{-}[rrrrr]&&&&&\ar@{-}[rrrrr]&&&&&\bullet\ar@{-}@/_/[uul]\ar@{-}[r]&\\
&&&&&&&&&&\\
&&\ar@{-}@/_/[drrrr]\ar@.[uuuurrrrrrrr]&&&&&&&&\ar@{-}@/_/[uur]\\
&&&&&&\bullet\ar@{-}@/_/[urrrr]\ar@{-}[d]\\
&&&&&&}$$

Thus, $P^1_\mathbb R$ corresponds to the upper semicircle of $\mathbb T$, with the endpoints identified, and so we obtain a circle, $P^1_\mathbb R=\mathbb T$, according to the following scheme: 
$$
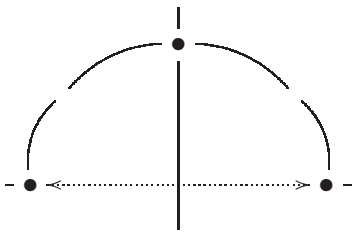
\xymatrix@R=7pt@C=2pt{
&&&&&&\ar@{-}[d]\\
&&&&&&\bullet\ar@{-}[dddd]\ar@{-}@/_/[dllll]&\\
&&\ar@{-}@/_/[ddl]&&&&&&&&\ar@{-}@/_/[ullll]\\
&&&&&&&&&&\\
\ar@{-}[r]&\bullet\ar@.[rrrrrrrrrr]&&&&&&&&&&\bullet\ar@{-}@/_/[uul]\ar@{-}[r]\ar@.[llllllllll]&\\
&&&&&&&&
}$$

(3) At $N=3$, the space $P^2_\mathbb R$ corresponds to the upper hemisphere of the sphere $S^2_\mathbb R\subset\mathbb R^3$, with the points on the equator identified via $x=-x$. Topologically speaking, we can deform if we want the hemisphere into a square, with the equator becoming the boundary of this square, and in this picture, the $x=-x$ identification corresponds to a ``identify opposite edges, with opposite orientations'' folding method for the square:
$$
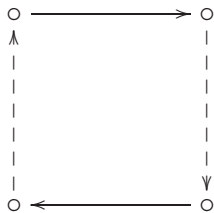
\xymatrix@R=55pt@C=55pt{
\circ\ar[r]&\circ\ar@{-->}[d]\\
\circ\ar@{-->}[u]&\circ\ar[l]}$$

(4) Thus, we have our space. In order to understand now what this beast is, let us look first at the other 3 possible methods of folding the square, which are as follows:
$$
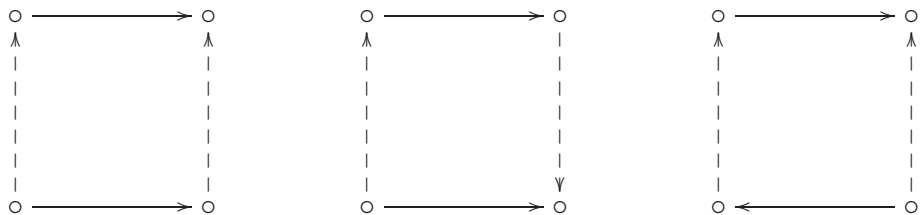
\xymatrix@R=55pt@C=55pt{
\circ\ar[r]&\circ\\
\circ\ar@{-->}[u]\ar[r]&\circ\ar@{-->}[u]}\qquad\qquad 
\xymatrix@R=55pt@C=55pt{
\circ\ar[r]&\circ\ar@{-->}[d]\\
\circ\ar[r]\ar@{-->}[u]&\circ}\qquad\qquad 
\xymatrix@R=55pt@C=55pt{
\circ\ar[r]&\circ\\
\circ\ar@{-->}[u]&\circ\ar[l]\ar@{-->}[u]}$$

Regarding the first space, the one on the left, things here are quite simple. Indeed, when identifying the solid edges we get a cylinder, and then when further identifying the dotted edges, what we get is some sort of closed cylinder, which is a torus.

\medskip

(5) Regarding the second space, the one in the middle, things here are more tricky. Indeed, when identifying the solid edges we get again a cylinder, but then when further identifying the dotted edges, we obtain some sort of ``impossible'' closed cylinder, called Klein bottle. This Klein bottle obviously cannot be drawn in 3 dimensions, but with a bit of imagination, you can see it, in its full splendor, in 4 dimensions.

\medskip

(6) Finally, regarding the third space, the one on the right, we know by symmetry that this must be the Klein bottle too. But we can see this as well via our standard folding method, namely identifying solid edges first, and dotted edges afterwards. Indeed, we first obtain in this way a M\"obius strip, and then, well, the Klein bottle.

\medskip

(7) With these preliminaries made, and getting back now to the projective space $P^2_\mathbb R$, we can see that this is something more complicated, of the same type, reminding the torus and the Klein bottle. So, we will call it ``sort of twisted sphere'', as in the statement, and exercise for you to figure out how this beast looks like, in 4 dimensions.
\end{proof}

All this is quite exciting, and reminds childhood and primary school, but is however a bit tiring for our neurons, guess that is pure mathematics. It is possible to come up with some explicit formulae for the embedding $P^2_\mathbb R\subset\mathbb R^4$, which are useful in practice, allowing us to do some analysis over $P^2_\mathbb R$, and we will leave this as an instructive exercise.

\bigskip

On the same topic, there is some linear algebra to be done too, as follows:

\index{rank 1 projection}

\begin{theorem}
The projective space $P^{N-1}_\mathbb R$ can be thought of as being the space of rank $1$ projections in the matrix algebra $M_N(\mathbb R)$, given by
$$P_x=\frac{1}{||x||^2}(x_ix_j)_{ij}$$
by identifying the lines in $\mathbb R^N$ passing through the origin with the corresponding rank $1$ projections in $M_N(\mathbb R)$, in the obvious way.
\end{theorem}

\begin{proof}
There are several things going on here, the idea being as follows:

\medskip

(1) The main assertion is more or less clear from definitions, the point being that the lines in $\mathbb R^N$ passing through the origin are obviously in bijection with the corresponding rank $1$ projections. Thus, we obtain the interpretation of $P^{N-1}_\mathbb R$ in the statement.

\medskip

(2) Regarding now the formula of the rank 1 projections, which is a must-know, for this, and in everyday life, consider a vector $y\in\mathbb R^N$. Its projection on $\mathbb Rx$ must be a certain multiple of $x$, and we are led in this way to the following formula:
$$P_xy
=\frac{<y,x>}{<x,x>}\,x
=\frac{1}{||x||^2}<y,x>x$$

(3) But with this in hand, we can now compute the entries of $P_x$, as follows:
\begin{eqnarray*}
(P_x)_{ij}
&=&<P_xe_j,e_i>\\
&=&\frac{1}{||x||^2}<e_j,x><x,e_i>\\
&=&\frac{x_jx_i}{||x||^2}
\end{eqnarray*}

Thus, we are led to the formula in the statement.
\end{proof}

All this is very interesting, and afterwards one can talk about projective algebraic manifolds, having many remarkable properties. However, we will pause our study here, because we still have many other things to say. Getting now to finite fields, we have:

\index{finite field}
\index{Fano plane}

\begin{theorem}
Given a field $F$, we can talk about the projective space $P^{N-1}_F$, as being the space of lines in $F^N$ passing through the origin. At $N=3$ we have
$$|P^2_F|=q^2+q+1$$
where $q=|F|$, in the case where our field $F$ is finite.
\end{theorem}

\begin{proof}
This is indeed clear from definitions, with the cardinality coming from:
$$|P^2_F|
=\frac{|F^3-\{0\}|}{|F-\{0\}|}
=\frac{q^3-1}{q-1}
=q^2+q+1$$

Thus, we are led to the conclusions in the statement.
\end{proof}

As an example, let us see what happens for the simplest finite field that we know, namely $F=\mathbb Z_2$. Here our projective plane, having $4+2+1=7$ points, and 7 lines, is a famous combinatorial object, called Fano plane, which is depicted as follows:
$$
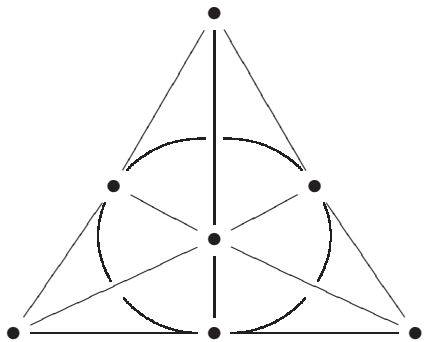
\xymatrix@R=6pt@C=6.5pt{
&&&&\bullet\ar@{-}[ddddd]\\
&&&&&&&\\
&&&&&&&\\
&&&&\ar@{-}@/^/[drr]\\
&&\bullet\ar@{-}[uuuurr]\ar@{-}@/^/[urr]\ar@{-}@/_/[dd]&&&&\bullet\ar@{-}[uuuull]&&\\
&&&&\bullet\ar@{-}[urr]\ar@{-}[ull]&&&&\\
&&\ar@{-}@/_/[drr]&&&&\ar@{-}@/^/[dll]\ar@{-}@/_/[uu]&&&\\
\bullet\ar@{-}[uuurr]\ar@{-}[rrrr]\ar@{-}[uurrrr]&&&&\bullet\ar@{-}[rrrr]\ar@{-}[uu]&&&&\bullet\ar@{-}[uuull]\ar@{-}[uullll]
}$$

Here the circle in the middle is by definition a line, and with this convention, the basic axioms in Definition 11.16 are satisfied, in the sense that any two points determine a line, and any two lines determine a point. And isn't this beautiful. So, let us record:

\begin{conclusion}
For getting started with geometry, all you need is $7$ points.
\end{conclusion}

Finally, no discussion about algebraic geometry would be complete without a look into algebraic topology. We have already seen, in the proof of Theorem 11.17, that ``shape'', taken in a basic topological sense, matters. So, let us further explore this.

\bigskip

Forgetting about manifolds, let us start with something that we know, namely:

\begin{definition}
A topological space $X$ is called connected when any two points $x,y\in X$ can be connected by a path. That is, given any two points $x,y\in X$, we can find a continuous function $f:[0,1]\to X$ such that $f(0)=x$ and $f(1)=y$.
\end{definition}

The problem is now, given a connected space $X$, how to count its ``holes''. And this is quite subtle problem, because as examples of such spaces we have:

\medskip

-- The sphere, the donut, the double-holed donut, the triple-holed donut, and so on. These spaces are quite simple, and intuition suggests to declare that the number of holes of the $N$-holed donut is, and you guessed right, $N$.

\medskip

-- However, we have as well as example the empty sphere, I mean just the crust of the sphere, and while this obviously falls into the class of ``one-holed spaces'', this is not the same thing as a donut, its hole being of different nature.

\medskip

-- As another example, consider again the sphere, but this time with two tunnels drilled into it, in the shape of a cross. Whether that missing cross should account for 1 hole, or for 2 holes, or for something in between, I will leave it up to you. 

\medskip

Summarizing, things are quite tricky, suggesting that the ``number of holes'' of a topological space $X$ is not an actual number, but rather something more complicated. Now with this in mind, let us formulate the following definition:

\index{homotopy group}
\index{hole}
\index{loop}

\begin{definition}
The homotopy group $\pi_1(X)$ of a connected space $X$ is the group of loops based at a given point $*\in X$, with the following conventions,
\begin{enumerate}
\item Two such loops are identified when one can pass continuously from one loop to the other, via a family of loops indexed by $t\in[0,1]$,

\item The composition of two such loops is the obvious one, namely is the loop obtaining by following the first loop, then the second loop, 

\item The unit loop is the null loop at $*$, which stays there, and the inverse of a given loop is the loop itself, followed backwards,
\end{enumerate}
with the remark that the group $\pi_1(X)$ defined in this way does not depend on the choice of the given point $*\in X$, where the loops are based.
\end{definition}

Here the fact that $\pi_1(X)$ defined in this way is indeed a group is obvious, and obvious as well is the fact that, since $X$ is assumed to be connected, this group does not depend on the choice of the given point $*\in X$, where the loops are based.

\bigskip

As basic examples, for spaces having ``no holes'', such as $\mathbb R$ itself, or $\mathbb R^N$, and so on, we have $\pi_1=\{1\}$. In fact, having no holes can only mean, by definition, that $\pi_1=\{1\}$. As further illustrations, here are now a few basic computations:

\index{free group}

\begin{theorem}
We have the following computations of homotopy groups:
\begin{enumerate}
\item For the circle, we obtain $\pi_1=\mathbb Z$.

\item For the torus, we obtain $\pi_1=\mathbb Z\times\mathbb Z$.

\item For the disk minus $2$ points, we have $\pi_1=\mathbb Z*\mathbb Z$.
\end{enumerate}
\end{theorem}

\begin{proof}
These results are all standard, as follows:

\medskip

(1) The first assertion is clear, because a loop on the circle must wind $n\in\mathbb Z$ times around the center, and this parameter $n\in\mathbb Z$ uniquely determines the loop, up to the identification in Definition 11.22. Thus, the homotopy group of the circle is the group of such parameters $n\in\mathbb Z$, which is of course the group $\mathbb Z$ itself.

\medskip

(2) In what regards now the second assertion, the torus being a product of two circles, we are led to the conclusion that its homotopy group must be some kind of product of $\mathbb Z$ with itself. But pictures show that the two standard generators of $\mathbb Z$, and so the two copies of $\mathbb Z$ themselves, commute, $gh=hg$, and so we obtain the product of $\mathbb Z$ with itself, subject to commutation, which is the usual product $\mathbb Z\times\mathbb Z$.

\medskip

(3) This is quite clear, because the homotopy group is generated by the 2 loops around the 2 missing points, which are obviously free, algebrically speaking. Thus, we obtain a free product of the group $\mathbb Z$ with itself, also known as free group on 2 generators.
\end{proof}

There are many other interesting things that can be said about homotopy groups. Also, another thing that can be done with the arbitrary spaces $X$, again in relation with studying their ``shape'', is that of looking at the fiber bundles over them, again up to continuous deformation. We are led in this way into a group, called $K_0(X)$. Moreover, both $\pi_1(X)$ and $K_0(X)$ have in fact higher analogues $\pi_n(X)$ and $K_n(X)$ as well, and the general goal of algebraic topology is that of understanding all these groups.

\bigskip

But all this, obviously, starts to become too complicated. So, leaving now the general manifolds and topological spaces aside, let us focus instead on the simplest objects of topology, namely the knots, with this meaning the smooth closed curves in $\mathbb R^3$:

\index{knot}

\begin{definition}
A knot is a smooth closed curve in $\mathbb R^3$, regarded modulo smooth transformations.
\end{definition}

And isn't this a beautiful definition. We are here at the core of everything that can be called geometry, and in fact, thinking a bit on how knots can be tied, in so many fascinating ways, we are led to the following philosophical conclusion:

\begin{conclusion}
Knots are to geometry and topology what prime numbers are to number theory.
\end{conclusion}

At the level of questions now, once we have a closed curve, say given via its algebraic equations, can we decide if is tied or not, and if tied, how complicated is it tied, how to untie it, and so on? But these are, obviously, quite difficult questions. 

\bigskip

Perhaps simpler now, experience with cables and ropes shows that a random closed curve is usually tied. But can we really prove this? Once again, difficult question. So, we will stop here, and exercise for you, to say something non-trivial about knots.

\section*{11d. Differential manifolds}

Back now to applied mathematics, involving analysis and physics, the situation here is a bit different. Although many interesting algebraic manifolds appear at the advanced level, making algebraic geometry a key tool in advanced physics, in what concerns the basics, here we are mostly in need of a different definition, as follows:

\index{differential manifold}
\index{smooth manifold}
\index{chart}
\index{coordinates}

\begin{definition}
A smooth manifold is a space $X$ which is locally isomorphic to $\mathbb R^N$. To be more precise, this space $X$ must be covered by charts, bijectively mapping open pieces of it to open pieces of $\mathbb R^N$, with the changes of charts being $C^\infty$ functions.
\end{definition}

As a basic example, we have $\mathbb R^N$ itself, or any open subset $X\subset\mathbb R^N$. Another example is the circle, or curves like ellipses and so on, for obvious reasons. To be more precise, the unit circle can be covered by 2 charts as above, by using polar coordinates, in the obvious way, and then by applying dilations, translations and other such transformations, namely bijections which are smooth, we obtain a whole menagery of circle-looking manifolds.  

\bigskip

In particular, we can see from this that Definition 11.26 is a serious rival to Definition 11.12, because both generalize, in a natural way, the conics that we know well. 

\bigskip

Going back now to Definition 11.26 as stated, let us first explore the basic examples. We have already talked about them in the above, and we have:

\begin{proposition}
The following are smooth manifolds, in the plane:
\begin{enumerate}
\item The circles.

\item The ellipses.

\item The non-degenerate conics.

\item Smooth deformations of these.
\end{enumerate}
\end{proposition}

\begin{proof}
All this is quite intuitive, the idea being as follows:

\medskip

(1) Consider the unit circle, $x^2+y^2=1$. We can write then $x=\cos t$, $y=\sin t$, with $t\in[0,2\pi)$, and we seem to have here the solution to our problem, just using 1 chart. But this is of course wrong, because $[0,2\pi)$ is not open, and we have a problem at $0$. In practice we need to use 2 such charts, say with the first one being with $t\in(0,3\pi/2)$, and the second one being with $t\in(\pi,5\pi/2)$. As for the fact that the change of charts is indeed smooth, this comes by writing down the formulae, or just thinking a bit, and arguing that this change of chart being actually a translation, it is automatically linear.

\medskip

(2) This follows from (1), by pulling the circle in both the $Ox$ and $Oy$ directions, and the  formulae here, based on Propositions 11.6 or 11.7, are left to you, reader.

\medskip

(3) We already have the ellipses, and the case of the parabolas and hyperbolas is elementary as well, and in fact simpler than the case of the ellipses. Indeed, a parablola is clearly homeomorphic to $\mathbb R$, and a hyperbola, to two copies of $\mathbb R$.

\medskip

(4) This is something which is clear too, depending of course on what exactly we mean by ``smooth deformation'', and by using a bit of multivariable calculus if needed.
\end{proof}

In higher dimensions now, as basic examples here, we have the unit sphere in $\mathbb R^N$, and smooth deformations of it, once again, somehow by obvious reasons. In case you are wondering on how to construct explicit charts for the sphere, the answer comes from:

\index{polar coordinates}
\index{spherical coordinates}

\begin{theorem}
We have spherical coordinates in $N$ dimensions,
$$\begin{cases}
x_1\!\!\!&=\ r\cos t_1\\
x_2\!\!\!&=\ r\sin t_1\cos t_2\\
\vdots\\
x_{N-1}\!\!\!&=\ r\sin t_1\sin t_2\ldots\sin t_{N-2}\cos t_{N-1}\\
x_N\!\!\!&=\ r\sin t_1\sin t_2\ldots\sin t_{N-2}\sin t_{N-1}
\end{cases}$$
with this guaranteeing that the sphere is indeed a smooth manifold.
\end{theorem}

\begin{proof}
There are several things going on here, the idea being as follows:

\medskip

(1) The fact that we have indeed spherical coordinates is clear, with the only point to be clarified being the identification of the precise ranges of the angles, which follows from some geometric thinking, first at $N=2,3$, and then in general. 

\medskip

(2) As for the last assertion, the manifold one, this can proved a bit like for the circle, as we did in the proof of Proposition 11.27 (1), basically by cutting the sphere into $2^N$ parts, and we will leave the details here as an instructive exercise. 
\end{proof}

In relation with the parametrization question for the spheres, we have the stereographic projection as well, which is something working more directly, as follows:

\index{stereographic projection}

\begin{theorem}
The stereographic projection is given by inverse maps
$$\Phi:\mathbb R^N\to S^N_\mathbb R-\{\infty\}\quad,\quad 
\Psi:S^N_\mathbb R-\{\infty\}\to\mathbb R^N$$
given by the following formulae,
$$\Phi(v)=(1,0)+\frac{2}{1+||v||^2}\,(-1,v)\quad,\quad 
\Psi(c,x)=\frac{x}{1-c}$$
with the convention $\mathbb R^{N+1}=\mathbb R\times\mathbb R^N$, and with the coordinate of $\mathbb R$ denoted $x_0$, and with the coordinates of $\mathbb R^N$ denoted $x_1,\ldots,x_N$.
\end{theorem}

\begin{proof}
We are looking for the formulae of the isomorphism $\mathbb R^N\simeq S^N_\mathbb R-\{\infty\}$, obtained by identifying $\mathbb R^N=\mathbb R^N\times\{0\}\subset\mathbb R^{N+1}$ with the unit sphere $S^N_\mathbb R\subset\mathbb R^{N+1}$, with the convention that the point which is added is $\infty=(1,0,\ldots,0)$, via the stereographic projection. That is, we need the precise formulae of two inverse maps, as follows:
$$\Phi:\mathbb R^N\to S^N_\mathbb R-\{\infty\}\quad,\quad 
\Psi:S^N_\mathbb R-\{\infty\}\to\mathbb R^N$$

In one sense, according to our conventions above, we must have a formula as follows for our map $\Phi$, with the parameter $t\in(0,1)$ being such that $||\Phi(v)||=1$:
$$\Phi(v)=t(0,v)+(1-t)(1,0)$$

The equation for the parameter $t\in(0,1)$ can be solved as follows:
$$(1-t)^2+t^2||v||^2=1
\iff t=\frac{2}{1+||v||^2}$$

Thus, we are led to the formula in the statement for $\Phi$. In the other sense now, we must have a formula as follows, for a certain $\alpha\in\mathbb R$:
$$(0,\Psi(c,x))=\alpha(c,x)+(1-\alpha)(1,0)$$

But from $\alpha c+1-\alpha=0$ we get the following formula for the parameter $\alpha$:
$$\alpha=\frac{1}{1-c}$$

Thus, we are led to the formula in the statement for $\Psi$.
\end{proof}

There are of course many other possible parametrizations of the sphere, such as the one using cylindrical coordinates. As an interesting question here, we have:

\index{cartography}

\begin{question}
What is the best parametrization of the unit sphere in $\mathbb R^3$, for purely mathematical reasons? What about for cartography reasons?
\end{question}

To be more precise, the first question, which is quite challenging, is that of finding the simplest proof ever for the fact that the sphere in $\mathbb R^3$ is indeed a smooth manifold. As for the second question, which is even more challenging, the problem here is to have your charts reflecting as nicely as possible useful things such as lengths, angles and areas.

\bigskip

By the way, speaking lengths, angles and areas, observe that the general differential manifold formalism from Definition 11.26 is obviously too broad for talking about these. In order to do so, several more axioms must be added, and we end up with something called Riemannian manifold, which is the main object of study of advanced differential geometry. And, regarding such manifolds, there are many deep theorems, including:

\index{Nash theorem}

\begin{theorem}[Nash]
Any Riemannian manifold embeds as $X\subset\mathbb R^N$.
\end{theorem}

\begin{proof}
This is something quite technical. Note however that, philosophically, this theorem suggests ignoring itself, and everything advanced, that is, pick for your geometry works good old manifolds $X\subset\mathbb R^N$, and study them via multivariable calculus. 
\end{proof}

So long for geometry, in a large sense. As a conclusion to all this, geometry comes in many flavors, algebraic or differential, Riemannian or not, over $\mathbb R$, or $\mathbb C$, or some other field $F$, and in addition to this we can talk about affine or projective geometry, or about discrete or continuous geometry, and so on. Many interesting things, and if excited by all this, orient yourself towards physics, where geometry in all its flavors is needed.

\section*{11e. Exercises}

This was a pure mathematics chapter, save for some old school computations coming from Kepler and Newton, and as exercises, all of pure mathematics type, we have: 

\begin{exercise}
Learn some algebraic geometry.
\end{exercise}

\begin{exercise}
Learn some differential geometry.
\end{exercise}

\begin{exercise}
Learn some Riemannian geometry.
\end{exercise}

\begin{exercise}
Learn some symplectic geometry too.
\end{exercise}

With apologies for this, but that is how things go in pure mathematics, ``learn'' being the keyword. As bonus exercise, more concrete, have some fun with cartography.

\chapter{Higher derivatives}

\section*{12a. Higher derivatives}

Welcome back to analysis. As another continuation of the basic multivariable calculus from chapter 10, we will go here for the real thing, namely the computation of maxima and minima of functions. Want to optimize the characteristics of your building, bridge, plane, engine, medicine, computer, or casino and other racketeering operations? It all comes down, you guessed right, to maximizing or minimizing a certain function.

\bigskip

Let us start with a reminder of what we know in 1 variable, from chapter 3:

\begin{theorem}
The one-variable functions are subject to the Taylor formula
$$f(x+t)\simeq f(x)+f'(x)t+\frac{f''(x)}{2}\,t^2+\ldots$$
which allows, via suitable truncations, to determine the local maxima and minima. 
\end{theorem}

\begin{proof}
This is something very compact, the idea being as follows:

\medskip

(1) To start with, in what regards the precise formulation of the Taylor formula, involving a function $f:X\to\mathbb R$, with $X\subset\mathbb R$ open, and then the differentiability of this function, at order 1, 2, 3 and so on, we refer to the material in chapter 3.

\medskip

(2) In order to compute the local maxima and minima of $f$, the idea is that we have an algorithm. At the first step, assuming that $f$ is differentiable, we can use:
$$f(x+t)\simeq f(x)+f'(x)t$$

Indeed, this formula shows that when $f'(x)\neq0$, the point $x$ cannot be a local maximum or minimum, due to the fact that $t\to-t$ will invert the growth. Thus, in order to find the local maxima and minima, we must compute first the points $x$ satisfying $f'(x)=0$, and then perform a more detailed study of each solution $x$ that we found. 

\medskip

(3) Next, assuming that we have $f'(x)=0$, as required by the study in (2), and that $f$ is twice differentiable at $x$, the second order Taylor formula there reads:
$$f(x+t)\simeq f(x)+\frac{f''(x)}{2}\,t^2$$

But this is something very useful, telling us that when $f''(x)<0$, what we have is a local maximum, and when $f''(x)>0$, what we have is a local minimum. As for the remaining case, that when $f''(x)=0$, things here remain open.

\medskip

(4) Assume now that we are in the case $f'(x)=f''(x)=0$, which is where our joint algorithm coming from (2) and (3) fails. In this case, we can go at order 3:
$$f(x+t)\simeq f(x)+\frac{f'''(x)}{6}\,t^3$$

But this solves the problem in the case $f'''(x)\neq0$, because here we cannot have a local minimum or maximum, due to $t\to-t$ which switches growth. As for the remaining case, $f'''(x)=0$, things here remain open, and we have to go at higher order.

\medskip

(5) Summarizing, we have a recurrence method for solving our problem, with the conclusion being that, with $n\in\mathbb N$ being minimal such that $f^{(n)}(x)\neq0$, we can use:
$$f(x+t)\simeq f(x)+\frac{f^{(n)}(x)}{n!}\,t^n$$

Indeed, local extremum needs this $n$ to be even, with in this case $f^{(n)}(x)<0$ corresponding to a local maximum, and $f^{(n)}(x)>0$ corresponding to a local minimum.
\end{proof}

Wiht this discussed, what about several variables? Browsing through what we did in chapter 10, first derivatives, we are at step 1 of the above algorithm, as shown by:

\begin{theorem}
In order for a differentiable function $f:\mathbb R^N\to\mathbb R$ to have a local maximum or minimum at $x\in\mathbb R^N$, the derivative must vanish there, $f'(x)=0$.
\end{theorem}

\begin{proof}
This comes, as in the 1-variable case, from the following formula:
$$f(x+t)\simeq f(x)+f'(x)t$$

Indeed, the error term $f'(x)t$ being a certain linear combination of the entries of $t\in\mathbb R^N$, this linear combination must vanish, $f'(x)=0$. Alternatively, and for being a bit more explicit, according to our theory from chapter 10, the above formula reads:
$$f(x+t)\simeq f(x)+\sum_{i=1}^N\frac{df}{dx_i}\,t_i$$

Thus, all the partial derivatives must vanish at $x$, telling us that $f'(x)=0$.
\end{proof}

Getting now to order 2 and higher, we have the following question, to be solved:

\begin{question}
What is the second derivative of $f:\mathbb R^N\to\mathbb R$, making the formula
$$f(x+t)\simeq f(x)+f'(x)t+\frac{f''(x)}{2}*t*t$$
work, with a suitable meaning for $*$? What about the higher order derivatives?
\end{question}

To be more precise here, this question is something quite self-explanatory, and with that $*$ symbols being used, instead of usual products, because on the right we can certainly not multiply a vector $t\in\mathbb R^N$ by itself, just like that, and on the left, who really knows in advance what that $f''(x)$ beast can be, and how can that multiply by $t*t$.

\bigskip

In answer now, we will follow the pedestrian way, a bit as in chapter 10, when talking about first derivatives. On the menu, second and higher order partial derivatives, that we actually already met, in chapter 8 when doing physics, then some algebra in order to make sense of that mysterious $*$ operation, then Taylor formula and applications.

\bigskip

Getting started now, we can talk about higher derivatives, in the obvious way, simply by performing the operation of taking derivatives recursively. As result here, we have:

\index{higher derivative}

\begin{theorem}
Given $f:\mathbb R^N\to\mathbb R$, we can talk about its higher derivatives
$$\frac{d^kf}{dx_{i_1}\ldots dx_{i_k}}=\frac{d}{dx_{i_1}}\cdots\frac{d}{dx_{i_k}}(f)$$
provided that these derivatives exist indeed. Moreover, due to the Clairaut formula,
$$\frac{d^2f}{dx_idx_j}=\frac{d^2f}{dx_jdx_i}$$
the order in which these higher derivatives are computed is irrelevant.
\end{theorem}

\begin{proof}
There are several things going on here, the idea being as follows:

\medskip

(1) First of all, we can talk about the quantities in the statement, with the remark however that at each step of our recursion, the corresponding partial derivative can exist of not. We will say in what follows that our function is $n$ times differentiable if the quantities in the statement exist at any $k\leq n$, and smooth, if this works with $n=\infty$.

\medskip

(2) Regarding now the second assertion, this is something more tricky. Let us first recall from chapter 8 that the second derivatives of a twice differentiable function of two variables $f:\mathbb R^2\to\mathbb R$ are subject to the Clairaut formula, namely:
$$\frac{d^2f}{dxdy}=\frac{d^2f}{dydx}$$

(3) But this result clearly extends to our function $f:\mathbb R^N\to\mathbb R$, simply by ignoring the unneeded variables, so we have the Clairaut formula in general, also called Schwarz formula, which is the one in the statement, namely:
$$\frac{d^2f}{dx_idx_j}=\frac{d^2f}{dx_jdx_i}$$

(4) Now observe that this tells us that the order in which the higher derivatives are computed is irrelevant. That is, we can permute the order of our partial derivative computations, and a standard way of doing this is by differentiating first with respect to $x_1$, as many times as needed, then with respect to $x_2$, and so on. Thus, the collection of partial derivatives can be written, in a more convenient form, as follows:
$$\frac{d^kf}{dx_1^{k_1}\ldots dx_N^{k_N}}=\frac{d^{k_1}}{dx_1^{k_1}}\cdots\frac{d^{k_N}}{dx_N^{k_N}}(f)$$

(5) To be more precise, here $k\in\mathbb N$ is as usual the global order of our derivatives, the exponents $k_1,\ldots,k_N\in\mathbb N$ are subject to the condition $k_1+\ldots+k_N=k$, and the operations on the right are the familiar one-variable higher derivative operations.
\end{proof}

All this is very nice, and as an illustration for the above, let us work out the case $k=2$. Here things are quite special, and we can formulate the following definition:

\begin{definition}
Given a twice differentiable function $f:\mathbb R^N\to\mathbb R$, we set
$$f''(x)=\left(\frac{d^2f}{dx_idx_j}\right)_{ij}$$
which is a symmetric matrix, called Hessian matrix of $f$ at the point $x\in\mathbb R^N$.
\end{definition}

To be more precise, we know that when $f:\mathbb R^N\to\mathbb R$ is twice differentiable, its order $k=2$ partial derivatives are the numbers in the statement. Now since these numbers naturally form a $N\times N$ matrix, the temptation is high to call this matrix $f''(x)$, and so we will do. And finally, we know from Clairaut that this matrix is symmetric:
$$f''(x)_{ij}=f''(x)_{ji}$$

Observe that at $N=1$ this is compatible with our previous definition of the second derivative $f''$, and this because in this case, the $1\times1$ matrix from Definition 12.5 is:
$$f''(x)=(f''(x))\in M_{1\times1}(\mathbb R)$$

As a word of warning, however, never use Definition 12.5 for functions $f:\mathbb R^N\to\mathbb R^M$, where the second derivative can only be something more complicated. Also, never attempt either to do something similar at $k=3$ or higher, for functions $f:\mathbb R^N\to\mathbb R$ with $N>1$, because again, that beast has too many indices, for being a true, honest matrix.

\bigskip

Back now to business, Taylor formula at order 2, we have here:

\index{second order derivative}
\index{Hessian}
\index{positive matrix}
\index{minimum}
\index{maximum}
\index{local minimum}
\index{local maximum}

\begin{theorem}
Given a twice differentiable function $f:\mathbb R^N\to\mathbb R$, we have
$$f(x+t)\simeq f(x)+f'(x)t+\frac{<f''(x)t,t>}{2}$$
where $f''(x)\in M_N(\mathbb R)$ stands as usual for the Hessian matrix.
\end{theorem}

\begin{proof}
This is something very standard, the idea being as follows:

\medskip

(1) As a first observation, at $N=1$ the Hessian matrix as constructed in Definition 12.5 is the $1\times1$ matrix having as entry the second derivative $f''(x)$, and the formula in the statement is something that we know well from chapter 3, namely:
$$f(x+t)\simeq f(x)+f'(x)t+\frac{f''(x)t^2}{2}$$

(2) In general now, this is in fact something which does not need a new proof, because it follows from the one-variable formula above, applied to the restriction of $f$ to the following segment in $\mathbb R^N$, which can be regarded as being a one-variable interval:
$$I=[x,x+t]$$

To be more precise, let $y\in\mathbb R^N$, and consider the following function, with $r\in\mathbb R$:
$$g(r)=f(x+ry)$$

We know from (1) that the Taylor formula for $g$, at the point $r=0$, reads:
$$g(r)\simeq g(0)+g'(0)r+\frac{g''(0)r^2}{2}$$

And our claim is that, with $t=ry$, this is precisely the formula in the statement.

\medskip

(3) So, let us see if our claim is correct. By using the chain rule, we have the following formula, with on the right, as usual, a row vector multiplied by a column vector:
$$g'(r)=f'(x+ry)\cdot y$$

By using again the chain rule, we can compute the second derivative as well:
\begin{eqnarray*}
g''(r)
&=&(f'(x+ry)\cdot y)'\\
&=&\left(\sum_i\frac{df}{dx_i}(x+ry)\cdot y_i\right)'\\
&=&\sum_i\sum_j\frac{d^2f}{dx_idx_j}(x+ry)\cdot\frac{d(x+ry)_j}{dr}\cdot y_i\\
&=&\sum_i\sum_j\frac{d^2f}{dx_idx_j}(x+ry)\cdot y_iy_j\\
&=&<f''(x+ry)y,y>
\end{eqnarray*}

(4) Time now to conclude. We know that we have $g(r)=f(x+ry)$, and according to our various computations above, we have the following formulae:
$$g(0)=f(x)\quad,\quad 
g'(0)=f'(x)\quad,\quad 
g''(0)=<f''(x)y,y>$$

Buit with this data in hand, the usual Taylor formula for our one variable function $g$, at order 2, at the point $r=0$, takes the following form, with $t=ry$:
\begin{eqnarray*}
f(x+ry)
&\simeq&f(x)+f'(x)ry+\frac{<f''(x)y,y>r^2}{2}\\
&=&f(x)+f'(x)t+\frac{<f''(x)t,t>}{2}
\end{eqnarray*}

Thus, we have obtained the formula in the statement.

\medskip

(5) Finally, for completness, let us record as well a more concrete formulation of what we found. According to our usual rules for matrix calculus, what we found is:
$$f(x+t)\simeq f(x)+\sum_{i=1}^N\frac{df}{dx_i}\,t_i+\frac{1}{2}\sum_{i=1}^N\sum_{j=1}^N\frac{d^2f}{dx_idx_j}\,t_it_j$$

Observe that, since the Hessian matrix $f''(x)$ is symmetric, most of the terms on the right will appear in pairs, making it clear what the $1/2$ rescaling is there for, namely avoiding redundancies. However, this is only true for the off-diagonal terms, so instead of further messing up our formula above, we will just leave it like this.
\end{proof}

We can now go back to local extrema, and we have, improving Theorem 12.2:

\begin{theorem}
In order for a twice differentiable function $f:\mathbb R^N\to\mathbb R$ to have a local minimum or maximum at $x\in\mathbb R^N$, the first derivative must vanish there,
$$f'(x)=0$$
and the Hessian must be positive or negative, in the sense that the quantities
$$<f''(x)t,t>\in\mathbb R$$
must keep a constant sign, positive or negative, when $t\in\mathbb R^N$ varies.
\end{theorem}

\begin{proof}
This is clear from Theorem 12.6. Consider indeed the formula there, namely:
$$f(x+t)\simeq f(x)+f'(x)t+\frac{<f''(x)t,t>}{2}$$

As explained in Theorem 12.2, in order for our function to have a local minimum or maximum at $x\in\mathbb R^N$, the first derivative must vanish there, $f'(x)=0$. And then, with this assumption made, the approximation that we have around $x$ becomes:
$$f(x+t)\simeq f(x)+\frac{<f''(x)t,t>}{2}$$

Now assuming that $x$ is a local minimum, we must have $<f''(x)t,t>\,\geq0$, and assuming that $x$ is a local maximum, we must have $<f''(x)t,t>\,\leq0$, as stated.
\end{proof}

The above result is not the end of the story at order 2, because we still have to talk about the converse, depending on how $f''(x)$ exactly acts. So, let us formulate:

\begin{definition}
Given a symmetric matrix $A\in M_N(\mathbb R)$, we write:
\begin{enumerate}
\item $A\geq0$, when $<At,t>\,\geq0$, for any $t\in\mathbb R^N$.

\item $A>0$, when $<At,t>\,>0$, for any $t\neq0$.

\item $A\leq0$, when $<At,t>\,\leq0$, for any $t\in\mathbb R^N$.

\item $A<0$, when $<At,t>\,<0$, for any $t\neq0$.
\end{enumerate}
\end{definition}

Obviously, some linear algebra is going on here, and you would probably say that, in terms of eigenvalues, the above cases should read $\lambda_i\geq0$, $\lambda_i>0$, $\lambda_i\leq0$, $\lambda_i<0$. Which is indeed the case, but for having all this clarified we still need to show that $A$ is diagonalizable, and with the basis change leaving invariant the quantities $<At,t>$. Summarizing, some work to be done here, and we will be back to this in a moment.   

\bigskip

In the meantime, with this in hand, we can complement Theorem 12.7 with:

\begin{theorem}[continuation]
Given a twice differentiable function $f:\mathbb R^N\to\mathbb R$ as above, and assuming that we have a point $x\in\mathbb R^N$ where $f'(x)=0$:
\begin{enumerate}
\item $f''(x)\geq0$ is needed for $x$ to be a local minimum.

\item $f''(x)>0$ guarantees that $x$ is a local minimum.

\item $f''(x)\leq0$ is needed for $x$ to be a local maximum.

\item $f''(x)<0$ guarantees that $x$ is a local maximum.
\end{enumerate}
\end{theorem}

\begin{proof}
This comes from Theorem 12.7, and the study from its proof, with the various positivity notions being those from Definition 12.8. To be more precise, consider the Taylor formula at order 2, under the critical point assumption $f'(x)=0$, namely:
$$f(x+t)\simeq f(x)+\frac{<f''(x)t,t>}{2}$$

(1) This is something clear, that we already know, from Theorem 12.7.

\medskip

(2) This comes as a variation of (1), because when the second order term is strictly positive, $<f''(x)t,t>\,>0$ for $t\neq0$, with this term overcoming anyway the error, we obtain a strict estimate, $f(x+t)>f(x)$ for $t\simeq0$ small, so we have a local minimum.

\medskip

(3-4) The discussion here is similar, with all inequalities being reversed.
\end{proof}

All the above looks quite good, nicely generalizing what we knew from 1-variable calculus, save for the following key fact, which is multivariable calculus specific:

\begin{warning}
Contrary to the real numbers, which can be positive, negative or zero, the symmetric matrices can be strictly positive $A>0$, strictly negative $A<0$, or many other things in between, and only the following cases stop our computation at order $2$:
\begin{enumerate}
\item $A>0$, where we have a local minimum, which is strict.

\item $A<0$, where we have a local maximum, which is strict.

\item $A\not\geq0$ and $A\not\leq0$, where we cannot have a local minimum or maximum.
\end{enumerate}
That is, the bad cases, where we need order $3$, are $A\geq0$, $A\not>0$, and $A\leq0$, $A\not<0$.
\end{warning}

And more on this in a moment, in terms of eigenvalues, after learning more linear algebra. Finally, at higher order things become more complicated, as follows:

\begin{theorem}
Given a higher order differentiable function $f:\mathbb R^N\to\mathbb R$, we have
$$f(x+t)\simeq f(x)+f'(x)t+\frac{<f''(x)t,t>}{2}+\ldots$$
and this helps in identifying the local extrema, when $f'(x)=0$ and $f''(x)=0$.
\end{theorem}

\begin{proof}
The study here is very similar to that at $k=2$, from the proof of Theorem 12.6, with everything coming from the usual Taylor formula, applied on:
$$I=[x,x+t]$$

And we will leave this as an exercise, for your long Summer nights. Actually, having the order 3 formula fully worked out would be already a great thing. Enjoy.
\end{proof}

\section*{12b. Matrices, positivity}

Time now for some advanced linear algebra, in relation with the various notions from Definition 12.8. As a main question that we would like to solve, we have:

\begin{question}
Which symmetric matrices $A\in M_N(\mathbb R)$ have the property
$$<At,t>\,\geq0$$
for any $t\in\mathbb R^N$? Also, when is a matrix $A\in M_N(\mathbb R)$ diagonalizable?
\end{question}

Observe that in the second question we have lifted the assumption that $A$ is symmetric, and this for good reason. Indeed, besides the local extremum problematics for the functions $f:\mathbb R^N\to\mathbb R$, discussed above, we have as well the quite interesting problem of studying the functions $f:\mathbb R^N\to\mathbb R^N$, whose derivatives are square matrices, $f'(x)\in M_N(\mathbb R)$, waiting to be diagonalized too. So, we are trying here two shoot 2 rabbits at the same time, with our Question 12.12, as formulated above. 

\bigskip

Obviously, hunting matters, so time to ask the expert. And the expert says:

\begin{cat}
In the lack of speed and claws, yes, develop as much linear algebra as you can. Without even caring for applications, these will come naturally.
\end{cat}

Thanks cat, so this will be our plan for this section, develop as much linear algebra as we can. And for applications, we will most likely leave them to you, reader, for later in life, depending on the precise physics and engineering questions that you will be interested in. And with the advice of course to follow the feline way, relax, and no mercy.

\bigskip

With this plan made, let us go back to the diagonalization question, from chapter 9. We will need here diagonalization results which are far more powerful. We first have:

\index{self-adjoint matrix}

\begin{theorem}
Any matrix $A\in M_N(\mathbb C)$ which is self-adjoint, $A=A^*$, is diagonalizable, with the diagonalization being of the following type,
$$A=UDU^*$$
with $U\in U_N$, meaning $U^*=U^{-1}$, and $D\in M_N(\mathbb R)$ diagonal. The converse holds too.
\end{theorem}

\begin{proof}
As a first remark, the converse trivially holds, because if we take a matrix of the form $A=UDU^*$, with $U$ unitary and $D$ diagonal and real, then we have:
$$A^*
=(UDU^*)^*
=UD^*U^*
=UDU^*
=A$$

In the other sense now, assume that $A$ is self-adjoint, $A=A^*$.  Our first claim is that the eigenvalues are real. Indeed, assuming $Av=\lambda v$, we have:
\begin{eqnarray*}
\lambda<v,v>
&=&<Av,v>\\
&=&<v,Av>\\
&=&\bar{\lambda}<v,v>
\end{eqnarray*}

Thus we obtain $\lambda\in\mathbb R$, as claimed. Our next claim now is that the eigenspaces corresponding to different eigenvalues are pairwise orthogonal. Assume indeed that:
$$Av=\lambda v\quad,\quad 
Aw=\mu w$$

We have then the following computation, using the fact that we have $\mu\in\mathbb R$:
\begin{eqnarray*}
\lambda<v,w>
&=&<Av,w>\\
&=&<v,Aw>\\
&=&\mu<v,w>
\end{eqnarray*}

Thus $\lambda\neq\mu$ implies $v\perp w$, as claimed. In order now to finish the proof, it remains to prove that the eigenspaces of $A$ span the whole space $\mathbb C^N$. For this purpose, we will use a recurrence method. Let us pick an eigenvector of our matrix:
$$Av=\lambda v$$

Assuming now that we have a vector $w$ orthogonal to it, $v\perp w$, we have:
\begin{eqnarray*}
<Aw,v>
&=&<w,Av>\\
&=&<w,\lambda v>\\
&=&\lambda<w,v>\\
&=&0
\end{eqnarray*}

Thus, if $v$ is an eigenvector, then the vector space $v^\perp$ is invariant under $A$. Moreover, since a matrix $A$ is self-adjoint precisely when $<Av,v>\in\mathbb R$ for any vector $v\in\mathbb C^N$, as one can see by expanding the scalar product, the restriction of $A$ to the subspace $v^\perp$ is self-adjoint. Thus, we can proceed by recurrence, and we obtain the result.
\end{proof}

As an important consequence of the above result, we have:

\begin{theorem}
Any matrix $A\in M_N(\mathbb R)$ which is symmetric, $A=A^t$, is diagonalizable, with the diagonalization being of the following type,
$$A=UDU^t$$
with $U\in O_N$, meaning $U^t=U^{-1}$, and $D\in M_N(\mathbb R)$ diagonal. The converse holds too.
\end{theorem}

\begin{proof}
As before, the converse trivially holds, because if we take a matrix of the form $A=UDU^t$, with $U$ orthogonal and $D$ diagonal and real, then we have $A^t=A$. In the other sense now, this follows from Theorem 12.14, and its proof.
\end{proof}

Getting now to positivity issues, the basic theory here is as follows:

\index{positive matrix}
\index{square root}
\index{positive eigenvalues}

\begin{theorem}
For a matrix $A\in M_N(\mathbb C)$ the following conditions are equivalent, and if they are satisfied, we say that $A$ is positive, and write $A\geq0$:
\begin{enumerate}
\item $A=B^2$, with $B=B^*$.

\item $A=CC^*$, for some $C\in M_N(\mathbb C)$.

\item $<Ax,x>\geq0$, for any vector $x\in\mathbb C^N$.

\item $A=A^*$, and the eigenvalues are positive, $\lambda_i\geq0$.

\item $A=UDU^*$, with $U\in U_N$ and with $D\in M_N(\mathbb R_+)$ diagonal.
\end{enumerate}
\end{theorem}

\begin{proof}
The idea is that the equivalences in the statement basically follow from some elementary computations, with only Theorem 12.14 being needed, at some point:

\medskip

$(1)\implies(2)$ This is clear, because we can take $C=B$.

\medskip

$(2)\implies(3)$ This follows from the following computation:
\begin{eqnarray*}
<Ax,x>
&=&<CC^*x,x>\\
&=&<C^*x,C^*x>\\
&\geq&0
\end{eqnarray*}

$(3)\implies(4)$ By using the fact that $<Ax,x>$ is real, we have:
$$<Ax,x>
=<x,A^*x>
=<A^*x,x>$$

Thus we have $A=A^*$, and the remaining assertion, regarding the eigenvalues, follows from the following computation, assuming $Ax=\lambda x$:
$$<Ax,x>
=<\lambda x,x>
=\lambda<x,x>
\geq0$$

$(4)\implies(5)$ This follows indeed by using Theorem 12.14.

\medskip

$(5)\implies(1)$ Assuming $A=UDU^*$ with $U\in U_N$, and with $D\in M_N(\mathbb R_+)$ diagonal, we can set $B=U\sqrt{D}U^*$. Then $B$ is self-adjoint, and its square is given by:
$$B^2
=U\sqrt{D}U^*\cdot U\sqrt{D}U^*
=UDU^*
=A$$

Thus, we are led to the conclusion in the statement.
\end{proof}

Let us record as well the following version of the above result:

\index{strictly positive matrix}

\begin{theorem}
For a matrix $A\in M_N(\mathbb C)$ the following conditions are equivalent, and if they are satisfied, we say that $A$ is strictly positive, and write $A>0$:
\begin{enumerate}
\item $A=B^2$, with $B=B^*$, invertible.

\item $A=CC^*$, for some $C\in M_N(\mathbb C)$ invertible.

\item $<Ax,x>>0$, for any nonzero vector $x\in\mathbb C^N$.

\item $A=A^*$, and the eigenvalues are strictly positive, $\lambda_i>0$.

\item $A=UDU^*$, with $U\in U_N$ and with $D\in M_N(\mathbb R_+^*)$ diagonal.
\end{enumerate}
\end{theorem}

\begin{proof}
This follows either from Theorem 12.16, by adding the above various extra assumptions, or from the proof of Theorem 12.16, by modifying where needed.
\end{proof}

And with this, good news, we can pause our linear algebra study, go back to our analysis questions, and formulate the following final result, regarding order 2:

\begin{theorem}
Given a twice differentiable function $f:\mathbb R^N\to\mathbb R$, assume that we have a point where $f'(x)=0$, and let $\lambda_1,\ldots,\lambda_N$ be the eigenvalues of $f''(x)$. Then:
\begin{enumerate}
\item $\lambda_i\geq0$ is needed for $x$ to be a local minimum.

\item $\lambda_i>0$ guarantees that $x$ is a local minimum.

\item $\lambda_i\leq0$ is needed for $x$ to be a local maximum.

\item $\lambda_i<0$ guarantees that $x$ is a local maximum.
\end{enumerate}
\end{theorem}

\begin{proof}
This is indeed a reformulation of Theorem 12.9, by using the above linear algebra material. In short, and for simplifying a bit, the idea is as follows:

\medskip

(1) We know from Theorem 12.15 that the Hessian matrix $f''(x)$ is diagonalized by a certain orthogonal matrix $U\in O_N$. But with this in hand, we can simply change the basis of $\mathbb R^N$, with the help of this matrix $U\in O_N$, and the Taylor formula becomes:
$$f(x+t)\simeq f(x)+\sum_{i=1}^N\lambda_it_i^2$$

And this latter formula, obviously, gives all the assertions in the statement.

\medskip

(2) This was for the theory, but in practice, there are some other things that can be useful. Consider for instance a function $f:\mathbb R^2\to\mathbb R$, whose Hessian looks as follows:
$$f''(x)=\begin{pmatrix}a&b\\ c&d\end{pmatrix}$$

The eigenvalues are then given by the following trace and determinant equations:
$$\lambda_1+\lambda_2=a+d\quad,\quad\lambda_1\lambda_2=ad-bc$$

Thus, without even computing the einegvalues, we can say right away, depending on the signs of $a+d$, $ad-bc$, if we are in one of the situations (1,2,3,4) in the statement.

\medskip

(3) In more dimensions things are more complicated, but there are still tricks, that can help, and the more you learn and know here, the better your analysis will be.
\end{proof}

With this discussed, shall we stop here, with our linear algebra study? In view of Question 12.12 and Cat 12.13, I would say no, remember indeed that we still have the first derivatives $f'(x)\in M_N(\mathbb R)$ of functions $f:\mathbb R^N\to\mathbb R^N$, to be better understood. 

\bigskip

So, as a continuation of our previous diagonalization results, known as ``spectral theorems'', let us discuss now the diagonalization of unitary matrices. We have here:

\index{unitary}

\begin{theorem}
Any matrix $U\in M_N(\mathbb C)$ which is unitary, $U^*=U^{-1}$, is diagonalizable, with the eigenvalues on $\mathbb T$. More precisely we have
$$U=VDV^*$$
with $V\in U_N$, and with $D\in M_N(\mathbb T)$ diagonal. The converse holds too.
\end{theorem}

\begin{proof}
As a first remark, the converse trivially holds, because given a matrix of type $U=VDV^*$, with $V\in U_N$, and with $D\in M_N(\mathbb T)$ being diagonal, we have:
$$U^*
=V\bar{D}V^*
=(VDV^*)^{-1}
=U^{-1}$$

Let us prove now the first assertion, stating that the eigenvalues of a unitary matrix $U\in U_N$ belong to $\mathbb T$. Indeed, assuming $Uv=\lambda v$, we have:
\begin{eqnarray*}
<v,v>
&=&<U^*Uv,v>\\
&=&<Uv,Uv>\\
&=&<\lambda v,\lambda v>\\
&=&|\lambda|^2<v,v>
\end{eqnarray*}

Thus we obtain $\lambda\in\mathbb T$, as claimed. Our next claim now is that the eigenspaces corresponding to different eigenvalues are pairwise orthogonal. Assume indeed that:
$$Uv=\lambda v\quad,\quad 
Uw=\mu w$$

We have then the following computation, using $U^*=U^{-1}$ and $\lambda,\mu\in\mathbb T$:
\begin{eqnarray*}
\lambda<v,w>
&=&<Uv,w>\\
&=&<v,U^*w>\\
&=&<v,U^{-1}w>\\
&=&<v,\mu^{-1}w>\\
&=&\mu<v,w>
\end{eqnarray*}

Thus $\lambda\neq\mu$ implies $v\perp w$, as claimed. In order now to finish the proof, it remains to prove that the eigenspaces of $U$ span the whole space $\mathbb C^N$. For this purpose, we will use a recurrence method. Let us pick an eigenvector of our matrix:
$$Uv=\lambda v$$

Assuming that we have a vector $w$ orthogonal to it, $v\perp w$, we have:
\begin{eqnarray*}
<Uw,v>
&=&<w,U^*v>\\
&=&<w,U^{-1}v>\\
&=&<w,\lambda^{-1}v>\\
&=&\lambda<w,v>\\
&=&0
\end{eqnarray*}

Thus, if $v$ is an eigenvector, then the vector space $v^\perp$ is invariant under $U$. Now since $U$ is an isometry, so is its restriction to this space $v^\perp$. Thus this restriction is a unitary, and so we can proceed by recurrence, and we obtain the result.
\end{proof}

Let us record as well the real version of the above result, in a weak form:

\begin{proposition}
Any matrix $U\in M_N(\mathbb R)$ which is orthogonal, $U^t=U^{-1}$,
is diagonalizable, with the eigenvalues on $\mathbb T$. More precisely we have
$$U=VDV^*$$
with $V\in U_N$, and with $D\in M_N(\mathbb T)$ being diagonal.
\end{proposition}

\begin{proof}
This comes indeed as a trivial consequence of Theorem 12.19.
\end{proof}

Observe that the above result does not provide us with a complete characterization of the matrices $U\in M_N(\mathbb R)$ which are orthogonal. To be more precise, the question left is that of understanding when the matrices of type $U=VDV^*$, with $V\in U_N$, and with $D\in M_N(\mathbb T)$ being diagonal, are real, and this is something non-trivial.

\bigskip

As an illustration, for the simplest unitaries that we know, namely the rotations in the real plane, we have the following formula, that we know well from chapter 9:
$$\begin{pmatrix}\cos t&-\sin t\\ \sin t&\cos t\end{pmatrix}
=\frac{1}{2}\begin{pmatrix}1&1\\i&-i\end{pmatrix}
\begin{pmatrix}e^{-it}&0\\0&e^{it}\end{pmatrix}
\begin{pmatrix}1&-i\\1&i\end{pmatrix}$$

Back to generalities, the self-adjoint matrices and the unitary matrices are particular cases of the general notion of a ``normal matrix'', and we have here:

\index{normal matrix}
\index{diagonalizable matrix}
\index{spectral theorem}

\begin{theorem}
Any matrix $A\in M_N(\mathbb C)$ which is normal, $AA^*=A^*A$, is diagonalizable, with the diagonalization being of the following type,
$$A=UDU^*$$
with $U\in U_N$, and with $D\in M_N(\mathbb C)$ diagonal. The converse holds too.
\end{theorem}

\begin{proof}
As a first remark, the converse trivially holds, because if we take a matrix of the form $A=UDU^*$, with $U$ unitary and $D$ diagonal, then we have:
$$AA^*
=UDD^*U^*
=UD^*DU^*
=A^*A$$

In the other sense now, this is something more technical. Our first claim is that a matrix $A$ is normal precisely when the following happens, for any vector $v$:
$$||Av||=||A^*v||$$

Indeed, the above equality can be written as follows:
$$<AA^*v,v>=<A^*Av,v>$$

But this is equivalent to $AA^*=A^*A$, by expanding the scalar products. Our claim now is that $A,A^*$ have the same eigenvectors, with conjugate eigenvalues:
$$Av=\lambda v\implies A^*v=\bar{\lambda}v$$

Indeed, this follows from the following computation, and from the trivial fact that if $A$ is normal, then so is any matrix of type $A-\lambda 1_N$:
\begin{eqnarray*}
||(A^*-\bar{\lambda}1_N)v||
&=&||(A-\lambda 1_N)^*v||\\
&=&||(A-\lambda 1_N)v||\\
&=&0
\end{eqnarray*}

Let us prove now, by using this, that the eigenspaces of $A$ are pairwise orthogonal. Assume that we have two eigenvectors, corresponding to different eigenvalues, $\lambda\neq\mu$:
$$Av=\lambda v\quad,\quad 
Aw=\mu w$$

We have the following computation, which shows that $\lambda\neq\mu$ implies $v\perp w$:
\begin{eqnarray*}
\lambda<v,w>
&=&<\lambda v,w>\\
&=&<Av,w>\\
&=&<v,A^*w>\\
&=&<v,\bar{\mu}w>\\
&=&\mu<v,w>
\end{eqnarray*}

In order to finish, it remains to prove that the eigenspaces of $A$ span the whole $\mathbb C^N$. This is something that we have already seen for the self-adjoint matrices, and for unitaries, and we will use here these results, in order to deal with the general normal case. As a first observation, given an arbitrary matrix $A$, the matrix $AA^*$ is self-adjoint:
$$(AA^*)^*=AA^*$$

Thus, we can diagonalize this matrix $AA^*$, as follows, with the passage matrix being a unitary, $V\in U_N$, and with the diagonal form being real, $E\in M_N(\mathbb R)$:
$$AA^*=VEV^*$$

Now observe that, for matrices of type $A=UDU^*$, which are those that we supposed to deal with, we have the following formulae:
$$V=U\quad,\quad 
E=D\bar{D}$$

In particular, the matrices $A$ and $AA^*$ have the same eigenspaces. So, this will be our idea, proving that the eigenspaces of $AA^*$ are eigenspaces of $A$. In order to do so, let us pick two eigenvectors $v,w$ of the matrix $AA^*$, corresponding to different eigenvalues, $\lambda\neq\mu$. The eigenvalue equations are then as follows:
$$AA^*v=\lambda v\quad,\quad 
AA^*w=\mu w$$

We have the following computation, using the normality condition $AA^*=A^*A$, and the fact that the eigenvalues of $AA^*$, and in particular $\mu$, are real:
\begin{eqnarray*}
\lambda<Av,w>
&=&<\lambda Av,w>\\
&=&<A\lambda v,w>\\
&=&<AAA^*v,w>\\
&=&<AA^*Av,w>\\
&=&<Av,AA^*w>\\
&=&<Av,\mu w>\\
&=&\mu<Av,w>
\end{eqnarray*}

We conclude that we have $<Av,w>=0$. But this reformulates as follows:
$$\lambda\neq\mu\implies A(E_\lambda)\perp E_\mu$$

Now since the eigenspaces of $AA^*$ are pairwise orthogonal, and span the whole $\mathbb C^N$, we deduce from this that these eigenspaces are invariant under $A$:
$$A(E_\lambda)\subset E_\lambda$$

But with this observation in hand, we can finish. Indeed, we can decompose the problem, and the matrix $A$ itself, following these eigenspaces of $AA^*$, which in practice amounts in saying that we can assume that we only have 1 eigenspace. But by rescaling, this is the same as assuming that we have $AA^*=1$, and so we are now into the unitary case, that we know how to solve, as explained in Theorem 12.19.
\end{proof}

As a first application of all this, we have the following result:

\index{absolute value}
\index{modulus of matrix}

\begin{theorem}
Given a matrix $A\in M_N(\mathbb C)$, we can construct a matrix $|A|$ as follows, by using the fact that $A^*A$ is diagonalizable, with positive eigenvalues:
$$|A|=\sqrt{A^*A}$$
This matrix $|A|$ is then positive, and its square is $|A|^2=A^*A$. In the case $N=1$, we obtain in this way the usual absolute value of the complex numbers.
\end{theorem}

\begin{proof}
Consider indeed the matrix $A^*A$, which is normal. According to Theorem 12.21, we can diagonalize this matrix as follows, with $U\in U_N$, and with $D$ diagonal:
$$A=UDU^*$$

From $A^*A\geq0$ we obtain $D\geq0$. But this means that the entries of $D$ are real, and positive. Thus we can extract the square root $\sqrt{D}$, and then set:
$$\sqrt{A^*A}=U\sqrt{D}U^*$$

Thus, we are basically done. Indeed, if we call this latter matrix $|A|$, we are led to the conclusions in the statement. Finally, the last assertion is clear from definitions.
\end{proof}

We can now formulate a useful polar decomposition result, as follows:

\index{polar decomposition}
\index{partial isometry}

\begin{theorem}
Any square matrix $A\in M_N(\mathbb C)$ decomposes as
$$A=U|A|$$
with $U$ being a partial isometry. When $A$ is invertible, $U$ is a unitary.
\end{theorem}

\begin{proof}
This is routine, and follows by comparing the actions of $A,|A|$ on the vectors $v\in\mathbb C^N$, and deducing from this the existence of a partial isometry $U$ as above.
\end{proof}

So long for advanced linear algebra. There are actually many other decomposition results for  the real matrices, quite often in relation with positivity, and the Jordan form too, which are all useful for questions in analysis, via derivatives and Hessians. 

\bigskip

Good luck of course in learning all this, when needed later in life, for your various math problems at that time. And always have in mind expert's advice, Cat 12.13.

\section*{12c. Harmonic functions} 

As another application of the second derivatives, we can now talk about harmonic functions in general, extending what we know from chapter 8. We first have:

\index{Laplace operator}
\index{harmonic function}

\begin{theorem}
We can talk about the Laplace operator in $N$ dimensions,
$$\Delta f=\sum_{i=1}^N\frac{d^2f}{dx_i^2}$$
and with this notion in hand, the following happen:
\begin{enumerate}
\item The Laplacian is the trace of the Hessian, $\Delta f=Tr(f'')$.

\item $\Delta f(x)$ expresses how far is $f(z)$, with $z\simeq x$, from $f(x)$.

\item The heat diffusion equation, namely $\dot{f}=\alpha\Delta f$, still holds.

\item The wave equation, namely $\ddot{f}=v^2\Delta f$, still holds too.

\item In particular, the light equation in $3D$ vacuum is $\ddot{f}=c^2\Delta f$.
\end{enumerate}
\end{theorem}

\begin{proof}
This is a straightforward remake of what we discussed in chapter 8, in two dimensions, the idea being that (1) is clear from definitions, (2) comes from the Taylor formula at order 2, then (3), which is supported by (2), comes from a lattice model, as in chapter 8, then (4) also comes from a lattice model, as in chapter 8, and finally (5) comes from (4), via some physics still to be discussed, and more on this in chapter 15.
\end{proof}

As before in chapter 8, save for some issues with the exact domain and image, we can think of the Laplace operator as being a linear operator of the space of functions $f:\mathbb R^N\to\mathbb R$, and linear algebra suggests to look at its eigenvalues and eigenvectors, and more specifically, to start with, at the $\lambda=0$ eigenspace. Which leads us to:

\begin{definition}
A function $f:X\to\mathbb R$ with $X\subset\mathbb R^N$ is called harmonic when
$$\Delta f=0$$
with this being called Laplace equation.
\end{definition}

Regarding now the harmonic functions, we have seen a lot of interesting mathematics regarding them, in chapter 8, in two dimensions. However, the case $N=2$ remains something quite special, due to the isomorphism $\mathbb R^2\simeq\mathbb C$, that we heavily used there. In general, in arbitrary $N$ dimensions, as a first interesting result, we have:

\index{Laplace equation}
\index{radial function}

\begin{theorem}
The fundamental radial solutions of $\Delta f=0$ are
$$f(x)=\begin{cases}
||x||^{2-N}&(N\neq 2)\\
\log||x||&(N=2)
\end{cases}$$
with the $\log$ at $N=2$ basically coming from $\log'=1/x$.
\end{theorem}

\begin{proof}
This is something that we actually stated in chapter 8, so time now for the proof.  Consider indeed a radial function, defined outside the origin $x=0$. This function can be written as follows, with $\varphi:(0,\infty)\to\mathbb C$ being a certain function:
$$f:\mathbb R^N-\{0\}\to\mathbb C\quad,\quad f(x)=\varphi(||x||)$$

Our first goal will be that of reformulating the Laplace equation $\Delta f=0$ in terms of the one-variable function $\varphi:(0,\infty)\to\mathbb C$. For this purpose, observe that we have:
\begin{eqnarray*}
\frac{d||x||}{dx_i}
&=&\frac{d\sqrt{\sum_{i=1}^Nx_i^2}}{dx_i}\\
&=&\frac{1}{2}\cdot\frac{1}{\sqrt{\sum_{i=1}^Nx_i^2}}\cdot\frac{d\left(\sum_{i=1}^Nx_i^2\right)}{dx_i}\\
&=&\frac{1}{2}\cdot\frac{1}{||x||}\cdot 2x_i\\
&=&\frac{x_i}{||x||}
\end{eqnarray*}

By using this formula, we have the following computation:
$$\frac{df}{dx_i}
=\frac{d\varphi(||x||)}{dx_i}
=\varphi'(||x||)\cdot\frac{d||x||}{dx_i}
=\varphi'(||x||)\cdot\frac{x_i}{||x||}$$

By differentiating one more time, we obtain the following formula:
\begin{eqnarray*}
\frac{d^2f}{dx_i^2}
&=&\frac{d}{dx_i}\left(\varphi'(||x||)\cdot\frac{x_i}{||x||}\right)\\
&=&\frac{d\varphi'(||x||)}{dx_i}\cdot\frac{x_i}{||x||}+\varphi'(||x||)\cdot\frac{d}{dx_i}\left(\frac{x_i}{||x||}\right)\\
&=&\left(\varphi''(||x||)\cdot\frac{x_i}{||x||}\right)\cdot\frac{x_i}{||x||}+\varphi'(||x||)\cdot\frac{||x||-x_i\cdot x_i/||x||}{||x||^2}\\
&=&\varphi''(||x||)\cdot\frac{x_i^2}{||x||^2}+\varphi'(||x||)\cdot\frac{||x||^2-x_i^2}{||x||^3}
\end{eqnarray*}

Now by summing over $i\in\{1,\ldots,N\}$, this gives the following formula:
\begin{eqnarray*}
\Delta f
&=&\sum_{i=1}^N\varphi''(||x||)\cdot\frac{x_i^2}{||x||^2}+\sum_{i=1}^N\varphi'(||x||)\cdot\frac{||x||^2-x_i^2}{||x||^3}\\
&=&\varphi''(||x||)\cdot\frac{||x||^2}{||x||^2}+\varphi'(||x||)\cdot\frac{(N-1)||x||^2}{||x||^3}\\
&=&\varphi''(||x||)+\varphi'(||x||)\cdot\frac{N-1}{||x||}
\end{eqnarray*}

Thus, with $r=||x||$, the Laplace equation $\Delta f=0$ can be reformulated as follows:
$$\varphi''(r)+\frac{(N-1)\varphi'(r)}{r}=0$$

Equivalently, the equation that we want to solve is as follows:
$$r\varphi''+(N-1)\varphi'=0$$

Now observe that we have the following formula:
\begin{eqnarray*}
(r^{N-1}\varphi')'
&=&(N-1)r^{N-2}\varphi'+r^{N-1}\varphi''\\
&=&r^{N-2}((N-1)\varphi'+r\varphi'')
\end{eqnarray*}

Thus, the equation to be solved can be simply written as follows:
$$(r^{N-1}\varphi')'=0$$

We conclude that $r^{N-1}\varphi'$ must be a constant $K$, and so, that we must have:
$$\varphi'=Kr^{1-N}$$

But the fundamental solutions of this latter equation are as follows:
$$\varphi(r)=\begin{cases}
r^{2-N}&(N\neq 2)\\
\log r&(N=2)
\end{cases}$$

Thus, we are led to the conclusion in the statement.
\end{proof}

Not bad for a start, all this. In analogy now with the one complex variable theory that we know from chapter 6 and chapter 8, we have the following result:

\index{maximum modulus}
\index{mean value formula}
\index{Liouville theorem}

\begin{theorem}
The harmonic functions in $N$ dimensions obey to the same general principles as the holomorphic functions, namely:
\begin{enumerate}
\item The plain mean value formula. 

\item The boundary mean value formula.

\item The maximum modulus principle.

\item The Liouville theorem.
\end{enumerate}
\end{theorem}

\begin{proof}
This follows exactly as in 2 dimensions, the idea being as follows:

\medskip

(1) Regarding the plain mean value formula, here the statement is that given an harmonic function $f:X\to\mathbb C$, and a ball $B$, the following happens:
$$f(x)=\int_Bf(y)dy$$

In order to prove this, we can assume that $B$ is centered at $0$, of radius $r>0$. If we denote by $\chi_r$ the characteristic function of this ball, nomalized as to integrate up to 1, in terms of the convolution operation from chapter 7, we want to prove that we have:
$$f=f*\chi_r$$

For doing so, let us pick a number $0<s<r$, and then a solution $w$ of the following equation on $B$, which can be constructed explicitly, say as a radial function:
$$\Delta w=\chi_r-\chi_s$$

By using the properties of the convolution operation $*$ from chapter 7, we have:
\begin{eqnarray*}
f*\chi_r-f*\chi_s
&=&f*(\chi_r-\chi_s)\\
&=&f*\Delta w\\
&=&\Delta f*w\\
&=&0
\end{eqnarray*}

Thus $f*\chi_r=f*\chi_s$, and by letting now $s\to0$, we get $f*\chi_r=f$, as desired.

\medskip

(2) Regarding the boundary mean value formula, here the statement is that given an harmonic function $f:X\to\mathbb C$, and a ball $B$, with boundary $\gamma$, the following happens:
$$f(x)=\int_\gamma f(y)dy$$

But this follows as a consequence of the plain mean value formula in (1), with our two mean value formulae, the one there and the one here, being in fact equivalent, by using annuli and radial integration for the proof of the equivalence, in the obvious way. 

\medskip

(3) Regarding the maximum modulus principle, the statement here is that any holomorphic function $f:X\to\mathbb C$ has the property that the maximum of $|f|$ over a domain is attained on its boundary. That is, given a domain $D$, with boundary $\gamma$, we have:
$$\exists x\in\gamma\quad,\quad |f(x)|=\max_{y\in D}|f(y)|$$

But this is something which follows again from the mean value formula in (1), first for the balls, and then in general, by using a standard division argument.

\medskip

(4) Finally, regarding the Liouville theorem, the statement here is that an entire, bounded harmonic function must be constant: 
$$f:\mathbb R^N\to\mathbb C\quad,\quad\Delta f=0\quad,\quad |f|\leq M\quad\implies\quad f={\rm constant}$$

As a slightly weaker statement, again called Liouville theorem, we have the fact that an entire harmonic function which vanishes at $\infty$ must vanish globally: 
$$f:\mathbb R^N\to\mathbb C\quad,\quad\Delta f=0\quad,\quad\lim_{x\to\infty}f(x)=0\quad\implies\quad f=0$$

But can view these as a consequence of the mean value formula in (1), because given two points $x\neq y$, we can view the values of $f$ at these points as averages over big balls centered at these points, say $B=B_x(R)$ and $C=B_y(R)$, with $R>>0$:
$$f(x)=\int_Bf(z)dz\quad,\quad f(y)=\int_Cf(z)dz$$

Indeed, the point is that when the radius goes to $\infty$, these averages tend to be equal, and so we have $f(x)\simeq f(y)$, which gives $f(x)=f(y)$ in the limit, as desired.
\end{proof}

Many other interesting things can be said about harmonic functions, and we will leave some learning here, both mathematics and physics, as an instructive exercise.

\section*{12d. Lagrange multipliers}

Let us discuss now some more specialized topics, in relation with optimization. First of all, thinking well, the functions that we have to minimize or maximize, in the real life, are often defined on a manifold, instead of being defined on the whole $\mathbb R^N$. Fortunately, the good old principle $f'(x)=0$ can be adapted to the manifold case, as follows:

\begin{principle}
In order for a function $f:X\to\mathbb R$ defined on a manifold $X$ to have a local extremum at $x\in X$, we must have, as usual 
$$f'(x)=0$$
but with this taking into account the fact that the equations defining the manifold count as well as ``zero'', and so must be incorporated into the formula $f'(x)=0$. 
\end{principle}

Obviously, we are punching here about our weight, because our discussion about manifolds from chapter 11 was quite introductory, and we have no tools in our bag for proving such things, or even for properly formulating them. However, we can certainly talk about all this, a bit like physicists do. So, our principle will be the one above, and in practice, the idea is that we must have a formula as follows, with $g_i$ being the constraint functions for our manifold $X$, and with $\lambda_i\in\mathbb R$ being certain scalars, called Lagrange multipliers:
$$f'(x)=\sum_i\lambda_ig_i'(x)$$

As a basic illustration for this, our claim is that, by using a suitable manifold, and a suitable function, and Lagrange multipliers, we can prove in this way the H\"older inequality, that we know well of course, but without any computation. Let us start with:

\begin{proposition}
For any exponent $p>1$, the following set
$$S_p=\left\{x\in\mathbb R^N\Big|\sum_i|x_i|^p=1\right\}$$
is a submanifold of $\mathbb R^N$.
\end{proposition}

\begin{proof}
We know from chapter 11 that the unit sphere in $\mathbb R^N$ is a manifold. In our terms, this solves our problem at $p=2$, because this unit sphere is:
$$S_2=\left\{x\in\mathbb R^N\Big|\sum_ix_i^2=1\right\}$$

Now observe that we have a bijection $S_p\simeq S_2$, at least on the part where all the coordinates are positive, $x_i>0$, given by the following function:
$$x_i\to x_i^{2/p}$$

Thus we obtain that $S_p$ is indeed a manifold, as claimed.
\end{proof}

We already know that the manifold $S_p$ constructed above is the unit sphere, in the case $p=2$. In order to have a better geometric picture of what is going on, in general, observe that $S_p$ can be constructed as well at $p=1$, as follows:
$$S_1=\left\{x\in\mathbb R^N\Big|\sum_i|x_i|=1\right\}$$

However, this is no longer a manifold, as we can see for instance at $N=2$, where we obtain a square. Now observe that we can talk as well about $p=\infty$, as follows:
$$S_\infty=\left\{x\in\mathbb R^N\Big|\sup_i|x_i|=1\right\}$$

This letter set is no longer a manifold either, as we can see for instance at $N=2$, where we obtain again a square, containing the previous square, the one at $p=1$.

\bigskip

With these limiting constructions in hand, we can have now a better geometric picture of what is going on, in the general context of Proposition 12.29. Indeed, let us draw, at $N=2$ for simplifying, our sets $S_p$ at the values $p=1,2,\infty$ of the exponent:
$$
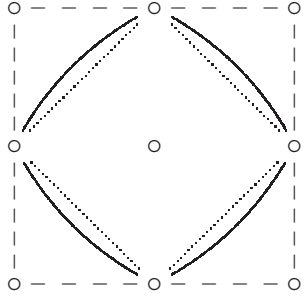
\xymatrix@R=40pt@C=44pt{
\circ\ar@{--}[r]\ar@{--}[d]&\circ\ar@{--}[r]\ar@{.}[dr]\ar@{.}[dl]\ar@{-}@/^/[dr]&\circ\ar@{--}[d]\\
\circ\ar@{--}[d]\ar@{.}[dr]\ar@{-}@/^/[ur]&\circ&\circ\ar@{--}[d]\ar@{.}[dl]\ar@{-}@/^/[dl]\\
\circ\ar@{--}[r]&\circ\ar@{--}[r]\ar@{-}@/^/[ul]&\circ
}$$

We can see that what we have is a small square, at $p=1$, becoming smooth and inflating towards the circle, in the parameter range $p\in(1,2]$, and then further inflating, in the parameter range $p\in[2,\infty)$, towards the big square appearing at $p=\infty$.

\bigskip

With these preliminaries in hand, we can formulate our result, as follows:

\begin{theorem}
The local extrema over $S_p$ of the function
$$f(x)=\sum_ix_iy_i$$ 
can be computed by using Lagrange multipliers, and this gives
$$\left|\sum_ix_iy_i\right|\leq\left(\sum_i|x_i|^p\right)^{1/p}\left(\sum_i|y_i|^q\right)^{1/q}$$
with $1/p+1/q=1$, that is, the H\"older inequality, with a purely geometric proof.
\end{theorem}

\begin{proof}
We can restrict the attention to the case where all the coordinates are positive, $x_i>0$ and  $y_i>0$. The derivative of the function in the statement is:
$$f'(x)=(y_1,\ldots,y_N)$$

On the other hand, we know that the manifold $S_p$ appears by definition as the set of zeroes of the function $\varphi(x)=\sum_ix_i^p-1$, having derivative as follows:
$$\varphi'(x)=p(x_1^{p-1},\ldots,x_N^{p-1})$$

Thus, by using Lagrange multipliers, the critical points of $f$ must satisfy:
$$(y_1,\ldots,y_N)\sim(x_1^{p-1},\ldots,x_N^{p-1})$$

In other words, the critical points must satisfy $x_i=\lambda y_i^{1/(p-1)}$, for some $\lambda>0$, and by using now $\sum_ix_i^p=1$ we can compute the precise value of $\lambda$, and we get:
$$\lambda=\left(\sum_iy_i^{p/(p-1)}\right)^{-1/p}$$ 

Now let us see what this means. Since the critical point is unique, this must be a maximum of our function, and we conclude that for any $x\in S_p$, we have:
$$\sum_ix_iy_i
\leq\sum_i\lambda y_i^{1/(p-1)}\cdot y_i
=\left(\sum_iy_i^{p/(p-1)}\right)^{1-1/p}
=\left(\sum_iy_i^q\right)^{1/q}$$

Thus we have H\"older, and the general case follows from this, by rescaling.
\end{proof}

As a second illustration for the method of Lagrange multipliers, this time in relation with certain questions from linear algebra, let us go back to the Hadamard matrices, that we met in chapter 7, when discussing discrete aspects of Fourier analysis.

\bigskip

In the real case, the basic theory of these matrices is as follows:

\begin{theorem}
The real Hadamard matrices, $H\in M_N(-1,1)$ having pairwise orthogonal rows, have the following properties:
\begin{enumerate}
\item The set of Hadamard matrices is $X_N=M_N(-1,1)\cap\sqrt{N}O_N$.

\item In order to have $X_N\neq\emptyset$, the matrix size must be $N\in\{2\}\cup 4\mathbb N$.

\item For $H\in M_N(-1,1)$ we have $|\det H|\leq N^{N/2}$, with equality when $H$ is Hadamard.

\item For $U\in O_N$ we have $||U||_1\leq N\sqrt{N}$, with equality when $H=\sqrt{N}U$ is Hadamard.
\end{enumerate}
\end{theorem}

\begin{proof}
Many things going on here, the idea being as follows:

\medskip

(1) This is just a reformulation of the Hadamard matrix condition.

\medskip

(2) This follows by playing with the first 3 rows, exercise for you.

\medskip

(3) This follows from our definition of the determinant, as a signed volume.

\medskip

(4) This follows from $||U||_2=\sqrt{N}$ and Cauchy-Schwarz, easy exercise for you. 
\end{proof}

All the above is quite interesting, and (1,2) raise the question of finding the correct generalizations of the Hadamard matrices, at $N\notin\{2\}\cup 4\mathbb N$. But the answer here comes from (3,4), which suggest looking either at the maximizers of $|\det|$ on $M_N(-1,1)$, or of $||.||_1$ on $\sqrt{N}O_N$. By following this latter way, we are led to the following question:

\begin{question}
What are the critical points of the $1$-norm on $O_N$?
\end{question}

And, good news, we can solve this latter question by using the theory of Lagrange multipliers developed in the above, the result here being as follows:

\index{Lagrange multipliers}

\begin{theorem}
An orthogonal matrix with nonzero entries is a critical point of
$$||U||_1=\sum_{ij}|U_{ij}|$$
precisely when $SU^t$ is symmetric, where $S_{ij}={\rm sgn}(U_{ij})$.
\end{theorem}

\begin{proof}
We regard $O_N$ as a real algebraic manifold, with coordinates $U_{ij}$. This manifold consists by definition of the zeroes of the following polynomials: 
$$A_{ij}=\sum_kU_{ik}U_{jk}-\delta_{ij}$$

Thus $U\in O_N$ is a critical point of $F(U)=||U||_1$ when the following is satisfied: 
$$dF\in span(dA_{ij})$$

Regarding the space $span(dA_{ij})$, this consists of the following quantities:
\begin{eqnarray*}
\sum_{ij}M_{ij}dA_{ij}
&=&\sum_{ijk}M_{ij}(U_{ik}dU_{jk}+U_{jk}dU_{ik})\\
&=&\sum_{ij}(M^tU)_{ij}dU_{ij}+\sum_{ij}(MU)_{ij}dU_{ij}
\end{eqnarray*}

In order to compute $dF$, observe first that, with $S_{ij}=sgn(U_{ij})$, we have:
$$d|U_{ij}|
=d\sqrt{U_{ij}^2}
=\frac{U_{ij}dU_{ij}}{|U_{ij}|}
=S_{ij}dU_{ij}$$

Thus $dF=\sum_{ij}S_{ij}dU_{ij}$, and so $U\in O_N$ is a critical point of $F$ precisely when there exists a matrix $M\in M_N(\mathbb R)$ such that the following two conditions are satisfied:
$$S=M^tU\quad,\quad 
S=MU$$

Now observe that these two equations can be written as follows:
$$M^t=SU^t\quad,\quad 
M=SU^t$$

Thus, the matrix $SU^t$ must be symmetric, as claimed.
\end{proof}

\section*{12e. Exercises}

This chapter was tough analysis as it comes, and as exercises, we have:

\begin{exercise}
Clarify the Taylor formula at order $3$, and its applications.
\end{exercise}

\begin{exercise}
Learn the Jordan form, as a complement to the linear algebra here.
\end{exercise}

\begin{exercise}
Harmonic functions in $3$ dimensions, and their applications.
\end{exercise}

\begin{exercise}
Clarify all the details, in relation with Lagrange multipliers.
\end{exercise}

As bonus exercise, learn the gradient method. You won't escape from that.

\part{Integration theory}

\ \vskip50mm

\begin{center}
{\em I'm the left eye

You're the right

Would it not be madness to fight

We come one}
\end{center}

\chapter{Multiple integrals}

\section*{13a. Multiple integrals}

Welcome to advanced calculus. We have kept the best for the end, and in this whole last part of the present book we will learn how to integrate functions of several variables. The general, fascinating question that we will be interested in is as follows:

\begin{question}
Given a function $f:\mathbb R^N\to\mathbb R$, how to compute
$$\int_{\mathbb R^N}f(x)\,dx_1\ldots dx_N$$
or at least, what are the rules satisfied by such integrals?
\end{question}

Here we adopt, somehow by definition, the convention that the above integral is constructed a bit like the one-variable integrals, via Riemann sums as in chapter 4, by using this time divisions of the space $\mathbb R^N$ into small $N$-dimensional cubes. There are of course some theory and details to be worked out here, not exactly trivial, but since we have not really done this in chapter 4, nor we will do this here. Let's focus on computations.

\bigskip

At the first glance, solving Question 13.1 looks like an easy task, because we can iterate one-variable integrations, which are something that we know well, from chapter 4. For instance the integral of a function $f:\mathbb R^2\to\mathbb R$ can be computed by using:
$$\int_{\mathbb R^2}f(z)dz=\int_\mathbb R\int_\mathbb Rf(x,y)dxdy$$

This being said, when doing so, we are faced right away with a dilemma. Indeed, we can use as well the following rival method, normally yielding the same answer:
$$\int_{\mathbb R^2}f(z)dz=\int_\mathbb R\int_\mathbb Rf(x,y)dydx$$

So, which method is the best? Depends on $f$, of course. However, things do not stop here, because in certain situations it is better to use polar coordinates, as follows:
$$\int_{\mathbb R^2}f(z)dz=\int_0^{2\pi}\int_0^\infty f(r\cos t,r\sin t)Jdrdt$$

Here the factor on the left is $J=dxdy/drdt$, which remains to be computed. And for the picture to be complete, we have as well the following fourth formula:
$$\int_{\mathbb R^2}f(z)dz=\int_0^\infty\int_0^{2\pi} f(r\cos t,r\sin t)Jdtdr$$

In short, you got my point, I hope, things are quite complicated in several variables, with the complications starting already in 2 variables. And actually, if you think a bit about 3 variables,  it is quite clear that the above complications can become true nightmares, leading to long nights spent in computing integrals, there in 3 variables.

\bigskip

So, what to do? Work and patience, of course, and here is our plan:

\bigskip

(1) In this chapter we will get used to the multiple integrals, with some general rules for their computation, including a rule for computing the above factor $J=dxdy/drdt$. Then, we will enjoy all this by computing some integrals over the spheres in $\mathbb R^N$. 

\bigskip

(2) In chapter 14 we will review probability theory, and further develop it, notably with the theory of normal variables. And finally, in chapters 15-16 we will get back to physics, with some sharp results in 3 dimensions, and then in infinite dimensions.

\bigskip

This sounds good, but as a matter of doublechecking what we are doing, make sure that it is wise indeed, let us ask the cat about what he thinks. And cat answers: 

\begin{cat}
Not quite sure about your formula 
$$\int_\mathbb R\int_\mathbb Rf(x,y)dxdy=\int_\mathbb R\int_\mathbb Rf(x,y)dydx$$
and I doubt too that you can properly compute $J=dxdy/drdt$. Read Rudin.
\end{cat}

Oh dear. What can I say. Sure I read Rudin, as a Gen X mathematician, but we are now well into the 21st century, and shall I go ahead with heavy measure theory, for having all this properly developed? Not quite sure, I'd rather stick to my plan.

\bigskip

So, ignoring what cat says, but get however a copy of Rudin's red book \cite{ru2}, and don't forget about Mao's too, and getting back to our plan, as a first goal, we would like to compute the factor $J=dxdy/drdt$. Let us start with something that we know, in 1D:

\index{change of variable}

\begin{proposition}
We have the change of variable formula
$$\int_a^bf(x)dx=\int_c^df(\varphi(t))\varphi'(t)dt$$
where $c=\varphi^{-1}(a)$ and $d=\varphi^{-1}(b)$.
\end{proposition}

\begin{proof}
This is something that we know well from chapter 4, obtained by integrating the standard differentiation rule $(F\varphi)'(t)=F'(\varphi(t))\varphi'(t)$, with $f=F'$.
\end{proof}

In several variables now, we have the following extension of the above result:

\index{change of variable}
\index{Jacobian}

\begin{theorem}
Given a transformation $\varphi=(\varphi_1,\ldots,\varphi_N)$, we have
$$\int_Ef(x)dx=\int_{\varphi^{-1}(E)}f(\varphi(t))|J_\varphi(t)|dt$$
with the $J_\varphi$ quantity, called Jacobian, being given by
$$J_\varphi(t)=\det\left[\left(\frac{d\varphi_i}{dx_j}(x)\right)_{ij}\right]$$ 
and with this generalizing the formula from Proposition 13.3.
\end{theorem}

\begin{proof}
This is something quite tricky, the idea being as follows:

\medskip

(1) Observe first that this generalizes Proposition 13.3, with the absolute value on the derivative appearing as to compensate for the lack of explicit bounds for the integral.

\medskip 

(2) As a second observation, we can assume if we want, by linearity, that we are dealing with the constant function $f=1$. And with $D={\varphi^{-1}(E)}$, our formula reads:
$$vol(\varphi(D))=\int_D|J_\varphi(t)|dt$$

Now since this latter formula is additive with respect to $D$, it is enough to prove it for small cubes $D$. And here, as a first remark, our formula is clear for the linear maps $\varphi$, by using the definition of the determinant of real matrices, as a signed volume.

\medskip

(3) However, the extension of this to the case of non-linear maps $\varphi$ is something which looks non-trivial, so we will not follow this path, in what follows.

\medskip

(4) In order to prove the theorem, as stated, let us rather focus on the transformations used $\varphi$, instead of the functions to be integrated $f$. Our first claim is that the validity of the theorem is stable under taking compositions of such transformations $\varphi$.

\medskip

(5) In order to prove this claim, consider a composition, as follows:
$$\varphi:E\to F\quad,\quad 
\psi:D\to E\quad,\quad 
\varphi\circ\psi:D\to F$$

Assuming that the theorem holds for $\varphi,\psi$, we have the following computation:
\begin{eqnarray*}
\int_Ff(x)dx
&=&\int_Ef(\varphi(s))|J_\varphi(s)|ds\\
&=&\int_Df(\varphi\circ\psi(t))|J_\varphi(\psi(t))|\cdot|J_\psi(t)|dt\\
&=&\int_Df(\varphi\circ\psi(t))|J_{\varphi\circ\psi}(t)|dt
\end{eqnarray*}

Thus, our theorem holds as well for $\varphi\circ\psi$, and we have proved our claim.

\medskip

(6) Next, as a key ingredient, let us examine the case where we are in $N=2$ dimensions, and our transformation $\varphi$ has one of the following special forms:
$$\varphi(x,y)=(\psi(x,y),y)\quad,\quad\varphi(x,y)=(x,\psi(x,y))$$

By symmetry, it is enough to deal with the first case. Here the Jacobian is $d\psi/dx$, and by replacing if needed $\psi\to-\psi$, we can assume that this Jacobian is positive, $d\psi/dx>0$. Now by assuming as before that $D=\varphi^{-1}(E)$ is a rectangle, $D=[a,b]\times[c,d]$, we can prove our formula by using the change of variables in 1 dimension, as follows:
\begin{eqnarray*}
\int_Ef(s)ds
&=&\int_{\varphi(D)}f(x,y)dxdy\\
&=&\int_c^d\int_{\psi(a,y)}^{\psi(b,y)}f(x,y)dxdy\\
&=&\int_c^d\int_a^bf(\psi(x,y),y)\frac{d\psi}{dx}\,dxdy\\
&=&\int_Df(\varphi(t))J_\varphi(t)dt
\end{eqnarray*}

(7) But with this, we can now prove the theorem, in $N=2$ dimensions. Indeed, given a transformation $\varphi=(\varphi_1,\varphi_2)$, consider the following two transformations:
$$\phi(x,y)=(\varphi_1(x,y),y)\quad,\quad \psi(x,y)=(x,\varphi_2\circ\phi^{-1}(x,y))$$

We have then $\varphi=\psi\circ\phi$, and by using (6) for $\psi,\phi$, which are of the special form there, and then (3) for composing, we conclude that the theorem holds for $\varphi$, as desired.

\medskip

(8) Thus, theorem proved in $N=2$ dimensions, at least in the generic situation, and we will leave the remaining details as an exercise. And the extension of the above proof to arbitrary $N$ dimensions is straightforward, that we will leave as an exercise too.
\end{proof}

Quite nice all this, we have some theory going on, and are on the good way for computing $J=dxdy/drdt$, and this regardless of what cat says. This being said, getting back to Cat 13.2, there are some deep truths there. Consider for instance:
\begin{eqnarray*}
\int_0^1\int_0^1\frac{y^2-x^2}{(x^2+y^2)^2}\,dxdy
&=&\int_0^1\left[\frac{x}{x^2+y^2}\right]_0^1dy\\
&=&\int_0^1\frac{1}{1+y^2}\,dy\\
&=&\frac{\pi}{4}
\end{eqnarray*}

Easy calculus you would say, but the point is that, when changing $dxdy\to dydx$, the integral will change its sign, and so become $-\pi/4$, different from the original integral.

\bigskip

Quite bizarre all this, hope you agree with me. So, what to do? In answer, there are 3 possible ways of dealing with the potential traps of multivariable integration:

\medskip

(1) Avoid them by reading Rudin \cite{ru2}, as cat says, where measure theory, and integration of the functions $f:\mathbb R^N\to\mathbb R$ is developed, in a fully rigorous way.

\medskip

(2) Know at least the main findings of measure theory. With these being, basically, ``things fine as long as your computations are bounded by an integrable function''.

\medskip

(3) Know nothing advanced, but be extremely careful with your computations, with checks and doublechecks, as per our Advice 1.3, from the beginning of this book. 

\medskip

And we will stop here with the philosophy. In what follows I will use a combination of (2) and (3), like everyone or almost in the business does, although we are surely familiar with (1) too, with all the formulae below being guaranteed to be correct.

\section*{13b. Gauss and Fresnel}

Time now do some exciting computations, with the technology that we have. In what regards the applications of Theorem 13.4, these often come via:

\index{polar coordinates}

\begin{proposition}
We have polar coordinates in $2$ dimensions,
$$\begin{cases}
x\!\!\!&=\ r\cos t\\
y\!\!\!&=\ r\sin t
\end{cases}$$
the corresponding Jacobian being $J=r$.
\end{proposition}

\begin{proof}
This is elementary, the Jacobian being:
\begin{eqnarray*}
J
&=&\begin{vmatrix}
\frac{d(r\cos t)}{dr}&&\frac{d(r\cos t)}{dt}\\
\\
\frac{d(r\sin t)}{dr}&&\frac{d(r\sin t)}{dt}
\end{vmatrix}\\
&=&\begin{vmatrix}
\cos t&-r\sin t\\
\sin t&r\cos t
\end{vmatrix}\\
&=&r\cos^2t+r\sin^2t\\
&=&r
\end{eqnarray*}

Thus, we have indeed the formula in the statement.
\end{proof}

We can now compute the Gauss integral, which is the best calculus formula ever:

\begin{theorem}
We have the following formula,
$$\int_\mathbb Re^{-x^2}dx=\sqrt{\pi}$$
called Gauss integral formula.
\end{theorem}

\begin{proof}
This is something truly magic, the idea being as follows:

\medskip

(1) To start with, we can certainly integrate $e^{-x^2}$ by using the formula of the exponential series, and the primitive which is worth 0 at $x=0$ is given by:
$$\int e^{-x^2}=\sum_{k=0}^\infty(-1)^k\frac{x^{2k+1}}{(2k+1)k!}$$

However, this series is not computable, in terms of the known, familiar series.

\medskip

(2) Thus, no primitive, but we can still ask for the computation of $\int_\mathbb Re^{-x^2}dx$, who knows. And here, another surprise awaits us, this is simply undoable, with bare hands, I mean all the formulae and tricks that we learned in chapter 4 fail, for this integral.

\medskip

(3) Which seems to send our problem to the trash can. However, and here comes the magic, the Gauss integral can be computed by using two dimensions, as follows:
\begin{eqnarray*}
\left(\int_\mathbb Re^{-x^2}dx\right)^2
&=&\int_\mathbb R\int_\mathbb Re^{-x^2-y^2}dxdy\\
&=&\int_0^{2\pi}\int_0^\infty e^{-r^2}rdrdt\\
&=&2\pi\int_0^\infty\left(-\frac{e^{-r^2}}{2}\right)'dr\\
&=&2\pi\left[0-\left(-\frac{1}{2}\right)\right]\\
&=&\pi
\end{eqnarray*}

(4) Amazing all this, isn't it. There are of course some other known proofs of the Gauss formula, but all using two dimensions, appearing as variations of the above.
\end{proof}

As a first application of the Gauss formula, we can fix some bugs with the physics from chapter 8, in relation with the normalization of the heat kernel, and we have:

\begin{theorem}
The heat diffusion equation, $\dot{f}=\alpha f''$ with $\alpha>0$, with initial condition $f(x,0)=g(x)$, has as solution the function
$$f(x,t)=\int_{\mathbb R^N}K_t(x-y)g(y)dy$$
where the function $K_t:\mathbb R^N\to\mathbb R$, called heat kernel, given by
$$K_t(x)=\frac{1}{\sqrt{(4\pi\alpha t)^N}}\,e^{-||x||^2/4\alpha t}$$
is the standard solution, coming from the initial data $g=\delta_0$, Dirac mass at $0$.
\end{theorem}

\begin{proof}
We already talked about the heat diffusion equation, on several occasions, and notably in chapter 8, with a quite detailed discussion of what happens in 1D. That discussion generalizes in a straightforward way in arbitrary $N$ dimensions, and we will leave the details here as an exercise. So, let us focus now on the thing which was not done in chapter 8, namely normalization of the heat kernel. In one dimension, we have:
\begin{eqnarray*}
\int_\mathbb RK_t(x)dx
&=&\frac{1}{\sqrt{4\pi\alpha t}}\int_\mathbb Re^{-x^2/4\alpha t}dx\\
&=&\frac{1}{\sqrt{4\pi\alpha t}}\int_\mathbb Re^{-y^2}\sqrt{4\alpha t}\,dy\\
&=&\frac{1}{\sqrt{4\pi\alpha t}}\cdot \sqrt{4\alpha t}\cdot\sqrt{\pi}\\
&=&1
\end{eqnarray*}

As for the mass 1 property in arbitrary $N$ dimensions, this comes from this, by obvious multiplicativity. Thus, we are led to the conclusions in the statement.
\end{proof}

And more on the applications of the Gauss formula in chapter 15, when talking about probability and normal variables, which are also called Gaussian variables.

\bigskip

We would like to end this section with a quick discussion of the Fresnel integrals. These integrals, coming as well from physics, and more specifically from the work on Fresnel on optics, are quite similar to the Gauss integral, the result being as follows:

\begin{theorem}
We have the following formulae,
$$\int_0^\infty\sin(t^2)dt=\int_0^\infty\cos(t^2)dt=\sqrt{\frac{\pi}{8}}$$
with these being called Fresnel integrals.
\end{theorem}

\begin{proof}
This is something quite tricky, a bit as before for the Gauss integral, that we discussed in Theorem 13.6, the idea with this being as follows:

\medskip

(1) As before with Gauss, we can certainly integrate $\sin(t^2)$ and $\cos(t^2)$ by using the series of $\sin$ and $\cos$, and the primitives which are worth 0 at $x=0$ are given by:
$$\int\sin(t^2)=\sum_{k=0}^\infty(-1)^k\frac{x^{4k+3}}{(4k+3)(2k+1)!}$$
$$\int\cos(t^2)=\sum_{k=0}^\infty(-1)^k\frac{x^{4k+1}}{(4k+1)(2k)!}$$

However, these series are not computable, in terms of the known, familiar series.

\medskip

(2) Next, and again as before with Gauss, we can try however to compute the integrals in the statement, and here the basic tools of one-variable calculus fail. But then, again as with Gauss, the two dimensions come to the rescue. In order to discuss this, observe first that, due to $e^{it^2}=\cos(t^2)+i\sin(t^2)$, the Fresnel formulae are equivalent to:
$$\int_0^\infty e^{it^2}dt=\sqrt{\frac{\pi}{2}}\cdot\frac{1+i}{2}$$

(3) So, let us prove this. For this purpose, consider the following function:
$$f(t)=\int_0^\infty\frac{e^{(i-u^2)t^2}}{i-u^2}\,du$$

The derivative of this function is then given by the following formula:
\begin{eqnarray*}
f'(t)
&=&2t\int_0^\infty e^{(i-u^2)t^2}du\\
&=&2te^{it^2}\int_0^\infty e^{-u^2t^2}du\\
&=&2e^{it^2}\int_0^\infty e^{-v^2}dv\\
&=&2e^{it^2}\times\frac{\sqrt{\pi}}{2}\\
&=&\sqrt{\pi}e^{it^2}
\end{eqnarray*}

(4) Now let us integrate this derivative, from 0 to $\infty$. We obtain in this way:
\begin{eqnarray*}
\sqrt{\pi}\int_0^\infty e^{it^2}dt
&=&\int_0^\infty f'(t)dt\\
&=&f(\infty)-f(0)\\
&=&0-\int_0^\infty\frac{1}{i-u^2}\,du\\
&=&\int_0^\infty\frac{1}{u^2-i}\,du
\end{eqnarray*}

Summarizing, we have obtained the following formula, for the Fresnel integral:
$$\int_0^\infty e^{it^2}dt=\frac{1}{\sqrt{\pi}}\int_0^\infty\frac{1}{u^2-i}\,du$$

(5) In order to compute the latter integral, set $w=e^{\pi i/4}$. Then $w^2=i$, and so:
\begin{eqnarray*}
\int_0^\infty\frac{1}{u^2-i}\,du
&=&\int_0^\infty\frac{1}{u^2-w^2}\,du\\
&=&\int_0^\infty\frac{1}{2w}\left(\frac{1}{u-w}-\frac{1}{u+w}\right)du\\
&=&\frac{1}{2w}\left[\log\left(\frac{u-w}{u+w}\right)\right]_0^\infty\\
&=&-\frac{1}{2w}\log(-1)
\end{eqnarray*}

Which brings us to the question, what is $\log(-1)$, in relation with this. In answer, and trust me here, $\log(-1)=-\pi i$, and we can finish our computation as follows:
$$\int_0^\infty\frac{1}{u^2-i}\,du
=-\frac{1}{2w}\times(-\pi i)
=\frac{\pi i}{2w}=\frac{\pi w}{2}
=\frac{\pi}{2}\cdot\frac{1+i}{\sqrt{2}}$$

(6) Alternatively, for avoiding this complex number mess, observe first that:
$$\int_0^\infty\frac{1}{u^2-i}\,du
=\int_0^\infty\frac{u^2}{u^4+1}\,du+i\int_0^\infty\frac{1}{u^4+1}\,du$$

Next, the two real integrals are equal, because with $u\to u^{-1}$ we obtain:
$$\int_0^\infty\frac{u^2}{u^4+1}\,du
=\int_0^\infty\frac{u^{-2}}{u^{-4}+1}\,u^{-2}du
=\int_0^\infty\frac{1}{u^4+1}\,du$$

Also, we can compute the sum of these integrals by using $t=u-u^{-1}$, as follows:
$$\int_0^\infty\frac{u^2+1}{u^4+1}\,du
=\int_0^\infty\frac{1+u^{-2}}{u^2+u^{-2}}\,du
=\int_0^\infty\frac{dt}{t^2+2}=\frac{\pi}{\sqrt{2}}$$

Thus, we are led to the following conclusion, exactly as we found in (5):
$$\int_0^\infty\frac{1}{u^2-i}\,du
=\frac{\pi}{2\sqrt{2}}+i\,\frac{\pi}{2\sqrt{2}}
=\frac{\pi}{2}\cdot\frac{1+i}{\sqrt{2}}$$

(7) Summarizing, computation done, one way or another, and this gives:
$$\int_0^\infty e^{it^2}dt=\frac{1}{\sqrt{\pi}}\times\frac{\pi}{2}\cdot\frac{1+i}{\sqrt{2}}
=\sqrt{\frac{\pi}{2}}\cdot\frac{1+i}{2}$$

But this is exactly what we wanted, and this ends the proof of our result.
\end{proof}

\section*{13c. Wallis and Stirling}

Moving now to 3 dimensions and higher, we already talked about spherical coordinates in chapter 11, when doing geometry. At the analytic level, we first have, in 3D:

\index{spherical coordinates}

\begin{proposition}
We have spherical coordinates in $3$ dimensions,
$$\begin{cases}
x\!\!\!&=\ r\cos s\\
y\!\!\!&=\ r\sin s\cos t\\
z\!\!\!&=\ r\sin s\sin t
\end{cases}$$
the corresponding Jacobian being $J(r,s,t)=r^2\sin s$.
\end{proposition}

\begin{proof}
The Jacobian is indeed given by the following formula:
\begin{eqnarray*}
&&J(r,s,t)\\
&=&\begin{vmatrix}
\cos s&-r\sin s&0\\
\sin s\cos t&r\cos s\cos t&-r\sin s\sin t\\
\sin s\sin t&r\cos s\sin t&r\sin s\cos t
\end{vmatrix}\\
&=&r^2\sin s\sin t
\begin{vmatrix}\cos s&-r\sin s\\ \sin s\sin t&r\cos s\sin t\end{vmatrix}
+r\sin s\cos t\begin{vmatrix}\cos s&-r\sin s\\ \sin s\cos t&r\cos s\cos t\end{vmatrix}\\
&=&r\sin s\sin^2 t
\begin{vmatrix}\cos s&-r\sin s\\ \sin s&r\cos s\end{vmatrix}
+r\sin s\cos^2 t\begin{vmatrix}\cos s&-r\sin s\\ \sin s&r\cos s\end{vmatrix}\\
&=&r\sin s(\sin^2t+\cos^2t)\begin{vmatrix}\cos s&-r\sin s\\ \sin s&r\cos s\end{vmatrix}\\
&=&r\sin s\times 1\times r\\
&=&r^2\sin s
\end{eqnarray*}

Thus, we have indeed the formula in the statement.
\end{proof}

In general, the result, which generalizes those at $N=2,3$, is as follows:

\index{spherical coordinates}

\begin{theorem}
We have spherical coordinates in $N$ dimensions,
$$\begin{cases}
x_1\!\!\!&=\ r\cos t_1\\
x_2\!\!\!&=\ r\sin t_1\cos t_2\\
\vdots\\
x_{N-1}\!\!\!&=\ r\sin t_1\sin t_2\ldots\sin t_{N-2}\cos t_{N-1}\\
x_N\!\!\!&=\ r\sin t_1\sin t_2\ldots\sin t_{N-2}\sin t_{N-1}
\end{cases}$$
the corresponding Jacobian being given by the following formula,
$$J(r,t)=r^{N-1}\sin^{N-2}t_1\sin^{N-3}t_2\,\ldots\,\sin^2t_{N-3}\sin t_{N-2}$$
and with this generalizing the known formulae at $N=2,3$.
\end{theorem}

\begin{proof}
As before, the fact that we have spherical coordinates is clear. Regarding now the Jacobian, also as before, by developing over the last column, we have:
\begin{eqnarray*}
J_N
&=&r\sin t_1\ldots\sin t_{N-2}\sin t_{N-1}\times \sin t_{N-1}J_{N-1}\\
&+&r\sin t_1\ldots \sin t_{N-2}\cos t_{N-1}\times\cos t_{N-1}J_{N-1}\\
&=&r\sin t_1\ldots\sin t_{N-2}(\sin^2 t_{N-1}+\cos^2 t_{N-1})J_{N-1}\\
&=&r\sin t_1\ldots\sin t_{N-2}J_{N-1}
\end{eqnarray*}

Thus, we obtain the formula in the statement, by recurrence.
\end{proof}

As a comment here, the above convention for spherical coordinates is one among many, designed to best work in arbitrary $N$ dimensions. Also, in what regards the precise range of the angles $t_1,\ldots,t_{N-1}$, we will leave this to you, as an instructive exercise.

\bigskip

As an application, let us compute the volumes of spheres. For this purpose, we must understand how the products of coordinates integrate over spheres. Let us start with the case $N=2$. Here the sphere is the unit circle $\mathbb T$, and with $z=e^{it}$ the coordinates are $\cos t,\sin t$. We can first integrate arbitrary powers of these coordinates, as follows:

\index{double factorials}
\index{trigonometric integral}

\begin{theorem}[Wallis]
We have the following formulae,
$$\int_0^{\pi/2}\cos^pt\,dt=\int_0^{\pi/2}\sin^pt\,dt=\left(\frac{\pi}{2}\right)^{\varepsilon(p)}\frac{p!!}{(p+1)!!}$$
where $\varepsilon(p)=1$ if $p$ is even, and $\varepsilon(p)=0$ if $p$ is odd, and where
$$m!!=(m-1)(m-3)(m-5)\ldots$$
with the product ending at $2$ if $m$ is odd, and ending at $1$ if $m$ is even.
\end{theorem}

\begin{proof}
Let us first compute the integral on the left in the statement:
$$I_p=\int_0^{\pi/2}\cos^pt\,dt$$

We can do this by partial integration. We have the following formula:
\begin{eqnarray*}
(\cos^pt\sin t)'
&=&p\cos^{p-1}t(-\sin t)\sin t+\cos^pt\cos t\\
&=&p\cos^{p+1}t-p\cos^{p-1}t+\cos^{p+1}t\\
&=&(p+1)\cos^{p+1}t-p\cos^{p-1}t
\end{eqnarray*}

By integrating between $0$ and $\pi/2$, we obtain the following formula:
$$(p+1)I_{p+1}=pI_{p-1}$$

Thus we can compute $I_p$ by recurrence, and we obtain:
\begin{eqnarray*}
I_p
&=&\frac{p-1}{p}\,I_{p-2}\\
&=&\frac{p-1}{p}\cdot\frac{p-3}{p-2}\,I_{p-4}\\
&=&\frac{p-1}{p}\cdot\frac{p-3}{p-2}\cdot\frac{p-5}{p-4}\,I_{p-6}\\
&&\vdots\\
&=&\frac{p!!}{(p+1)!!}\,I_{1-\varepsilon(p)}
\end{eqnarray*}

But $I_0=\frac{\pi}{2}$ and $I_1=1$, so we get the result. As for the second formula, this follows from the first one, with $t=\frac{\pi}{2}-s$. Thus, we have proved both formulae in the statement.
\end{proof}

We can now compute the volumes of the $N$-dimensional spheres, as follows:

\index{volume of sphere}
\index{double factorial}

\begin{theorem}
The volume of the unit sphere in $\mathbb R^N$ is given by
$$V=\left(\frac{\pi}{2}\right)^{[N/2]}\frac{2^N}{(N+1)!!}$$
with our usual convention $N!!=(N-1)(N-3)(N-5)\ldots$
\end{theorem}

\begin{proof}
Let us denote by $B^+$ the positive part of the unit sphere, or rather unit ball $B$, obtained by cutting this unit ball in $2^N$ parts. At the level of volumes, we have:
$$V=2^NV^+$$

We have the following computation, using spherical coordinates:
\begin{eqnarray*}
V^+
&=&\int_{B^+}1\\
&=&\int_0^1\int_0^{\pi/2}\ldots\int_0^{\pi/2}r^{N-1}\sin^{N-2}t_1\ldots\sin t_{N-2}\,drdt_1\ldots dt_{N-1}\\
&=&\int_0^1r^{N-1}\,dr\int_0^{\pi/2}\sin^{N-2}t_1\,dt_1\ldots\int_0^{\pi/2}\sin t_{N-2}dt_{N-2}\int_0^{\pi/2}1dt_{N-1}\\
&=&\frac{1}{N}\times\left(\frac{\pi}{2}\right)^{[N/2]}\times\frac{(N-2)!!}{(N-1)!!}\cdot\frac{(N-3)!!}{(N-2)!!}\ldots\frac{2!!}{3!!}\cdot\frac{1!!}{2!!}\cdot1\\
&=&\frac{1}{N}\times\left(\frac{\pi}{2}\right)^{[N/2]}\times\frac{1}{(N-1)!!}\\
&=&\left(\frac{\pi}{2}\right)^{[N/2]}\frac{1}{(N+1)!!}
\end{eqnarray*}

Thus, we obtain the formula in the statement.
\end{proof}

Let us record as well the formula of the area of the sphere, as follows:

\index{area of sphere}

\begin{theorem}
The area of the unit sphere in $\mathbb R^N$ is given by
$$A=\left(\frac{\pi}{2}\right)^{[N/2]}\frac{2^N}{(N-1)!!}$$
with our usual convention $N!!=(N-1)(N-3)(N-5)\ldots$
\end{theorem}

\begin{proof}
As shown by the pizza argument from chapter 1, when talking about $\pi$, which extends to $N$ dimensions, the area and volume of the sphere in $\mathbb R^N$ are related by:
$$A=N\cdot V$$

Together with the formula in Theorem 13.12 for $V$, this gives the result.
\end{proof}

The formula in Theorem 13.12 is certainly nice, but in practice, we would like to have estimates for that sphere volumes too. For this purpose, we will need:

\index{Stirling formula}
\index{Riemann sum}

\begin{theorem}
We have the Stirling formula
$$N!\simeq\left(\frac{N}{e}\right)^N\sqrt{2\pi N}$$
valid in the $N\to\infty$ limit.
\end{theorem}

\begin{proof}
This is something quite tricky, the idea being as follows:

\medskip

(1) Let us first see what we can get with Riemann sums. We have:
$$\log(N!)
=\sum_{k=1}^N\log k
\approx\int_1^N\log x\,dx
=N\log N-N+1$$

By exponentiating, this gives the following estimate, which is not bad:
$$N!\approx\left(\frac{N}{e}\right)^N\cdot e$$

(2) We can improve our estimate by replacing the rectangles from the Riemann sum approach to the integrals by trapezoids. In practice, this gives the following estimate:
$$\log(N!)
\approx\int_1^N\log x\,dx+\frac{\log 1+\log N}{2}
=N\log N-N+1+\frac{\log N}{2}$$

By exponentiating, this gives the following estimate, which gets us closer:
$$N!\approx\left(\frac{N}{e}\right)^N\cdot e\cdot\sqrt{N}$$

(3) In order to conclude, we must take some kind of mathematical magnifier, and carefully estimate the error made in (2). Fortunately, this mathematical magnifier exists, called Euler-Maclaurin formula, and after some computations, this leads to:
$$N!\simeq\left(\frac{N}{e}\right)^N\sqrt{2\pi N}$$

(4) However, all this remains a bit complicated, so we would like to present now an alternative approach to (3), which also misses some details, but better does the job, explaining where the $\sqrt{2\pi}$ factor comes from. First, by partial integration we have:
$$N!=\int_0^\infty x^Ne^{-x}dx$$

Since the integrand is sharply peaked at $x=N$, as you can see by computing the derivative of $\log(x^Ne^{-x})$, this suggests writing $x=N+y$, and we obtain:
\begin{eqnarray*}
\log(x^Ne^{-x})
&=&N\log x-x\\
&=&N\log(N+y)-(N+y)\\
&=&N\log N+N\log\left(1+\frac{y}{N}\right)-(N+y)\\
&\simeq&N\log N+N\left(\frac{y}{N}-\frac{y^2}{2N^2}\right)-(N+y)\\
&=&N\log N-N-\frac{y^2}{2N}
\end{eqnarray*}

By exponentiating, we obtain from this the following estimate:
$$x^Ne^{-x}\simeq\left(\frac{N}{e}\right)^Ne^{-y^2/2N}$$

(5) Now by integrating, and using the Gauss formula, we obtain from this:
\begin{eqnarray*}
N!
&=&\int_0^\infty x^Ne^{-x}dx\\
&\simeq&\int_{-N}^N\left(\frac{N}{e}\right)^Ne^{-y^2/2N}\,dy\\
&\simeq&\left(\frac{N}{e}\right)^N\int_\mathbb Re^{-y^2/2N}\,dy\\
&=&\left(\frac{N}{e}\right)^N\sqrt{2\pi N}
\end{eqnarray*}

Thus, we have proved the Stirling formula, as formulated in the statement.
\end{proof}

We can now estimate the volumes of the spheres, as follows:

\begin{theorem}
The volume of the unit sphere in $\mathbb R^N$ is given by
$$V\simeq\left(\frac{2\pi e}{N}\right)^{N/2}\frac{1}{\sqrt{\pi N}}$$
in the $N\to\infty$ limit.
\end{theorem}

\begin{proof}
This is very standard, using the formula in Theorem 13.12,  as follows:

\medskip

(1) The double factorials can be estimated by using the Stirling formula. Indeed, in the case where $N=2K$ is even, we have the following computation:
$$(N+1)!!
=2^KK!
\simeq\left(\frac{2K}{e}\right)^K\sqrt{2\pi K}
=\left(\frac{N}{e}\right)^{N/2}\sqrt{\pi N}$$

(2) As for the case where $N=2K-1$ is odd, here the estimate goes as follows:
\begin{eqnarray*}
(N+1)!!
&=&\frac{(2K)!}{2^KK!}\\
&\simeq&\frac{1}{2^K}\left(\frac{2K}{e}\right)^{2K}\sqrt{4\pi K}\left(\frac{e}{K}\right)^K\frac{1}{\sqrt{2\pi K}}\\
&=&\left(\frac{2K}{e}\right)^K\sqrt{2}\\
&=&\left(\frac{N+1}{e}\right)^{(N+1)/2}\sqrt{2}\\
&=&\left(\frac{N}{e}\right)^{N/2}\left(\frac{N+1}{N}\right)^{N/2}\sqrt{\frac{N+1}{e}}\cdot\sqrt{2}\\
&\simeq&\left(\frac{N}{e}\right)^{N/2}\sqrt{e}\cdot\sqrt{\frac{N}{e}}\cdot\sqrt{2}\\
&=&\left(\frac{N}{e}\right)^{N/2}\sqrt{2N}
\end{eqnarray*}

(3) Now back to the spheres, when $N$ is even, the estimate goes as follows:
\begin{eqnarray*}
V
&=&\left(\frac{\pi}{2}\right)^{N/2}\frac{2^N}{(N+1)!!}\\
&\simeq&\left(\frac{\pi}{2}\right)^{N/2}2^N\left(\frac{e}{N}\right)^{N/2}\frac{1}{\sqrt{\pi N}}\\
&=&\left(\frac{2\pi e}{N}\right)^{N/2}\frac{1}{\sqrt{\pi N}}
\end{eqnarray*}

(4) As for the case where $N$ is odd, here the estimate goes as follows:
\begin{eqnarray*}
V
&=&\left(\frac{\pi}{2}\right)^{(N-1)/2}\frac{2^N}{(N+1)!!}\\
&\simeq&\left(\frac{\pi}{2}\right)^{(N-1)/2}2^N\left(\frac{e}{N}\right)^{N/2}\frac{1}{\sqrt{2N}}\\
&=&\sqrt{\frac{2}{\pi}}\left(\frac{2\pi e}{N}\right)^{N/2}\frac{1}{\sqrt{2N}}\\
&=&\left(\frac{2\pi e}{N}\right)^{N/2}\frac{1}{\sqrt{\pi N}}
\end{eqnarray*}

Thus, we are led to the uniform formula in the statement.
\end{proof}

So long for high dimensional spheres and their volumes. Needless to say, all this is very useful when dealing with harmonic functions are related equations, such as the wave and heat ones, and exercise of course for you, to learn more about all this.

\section*{13d. Spherical integrals} 

Moving on, as yet another useful piece of mathematics, coming from the spherical coordinates, let us discuss now the computation of the arbitrary integrals over the sphere. We will need a technical result extending Theorem 13.11, as follows:

\index{trigonometric integral}

\begin{theorem}[Wallis 2]
We have the following formula,
$$\int_0^{\pi/2}\cos^pt\sin^qt\,dt=\left(\frac{\pi}{2}\right)^{\varepsilon(p)\varepsilon(q)}\frac{p!!q!!}{(p+q+1)!!}$$
where $\varepsilon(p)=1$ if $p$ is even, and $\varepsilon(p)=0$ if $p$ is odd, and where
$$m!!=(m-1)(m-3)(m-5)\ldots$$
with the product ending at $2$ if $m$ is odd, and ending at $1$ if $m$ is even.
\end{theorem}

\begin{proof}
We use the same idea as in Theorem 13.11. Let $I_{pq}$ be the integral in the statement. In order to do the partial integration, observe that we have:
\begin{eqnarray*}
(\cos^pt\sin^qt)'
&=&p\cos^{p-1}t(-\sin t)\sin^qt\\
&+&\cos^pt\cdot q\sin^{q-1}t\cos t\\
&=&-p\cos^{p-1}t\sin^{q+1}t+q\cos^{p+1}t\sin^{q-1}t
\end{eqnarray*}

By integrating between $0$ and $\pi/2$, we obtain, for $p,q>0$:
$$pI_{p-1,q+1}=qI_{p+1,q-1}$$

Thus, we can compute $I_{pq}$ by recurrence. When $q$ is even we have:
\begin{eqnarray*}
I_{pq}
&=&\frac{q-1}{p+1}\,I_{p+2,q-2}\\
&=&\frac{q-1}{p+1}\cdot\frac{q-3}{p+3}\,I_{p+4,q-4}\\
&=&\frac{q-1}{p+1}\cdot\frac{q-3}{p+3}\cdot\frac{q-5}{p+5}\,I_{p+6,q-6}\\
&=&\vdots\\
&=&\frac{p!!q!!}{(p+q)!!}\,I_{p+q}
\end{eqnarray*}

But the last term comes from Theorem 13.11, and we obtain the result:
\begin{eqnarray*}
I_{pq}
&=&\frac{p!!q!!}{(p+q)!!}\,I_{p+q}\\
&=&\frac{p!!q!!}{(p+q)!!}\left(\frac{\pi}{2}\right)^{\varepsilon(p+q)}\frac{(p+q)!!}{(p+q+1)!!}\\
&=&\left(\frac{\pi}{2}\right)^{\varepsilon(p)\varepsilon(q)}\frac{p!!q!!}{(p+q+1)!!}
\end{eqnarray*}

Observe that this gives the result for $p$ even as well, by symmetry. In the remaining case now, where both the exponents $p,q$ are odd, we can use once again the formula $pI_{p-1,q+1}=qI_{p+1,q-1}$ established above, and the recurrence goes this time as follows:
\begin{eqnarray*}
I_{pq}
&=&\frac{q-1}{p+1}\,I_{p+2,q-2}\\
&=&\frac{q-1}{p+1}\cdot\frac{q-3}{p+3}\,I_{p+4,q-4}\\
&=&\frac{q-1}{p+1}\cdot\frac{q-3}{p+3}\cdot\frac{q-5}{p+5}\,I_{p+6,q-6}\\
&=&\vdots\\
&=&\frac{p!!q!!}{(p+q-1)!!}\,I_{p+q-1,1}
\end{eqnarray*}

In order to compute the last term, observe that we have:
$$I_{p1}
=\int_0^{\pi/2}\cos^pt\sin t\,dt
=-\frac{1}{p+1}\int_0^{\pi/2}(\cos^{p+1}t)'\,dt
=\frac{1}{p+1}$$

Thus, we can finish our computation in the case $p,q$ odd, as follows:
\begin{eqnarray*}
I_{pq}
&=&\frac{p!!q!!}{(p+q-1)!!}\,I_{p+q-1,1}\\
&=&\frac{p!!q!!}{(p+q-1)!!}\cdot\frac{1}{p+q}\\
&=&\frac{p!!q!!}{(p+q+1)!!}
\end{eqnarray*}

Thus, we obtain the formula in the statement, the exponent of $\pi/2$ appearing there being $\varepsilon(p)\varepsilon(q)=0\cdot 0=0$ in the present case, and this finishes the proof.
\end{proof}

We can now integrate in full generality over the spheres, as follows:

\index{spherical integral}
\index{double factorials}

\begin{theorem}
The polynomial integrals over the unit sphere $S^{N-1}_\mathbb R\subset\mathbb R^N$, with respect to the normalized, mass $1$ measure, are given by the following formula,
$$\int_{S^{N-1}_\mathbb R}x_1^{k_1}\ldots x_N^{k_N}\,dx=\frac{(N-1)!!k_1!!\ldots k_N!!}{(N+\Sigma k_i-1)!!}$$
valid when all exponents $k_i$ are even. If an exponent $k_i$ is odd, the integral vanishes.
\end{theorem}

\begin{proof}
Assume first that one of the exponents $k_i$ is odd. We can make then the following change of variables, which shows that the integral in the statement vanishes:
$$x_i\to-x_i$$

Assume now that all exponents $k_i$ are even. As a first observation, the result holds indeed at $N=2$, due to the formula from Theorem 13.16, which reads:
$$\int_0^{\pi/2}\cos^pt\sin^qt\,dt
=\left(\frac{\pi}{2}\right)^{\varepsilon(p)\varepsilon(q)}\frac{p!!q!!}{(p+q+1)!!}
=\frac{p!!q!!}{(p+q+1)!!}$$

In the general case now, where the dimension $N\in\mathbb N$ is arbitrary, the integral in the statement can be written in spherical coordinates, as follows:
$$I=\frac{2^N}{A}\int_0^{\pi/2}\ldots\int_0^{\pi/2}x_1^{k_1}\ldots x_N^{k_N}J\,dt_1\ldots dt_{N-1}$$

Here $A$ is the area of the sphere, $J$ is the Jacobian, and the $2^N$ factor comes from the restriction to the $1/2^N$ part of the sphere where all the coordinates are positive. According to Theorem 13.13, the normalization constant in front of the integral is:
$$\frac{2^N}{A}=\left(\frac{2}{\pi}\right)^{[N/2]}(N-1)!!$$

As for the unnormalized integral, this is given by:
\begin{eqnarray*}
I'=\int_0^{\pi/2}\ldots\int_0^{\pi/2}
&&(\cos t_1)^{k_1}
(\sin t_1\cos t_2)^{k_2}\\
&&\vdots\\
&&(\sin t_1\sin t_2\ldots\sin t_{N-2}\cos t_{N-1})^{k_{N-1}}\\
&&(\sin t_1\sin t_2\ldots\sin t_{N-2}\sin t_{N-1})^{k_N}\\
&&\sin^{N-2}t_1\sin^{N-3}t_2\ldots\sin^2t_{N-3}\sin t_{N-2}\\
&&dt_1\ldots dt_{N-1}
\end{eqnarray*}

By rearranging the terms, we obtain the following formula:
\begin{eqnarray*}
I'
&=&\int_0^{\pi/2}\cos^{k_1}t_1\sin^{k_2+\ldots+k_N+N-2}t_1\,dt_1\\
&&\int_0^{\pi/2}\cos^{k_2}t_2\sin^{k_3+\ldots+k_N+N-3}t_2\,dt_2\\
&&\vdots\\
&&\int_0^{\pi/2}\cos^{k_{N-2}}t_{N-2}\sin^{k_{N-1}+k_N+1}t_{N-2}\,dt_{N-2}\\
&&\int_0^{\pi/2}\cos^{k_{N-1}}t_{N-1}\sin^{k_N}t_{N-1}\,dt_{N-1}
\end{eqnarray*}

Now by using the above-mentioned formula at $N=2$, this gives:
\begin{eqnarray*}
I'
&=&\frac{k_1!!(k_2+\ldots+k_N+N-2)!!}{(k_1+\ldots+k_N+N-1)!!}\left(\frac{\pi}{2}\right)^{\varepsilon(N-2)}\\
&&\frac{k_2!!(k_3+\ldots+k_N+N-3)!!}{(k_2+\ldots+k_N+N-2)!!}\left(\frac{\pi}{2}\right)^{\varepsilon(N-3)}\\
&&\vdots\\
&&\frac{k_{N-2}!!(k_{N-1}+k_N+1)!!}{(k_{N-2}+k_{N-1}+l_N+2)!!}\left(\frac{\pi}{2}\right)^{\varepsilon(1)}\\
&&\frac{k_{N-1}!!k_N!!}{(k_{N-1}+k_N+1)!!}\left(\frac{\pi}{2}\right)^{\varepsilon(0)}
\end{eqnarray*}

By simplifying and putting everything together, we obtain from this:
$$I
=\left(\frac{2}{\pi}\right)^{[N/2]}(N-1)!!\times\frac{k_1!!\ldots k_N!!}{(\Sigma k_i+N-1)!!}$$

Thus, we are led to the conclusion in the statement.
\end{proof}

We have the following useful generalization of the above formula:

\index{spherical integral}

\begin{theorem}
We have the following integration formula over the sphere $S^{N-1}_\mathbb R\subset\mathbb R^N$, with respect to the normalized, mass $1$ measure, valid for any exponents $k_i\in\mathbb N$,
$$\int_{S^{N-1}_\mathbb R}|x_1^{k_1}\ldots x_N^{k_N}|\,dx=\left(\frac{2}{\pi}\right)^{\Sigma(k_1,\ldots,k_N)}\frac{(N-1)!!k_1!!\ldots k_N!!}{(N+\Sigma k_i-1)!!}$$
with $\Sigma=[odds/2]$ if $N$ is odd and $\Sigma=[(odds+1)/2]$ if $N$ is even, where ``odds'' denotes the number of odd numbers in the sequence $k_1,\ldots,k_N$.
\end{theorem}

\begin{proof}
As before, the formula holds at $N=2$, due to Theorem 13.16. In general, the integral in the statement can be written in spherical coordinates, as in the proof of Theorem 13.17, with the unnormalized integral being given by the following formula:
\begin{eqnarray*}
I'
&=&\int_0^{\pi/2}\cos^{k_1}t_1\sin^{k_2+\ldots+k_N+N-2}t_1\,dt_1\\
&&\int_0^{\pi/2}\cos^{k_2}t_2\sin^{k_3+\ldots+k_N+N-3}t_2\,dt_2\\
&&\vdots\\
&&\int_0^{\pi/2}\cos^{k_{N-2}}t_{N-2}\sin^{k_{N-1}+k_N+1}t_{N-2}\,dt_{N-2}\\
&&\int_0^{\pi/2}\cos^{k_{N-1}}t_{N-1}\sin^{k_N}t_{N-1}\,dt_{N-1}
\end{eqnarray*}

Now by using the Wallis formula at $N=2$, we get:
\begin{eqnarray*}
I'
&=&\frac{\pi}{2}\cdot\frac{k_1!!(k_2+\ldots+k_N+N-2)!!}{(k_1+\ldots+k_N+N-1)!!}\left(\frac{2}{\pi}\right)^{\delta(k_1,k_2+\ldots+k_N+N-2)}\\
&&\frac{\pi}{2}\cdot\frac{k_2!!(k_3+\ldots+k_N+N-3)!!}{(k_2+\ldots+k_N+N-2)!!}\left(\frac{2}{\pi}\right)^{\delta(k_2,k_3+\ldots+k_N+N-3)}\\
&&\vdots\\
&&\frac{\pi}{2}\cdot\frac{k_{N-2}!!(k_{N-1}+k_N+1)!!}{(k_{N-2}+k_{N-1}+k_N+2)!!}\left(\frac{2}{\pi}\right)^{\delta(k_{N-2},k_{N-1}+k_N+1)}\\
&&\frac{\pi}{2}\cdot\frac{k_{N-1}!!k_N!!}{(k_{N-1}+k_N+1)!!}\left(\frac{2}{\pi}\right)^{\delta(k_{N-1},k_N)}
\end{eqnarray*}

In order to compute this latter quantity, let us denote by $F$ the part involving the double factorials, and by $P$ the part involving the powers of $\pi/2$. We have:
$$F=\frac{k_1!!\ldots k_N!!}{(\Sigma k_i+N-1)!!}$$

As in what regards $P$, the above $\delta$ exponents sum up to the following number:
$$\Delta(k_1,\ldots,k_N)=\sum_{i=1}^{N-1}\delta(k_i,k_{i+1}+\ldots+k_N+N-i-1)$$

In other words, with this notation, the formula of $I'$ found before reads:
\begin{eqnarray*}
I'
&=&\left(\frac{2}{\pi}\right)^{\Delta(k_1,\ldots,k_N)-N+1}\frac{k_1!!k_2!!\ldots k_N!!}{(k_1+\ldots+k_N+N-1)!!}\\
&=&\left(\frac{2}{\pi}\right)^{\Sigma(k_1,\ldots,k_N)-[N/2]}\frac{k_1!!k_2!!\ldots k_N!!}{(k_1+\ldots+k_N+N-1)!!}
\end{eqnarray*}

Together with $I=(2^N/V)I'$, this gives the formula in the statement.
\end{proof}

Finally, we have the following related result, dealing with the complex sphere:

\index{spherical integral}

\begin{theorem}
We have the following integration formula over the complex sphere $S^{N-1}_\mathbb C\subset\mathbb C^N$, with respect to the normalized uniform measure, 
$$\int_{S^{N-1}_\mathbb C}|z_1|^{2k_1}\ldots|z_N|^{2k_N}\,dz=\frac{(N-1)!k_1!\ldots k_n!}{(N+\sum k_i-1)!}$$
valid for any exponents $k_i\in\mathbb N$. As for the other polynomial integrals in $z_1,\ldots,z_N$ and their conjugates $\bar{z}_1,\ldots,\bar{z}_N$, these all vanish.
\end{theorem}

\begin{proof}
Consider an arbitrary polynomial integral over $S^{N-1}_\mathbb C$, containing the same number of plain and conjugated variables, as to not vanish trivially, written as follows:
$$I=\int_{S^{N-1}_\mathbb C}z_{i_1}\bar{z}_{i_2}\ldots z_{i_{2k-1}}\bar{z}_{i_{2k}}\,dz$$

By using transformations of type $p\to\lambda p$ with $|\lambda|=1$, we see that this integral $I$ vanishes, unless each $z_a$ appears as many times as $\bar{z}_a$ does, and this gives the last assertion. So, assume now that we are in the non-vanishing case. Then the $k_a$ copies of $z_a$ and the $k_a$ copies of $\bar{z}_a$ produce by multiplication a factor $|z_a|^{2k_a}$, so we have:
$$I=\int_{S^{N-1}_\mathbb C}|z_1|^{2k_1}\ldots|z_N|^{2k_N}\,dz$$

Now by using the standard identification $S^{N-1}_\mathbb C\simeq S^{2N-1}_\mathbb R$, we obtain:
\begin{eqnarray*}
I
&=&\int_{S^{2N-1}_\mathbb R}(x_1^2+y_1^2)^{k_1}\ldots(x_N^2+y_N^2)^{k_N}\,d(x,y)\\
&=&\sum_{r_1\ldots r_N}\binom{k_1}{r_1}\ldots\binom{k_N}{r_N}\int_{S^{2N-1}_\mathbb R}x_1^{2k_1-2r_1}y_1^{2r_1}\ldots x_N^{2k_N-2r_N}y_N^{2r_N}\,d(x,y)
\end{eqnarray*}

By using the formula in Theorem 13.17, we obtain:
\begin{eqnarray*}
&&I\\
&=&\sum_{r_1\ldots r_N}\binom{k_1}{r_1}\ldots\binom{k_N}{r_N}\frac{(2N-1)!!(2r_1)!!\ldots(2r_N)!!(2k_1-2r_1)!!\ldots (2k_N-2r_N)!!}{(2N+2\sum k_i-1)!!}\\
&=&\sum_{r_1\ldots r_N}\binom{k_1}{r_1}\ldots\binom{k_N}{r_N}\frac{2^{N-1}(N-1)!\prod(2r_i)!/(2^{r_i}r_i!)\prod(2k_i-2r_i)!/(2^{k_i-r_i}(k_i-r_i)!)}{2^{N+\sum k_i-1}(N+\sum k_i-1)!}\\
&=&\sum_{r_1\ldots r_N}\binom{k_1}{r_1}\ldots\binom{k_N}{r_N}
\frac{(N-1)!(2r_1)!\ldots (2r_N)!(2k_1-2r_1)!\ldots (2k_N-2r_N)!}{4^{\sum k_i}(N+\sum k_i-1)!r_1!\ldots r_N!(k_1-r_1)!\ldots (k_N-r_N)!}
\end{eqnarray*}

Now observe that can rewrite this quantity in the following way:
\begin{eqnarray*}
&&I\\
&=&\sum_{r_1\ldots r_N}\frac{k_1!\ldots k_N!(N-1)!(2r_1)!\ldots (2r_N)!(2k_1-2r_1)!\ldots (2k_N-2r_N)!}{4^{\sum k_i}(N+\sum k_i-1)!(r_1!\ldots r_N!(k_1-r_1)!\ldots (k_N-r_N)!)^2}\\
&=&\sum_{r_1}\binom{2r_1}{r_1}\binom{2k_1-2r_1}{k_1-r_1}\ldots\sum_{r_N}\binom{2r_N}{r_N}\binom{2k_N-2r_N}{k_N-r_N}\frac{(N-1)!k_1!\ldots k_N!}{4^{\sum k_i}(N+\sum k_i-1)!}\\
&=&4^{k_1}\times\ldots\times 4^{k_N}\times\frac{(N-1)!k_1!\ldots k_N!}{4^{\sum k_i}(N+\sum k_i-1)!}\\
&=&\frac{(N-1)!k_1!\ldots k_N!}{(N+\sum k_i-1)!}
\end{eqnarray*}

Thus, we are led to the formula in the statement.
\end{proof}

We will see applications of all this in the next chapter, with some quite conceptual results regarding the spherical coordinates, in the $N\to\infty$ limit, obtained by processing the above results, by using standard tools from probability.

\section*{13e. Exercises}

Here are some exercises on the above, often insisting on missing details:

\begin{exercise}
Clarify the details of the change of variable formula.
\end{exercise}

\begin{exercise}
Find the ranges of angles in the spherical coordinate formula.
\end{exercise}

\begin{exercise}
Further refine the Stirling formula, with more terms.
\end{exercise}

\begin{exercise}
Prove $\sum_r\binom{2r}{r}\binom{2k-2r}{k-r}=4^k$, used in the above, at the end.
\end{exercise}

As a bonus exercise, compute areas and volumes, as many as you can.

\chapter{Normal variables}

\section*{14a. Probability basics}

Good news, with the calculus that we know, we can have a crash course in advanced probability. We have already met some probability in this book, scattered at the end of chapters 4, 6, 7, and we will first review that material. Then, helped by the Gauss integral formula $\int_\mathbb R e^{-x^2}=\sqrt{\pi}$ that we just learned, and by multivariable integration in general, we will develop the theory of normal variables, real and complex, which is central to advanced probability. And then, there will certainly be more, all interesting things. 

\bigskip

Sounds exciting, doesn't it. Cat however seems unfazed, and declares:

\begin{cat}
Probability is the same thing as measure theory. Read Rudin.
\end{cat}

Damn cat, looks like we disagree more and more as time goes by, especially on these multivariable integration topics. So, let me ask you cat, how many times did you come upon a function which is not integrable? How many times did you fail catching a mouse, due to a failure of Fubini? What about catching birds, does the Zorn lemma really help there? Things are nice and smooth in life, or at least that's my belief.

\begin{cat}
Yes for mice and birds, but electrons can be quite tricky.
\end{cat}

Humm, good point, but let's leave electrons for later, for chapter 16. So, forgetting now about philosophy, and pedagogy matters, but dear reader feel free to have your own opinion here, and why not agreeing with cat, and going ahead with our plan, let us first review the few things that we know about probability, from chapters 4, 6, 7. 

\bigskip

There will be of course a bit of redundancy, but always good to talk again about that things, with a bit less details of course, this time. As a starting point, we have:

\index{probability space}
\index{random variable}
\index{moments}
\index{law}
\index{distribution}

\begin{definition}
Let $X$ be a probability space, that is, a space with a probability measure, and with the corresponding integration denoted $E$, and called expectation.
\begin{enumerate}
\item The random variables are the real functions $f\in L^\infty(X)$.

\item The moments of such a variable are the numbers $M_k(f)=E(f^k)$.

\item The law of such a variable is the measure given by $M_k(f)=\int_\mathbb Rx^kd\mu_f(x)$.
\end{enumerate}
\end{definition}

Here, as explained in chapter 7, the fact that $\mu_f$ as above exists indeed is not exactly trivial. But we can do this by looking at formulae of the following type:
$$E(\varphi(f))=\int_\mathbb R\varphi(x)d\mu_f(x)$$

Indeed, having this for monomials $\varphi(x)=x^n$, as above, is the same as having it for polynomials $\varphi\in\mathbb R[X]$, which in turn is the same as having it for the characteristic functions $\varphi=\chi_I$ of measurable sets $I\subset\mathbb R$. Thus, in the end, what we need is:
$$P(f\in I)=\mu_f(I)$$

But this latter formula can serve as a definition for $\mu_f$, and we are done. 

\bigskip

Regarding now independence, as explained in chapter 7, we can formulate here:

\index{independence}

\begin{definition}
Two variables $f,g\in L^\infty(X)$ are called independent when
$$E(f^kg^l)=E(f^k)\,E(g^l)$$
happens, for any $k,l\in\mathbb N$.
\end{definition}

Again, this definition hides some non-trivial things, the idea being a bit as before, namely that of looking at formulae of the following type:
$$E[\varphi(f)\psi(g)]=E[\varphi(f)]\,E[\psi(g)]$$

To be more precise, passing as before from monomials to polynomials, then to characteristic functions, we are led to the usual definition of independence, namely:
$$P(f\in I,g\in J)=P(f\in I)\,P(g\in J)$$

As a first result now, that we know well from chapter 7, we have:

\index{convolution}

\begin{theorem}
Assuming that $f,g\in L^\infty(X)$ are independent, we have
$$\mu_{f+g}=\mu_f*\mu_g$$
where $*$ is the convolution of real probability measures.
\end{theorem}

\begin{proof}
We have the following computation, using the independence of $f,g$:
$$\int_\mathbb Rx^kd\mu_{f+g}(x)
=E((f+g)^k)
=\sum_r\binom{k}{r}M_r(f)M_{k-r}(g)$$

On the other hand, we have as well the following computation:
\begin{eqnarray*}
\int_\mathbb Rx^kd(\mu_f*\mu_g)(x)
&=&\int_{\mathbb R\times\mathbb R}(x+y)^kd\mu_f(x)d\mu_g(y)\\
&=&\sum_r\binom{k}{r}M_r(f)M_{k-r}(g)
\end{eqnarray*}

Thus $\mu_{f+g}$ and $\mu_f*\mu_g$ have the same moments, so they coincide, as claimed.
\end{proof}

As a second result on independence, which is more advanced, we have:

\index{independence}
\index{Fourier transform}

\begin{theorem}
Assuming that $f,g\in L^\infty(X)$ are independent, we have
$$F_{f+g}=F_fF_g$$
where $F_f(x)=E(e^{ixf})$ is the Fourier transform.
\end{theorem}

\begin{proof}
This is something that we know too from chapter 7, coming from:
\begin{eqnarray*}
F_{f+g}(x)
&=&\int_\mathbb Re^{ixz}d(\mu_f*\mu_g)(z)\\
&=&\int_{\mathbb R\times\mathbb R}e^{ix(z+t)}d\mu_f(z)d\mu_g(t)\\
&=&\int_\mathbb Re^{ixz}d\mu_f(z)\int_\mathbb Re^{ixt}d\mu_g(t)\\
&=&F_f(x)F_g(x)
\end{eqnarray*}

Thus, we are led to the conclusion in the statement.
\end{proof}

Getting now to more advanced theory, a key problem is that of recovering a probability measure out of its moments. And here, we know from chapter 6 that we have:

\begin{theorem}
The density of a real probability measure $\mu$ can be recaptured from the sequence of moments $\{M_k\}_{k\geq0}$ via the Stieltjes inversion formula
$$d\mu (x)=\lim_{t\searrow 0}-\frac{1}{\pi}\,Im\left(G(x+it)\right)\cdot dx$$
where the function on the right, given in terms of moments by
$$G(\xi)=\xi^{-1}+M_1\xi^{-2}+M_2\xi^{-3}+\ldots$$
is the Cauchy transform of the measure $\mu$.
\end{theorem}

\begin{proof}
This is something quite subtle and heavy, and for the full proof, along with some basic applications, we refer to chapter 6, the idea being as follows:

\medskip

(1) Regarding the proof, the Cauchy transform of our measure $\mu$ is given by:
$$G(\xi)
=\xi^{-1}\sum_{k=0}^\infty M_k\xi^{-k}
=\int_\mathbb R\frac{1}{\xi-y}\,d\mu(y)$$

Now with $\xi=x+it$, we obtain from this the following formula:
\begin{eqnarray*}
Im(G(x+it))
&=&\int_\mathbb RIm\left(\frac{1}{x-y+it}\right)d\mu(y)\\
&=&-\int_\mathbb R\frac{t}{(x-y)^2+t^2}\,d\mu(y)
\end{eqnarray*}

By integrating over $[a,b]$ we obtain, with the change of variables $x=y+tz$:
\begin{eqnarray*}
\int_a^bIm(G(x+it))dx
&=&-\int_\mathbb R\int_{(a-y)/t}^{(b-y)/t}\frac{t}{(tz)^2+t^2}\,t\,dz\,d\mu(y)\\
&=&-\int_\mathbb R\int_{(a-y)/t}^{(b-y)/t}\frac{1}{1+z^2}\,dz\,d\mu(y)\\
&=&-\int_\mathbb R\left(\arctan\frac{b-y}{t}-\arctan\frac{a-y}{t}\right)d\mu(y)
\end{eqnarray*}

(2) The point now is that with $t\searrow0$ we have the following estimates:
$$\lim_{t\searrow0}\left(\arctan\frac{b-y}{t}-\arctan\frac{a-y}{t}\right)
=\begin{cases}
\frac{\pi}{2}-\frac{\pi}{2}=0& (y<a)\\
\frac{\pi}{2}-0=\frac{\pi}{2}& (y=a)\\
\frac{\pi}{2}-(-\frac{\pi}{2})=\pi& (a<y<b)\\
0-(-\frac{\pi}{2})=\frac{\pi}{2}& (y=b)\\
-\frac{\pi}{2}-(-\frac{\pi}{2})=0& (y>b)
\end{cases}$$

We therefore obtain the following formula, which proves our result:
$$\lim_{t\searrow0}\int_a^bIm(G(x+it))dx=-\pi\left(\mu(a,b)+\frac{\mu(a)+\mu(b)}{2}\right)$$

(3) As applications, the first thing that we did in chapter 6 was to find the measure having as even moments the Catalan numbers, $C_k=\frac{1}{k+1}\binom{2k}{k}$, and having all odd moments $0$. We found the following measure, called Wigner semicircle law on $[-2,2]$:
$$\gamma_1=\frac{1}{2\pi}\sqrt{4-x^2}dx$$

(4) We also computed the measure having as moments the Catalan numbers, $M_k=C_k$. We found the following measure, called Marchenko-Pastur law on $[0,4]$:
$$\pi_1=\frac{1}{2\pi}\sqrt{4x^{-1}-1}\,dx$$

(5) Next, we found that the measure having as moments the central binomial coefficients, $D_k=\binom{2k}{k}$, is the following measure, called arcsine law on $[0,4]$:
$$\alpha_1=\frac{1}{\pi\sqrt{x(4-x)}}\,dx$$

(6) Finally, for the middle binomial coefficients, $E_k=\binom{k}{[k/2]}$, we found the following law on $[-2,2]$, called modified arcsine law on $[-2,2]$:
$$\sigma_1=\frac{1}{2\pi}\sqrt{\frac{2+x}{2-x}}\,dx$$

(7) All this is very nice, and as already mentioned in chapter 6, these measures are those appearing via random walks on basic graphs. In addition, these measures are as well the main laws in Random Matrix Theory (RMT). More on them later.
\end{proof}

We have as well the following result, also from chapter 6, which can be useful too:

\begin{theorem}
A sequence of numbers $M_0,M_1,M_2,M_3,\ldots\in\mathbb R$, with $M_0=1$, is the series of moments of a real probability measure $\mu$ precisely when:
$$\begin{vmatrix}M_0\end{vmatrix}\geq0\quad,\quad 
\begin{vmatrix}
M_0&M_1\\
M_1&M_2
\end{vmatrix}\geq0\quad,\quad 
\begin{vmatrix}
M_0&M_1&M_2\\
M_1&M_2&M_3\\
M_2&M_3&M_4\\
\end{vmatrix}\geq0\quad,\quad 
\ldots$$
That is, the associated Hankel determinants must be all positive.
\end{theorem}

\begin{proof}
This is something a bit heavier, and as a first observation, the positivity conditions in the statement tell us that the following linear forms must be positive:
$$\sum_{i,j=1}^nc_i\bar{c}_jM_{i+j}\geq0$$

But this is something very classical, in one sense the result being elementary, coming from the following computation, which shows that we have positivity indeed:
$$\int_\mathbb R\left|\sum_{i=1}^nc_ix^i\right|^2d\mu(x)
=\int_\mathbb R\sum_{i,j=1}^nc_i\bar{c}_jx^{i+j}d\mu(x)
=\sum_{i,j=1}^nc_i\bar{c}_jM_{i+j}$$

As for the other sense, here the result comes once again from the above formula, this time via some standard study, inspired from the positivity results from chapter 12.
\end{proof}

Finally, still in relation with this, moment problem, let us formulate:

\begin{definition}
The orthogonal polynomials with respect to $d\mu(x)=f(x)dx$ are polynomials $P_k\in\mathbb R[x]$ of degree $k\in\mathbb N$, which are orthogonal inside $H=L^2(\mathbb R,\mu)$:
$$\int_\mathbb RP_k(x)P_l(x)f(x)dx=0\quad,\quad\forall k\neq l$$
Equivalently, these orthogonal polynomials $\{P_k\}_{k\in\mathbb N}$, which are each unique modulo scalars, appear from the Weierstrass basis $\{x^k\}_{k\in\mathbb N}$, by doing Gram-Schmidt.
\end{definition}

Observe that these polynomials exist indeed for any real measure $d\mu(x)=f(x)dx$, as a consequence of our general discussion from chapter 7, involving Hilbert spaces, and the Weierstrass theorem. It is possible to be a bit more explicit here, as follows:

\index{moments of measure}

\begin{theorem}
The orthogonal polynomials with respect to $\mu$ are given by
$$P_k=c_k\begin{vmatrix}
M_0&M_1&\ldots&M_k\\
M_1&M_2&\ldots&M_{k+1}\\
\vdots&\vdots&&\vdots\\
M_{k-1}&M_k&\ldots&M_{2k-1}\\
1&x&\ldots&x^k
\end{vmatrix}$$
where $M_k=\int_\mathbb Rx^kd\mu(x)$ are the moments of $\mu$, and $c_k\in\mathbb R^*$ can be any numbers.
\end{theorem}

\begin{proof}
Let us first see what happens at small values of $k\in\mathbb N$. At $k=0$ our formula is as follows, stating that the first polynomial $P_0$ must be a constant, as it should:
$$P_0=c_0|M_0|=c_0$$

At $k=1$ now, again by using $M_0=1$, the formula is as follows:
$$P_1=c_1\begin{vmatrix}M_0&M_1\\ 1&x\end{vmatrix}=c_1(x-M_1)$$

But this is again the good formula, because the degree is 1, and we have:
\begin{eqnarray*}
<1,P_1>
&=&c_1<1,x-M_1>\\
&=&c_1(<1,x>-<1,M_1>)\\
&=&c_1(M_1-M_1)\\
&=&0
\end{eqnarray*}

At $k=2$ now, things get more complicated, with the formula being as follows:
$$P_2=c_2\begin{vmatrix}
M_0&M_1&M_2\\
M_1&M_2&M_3\\
1&x&x^2
\end{vmatrix}$$

However, no need for big computations here, in order to check the orthogonality, because by using the fact that $x^k$ integrates up to $M_k$, we obtain:
$$<1,P_2>=\int_\mathbb RP_2(x)d\mu(x)=c_2\begin{vmatrix}
M_0&M_1&M_2\\
M_1&M_2&M_3\\
M_0&M_1&M_2
\end{vmatrix}=0$$

Similarly, again by using the fact that $x^k$ integrates up to $M_k$, we have as well:
$$<x,P_2>=\int_\mathbb RxP_2(x)d\mu(x)=c_2\begin{vmatrix}
M_0&M_1&M_2\\
M_1&M_2&M_3\\
M_1&M_2&M_3
\end{vmatrix}=0$$

Thus, result proved at $k=0,1,2$, and the proof in general is similar.
\end{proof}

We will be back to all this later, in chapter 16, when doing quantum mechanics.

\section*{14b. Discrete laws}

All the above is very nice, we have some interesting theory going on. Let us discuss now some illustrations. We will first talk about discrete probability. First, we have:

\index{Bernoulli law}

\begin{definition}
The Bernoulli law of parameter $x\in[0,1]$ is the law
$$\rho_x=(1-x)\delta_0+x\delta_1$$
appearing when flipping a biased coin, $P({\rm heads})=x$, $P({\rm tails})=1-x$.
\end{definition}

To be more precise, when flipping a biased coin as above, and betting heads, your winning law is $\rho_x$. Next, let us flip the biased coin several times in a row. This leads to:

\index{binomial law}
\index{biased coin}
\index{convolution}
\index{independence}

\begin{theorem}
When flipping a $x$-biased coin $n$ times in a row, the law is
$$\rho_{xn}=\sum_{k=0}^n\binom{n}{k}x^k(1-x)^{n-k}\delta_k$$
called binomial law of parameters $x\in[0,1]$ and $n\in\mathbb N$.
\end{theorem}

\begin{proof}
This is something very standard, the idea being as follows:

\medskip

(1) Observe first that at $n=1$ we have indeed the Bernoulli law $\rho_x$. 

\medskip

(2) In general, we can argue that when flipping the coin $n$ times in a row, and betting heads, the probability of winning $k$ times, among our $n$ attempts, is given by:
$$P(k\ {\rm wins})=\binom{n}{k}P({\rm heads})^kP({\rm tails})^{n-k}=\binom{n}{k}x^k(1-x)^{n-k}$$

Thus, we are led to the formula of $\rho_{xn}$ in the statement.

\medskip

(3) Alternatively, and being a bit more formal, since our $n$ coin tosses are independent, and independence corresponds to convolution, at the level of laws, we have:
\begin{eqnarray*}
\rho_{xn}
&=&\rho_x^{*n}\\
&=&\Big[(1-x)\delta_0+x\delta_1\Big]^{*n}\\
&=&\sum_{k=0}^n\binom{n}{k}x^k(1-x)^{n-k}\delta_k
\end{eqnarray*}

(4) Thus, one way or another, we are led to the formula in the statement.
\end{proof}

Getting now to the study of the binomial laws, we have here:

\index{mean}
\index{variance}
\index{moments}

\begin{theorem}
The binomial law $\rho_{xn}$ has the following properties:
\begin{enumerate}
\item The mean is $E=nx$.

\item The variance is $V=nx(1-x)$.
\end{enumerate}
\end{theorem}

\begin{proof}
In what regards the mean, the computation is as follows:
\begin{eqnarray*}
E
&=&\sum_{k=1}^nk\binom{n}{k}x^k(1-x)^{n-k}\\
&=&\sum_{k=1}^n\frac{n!}{(k-1)!(n-k)!}x^k(1-x)^{n-k}\\
&=&nx\sum_{k=1}^n\frac{(n-1)!}{(k-1)!(n-k)!}x^{k-1}(1-x)^{n-k}\\
&=&nx\sum_{t=0}^{n-1}\binom{n-1}{t}x^t(1-x)^{n-t-1}\\
&=&nx(x+1-x)^{n-1}\\
&=&nx
\end{eqnarray*}

With the same trick, we can compute the difference of the first two moments:
\begin{eqnarray*}
M_2-M_1
&=&\sum_{k=2}^n(k^2-k)\binom{n}{k}x^k(1-x)^{n-k}\\
&=&\sum_{k=2}^n\frac{n!}{(k-2)!(n-k)!}x^k(1-x)^{n-k}\\
&=&n(n-1)x^2\sum_{k=2}^n\frac{(n-2)!}{(k-2)!(n-k)!}x^{k-2}(1-x)^{n-k}\\
&=&n(n-1)x^2\sum_{t=0}^{n-2}\binom{n-2}{t}x^t(1-x)^{n-t-2}\\
&=&n(n-1)x^2(x+1-x)^{n-2}\\
&=&n(n-1)x^2
\end{eqnarray*}

We conclude that the second moment is given by the following formula:
$$M_2=n(n-1)x^2+nx=nx((n-1)x+1)$$

As for the variance $V=M_2-M_1^2$, this is given by the following formula:
$$V=nx((n-1)x+1)-(nx)^2=nx(1-x)$$

Thus, we are led to the conclusions in the statement.
\end{proof}

Many other things can be said about the binomial laws, and among others, as mentioned in chapter 7, these can be used in order to prove the Weierstrass approximation theorem, with the computations there crucially using the formulae in Theorem 14.13.

\bigskip

Moving on, the central objects in discrete probability theory are the Poisson laws:

\index{Poisson law}

\begin{definition}
The Poisson law of parameter $t>0$ is the measure
$$p_t=e^{-t}\sum_{k\in\mathbb N}\frac{t^k}{k!}\,\delta_k$$
with the letter ``p'' standing for Poisson.
\end{definition}

We already talked about these laws, first in chapter 4, with an elementary discussion, using binomials and factorials, and then in chapter 7 too. So, let us quickly review what we know. Going directly for the kill, Fourier transform computation, we have:

\index{Fourier transform}
\index{convolution semigroup}

\begin{theorem}
The Fourier transform of $p_t$ is given by:
$$F_{p_t}(y)=\exp\left((e^{iy}-1)t\right)$$
In particular we have $p_s*p_t=p_{s+t}$, called convolution semigroup property.
\end{theorem}

\begin{proof}
We have indeed the following computation, that we know from chapter 7:
\begin{eqnarray*}
F_{p_t}(y)
&=&e^{-t}\sum_k\frac{t^k}{k!}F_{\delta_k}(y)\\
&=&e^{-t}\sum_k\frac{t^k}{k!}\,e^{iky}\\
&=&e^{-t}\sum_k\frac{(e^{iy}t)^k}{k!}\\
&=&\exp\left((e^{iy}-1)t\right)
\end{eqnarray*}

As for the second assertion, this follows from the fact that $\log F_{p_t}$ is linear in $t$, via the linearization property for the convolution from Theorem 14.6.
\end{proof}

We can now establish the Poisson Limit Theorem, as follows:

\index{PLT}
\index{Poisson Limit Theorem}
\index{Bernoulli laws}

\begin{theorem}[PLT]
We have the following convergence, in moments,
$$\left(\left(1-\frac{t}{n}\right)\delta_0+\frac{t}{n}\delta_1\right)^{*n}\to p_t$$
for any $t>0$.
\end{theorem}

\begin{proof}
Indeed, if we denote by $\nu_n$ the measure under the convolution sign, we have the following computation, for the Fourier transform of the limit: 
\begin{eqnarray*}
F_{\delta_r}(y)=e^{iry}
&\implies&F_{\nu_n}(y)=\left(1-\frac{t}{n}\right)+\frac{t}{n}e^{iy}\\
&\implies&F_{\nu_n^{*n}}(y)=\left(\left(1-\frac{t}{n}\right)+\frac{t}{n}e^{iy}\right)^n\\
&\implies&F_{\nu_n^{*n}}(y)=\left(1+\frac{(e^{iy}-1)t}{n}\right)^n\\
&\implies&F(y)=\exp\left((e^{iy}-1)t\right)
\end{eqnarray*}

Thus, we obtain indeed the Fourier transform of $p_t$, as desired.
\end{proof}

Finally, one more thing that we know about the Poisson laws are some interesting formulae for their moments, from chapter 4. The result there was as follows:

\index{Bell numbers}
\index{partitions}

\begin{theorem}
The moments of $p_1$ are the Bell numbers,
$$M_k(p_1)=|P(k)|$$
where $P(k)$ is the set of partitions of $\{1,\ldots,k\}$. More generally, we have
$$M_k(p_t)=\sum_{\pi\in P(k)}t^{|\pi|}$$
for any $t>0$, where $|.|$ is the number of blocks.
\end{theorem}

\begin{proof}
We know that the moments of $p_1$ are given by the following formula:
$$M_k=\frac{1}{e}\sum_r\frac{r^k}{r!}$$

We therefore have the following recurrence formula for these moments:
\begin{eqnarray*}
M_{k+1}
&=&\frac{1}{e}\sum_r\frac{r^k}{r!}\left(1+\frac{1}{r}\right)^k\\
&=&\frac{1}{e}\sum_r\frac{r^k}{r!}\sum_s\binom{k}{s}r^{-s}\\
&=&\sum_s\binom{k}{s}M_{k-s}
\end{eqnarray*}

But the Bell numbers $B_k=|P(k)|$ satisfy the same recurrence, so we have $M_k=B_k$, as claimed. As for the proof of the formula at $t>0$ arbitrary, this is similar. 
\end{proof}

\section*{14c. Normal variables}

Getting now to the continuous case, as a key application of the Gauss integral formula, established in chapter 13, we can introduce the normal laws, as follows:

\begin{definition}
The normal law of parameter $1$ is the following measure:
$$g_1=\frac{1}{\sqrt{2\pi}}e^{-x^2/2}dx$$
More generally, the normal law of parameter $t>0$ is the following measure:
$$g_t=\frac{1}{\sqrt{2\pi t}}e^{-x^2/2t}dx$$
These are also called Gaussian distributions, with ``g'' standing for Gauss.
\end{definition}

Observe that the above laws have indeed mass 1, as they should. This follows indeed from the Gauss formula, which gives, with $x=\sqrt{2t}\,y$:
\begin{eqnarray*}
\int_\mathbb R e^{-x^2/2t}dx
&=&\int_\mathbb R e^{-y^2}\sqrt{2t}\,dy\\
&=&\sqrt{2t}\int_\mathbb R e^{-y^2}dy\\
&=&\sqrt{2\pi t}
\end{eqnarray*}

Generally speaking, the normal laws appear as bit everywhere, in real life. The reasons behind this phenomenon come from the Central Limit Theorem (CLT), that we will explain in a moment, after developing some general theory. As a first result, we have:

\begin{proposition}
We have the variance formula
$$V(g_t)=t$$
valid for any $t>0$.
\end{proposition}

\begin{proof}
The first moment is 0, because our normal law $g_t$ is centered. As for the second moment, this can be computed as follows:
\begin{eqnarray*}
M_2
&=&\frac{1}{\sqrt{2\pi t}}\int_\mathbb Rx^2e^{-x^2/2t}dx\\
&=&\frac{1}{\sqrt{2\pi t}}\int_\mathbb R(tx)\left(-e^{-x^2/2t}\right)'dx\\
&=&\frac{1}{\sqrt{2\pi t}}\int_\mathbb Rte^{-x^2/2t}dx\\
&=&t
\end{eqnarray*}

We conclude from this that the variance is $V=M_2=t$.
\end{proof}

Here is another result, which is the key one for the study of the normal laws:

\begin{theorem}
We have the following formula, valid for any $t>0$:
$$F_{g_t}(x)=e^{-tx^2/2}$$
In particular, the normal laws satisfy $g_s*g_t=g_{s+t}$, for any $s,t>0$.
\end{theorem}

\begin{proof}
The Fourier transform formula can be established as follows:
\begin{eqnarray*}
F_{g_t}(x)
&=&\frac{1}{\sqrt{2\pi t}}\int_\mathbb Re^{-y^2/2t+ixy}dy\\
&=&\frac{1}{\sqrt{2\pi t}}\int_\mathbb Re^{-(y/\sqrt{2t}-\sqrt{t/2}ix)^2-tx^2/2}dy\\
&=&\frac{1}{\sqrt{2\pi t}}\int_\mathbb Re^{-z^2-tx^2/2}\sqrt{2t}dz\\
&=&\frac{1}{\sqrt{\pi}}e^{-tx^2/2}\int_\mathbb Re^{-z^2}dz\\
&=&\frac{1}{\sqrt{\pi}}e^{-tx^2/2}\cdot\sqrt{\pi}\\
&=&e^{-tx^2/2}
\end{eqnarray*}

As for the last assertion, this follows from the fact that $\log F_{g_t}$ is linear in $t$, via the linearization property for the convolution from Theorem 14.6.
\end{proof}

We are now ready to state and prove the CLT, as follows:

\begin{theorem}[CLT]
Given random variables $f_1,f_2,f_3,\ldots\in L^\infty(X)$ which are i.i.d., centered, and with variance $t>0$, we have, with $n\to\infty$, in moments,
$$\frac{1}{\sqrt{n}}\sum_{i=1}^nf_i\sim g_t$$
where $g_t$ is the Gaussian law of parameter $t$, having as density $\frac{1}{\sqrt{2\pi t}}e^{-y^2/2t}dy$.
\end{theorem}

\begin{proof}
In terms of moments, the Fourier transform is given by:
\begin{eqnarray*}
F_f(x)
&=&E\left(\sum_{k=0}^\infty\frac{(ixf)^k}{k!}\right)\\
&=&\sum_{k=0}^\infty\frac{(ix)^kE(f^k)}{k!}\\
&=&\sum_{k=0}^\infty\frac{i^kM_k(f)}{k!}\,x^k
\end{eqnarray*}

We conclude that the Fourier transform of the variable in the statement is:
\begin{eqnarray*}
F(x)
&=&\left[F_f\left(\frac{x}{\sqrt{n}}\right)\right]^n\\
&=&\left[1-\frac{tx^2}{2n}+O(n^{-2})\right]^n\\
&\simeq&\left[1-\frac{tx^2}{2n}\right]^n\\
&\simeq&e^{-tx^2/2}
\end{eqnarray*}

But this latter function being the Fourier transform of $g_t$, we obtain the result.
\end{proof}

Let us discuss now some further properties of the normal law. We first have:

\begin{proposition}
The even moments of the normal law are the numbers
$$M_k(g_t)=t^{k/2}\times k!!$$
where $k!!=(k-1)(k-3)(k-5)\ldots\,$, and the odd moments vanish. 
\end{proposition}

\begin{proof}
We have the following computation, valid for any integer $k\in\mathbb N$:
\begin{eqnarray*}
M_k
&=&\frac{1}{\sqrt{2\pi t}}\int_\mathbb Ry^ke^{-y^2/2t}dy\\
&=&\frac{1}{\sqrt{2\pi t}}\int_\mathbb R(ty^{k-1})\left(-e^{-y^2/2t}\right)'dy\\
&=&\frac{1}{\sqrt{2\pi t}}\int_\mathbb Rt(k-1)y^{k-2}e^{-y^2/2t}dy\\
&=&t(k-1)\times\frac{1}{\sqrt{2\pi t}}\int_\mathbb Ry^{k-2}e^{-y^2/2t}dy\\
&=&t(k-1)M_{k-2}
\end{eqnarray*}

Thus by recurrence, we are led to the formula in the statement.
\end{proof}

We have the following alternative formulation of the above result:

\begin{proposition}
The moments of the normal law are the numbers
$$M_k(g_t)=t^{k/2}|P_2(k)|$$
where $P_2(k)$ is the set of pairings of $\{1,\ldots,k\}$.
\end{proposition}

\begin{proof}
Let us count the pairings of $\{1,\ldots,k\}$. In order to have such a pairing, we must pair $1$ with one of the numbers $2,\ldots,k$, and then use a pairing of the remaining $k-2$ numbers. Thus, we have the following recurrence formula:
$$|P_2(k)|=(k-1)|P_2(k-2)|$$

As for the initial data, this is $P_1=0$, $P_2=1$. Thus, we are led to the result.
\end{proof}

We are not done yet, and here is one more improvement of the above:

\begin{theorem}
The moments of the normal law are the numbers
$$M_k(g_t)=\sum_{\pi\in P_2(k)}t^{|\pi|}$$
where $P_2(k)$ is the set of pairings of $\{1,\ldots,k\}$, and $|.|$ is the number of blocks.
\end{theorem}

\begin{proof}
This follows indeed from Proposition 14.23, because the number of blocks of a pairing of $\{1,\ldots,k\}$ is trivially $k/2$, independently of the pairing.
\end{proof}

Let us discuss now the complex analogues of all this, with a notion of complex normal, or Gaussian law. To start with, we have the following definition:

\begin{definition}
The complex normal, or Gaussian law of parameter $t>0$ is
$$G_t=law\left(\frac{1}{\sqrt{2}}(a+ib)\right)$$
where $a,b$ are independent, each following the law $g_t$.
\end{definition}

In short, the complex normal laws appear as natural complexifications of the real normal laws. As in the real case, these measures form convolution semigroups:

\begin{proposition}
The complex Gaussian laws have the property
$$G_s*G_t=G_{s+t}$$
for any $s,t>0$, and so they form a convolution semigroup.
\end{proposition}

\begin{proof}
This follows indeed from the real result, namely $g_s*g_t=g_{s+t}$, established in Theorem 14.20, simply by taking the real and imaginary parts.
\end{proof}

We have as well the following complex analogue of the CLT:

\index{CCLT}
\index{Complex CLT}

\begin{theorem}[CCLT]
Given complex variables $f_1,f_2,f_3,\ldots\in L^\infty(X)$ which are i.i.d., centered, and with common variance $t>0$, we have
$$\frac{1}{\sqrt{n}}\sum_{i=1}^nf_i\sim G_t$$
with $n\to\infty$, in moments.
\end{theorem}

\begin{proof}
This follows indeed from the real CLT, established in Theorem 14.21, simply by taking the real and imaginary parts of all the variables involved.
\end{proof}

Regarding now the moments, the situation here is more complicated than in the real case, because in order to have good results, we have to deal with both the complex variables, and their conjugates. Let us formulate the following definition:

\index{colored integers}
\index{colored moments}

\begin{definition}
The moments a complex variable $f\in L^\infty(X)$ are the numbers
$$M_k=E(f^k)$$
depending on colored integers $k=\circ\bullet\bullet\circ\ldots\,$, with the conventions
$$f^\emptyset=1\quad,\quad f^\circ=f\quad,\quad f^\bullet=\bar{f}$$
and multiplicativity, in order to define the colored powers $f^k$.
\end{definition}

Observe that, since the variables $f,\bar{f}$ commute, we can permute terms, and restrict the attention to exponents of type $k=\ldots\circ\circ\circ\bullet\bullet\bullet\ldots\,$, if we want to. However, our results about the complex Gaussian laws, and about some other complex laws too, later, will actually look better without doing this, so we will use Definition 14.28 as stated. 

\bigskip

As an illustration, the moments indexed by low length colored integers are as follows:
$$M_\emptyset=1$$
$$M_\circ=E(f)\quad,\quad M_\bullet=E(\bar{f})$$
$$M_{\circ\circ}=E(f^2)\quad,\quad 
M_{\circ\bullet}=M_{\bullet\circ}=E(f\bar{f})\quad,\quad
M_{\bullet\bullet}=E(\bar{f}^2)$$
$$\vdots$$

Regarding now the moments of the complex normal laws, we have here:

\begin{theorem}
The moments of the complex normal law are given by
$$M_k(G_t)=\begin{cases}
t^pp!&(k\ {\rm uniform, of\ length}\ 2p)\\
0&(k\ {\rm not\ uniform})
\end{cases}$$
where $k=\circ\bullet\bullet\circ\ldots$ is called uniform when it contains the same number of $\circ$ and $\bullet$.
\end{theorem}

\begin{proof}
We must compute the moments, with respect to colored integer exponents $k=\circ\bullet\bullet\circ\ldots$\,, of the variable from Definition 14.25, namely:
$$f=\frac{1}{\sqrt{2}}(a+ib)$$

We can assume that we are in the case $t=1$, and the proof here goes as follows:

\medskip

(1) As a first observation, in the case where our exponent $k=\circ\bullet\bullet\circ\ldots$ is not uniform, a standard rotation argument shows that the corresponding moment of $f$ vanishes. To be more precise, the variable $f'=wf$ is complex Gaussian too, for any complex number $w\in\mathbb T$, and from $M_k(f)=M_k(f')$ we obtain $M_k(f)=0$, in this case.

\medskip

(2) In the uniform case now, where the exponent $k=\circ\bullet\bullet\circ\ldots$ consists of $p$ copies of $\circ$ and $p$ copies of $\bullet$\,, the corresponding moment can be computed as follows:
\begin{eqnarray*}
M_k
&=&\int(f\bar{f})^p\\
&=&\frac{1}{2^p}\int(a^2+b^2)^p\\
&=&\frac{1}{2^p}\sum_r\binom{p}{r}\int a^{2r}\int b^{2p-2r}\\
&=&\frac{1}{2^p}\sum_r\binom{p}{r}(2r)!!(2p-2r)!!\\
&=&\frac{1}{2^p}\sum_r\frac{p!}{r!(p-r)!}\cdot\frac{(2r)!}{2^rr!}\cdot\frac{(2p-2r)!}{2^{p-r}(p-r)!}\\
&=&\frac{p!}{4^p}\sum_r\binom{2r}{r}\binom{2p-2r}{p-r}
\end{eqnarray*}

(3) In order to finish now the computation, let us recall that we have the following formula, coming from the generalized binomial formula, or from the Taylor formula:
$$\frac{1}{\sqrt{1+t}}=\sum_{q=0}^\infty\binom{2q}{q}\left(\frac{-t}{4}\right)^q$$

By taking the square of this series, we obtain the following formula:
$$\frac{1}{1+t}
=\sum_p\left(\frac{-t}{4}\right)^p\sum_r\binom{2r}{r}\binom{2p-2r}{p-r}$$

Now by looking at the coefficient of $t^p$ on both sides, we conclude that the sum on the right equals $4^p$. Thus, we can finish the moment computation in (2), as follows:
$$M_k=\frac{p!}{4^p}\times 4^p=p!$$

We are therefore led to the conclusion in the statement.
\end{proof}

As before with the real Gaussian laws, a better-looking statement  is in terms of partitions. Given a colored integer $k=\circ\bullet\bullet\circ\ldots\,$, we say that a pairing $\pi\in P_2(k)$ is matching when it pairs $\circ-\bullet$ symbols. With this convention, we have the following result:

\index{matching pairings}

\begin{theorem}
The moments of the complex normal law are the numbers
$$M_k(G_t)=\sum_{\pi\in\mathcal P_2(k)}t^{|\pi|}$$
where $\mathcal P_2(k)$ are the matching pairings of $\{1,\ldots,k\}$, and $|.|$ is the number of blocks.
\end{theorem}

\begin{proof}
This is a reformulation of Theorem 14.29. Indeed, we can assume that we are in the case $t=1$, and here we know from Theorem 14.29 that the moments are:
$$M_k=\begin{cases}
(|k|/2)!&(k\ {\rm uniform})\\
0&(k\ {\rm not\ uniform})
\end{cases}$$

On the other hand, the numbers $|\mathcal P_2(k)|$ are given by exactly the same formula. Indeed, in order to have a matching pairing of $k$, our exponent $k=\circ\bullet\bullet\circ\ldots$ must be uniform, consisting of $p$ copies of $\circ$ and $p$ copies of $\bullet$, with $p=|k|/2$. But then the matching pairings of $k$ correspond to the permutations of the $\bullet$ symbols, as to be matched with $\circ$ symbols, and so we have $p!$ such pairings. Thus, we have the same formula as for the moments of $f$, and we are led to the conclusion in the statement.
\end{proof}

In practice, we also need to know how to compute joint moments. We have here:

\begin{theorem}[Wick formula]
Given independent variables $f_i$, each following the complex normal law $G_t$, with $t>0$ being a fixed parameter, we have the formula
$$E\left(f_{i_1}^{k_1}\ldots f_{i_s}^{k_s}\right)=t^{s/2}\#\left\{\pi\in\mathcal P_2(k)\Big|\pi\leq\ker i\right\}$$
where $k=k_1\ldots k_s$ and $i=i_1\ldots i_s$, for the joint moments of these variables, where $\pi\leq\ker i$ means that the indices of $i$ must fit into the blocks of $\pi$, in the obvious way.
\end{theorem}

\begin{proof}
This is something well-known, which can be proved as follows:

\medskip

(1) Let us first discuss the case where we have a single variable $f$, which amounts in taking $f_i=f$ for any $i$ in the formula in the statement. What we have to compute here are the moments of $f$, with respect to colored integer exponents $k=\circ\bullet\bullet\circ\ldots\,$, and the formula in the statement tells us that these moments must be:
$$E(f^k)=t^{|k|/2}|\mathcal P_2(k)|$$

But this is the formula in Theorem 14.30, so we are done with this case.

\medskip

(2) In general now, when expanding the product $f_{i_1}^{k_1}\ldots f_{i_s}^{k_s}$ and rearranging the terms, we are left with doing a number of computations as in (1), and then making the product of the expectations that we found. But this amounts in counting the partitions in the statement, with the condition $\pi\leq\ker i$ there standing for the fact that we are doing the various type (1) computations independently, and then making the product.
\end{proof}

So long for the normal laws, real and complex. These appear pretty much everywhere in mathematics, and with the above, you are fully armed for dealing with them. In particular the Wick formula is what you need for Random Matrix Theory (RMT), shall you ever get interested in that, and with a good reference here being Mehta \cite{meh}.

\section*{14d. Geometric aspects}

We can do many things with the probability theory that we learned so far, and especially with the real and complex normal laws. As a first application, let us go back to the spherical integrals, studied in chapter 13. In the real case, we have:

\begin{theorem}
The moments of the hyperspherical variables are
$$\int_{S^{N-1}_\mathbb R}x_i^pdx=\frac{(N-1)!!p!!}{(N+p-1)!!}$$
and the rescaled variables $y_i=\sqrt{N}x_i$ become normal and independent with $N\to\infty$.
\end{theorem}

\begin{proof}
The moment formula in the statement follows from the general formulae in chapter 13. As a consequence, with $N\to\infty$ we have the following estimate:
\begin{eqnarray*}
\int_{S^{N-1}_\mathbb R}x_i^pdx
&\simeq&N^{-p/2}\times p!!\\
&=&N^{-p/2}M_p(g_1)
\end{eqnarray*}

Thus, the rescaled variables $\sqrt{N}x_i$ become normal with $N\to\infty$, as claimed. As for the proof of the asymptotic independence, this is standard too, once again by using the formulae in chapter 13. Indeed, the joint moments of $x_1,\ldots,x_N$ are given by:
\begin{eqnarray*}
\int_{S^{N-1}_\mathbb R}x_1^{k_1}\ldots x_N^{k_N}\,dx
&=&\frac{(N-1)!!k_1!!\ldots k_N!!}{(N+\Sigma k_i-1)!!}\\
&\simeq&N^{-\Sigma k_i}\times k_1!!\ldots k_N!!
\end{eqnarray*}

By rescaling, the joint moments of the variables $y_i=\sqrt{N}x_i$ are given by:
$$\int_{S^{N-1}_\mathbb R}y_1^{k_1}\ldots y_N^{k_N}\,dx\simeq k_1!!\ldots k_N!!$$

Thus, we have multiplicativity, and so independence with $N\to\infty$, as claimed.
\end{proof}

Importantly, we can recover the normal laws as well in connection with the rotation groups. Indeed, we have the following reformulation of Theorem 14.32:

\begin{theorem}
We have the integration formula
$$\int_{O_N}U_{ij}^p\,dU=\frac{(N-1)!!p!!}{(N+p-1)!!}$$
and the rescaled variables $V_{ij}=\sqrt{N}U_{ij}$ become normal and independent with $N\to\infty$.
\end{theorem}

\begin{proof}
We use the basic fact that the rotations $U\in O_N$ act on the points of the real sphere $z\in S^{N-1}_\mathbb R$, with the stabilizer of $z=(1,0,\ldots,0)$ being the subgroup $O_{N-1}\subset O_N$. In algebraic terms, this gives an identification as follows:
$$S^{N-1}_\mathbb R=O_N/O_{N-1}$$

In functional analytic terms, this result provides us with an embedding as follows, for any $i$, which makes correspond the respective integration functionals:
$$C(S^{N-1}_\mathbb R)\subset C(O_N)\quad,\quad 
x_i\to U_{1i}$$

With this identification made, the result follows from Theorem 14.32.
\end{proof}

In the complex case, the analogues of the above results are as follows:

\begin{theorem}
The rescalings $\sqrt{N}z_i$ of the unit complex sphere coordinates
$$z_i:S^{N-1}_\mathbb C\to\mathbb C$$
as well as the rescalings $\sqrt{N}U_{ij}$ of the unitary group coordinates
$$U_{ij}:U_N\to\mathbb C$$
become complex Gaussian and independent with $N\to\infty$. 
\end{theorem}

\begin{proof}
We have several assertions to be proved, the idea being as follows:

\medskip

(1) According to the formulae in chapter 13, the polynomials integrals in $z_i,\bar{z}_i$ vanish, unless the number of $z_i,\bar{z}_i$ is the same. In this latter case these terms can be grouped together, by using $z_i\bar{z}_i=|z_i|^2$, and the relevant integration formula is:
$$\int_{S^{N-1}_\mathbb C}|z_i|^{2k}\,dz=\frac{(N-1)!k!}{(N+k-1)!}$$

Now with $N\to\infty$, we obtain from this the following estimate:
$$\int_{S^{N-1}_\mathbb C}|z_i|^{2k}dx
\simeq N^{-k}\times k!$$

Thus, the rescaled variables $\sqrt{N}z_i$ become normal with $N\to\infty$, as claimed. 

\medskip

(2) As for the proof of the asymptotic independence, this is standard too, again by using the formulae in chapter 13. Indeed, the joint moments of $z_1,\ldots,z_N$ are given by:
\begin{eqnarray*}
\int_{S^{N-1}_\mathbb R}|z_1|^{2k_1}\ldots|z_N|^{2k_N}\,dx
&=&\frac{(N-1)!k_1!\ldots k_n!}{(N+\sum k_i-1)!}\\
&\simeq&N^{-\Sigma k_i}\times k_1!\ldots k_N!
\end{eqnarray*}

By rescaling, the joint moments of the variables $y_i=\sqrt{N}z_i$ are given by:
$$\int_{S^{N-1}_\mathbb R}|y_1|^{2k_1}\ldots|y_N|^{2k_N}\,dx\simeq k_1!\ldots k_N!$$

Thus, we have multiplicativity, and so independence with $N\to\infty$, as claimed.

\medskip

(3) Regarding the last assertion, we can use here the basic fact that the rotations $U\in U_N$ act on the points of the sphere $z\in S^{N-1}_\mathbb C$, with the stabilizer of $z=(1,0,\ldots,0)$ being the subgroup $U_{N-1}\subset U_N$. In algebraic terms, this gives an equality as follows:
$$S^{N-1}_\mathbb C=U_N/U_{N-1}$$

In functional analytic terms, this result provides us with an embedding as follows, for any $i$, which makes correspond the respective integration functionals:
$$C(S^{N-1}_\mathbb C)\subset C(U_N)\quad,\quad 
x_i\to U_{1i}$$

With this identification made, the result follows from (1,2).
\end{proof}

All the above is quite nice, and suggests that there are deep relations between group theory and probability. So, let us discuss this now, by starting with the discrete case, which is the simplest one. We have here the following beautiful, well-known result:

\begin{theorem}
The probability for a random $\sigma\in S_N$ to have no fixed points is
$$P\simeq\frac{1}{e}$$
in the $N\to\infty$ limit, where $e=2.71828\ldots$ is the usual constant from analysis.
\end{theorem}

\begin{proof}
This is best viewed by using the inclusion-exclusion principle. Let us set:
$$S_N^k=\left\{\sigma\in S_N\Big|\sigma(k)=k\right\}$$

The set of permutations having no fixed points, called derangements, is then:
$$X_N=\left(\bigcup_kS_N^k\right)^c$$

Now the inclusion-exclusion principle tells us that we have:
\begin{eqnarray*}
|X_N|
&=&\left|\left(\bigcup_kS_N^k\right)^c\right|\\
&=&|S_N|-\sum_k|S_N^k|+\sum_{k<l}|S_N^k\cap S_N^l|-\ldots+(-1)^N\sum_{k_1<\ldots<k_N}|S_N^{k_1}\cap\ldots\cap S_N^{k_N}|\\
&=&N!-N(N-1)!+\binom{N}{2}(N-2)!-\ldots+(-1)^N\binom{N}{N}(N-N)!\\
&=&\sum_{r=0}^N(-1)^r\binom{N}{r}(N-r)!
\end{eqnarray*}

Thus, the probability that we are interested in, for a random permutation $\sigma\in S_N$ to have no fixed points, is given by the following formula:
$$P
=\frac{|X_N|}{N!}
=\sum_{r=0}^N\frac{(-1)^r}{r!}$$

Since on the right we have the expansion of $1/e$, this gives the result.
\end{proof}

The above result looks quite exciting. In order to further explore and refine it, we will need some notions from group theory, and more specifically, the following definition:

\begin{definition} 
Given a closed subgroup $G\subset U_N$, the function
$$\chi:G\to\mathbb C\quad,\quad 
\chi(g)=\sum_{i=1}^Ng_{ii}$$
is called main character of $G$.
\end{definition}

We should mention that there is a long story with these characters, the idea being that a compact group $G$ can have several representations $\pi:G\subset U_N$, which can be studied via their characters $\chi_\pi:G\to\mathbb C$, and with all this coming from Weyl \cite{wy1}, \cite{wy2}. 

\bigskip

All this is quite heavy, but technically, we will not need all this here. Indeed, in relation with our questions, we can formulate the following nice, elementary result:

\begin{theorem}
Consider the symmetric group $S_N$, regarded as the permutation group, $S_N\subset O_N$, of the $N$ coordinate axes of $\mathbb R^N$.
\begin{enumerate}
\item The main character $\chi\in C(S_N)$ counts the number of fixed points.

\item The law of $\chi\in C(S_N)$ becomes Poisson $(1)$, in the $N\to\infty$ limit.
\end{enumerate}
\end{theorem}

\begin{proof}
We have two things to be proved here, the idea being as follows:

\medskip

(1) The permutation matrices $\sigma\in O_N$, which give the embedding $S_N\subset O_N$ in the statement, being given by $\sigma_{ij}=\delta_{i\sigma(j)}$, we have the following computation:
$$\chi(\sigma)
=\sum_i\delta_{\sigma(i)i}
=\#\left\{i\in\{1,\ldots,N\}\Big|\sigma(i)=i\right\}$$

(2) In order to establish now the asymptotic result in the statement, we must prove the following formula, for any $r\in\mathbb N$, in the $N\to\infty$ limit:
$$P(\chi=r)\simeq\frac{1}{r!e}$$

We already know, from Theorem 14.35, that this formula holds at $r=0$. In the general case now, we have to count the permutations $\sigma\in S_N$ having exactly $r$ points. Now since having such a permutation amounts in choosing $r$ points among $1,\ldots,N$, and then permuting the $N-r$ points left, without fixed points allowed, we have:
\begin{eqnarray*}
\#\left\{\sigma\in S_N\Big|\chi(\sigma)=r\right\}
&=&\binom{N}{r}\#\left\{\sigma\in S_{N-r}\Big|\chi(\sigma)=0\right\}\\
&=&\frac{N!}{r!(N-r)!}\#\left\{\sigma\in S_{N-r}\Big|\chi(\sigma)=0\right\}\\
&=&N!\times\frac{1}{r!}\times\frac{\#\left\{\sigma\in S_{N-r}\Big|\chi(\sigma)=0\right\}}{(N-r)!}
\end{eqnarray*}

By dividing everything by $N!$, we obtain from this the following formula:
$$\frac{\#\left\{\sigma\in S_N\Big|\chi(\sigma)=r\right\}}{N!}=\frac{1}{r!}\times\frac{\#\left\{\sigma\in S_{N-r}\Big|\chi(\sigma)=0\right\}}{(N-r)!}$$

Now by using the computation at $r=0$, that we already have, from Theorem 14.35, it follows that with $N\to\infty$ we have the following estimate:
$$P(\chi=r)
\simeq\frac{1}{r!}\cdot P(\chi=0)
\simeq\frac{1}{r!}\cdot\frac{1}{e}$$

Thus, we obtain as limiting measure the Poisson law of parameter 1, as stated.
\end{proof}

As a next step, let us try now to generalize what we have, namely Theorem 14.37, as to reach to the Poisson laws of arbitrary parameter $t>0$. Again we will need some notions from group theory, and more specifically, the following definition:

\begin{definition}
Given a closed subgroup $G\subset U_N$, the function
$$\chi_t:G\to\mathbb C\quad,\quad 
\chi_t(g)=\sum_{i=1}^{[tN]}g_{ii}$$
is called main truncated character of $G$, of parameter $t\in(0,1]$.
\end{definition}

As before with plain characters, there is some theory behind this definition, involving this time Random Matrix Theory (RMT), but technically, we will not need this here. Indeed, we can now formulate the following result, which is nice and elementary, generalizing Theorem 14.37, and which will be our final saying on the subject:

\begin{theorem}
Consider the symmetric group $S_N$, regarded as the permutation group, $S_N\subset O_N$, of the $N$ coordinate axes of $\mathbb R^N$.
\begin{enumerate}
\item The variable $\chi_t$ counts the number of fixed points among $1,\ldots,[tN]$.

\item The law of this variable $\chi_t$ becomes Poisson $(t)$, in the $N\to\infty$ limit.
\end{enumerate}
\end{theorem}

\begin{proof}
We already know from Theorem 14.37 that the results hold at $t=1$. In general, the proof is similar, the idea being as follows:

\medskip

(1) We have indeed the following computation, coming from definitions:
$$\chi_t(\sigma)
=\sum_{i=1}^{[tN]}\delta_{\sigma(i)i}
=\#\left\{i\in\{1,\ldots,[tN]\}\Big|\sigma(i)=i\right\}$$

(2) Consider indeed the following sets, as in the proof of Theorem 14.35:
$$S_N^k=\left\{\sigma\in S_N\Big|\sigma(k)=k\right\}$$

The set of permutations having no fixed points among $1,\ldots,[tN]$ is then:
$$X_N=\left(\bigcup_{k\leq[tN]}S_N^k\right)^c$$

As before in the proof of Theorem 14.37, we obtain by inclusion-exclusion that:
\begin{eqnarray*}
P(\chi_t=0)
&=&\frac{1}{N!}\sum_{r=0}^{[tN]}(-1)^r\sum_{k_1<\ldots<k_r<[tN]}|S_N^{k_1}\cap\ldots\cap S_N^{k_r}|\\
&=&\frac{1}{N!}\sum_{r=0}^{[tN]}(-1)^r\binom{[tN]}{r}(N-r)!\\
&=&\sum_{r=0}^{[tN]}\frac{(-1)^r}{r!}\cdot\frac{[tN]!(N-r)!}{N!([tN]-r)!}
\end{eqnarray*}

Now with $N\to\infty$, we obtain from this the following estimate:
$$P(\chi_t=0)
\simeq\sum_{r=0}^{[tN]}\frac{(-1)^r}{r!}\cdot t^r
=\sum_{r=0}^{[tN]}\frac{(-t)^r}{r!}
\simeq e^{-t}$$

More generally, by counting the permutations $\sigma\in S_N$ having exactly $r$ fixed points among $1,\ldots,[tN]$, as in the proof of Theorem 14.37, we obtain:
$$P(\chi_t=r)\simeq\frac{t^r}{r!e^t}$$

Thus, we obtain in the limit a Poisson law of parameter $t$, as stated.
\end{proof}

It is possible to establish as well continuous analogues of the above results, involving the groups $O_N,U_N$, and the real and complex normal laws. However, the proofs here are more technical. For an introduction to this, you can have a look at my book \cite{ba2}.

\bigskip

Finally, let us end all this with some low-dimensional magic, involving physicists' favorite group, $SU_2$, and one of their favorite laws, the Wigner law $\gamma_1$. We have: 

\begin{theorem}
The main character of $SU_2$ follows the following law,
$$\gamma_1=\frac{1}{2\pi}\sqrt{4-x^2}dx$$
which is the Wigner law of parameter $1$.
\end{theorem}

\begin{proof}
The group $SU_2$ is by definition the group of unitary rotations $U\in U_2$ of determinant one, and by solving the equation $U^*=U^{-1}$, we are led to:
$$SU_2=\left\{\begin{pmatrix}a+ib&c+id\\ -c+id&a-ib\end{pmatrix}\ \Big|\ a^2+b^2+c^2+d^2=1\right\}$$

In this picture, the main character is given by the following formula:
$$\chi\begin{pmatrix}a+ib&c+id\\ -c+id&a-ib\end{pmatrix}=2a$$

We are therefore left with computing the law of the following variable:
$$a\in C(S^3_\mathbb R)$$

But this is something very familiar, namely a hyperspherical variable at $N=4$, so we can use here Theorem 14.32. We obtain the following moment formula:
$$\int_{S^3_\mathbb R}a^{2k}
=\frac{3!!(2k)!!}{(2k+3)!!}
=\frac{1}{4^k}\cdot\frac{1}{k+1}\binom{2k}{k}
=\frac{C_k}{4^k}$$

Thus the variable $2a\in C(S^3_\mathbb R)$ follows the Wigner semicircle law $\gamma_1$, as claimed.
\end{proof}

There are far more things that can be said, in relation with this, all beautiful mathematics, useful in relation with physics. If interested, have a look at my book \cite{ba2}.

\section*{14e. Exercises}

Probability theory is a big business, and many other things can be said, in relation with the above, which are all useful, and good to know. As exercises, we have:

\begin{exercise}
Learn more about the moment problem, and related topics.
\end{exercise}

\begin{exercise}
Learn about the compound Poisson laws, and their properties.
\end{exercise}

\begin{exercise}
Study more in detail the CLT, and the convergence there.
\end{exercise}

\begin{exercise}
Recover the Marchenko-Pastur law $\pi_1$, in relation with $SO_3$.
\end{exercise}

As bonus exercise, learn some compact group theory, and more specifically the Peter-Weyl theory for such groups. This is quite tough, but is an excellent investment.

\chapter{Partial integration}

\section*{15a. Vector products}

Good news, done with mathematics, and in the remainder of this book, last 50 pages, we will get into physics. We already met some physics in this book, namely gravity, waves and heat in one dimension, in chapters 3-4, then waves and heat again, in two dimensions, along with some basic electrodynamics, light and early atomic theory, in chapter 8, and finally harmonic functions, this time in arbitrary $N$ dimensions, in chapter 12.

\bigskip

Obviously, time to have some cleanup here, review what we know, and do more, in a more systematic way. We have divided what we have two say in two parts, as follows:

\bigskip

(1) In the present chapter we will discuss physics in low dimensions, $N=2,3,4$, with on the menu more classical mechanics, notably with the Coriolis force, then some relativity following Einstein, and finally the basics of electrostatics, following Gauss and others, which are quite interesting, mathematically speaking. We will be mostly interested in $N=3$, and the common feature of everything that we want to talk about will be the vector product $\times$, which exists only in 3 dimensions, and is something quite magic. 

\bigskip

(2) In the next and final chapter we will discuss physics in infinite dimensions, $N=\infty$. Our goal here will be that of having some quantum mechanics theory started, along the lines suggested at the end of chapter 8, and more specifically, solving the hydrogen atom. There are actually several ways of proceeding here, following Heisenberg, Schr\"odinger and others, and matter of annoying my cat, who seems to be a big fan of measure theory and Hilbert spaces, we will opt here for the Schr\"odinger approach, which is elementary.

\bigskip

Getting now to what we want to do in this chapter, namely 3 dimensions, vector products $\times$, and all sorts of related physics, there is some important mathematical interest in doing this, because we will reach in this way to answers to the following question:

\begin{question}
What are the higher dimensional analogues of the formula
$$\int_a^bF'(x)dx=F(b)-F(a)$$
that is, of the fundamental theorem of calculus?
\end{question}

And isn't this a fundamental question for us, mathematicians, because we have so far in our bag multivariable extensions of all the main theorems of one-variable calculus, except for this. So, definitely something to be solved, before the end of this book.

\bigskip

So, let us discuss this first. The fundamental theorem of calculus tells us that the integral of a function on an interval $[a,b]$ can be suitably recaptured from what happens on the boundary $\{a,b\}$ of this interval. Obviously, this is something quite magic, and thinking now at what we can expect in $N=2,3$ or more dimensions, that can only be quite complicated, involving curves, surfaces, solid bodies and so on, and with all this vaguely reminding all sorts of physics things, such as field lines for gravity, magnets and so on. In short, despite having no formal proof yet for all this, let us formulate:

\begin{answer}
Partial integration in several dimensions most likely means physics, and we will probably only get some partial results here, at $N=2,3$.
\end{answer}

Of course, all this remains to be confirmed, but assuming that you trust me a bit, here we are now at the plan that we made before, for this chapter. That is, do physics in low dimensions, guided by the beauty of the world surrounding us, and once this physics done, record some mathematical corollaries too, in relation with Question 15.1.

\bigskip

Before getting started, however, as usual when struggling with pedagogical matters, and other delicate dillemas, let us ask the cat. But here, unfortunately, no surprise:

\begin{cat}
Read Rudin.
\end{cat}

Thanks cat, and guess this sort of discussion, that we started in chapter 13, looks more and more cyclic. So, I'll just go my way, on your side have a good hunt, and by the way make sure to catch enough mice and birds, using your measure theory and differential forms techniques, because there is a bit of shortage of cat food today, sorry for that.

\bigskip

Getting started now, here is what we need:

\index{vector product}
\index{right-hand rule}

\begin{definition}
The vector product of two vectors in $\mathbb R^3$ is given by
$$x\times y=||x||\cdot||y||\cdot\sin\theta\cdot n$$
where $n\in\mathbb R^3$ with $n\perp x,y$ and $||n||=1$ is constructed using the right-hand rule:
$$\begin{matrix}
\ \ \ \ \ \ \ \ \ \uparrow_{x\times y}\\
\leftarrow_x\\
\swarrow_y
\end{matrix}$$
Alternatively, in usual vertical linear algebra notation for all vectors,
$$\begin{pmatrix}x_1\\ x_2\\x_3\end{pmatrix}
\times\begin{pmatrix}y_1\\ y_2\\y_3\end{pmatrix}
=\begin{pmatrix}x_2y_3-x_3y_2\\ x_3y_1-x_1y_3\\x_1y_2-x_2y_1\end{pmatrix}$$
the rule being that of computing $2\times2$ determinants, and adding a middle sign.
\end{definition}

Obviously, this definition is something quite subtle, and also something very annoying, because you always need this, and always forget the formula. Here are my personal methods. With the first definition, what I always remember is that:
$$||x\times y||\sim||x||,||y||\quad,\quad 
x\times x=0\quad,\quad 
e_1\times e_2=e_3$$

So, here's how it works. We are looking for a vector $x\times y$ whose length is proportional to those of $x,y$. But the second formula tells us that the angle $\theta$ between $x,y$ must be involved via $0\to0$, and so the factor can only be $\sin\theta$. And with this we are almost there, it's just a matter of choosing the orientation, and this comes from $e_1\times e_2=e_3$.

\bigskip

As with the second definition, that I like the most, what I remember here is simply:
$$\begin{vmatrix}
1&x_1&y_1\\
1&x_2&y_2\\
1&x_3&y_3
\end{vmatrix}=?$$

Indeed, when trying to compute this determinant, by developing over the first column, what you get as coefficients are the entries of $x\times y$. And with the good middle sign.

\bigskip

In practice now, in order to get familiar with the vector products, nothing better than doing some classical mechanics. We have here the following key result:

\begin{theorem}
In the gravitational $2$-body problem, the angular momentum
$$J=x\times p$$
with $p=mv$ being the usual momentum, is conserved.
\end{theorem}

\begin{proof}
There are several things to be said here, the idea being as follows:

\medskip

(1) First of all the usual momentum, $p=mv$, is not conserved, because the simplest solution is the circular motion, where the moment gets turned around. But this suggests precisely that, in order to fix the lack of conservation of the momentum $p$, what we have to do is to make a vector product with the position $x$. Leading to $J$, as above.

\medskip

(2) Regarding now the proof, consider indeed a particle $m$ moving under the gravitational force of a particle $M$, assumed, as usual, to be fixed at 0. By using the fact that for two proportional vectors, $p\sim q$, we have $p\times q=0$, we obtain:
\begin{eqnarray*}
\dot{J}
&=&\dot{x}\times p+x\times\dot{p}\\
&=&v\times mv+x\times ma\\
&=&m(v\times v+x\times a)\\
&=&m(0+0)\\
&=&0
\end{eqnarray*}

Now since the derivative of $J$ vanishes, this quantity is constant, as stated.
\end{proof}

While the above principle looks like something quite trivial, the mathematics behind it is quite interesting, and has several notable consequences, as follows:

\begin{theorem}
In the context of a $2$-body problem, the following happen:
\begin{enumerate}
\item The fact that the direction of $J$ is fixed tells us that the trajectory of one body with respect to the other lies in a plane.

\item The fact that the magnitude of $J$ is fixed tells us that the Kepler 2 law holds, namely that we have same areas sweeped by $Ox$ over the same times. 
\end{enumerate}
\end{theorem}

\begin{proof}
This follows indeed from Theorem 15.5, as follows:

\medskip

(1) We have by definition $J=m(x\times v)$, and since a vector product is orthogonal on both the vectors it comes from, we deduce from this that we have:
$$J\perp x,v$$

But this can be written as follows, with $J^\perp$ standing for the plane orthogonal to $J$:
$$x,v\in J^\perp$$

Now since $J$ is fixed by Theorem 15.5, we conclude that both $x,v$, and in particular the position $x$, and so the whole trajectory, lie in this fixed plane $J^\perp$, as claimed.

\medskip 

(2) Conversely now, forget about Theorem 15.5, and assume that the trajectory lies in a certain plane $E$. Thus $x\in E$, and by differentiating we have $v\in E$ too, and so $x,v\in E$. Thus $E=J^\perp$, and so $J=E^\perp$, so the direction of $J$ is fixed, as claimed.

\medskip

(3) Regarding now the last assertion, we already know from the various formulae from chapter 11 that the Kepler 2 law is more or less equivalent to the formula $\dot{\theta}=\lambda/r^2$. However, the derivation of $\dot{\theta}=\lambda/r^2$ was something tricky, and what we want to prove now is that this appears as a simple consequence of $||J||$ = constant.

\medskip

(4) In order to to so, let us compute $J$, according to its definition $J=x\times p$, but in polar coordinates, which will change everything. Since $p=m\dot{x}$, we have:
$$J=r\begin{pmatrix}\cos\theta\\ \sin\theta\\ 0\end{pmatrix}
\times m\begin{pmatrix}
\dot{r}\cos\theta-r\sin\theta\cdot\dot{\theta}\\
\dot{r}\sin\theta+r\cos\theta\cdot\dot{\theta}\\
0\end{pmatrix}$$

Now recall from the definition of the vector product that we have:
$$\begin{pmatrix}a\\b\\0\end{pmatrix}
\times\begin{pmatrix}c\\d\\0\end{pmatrix}
=\begin{pmatrix}0\\0\\ad-bc\end{pmatrix}$$

Thus $J$ is a vector of the above form, with its last component being:
\begin{eqnarray*}
J_z
&=&rm\begin{vmatrix}
\cos\theta&&\dot{r}\cos\theta-r\sin\theta\cdot\dot{\theta}\\
\sin\theta&&\dot{r}\sin\theta+r\cos\theta\cdot\dot{\theta}
\end{vmatrix}\\
&=&rm\cdot r(\cos^2\theta+\sin^2\theta)\dot{\theta}\\
&=&r^2m\cdot\dot{\theta}
\end{eqnarray*}

(5) Now with the above formula in hand, our claim is that the magnitude $||J||$ is constant precisely when $\dot{\theta}=\lambda/r^2$, for some $\lambda\in\mathbb R$. Indeed, up to the obvious fact that the orientation of $J$ is a binary parameter, who cannot just switch like that, let us just agree on this, knowing $J$ is the same as knowing $J_z$, and is also the same as knowing $||J||$. Thus, our claim is proved, and this leads to the conclusion in the statement.
\end{proof}

As another basic application of the vector products, still staying with classical mechanics, we have all sorts of useful formulae regarding rotating frames. We first have:

\begin{theorem}
Assume that a $3D$ body rotates along an axis, with angular speed $w$. For a fixed point of the body, with position vector $x$, the usual $3D$ speed is 
$$v=\omega\times x$$ 
where $\omega=wn$, with $n$ unit vector pointing North. When the point moves on the body
$$V=\dot{x}+\omega\times x$$
is its speed computed by an inertial observer $O$ on the rotation axis.
\end{theorem} 

\begin{proof}
We have two assertions here, both requiring some 3D thinking, as follows:

\medskip

(1) Assuming that the point is fixed, the magnitude of $\omega\times x$ is the good one, due to the following computation, with $r$ being the distance from the point to the axis:
$$||\omega\times x||=w||x||\sin t=wr=||v||$$

As for the orientation of $\omega\times x$, this is the good one as well, because the North pole rule used above amounts in applying the right-hand rule for finding $n$, and so $\omega$, and this right-hand rule was precisely the one used in defining the vector products $\times$.

\medskip

(2) Next, when the point moves on the body, the inertial observer $O$ can compute its speed by using a frame $(u_1,u_2,u_3)$ which rotates with the body, as follows:
\begin{eqnarray*}
V
&=&\dot{x}_1u_1+\dot{x}_2u_2+\dot{x}_3u_3+x_1\dot{u}_1+x_2\dot{u}_2+x_3\dot{u}_3\\
&=&\dot{x}+(x_1\cdot\omega\times u_1+x_2\cdot\omega\times u_2+x_3\cdot\omega\times u_3)\\
&=&\dot{x}+w\times(x_1u_1+x_2u_2+x_3u_3)\\
&=&\dot{x}+\omega\times x
\end{eqnarray*}

Thus, we are led to the conclusions in the statement.
\end{proof}

In what regards now the acceleration, the result, which is famous, is as follows:

\begin{theorem}
Assuming as before that a $3D$ body rotates along an axis, the acceleration of a moving point on the body, computed by $O$ as before, is given by
$$A=a+2\omega\times v+\omega\times(\omega\times x)$$
with $\omega=wn$ being as before. In this formula the second term is called Coriolis acceleration, and the third term is called centripetal acceleration.
\end{theorem}

\begin{proof}
This comes by using twice the formulae in Theorem 15.7, as follows:
\begin{eqnarray*}
A
&=&\dot{V}+\omega\times V\\
&=&(\ddot{x}+\dot{\omega}\times x+\omega\times\dot{x})+(\omega\times\dot{x}+\omega\times(\omega\times x))\\
&=&\ddot{x}+\omega\times\dot{x}+\omega\times\dot{x}+\omega\times(\omega\times x)\\
&=&a+2\omega\times v+\omega\times(\omega\times x)
\end{eqnarray*}

Thus, we are led to the conclusion in the statement.
\end{proof}

The truly famous result is actually the one regarding forces, obtained by multiplying everything by a mass $m$, and writing things the other way around, as follows: 
$$ma=mA-2m\omega\times v-m\omega\times(\omega\times x)$$

Here the second term is called Coriolis force, and the third term is called centrifugal force. These forces are both called apparent, or fictious, because they do not exist in the inertial frame, but they exist however in the non-inertial frame of reference, as explained above. And with of course the terms centrifugal and centripetal not to be messed up.

\bigskip

In fact, even more famous is the terrestrial application of all this, as follows:

\begin{theorem}
The acceleration of an object $m$ subject to a force $F$ is given by
$$ma=F-mg-2m\omega\times v-m\omega\times(\omega\times x)$$
with $g$ pointing upwards, and with the last terms being the Coriolis and centrifugal forces.
\end{theorem}

\begin{proof}
This follows indeed from the above discussion, by assuming that the acceleration $A$ there comes from the combined effect of a force $F$, and of the usual $g$.
\end{proof}

We refer to any standard undergraduate mechanics book, such as Feynman \cite{fe1}, Kibble \cite{kbe} or Taylor \cite{tay} for more on the above, including various numerics on what happens here on Earth, the Foucault pendulum, history of all this, and many other things. Let us just mention here, as a basic illustration for all this, that a rock dropped from 100m deviates about 1cm from its intended target, due to the formula in Theorem 15.9.

\section*{15b. Einstein addition}

As another application of the vector products, let us discuss now the speed addition in  relativity. Based on experiments by Fizeau, then Michelson-Morley and others, and some physics by Maxwell and Lorentz, Einstein came upon the following principles:

\begin{fact}[Einstein principles]
The following happen:
\begin{enumerate}
\item Light travels in vacuum at a finite speed, $c<\infty$.

\item This speed $c$ is the same for all inertial observers.

\item In non-vacuum, the light speed is lower, $v<c$.

\item Nothing can travel faster than light, $v\not>c$.
\end{enumerate}
\end{fact} 

The point now is that, obviously, something is wrong here. Indeed, assuming for instance that we have a train, running in vacuum at speed $v>0$, and someone on board lights a flashlight $\ast$ towards the locomotive, an observer $\circ$ on the ground will see the light traveling at speed $c+v>c$, which is a contradiction:
$$
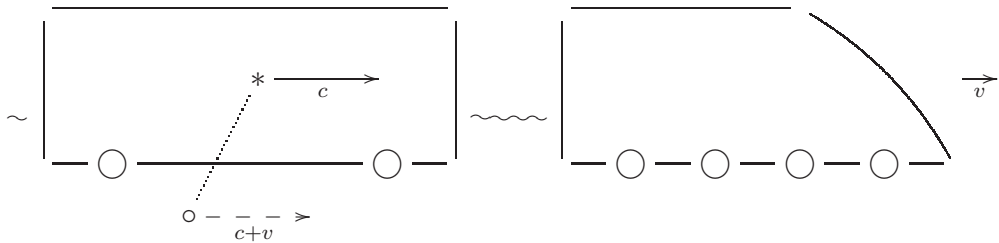
\xymatrix@R=5pt@C=13pt{
&\ar@{-}[rrrrrr]&&&&&&&&\ar@{-}[rrr]&&&\\
&&&&&&&&&\\
&&&&\ast\ar[rr]_c&&&&&&&&&&\ar[r]_v&\\
\ar@{~}[r]&&&&&&&\ar@{~}[rr]&&\\
&\ar@{-}[r]\ar@{-}[uuuu]&\bigcirc\ar@{-}[rrrr]&&&&\bigcirc\ar@{-}[r]&\ar@{-}[uuuu]&&\ar@{-}[r]\ar@{-}[uuuu]&\bigcirc\ar@{-}[r]&\bigcirc\ar@{-}[r]&\bigcirc\ar@{-}[r]&\bigcirc\ar@{-}[r]&\ar@/_/@{-}[uuuull]\\
&&&\circ\ar@{-->}[rr]_{c+v}\ar@{.}[uuur]&&
}$$

Equivalently, with the same train running, in vacuum at speed $v>0$, if the observer on the ground lights a flashlight $\ast$ towards the back of the train, then viewed from the train, that light will travel at speed $c+v>c$, which is a contradiction again:
$$
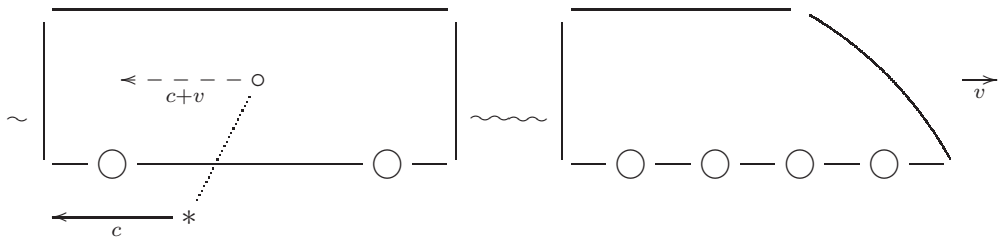
\xymatrix@R=5pt@C=13pt{
&\ar@{-}[rrrrrr]&&&&&&&&\ar@{-}[rrr]&&&\\
&&&&&&&&&\\
&&&&\circ\ar@{-->}[ll]^{c+v}&&&&&&&&&&\ar[r]_v&\\
\ar@{~}[r]&&&&&&&\ar@{~}[rr]&&\\
&\ar@{-}[r]\ar@{-}[uuuu]&\bigcirc\ar@{-}[rrrr]&&&&\bigcirc\ar@{-}[r]&\ar@{-}[uuuu]&&\ar@{-}[r]\ar@{-}[uuuu]&\bigcirc\ar@{-}[r]&\bigcirc\ar@{-}[r]&\bigcirc\ar@{-}[r]&\bigcirc\ar@{-}[r]&\ar@/_/@{-}[uuuull]\\
&&&\ast\ar[ll]^c\ar@{.}[uuur]&&
}$$

Summarizing, Fact 15.10, while physically true, implies $c+v=c$, so contradicts classical mechanics, which needs a fix. In the classical case, to start with, we have:

\begin{proposition}
The classical speeds add according to the Galileo formula
$$v_{AC}=v_{AB}+v_{BC}$$
where $v_{AB}$ denotes the relative speed of $A$ with respect to $B$.
\end{proposition}

\begin{proof}
This is clear indeed from the definition of speed, and very intuitive.
\end{proof}

In order to find the fix, we will first discuss the 1D case, and leave the 3D case, which is a bit more complicated, for later. We will use two tricks. First, let us forget about absolute speeds, with respect to a given frame, and talk about relative speeds only. In this case we are allowed to sum only quantities of type $v_{AB},v_{BC}$, and we denote by $v_{AB}+_gv_{BC}$ the corresponding sum $v_{AC}$. With this convention, the Galileo formula becomes:
$$u+_gv=u+v$$

As a second trick now, observe that this Galileo formula holds in any system of units. In order now to deal with our problems, basically involving high speeds, it is convenient to change the system of units, as to have $c=1$. With this convention our $c+v=c$ problem becomes $1+v=1$, and the solution to it is quite obvious, as follows:

\begin{theorem}
If we define the Einstein sum $+_e$ of relative speeds by
$$u+_ev=\frac{u+v}{1+uv}$$
in $c=1$ units, then we have the formula $1+_ev=1$, valid for any $v$.
\end{theorem}

\begin{proof}
This is obvious indeed from our definition of $+_e$, because if we plug in $u=1$ in the above formula, we obtain as result:
$$1+_ev=\frac{1+v}{1+v}=1$$

Thus, we are led to the conclusion in the statement.
\end{proof}

Summarizing, we have solved our problem. In order now to formulate a final result, we must do some reverse engineering, by waiving the above two tricks. First, by getting back to usual units, $v\to v/c$, our new addition formula becomes:
$$\frac{u}{c}+_e\frac{v}{c}=\frac{\frac{u}{c}+\frac{v}{c}}{1+\frac{u}{c}\cdot\frac{v}{c}}$$

By multiplying by $c$, we can write this formula in a better way, as follows:
$$u+_ev=\frac{u+v}{1+uv/c^2}$$

In order now to finish, it remains to get back to absolute speeds, as in Proposition 15.11. And by doing so, we are led to the following result:

\begin{theorem}
If we sum the speeds according to the Einstein formula
$$v_{AC}=\frac{v_{AB}+v_{BC}}{1+v_{AB}v_{BC}/c^2}$$
then the Galileo formula still holds, approximately, for low speeds
$$v_{AC}\simeq v_{AB}+v_{BC}$$
and if we have $v_{AB}=c$ or $v_{BC}=c$, the resulting sum is $v_{AC}=c$.
\end{theorem}

\begin{proof}
We have two assertions here, which are both clear, as follows:

\medskip

(1) Regarding the first assertion, if we are at low speeds, $v_{AB},v_{BC}<<c$, the correction term $v_{AB}v_{BC}/c^2$ dissapears, and we are left with the Galileo formula, as claimed.

\medskip

(2) As for the second assertion, this follows from the above discussion.
\end{proof}

All the above is very nice, but remember, takes place in 1D. So, time now to get seriously to work, and see what all this becomes in 3D. We have here:

\begin{question}
What is the correct analogue of the Einstein summation formula
$$u+_ev=\frac{u+v}{1+uv}$$
in $2$ and $3$ dimensions?
\end{question}

In order to discuss this question, let us attempt to construct $u+_ev$ in arbitrary dimensions, just by using our common sense and intuition. When the vectors $u,v\in\mathbb R^N$ are proportional, we are basically in 1D, and so our addition formula must satisfy:
$$u\sim v\implies u+_ev=\frac{u+v}{1+<u,v>}$$

However, the formula on the right will not work as such in general, for arbitrary speeds $u,v\in\mathbb R^N$, and this because we have, as main requirement for our operation, in analogy with the $1+_ev=1$ formula from 1D, the following condition:
$$||u||=1\implies u+_ev=u$$

Equivalently, in analogy with $u+_e1=1$ from 1D, we would like to have:
$$||v||=1\implies u+_ev=v$$

Summarizing, our $u\sim v$ formula above is not bad, as a start, but we must add a correction term to it, for the above requirements to be satisfied, and of course with the correction term vanishing when $u\sim v$. So, we are led to a math puzzle:

\begin{puzzle}
What vanishes when $u\sim v$, and then how to correctly define
$$u+_ev=\frac{u+v+\gamma_{uv}}{1+<u,v>}$$
as for the correction term $\gamma_{uv}$ to vanish when $u\sim v$?
\end{puzzle}

But the solution to the first question is well-known in 3D. Indeed, here we can use the vector product $u\times v$, that we met before, which notoriously satisfies:
$$u\sim v\implies u\times v=0$$

Thus, our correction term $\gamma_{uv}$ must be something containing $w=u\times v$, which vanishes when this vector $w$ vanishes, and in addition arranged such that $||u||=1$ produces a simplification, with $u+_ev=u$ as end result, and with $||v||=1$ producing a simplification too, with $u+_ev=v$ as end result. Thus, our vector calculus puzzle becomes:

\begin{puzzle}
How to correctly define the Einstein summation in $3$ dimensions,
$$u+_ev=\frac{u+v+\gamma_{uvw}}{1+<u,v>}$$
with $w=u\times v$, in such a way as for the correction term $\gamma_{uvw}$ to satisfy 
$$w=0\implies\gamma_{uvw}=0$$
and also such that $||u||=1\implies u+_ev=u$, and $||v||\implies u+_ev=v$?
\end{puzzle}

In order to solve this latter puzzle, we must ``transport'' the vector $w$ to the plane spanned by $u,v$. But this is simplest done by taking the vector product with any vector in this plane, and so as a reasonable candidate for our correction term, we have:
$$\gamma_{uvw}=(\alpha u+\beta v)\times w$$

Here $\alpha,\beta\in\mathbb R$ are some scalars to be determined, but let us take a break, and leave the computations for later. We did some good work, time to update our puzzle:

\begin{puzzle}
How to define the Einstein summation in $3$ dimensions,
$$u+_ev=\frac{u+v+\gamma_{uvw}}{1+<u,v>}$$
with the correction term being of the following form, with $w=u\times v$, and $\alpha,\beta\in\mathbb R$,
$$\gamma_{uvw}=(\alpha u+\beta v)\times w$$
in such a way as to have $||u||=1\implies u+_ev=u$, and $||v||\implies u+_ev=v$?
\end{puzzle}

In order to investigate what happens when $||u||=1$ or $||v||=1$, we must compute the vector products $u\times w$ and $v\times w$. So, pausing now our study for consulting the vector calculus database, and then coming back, here is the formula that we need:
$$u\times(u\times v)=<u,v>u-<u,u>v$$

As for the formula of $v\times w$, that I forgot to record, we can recover it from the one above of $u\times w$, by using the basic properties of the vector products, as follows:
\begin{eqnarray*}
v\times(u\times v)
&=&-v\times(v\times u)\\
&=&-(<v,u>v-<v,v>u)\\
&=&<v,v>u-<u,v>v
\end{eqnarray*}

With these formulae in hand, we can now compute the correction term, with the result here, that we will need several times in what comes next, being as follows:

\begin{proposition}
The correction term $\gamma_{uvw}=(\alpha u+\beta v)\times w$ is given by
$$\gamma_{uvw}=(\alpha<u,v>+\beta<v,v>)u-(\alpha<u,u>+\beta<u,v>)v$$
for any values of the scalars $\alpha,\beta\in\mathbb R$.
\end{proposition}

\begin{proof}
According to our vector product formulae above, we have:
\begin{eqnarray*}
\gamma_{uvw}
&=&(\alpha u+\beta v)\times w\\
&=&\alpha(<u,v>u-<u,u>v)+\beta(<v,v>u-<u,v>v)\\
&=&(\alpha<u,v>+\beta<v,v>)u-(\alpha<u,u>+\beta<u,v>)v
\end{eqnarray*}

Thus, we are led to the conclusion in the statement.
\end{proof}

Time now to get into the real thing, see what happens when $||u||=1$ and $||v||=1$, if we can get indeed $u+_ev=u$ and $u+_ev=v$. It is convenient here to do some reverse engineering. Regarding the first desired formula, namely $u+_ev=u$, we have:
\begin{eqnarray*}
u+_ev=u
&\iff&u+v+\gamma_{uvw}=(1+<u,v>)u\\
&\iff&\gamma_{uvw}=<u,v>u-v\\
&\iff&\alpha=1,\ \beta=0,\ ||u||=1
\end{eqnarray*}

Thus, with the parameter choice $\alpha=1,\beta=0$, we will have, as desired:
$$||u||=1\implies u+_ev=u$$

In what regards now the second desired formula, namely $u+_ev=v$, here the computation is almost identical, save for a sign switch, which after some thinking comes from our choice $w=u\times v$ instead of $w=v\times u$, clearly favoring $u$, as follows:
\begin{eqnarray*}
u+_ev=v
&\iff&u+v+\gamma_{uvw}=(1+<u,v>)v\\
&\iff&\gamma_{uvw}=-u+<u,v>v\\
&\iff&\alpha=0,\ \beta=-1,\ ||v||=1
\end{eqnarray*}

Thus, with the parameter choice $\alpha=0,\beta=-1$, we will have, as desired:
$$||v||=1\implies u+_ev=v$$

All this is mixed news, because we managed to solve both our problems, at $||u||=1$ and at $||v||=1$, but our solutions are different. So, time to breathe, decide that we did enough interesting work for the day, and formulate our conclusion as follows:

\begin{proposition}
When defining the Einstein speed summation in $3{\rm D}$ as 
$$u+_ev=\frac{u+v+u\times(u\times v)}{1+<u,v>}$$
in $c=1$ units, the following happen:
\begin{enumerate}
\item When $u\sim v$, we recover the previous $1{\rm D}$ formula.

\item When $||u||=1$, speed of light, we have $u+_ev=u$.

\item However, $||v||=1$ does not imply $u+_ev=v$.

\item Also, the formula $u+_ev=v+_eu$ fails.
\end{enumerate}
\end{proposition}

\begin{proof}
Here (1) and (2) follow from the above discussion, with the following choice for the correction term, by favoring the $||u||=1$ problem over the $||v||=1$ one:
$$\gamma_{uvw}=u\times w$$

In fact, with this choice made, the computation is very simple, as follows:
\begin{eqnarray*}
||u||=1
&\implies&\gamma_{uvw}=<u,v>u-v\\
&\implies&u+v+\gamma_{uvw}=u+<u,v>u\\
&\implies&\frac{u+v+\gamma_{uvw}}{1+<u,v>}=u
\end{eqnarray*}

As for (3) and (4), these are also clear from the above discussion.
\end{proof}

Looking now at Proposition 15.19 from an abstract, mathematical perspective, there are still many things missing from there, which can be summarized as follows:

\begin{question}
Can we fine-tune the Einstein speed summation in $3{\rm D}$ into
$$u+_ev=\frac{u+v+\lambda\cdot u\times(u\times v)}{1+<u,v>}$$
with $\lambda\in\mathbb R$, chosen such that $||u||=1\implies\lambda=1$, as to have:
\begin{enumerate}
\item $||u||,||v||<1\implies||u+_ev||<1$.

\item $||v||=1\implies||u+_ev||=1$.
\end{enumerate}
\end{question}

Obviously, as simplest answer, $\lambda$ must be some well-chosen function of $||u||$, or rather of $||u||^2$, because it is always better to use square norms, when possible. But then, with this idea in mind, after a few computations we are led to the following solution:
$$\lambda=\frac{1}{1+\sqrt{1-||u||^2}}$$

Summarizing, final correction done, and with this being the end of mathematics, we did a nice job, and we can now formulate our findings as a theorem, as follows:

\begin{theorem}
When defining the Einstein speed summation in $3{\rm D}$ as 
$$u+_ev=\frac{1}{1+<u,v>}\left(u+v+\frac{u\times(u\times v)}{1+\sqrt{1-||u||^2}}\right)$$
in $c=1$ units, the following happen:
\begin{enumerate}
\item When $u\sim v$, we recover the previous $1{\rm D}$ formula.

\item We have $||u||,||v||<1\implies||u+_ev||<1$.

\item When $||u||=1$, we have $u+_ev=u$.

\item When $||v||=1$, we have $||u+_ev||=1$.

\item However, $||v||=1$ does not imply $u+_ev=v$.

\item Also, the formula $u+_ev=v+_eu$ fails.
\end{enumerate}
\end{theorem}

\begin{proof}
This follows from the above discussion, as follows:

\medskip

(1) This is something that we know from Proposition 15.19.

\medskip

(2) In order to simplify notation, let us set $\delta=\sqrt{1-||u||^2}$, which is the inverse of the quantity $\gamma=1/\sqrt{1-||u||^2}$. With this convention, we have:
\begin{eqnarray*}
u+_ev
&=&\frac{1}{1+<u,v>}\left(u+v+\frac{<u,v>u-||u||^2v}{1+\delta}\right)\\
&=&\frac{(1+\delta+<u,v>)u+(1+\delta-||u||^2)v}{(1+<u,v>)(1+\delta)}
\end{eqnarray*}

Taking the squared norm and computing gives the following formula:
$$||u+_ev||^2
=\frac{(1+\delta)^2||u+v||^2+(||u||^2-2(1+\delta))(||u||^2||v||^2-<u,v>^2)}{(1+<u,v>)^2(1+\delta)^2}$$

Now this formula can be further processed by using $\delta=\sqrt{1-||u||^2}$, and we get:
$$||u+_ev||^2=\frac{||u+v||^2-||u||^2||v||^2+<u,v>^2}{(1+<u,v>)^2}$$

But this type of formula is exactly what we need, for what we want to do. Indeed, by assuming $||u||,||v||<1$, we have the following estimate:
\begin{eqnarray*}
||u+_ev||^2<1
&\iff&||u+v||^2-||u||^2||v||^2+<u,v>^2<(1+<u,v>)^2\\
&\iff&||u+v||^2-||u||^2||v||^2<1+2<u,v>\\
&\iff&||u||^2+||v||^2-||u||^2||v||^2<1\\
&\iff&(1-||u||^2)(1-||v||^2)>0
\end{eqnarray*}

Thus, we are led to the conclusion in the statement.

\medskip

(3) This is something that we know from Proposition 15.19.

\medskip

(4) This comes from the squared norm formula established in the proof of (2) above, because when assuming $||v||=1$, we obtain:
\begin{eqnarray*}
||u+_ev||^2
&=&\frac{||u+v||^2-||u||^2+<u,v>^2}{(1+<u,v>)^2}\\
&=&\frac{||u||^2+1+2<u,v>-||u||^2+<u,v>^2}{(1+<u,v>)^2}\\
&=&\frac{1+2<u,v>+<u,v>^2}{(1+<u,v>)^2}\\
&=&1
\end{eqnarray*}

(5) This is clear, from the obvious lack of symmetry of our formula.

\medskip

(6) This is again clear, from the obvious lack of symmetry of our formula.
\end{proof}

Still with me I hope, after all these computations. And the good news is that, we did indeed some first-class work, because the formula in Theorem 15.21 is the good one, confirmed by experimental physics. For more on all this, and for the continuation of the story, special and then general relativity, you have the book of Einstein \cite{ein}.

\section*{15c. Charges and flux}

Time for electricity. In analogy with the usual study of gravity, let us start with:

\index{electric field}
\index{charge}

\begin{definition}
Given charges $q_1,\ldots,q_k\in\mathbb R$ located at positions $x_1,\ldots,x_k\in\mathbb R^3$, we define their electric field to be the vector function
$$E(x)=K\sum_i\frac{q_i(x-x_i)}{||x-x_i||^3}$$
so that their force applied to a charge $Q\in\mathbb R$ positioned at $x\in\mathbb R^3$ is given by $F=QE$.
\end{definition}

More generally, we will be interested in electric fields of various non-discrete configurations of charges, such as charged curves, surfaces and solid bodies. Indeed, things like wires or metal sheets or solid bodies coming in all sorts of shapes, tailored for their purpose, play a key role, so this extension is essential. So, let us go ahead with:

\begin{definition}
The electric field of a charge configuration $L\subset\mathbb R^3$, with charge density function $\rho:L\to\mathbb R$, is the vector function
$$E(x)=K\int_L\frac{\rho(z)(x-z)}{||x-z||^3}\,dz$$
so that the force of $L$ applied to a charge $Q$ positioned at $x$ is given by $F=QE$.
\end{definition}

With the above definitions in hand, it is most convenient now to forget about the charges, and focus on the study of the corresponding electric fields $E$. 

\bigskip

These fields are by definition vector functions $E:\mathbb R^3\to\mathbb R^3$, with the convention that they take $\pm\infty$ values at the places where the charges are located, and intuitively, are best represented by their field lines, which are constructed as follows:

\index{field lines}

\begin{definition}
The field lines of an electric field $E:\mathbb R^3\to\mathbb R^3$ are the oriented curves $\gamma\subset\mathbb R^3$ pointing at every point $x\in\mathbb R^3$ at the direction of the field, $E(x)\in\mathbb R^3$.
\end{definition}

As a basic example here, for one charge the field lines are the half-lines emanating from its position, oriented according to the sign of the charge:
$$\begin{matrix}
\nwarrow&\uparrow&\nearrow\\
\leftarrow&\oplus&\rightarrow\\
\swarrow&\downarrow&\searrow
\end{matrix}\qquad\qquad\qquad\qquad
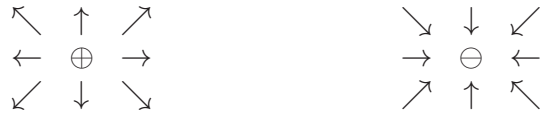
\begin{matrix}
\searrow&\downarrow&\swarrow\\
\rightarrow&\ominus&\leftarrow\\
\nearrow&\uparrow&\nwarrow
\end{matrix}$$

For two charges now, if these are of opposite signs, $+$ and $-$, you get a picture that you are very familiar with, namely that of the field lines of a bar magnet:
$$\begin{matrix}
\nearrow&\ \ \nearrow&\rightarrow&\rightarrow&\rightarrow&\rightarrow&\searrow\ \ \ &\!\!\!\searrow\\
\nwarrow&\!\!\!\uparrow&\nearrow&\rightarrow&\rightarrow&\searrow&\downarrow&\!\!\!\swarrow\\
\leftarrow&\!\!\!\oplus&\rightarrow&\rightarrow&\rightarrow&\rightarrow&\ominus&\!\!\!\leftarrow\\
\swarrow&\!\!\!\downarrow&\searrow&\rightarrow&\rightarrow&\nearrow&\uparrow&\!\!\!\nwarrow\\
\searrow&\ \ \searrow&\rightarrow&\rightarrow&\rightarrow&\rightarrow&\nearrow\ \ \ &\!\!\!\nearrow
\end{matrix}$$

If the charges are $+,+$ or $-,-$, you get something of similar type, but repulsive this time, with the field lines emanating from the charges being no longer shared:
$$\begin{matrix}
\leftarrow\ &\!\!\!\!\!\nwarrow&\nwarrow&&&\nearrow&\ \ \ \nearrow&\rightarrow\\
&\uparrow&\nearrow&&&\nwarrow&\uparrow&\\
\leftarrow&\oplus&\ &\ &\ \ \ \ \ \ \ &\ &\oplus&\rightarrow\\
&\downarrow&\searrow&&&\swarrow&\downarrow&\\
\leftarrow\ &\!\!\!\!\!\swarrow&\swarrow&&&\searrow&\ \ \ \searrow&\rightarrow
\end{matrix}$$

These pictures, and notably the last one, with $+,+$ charges, are quite interesting, because the repulsion situation does not appear in the context of gravity. Thus, we can only expect our geometry here to be far more complicated than that of gravity.

\bigskip

The field lines, as constructed in Definition 15.24, obviously do not encapsulate the whole information about the field, with the direction of each vector $E(x)\in\mathbb R^3$ being there, but with the magnitude $||E(x)||\geq0$ of this vector missing. However, say when drawing, when picking up uniformly radially spaced field lines around each charge, and with the number of these lines proportional to the magnitude of the charge, and then completing the picture, the density of the field lines around each point $x\in\mathbb R$ will give you then the magnitude $||E(x)||\geq0$ of the field there, up to a scalar.

\bigskip

Let us summarize these observations as follows:

\begin{proposition}
Given an electric field $E:\mathbb R^3\to\mathbb R^3$, the knowledge of its field lines is the same as the knowledge of the composition
$$nE:\mathbb R^3\to\mathbb R^3\to S$$
where $S\subset\mathbb R^3$ is the unit sphere, and $n:\mathbb R^3\to S$ is the rescaling map, namely:
$$n(x)=\frac{x}{||x||}$$
However, in practice, when the field lines are accurately drawn, the density of the field lines gives you the magnitude of the field, up to a scalar.
\end{proposition}

\begin{proof}
We have two assertions here, the idea being as follows:

\medskip

(1) The first assertion is clear from definitions, with of course our usual convention that the electric field and its problematics take place outside the locations of the charges, which makes everything in the statement to be indeed well-defined. 

\medskip

(2) Regarding now the last assertion, which is of course a bit informal, this follows from the above discussion. It is possible to be a bit more mathematical here, with a definition, formula and everything, but we will not need this, in what follows.
\end{proof}

Let us introduce now a key definition, as follows:

\index{flux}

\begin{definition}
The flux of an electric field $E:\mathbb R^3\to\mathbb R^3$ through a surface $S\subset\mathbb R^3$, assumed to be oriented, is the quantity
$$\Phi_E(S)=\int_S<E(x),n(x)>dx$$
with $n(x)$ being unit vectors orthogonal to $S$, following the orientation of $S$. Intuitively, the flux measures the signed number of field lines crossing $S$.
\end{definition}

Here by orientation of $S$ we mean precisely the choice of unit vectors $n(x)$ as above, orthogonal to $S$, which must vary continuously with $x$. For instance a sphere has two possible orientations, one with all these vectors $n(x)$ pointing inside, and one with all these vectors $n(x)$ pointing outside. More generally, any surface has locally two possible orientations, so if it is connected, it has two possible orientations. In what follows the convention is that the closed surfaces are oriented with each $n(x)$ pointing outside.

\bigskip

As a first illustration, let us do a basic computation, as follows:

\index{point charge}

\begin{proposition}
For a point charge $q\in\mathbb R$ at the center of a sphere $S$,
$$\Phi_E(S)=\frac{q}{\varepsilon_0}$$
where the constant is $\varepsilon_0=1/(4\pi K)$, independently of the radius of $S$.
\end{proposition}

\begin{proof}
Assuming that $S$ has radius $r$, we have the following computation:
\begin{eqnarray*}
\Phi_E(S)
&=&\int_S<E(x),n(x)>dx\\
&=&\int_S\left<\frac{Kqx}{r^3},\frac{x}{r}\right>dx\\
&=&\int_S\frac{Kq}{r^2}\,dx\\
&=&\frac{Kq}{r^2}\times 4\pi r^2\\
&=&4\pi Kq
\end{eqnarray*} 

Thus with $\varepsilon_0=1/(4\pi K)$ as above, we obtain the result.
\end{proof}

More generally now, we have the following result:

\index{flux}
\index{Jacobian}

\begin{theorem}
The flux of a field $E$ through a sphere $S$ is given by
$$\Phi_E(S)=\frac{Q_{enc}}{\varepsilon_0}$$
where $Q_{enc}$ is the total charge enclosed by $S$, and $\varepsilon_0=1/(4\pi K)$.
\end{theorem}

\begin{proof}
This can be done in several steps, as follows:

\medskip

(1) Before jumping into computations, let us do some manipulations. First, by discretizing the problem, we can assume that we are dealing with a system of point charges. Moreover, by additivity, we can assume that we are dealing with a single charge. And if we denote by $q\in\mathbb R$ this charge, located at $v\in\mathbb R^3$, we want to prove that we have the following formula, where $B\subset\mathbb R^3$ denotes the ball enclosed by $S$:
$$\Phi_E(S)=\frac{q}{\varepsilon_0}\,\delta_{v\in B}$$

(2) By linearity we can assume that we are dealing with the unit sphere $S$. Moreover, by rotating we can assume that our charge $q$ lies on the $Ox$ axis, that is, that we have $v=(r,0,0)$ with $r\geq0$, $r\neq1$. The formula that we want to prove becomes:
$$\Phi_E(S)=\frac{q}{\varepsilon_0}\,\delta_{r<1}$$

(3) Let us start now the computation. With $u=(x,y,z)$, we have: 
\begin{eqnarray*}
\Phi_E(S)
&=&\int_S<E(u),u>du\\
&=&\int_S\left<\frac{Kq(u-v)}{||u-v||^3},u\right>du\\
&=&Kq\int_S\frac{<u-v,u>}{||u-v||^3}\,du\\
&=&Kq\int_S\frac{1-<v,u>}{||u-v||^3}\,du\\
&=&Kq\int_S\frac{1-rx}{(1-2xr+r^2)^{3/2}}\,du
\end{eqnarray*}

(4) In order to compute the above integral, we will use spherical coordinates for the unit sphere $S$, which are as follows, with $s\in[0,\pi]$ and $t\in[0,2\pi]$:
$$\begin{cases}
x\!\!\!&=\ \cos s\\
y\!\!\!&=\ \sin s\cos t\\
z\!\!\!&=\ \sin s\sin t
\end{cases}$$

We recall that the corresponding Jacobian, computed before, is given by:
$$J=\sin s$$

(5) With the above change of coordinates, our integral from (3) becomes:
\begin{eqnarray*}
\Phi_E(S)
&=&Kq\int_S\frac{1-rx}{(1-2xr+r^2)^{3/2}}\,du\\
&=&Kq\int_0^{2\pi}\int_0^\pi\frac{1-r\cos s}{(1-2r\cos s+r^2)^{3/2}}\cdot\sin s\,ds\,dt\\
&=&2\pi Kq\int_0^\pi\frac{(1-r\cos s)\sin s}{(1-2r\cos s+r^2)^{3/2}}\,ds\\
&=&\frac{q}{2\varepsilon_0}\int_0^\pi\frac{(1-r\cos s)\sin s}{(1-2r\cos s+r^2)^{3/2}}\,ds
\end{eqnarray*}

(6) The point now is that the integral on the right can be computed with the change of variables $x=\cos s$. Indeed, we have $dx=-\sin s\,ds$, and we obtain:
\begin{eqnarray*}
\int_0^\pi\frac{(1-r\cos s)\sin s}{(1-2r\cos s+r^2)^{3/2}}\,ds
&=&\int_{-1}^1\frac{1-rx}{(1-2rx+r^2)^{3/2}}\,dx\\
&=&\left[\frac{x-r}{\sqrt{1-2rx+r^2}}\right]_{-1}^1\\
&=&\frac{1-r}{\sqrt{1-2r+r^2}}-\frac{-1-r}{\sqrt{1+2r+r^2}}\\
&=&\frac{1-r}{|1-r|}+1\\
&=&2\delta_{r<1}
\end{eqnarray*}

Thus, we are led to the formula in the statement.
\end{proof}

More generally now, we have the following key result, due to Gauss:

\index{Gauss law}
\index{flux}

\begin{theorem}[Gauss law]
The flux of a field $E$ through a surface $S$ is given by
$$\Phi_E(S)=\frac{Q_{enc}}{\varepsilon_0}$$
where $Q_{enc}$ is the total charge enclosed by $S$, and $\varepsilon_0=1/(4\pi K)$.
\end{theorem}

\begin{proof}
This basically follows from Theorem 15.28, or even from Proposition 15.27, by adding to the results there a number of new ingredients, as follows:

\medskip

(1) Our first claim is that given a closed surface $S$, with no charges inside, the flux through it of any choice of external charges vanishes:
$$\Phi_E(S)=0$$

This claim is indeed supported by the intuitive interpretation of the flux, as corresponding to the signed number of field lines crossing $S$. Indeed, any field line entering as $+$ must exit somewhere as $-$, and vice versa, so when summing we get $0$.

\medskip

(2) In practice now, in order to prove this rigorously, there are several ways. A standard argument, which is quite elementary, is the one used by Feynman in \cite{fe2}, based on the fact that, due to $F\sim 1/d^2$, local deformations of $S$ will leave invariant the flux, and so in the end we are left with a rotationally invariant surface, where the result is clear.

\medskip

(3) The point now is that, with this and Proposition 15.27 in hand, we can finish by using a standard math trick. Let us assume indeed, by discretizing, that our system of charges is discrete, consisting of enclosed charges $q_1,\ldots,q_k\in\mathbb R$, and an exterior total charge $Q_{ext}$. We can surround each of $q_1,\ldots,q_k$ by small disjoint spheres $U_1,\ldots,U_k$, chosen such that their interiors do not touch $S$, and we have:
\begin{eqnarray*}
\Phi_E(S)
&=&\Phi_E(S-\cup U_i)+\Phi_E(\cup U_i)\\
&=&0+\Phi_E(\cup U_i)\\
&=&\sum_i\Phi_E(U_i)\\
&=&\sum_i\frac{q_i}{\varepsilon_0}\\
&=&\frac{Q_{enc}}{\varepsilon_0}
\end{eqnarray*}

(4) To be more precise, in the above the union $\cup U_i$ is a usual disjoint union, and the flux is of course additive over components. As for the difference $S-\cup U_i$, this is by definition the disjoint union of $S$ with the disjoint union $\cup(-U_i)$, with each $-U_i$ standing for $U_i$ with orientation reversed, and since this difference has no enclosed charges, the flux through it vanishes by (2). Finally, the end makes use of Proposition 15.27.
\end{proof}

So long for the Gauss law. We will be back to it in a moment, with a new, better proof, by using some advanced 3D calculus, that we will have to learn first. 

\section*{15d. Gauss, Green, Stokes}

Getting back now to calculus tools, what was missing from our picture was the higher dimensional analogue of the fundamental theorem of calculus, and more generally of the partial integration formula. In 2 dimensions, to start with, we have:

\index{Green formula}
\index{Gauss formula}

\begin{theorem}[Green]
Given a plane curve $C\subset\mathbb R^2$, we have the formula
$$\int_CPdx+Qdy=\int_D\left(\frac{dQ}{dx}-\frac{dP}{dy}\right)dxdy$$
where $D\subset\mathbb R^2$ is the domain enclosed by $C$.
\end{theorem}

\begin{proof}
Assume indeed that we have a plane curve $C\subset\mathbb R^2$, without self-intersections, which is piecewise $C^1$, and is assumed to be counterclockwise oriented. In order to prove the formula regarding $P$, we can parametrize the enclosed domain $D$ as follows:
$$D=\left\{(x,y)\Big|a\leq x\leq b,f(x)\leq y\leq g(x)\right\}$$

We have then the following computation, which gives the result for $P$:
\begin{eqnarray*}
\int_D-\frac{dP}{dy}\,dxdy
&=&\int_a^b\left(\int_{f(x)}^{g(x)}-\frac{dP}{dy}(x,y)dy\right)dx\\
&=&\int_a^bP(x,f(x))-P(x,g(x))dx\\
&=&\int_a^bP(x,f(x))dx+\int_b^aP(x,g(x))dx\\
&=&\int_CPdx
\end{eqnarray*}

As for the result for the $Q$ term, the computation here is similar.
\end{proof}

Getting now to 3D, let us begin with a standard definition, as follows:

\index{gradient}
\index{divergence}
\index{curl}

\begin{definition}
Given a function $f:\mathbb R^3\to\mathbb R$, its usual derivative $f'(u)\in\mathbb R^3$ can be written as $f'(u)=\nabla f(u)$, where the gradient operator $\nabla$ is given by:
$$\nabla=\begin{pmatrix}
\frac{d}{dx}\\
\frac{d}{dy}\\
\frac{d}{dz}
\end{pmatrix}$$
By using $\nabla$, we can talk about the divergence of a function $\varphi:\mathbb R^3\to\mathbb R^3$, as being
$$<\nabla,\varphi>=
\left<\begin{pmatrix}
\frac{d}{dx}\\
\frac{d}{dy}\\
\frac{d}{dz}
\end{pmatrix},\begin{pmatrix}
\varphi_x\\
\varphi_y\\
\varphi_z
\end{pmatrix}\right>=
\frac{d\varphi_x}{dx}+\frac{d\varphi_y}{dy}+\frac{d\varphi_z}{dz}$$
as well as about the curl of the same function $\varphi:\mathbb R^3\to\mathbb R^3$, as being
$$\nabla\times\varphi=
\begin{vmatrix}
u_x&\frac{d}{dx}&\varphi_x\\
u_y&\frac{d}{dy}&\varphi_y\\
u_z&\frac{d}{dz}&\varphi_z
\end{vmatrix}
=\begin{pmatrix}
\frac{d\varphi_z}{dy}-\frac{d\varphi_y}{dz}\\
\frac{d\varphi_x}{dz}-\frac{d\varphi_z}{dx}\\
\frac{d\varphi_y}{dx}-\frac{d\varphi_x}{dy}
\end{pmatrix}
$$
where $u_x,u_y,u_z$ are the unit vectors along the coordinate directions $x,y,z$.
\end{definition}

All this might seem a bit abstract, but is in fact very intuitive. The gradient $\nabla f$ points in the direction of the maximal increase of $f$, with $|\nabla f|$ giving you the rate of increase of $f$, in that direction. As for the divergence and curl, these measure the divergence and curl of the vectors $\varphi(u+v)$ around a given point $u\in\mathbb R^3$, in a usual, real-life sense.

\bigskip

As a continuation now of our calculus work, we have the following result:

\index{Stokes formula}

\begin{theorem}[Stokes]
Given a smooth oriented surface $S\subset\mathbb R^3$, with boundary $C\subset\mathbb R^3$, and a vector field $F$, we have the following formula:
$$\int_S<(\nabla\times F)(x),n(x)>dx=\int_C<F(x),dx>$$
In other words, the line integral of a vector field $F$ over a loop $C$ equals the surface integral of the curl of the vector field $\varphi$, over the enclosed surface $S$. 
\end{theorem}

\begin{proof}
This basically follows from the Green theorem, the idea being as follows:

\medskip

(1) Let us first parametrize our surface $S$, and its boundary $C$. We can assume that we are in the situation where we have a closed oriented curve $\gamma:[a,b]\to\mathbb R^2$, with interior $D\subset\mathbb R^2$, and where the surface appears as $S=\psi(D)$, with $\psi:D\to\mathbb R^3$. In this case, the function $\delta=\psi\circ\gamma$ parametrizes the boundary of our surface, $C=\delta[a,b]$.

\medskip

(2) Let us first look at the integral on the right in the statement. We have the following formula, with $J_y(\psi)$ standing for the Jacobian of $\psi$ at the point $y=\gamma(t)$:
\begin{eqnarray*}
\int_C<F(x),dx>
&=&\int_\gamma<F(\psi(\gamma)),d\psi(\gamma)>\\
&=&\int_\gamma<F(\psi(y)),J_y(\psi)d\gamma>
\end{eqnarray*}

In order to further process this formula, let us introduce the following function:
$$P(u,v)=\left<F(\psi(u,v)),\frac{d\psi}{du}(u,v)\right>e_u+\left<F(\psi(u,v)),\frac{d\psi}{dv}(u,v)\right>e_v$$

In terms of this function, we have the following formula for our line integral:
$$\int_C<F(x),dx>=\int_\gamma<P(y),dy>$$

(3) In order to compute now the other integral in the statement, we first have:
\begin{eqnarray*}
&&\frac{dP_v}{du}-\frac{dP_u}{dv}\\
&=&\left<\frac{d(F\psi)}{du},\frac{d\psi}{dv}\right>+\left<F\psi,\frac{d^2\psi}{dudv}\right>
-\left<\frac{d(F\psi)}{dv},\frac{d\psi}{du}\right>-\left<F\psi,\frac{d^2\psi}{dvdu}\right>\\
&=&\left<\frac{d(F\psi)}{du},\frac{d\psi}{dv}\right>-\left<\frac{d(F\psi)}{dv},\frac{d\psi}{du}\right>\\
&=&\left<\frac{d\psi}{dv},\left(J_{\psi(u,v)}F-(J_{\psi(u,v)}F)^t\right)\frac{d\psi}{du}\right>\\
&=&\left<\frac{d\psi}{dv},(\nabla\times F)\times\frac{d\psi}{du}\right>\\
&=&\left<\nabla\times F,\frac{d\psi}{du}\times\frac{d\psi}{dv}\right>
\end{eqnarray*}

We conclude that the integral on the left in the statement is given by:
\begin{eqnarray*}
&&\int_S<(\nabla\times F)(x),n(x)>dx\\
&=&\int_D\left<(\nabla\times F)(\psi(u,v)),\frac{d\psi}{du}(u,v)\times\frac{d\psi}{dv}(u,v)\right>dudv\\
&=&\int_D\left(\frac{dP_v}{du}-\frac{dP_u}{dv}\right)dudv
\end{eqnarray*}

(4) But with this, we are done, because the integrals computed in (2) and (3) are indeed equal, due to the Green theorem. Thus, the Stokes formula holds indeed.
\end{proof}

As a conclusion to what we did so far, we have the following statement:

\begin{theorem}
The following results hold, in $3$ dimensions:
\begin{enumerate}
\item Fundamental theorem for gradients, namely
$$\int_a^b<\nabla f,dx>=f(b)-f(a)$$

\item Fundamental theorem for divergences, or Gauss or Green formula,
$$\int_B<\nabla,\varphi>=\int_S<\varphi(x),n(x)>dx$$

\item Fundamental theorem for curls, or Stokes formula,
$$\int_A<(\nabla\times\varphi)(x),n(x)>dx=\int_P<\varphi(x),dx>$$
\end{enumerate}
where $S$ is the boundary of the body $B$, and $P$ is the boundary of the area $A$.
\end{theorem}

\begin{proof}
This is a mixture of trivial and non-trivial results, with (1) being something that we know well in 1D, namely fundamental theorem of calculus, and the general, $N$-dimensional formula following from that, and with (2) and (3) being established above.
\end{proof}

Getting back now to electrostatics, as a main application of the above, we have the following new point of view on the Gauss formula, which is more conceptual:

\index{Gauss formula}

\begin{theorem}[Gauss]
Given an electric potential $E$, its divergence is given by
$$<\nabla,E>=\frac{\rho}{\varepsilon_0}$$
where $\rho$ denotes as usual the charge distribution. Also, we have
$$\nabla\times E=0$$
meaning that the curl of $E$ vanishes.
\end{theorem}

\begin{proof}
We have several assertions here, the idea being as follows:

\medskip

(1) The first formula, called Gauss law in differential form, follows from:
\begin{eqnarray*}
\int_B<\nabla,E>
&=&\int_S<E(x),n(x)>dx\\
&=&\Phi_E(S)\\
&=&\frac{Q_{enc}}{\varepsilon_0}\\
&=&\int_B\frac{\rho}{\varepsilon_0}
\end{eqnarray*}

Now since this must hold for any $B$, this gives the formula in the statement.

\medskip

(2) As a side remark, the Gauss law in differential form can be established as well directly, with the computation, involving a Dirac mass, being as follows:
\begin{eqnarray*}
<\nabla,E>(x)
&=&\left<\nabla,K\int_{\mathbb R^3}\frac{\rho(z)(x-z)}{||x-z||^3}\,dz\right>\\
&=&K\int_{\mathbb R^3}\left<\nabla,\frac{x-z}{||x-z||^3}\right>\,\rho(z)\,dz\\
&=&K\int_{\mathbb R^3} 4\pi\delta_x\cdot\rho(z)dz\\
&=&4\pi K\int_{\mathbb R^3}\delta_x\,\rho(z)dz\\
&=&\frac{\rho(x)}{\varepsilon_0}
\end{eqnarray*}

And with this in hand, we have via (1) a new proof of the usual Gauss law.

\medskip

(3) Regarding the curl, by discretizing and linearity we can assume that we are dealing with a single charge $q$, positioned at $0$. We have, by using spherical coordinates $r,s,t$:
\begin{eqnarray*}
\int_a^b<E(x),dx>
&=&\int_a^b\left<\frac{Kqx}{||x||^3},dx\right>\\
&=&\int_a^b\left<\frac{Kq}{r^2}\cdot\frac{x}{||x||},dx\right>\\
&=&\int_a^b\frac{Kq}{r^2}\,dr\\
&=&\left[-\frac{Kq}{r}\right]_a^b\\
&=&Kq\left(\frac{1}{r_a}-\frac{1}{r_b}\right)
\end{eqnarray*}

In particular the integral of $E$ over any closed loop vanishes, and by using now Stokes' theorem, we conclude that the curl of $E$ vanishes, as stated.

\medskip

(4) Finally, as a side remark, both the formula of the divergence and the vanishing of the curl are somewhat clear by looking at the field lines of $E$. However, as all the above mathematics shows, there is certainly something to be understood, in all this.
\end{proof}

So long for electrostatics, which provide a good motivation and illustration for our mathematics. When upgrading to electrodynamics, things become even more interesting, because our technology can be used in order to understand the Maxwell equations:
$$<\nabla,E>=\frac{\rho}{\varepsilon_0}\quad,\quad<\nabla,B>=0$$
$$\nabla\times E=-\dot{B}\quad,\quad\nabla\times B=\mu_0J+\mu_0\varepsilon_0\dot{E}$$

And we will leave some learning here, as a continuation of our previous discussion about these equations, from chapter 8, as a long, interesting, and must-do exercise.

\section*{15e. Exercises}

We have a tough physics chapter here, and as exercises, we have:

\begin{exercise}
Work out some numerics for the Coriolis force, on Earth.
\end{exercise}

\begin{exercise}
Learn about Einstein summation via Lorentz transformation.
\end{exercise}

\begin{exercise}
Prove the Gauss law by counting the flux lines. Can you?
\end{exercise}

\begin{exercise}
Find some further applications of Gauss, Green, Stokes.
\end{exercise}

As bonus exercise, learn full electrodynamics, say from Griffiths \cite{gr1}.

\chapter{Infinite dimensions}

\section*{16a. Operators, matrices}

Welcome to quantum mechanics, and in the hope that we will survive. We already talked in chapter 8 about the main idea of Heisenberg, namely using infinite matrices in order to axiomatize quantum mechanics, based on the following key fact:

\begin{fact}[Rydberg, Ritz]
The spectral lines of the hydrogen atom are given by the Rydberg formula, as follows, depending on integer parameters $n_1<n_2$:
$$\frac{1}{\lambda_{n_1n_2}}=R\left(\frac{1}{n_1^2}-\frac{1}{n_2^2}\right)$$
These spectral lines combine according to the Ritz-Rydberg principle, as follows:
$$\frac{1}{\lambda_{n_1n_2}}+\frac{1}{\lambda_{n_2n_3}}=\frac{1}{\lambda_{n_1n_3}}$$
Similar formulae hold for other atoms, with suitable fine-tunings of the constant $R$.
\end{fact}

We refer to chapter 8 for the full story with all this, which is theory based on some key observations of Lyman, Balmer, Paschen, around 1890-1900. The point now is that the above combination principle reminds the multiplication formula $e_{n_1n_2}e_{n_2n_3}=e_{n_1n_3}$ for the elementary matrices $e_{ij}:e_j\to e_i$, which leads to the following principle:

\begin{principle}[Heisenberg]
Observables in quantum mechanics should be some sort of infinite matrices, generalizing  the Lyman, Balmer, Paschen lines of the hydrogen atom, and multiplying between them as the matrices do, as to produce further observables.
\end{principle}

All this is quite deep, and needs a number of comments, as follows:

\bigskip

(1) First of all, our matrices must be indeed infinite, because so are the series observed by Lyman, Balmer, Paschen, corresponding to $n_1=1,2,3$ in the Rydberg formula, and making it clear that the range of the second parameter $n_2>n_1$ is up to $\infty$. 

\bigskip

(2) Although this was not known to Ritz-Rydberg and Heisenberg, let us mention too that some later results of Brackett, Pfund, Humphreys and others, at $n_1=4,5,6,\ldots\,$, confirmed the fact that the range of the first parameter $n_1$ is up to $\infty$ too.

\bigskip

(3) As a more tricky comment now, going beyond what Principle 16.2 says, our infinite matrices must be in fact complex. This was something known to Heisenberg, and later Schr\"odinger came with proof that quantum mechanics naturally lives over $\mathbb C$.

\bigskip

(4) But all this leads us into some tricky mathematics, because the infinite matrices $A\in M_\infty(\mathbb C)$ do not act on the vectors $v\in\mathbb C^\infty$ just like that. For instance the all-one matrix $A_{ij}=1$ does not act on the all-one vector $v_i=1$, for obvious reasons. 

\bigskip

Summarizing, in order to get to some mathematical theory going, out of Principle 16.2, we must assume that our matrices $A\in M_\infty(\mathbb C)$ must be ``bounded'' in some sense. Or perhaps the vectors $v\in\mathbb C^\infty$ must be bounded. Or perhaps, both.

\bigskip

In order to fix all this, let us start with $\mathbb C^\infty$. We would like to replace it with its subspace $H=l^2(\mathbb N)$ consisting of vectors having finite norm, as for our computations to converge. This being said, taking a look at what Schr\"odinger was saying too, a bit later, why not including right away in our theory spaces like $H=L^2(\mathbb R^3)$ too, which are perhaps a bit more relevant than Heisenberg's $l^2(\mathbb N)$. We are led in this way into the following notion, that we already met in chapter 7, and that we intend to discuss now in detail:

\index{Hilbert space}
\index{scalar product}

\begin{definition}
A Hilbert space is a complex vector space $H$ with a scalar product $<x,y>$, which will be linear at left and antilinear at right,
$$<\lambda x,y>=\lambda<x,y>\quad,\quad <x,\lambda y>=\bar{\lambda}<x,y>$$
and which is complete with respect to corresponding norm
$$||x||=\sqrt{<x,x>}$$
in the sense that any sequence $\{x_n\}$ which is a Cauchy sequence, having the property $||x_n-x_m||\to0$ with $n,m\to\infty$, has a limit, $x_n\to x$.
\end{definition}

Here our convention for the scalar products, written $<x,y>$ and being linear at left, is one among others, often used by mathematicians, and we will just use this, in the lack of a physicist with an axe around. As further comments now on Definition 16.3, there is some mathematics encapsulated there, needing some discussion. First, we have:

\begin{theorem}
Given an index set $I$, which can be finite or not, the space of square-summable vectors having indices in $I$, namely
$$l^2(I)=\left\{(x_i)_{i\in I}\Big|\sum_i|x_i|^2<\infty\right\}$$
is a Hilbert space, with scalar product as follows:
$$<x,y>=\sum_ix_i\bar{y}_i$$
When $I$ is finite, $I=\{1,\ldots,N\}$, we obtain in this way the usual space $H=\mathbb C^N$.
\end{theorem}

\begin{proof}
We already met such things in chapter 7, but let us recall all this:

\medskip

(1) We know that $l^2(I)\subset\mathbb C^I$ is the space of vectors satisfying $||x||<\infty$. We want to prove that $l^2(I)$ is a vector space, that $<x,y>$ is a scalar product on it, that $l^2(I)$ is complete with respect to $||.||$, and finally that for $|I|<\infty$ we have $l^2(I)=\mathbb C^{|I|}$.

\medskip

(2) The last assertion, $l^2(I)=\mathbb C^{|I|}$ for $|I|<\infty$, is clear, because in this case the sums are finite, so the condition $||x||<\infty$ is automatic. So, we know at least one thing.

\medskip

(3) Regarding the rest, our claim here, which will more or less prove everything, is that for any two vectors $x,y\in l^2(I)$ we have the Cauchy-Schwarz inequality:
$$|<x,y>|\leq||x||\cdot||y||$$

But this follows from the positivity of the following degree 2 quantity, depending on a real variable $t\in\mathbb R$, and on a variable on the unit circle, $w\in\mathbb T$:
$$f(t)=||twx+y||^2$$

(4) Now with Cauchy-Schwarz proved, everything is straightforward. We first obtain, by raising to the square and expanding, that for any $x,y\in l^2(I)$ we have:
$$||x+y||\leq||x||+||y||$$

Thus $l^2(I)$ is indeed a vector space, the other vector space conditions being trivial.

\medskip

(5) Also, $<x,y>$ is surely a scalar product on this vector space, because all the conditions for a scalar product are trivially satisfied.

\medskip

(6) Finally, the fact that our space $l^2(I)$ is indeed complete with respect to its norm $||.||$ follows in the obvious way, the limit of a Cauchy sequence $\{x_n\}$ being the vector $y=(y_i)$ given by $y_i=\lim_{n\to\infty}x_{ni}$, with all the verifications here being trivial.
\end{proof}

Going now a bit abstract, we have, more generally, the following result, which shows that our formalism covers as well the Schr\"odinger spaces of type $L^2(\mathbb R^3)$:

\begin{theorem}
Given an arbitrary space $X$ with a positive measure $\mu$ on it, the space of square-summable complex functions on it, namely
$$L^2(X)=\left\{f:X\to\mathbb C\Big|\int_X|f(x)|^2\,d\mu(x)<\infty\right\}$$
is a Hilbert space, with scalar product as follows:
$$<f,g>=\int_Xf(x)\overline{g(x)}\,d\mu(x)$$
When $X=I$ is discrete, meaning that the measure $\mu$ on it is the counting measure, $\mu(\{x\})=1$ for any $x\in X$, we obtain in this way the previous spaces $l^2(I)$.
\end{theorem}

\begin{proof}
This is something routine, remake of Theorem 16.4, as follows:

\medskip

(1) The proof of the first, and main assertion is something perfectly similar to the proof of Theorem 16.4, by replacing everywhere the sums by integrals. 

\medskip

(2) With the remark that we forgot to say in the statement that the $L^2$ functions are by definition taken up to equality almost everywhere, $f=g$ when $||f-g||=0$.

\medskip

(3) As for the last assertion, when $\mu$ is the counting measure all our integrals here become usual sums, and so we recover in this way Theorem 16.4.
\end{proof}

As a third and last theorem about Hilbert spaces, that we will need, we have:

\index{orthonormal basis}
\index{Gram-Schmidt}
\index{separable space}

\begin{theorem}
Any Hilbert space $H$ has an orthonormal basis $\{e_i\}_{i\in I}$, which is by definition a set of vectors whose span is dense in $H$, and which satisfy
$$<e_i,e_j>=\delta_{ij}$$
with $\delta$ being a Kronecker symbol. The cardinality $|I|$ of the index set, which can be finite, countable, or worse, depends only on $H$, and is called dimension of $H$. We have
$$H\simeq l^2(I)$$
in the obvious way, mapping $\sum\lambda_ie_i\to(\lambda_i)$. The Hilbert spaces with $\dim H=|I|$ being countable, including $l^2(\mathbb N)$ and $L^2(\mathbb R)$, are all isomorphic, and are called separable.
\end{theorem}

\begin{proof}
We have many assertions here, the idea being as follows:

\medskip

(1) In finite dimensions an orthonormal basis $\{e_i\}_{i\in I}$ can be constructed by starting with any vector space basis $\{x_i\}_{i\in I}$, and using the Gram-Schmidt procedure. As for the other assertions, these are all clear, from basic linear algebra.

\medskip

(2) In general, the same method works, namely Gram-Schmidt, with a subtlety coming from the fact that the basis $\{e_i\}_{i\in I}$ will not span in general the whole $H$, but just a dense subspace of it, as it is in fact obvious by looking at the standard basis of $l^2(\mathbb N)$. 

\medskip

(3) And there is a second subtlety as well, coming from the fact that the recurrence procedure needed for Gram-Schmidt must be replaced by some sort of ``transfinite recurrence'', using scary tools from logic, and more specifically the Zorn lemma.

\medskip

(4) Finally, everything at the end is clear from definitions, except perhaps for the fact that $L^2(\mathbb R)$ is separable. But here we can argue that, since functions can be approximated by polynomials, we have a countable algebraic basis, namely $\{x^n\}_{n\in\mathbb N}$, called the Weierstrass basis, that we can orthogonalize afterwards by using Gram-Schmidt.
\end{proof}

Moving ahead, now that we know what our vector spaces are, we can talk about infinite matrices with respect to them. And the situation here is as follows:

\index{infinite matrix}
\index{linear operator}
\index{bounded operator}

\begin{theorem}
Given a Hilbert space $H$, consider the linear operators $T:H\to H$, and for each such operator define its norm by the following formula:
$$||T||=\sup_{||x||=1}||Tx||$$
The operators which are bounded, $||T||<\infty$, form then a complex algebra $B(H)$, which is complete with respect to $||.||$. When $H$ comes with a basis $\{e_i\}_{i\in I}$, we have
$$B(H)\subset M_I(\mathbb C)$$
with the correspondence $T\leftrightarrow M$ obtained via the usual linear algebra formulae, namely:
$$T(x)=Mx\quad,\quad M_{ij}=<Te_j,e_i>$$
In infinite dimensions, the above inclusion is not an equality.
\end{theorem}

\begin{proof}
This is something straightforward, the idea being as follows:

\medskip

(1) The fact that we have indeed an algebra, satisfying the product condition in the statement, follows from the following estimates, which are all elementary:
$$||S+T||\leq||S||+||T||\quad,\quad 
||\lambda T||=|\lambda|\cdot||T||\quad,\quad 
||ST||\leq||S||\cdot||T||$$

(2) Regarding now the completness assertion, if $\{T_n\}\subset B(H)$ is Cauchy then $\{T_nx\}$ is Cauchy for any $x\in H$, so we can define the limit $T=\lim_{n\to\infty}T_n$ by setting:
$$Tx=\lim_{n\to\infty}T_nx$$

Let us first check that the application $x\to Tx$ is linear. We have:
\begin{eqnarray*}
T(x+y)
&=&\lim_{n\to\infty}T_n(x+y)\\
&=&\lim_{n\to\infty}T_n(x)+T_n(y)\\
&=&\lim_{n\to\infty}T_n(x)+\lim_{n\to\infty}T_n(y)\\
&=&T(x)+T(y)
\end{eqnarray*}

Similarly, we have $T(\lambda x)=\lambda T(x)$, and we conclude that $T\in\mathcal L(H)$.

\medskip

(3) With this done, it remains to prove now that we have $T\in B(H)$, and that $T_n\to T$ in norm. For this purpose, observe that we have:
\begin{eqnarray*}
||T_n-T_m||\leq\varepsilon\ ,\ \forall n,m\geq N
&\implies&||T_nx-T_mx||\leq\varepsilon\ ,\ \forall||x||=1\ ,\ \forall n,m\geq N\\
&\implies&||T_nx-Tx||\leq\varepsilon\ ,\ \forall||x||=1\ ,\ \forall n\geq N\\
&\implies&||T_Nx-Tx||\leq\varepsilon\ ,\ \forall||x||=1\\
&\implies&||T_N-T||\leq\varepsilon
\end{eqnarray*}

But this gives both $T\in B(H)$, and $T_N\to T$ in norm, and we are done.

\medskip

(4) Regarding the embedding, the correspondence $T\to M$ in the statement is indeed linear, and its kernel is $\{0\}$, so we have indeed an embedding as follows, as claimed:
$$B(H)\subset M_I(\mathbb C)$$

In finite dimensions we have an isomorphism, because any $M\in M_N(\mathbb C)$ determines an operator $T:\mathbb C^N\to\mathbb C^N$, given by $<Te_j,e_i>=M_{ij}$. However, in infinite dimensions, we have matrices not producing operators, as for instance the all-one matrix.
\end{proof}

Finally, as a second and last result regarding the operators, we will need:

\index{adjoint operator}

\begin{theorem}
Each operator $T\in B(H)$ has an adjoint $T^*\in B(H)$, given by: 
$$<Tx,y>=<x,T^*y>$$
The operation $T\to T^*$ is antilinear, antimultiplicative, involutive, and satisfies:
$$||T||=||T^*||\quad,\quad ||TT^*||=||T||^2$$
When $H$ comes with a basis $\{e_i\}_{i\in I}$, the operation $T\to T^*$ corresponds to
$$(M^*)_{ij}=\overline{M}_{ji}$$ 
at the level of the associated matrices $M\in M_I(\mathbb C)$.
\end{theorem}

\begin{proof}
This is standard too, and can be proved in 3 steps, as follows:

\medskip

(1) The existence of the adjoint operator $T^*$, given by the formula in the statement, comes from the fact that the function $\varphi(x)=<Tx,y>$ being a linear map $H\to\mathbb C$, we must have a formula as follows, for a certain vector $T^*y\in H$:
$$\varphi(x)=<x,T^*y>$$

Moreover, since this vector is unique, $T^*$ is unique too, and we have as well:
$$(S+T)^*=S^*+T^*\quad,\quad
(\lambda T)^*=\bar{\lambda}T^*\quad,\quad 
(ST)^*=T^*S^*\quad,\quad 
(T^*)^*=T$$

Observe also that we have indeed $T^*\in B(H)$, because:
\begin{eqnarray*}
||T||
&=&\sup_{||x||=1}\sup_{||y||=1}<Tx,y>\\
&=&\sup_{||y||=1}\sup_{||x||=1}<x,T^*y>\\
&=&||T^*||
\end{eqnarray*}

(2) Regarding now $||TT^*||=||T||^2$, which is a key formula, observe that we have:
$$||TT^*||
\leq||T||\cdot||T^*||
=||T||^2$$

On the other hand, we have as well the following estimate:
\begin{eqnarray*}
||T||^2
&=&\sup_{||x||=1}|<Tx,Tx>|\\
&=&\sup_{||x||=1}|<x,T^*Tx>|\\
&\leq&||T^*T||
\end{eqnarray*}

By replacing $T\to T^*$ we obtain from this $||T||^2\leq||TT^*||$, as desired.

\medskip

(3) Finally, when $H$ comes with a basis, the formula $<Tx,y>=<x,T^*y>$ applied with $x=e_i$, $y=e_j$ translates into the formula $(M^*)_{ij}=\overline{M}_{ji}$, as desired.
\end{proof}

So long for Hilbert spaces and operators. For more, you can check my book \cite{ba2}.

\section*{16b. Schr\"odinger equation}

Time now to get into quantum mechanics. This is a bit of a nightmare, and I am afraid that I will have to ask again the cat. With him, things going better since last time, we sort of agree now on multivariable calculus, following a long philosophical debate, and my willingness to start sharing his diet, namely raw mice and birds. But let's see now in relation with this highly divisive topic which is quantum mechanics:

\begin{cat}
Schr\"odinger.
\end{cat}

Thanks cat, this is actually quite surprising, and I am particularly pleased to hear this. So, before getting back to what Heisenberg was saying, based on Lyman, Balmer, Paschen, namely developing some sort of ``matrix mechanics'', let us hear as well the point of view of Schr\"odinger, which came a few years later. His idea was to forget about exact things, and try to investigate the hydrogen atom statistically. Let us start with:

\begin{question}
In the context of the hydrogen atom, assuming that the proton is fixed, what is the probability density $\varphi_t(x)$ of the position of the electron $e$, at time $t$,
$$P_t(e\in V)=\int_V\varphi_t(x)dx$$
as function of an intial probability density $\varphi_0(x)$? Moreover, can the corresponding equation be solved, and will this prove the Bohr claims for hydrogen, statistically?
\end{question}

In order to get familiar with this question, let us first look at examples coming from classical mechanics. In the context of a particle whose position at time $t$ is given by $x_0+\gamma(t)$, the evolution of the probability density will be given by: 
$$\varphi_t(x)=\varphi_0(x)+\gamma(t)$$

However, such examples are somewhat trivial, of course not in relation with the computation of $\gamma$, usually a difficult question, but in relation with our questions, and do not apply to the electron. The point indeed is that, in what regards the electron, we have:

\begin{fact}
In respect with various simple interference experiments:
\begin{enumerate}
\item The electron is definitely not a particle in the usual sense.

\item But in most situations it behaves exactly like a wave.

\item But in other situations it behaves like a particle.
\end{enumerate}
\end{fact}

Getting back now to the Schr\"odinger question, all this suggests to use, as for the waves, an amplitude function $\psi_t(x)\in\mathbb C$, related to the density $\varphi_t(x)>0$ by the formula $\varphi_t(x)=|\psi_t(x)|^2$. Not that a big deal, you would say, because the two are related by simple formulae as follows, with $\theta_t(x)$ being an arbitrary phase function:
$$\varphi_t(x)=|\psi_t(x)|^2\quad,\quad \psi_t(x)=e^{i\theta_t(x)}\sqrt{\varphi_t(x)}$$

However, such manipulations can be crucial, raising for instance the possibility that the amplitude function satisfies some simple equation, while the density itself, maybe not. And this is what happens indeed. Schr\"odinger was led in this way to:

\index{Schr\"odinger equation}
\index{wave function}
\index{Planck constant}
\index{electron}

\begin{claim}[Schr\"odinger]
In the context of the hydrogen atom, the amplitude function of the electron $\psi=\psi_t(x)$  is subject to the Schr\"odinger equation
$$ih\dot{\psi}=-\frac{h^2}{2m}\Delta\psi+V\psi$$
$m$ being the mass, $h=h_0/2\pi$ the reduced Planck constant, and $V$ the Coulomb potential of the proton. The same holds for movements of the electron under any potential $V$.
\end{claim}

Observe the similarity with the wave equation $\ddot{f}=v^2\Delta f$, and with the heat equation $\dot{f}=\alpha\Delta f$ too. Many things can be said here. Following now Heisenberg and Schr\"odinger, and then especially Dirac, who did the axiomatization work, we have:

\index{observable}

\begin{definition}
In quantum mechanics the states of the system are vectors of a Hilbert space $H$, and the observables of the system are linear operators 
$$T:H\to H$$
which can be densely defined, and are taken self-adjoint, $T=T^*$. The average value of such an observable $T$, evaluated on a state $\xi\in H$, is given by:
$$<T>=<T\xi,\xi>$$
In the context of the Schr\"odinger mechanics of the hydrogen atom, the Hilbert space is the space $H=L^2(\mathbb R^3)$ where the wave function $\psi$ lives, and we have
$$<T>=\int_{\mathbb R^3}T(\psi)\cdot\bar{\psi}\,dx$$
which is called ``sandwiching'' formula, with the operators 
$$x\quad,\quad-\frac{ih}{m}\nabla\quad,\quad-ih\nabla\quad,\quad -\frac{h^2\Delta}{2m}\quad,\quad -\frac{h^2\Delta}{2m}+V$$
representing the position, speed, momentum, kinetic energy, and total energy.
\end{definition}

In other words, we are doing here two things. First, we are declaring by axiom that various ``sandwiching'' formulae found before by Heisenberg, involving the operators at the end, that we will not get into in detail here, hold true. And second, we are raising the possibility for other quantum mechanical systems, more complicated, to be described as well by the mathematics of the operators on a certain Hilbert space $H$, as above.

\bigskip

So, this was the story of early quantum mechanics, over-simplified as to fit here in a few pages. For more, you can check Feynman \cite{fe3} for foundations, and everything, including for some nice pictures and explanations regarding Fact 16.11. You have as well Griffiths \cite{gr2} or  Weinberg \cite{we2}, for further explanations on Definition 16.13, not to forget Dirac's original text \cite{dir}, and all this is discussed as well in my book \cite{ba3}.

\section*{16c. Spherical coordinates}

In order to solve now the hydrogen atom, the idea will be that of reformulating the Schr\"odinger equation in spherical coordinates. And for this purpose, we will need:

\index{Laplace operator}
\index{spherical coordinates}

\begin{theorem}
The Laplace operator in spherical coordinates is
$$\Delta=\frac{1}{r^2}\cdot\frac{d}{dr}\left(r^2\cdot\frac{d}{dr}\right)
+\frac{1}{r^2\sin s}\cdot\frac{d}{ds}\left(\sin s\cdot\frac{d}{ds}\right)
+\frac{1}{r^2\sin^2s}\cdot\frac{d^2}{dt^2}$$
with our standard conventions for these coordinates, in 3D.
\end{theorem}

\begin{proof}
There are several proofs here, a short, elementary one being as follows:

\medskip

(1) Let us first see how $\Delta$ behaves under a change of coordinates $\{x_i\}\to\{y_i\}$, in arbitrary $N$ dimensions. Our starting point is the chain rule for derivatives:
$$\frac{d}{dx_i}=\sum_j\frac{d}{dy_j}\cdot\frac{dy_j}{dx_i}$$

By using this rule, then Leibnitz for products, then again this rule, we obtain:
\begin{eqnarray*}
\frac{d^2f}{dx_i^2}
&=&\sum_j\frac{d}{dx_i}\left(\frac{df}{dy_j}\cdot\frac{dy_j}{dx_i}\right)\\
&=&\sum_j\frac{d}{dx_i}\left(\frac{df}{dy_j}\right)\cdot\frac{dy_j}{dx_i}+\frac{df}{dy_j}\cdot\frac{d}{dx_i}\left(\frac{dy_j}{dx_i}\right)\\
&=&\sum_j\left(\sum_k\frac{d}{dy_k}\cdot\frac{dy_k}{dx_i}\right)
\left(\frac{df}{dy_j}\right)\cdot\frac{dy_j}{dx_i}+\frac{df}{dy_j}\cdot\frac{d^2y_j}{dx_i^2}\\
&=&\sum_{jk}\frac{d^2f}{dy_kdy_j}\cdot\frac{dy_k}{dx_i}\cdot\frac{dy_j}{dx_i}+\sum_j\frac{df}{dy_j}\cdot\frac{d^2y_j}{dx_i^2}
\end{eqnarray*}

(2) Now by summing over $i$, we obtain the following formula, with $A$ being the derivative of $x\to y$, that is to say, the matrix of partial derivatives $dy_i/dx_j$:
\begin{eqnarray*}
\Delta f
&=&\sum_{ijk}\frac{d^2f}{dy_kdy_j}\cdot\frac{dy_k}{dx_i}\cdot\frac{dy_j}{dx_i}+\sum_{ij}\frac{df}{dy_j}\cdot\frac{d^2y_j}{dx_i^2}\\
&=&\sum_{ijk}A_{ki}A_{ji}\frac{d^2f}{dy_kdy_j}+\sum_{ij}\frac{d^2y_j}{dx_i^2}\cdot\frac{df}{dy_j}\\
&=&\sum_{jk}(AA^t)_{jk}\frac{d^2f}{dy_kdy_j}+\sum_j\Delta(y_j)\frac{df}{dy_j}
\end{eqnarray*}

(3) So, this will be the formula that we will need. Observe that this formula can be further compacted as follows, with all the notations being self-explanatory:
$$\Delta f=Tr(AA^tH_y(f))+<\Delta(y),\nabla_y(f)>$$

(4) Getting now to spherical coordinates, $(x,y,z)\to(r,s,t)$, the derivative of the inverse, obtained by differentiating $x,y,z$ with respect to $r,s,t$, is given by:
$$A^{-1}=\begin{pmatrix}
\cos s&-r\sin s&0\\
\sin s\cos t&r\cos s\cos t&-r\sin s\sin t\\
\sin s\sin t&r\cos s\sin t&r\sin s\cos t
\end{pmatrix}$$

The product $(A^{-1})^tA^{-1}$ of the transpose of this matrix with itself is then:
$$\begin{pmatrix}
\cos s&\sin s\cos t&\sin s\sin t\\
-r\sin s&r\cos s\cos t&r\cos s\sin t\\
0&-r\sin s\sin t&r\sin s\cos t
\end{pmatrix}
\begin{pmatrix}
\cos s&-r\sin s&0\\
\sin s\cos t&r\cos s\cos t&-r\sin s\sin t\\
\sin s\sin t&r\cos s\sin t&r\sin s\cos t
\end{pmatrix}$$

But everything simplifies here, and we have the following remarkable formula, which by the way is something very useful, worth to be memorized:
$$(A^{-1})^tA^{-1}=\begin{pmatrix}1&0&0\\ 0&r^2&0\\ 0&0&r^2\sin^2s\end{pmatrix}$$

Now by inverting, we obtain the following formula, in relation with the above:
$$AA^t=\begin{pmatrix}1&0&0\\ 0&1/r^2&0\\ 0&0&1/(r^2\sin^2s)\end{pmatrix}$$

(5) Let us compute now the Laplacian of $r,s,t$. We first have the following formula, that we will use many times in what follows, and is worth to be memorized:
\begin{eqnarray*}
\frac{dr}{dx}
&=&\frac{d}{dx}\sqrt{x^2+y^2+z^2}\\
&=&\frac{1}{2}\cdot\frac{2x}{\sqrt{x^2+y^2+z^2}}\\
&=&\frac{x}{r}
\end{eqnarray*}

Of course the same computation works for $y,z$ too, and we therefore have:
$$\frac{dr}{dx}=\frac{x}{r}\quad,\quad\frac{dr}{dy}=\frac{y}{r}\quad,\quad\frac{dr}{dz}=\frac{z}{r}$$

(6) By using the above formulae, twice, we can compute the Laplacian of $r$:
\begin{eqnarray*}
\Delta(r)
&=&\Delta\left(\sqrt{x^2+y^2+z^2}\right)\\
&=&\frac{d}{dx}\left(\frac{x}{r}\right)+\frac{d}{dy}\left(\frac{y}{r}\right)+\frac{d}{dz}\left(\frac{z}{r}\right)\\
&=&\frac{r^2-x^2}{r^3}+\frac{r^2-y^2}{r^3}+\frac{r^2-z^2}{r^3}\\
&=&\frac{2}{r}
\end{eqnarray*}

(7) In what regards now $s$, the computation here goes as follows:
\begin{eqnarray*}
\Delta(s)
&=&\Delta\left(\arccos\left(\frac{x}{r}\right)\right)\\
&=&\frac{d}{dx}\left(-\frac{\sqrt{r^2-x^2}}{r^2}\right)
+\frac{d}{dy}\left(\frac{xy}{r^2\sqrt{r^2-x^2}}\right)
+\frac{d}{dz}\left(\frac{xz}{r^2\sqrt{r^2-x^2}}\right)\\
&=&\frac{2x\sqrt{r^2-x^2}}{r^4}
+\frac{r^2(z^2-2y^2)+2x^2y^2}{r^4\sqrt{r^2-x^2}}
+\frac{r^2(y^2-2z^2)+2x^2z^2}{r^4\sqrt{r^2-x^2}}\\
&=&\frac{2x\sqrt{r^2-x^2}}{r^4}
+\frac{x(2x^2-r^2)}{r^4\sqrt{r^2-x^2}}\\
&=&\frac{x}{r^2\sqrt{r^2-x^2}}\\
&=&\frac{\cos s}{r^2\sin s}
\end{eqnarray*}

(8) Finally, in what regards $t$, the computation here goes as follows:
\begin{eqnarray*}
\Delta(t)
&=&\Delta\left(\arctan\left(\frac{z}{y}\right)\right)\\
&=&\frac{d}{dx}(0)+\frac{d}{dy}\left(-\frac{z}{y^2+z^2}\right)+\frac{d}{dz}\left(\frac{y}{y^2+z^2}\right)\\
&=&0-\frac{2yz}{(y^2+z^2)^2}+\frac{2yz}{(y^2+z^2)^2}\\
&=&0
\end{eqnarray*}

(9) We can now plug the data from (4) and (6,7,8) in the general formula that we found in (2) above, and we obtain in this way:
\begin{eqnarray*}
\Delta f
&=&\frac{d^2f}{dr^2}+\frac{1}{r^2}\cdot\frac{d^2f}{ds^2}+\frac{1}{r^2\sin^2s}\cdot\frac{d^2f}{dt^2}+\frac{2}{r}\cdot\frac{df}{dr}+\frac{\cos s}{r^2\sin s}\cdot\frac{df}{ds}\\
&=&\frac{2}{r}\cdot\frac{df}{dr}+\frac{d^2f}{dr^2}+\frac{\cos s}{r^2\sin s}\cdot\frac{df}{ds}+\frac{1}{r^2}\cdot\frac{d^2f}{ds^2}+\frac{1}{r^2\sin^2s}\cdot\frac{d^2f}{dt^2}\\
&=&\frac{1}{r^2}\cdot\frac{d}{dr}\left(r^2\cdot\frac{df}{dr}\right)
+\frac{1}{r^2\sin s}\cdot\frac{d}{ds}\left(\sin s\cdot\frac{df}{ds}\right)
+\frac{1}{r^2\sin^2s}\cdot\frac{d^2f}{dt^2}
\end{eqnarray*}

Thus, we are led to the formula in the statement.
\end{proof}

Still with me, I hope, and do not worry, one day you will have such computations for breakfast. We can now reformulate the Schr\"odinger equation in spherical coordinates, and separate the variables, which leads to a radial and angular equation, as follows:

\index{time-independent equation}
\index{separation of variables}
\index{radial equation}
\index{angular equation}

\begin{theorem}
The time-independent Schr\"odinger equation in spherical coordinates separates, for solutions of type $\phi=\rho(r)\alpha(s,t)$, into two equations, as follows,
$$\frac{d}{dr}\left(r^2\cdot\frac{d\rho}{dr}\right)
-\frac{2mr^2}{h^2}(V-E)\rho=K\rho$$
$$\sin s\cdot\frac{d}{ds}\left(\sin s\cdot\frac{d\alpha}{ds}\right)
+\frac{d^2\alpha}{dt^2}=-K\sin^2s\cdot\alpha$$
with $K$ being a constant, called radial equation, and angular equation.
\end{theorem}

\begin{proof}
By using the formula in Theorem 16.14, the time-independent Schr\"odinger equation reformulates in spherical coordinates as follows:
$$(V-E)\phi=\frac{h^2}{2m}\left[\frac{1}{r^2}\cdot\frac{d}{dr}\left(r^2\cdot\frac{d\phi}{dr}\right)
+\frac{1}{r^2\sin s}\cdot\frac{d}{ds}\left(\sin s\cdot\frac{d\phi}{ds}\right)
+\frac{1}{r^2\sin^2s}\cdot\frac{d^2\phi}{dt^2}\right]$$

Let us look now for separable solutions for this latter equation, consisting of a radial part and an angular part, as in the statement, namely:
$$\phi(r,s,t)=\rho(r)\alpha(s,t)$$

By plugging this function into our equation, we obtain:
$$(V-E)\rho\alpha=\frac{h^2}{2m}\left[\frac{\alpha}{r^2}\cdot\frac{d}{dr}\left(r^2\cdot\frac{d\rho}{dr}\right)
+\frac{\rho}{r^2\sin s}\cdot\frac{d}{ds}\left(\sin s\cdot\frac{d\alpha}{ds}\right)
+\frac{\rho}{r^2\sin^2s}\cdot\frac{d^2\alpha}{dt^2}\right]$$

In order to solve this equation, we will do two manipulations. First, by multiplying everything by $2mr^2/(h^2\rho\alpha)$, this equation takes the following more convenient form:
$$\frac{2mr^2}{h^2}(V-E)=\frac{1}{\rho}\cdot\frac{d}{dr}\left(r^2\cdot\frac{d\rho}{dr}\right)
+\frac{1}{\alpha\sin s}\cdot\frac{d}{ds}\left(\sin s\cdot\frac{d\alpha}{ds}\right)
+\frac{1}{\alpha\sin^2s}\cdot\frac{d^2\alpha}{dt^2}$$

Now observe that by moving the radial terms to the left, and the angular terms to the right, this latter equation can be written as follows:
$$\frac{2mr^2}{h^2}(V-E)-\frac{1}{\rho}\cdot\frac{d}{dr}\left(r^2\cdot\frac{d\rho}{dr}\right)
=\frac{1}{\alpha\sin^2s}\left[\sin s\cdot\frac{d}{ds}\left(\sin s\cdot\frac{d\alpha}{ds}\right)
+\frac{d^2\alpha}{dt^2}\right]$$

Since this latter equation is now separated between radial and angular variables, both sides must be equal to a certain constant $-K$, as follows:
$$\frac{2mr^2}{h^2}(V-E)-\frac{1}{\rho}\cdot\frac{d}{dr}\left(r^2\cdot\frac{d\rho}{dr}\right)
=-K$$
$$\frac{1}{\alpha\sin^2s}\left[\sin s\cdot\frac{d}{ds}\left(\sin s\cdot\frac{d\alpha}{ds}\right)
+\frac{d^2\alpha}{dt^2}\right]=-K$$

But this leads to the conclusion in the statement.
\end{proof}

Let us first study the angular equation. We first have the following result:

\index{polar equation}
\index{azimuthal equation}

\begin{proposition}
The angular equation that we found before, namely
$$\sin s\cdot\frac{d}{ds}\left(\sin s\cdot\frac{d\alpha}{ds}\right)
+\frac{d^2\alpha}{dt^2}=-K\sin^2s\cdot\alpha$$
separates, for solutions of type $\alpha=\sigma(s)\theta(t)$, into two equations, as follows,
$$\frac{1}{\theta}\cdot\frac{d^2\theta}{dt^2}=-m^2$$
$$\frac{\sin s}{\sigma}\cdot\frac{d}{ds}\left(\sin s\cdot\frac{d\sigma}{ds}\right)
+K\sin^2s=m^2$$
with $m$ being a constant, called azimuthal equation, and polar equation.
\end{proposition}

\begin{proof}
This is something elementary, the idea being as follows:

\medskip

(1) Let us first recall that $r\in[0,\infty)$ is the radius, $s\in[0,\pi]$ is the polar angle, and $t\in[0,2\pi]$ is the azimuthal angle. Be said in passing, there are several conventions and notations here, and the above ones, that we use here, come from the general ones in $N$ dimensions, because further coordinates can be easily added, in the obvious way.

\medskip

(2) Getting back now to our question, by plugging $\alpha=\sigma(s)\theta(t)$ into the angular equation, we obtain:
$$\sin s\cdot\theta\cdot\frac{d}{ds}\left(\sin s\cdot\frac{d\sigma}{ds}\right)
+\sigma\cdot\frac{d^2\theta}{dt^2}=-K\sin^2s\cdot\sigma\theta$$

By dividing everything by $\sigma\theta$, this equation can be written as follows:
$$-\frac{1}{\theta}\cdot\frac{d^2\theta}{dt^2}=\frac{\sin s}{\sigma}\cdot\frac{d}{ds}\left(\sin s\cdot\frac{d\sigma}{ds}\right)
+K\sin^2s$$

Since the variables are separated, we must have, for a certain constant $m$:
$$\frac{1}{\theta}\cdot\frac{d^2\theta}{dt^2}=-m^2$$
$$\frac{\sin s}{\sigma}\cdot\frac{d}{ds}\left(\sin s\cdot\frac{d\sigma}{ds}\right)
+K\sin^2s=m^2$$

Thus, we are led to the conclusion in the statement.
\end{proof}

Regarding the azimuthal equation, things here are quickly settled, as follows:

\begin{proposition}
The solutions of the azimuthal equation, namely
$$\frac{1}{\theta}\cdot\frac{d^2\theta}{dt^2}=-m^2$$
are the functions as follows, with $a,b\in\mathbb C$ being parameters,
$$\theta(t)=ae^{imt}+be^{-imt}$$
and with only the case $m\in\mathbb Z$ being acceptable, on physical grounds.
\end{proposition}

\begin{proof}
The first assertion is clear, because we have a second order equation, and two obvious solutions for it, $e^{\pm imt}$, and then their linear combinations, and that's all. Regarding the last assertion, the point here is that by using $\theta(t)=\theta(t+2\pi)$, which is a natural physical assumption on the wave function, we are led to $m\in\mathbb Z$, as stated.
\end{proof}

We are now about to solve the angular equation, with only the polar equation remaining to be studied. However, in practice, this polar equation is 10 times more difficult that everything what we did so far, so be patient. We first have:

\begin{proposition}
The polar equation that we found before, namely
$$\frac{\sin s}{\sigma}\cdot\frac{d}{ds}\left(\sin s\cdot\frac{d\sigma}{ds}\right)
+K\sin^2s=m^2$$
with $m\in\mathbb Z$, translates via $\sigma(s)=f(\cos s)$ into the following equation,
$$(1-x^2)f''(x)-2xf'(x)=\left(\frac{m^2}{1-x^2}-K\right)f(x)$$
where $x=\cos s$, called Legendre equation.
\end{proposition}

\begin{proof}
Let us first do a number of manipulations on our equation, before making the change of variables. By multiplying by $\sigma$, our equation becomes:
$$\sin s\cdot\frac{d}{ds}\left(\sin s\cdot\frac{d\sigma}{ds}\right)=\left(m^2-K\sin^2s\right)\sigma$$

By differentiating at left, this equation becomes:
$$\sin s\left(\cos s\cdot\sigma'+\sin s\cdot\sigma''\right)=\left(m^2-K\sin^2s\right)\sigma$$

Finally, by dividing everything by $\sin^2s$, our equation becomes:
$$\sigma''+\frac{\cos s}{\sin s}\cdot\sigma'=\left(\frac{m^2}{\sin^2s}-K\right)\sigma$$

Now let us set $\sigma(s)=f(\cos s)$. With this change of variables, we have:
$$\sigma=f(\cos s)$$
$$\sigma'=-\sin s\cdot f'(\cos s)$$
$$\sigma''=-\cos s\cdot f'(\cos s)+\sin^2s\cdot f''(\cos s)$$

By plugging this data, our radial equation becomes:
$$\sin^2s\cdot f''(\cos s)-2\cos s\cdot f'(\cos s)=\left(\frac{m^2}{\sin^2s}-K\right)f(\cos s)$$

Now with $x=\cos s$, which is our new variable, this equation reads:
$$(1-x^2)f''(x)-2xf'(x)=\left(\frac{m^2}{1-x^2}-K\right)f(x)$$

But this is the Legendre equation, as stated.
\end{proof}

Here comes now the difficult point. We have the following non-trivial result:

\begin{theorem}
The solutions of the Legendre equation, namely
$$(1-x^2)f''(x)-2xf'(x)=\left(\frac{m^2}{1-x^2}-K\right)f(x)$$
can be explicitely computed, via complicated math, and only the case
$$K=l(l+1)\quad:\quad l\in\mathbb N$$
is acceptable, on physical grounds. 
\end{theorem}

\begin{proof}
The first part is something quite complicated, involving the hypergeometric functions $\!\!{\ }_2F_1$, that you don't want to hear about, believe me. As for the second part, analysis and physics, this is something not trivial either. See Griffiths \cite{gr2}.
\end{proof}

In order to construct the solutions, we will need:

\begin{theorem}
The orthonormal basis of $L^2[-1,1]$ obtained by starting with the Weierstrass basis $\{x^l\}$, and doing Gram-Schmidt, is the family of polynomials $\{P_l\}$, with each $P_l$ being of degree $l$, and with positive leading coefficient, subject to:
$$\int_{-1}^1P_k(x)P_l(x)\,dx=\delta_{kl}$$
These polynomials, called Legendre polynomials, satisfy the equation
$$(1-x^2)P_l''(x)-2xP_l'(x)+l(l+1)P_l(x)=0$$
which is the Legendre equation at $m=0$, and with $K=l(l+1)$. Moreover,
$$P_l(x)=\frac{1}{2^ll!}\left(\frac{d}{dx}\right)^l(x^2-1)^l$$
which is called the Rodrigues formula for Legendre polynomials.
\end{theorem}

\begin{proof}
As a first observation, we are not lost somewhere in abstract math, because of the occurrence of the Legendre equation. As for the proof, this goes as follows:

\medskip

(1) The first assertion is clear, because the Gram-Schmidt procedure applied to the Weierstrass basis $\{x^l\}$ can only lead to a certain family of polynomials $\{P_l\}$, with each $P_l$ being of degree $l$, and also unique, if we assume that it has positive leading coefficient, with this $\pm$ choice being needed, as usual, at each step of Gram-Schmidt.

\medskip

(2) In order to have now an idea about these beasts, here are the first few of them, which can be obtained say via a straightforward application of Gram-Schmidt:
\begin{eqnarray*}
P_0&=&1\\
P_1&=&x\\
P_2&=&(3x^2-1)/2\\
P_3&=&(5x^3-3x)/2\\
P_4&=&(35x^4-30x^2+3)/8\\
P_5&=&(63x^5-70x^3+15x)/8
\end{eqnarray*}

(3) Now thinking about what Gram-Schmidt does, this is certainly something by recurrence. And examining the recurrence leads to the Legendre equation, as stated.

\medskip

(4) As for the Rodrigues formula, by uniqueness no need to try to understand where this formula comes from, and we have two choices here, either by verifying that $\{P_l\}$ is orthonormal, or by verifying the Legendre equation. And both methods work.
\end{proof}

Going ahead now, we can solve in fact the Legendre equation at any $m$, as follows:

\begin{proposition}
The general Legendre equation, with parameters $m\in\mathbb N$ and $K=l(l+1)$ with $l\in\mathbb N$, namely
$$(1-x^2)f''(x)-2xf'(x)=\left(\frac{m^2}{1-x^2}-l(l+1)\right)f(x)$$
is solved by the following functions, called Legendre functions,
$$P_l^m(x)=(-1)^m(1-x^2)^{m/2}\left(\frac{d}{dx}\right)^mP_l(x)$$
where $P_l$ are as before the Legendre polynomials. Also, we have
$$P_l^m(x)=(-1)^m\frac{(1-x^2)^{m/2}}{2^ll!}\left(\frac{d}{dx}\right)^{l+m}(x^2-1)^l$$
called Rodrigues formula for Legendre functions.
\end{proposition}

\begin{proof}
The first assertion is something elementary, coming by differentiating $m$ times the Legendre equation, which leads to the general Legendre equation. As for the second assertion, this follows from the Rodrigues formula for Legendre polynomials.
\end{proof}

And this is the end of our study. Eventually. By putting together all the above results, we are led to the following conclusion:

\index{separation of variables}
\index{Legendre function}
\index{spherical harmonics}

\begin{theorem}
The separated solutions $\alpha=\sigma(s)\theta(t)$ of the angular equation,
$$\sin s\cdot\frac{d}{ds}\left(\sin s\cdot\frac{d\alpha}{ds}\right)
+\frac{d^2\alpha}{dt^2}=-K\sin^2s\cdot\alpha$$
are given by the following formulae, where $l\in\mathbb N$ is such that $K=l(l+1)$,
$$\sigma(s)=P_l^m(\cos s)\quad,\quad\theta(t)=e^{imt}$$
and where $m\in\mathbb Z$ is a constant, and with $P_l^m$ being the Legendre function,
$$P_l^m(x)=(-1)^m(1-x^2)^{m/2}\left(\frac{d}{dx}\right)^mP_l(x)$$
where $P_l$ are the Legendre polynomials, given by the following formula:
$$P_l(x)=\frac{1}{2^ll!}\left(\frac{d}{dx}\right)^l(x^2-1)^l$$
These solutions $\alpha=\sigma(s)\theta(t)$ are called spherical harmonics.
\end{theorem}

\begin{proof}
This follows indeed from all the above, and with the comment that everything is taken up to linear combinations. We will normalize the wave function later.
\end{proof}

\section*{16d. The hydrogen atom}

Hydrogen, eventually. In order now to finish our study, and eventually get to conclusions about hydrogen, it remains to solve the radial equation, for the Coulomb potential $V$ of the proton. Let us begin with some generalities, valid for any time-independent potential $V$. As a first manipulation on the radial equation, we have:

\index{radial equation}

\begin{proposition}
The radial equation, written with $K=l(l+1)$, 
$$(r^2\rho')'-\frac{2mr^2}{h^2}(V-E)\rho=l(l+1)\rho$$
takes with $\rho=u/r$ the following form, called modified radial equation,
$$Eu=-\frac{h^2}{2m}\cdot u''+\left(V+\frac{h^2l(l+1)}{2mr^2}\right)u$$
which is a time-independent 1D Schr\"odinger equation.
\end{proposition}

\begin{proof}
With $\rho=u/r$ as in the statement, we have:
$$\rho=\frac{u}{r}\quad,\quad \rho'=\frac{u'r-u}{r^2}\quad,\quad (r^2\rho')'=u''r$$

By plugging this data into the radial equation, this becomes:
$$u''r-\frac{2mr}{h^2}(V-E)u=\frac{l(l+1)}{r}\cdot u$$

By multiplying everything by $h^2/(2mr)$, this latter equation becomes:
$$\frac{h^2}{2m}\cdot u''-(V-E)u=\frac{h^2l(l+1)}{2mr^2}\cdot u$$

But this gives the formula in the statement. As for the interpretation, as time-independent 1D Schr\"odinger equation, this is clear as well, and with the comment here that the term added to the potential $V$ is some sort of centrifugal term.
\end{proof}

Getting back now to the Coulomb potential of the proton, we have here:

\begin{fact}
The Coulomb potential of the hydrogen atom proton, acting on the electron by attraction, is given according to the Coulomb law by
$$V=-\frac{Kep}{r}$$
where $p$ is the charge of the proton, and $K$ is the Coulomb constant. In practice however we have $p\simeq e$ up to order $10^{-7}$, and so our formula can be written as
$$V\simeq-\frac{Ke^2}{r}$$
and we will use this latter formula, and with $=$ sign, for simplifying.
\end{fact} 

Getting back now to math, it remains to solve the modified radial equation, for the above potential $V$. And we have here the following result, which does not exactly solve this radial equation, but provides us instead with something far better, namely the proof of the original claim by Bohr, which was at the origin of everything:

\index{Bohr formula}

\begin{theorem}[Schr\"odinger]
In the case of the hydrogen atom, where $V$ is the Coulomb potential of the proton, the modified radial equation, which reads
$$Eu=-\frac{h^2}{2m}\cdot u''+\left(-\frac{Ke^2}{r}+\frac{h^2l(l+1)}{2mr^2}\right)u$$
leads to the Bohr formula for allowed energies,
$$E_n=-\frac{m}{2}\left(\frac{Ke^2}{h}\right)^2\cdot\frac{1}{n^2}$$
with $n\in\mathbb N$, the binding energy being 
$$E_1\simeq-2.177\times 10^{-18}$$
with means $E_1\simeq-13.591\ {\rm eV}$.
\end{theorem}

\begin{proof}
This is again something non-trivial, and we will be following Griffiths \cite{gr2}, with some details missing. The idea is as follows:

\medskip

(1) By dividing our modified radial equation by $E$, this becomes:
$$-\frac{h^2}{2mE}\cdot u''=\left(1+\frac{Ke^2}{Er}-\frac{h^2l(l+1)}{2mEr^2}\right)u$$

In terms of $\gamma=\sqrt{-2mE}/h$, this equation takes the following form:
$$\frac{u''}{\gamma^2}=\left(1+\frac{Ke^2}{Er}+\frac{l(l+1)}{(\gamma r)^2}\right)u$$

In terms of the new variable $p=\gamma r$, this latter equation reads:
$$u''=\left(1+\frac{\gamma Ke^2}{Ep}+\frac{l(l+1)}{p^2}\right)u$$

Now let us introduce a new constant $S$ for our problem, as follows:
$$S=-\frac{\gamma Ke^2}{E}$$

In terms of this new constant, our equation reads:
$$u''=\left(1-\frac{S}{p}+\frac{l(l+1)}{p^2}\right)u$$

(2) The idea will be that of looking for a solution written as a power series, but before that, we must ``peel off'' the asymptotic behavior. Which is something that can be done, of course, heuristically. With $p\to\infty$ we are led to $u''=u$, and ignoring the solution $u=e^p$ which blows up, our approximate asymptotic solution is:
$$u\sim e^{-p}$$

Similarly, with $p\to0$ we are led to $u''=l(l+1)u/p^2$, and ignoring the solution $u=p^{-l}$ which blows up, our approximate asymptotic solution is:
$$u\sim p^{l+1}$$

(3) The above heuristic considerations suggest writing our function $u$ as follows:
$$u=p^{l+1}e^{-p}v$$

So, let us do this. In terms of $v$, we have the following formula:
$$u'=p^le^{-p}\left[(l+1-p)v+pv'\right]$$

Differentiating a second time gives the following formula:
$$u''=p^le^{-p}\left[\left(\frac{l(l+1)}{p}-2l-2+p\right)v+2(l+1-p)v'+pv''\right]$$

Thus the radial equation, as modified in (1) above, reads:
$$pv''+2(l+1-p)v'+(S-2(l+1))v=0$$

(4) We will be looking for a solution $v$ appearing as a power series:
$$v=\sum_{j=0}^\infty c_jp^j$$

But our equation leads to the following recurrence formula for the coefficients:
$$c_{j+1}=\frac{2(j+l+1)-S}{(j+1)(j+2l+2)}\cdot c_j$$

(5) We are in principle done, but we still must check that, with this choice for the coefficients $c_j$, our solution $v$, or rather our solution $u$, does not blow up. And the whole point is here. Indeed, at $j>>0$ our recurrence formula reads, approximately:
$$c_{j+1}\simeq\frac{2c_j}{j}$$

But, surprisingly, this leads to $v\simeq c_0e^{2p}$, and so to $u\simeq c_0p^{l+1}e^p$, which blows up.

\medskip

(6) As a conclusion, the only possibility for $u$ not to blow up is that where the series defining $v$ terminates at some point. Thus, we must have for a certain index $j$:
$$2(j+l+1)=S$$

In other words, we must have, for a certain integer $n>l$:
$$S=2n$$

(7) We are almost there. Recall from (1) above that $S$ was defined as follows:
$$S=-\frac{\gamma Ke^2}{E}\quad:\quad\gamma=\frac{\sqrt{-2mE}}{h}$$

Thus, we have the following formula for the square of $S$:
$$S^2=\frac{\gamma^2K^2e^4}{E^2}=-\frac{2mE}{h^2}\cdot\frac{K^2e^4}{E^2}=-\frac{2mK^2e^4}{h^2E}$$

Now by using the formula $S=2n$ from (6), the energy $E$ must be of the form:
$$E=-\frac{2mK^2e^4}{h^2S^2}=-\frac{mK^2e^4}{2h^2n^2}$$

Calling this energy $E_n$, depending on $n\in\mathbb N$, we have, as claimed:
$$E_n=-\frac{m}{2}\left(\frac{Ke^2}{h}\right)^2\cdot\frac{1}{n^2}$$

(8) Thus, we proved the Bohr formula. Regarding numerics, the data is as follows:
$$K=8.988\times10^9\quad,\quad
e=1.602\times10^{-19}$$
$$h=1.055\times10^{-34}\quad,\quad 
m=9.109\times10^{-31}$$

But this gives the formula of $E_1$ in the statement.
\end{proof}

As a first remark, all this agrees with the Rydberg formula, due to:

\index{Rydberg formula}
\index{Planck formula}

\begin{theorem}
The Rydberg constant for hydrogen is given by
$$R=-\frac{E_1}{h_0c}$$
where $E_1$ is the Bohr binding energy, and the Rydberg formula itself, namely
$$\frac{1}{\lambda_{n_1n_2}}=R\left(\frac{1}{n_1^2}-\frac{1}{n_2^2}\right)$$
simply reads, via the energy formula in Theorem 16.25,
$$\frac{1}{\lambda_{n_1n_2}}=\frac{E_{n_2}-E_{n_1}}{h_0c}$$
which is in agreement with the Planck formula $E=h_0c/\lambda$.
\end{theorem}

\begin{proof}
Here the first assertion is something numeric, coming from the fact that the formula in the statement gives, when evaluated, the Rydberg constant:
$$R=\frac{-E_1}{h_0c}=\frac{2.177\times10^{-18}}{6.626\times10^{-34}\times 2.998\times 10^8}=1.096\times10^7$$

As a consequence, and passed now what the experiments exactly say, we can define the Rydberg constant of hydrogen abstractly, by the following formula:
$$R=\frac{m}{2h_0c}\left(\frac{Ke^2}{h}\right)^2$$

Regarding now the second assertion, by dividing $R=-E_1/(h_0c)$ by any number of type $n^2$ we obtain, according to the energy convention in Theorem 16.25:
$$\frac{R}{n^2}=-\frac{E_n}{h_0c}$$

But these are exactly the numbers which are subject to substraction in the Rydberg formula, and so we are led to the conclusion in the statement.
\end{proof}

Let us go back now to our study of the Schr\"odinger equation. Our conclusions are:

\index{spherical harmonics}

\begin{theorem}
The wave functions of the hydrogen atom are the following functions, labelled by three quantum numbers, $n,l,m$,
$$\phi_{nlm}(r,s,t)=\rho_{nl}(r)\alpha_l^m(s,t)$$
where $\rho_{nl}(r)=p^{l+1}e^{-p}v(p)/r$ with $p=\gamma r$ as before, with the coefficients of $v$ subject to
$$c_{j+1}=\frac{2(j+l+1-n)}{(j+1)(j+2l+2)}\cdot c_j$$
and $\alpha_l^m(s,t)$ being the spherical harmonics found before.
\end{theorem}

\begin{proof}
This follows indeed by putting together all the results obtained so far, and with the remark that everything is up to the normalization of the wave function.
\end{proof}

In what regards the main wave function, that of the ground state, we have:

\index{ground state}
\index{Bohr radius}

\begin{theorem}
With the hydrogen atom in its ground state, the wave function is
$$\phi_{100}(r,s,t)=\frac{1}{\sqrt{\pi a^3}}\,e^{-r/a}$$
where $a=1/\gamma$ is the inverse of the parameter appearing in our computations above,
$$\gamma=\frac{\sqrt{-2mE}}{h}$$
called Bohr radius of the hydrogen atom. This Bohr radius is the mean distance between the electron and the proton, in the ground state, and is given by the formula
$$a=\frac{h^2}{mKe^2}$$
which numerically means $a\simeq5.291\times10^{-11}$.
\end{theorem}

\begin{proof}
There are several things going on here, as follows:

\medskip

(1) According to the various formulae in the proof of Theorem 16.25, taken at $n=1$, the parameter $\gamma$ appearing in the computations there is given by:
$$\gamma=\frac{\sqrt{-2mE}}{h}=\frac{1}{h}\cdot m\cdot\frac{Ke^2}{h}=\frac{mKe^2}{h^2}$$

Thus, the inverse $a=1/\gamma$ is indeed given by the formula in the statement.

\medskip

(2) Regarding the wave function, according to Theorem 16.27 this consists of:
$$\rho_{10}(r)=\frac{2e^{-r/a}}{\sqrt{a^3}}\quad,\quad \alpha_0^0(s,t)=\frac{1}{2\sqrt{\pi}}$$

By making the product, we obtain the formula of $\phi_{100}$ in the statement.

\medskip

(3) But this formula of $\phi_{100}$ shows in particular that the Bohr radius $a$ is indeed the mean distance between the electron and the proton, in the ground state.

\medskip

(4) Finally, in what regards the numerics, these are as follows:
$$a=\frac{1.055^2\times10^{-68}}{9.109\times10^{-31}\times8.988\times10^9\times1.602^2\times10^{-38}}=5.297\times10^{-11}$$

Thus, we are led to the conclusions in the statement.
\end{proof}

Getting back now to the general setting of Theorem 16.25, the point is that the polynomials $v(p)$ appearing there are well-known objects in mathematics, as follows:

\begin{proposition}
The polynomials $v(p)$ are given by the formula
$$v(p)=L_{n-l-1}^{2l+1}(p)$$
where the polynomials on the right, called associated Laguerre polynomials, are given by
$$L_q^p(x)=(-1)^p\left(\frac{d}{dx}\right)^pL_{p+q}(x)$$
with $L_{p+q}$ being the Laguerre polynomials, given by the following formula:
$$L_q(x)=\frac{e^x}{q!}\left(\frac{d}{dx}\right)^q(e^{-x}x^q)$$
\end{proposition}

\begin{proof}
The story here is very similar to that of the Legendre polynomials. Consider the Hilbert space $H=L^2[0,\infty)$, with the following scalar product on it:
$$<f,g>=\int_0^\infty f(x)g(x)e^{-x}\,dx$$

(1) The orthogonal basis obtained by applying Gram-Schmidt to the Weierstrass basis $\{x^q\}$ is then the basis formed by the Laguerre polynomials $\{L_q\}$.

\medskip

(2) We have the explicit formula for $L_q$ in the statement, which is analogous to the Rodrigues formula for the Legendre polynomials.

\medskip

(3) The first assertion follows from the fact that the coefficients of the associated Laguerre polynomials satisfy the equation for the coefficients of $v(p)$.

\medskip

(4) Alternatively, the first assertion follows as well by using an equation for the Laguerre polynomials, which is very similar to the Legendre equation.
\end{proof}

With the above result in hand, we can now improve Theorem 16.25, as follows:

\index{spherical harmonics}

\begin{theorem}
The wave functions of the hydrogen atom are given by
$$\phi_{nlm}(r,s,t)=\sqrt{\left(\frac{2}{na}\right)^3\frac{(n-l-1)!}{2n(n+l)!}}e^{-r/na}\left(\frac{2r}{na}\right)^lL_{n-l-1}^{2l+1}\left(\frac{2r}{na}\right)\alpha_l^m(s,t)$$
with $\alpha_l^m(s,t)$ being the spherical harmonics found before.
\end{theorem}

\begin{proof}
This follows indeed by putting together what we have, namely Theorem 16.25 and Proposition 16.29, and then doing some remaining work, concerning the normalization of the wave function, which leads to the normalization factor appearing above.
\end{proof}

And good news, that is all. The above formula is all you need, in everyday life.

\section*{16e. Exercises}

Congratulations for having read this book, and no exercises for this final chapter. However, if you enjoyed this book, and looking for more to read, have a look at the various books referenced below. Normally these are all good books, which all sort of go to the point, without bothering much with annoying details, and normally you should like them too. So, have a look at them, and start with the one that you like the most.

\bigskip

Before everything, however, learn more calculus, and especially do many exercises, as many as needed, first in order to be at ease with calculus, and then in order to really love calculus. In fact, and in the hope that you got it, the present book was just an introduction to calculus, not calculus itself. Just a beginning.

\baselineskip=14pt

\printindex

\end{document}